%% file: thesis.tex
  \documentclass[utopia]{ucalgarythesis}

\input{preamble.tex}

\begin{document}



  \title{
Dagger Linear Logic and Categorical Quantum Mechanics   \\ \bigskip
   }
   
  \author{Priyaa Varshinee Srinivasan}
  \thesis{Thesis}
  \dept{Graduate Program in Computer Science}
  \degree{Doctor of Philosophy}
  \gradyear{2021}
  \monthname{September}

  \frontmatter           
  \makethesistitle       



\begin{thesisabstract}  
 
  This thesis develops the categorical proof theory for the non-compact multiplicative dagger ($\dagger$) linear logic, 
  and investigates its applications to Categorical Quantum Mechanics (CQM). 
   
  The existing frameworks of CQM are categorical proof theories of compact $\dagger$-linear logic, and
  are motivated by the interpretation of quantum systems in the category of finite dimensional Hilbert spaces.
  This thesis describes a new non-compact framework called Mixed Unitary Categories which can 
  accommodate infinite dimensional systems, and develops models for the framework. To this end, it builds 
  on linearly distributive categories, and $*$-autonomous categories which are categorical proof theories 
  of (non-compact) multiplicative linear logic. The proof theory of non-compact $\dagger$-linear logic is 
  obtained from the basic setting of an LDC by adding a dagger functor satisfying appropriate coherences 
  to give a $\dagger$-LDC.
   
  From every (isomix) $\dagger$-LDC one can extract a canonical ``unitary core'' which up to equivalence is the traditional 
  CQM framework of $\dagger$-monoidal categories. This leads to the framework of Mixed Unitary Categories (MUCs):  
  every MUC contains a (compact) unitary core which is extended by a (non-compact) isomix $\dagger$-LDC. 
  Various models of MUCs based on Finiteness Spaces, Chu spaces, Hopf modules, etc., are developed in this thesis. 
  This thesis also generalizes the key algebraic structures of CQM, such as observables, measurement, and 
  complementarity, to MUC framework. Furthermore, using the MUC framework, this thesis establishes a 
  connection between the complementary observables of quantum mechanics and the exponential modalities of linear logic. 
   
\end{thesisabstract}



 
 


  \chapter{Acknowledgements}  
  
  I thank my supervisor Dr. Robin Cockett for the trust he placed in me even though 
  I was not a mathematician and had no clue about Category Theory when I joined his group. 
  I am eternally grateful for his unwavering support and patience in helping me refine my {\em art of learning}, 
  and, dually, {\em my art of teaching}. I am also grateful to him for providing an inclusive, safe and a comfortable
  environment for me to learn Mathematics and to do research. Thank you, Robin, 
  for all the proofs we worked out over a cup of tea and the group hikes.

  I am thankful to Jonathan, my friend, for the countless hours he spent in teaching 
  me Category Theory and helping me write proofs  in the initial days of my Ph.D., without 
  which this journey would not have been possible. Thank you, Jonathan!
  
 I thank my co-supervisor Dr. Gilad Gour for his patience and kind support 
 throughout my degree. I thank Dr. Kristine Bauer for her support as my interim supervisor 
 during Dr. Cockett's sabbatical and for all the fun conversations 
 we shared during the weekly seminars and symposiums. 
 I thank Dr. Renate Scheidler for her excellent mentorship during 
 the practicum of my SAGES teaching certification. I thank Dr. Jalal Kawash 
 for the kind guidance I received as his teaching assistant.  I thank Dr. Barry Sanders 
 for the financial support he provided me during the first three years of my degree. 
  I thank Dr. Bob Coecke, Dr. Jamie Vicary,  and Dr. Chris Heunen for hosting me at  the 
  University of Oxford and at the University of Edinburgh during my Mitacs Globallink internship. 
 I thank my examiners Dr. Peter Selinger, Dr. Philipp Woelfel, Dr. Carlo Maria Scandolo, and
 Dr. Kristine Bauer immensely for their time and effort in reading my thesis, and their feedback. 
 
 I thank my friend and colleague Jean Simon Lemay (JS) for many examples 
 in this thesis, and for making my research visit to Oxford quite enjoyable. Thank you, JS, 
 for lending me your bicycle during the visit! Thank you, Jonathan, Amolak, JS, and Ben for proof-reading 
 my thesis. Thank you, Cole, for the fun and successful collaborations on the CNOT and the dagger linear logic papers.
 A big thanks to all my colleagues, in the past and in the present, from Dr. Cockett's and Dr. Gour's 
 group for all the enlightening conversations and support. It was an honour for me to work with you all. 

A ton of thanks to my friends, for making me feel home in Calgary and for all the beautiful memories - 
 Sunday afternoons at coffee shops, night walks at the Confederation park, dumpling parties, 
 road trip to Okanagan, hiking in the Rockies, hot tea, home-made delicacies, shared lunch, 
 weeding bees, and many more. I thank all my overseas friends for being there for me over Skype 
 and other mediums whenever I was in need of you. 

 Many thanks to Maryam, Jackie, and Katie at the front office of the Computer Science Department  
 for providing excellent administrative support throughout my degree, and for their gentle smile 
 whenever I turned up at their desk with a request. 

Last, but most of all, I thank my family for their extra-ordinary support, 
and encouragement in all my endeavors, and for making possible this 
first ever doctorate degree in my family tree.
    

  \chapter[Dedication]{}
  
  \begin{dedication}
     {\em To my parents Shanthi and Srinivasan, 
     
     \vspace{0.5em}
     and to all the front-line workers in this world and in heaven (COVID-19).}
  \end{dedication}

 \chapter*{Symbols, abbreviations and nomenclature}

\medskip

{ \def\arraystretch{1.35}
 \begin{tabular*}{\textwidth}{c @{\extracolsep{\fill}} lllll}
  Symbol & Definition \\
  $\cup$, $\bigcup$ & Set union \\
  $\sqcup$ & Disjoint union of sets \\
  $\cap$, $\bigcap$ & Set intersection \\
  $\in$ &  A member of \\
  $\llbracket x_1, x_2, \cdots, x_n \rrbracket$ & finite multiset \\
  $\emptyset$ & Empty set \\
  $\{ * \}$ & Singleton set \\
  $\{ x \mid $ Condition on $x \}$ & defining a set by a condition \\
  $\X, \Y$ & Categories except for $\C$ and $\R$ which sometimes mean  \\ 
                 & complex numbers, and real numbers respectively \\
  $A, B, C, \cdots$ & Objects in a category \\
  $A \simeq B$ & Object $A$ is isomorphic to object $B$ \\
  $f, g, h$ & Maps in a category \\
  $F, G $ & Functors between categories \\
  $\ox$ & Tensor product \\ 
  $\oa$ & Par product (in linear logic)
 \end{tabular*}}

\medskip 

{   \def\arraystretch{1.35} 
\begin{tabular*}{\textwidth}{c @{\extracolsep{\fill}} lllll}
  Symbol & Definition \\
  $\prod$ & Product \\
  $\coprod$ & Coproduct \\
  $a_\ox$ & Associator natural isomorphism  for the tensor product \\ 
  $a_\oa$ & Associator natural isomorphism for the par product \\
  $u_\ox^l $ & Left unitor natural isomorphism for the tensor product \\
  $\top$ & Unit of the tensor product \\
  $\bot$ & Unit of the par product \\ 
  $\partial^L$ & Left linear distributor \\ 
  $\partial^R$ & Right linear distributor \\ 
  $m$, $\nabla$ & Multiplication for a monoid \\
  $u$, $\triunit{0.65}$ & Unit for a monoid \\
  $d$, $\Delta$ & Comultiplication for a comonoid \\
  $e$, $\tricounit{0.65}$ & Counit for a comonoid \\
  $\dagger$, $\ddagger$ & Dagger functor \\
  $(-)^*$ & Dual functor \\ 
  $\overline{(-)}$ & Conjugation functor \\ 
  $:=$ & The left-hand side is defined to the right-hand side of the equation \\
  $\implies$ & Implication \\
  $\Leftrightarrow$ & If and only if \\
  $\exists$ & there exists \\
  $\forall $ & For all \\ 
  $\N$ & Natural numbers  \\
  $\Z$ & Integers 
  \end{tabular*} }


  \begin{singlespace}   

  \renewcommand\contentsname{Table of Contents}
  \cleardoublepage\phantomsection
  \addcontentsline{toc}{chapter}{\contentsname}
  \tableofcontents

  \renewcommand{\listfigurename}{List of Figures and Illustrations}
  \cleardoublepage\phantomsection
  \addcontentsline{toc}{chapter}{\listfigurename}
  \listoffigures

  

    %
    %
    
  
  \end{singlespace}     
  
%

  \chapter{Epigraph}

  \begin{center}
  \includegraphics[scale=0.55]{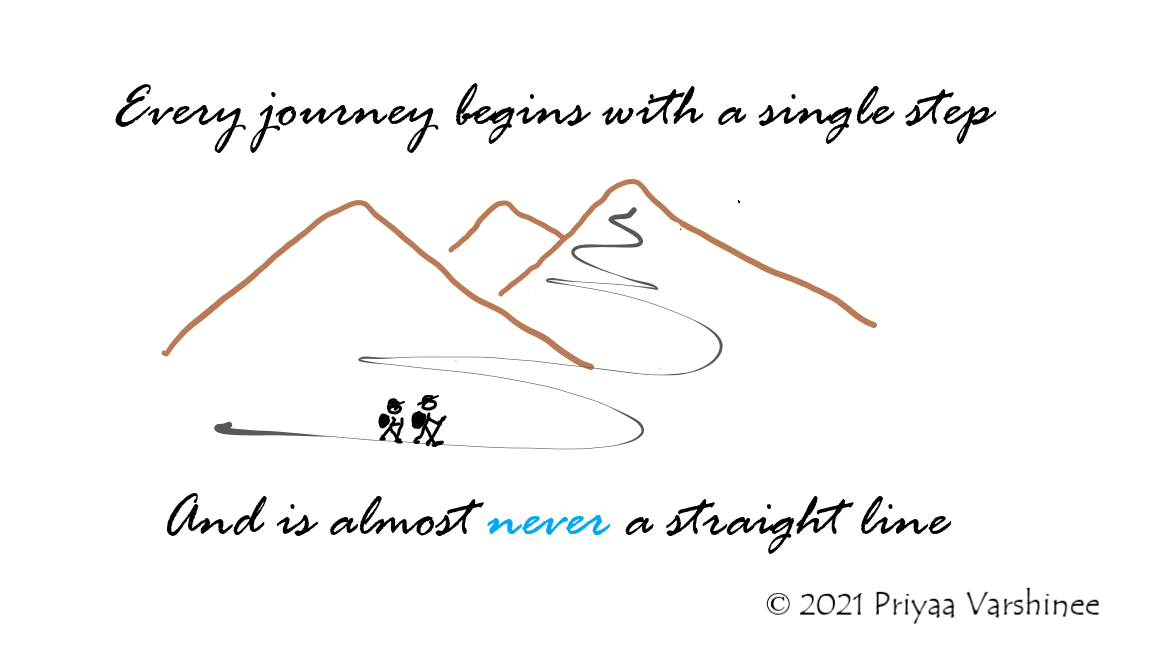}
  \end{center}

   


  \mainmatter           
    

\include{introduction}

\part{Dagger Linear Logic}
\label{Part: dll}
\include{chapter1}   
\include{chapter3} 	
\include{chapter4}  
\include{chapter5}

\part{Application of dagger linear logic to categorical quantum mechanics} 
\include{chapter2}  
\include{chapter-CP} 
\include{chapter-compaction} 
\include{chapter-compli}
\include{chapter-conclusion}



  \cleardoublepage\phantomsection
  \addcontentsline{toc}{chapter}{Bibliography}
  


    
   
 

   
\printbibliography 


  \appendix             

\end{document}

%% file: preamble.tex
\usepackage[pdftex,dvipsnames]{xcolor}  
\usepackage[utf8]{inputenc}
\usepackage[colorlinks, allcolors=blue]{hyperref} 
\usepackage{enumerate}
\usepackage{graphicx}
 
\usepackage{pdflscape}

\usepackage{bussproofs}
\usepackage{placeins}

\usepackage{amsmath,amsthm,amsfonts,graphicx}
\usepackage{amssymb}
\usepackage{stmaryrd}
\usepackage{tikz}
\usepackage{proof}
\usepackage{enumerate}
\usepackage{mathtools}
\usepackage{xypic}
\usepackage{lscape}
\usepackage{multicol}
\usepackage{leftidx}
\usepackage{bbm}
\usepackage{rotating}
\usepackage{extarrows}
\usepackage{cmll}

\usepackage{multirow}

\usepackage[english]{babel}
\usepackage[style=numeric, sorting=nty]{biblatex}
\usepackage{float}

\usepackage{array}   
\newcolumntype{L}{>{$}l<{$}} 

  \theoremstyle{plain}
  \newtheorem{theorem}{Theorem}[chapter]
  \newtheorem{lemma}[theorem]{Lemma}
  \newtheorem{corollary}[theorem]{Corollary}
  \newtheorem{proposition}[theorem]{Proposition}
  \newtheorem{remark}[theorem]{Remark}
  \newtheorem{example}[theorem]{Example}

  \newtheorem{definition}[theorem]{Definition}
   
  \newcommand{\invamalg}{\mathbin{\rotatebox[origin=c]{180}{$\amalg$}}}

  \newcommand{\s}{{\sf s}}
  \renewcommand{\t}{{\sf t}}
  \renewcommand{\u}{{\sf{u}}}
  \renewcommand{\v}{{\sf{v}}}

  \newcommand{\X}{\mathbb{X}}

  \newcommand{\C}{\mathbb{C}}
  \newcommand{\D}{\mathbb{D}}

  \newcommand{\N}{\mathbb{N}}
  \newcommand{\U}{\mathbb{U}}
  \newcommand{\V}{\mathbb{V}}
  \newcommand{\R}{\mathbb{R}}
  \newcommand{\Z}{\mathbb{Z}}
  \newcommand{\Y}{\mathbb{Y}}

  \newcommand{\m}{{\sf m}}

  \newcommand{\CP}{\mathsf{CP}}
  \newcommand{\ox}{\otimes}
  \newcommand{\pr}{\oplus}
  \newcommand{\oa}{\oplus}
  \newcommand{\op}{\mathsf{op}}
  \newcommand{\rev}{\mathsf{rev}}
  \newcommand{\mx}{\mathsf{mx}}
  \newcommand{\Mx}{\mathsf{Mx}}
  \newcommand{\Chu}{\mathsf{Chu}}
  \newcommand{\Chus}{\mathsf{Chus}}
  \newcommand{\FRel}{\mathsf{FRel}}
  \newcommand{\FMat}{\mathsf{FMat}}
  \newcommand{\Rel}{\mathsf{Rel}}
  \newcommand{\Mat}{\mathsf{Mat}}
  \newcommand{\Hilb}{\mathsf{Hilb}}
  \newcommand{\FHilb}{\mathsf{FHilb}}
  \newcommand{\Core}{\mathsf{Core}}
  \newcommand{\Unitary}{\mathsf{Unitary}}

  \newcommand{\lollipop}{\ensuremath{\!-\!\!\circ}}
  \renewcommand{\bar}[1]{\overline{#1}}
  \newcommand{\x}{\times}
  \newcommand{\poppilol} {\reflectbox{$\multimap$}}
  
  \newcommand{\dashvv}{\dashv \!\!\!\!\! \dashv}

  \newcommand{\Asp}{\mathsf{Asp}}

\DeclarePairedDelimiter\ceil{\lceil}{\rceil}
\DeclarePairedDelimiter\floor{\lfloor}{\rfloor}

\makeatletter

\newdimen\w@dth

\def\setw@dth#1#2{\setbox\z@\hbox{\scriptsize $#1$}\w@dth=\wd\z@
\setbox\@ne\hbox{\scriptsize $#2$}\ifnum\w@dth<\wd\@ne \w@dth=\wd\@ne \fi
\advance\w@dth by 1.2em}

\def\t@^#1_#2{\allowbreak\def\n@one{#1}\def\n@two{#2}\mathrel
{\setw@dth{#1}{#2}
\mathop{\hbox to \w@dth{\rightarrowfill}}\limits
\ifx\n@one\empty\else ^{\box\z@}\fi
\ifx\n@two\empty\else _{\box\@ne}\fi}}
\def\t@@^#1{\@ifnextchar_ {\t@^{#1}}{\t@^{#1}_{}}}

\def\t@left^#1_#2{\def\n@one{#1}\def\n@two{#2}\mathrel{\setw@dth{#1}{#2}
\mathop{\hbox to \w@dth{\leftarrowfill}}\limits
\ifx\n@one\empty\else ^{\box\z@}\fi
\ifx\n@two\empty\else _{\box\@ne}\fi}}
\def\t@@left^#1{\@ifnextchar_ {\t@left^{#1}}{\t@left^{#1}_{}}}

\def\two@^#1_#2{\def\n@one{#1}\def\n@two{#2}\mathrel{\setw@dth{#1}{#2}
\mathop{\vcenter{\hbox to \w@dth{\rightarrowfill}\kern-1.7ex
                 \hbox to \w@dth{\rightarrowfill}}%
       }\limits
\ifx\n@one\empty\else ^{\box\z@}\fi
\ifx\n@two\empty\else _{\box\@ne}\fi}}
\def\tw@@^#1{\@ifnextchar_ {\two@^{#1}}{\two@^{#1}_{}}}

\def\tofr@^#1_#2{\def\n@one{#1}\def\n@two{#2}\mathrel{\setw@dth{#1}{#2}
\mathop{\vcenter{\hbox to \w@dth{\rightarrowfill}\kern-1.7ex
                 \hbox to \w@dth{\leftarrowfill}}%
       }\limits
\ifx\n@one\empty\else ^{\box\z@}\fi
\ifx\n@two\empty\else _{\box\@ne}\fi}}
\def\t@fr@^#1{\@ifnextchar_ {\tofr@^{#1}}{\tofr@^{#1}_{}}}

\newdimen\W@dth
\def\setW@dth#1#2{\setbox\z@\hbox{$#1$}\W@dth=\wd\z@
\setbox\@ne\hbox{$#2$}\ifnum\W@dth<\wd\@ne \W@dth=\wd\@ne \fi
\advance\W@dth by 1.2em}

\def\T@^#1_#2{\allowbreak\def\N@one{#1}\def\N@two{#2}\mathrel
{\setW@dth{#1}{#2}
\mathop{\hbox to \W@dth{\rightarrowfill}}\limits
\ifx\N@one\empty\else ^{\box\z@}\fi
\ifx\N@two\empty\else _{\box\@ne}\fi}}
\def\T@@^#1{\@ifnextchar_ {\T@^{#1}}{\T@^{#1}_{}}}

\def\T@left^#1_#2{\def\N@one{#1}\def\N@two{#2}\mathrel{\setW@dth{#1}{#2}
\mathop{\hbox to \W@dth{\leftarrowfill}}\limits
\ifx\N@one\empty\else ^{\box\z@}\fi
\ifx\N@two\empty\else _{\box\@ne}\fi}}
\def\T@@left^#1{\@ifnextchar_ {\T@left^{#1}}{\T@left^{#1}_{}}}

\def\Tofr@^#1_#2{\def\N@one{#1}\def\N@two{#2}\mathrel{\setW@dth{#1}{#2}
\mathop{\vcenter{\hbox to \W@dth{\rightarrowfill}\kern-1.7ex
                 \hbox to \W@dth{\leftarrowfill}}%
       }\limits
\ifx\N@one\empty\else ^{\box\z@}\fi
\ifx\N@two\empty\else _{\box\@ne}\fi}}
\def\T@fr@^#1{\@ifnextchar_ {\Tofr@^{#1}}{\Tofr@^{#1}_{}}}

\def\Two@^#1_#2{\def\N@one{#1}\def\N@two{#2}\mathrel{\setW@dth{#1}{#2}
\mathop{\vcenter{\hbox to \W@dth{\rightarrowfill}\kern-1.7ex
                 \hbox to \W@dth{\rightarrowfill}}%
       }\limits
\ifx\N@one\empty\else ^{\box\z@}\fi
\ifx\N@two\empty\else _{\box\@ne}\fi}}
\def\Tw@@^#1{\@ifnextchar_ {\Two@^{#1}}{\Two@^{#1}_{}}}

\def\to{\@ifnextchar^ {\t@@}{\t@@^{}}}
\def\from{\@ifnextchar^ {\t@@left}{\t@@left^{}}}
\def\tofro{\@ifnextchar^ {\t@fr@}{\t@fr@^{}}}
\def\To{\@ifnextchar^ {\T@@}{\T@@^{}}}
\def\From{\@ifnextchar^ {\T@@left}{\T@@left^{}}}
\def\Two{\@ifnextchar^ {\Tw@@}{\Tw@@^{}}}
\def\Tofro{\@ifnextchar^ {\T@fr@}{\T@fr@^{}}}

\makeatother

\tikzstyle{strings}=[baseline={([yshift=-.5ex]current bounding box.center)}]

\tikzset{every picture/.append style={scale=.5}, transform shape, strings}

\tikzset{%
symbol/.style={%
draw=none,
every to/.append style={%
edge node={node [sloped, allow upside down, auto=false]{$#1$}}}
}
}

\usetikzlibrary{shapes.geometric}
\usetikzlibrary{patterns}
\usetikzlibrary{fit}
\usetikzlibrary{positioning}
\usetikzlibrary{calc}
\usetikzlibrary{arrows}
\usetikzlibrary{decorations.markings}
\usetikzlibrary{decorations.pathreplacing}
\usetikzlibrary{shapes}

\pgfdeclarelayer{nodelayer}
\pgfdeclarelayer{edgelayer}
\pgfsetlayers{edgelayer,nodelayer,main}

\tikzset{simple/.style={}}
\tikzset{nothing/.style={outer sep=-3.4pt}}
\tikzset{map/.style={draw,fill=white, rectangle}}
\tikzstyle{filled}=[-, fill=black]

\tikzset{dot/.style={thick, fill=black, circle, scale=1, inner sep = .05cm}}

\tikzset{oa/.style={draw, scale=0.9,minimum height=.1cm,circle,append after command={
[shorten >=\pgflinewidth, shorten <=\pgflinewidth,]
(\tikzlastnode.north) edge (\tikzlastnode.south)
(\tikzlastnode.east) edge (\tikzlastnode.west)
} } }

\tikzset{ox/.style={draw, scale=0.9,minimum height=.1cm,circle,append after command={
[shorten >=\pgflinewidth, shorten <=\pgflinewidth,]
(\tikzlastnode.north west) edge (\tikzlastnode.south east)
(\tikzlastnode.north east) edge (\tikzlastnode.south west) } } }

\tikzset{circ/.style={
shape=circle, inner sep=1pt, draw}}

\tikzstyle{none}=[inner sep=-1pt]
\tikzstyle{circle}=[shape=circle,draw]

\tikzstyle{onehalfcircle}=[shape=circle, scale=1.5, draw]
\tikzstyle{twocircle}=[shape=circle, scale=2, draw]
\tikzstyle{black}=[shape=circle, fill=black, draw]

\tikzset{wires/.style={}}

\tikzset{box/.style={inner sep=0pt, thick, draw=black, text height=1.5ex, text depth=.25ex, 
text centered, minimum height=3em, anchor=center}}

\newcommand{\envmap}{
\begin{tikzpicture}[scale=1.5]
\draw (0,0.25) -- (0,-0.05);
\draw (-0.15,-0.05) -- (0.15,-0.05);
\draw (-0.10,-0.1) -- (0.10,-0.1);
\draw (-0.05,-0.15) -- (0.05,-0.15);
\end{tikzpicture}
}

\newcommand{\anotherenvmap}[1]{
\begin{tikzpicture}[scale=#1]
	\begin{pgfonlayer}{nodelayer}
		\node [style=none] (0) at (-2, 0.7) {};
		\node [style=none] (1) at (-2, 0.5) {};
		\node [style=none] (2) at (-2.12, 0.5) {};
		\node [style=none] (3) at (-1.88, 0.5) {};
	\end{pgfonlayer}
	\begin{pgfonlayer}{edgelayer}
		\draw [style=none] (0.center) to (1.center);
		\draw [style=none, bend right=90, looseness=1.75] (2.center) to (3.center);
		\draw [style=none] (2.center) to (3.center);
	\end{pgfonlayer}
\end{tikzpicture}
}

\newcommand{\mulmap}[2]{
	\begin{tikzpicture}[scale={#1}]
		\begin{pgfonlayer}{nodelayer}
			\node [style=circle, scale=0.4, fill={#2}] (5) at (0.32, 0.25) {};
			\node [style=none] (6) at (0.07, 0.5) {};
			\node [style=none] (7) at (0.57, 0.5) {};
			\node [style=none] (8) at (0.32, 0) {};
			\node [style=none] (9) at (0.64, 0.5) {};
		\end{pgfonlayer}
		\begin{pgfonlayer}{edgelayer}
			\draw [style=none] (8.center) to (5);
			\draw [style=none, bend left, looseness=1.25] (5) to (6.center);
			\draw [style=none, bend right, looseness=1.25] (5) to (7.center);
		\end{pgfonlayer}
	\end{tikzpicture}	
}

\newcommand{\leftaction}[2]{
\begin{tikzpicture}[scale=#1]
	\begin{pgfonlayer}{nodelayer}
		\node [style=none] (0) at (0.5, 0.25) {};
		\node [style=none] (1) at (0.5, 0.75) {};
		\node [style=none] (2) at (0, 0.75) {};
		\node [style=none] (3) at (0.5, -0.75) {};
		\node [style=none] (4) at (0.5, 2) {};
		\node [style=none] (5) at (-0.75, 2) {};
	\end{pgfonlayer}
	\begin{pgfonlayer}{edgelayer}
		\draw[fill=#2] (0.center) -- (1.center) -- (2.center) -- (0.center);
		\draw (3.center) to (0.center);
		\draw (1.center) to (4.center);
		\draw [bend left, looseness=1.00] (2.center) to (5.center);
	\end{pgfonlayer}
\end{tikzpicture}
}

\newcommand{\rightaction}[2]{
	\begin{tikzpicture}[scale=#1, xscale=-1]
	\begin{pgfonlayer}{nodelayer}
	\node [style=none] (0) at (0.5, 0.25) {};
	\node [style=none] (1) at (0.5, 0.75) {};
	\node [style=none] (2) at (0, 0.75) {};
	\node [style=none] (3) at (0.5, -0.75) {};
	\node [style=none] (4) at (0.5, 2) {};
	\node [style=none] (5) at (-0.75, 2) {};
	\end{pgfonlayer}
	\begin{pgfonlayer}{edgelayer}
	\draw[fill=#2] (0.center) -- (1.center) -- (2.center) -- (0.center);
	\draw (3.center) to (0.center);
	\draw (1.center) to (4.center);
	\draw [bend left, looseness=1.00] (2.center) to (5.center);
	\end{pgfonlayer}
	\end{tikzpicture} }

\newcommand{\leftcoaction}[2]{
	\begin{tikzpicture}[scale=#1, yscale=-1]
	\begin{pgfonlayer}{nodelayer}
	\node [style=none] (0) at (0.5, 0.25) {};
	\node [style=none] (1) at (0.5, 0.75) {};
	\node [style=none] (2) at (0, 0.75) {};
	\node [style=none] (3) at (0.5, -0.75) {};
	\node [style=none] (4) at (0.5, 2) {};
	\node [style=none] (5) at (-0.75, 2) {};
	\end{pgfonlayer}
	\begin{pgfonlayer}{edgelayer}
	\draw[fill=#2] (0.center) -- (1.center) -- (2.center) -- (0.center);
	\draw (3.center) to (0.center);
	\draw (1.center) to (4.center);
	\draw [bend left, looseness=1.00] (2.center) to (5.center);
	\end{pgfonlayer}
	\end{tikzpicture} }

\newcommand{\rightcoaction}[2]{
	\begin{tikzpicture}[scale=#1, yscale=-1, xscale=-1]
	\begin{pgfonlayer}{nodelayer}
	\node [style=none] (0) at (0.5, 0.25) {};
	\node [style=none] (1) at (0.5, 0.75) {};
	\node [style=none] (2) at (0, 0.75) {};
	\node [style=none] (3) at (0.5, -0.75) {};
	\node [style=none] (4) at (0.5, 2) {};
	\node [style=none] (5) at (-0.75, 2) {};
	\end{pgfonlayer}
	\begin{pgfonlayer}{edgelayer}
	\draw[fill=#2] (0.center) -- (1.center) -- (2.center) -- (0.center);
	\draw (3.center) to (0.center);
	\draw (1.center) to (4.center);
	\draw [bend left, looseness=1.00] (2.center) to (5.center);
	\end{pgfonlayer}
	\end{tikzpicture}
}

\newcommand{\unitmap}[2]{
\begin{tikzpicture}[scale=#1]
	\begin{pgfonlayer}{nodelayer}
		\node [style=circle, scale=0.4, fill=#2] (0) at (0, 0) {};
		\node [style=none] (1) at (0, -0.4) {};
		\node [style=none] (4) at (0.13, 0) {};
	\end{pgfonlayer}
	\begin{pgfonlayer}{edgelayer}
		\draw [style=none] (0) to (1.center);
	\end{pgfonlayer}
\end{tikzpicture} }

\newcommand{\counitmap}[2]{
\begin{tikzpicture}[scale=#1, rotate=180]
	\begin{pgfonlayer}{nodelayer}
		\node [style=circle, scale=0.4, fill=#2] (0) at (0, 0) {};
		\node [style=none] (1) at (0, -0.4) {};
		\node [style=none] (4) at (0.13, 0) {};
	\end{pgfonlayer}
	\begin{pgfonlayer}{edgelayer}
		\draw [style=none] (0) to (1.center);
	\end{pgfonlayer}
\end{tikzpicture}
}

\newcommand{\bialgunitmap}[1]{
\begin{tikzpicture}[scale=#1]
	\begin{pgfonlayer}{nodelayer}
		\node [style=none] (0) at (-1, 2) {};
		\node [style=none] (1) at (-0.5, 2) {};
		\node [style=none] (2) at (-0.75, 1.75) {};
		\node [style=none] (3) at (-0.75, 1.25) {};
	\end{pgfonlayer}
	\begin{pgfonlayer}{edgelayer}
		\draw (0.center) to (1.center);
		\draw (1.center) to (2.center);
		\draw (2.center) to (0.center);
		\draw (2.center) to (3.center);
		\draw[fill=white] (0) -- (1) -- (2) -- (0);
	\end{pgfonlayer}
\end{tikzpicture}
}

\newcommand{\bialgcounitmap}[1]{
\begin{tikzpicture}[scale=#1]
	\begin{pgfonlayer}{nodelayer}
		\node [style=none] (0) at (-1, 1) {};
		\node [style=none] (1) at (-0.5, 1) {};
		\node [style=none] (2) at (-0.75, 1.25) {};
		\node [style=none] (3) at (-0.75, 1.75) {};
	\end{pgfonlayer}
	\begin{pgfonlayer}{edgelayer}
		\draw (0.center) to (1.center);
		\draw (1.center) to (2.center);
		\draw (2.center) to (0.center);
		\draw (3.center) to (2.center);
		\draw[fill=white] (0) -- (1) -- (2) -- (0);
	\end{pgfonlayer}
\end{tikzpicture}
}

\newcommand{\twincomul}[2] {
\begin{tikzpicture}[scale={#1}]
	\begin{pgfonlayer}{nodelayer}
		\node [style=circle, scale=0.6, fill=#2] (0) at (-2, 1) {};
		\node [style=circle, scale=0.6, fill=#2] (1) at (-1.25, 1) {};
		\node [style=none] (2) at (-2, 1.5) {};
		\node [style=none] (3) at (-1.25, 1.5) {};
		\node [style=none] (4) at (-1.75, 0.5) {};
		\node [style=none] (5) at (-0.75, 0.5) {};
		\node [style=none] (6) at (-2.5, 0.5) {};
		\node [style=none] (7) at (-1.5, 0.5) {};
	\end{pgfonlayer}
	\begin{pgfonlayer}{edgelayer}
		\draw [bend right=45] (1) to (4.center);
		\draw [bend left=45] (0) to (7.center);
		\draw [bend right=45] (0) to (6.center);
		\draw (0) to (2.center);
		\draw (3.center) to (1);
		\draw [bend left=45] (1) to (5.center);
	\end{pgfonlayer}
\end{tikzpicture}}

\newcommand{\twinmul}[2] {
\begin{tikzpicture}[scale={#1}, yscale=-1]
	\begin{pgfonlayer}{nodelayer}
		\node [style=circle, scale=0.6, fill=#2] (0) at (-2, 1) {};
		\node [style=circle, scale=0.6, fill=#2] (1) at (-1.25, 1) {};
		\node [style=none] (2) at (-2, 1.5) {};
		\node [style=none] (3) at (-1.25, 1.5) {};
		\node [style=none] (4) at (-1.75, 0.5) {};
		\node [style=none] (5) at (-0.75, 0.5) {};
		\node [style=none] (6) at (-2.5, 0.5) {};
		\node [style=none] (7) at (-1.5, 0.5) {};
	\end{pgfonlayer}
	\begin{pgfonlayer}{edgelayer}
		\draw [bend right=45] (1) to (4.center);
		\draw [bend left=45] (0) to (7.center);
		\draw [bend right=45] (0) to (6.center);
		\draw (0) to (2.center);
		\draw (3.center) to (1);
		\draw [bend left=45] (1) to (5.center);
	\end{pgfonlayer}
\end{tikzpicture}}

\newcommand{\twincounit}[2] {
\begin{tikzpicture}[scale={#1}]
	\begin{pgfonlayer}{nodelayer}
		\node [style=circle, scale=0.6, fill=#2] (0) at (-2, 0.5) {};
		\node [style=circle, scale=0.6, fill=#2] (1) at (-1.5, 0.5) {};
		\node [style=none] (2) at (-2, 1.25) {};
		\node [style=none] (3) at (-1.5, 1.25) {};
	\end{pgfonlayer}
	\begin{pgfonlayer}{edgelayer}
		\draw (0) to (2.center);
		\draw (3.center) to (1);
	\end{pgfonlayer}
\end{tikzpicture} }

\newcommand{\twinunit}[2] {
\begin{tikzpicture}[scale={#1}, yscale=-1]
	\begin{pgfonlayer}{nodelayer}
		\node [style=circle, scale=0.6, fill=#2] (0) at (-2, 0.5) {};
		\node [style=circle, scale=0.6, fill=#2] (1) at (-1.5, 0.5) {};
		\node [style=none] (2) at (-2, 1.25) {};
		\node [style=none] (3) at (-1.5, 1.25) {};
	\end{pgfonlayer}
	\begin{pgfonlayer}{edgelayer}
		\draw (0) to (2.center);
		\draw (3.center) to (1);
	\end{pgfonlayer}
\end{tikzpicture} }

\newcommand{\comulmap}[2]{
	\begin{tikzpicture}[scale={#1}]
		\begin{pgfonlayer}{nodelayer}
			\node [style=circle, scale=0.4, fill={#2}] (5) at (0.32, 0.25) {};
			\node [style=none] (6) at (0.07, 0) {};
			\node [style=none] (7) at (0.57, 0) {};
			\node [style=none] (8) at (0.32, 0.5) {};
			\node [style=none] (9) at (0.64, 0) {};
		\end{pgfonlayer}
		\begin{pgfonlayer}{edgelayer}
			\draw [style=none] (8.center) to (5);
			\draw [style=none, bend right, looseness=1.25] (5) to (6.center);
			\draw [style=none, bend left, looseness=1.25] (5) to (7.center);
		\end{pgfonlayer}
	\end{tikzpicture}
}

\newcommand{\trianglemult}[1]{
\begin{tikzpicture}[scale=#1]
	\begin{pgfonlayer}{nodelayer}
		\node [style=none] (0) at (-0.25, 3.5) {};
		\node [style=none] (1) at (-0.5, 3.75) {};
		\node [style=none] (2) at (0, 3.75) {};
		\node [style=none] (3) at (-0.25, 3) {};
		\node [style=none] (4) at (0.25, 4.25) {};
		\node [style=none] (5) at (-0.75, 4.25) {};
	\end{pgfonlayer}
	\begin{pgfonlayer}{edgelayer}
		\draw [bend right, looseness=1.00] (2.center) to (4.center);
		\draw (0.center) to (1.center);
		\draw (0.center) to (2.center);
		\draw (2.center) to (1.center);
		\draw [in=-90, out=165, looseness=0.75] (1.center) to (5.center);
		\draw (0.center) to (3.center);
	\end{pgfonlayer}
\end{tikzpicture} }

\newcommand{\trianglecomult}[1]{
\begin{tikzpicture}[scale=#1]
	\begin{pgfonlayer}{nodelayer}
		\node [style=none] (0) at (-0.25, 3.75) {};
		\node [style=none] (1) at (-0.5, 3.5) {};
		\node [style=none] (2) at (0, 3.5) {};
		\node [style=none] (3) at (-0.25, 4.25) {};
		\node [style=none] (4) at (0.25, 3) {};
		\node [style=none] (5) at (-0.75, 3) {};
	\end{pgfonlayer}
	\begin{pgfonlayer}{edgelayer}
		\draw [bend left, looseness=1.00] (2.center) to (4.center);
		\draw (0.center) to (1.center);
		\draw (0.center) to (2.center);
		\draw (2.center) to (1.center);
		\draw [in=90, out=-165, looseness=0.75] (1.center) to (5.center);
		\draw (0.center) to (3.center);
	\end{pgfonlayer}
\end{tikzpicture}
}

\newcommand{\trianglecomultblack}[1]{
\begin{tikzpicture}[scale=#1]
	\begin{pgfonlayer}{nodelayer}
		\node [style=none] (0) at (-0.25, 3.75) {};
		\node [style=none] (1) at (-0.5, 3.5) {};
		\node [style=none] (2) at (0, 3.5) {};
		\node [style=none] (3) at (-0.25, 4.25) {};
		\node [style=none] (4) at (0.25, 3) {};
		\node [style=none] (5) at (-0.75, 3) {};
	\end{pgfonlayer}
	\begin{pgfonlayer}{edgelayer}
		\draw [bend left, looseness=1.00] (2.center) to (4.center);
		\draw [fill=black] (0.center) -- (1.center) --  (2.center) -- (0.center);
		\draw [in=90, out=-165, looseness=0.75] (1.center) to (5.center);
		\draw (0.center) to (3.center);
	\end{pgfonlayer}
\end{tikzpicture}
}

\newcommand{\trianglecounit}[1]{
\begin{tikzpicture}[scale=#1]
	\begin{pgfonlayer}{nodelayer}
		\node [style=none] (0) at (-0.25, 3.5) {};
		\node [style=none] (1) at (-0.5, 3.25) {};
		\node [style=none] (2) at (0, 3.25) {};
		\node [style=none] (3) at (-0.25, 4.25) {};
		\node [style=none] (4) at (-0.25, 2.8) {};
	\end{pgfonlayer}
	\begin{pgfonlayer}{edgelayer}
		\draw (0.center) to (1.center);
		\draw (0.center) to (2.center);
		\draw (2.center) to (1.center);
		\draw (0.center) to (3.center);
	\end{pgfonlayer}
\end{tikzpicture}~\!\!}

\newcommand{\trianglecounitblack}[1]{
\begin{tikzpicture}[scale=#1]
	\begin{pgfonlayer}{nodelayer}
		\node [style=none] (0) at (-0.25, 3.5) {};
		\node [style=none] (1) at (-0.5, 3.25) {};
		\node [style=none] (2) at (0, 3.25) {};
		\node [style=none] (3) at (-0.25, 4.25) {};
		\node [style=none] (4) at (-0.25, 2.8) {};
	\end{pgfonlayer}
	\begin{pgfonlayer}{edgelayer}
		\draw [fill=black] (0.center) -- (1.center) --  (2.center) -- (0.center);
		\draw (0.center) to (3.center);
	\end{pgfonlayer}
\end{tikzpicture}~\!\!}

\newcommand{\triangleunit}[1]{
\begin{tikzpicture}[scale=#1]
	\begin{pgfonlayer}{nodelayer}
		\node [style=none] (0) at (-0.25, 4) {};
		\node [style=none] (1) at (-0.5, 4.25) {};
		\node [style=none] (2) at (0, 4.25) {};
		\node [style=none] (3) at (-0.25, 3.25) {};
		\node [style=none] (4) at (-0.25, 3) {};
	\end{pgfonlayer}
	\begin{pgfonlayer}{edgelayer}
		\draw (0.center) to (1.center);
		\draw (0.center) to (2.center);
		\draw (2.center) to (1.center);
		\draw (0.center) to (3.center);
	\end{pgfonlayer}
\end{tikzpicture}~\!\!}

\newcommand{\linmonwtik} {\begin{tikzpicture}
	\begin{pgfonlayer}{nodelayer}
		\node [style=none] (0) at (-2.7, 1.17) {};
		\node [style=none] (1) at (-1.85, 1.17) {};
		\node [style=none] (2) at (-2, 1.35) {};
		\node [style=none] (3) at (-2, 1) {};
		\node [style=none] (4) at (-1.85, 1) {};
		\node [style=none] (5) at (-1.85, 1.35) {};
		\node [style=circle, scale=0.6] (6) at (-2.35, 1.45) {};
		\node [style=none] (7) at (-1.6, 1.17) {};
		\node [style=none] (8) at (-2.95, 1.17) {};
	\end{pgfonlayer}
	\begin{pgfonlayer}{edgelayer}
		\draw (2.center) to (3.center);
		\draw (5.center) to (4.center);
		\draw (0.center) to (1.center);
	\end{pgfonlayer}
\end{tikzpicture}}

\newcommand{\expmonwtik} {\begin{tikzpicture}
	\begin{pgfonlayer}{nodelayer}
		\node [style=none] (0) at (-2.8, 1.17) {};
		\node [style=none] (1) at (-1.85, 1.17) {};
		\node [style=none] (2) at (-2, 1.35) {};
		\node [style=none] (3) at (-2, 1) {};
		\node [style=none] (4) at (-1.85, 1) {};
		\node [style=none] (5) at (-1.85, 1.35) {};
		\node [style=circle, scale=0.6] (6) at (-2.35, 1.55) {};
		\node [style=map, scale=1.7, fill opacity=0] (9) at (-2.35, 1.55) {};
		\node [style=none] (7) at (-1.6, 1.17) {};
		\node [style=none] (8) at (-3.05, 1.17) {};
	\end{pgfonlayer}
	\begin{pgfonlayer}{edgelayer}
		\draw (2.center) to (3.center);
		\draw (5.center) to (4.center);
		\draw (0.center) to (1.center);
	\end{pgfonlayer}
\end{tikzpicture}}

\newcommand{\dagmonwtik} {\begin{tikzpicture}
	\begin{pgfonlayer}{nodelayer}
		\node [style=none] (0) at (-2.7, 1.17) {};
		\node [style=none] (1) at (-1.85, 1.17) {};
		\node [style=none] (2) at (-2, 1.35) {};
		\node [style=none] (3) at (-2, 1) {};
		\node [style=none] (4) at (-1.85, 1) {};
		\node [style=none] (5) at (-1.85, 1.35) {};
		\node [style=circle, scale=0.6] (6) at (-2.25, 1.45) {};
		\node [style=none] (7) at (-1.6, 1.17) {};
		\node [style=none] (8) at (-2.95, 1.17) {};
		\node [style=none, scale=1.5] (9) at (-2.55, 1.45) {$\dag$};
	\end{pgfonlayer}
	\begin{pgfonlayer}{edgelayer}
		\draw (2.center) to (3.center);
		\draw (5.center) to (4.center);
		\draw (0.center) to (1.center);
	\end{pgfonlayer}
\end{tikzpicture}}

\newcommand{\lincomonwtik} {\begin{tikzpicture} 
	\begin{pgfonlayer}{nodelayer}
		\node [style=none] (0) at (-2.7, 1.17) {};
		\node [style=none] (1) at (-1.85, 1.17) {};
		\node [style=none] (2) at (-2, 1.35) {};
		\node [style=none] (3) at (-2, 1) {};
		\node [style=none] (4) at (-1.85, 1) {};
		\node [style=none] (5) at (-1.85, 1.35) {};
		\node [style=circle, scale=0.6] (6) at (-2.35, 0.9) {};
		\node [style=none] (7) at (-1.6, 1.17) {};
		\node [style=none] (8) at (-2.95, 1.17) {};
	\end{pgfonlayer}
	\begin{pgfonlayer}{edgelayer}
		\draw (2.center) to (3.center);
		\draw (5.center) to (4.center);
		\draw (0.center) to (1.center);
	\end{pgfonlayer}
\end{tikzpicture}}

\newcommand{\dagcomonwtik} {\begin{tikzpicture} 
	\begin{pgfonlayer}{nodelayer}
		\node [style=none] (0) at (-2.7, 1.17) {};
		\node [style=none] (1) at (-1.85, 1.17) {};
		\node [style=none] (2) at (-2, 1.35) {};
		\node [style=none] (3) at (-2, 1) {};
		\node [style=none] (4) at (-1.85, 1) {};
		\node [style=none] (5) at (-1.85, 1.35) {};
		\node [style=circle, scale=0.5] (6) at (-2.2, 0.9) {};
		\node [style=none, scale=1.4] (9) at (-2.5, 0.8) {$\dag$};
		\node [style=none] (7) at (-1.6, 1.17) {};
		\node [style=none] (8) at (-2.95, 1.17) {};
	\end{pgfonlayer}
	\begin{pgfonlayer}{edgelayer}
		\draw (2.center) to (3.center);
		\draw (5.center) to (4.center);
		\draw (0.center) to (1.center);
	\end{pgfonlayer}
\end{tikzpicture}}

\newcommand{\lincomonbtik} {\begin{tikzpicture} 
	\begin{pgfonlayer}{nodelayer}
		\node [style=none] (0) at (-2.7, 1.17) {};
		\node [style=none] (1) at (-1.85, 1.17) {};
		\node [style=none] (2) at (-2, 1.35) {};
		\node [style=none] (3) at (-2, 1) {};
		\node [style=none] (4) at (-1.85, 1) {};
		\node [style=none] (5) at (-1.85, 1.35) {};
		\node [style=circle, scale=0.6, fill=black] (6) at (-2.35, 0.9) {};
		\node [style=none] (7) at (-1.6, 1.17) {};
		\node [style=none] (8) at (-2.95, 1.17) {};
	\end{pgfonlayer}
	\begin{pgfonlayer}{edgelayer}
		\draw (2.center) to (3.center);
		\draw (5.center) to (4.center);
		\draw (0.center) to (1.center);
	\end{pgfonlayer}
\end{tikzpicture}}

\newcommand{\lincomonwtritik} {\begin{tikzpicture} 
	\begin{pgfonlayer}{nodelayer}
		\node [style=none] (0) at (-2.7, 1.17) {};
		\node [style=none] (1) at (-1.85, 1.17) {};
		\node [style=none] (2) at (-2, 1.35) {};
		\node [style=none] (3) at (-2, 1) {};
		\node [style=none] (4) at (-1.85, 1) {};
		\node [style=none] (5) at (-1.85, 1.35) {};
		\node [style=none] (6) at (-2.2, 1) {};
		\node [style=none] (7) at (-2.5, 1) {};
		\node [style=none] (8) at (-2.35, 0.82) {};
		\node [style=none] (9) at (-1.6, 1.17) {};
		\node [style=none] (10) at (-2.95, 1.17) {};
	\end{pgfonlayer}
	\begin{pgfonlayer}{edgelayer}
		\draw (2.center) to (3.center);
		\draw (5.center) to (4.center);
		\draw (0.center) to (1.center);
		\draw (6.center) -- (7.center) -- (8.center) -- (6.center);
	\end{pgfonlayer}
\end{tikzpicture}}

\newcommand{\dagcomonwtritik} {\begin{tikzpicture} 
	\begin{pgfonlayer}{nodelayer}
		\node [style=none] (0) at (-2.7, 1.17) {};
		\node [style=none] (1) at (-1.85, 1.17) {};
		\node [style=none] (2) at (-2, 1.35) {};
		\node [style=none] (3) at (-2, 1) {};
		\node [style=none] (4) at (-1.85, 1) {};
		\node [style=none] (5) at (-1.85, 1.35) {};
		\node [style=none] (6) at (-2.1, 0.85) {};
		\node [style=none] (7) at (-2.4, 0.85) {};
		\node [style=none] (8) at (-2.25, 1) {};
		\node [style=none] (9) at (-1.6, 1.17) {};
		\node [style=none] (10) at (-2.95, 1.17) {};
		\node [style=none, scale=1.5] (11) at (-2.6, 0.85) {$\dag$};
	\end{pgfonlayer}
	\begin{pgfonlayer}{edgelayer}
		\draw (2.center) to (3.center);
		\draw (5.center) to (4.center);
		\draw (0.center) to (1.center);
		\draw (6.center) -- (7.center) -- (8.center) -- (6.center);
	\end{pgfonlayer}
\end{tikzpicture}}

\newcommand{\lincomonbtritik} {\begin{tikzpicture} 
	\begin{pgfonlayer}{nodelayer}
		\node [style=none] (0) at (-2.7, 1.17) {};
		\node [style=none] (1) at (-1.85, 1.17) {};
		\node [style=none] (2) at (-2, 1.35) {};
		\node [style=none] (3) at (-2, 1) {};
		\node [style=none] (4) at (-1.85, 1) {};
		\node [style=none] (5) at (-1.85, 1.35) {};
		\node [style=none] (6) at (-2.2, 1) {};
		\node [style=none] (7) at (-2.5, 1) {};
		\node [style=none] (8) at (-2.35, 0.82) {};
		\node [style=none] (9) at (-1.6, 1.17) {};
		\node [style=none] (10) at (-2.95, 1.17) {};
	\end{pgfonlayer}
	\begin{pgfonlayer}{edgelayer}
		\draw (2.center) to (3.center);
		\draw (5.center) to (4.center);
		\draw (0.center) to (1.center);
		\draw[fill=black] (6.center) -- (7.center) -- (8.center) -- (6.center);
	\end{pgfonlayer}
\end{tikzpicture}}

\newcommand{\linbialgwtik} {\begin{tikzpicture}
	\begin{pgfonlayer}{nodelayer}
		\node [style=none] (0) at (-2.8, 1.17) {};
		\node [style=none] (1) at (-1.85, 1.17) {};
		\node [style=none] (2) at (-2, 1.35) {};
		\node [style=none] (3) at (-2, 1) {};
		\node [style=none] (4) at (-1.85, 1) {};
		\node [style=none] (5) at (-1.85, 1.35) {};
		\node [style=none] (6) at (-2.2, 1) {};
		\node [style=none] (7) at (-2.5, 1) {};
		\node [style=none] (8) at (-2.35, 0.82) {};
		\node [style=none] (9) at (-1.6, 1.17) {};
		\node [style=none] (10) at (-3.05, 1.17) {};
		\node [style=circle, scale=0.6] (11) at (-2.35, 1.45) {};
	\end{pgfonlayer}
	\begin{pgfonlayer}{edgelayer}
		\draw (2.center) to (3.center);
		\draw (5.center) to (4.center);
		\draw (0.center) to (1.center);
		\draw (6.center) -- (7.center) -- (8.center) -- (6.center);
	\end{pgfonlayer}
\end{tikzpicture}}

\newcommand{\expbialgwtik} {\begin{tikzpicture}
	\begin{pgfonlayer}{nodelayer}
		\node [style=none] (0) at (-2.7, 1.17) {};
		\node [style=none] (1) at (-1.85, 1.17) {};
		\node [style=none] (2) at (-2, 1.35) {};
		\node [style=none] (3) at (-2, 1) {};
		\node [style=none] (4) at (-1.85, 1) {};
		\node [style=none] (5) at (-1.85, 1.35) {};
		\node [style=none] (6) at (-2.2, 1) {};
		\node [style=none] (7) at (-2.5, 1) {};
		\node [style=none] (8) at (-2.35, 0.82) {};
		\node [style=none] (9) at (-1.6, 1.17) {};
		\node [style=none] (10) at (-2.95, 1.17) {};
		\node [style=circle, scale=0.6 ] (12) at (-2.35, 1.55) {};
		\node [style=map, scale=1.7, fill opacity=0] (11) at (-2.35, 1.55) {};
	\end{pgfonlayer}
	\begin{pgfonlayer}{edgelayer}
		\draw (2.center) to (3.center);
		\draw (5.center) to (4.center);
		\draw (0.center) to (1.center);
		\draw (6.center) -- (7.center) -- (8.center) -- (6.center);
	\end{pgfonlayer}
\end{tikzpicture}}

\newcommand{\dagbialgwtik} {\begin{tikzpicture}
	\begin{pgfonlayer}{nodelayer}
		\node [style=none] (0) at (-2.7, 1.17) {};
		\node [style=none] (1) at (-1.85, 1.17) {};
		\node [style=none] (2) at (-2, 1.35) {};
		\node [style=none] (3) at (-2, 1) {};
		\node [style=none] (4) at (-1.85, 1) {};
		\node [style=none] (5) at (-1.85, 1.35) {};
		\node [style=none] (6) at (-2.2, 1) {};
		\node [style=none] (7) at (-2.5, 1) {};
		\node [style=none] (8) at (-2.35, 0.82) {};
		\node [style=none] (9) at (-1.6, 1.17) {};
		\node [style=none] (10) at (-2.95, 1.17) {};
		\node [style=circle, scale=0.5] (11) at (-2.25, 1.45) {};		
		\node [style=none, scale=1.5] (12) at (-2.55, 1.45) {$\dagger$};
		\node [style=none] (13) at (-2.25, 0.5) {}; 
	\end{pgfonlayer}
	\begin{pgfonlayer}{edgelayer}
		\draw (2.center) to (3.center);
		\draw (5.center) to (4.center);
		\draw (0.center) to (1.center);
		\draw (6.center) -- (7.center) -- (8.center) -- (6.center);
	\end{pgfonlayer}
\end{tikzpicture}}

\newcommand{\linbialgbtik} {\begin{tikzpicture}
	\begin{pgfonlayer}{nodelayer}
		\node [style=none] (0) at (-2.7, 1.17) {};
		\node [style=none] (1) at (-1.85, 1.17) {};
		\node [style=none] (2) at (-2, 1.35) {};
		\node [style=none] (3) at (-2, 1) {};
		\node [style=none] (4) at (-1.85, 1) {};
		\node [style=none] (5) at (-1.85, 1.35) {};
		\node [style=none] (6) at (-2.2, 1) {};
		\node [style=none] (7) at (-2.5, 1) {};
		\node [style=none] (8) at (-2.35, 0.82) {};
		\node [style=none] (9) at (-1.6, 1.17) {};
		\node [style=none] (10) at (-2.95, 1.17) {};
		\node [style=circle, scale=0.6, fill=black] (11) at (-2.35, 1.45) {};
	\end{pgfonlayer}
	\begin{pgfonlayer}{edgelayer}
		\draw (2.center) to (3.center);
		\draw (5.center) to (4.center);
		\draw (0.center) to (1.center);
		\draw[fill=black] (6.center) -- (7.center) -- (8.center) -- (6.center);
	\end{pgfonlayer}
\end{tikzpicture}}

\newcommand{\expbialgbtik} {\begin{tikzpicture}
	\begin{pgfonlayer}{nodelayer}
		\node [style=none] (0) at (-2.7, 1.17) {};
		\node [style=none] (1) at (-1.85, 1.17) {};
		\node [style=none] (2) at (-2, 1.35) {};
		\node [style=none] (3) at (-2, 1) {};
		\node [style=none] (4) at (-1.85, 1) {};
		\node [style=none] (5) at (-1.85, 1.35) {};
		\node [style=none] (6) at (-2.2, 1) {};
		\node [style=none] (7) at (-2.5, 1) {};
		\node [style=none] (8) at (-2.35, 0.82) {};
		\node [style=none] (9) at (-1.6, 1.17) {};
		\node [style=none] (10) at (-2.95, 1.17) {};
		\node [style=circle, scale=0.6, fill=black ] (12) at (-2.35, 1.55) {};
		\node [style=map, scale=1.7, fill opacity=0] (11) at (-2.35, 1.55) {};
	\end{pgfonlayer}
	\begin{pgfonlayer}{edgelayer}
		\draw (2.center) to (3.center);
		\draw (5.center) to (4.center);
		\draw (0.center) to (1.center);
		\draw (6.center) -- (7.center) -- (8.center) -- (6.center);
	\end{pgfonlayer}
\end{tikzpicture}}

\newcommand{\tricomul}[1]{\trianglecomult{#1}}

\newcommand{\tricounit}[1]{\trianglecounit{#1}}
\newcommand{\triunit}[1]{\triangleunit{#1}}

\newcommand{\tricomulb}[1]{\trianglecomultblack{#1}}

\newcommand{\tricounitb}[1]{\trianglecounitblack{#1}}

\newcommand{\dagdual}{\xymatrixcolsep{4mm} \xymatrix{ \ar@{-||}[r]^{\dagger} & }}

\newcommand{\whitelin}{\xymatrixcolsep{3mm} \xymatrix{  \ar@{-||}[r]^{\circ} &  }}
\newcommand{\blacklin}{\xymatrixcolsep{3mm} \xymatrix{ \ar@{-||}[r]^{\bullet} & }}

\newcommand{\whitedag}[1]\dagmonw 
\newcommand{\blackdag}[1]\dagmonb 

\newcommand{\linmonw} {\xymatrixcolsep{4mm} \xymatrix{ \ar@{-||}[r]^{\circ} & }}
\newcommand{\linmonb} {\xymatrixcolsep{4mm} \xymatrix{  \ar@{-||}[r]^{\bullet} & }}
\newcommand{\dagmonw} {\xymatrixcolsep{4mm} \xymatrix{ \ar@{-||}[r]_{}^{\dagger\circ} & }}
\newcommand{\dagmonb} {\xymatrixcolsep{4mm} \xymatrix{  \ar@{-||}[r]_{}^{\dagger\bullet} & }}

\newcommand{\lincomonw} {\xymatrixcolsep{4mm} \xymatrix{ \ar@{-||}[r]_{\circ} & }}
\newcommand{\lincomonb} {\xymatrixcolsep{4mm} \xymatrix{  \ar@{-||}[r]_{\bullet} & }}
\newcommand{\dagcomonw} {\xymatrixcolsep{4mm} \xymatrix{ \ar@{-||}[r]_{\dagger\circ} & }}
\newcommand{\dagcomonb} {\xymatrixcolsep{4mm} \xymatrix{  \ar@{-||}[r]_{\dagger\bullet} & }}

\newcommand{\monoid}[1]{(#1, \mulmap{1.5}{white}, \unitmap{1.5}{white})}
\newcommand{\comonoid}[1]{(#1, \comulmap{1.5}{white}, \counitmap{1.5}{white})}
\newcommand{\comonoidb}[1]{(#1, \comulmap{1.5}{black}, \counitmap{1.5}{black})}
\newcommand{\Frob}[1]{(#1, \mulmap{1.5}{white}, \unitmap{1.5}{white}, \comulmap{1.5}{white}, \counitmap{1.5}{white})}
\newcommand{\bFrob}[1]{(#1, \mulmap{1.5}{black}, \unitmap{1.5}{black}, \comulmap{1.5}{black}, \counitmap{1.5}{black})}
\newcommand{\bialg}[1]{(#1, \mulmap{1.5}{white}, \unitmap{1.5}{white}, \comulmap{1.5}{black}, \counitmap{1.5}{black})}
\newcommand{\bialgb}[1]{(#1, \mulmap{1.5}{black}, \unitmap{1.5}{black}, \comulmap{1.5}{white}, \counitmap{1.5}{white})}

\newcommand{\tr}{\triangleright}
\newcommand{\tl}{\triangleleft}

\xymatrixrowsep{5mm}
\newdir{(>}{{}*!/-5pt/\dir{>}}  

\newcommand{\Lim}[1]{\underset{#1}{\sf Lim~}}
\newcommand{\Colim}[1]{\underset{#1}{\sf Colim~}}

\newcommand{\tri}[1]{
\begin{tikzpicture}[scale=#1]
	\begin{pgfonlayer}{nodelayer}
		\node [style=none] (0) at (2, 3.5) {};
		\node [style=none] (1) at (1.75, 3.25) {};
		\node [style=none] (2) at (2.25, 3.25) {};
	\end{pgfonlayer}
	\begin{pgfonlayer}{edgelayer}
		\draw (0.center) to (1.center);
		\draw (0.center) to (2.center);
		\draw (2.center) to (1.center);
	\end{pgfonlayer}
\end{tikzpicture} } 

\newcommand{\btri}[1]{
\begin{tikzpicture}[scale=#1]
	\begin{pgfonlayer}{nodelayer}
		\node [style=none] (0) at (2, 3.5) {};
		\node [style=none] (1) at (1.75, 3.25) {};
		\node [style=none] (2) at (2.25, 3.25) {};
	\end{pgfonlayer}
	\begin{pgfonlayer}{edgelayer}
		\draw [fill=black] (0.center) -- (1.center) -- (2.center) -- (0.center);
	\end{pgfonlayer}
\end{tikzpicture} } 

\newcommand{\linbialg}{\xymatrixcolsep{5mm} \xymatrix{ \ar@{-||}[r]^{\tri{0.75}} & }}
\newcommand{\linbialgblack}{\xymatrixcolsep{5mm} \xymatrix{ \ar@{-||}[r]^{\btri{0.75}} & }}
\newcommand{\dagbialg}{\xymatrixcolsep{6mm} \xymatrix{ \ar@{-||}[r]^{\dagger \tri{0.75}} & }}
\newcommand{\dagbialgb}{\xymatrixcolsep{6mm} \xymatrix{ \ar@{-||}[r]^{\dagger \btri{0.75}} & }}
 
\newlength{\llcfoo}



%% file: introduction.tex

\chapter{Introduction}

The programme of Categorical Quantum Mechanics (CQM) \cite{CoK17,HeV19} initiated by 
Abramsky and Coecke's seminal paper \cite{AC04} employs a graphical calculus 
to study quantum processes within the $\dag$-compact closed category ($\dag$-KCC)
of finite dimensional Hilbert spaces ($\FHilb$). From the perspective of logic, 
the graphical calculus is the proof theory of a compact fragment of multiplicative $\dagger$-linear logic \cite{Dun06}.
Thus, CQM provides a novel approach to quantum mechanics with its use of a graphical calculus backed by the 
rigor of Categorical Proof Theory.

 A well-known limitation of compact closed categories is that they model finite dimensional 
 Hilbert spaces, but they do not model infinite dimensional spaces \cite{Heunen16}.  
A categorical generalization of compact closed categories, in which infinite dimensional spaces 
 can be modelled are $*$-autonomous categories. 
 This can be taken a step further by generalizing to linearly 
 distributive categories (LDCs) in which the existence of dual objects is not assumed. 
 These linear settings come with a proof theory (for non-compact multiplicative linear logic) which is
 a graphical calculus. Thus, one does not abandon this attractive feature of CQM in these more general settings. In fact, 
 the graphical calculus of linearly distributive and $*$-autonomous categories subsume the 
 graphical calculus of compact closed categories. Thus, from a categorical perspective, an inviting direction to  
 accommodate infinite dimensional  quantum processes is to utilize $*$-autonomous and 
 linearly distributive categories as extensions of the compact closed settings used in CQM. The advantage 
of such a foundational approach is that it encompasses many of the earlier approaches and allows for a structural 
positioning of the different approaches.

 The aim of this thesis is to lay the foundations of generalization in this direction. Towards this aim, 
the first part of this thesis develops the categorical proof theory of non-compact $\dagger$-linear logic 
using linearly distributive and $*$-autonomous categories, and shows that one can always 
extract the usual settings for compact $\dagger$-linear logic from this new framework 
of $\dagger$-isomix categories. The first part also describes models for the new 
framework. The second part of this thesis explores the applicability of $\dagger$-isomix 
categories to CQM and studying quantum processes of arbitrary dimensions. 


\section{From linear logic to quantum mechanics}

\subsection{Linear logic}
Linear logic was introduced by Girard in \cite{Gir87}  as a logic of resources manipulation. Unlike classical 
logic which treats logical statements as truth values, linear logic treats logical statements as resources 
which cannot be duplicated or destroyed.  For example, consider the following statements, $p$ and $q$:
\begin{align*}
    p &:  \text{to spend a dollar} \\
    q &: \text{to buy an apple} 
\end{align*}

Then, in linear logic, the compound statement ``$p \Rightarrow q$" has the meaning that if a dollar is spent then an apple 
can be bought. This means that a person can either have {\em a dollar} or {\em an apple} at a given time but not both. 
The word ``linear" refers to this resource sensitivity of the 
logic: thus, a proof of a statement in linear logic may be regarded as a series of resource transformations.

The types (formulae) of linear logic can be defined inductively as follows (``$|$" is to be read as ``or''):
\begin{align*}
A :=&  ~p ~|~ A^{\perp} \\
&|~ A \ox  A ~|~ 1 ~|~ A \parr A ~|~ \bot \\ 
&|~ A \with A ~|~ \top ~|~ A  \oa A ~|~ 0 \\
&|~ !A ~|~ ?A     
\end{align*}

In the inductive definition of $A$, $p$ is an atomic formula such as 
``to spend a dollar" and $p^\perp$ is the negation of the formula which in this case will be ``to 
receive a dollar".  The {\bf negation} is an involution, hence, $(p^\perp)^\perp = p$.

The connectives $\ox$, $1$, $\parr$, and $\bot$ are called the {\bf multiplicatives}. A statement, 
$A \otimes B$ (read as A {\em tensor} B) 
allows the resources $A$ and $B$ to be available at the same time.  For example, consider the statement $r : $ {\em an orange}. Then, $p \Rightarrow q \otimes r$ refers to the fact that 
spending a dollar buys an apple and an orange at the same time. The formula $1$ is the multiplicative truth, hence, $A \otimes 1 = 1 \otimes A = A$. The connective, $\parr$, read as {\em par}, is the multiplicative disjunction and is dual to tensor. This means that,
$(A^\perp \otimes B^\perp)^\perp = A \parr B$. Smililary, $1$ and $\top$ are duals to one another $(1^\perp = \top)$.

The connectives $\with$, $\top$, $\oa$, and $0$ are called the {\bf additives}. The statement, $A \with B$ 
(read as $A$ {\em with} B), means that either $A$ or $B$ is available at a time.  For example, the statement, $p \Rightarrow q \with r$, refers to the fact that 
spending a dollar buys either an apple or an orange with the freedom to choose either. In computer science $\with$, 
represents non-determinism. The formula $0$ is the additive false: $ A \parr 0 = 0 \parr A = A$.
The additive disjunction $\parr$ (read as the {\em plus}) and $\with$ are dual to one another. Similary, 
 $0$ and $\bot$ are duals. 

Finally, the {\bf exponential operators} - $!$ read as {\em of course} or {\em bang} and $?$ read as {\em why not} or {\em whimper} - allow 
for the duplication and destruction of resources. For any resource $A$, $!A$ models an {\em unlimited} store from which the resource $A$ 
can be extracted 0 or more times. Morever, the storage itself can be duplicated or even destroyed.  The {\em why not} operator, $?$ encodes the notion of infinite demand, and is dual to the $!$, 
that is, $(?A)^\perp = !(A^\perp)$. 

Interested readers can refer to \cite{Gir95, Gir87} for more details on linear logic.

\subsection{Categorical semantics for linear logic}

Linear logic, being a logic of resources, emphasizes the structure of proofs rather than provability, that is, one is more 
interested to know how a statement can be proved, rather than, merely if the statement is provable. Proofs in sequent calculus 
often contain extraneous details due the sequential nature of the calculus. For example, in a sequent proof, in order to apply a set of rules 
to a sequent, one must choose an order for the application of the rules, even if the rules are independent (that is, the order in which 
the rules are applied does not affect the final result). 

In order to remove such unuseful information,  Girard \cite{Gir87} introduced proof-nets to 
represent proofs of linear logic, specifically for the multiplicative fragment without units 
i.e., only the $\ox$ and  $\parr$. Proof-nets are formalized as circuits. The study of algorithms to decide if a given circuit corresponds to a valid 
proof-net paved the way to the study of the categorical semantics for multiplicative linear logic (MLL). 
An overview of the categorical proof theories for different fragments of MLL is provided in the table below: 

  \begin{table}[h]
    \centering
	\begin{tabular}{ | l | l | }
    \hline
	 {\bf Linear logic fragment} & {\bf Categorical proof theory}  \\  
   \hline 
   \hline
	MLL & Linearly distributive categories  \cite{CS97} \\  
  \hline
	MLL with negation  & $*$-autonomous categories  \cite{Bar91b} \\ 
  \hline
	Compact MLL ($\ox = \parr$, $1 = \bot$) & Monoidal categories \cite{Kel64} \\ 
  \hline
	Compact MLL with negation & Compact closed categories \cite{Kel80} \\ 
  \hline
	\end{tabular}
  \caption{Categorical semantics for multiplicative linear logic}
  \label{Table: MLL} 
\end{table}

\FloatBarrier

The above listed semantics are sound and complete in the sense that  
there exists a one-to-one correspondence between the proof-nets of a fragment and the morphisms 
of its corresponding categorical setting.  The proof-nets provide a graphical calculus to these categorical settings, 
thereby, enabling a string diagrammatic presentation of the morphisms in these categories. 

The proof theory of compact MLL based on the graphical calculus 
of monoidal categories \cite{Sel10}  is used to derive an elegant description 
of quantum mechanics \cite{CoK17, HeV19} .

\subsection{Categorical quantum mechanics (CQM)}

Traditionally, Hilbert spaces \cite{NiC10} or more generally von Neumann Algebras \cite{Dix11} are used as the mathematical 
framework of quantum mechanics \cite{Neu18} . While these frameworks support detailed computation, they 
do not support an intuition for the problem: this leads to an approach described as ``shut up and calculate'' \cite{Mer89}. 
The programme of Categorical Quantum Mechanics (CQM) \cite{AC04, CoK17, HeV19}  emerged out of the desire 
to develop a more intuitive framework to represent and reason about quantum processes.

Linear logic captures \cite{Pra93, AbD04, AC04, Dun06, Bol19} the essence of quantum mechanics owing to its 
resource-sensitive character. In particular, linear logic does not allow duplication of an arbitrary type: in quantum mechanics, this is 
referred to as the {\em no-cloning theorem} \cite{NiC10} which states that it is impossible to duplicate an arbitrary quantum state. 
Motivated by this connection, CQM uses the compact multiplicative linear logic as the base framework for its purpose, 
and added the notion of `dagger' to this fragment, thus giving rise to compact $\dagger$-linear logic.  The `dagger' abstracts the notion of `adjoint' 
which is crucial to quantum mechanics: measurable properties of a quantum system are 
given by self-adjoint operators on a separable Hilbert space, that is, a Hilbert space with 
countable orthonormal basis.  

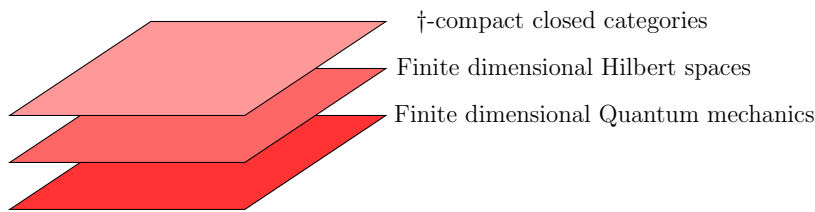
\begin{figure}[h]
\centering
	\begin{tikzpicture}[scale=1.25]
	\begin{pgfonlayer}{nodelayer}
		\node [style=none] (0) at (-4, -2) {};
		\node [style=none] (1) at (-1, 0) {};
		\node [style=none] (2) at (4, 0) {};
		\node [style=none] (3) at (1, -2) {};
		\node [style=none] (4) at (-4, -1) {};
		\node [style=none] (5) at (-1, 1) {};
		\node [style=none] (6) at (4, 1) {};
		\node [style=none] (7) at (1, -1) {};
		\node [style=none] (8) at (-4, 0) {};
		\node [style=none] (9) at (-1, 2) {};
		\node [style=none] (10) at (4, 2) {};
		\node [style=none] (11) at (1, 0) {};
		\node [style=none, scale=1.25] (12) at (7.75, 2) {$\dagger$-compact closed categories};
		\node [style=none, scale=1.25] (13) at (8, 1) { Finite dimensional Hilbert spaces};
		\node [style=none, scale=1.25] (14) at (8.65, 0) {Finite dimensional Quantum mechanics};
	\end{pgfonlayer}
	\begin{pgfonlayer}{edgelayer}
		\draw [style=filled, fill=red!80] (3.center)
			 to (0.center)
			 to (1.center)
			 to (2.center)
			 to cycle;
		\draw [style=filled, fill=red!60] (7.center)
			 to (4.center)
			 to (5.center)
			 to (6.center)
			 to cycle;
		\draw [style=filled, fill=red!40] (11.center)
			 to (8.center)
			 to (9.center)
			 to (10.center)
			 to cycle;
	\end{pgfonlayer}
\end{tikzpicture} 
\caption{$\dagger$-compact closed categories for quantum mechanics}
\label{Fig: layers}
\end{figure}

While monoidal categories provide the semantics for compact MLL, $\dagger$-monoidal categories provide the categorical semantics for 
compact $\dagger$-linear logic.  CQM uses $\dagger$-monoidal categories,  specifically, $\dagger$-compact closed 
categories (See Chapter \ref{Chap: CQM}) to develop a high-level, intuitive, formal language for quantum mechanics 
by abstracting the standard, traditionally-used, analytical framework of Hilbert spaces, as illustrated in Figure \ref{Fig: layers}.

In 2004, Abramsky and Coecke \cite{AC04} described the fundamental axioms of quantum mechanics 
within the framework of $\dagger$-compact closed categories ($\dagger$-KCCs). 
This was quite significant as it meant that the proof theory 
based on string diagrams of monoidal categories \cite{Sel10}, 
could be deployed to reason about quantum processes. For example, in CQM, 
physical systems are represented as wires and processes as circles.  
The label of a wire represents its type. Diagram $(a)$ represents a system $A$, and 
diagram $(b)$ represents a transformation from system $A$ to system $B$.
 Processes can composed sequentially by connecting the wires 
with matching types. Note that the string diagrams are to be read 
from top to bottom (following the direction of gravity), and from left to right. 

\[ (a) ~~~ \begin{tikzpicture}
	\begin{pgfonlayer}{nodelayer}
		\node [style=none] (4) at (2.5, 7.25) {};
		\node [style=none] (5) at (2.5, 3) {};
		\node [style=none] (7) at (2, 5.25) {$A$};
	\end{pgfonlayer}
	\begin{pgfonlayer}{edgelayer}
		\draw (5.center) to (4.center);
	\end{pgfonlayer}
\end{tikzpicture} ~~~~~~~~~~
 (b) ~~~ \begin{tikzpicture}
	\begin{pgfonlayer}{nodelayer}
		\node [style=circle, scale=2] (3) at (2.5, 5.25) {};
		\node [style=none] (4) at (2.5, 7.25) {};
		\node [style=none] (5) at (2.5, 3) {};
		\node [style=none] (6) at (2.5, 5.25) {$f$};
		\node [style=none] (7) at (2, 7) {$A$};
		\node [style=none] (8) at (2, 3.5) {$B$};
	\end{pgfonlayer}
	\begin{pgfonlayer}{edgelayer}
		\draw (4.center) to (3);
		\draw (3) to (5.center);
	\end{pgfonlayer}
\end{tikzpicture}
~~~~~~~~~~ (c) ~~~ \begin{tikzpicture}
	\begin{pgfonlayer}{nodelayer}
		\node [style=circle, scale=2] (0) at (0, 5.5) {};
		\node [style=none] (1) at (0, 7) {};
		\node [style=none] (2) at (0, 2.5) {};
		\node [style=none] (3) at (0, 5.5) {$f$};
		\node [style=none] (4) at (-0.5, 6.75) {$A$};
		\node [style=none] (5) at (-0.5, 4.75) {$B$};
		\node [style=circle, scale=2] (7) at (0, 4) {};
		\node [style=none] (6) at (0, 4) {$g$};
		\node [style=none] (8) at (-0.5, 2.75) {$C$};
	\end{pgfonlayer}
	\begin{pgfonlayer}{edgelayer}
		\draw (1.center) to (0);
		\draw (7) to (0);
		\draw (2.center) to (7);
	\end{pgfonlayer}
\end{tikzpicture}
  \]

Moreover, the wires and the boxes can be composed in parallel leading to processes as shown in diagrams 
$(e)$ and $(f)$. Morever, the wires are allowed to cross one another another as shown in diagram $(g)$.
\[ (e)~~~ \begin{tikzpicture}
	\begin{pgfonlayer}{nodelayer}
		\node [style=circle, scale=2] (0) at (0, 4.75) {};
		\node [style=none] (1) at (-0.75, 6.5) {};
		\node [style=none] (2) at (-0.75, 3) {};
		\node [style=none] (4) at (-0.75, 6.75) {$A_1$};
		\node [style=none] (5) at (0, 6.5) {...};
		\node [style=none] (6) at (0.75, 6.5) {};
		\node [style=none] (7) at (0.75, 3) {};
		\node [style=none] (8) at (0, 3) {...};
		\node [style=none] (9) at (0.75, 6.75) {$A_n$};
		\node [style=none] (10) at (0, 4.75) {$f$};
		\node [style=none] (11) at (-0.75, 2.75) {$B_1$};
		\node [style=none] (12) at (0.75, 2.75) {$B_m$};
	\end{pgfonlayer}
	\begin{pgfonlayer}{edgelayer}
		\draw [in=150, out=-90] (1.center) to (0);
		\draw [in=-150, out=90] (2.center) to (0);
		\draw [in=-30, out=90] (7.center) to (0);
		\draw [in=30, out=-90] (6.center) to (0);
	\end{pgfonlayer}
\end{tikzpicture}
~~~~~~~~~~ (g) ~~~ \begin{tikzpicture}
	\begin{pgfonlayer}{nodelayer}
		\node [style=none] (1) at (-1, 6.5) {};
		\node [style=none] (2) at (-1, 3) {};
		\node [style=none] (4) at (-1, 6.75) {$A$};
		\node [style=none] (6) at (0.25, 6.5) {};
		\node [style=none] (7) at (0.25, 3) {};
		\node [style=none] (9) at (0.25, 6.75) {$C$};
		\node [style=none] (11) at (-1, 2.75) {$B$};
		\node [style=none] (12) at (0.25, 2.75) {$D$};
		\node [style=circle, scale=2] (13) at (-1, 4.75) {};
		\node [style=circle, scale=2] (14) at (0.25, 4.75) {};
		\node [style=none] (15) at (-1, 4.75) {$f$};
		\node [style=none] (16) at (0.25, 4.75) {$g$};
	\end{pgfonlayer}
	\begin{pgfonlayer}{edgelayer}
		\draw (6.center) to (14);
		\draw (14) to (7.center);
		\draw (13) to (2.center);
		\draw (13) to (1.center);
	\end{pgfonlayer}
\end{tikzpicture}
~~~~~~~~~~ (h) ~~~
\begin{tikzpicture}
	\begin{pgfonlayer}{nodelayer}
		\node [style=none] (1) at (-1.25, 6.5) {};
		\node [style=none] (2) at (-1.25, 3) {};
		\node [style=none] (4) at (-1.25, 6.75) {$A$};
		\node [style=none] (6) at (0.25, 6.5) {};
		\node [style=none] (7) at (0.25, 3) {};
		\node [style=none] (9) at (0.25, 6.75) {$B$};
		\node [style=none] (11) at (-1.25, 2.75) {$B$};
		\node [style=none] (12) at (0.25, 2.75) {$A$};
	\end{pgfonlayer}
	\begin{pgfonlayer}{edgelayer}
		\draw [in=90, out=-90, looseness=1.25] (1.center) to (7.center);
		\draw [in=90, out=-90, looseness=1.25] (6.center) to (2.center);
	\end{pgfonlayer}
\end{tikzpicture} \]


It is far simpler to reason about processes using string diagrams since the human brain is good at processing visual information. 
For example, it is quite easy to see that the diagrams below represent the same process: one can prove the diagrams equal by 
fixing the ends of the wires and moving the circles. 

\[ \begin{tikzpicture}
	\begin{pgfonlayer}{nodelayer}
		\node [style=circle, scale=2] (0) at (-0.75, 4.75) {};
		\node [style=none] (1) at (1.25, 6) {};
		\node [style=none] (2) at (1.25, 2.25) {};
		\node [style=none] (10) at (-0.75, 4.75) {$f$};
		\node [style=circle, scale=2] (11) at (0.25, 4.75) {};
		\node [style=none] (12) at (0.25, 7) {};
		\node [style=none] (13) at (-1, 2.25) {};
		\node [style=circle, scale=2] (14) at (1.25, 4.75) {};
		\node [style=none] (15) at (0, 2.25) {};
		\node [style=none] (16) at (-0.75, 7) {};
		\node [style=none] (17) at (1.25, 7) {};
		\node [style=none] (18) at (0.25, 4.75) {$g$};
		\node [style=none] (19) at (1.25, 4.75) {$h$};
	\end{pgfonlayer}
	\begin{pgfonlayer}{edgelayer}
		\draw [in=90, out=-90, looseness=1.25] (1.center) to (0);
		\draw [in=-90, out=90] (2.center) to (0);
		\draw [in=90, out=-90] (12.center) to (11);
		\draw [in=-90, out=90, looseness=0.75] (13.center) to (11);
		\draw [in=-90, out=90, looseness=1.25] (15.center) to (14);
		\draw [in=-90, out=90] (14) to (16.center);
		\draw (17.center) to (1.center);
	\end{pgfonlayer}
\end{tikzpicture} = \begin{tikzpicture}
	\begin{pgfonlayer}{nodelayer}
		\node [style=circle, scale=2] (0) at (1, 4.75) {};
		\node [style=none] (1) at (1, 6) {};
		\node [style=none] (2) at (1, 2.25) {};
		\node [style=none] (10) at (1, 4.75) {$f$};
		\node [style=circle, scale=2] (11) at (-1, 4.75) {};
		\node [style=none] (12) at (0, 7) {};
		\node [style=none] (13) at (-1, 2.25) {};
		\node [style=circle, scale=2] (14) at (0, 4.75) {};
		\node [style=none] (15) at (0, 2.25) {};
		\node [style=none] (16) at (-1, 7) {};
		\node [style=none] (17) at (1, 7) {};
		\node [style=none] (18) at (-1, 4.75) {$g$};
		\node [style=none] (19) at (0, 4.75) {$h$};
	\end{pgfonlayer}
	\begin{pgfonlayer}{edgelayer}
		\draw [in=90, out=-90, looseness=1.25] (1.center) to (0);
		\draw [in=-90, out=90] (2.center) to (0);
		\draw [in=90, out=-90] (12.center) to (11);
		\draw [in=-90, out=90, looseness=0.75] (13.center) to (11);
		\draw [in=-90, out=90, looseness=1.25] (15.center) to (14);
		\draw [in=-90, out=90] (14) to (16.center);
		\draw (17.center) to (1.center);
	\end{pgfonlayer}
\end{tikzpicture} = \begin{tikzpicture}
	\begin{pgfonlayer}{nodelayer}
		\node [style=circle, scale=2] (0) at (1, 4.75) {};
		\node [style=none] (1) at (1, 6) {};
		\node [style=none] (2) at (1, 2.25) {};
		\node [style=none] (10) at (1, 4.75) {$f$};
		\node [style=circle, scale=2] (11) at (0, 4.75) {};
		\node [style=none] (12) at (0, 7) {};
		\node [style=none] (13) at (-1, 2.25) {};
		\node [style=circle, scale=2] (14) at (-1, 4.75) {};
		\node [style=none] (15) at (0, 2.25) {};
		\node [style=none] (16) at (-1, 7) {};
		\node [style=none] (17) at (1, 7) {};
		\node [style=none] (18) at (0, 4.75) {$g$};
		\node [style=none] (19) at (-1, 4.75) {$h$};
	\end{pgfonlayer}
	\begin{pgfonlayer}{edgelayer}
		\draw [in=90, out=-90, looseness=1.25] (1.center) to (0);
		\draw [in=-90, out=90] (2.center) to (0);
		\draw [in=90, out=-90] (12.center) to (11);
		\draw [in=-90, out=90] (13.center) to (11);
		\draw [in=-90, out=90, looseness=1.25] (15.center) to (14);
		\draw [in=-90, out=90] (14) to (16.center);
		\draw (17.center) to (1.center);
	\end{pgfonlayer}
\end{tikzpicture} \]

In compact closed categories, one can additionally bend wires into a cap and a cup as follows, thus adding to the expressive power of the language: 
\[  	\begin{tikzpicture}
		\begin{pgfonlayer}{nodelayer}
			\node [style=none] (0) at (-6, 5) {};
			\node [style=none] (1) at (-4.5, 5) {};
			\node [style=none] (2) at (-6, 6.25) {};
			\node [style=none] (3) at (-4.5, 6.25) {};
		\end{pgfonlayer}
		\begin{pgfonlayer}{edgelayer}
			\draw (0.center) to (2.center);
			\draw [bend left=90, looseness=1.50] (2.center) to (3.center);
			\draw (3.center) to (1.center);
		\end{pgfonlayer}
	\end{tikzpicture} ~~~~~~~~ \begin{tikzpicture}
		\begin{pgfonlayer}{nodelayer}
			\node [style=none] (0) at (-6, 6.25) {};
			\node [style=none] (1) at (-4.5, 6.25) {};
			\node [style=none] (2) at (-6, 5) {};
			\node [style=none] (3) at (-4.5, 5) {};
		\end{pgfonlayer}
		\begin{pgfonlayer}{edgelayer}
			\draw (0.center) to (2.center);
			\draw [bend right=90, looseness=1.50] (2.center) to (3.center);
			\draw (3.center) to (1.center);
		\end{pgfonlayer}
	\end{tikzpicture} \]


Around the same time that Coecke and Absramsky applied the graphical calculus of monoidal categories to study quantum mechanics, 
Selinger used structures in monoidal categories for designing a quantum programming language \cite{Sel04}. 
In 2007, Selinger \cite{Sel07} refined and extended  
the framework of Abramsky and Coecke to account for mixed quantum states, 
and coined the term `Dagger($\dagger$) compact closed categories' for the categorical framework of quantum mechanics.  
The $\dagger$-compact closed categories faithfully abstract the structure of finite-dimensional Hilbert spaces, 
thereby enabling a diagrammatic but rigorous reasoning technique for quantum processes and protocols
within the category of finite-dimensional Hilbert spaces and linear maps, $\FHilb$. The category of 
all Hilbert spaces and linear maps, ${\sf Hilb}$, is $\dagger$-monoidal but not compact closed. 

\subsection{Towards arbitrary dimensions in CQM}

CQM has been applied to study problems in areas such as quantum foundations, 
quantum information theory and quantum computing. CQM techniques have been used to study 
causality \cite{KiU19, Coe14} and non-locality \cite{CDKW12, CBS11} in quantum foundations. 
In quantum information theory, it has been used to construct structures 
like quantum Latin squares \cite{MuV16} which are crucial to many protocols in quantum information theory. Perhaps, the most significant 
outcome of CQM is in the field of quantum computation namely the ZX-calculus \cite{CoD11} which is a fine-grained 
diagrammatic calculus providing a set of generators and rewrite rules for designing and optimizing quantum circuits.  

Because the compact closed setting is restricted to finite dimensional Hilbert spaces, practical applications of CQM, 
so far, have been limited to the areas such as quantum computing and quantum information theory 
which are off-shoots of finite dimensional quantum mechanics. This is because the base framework of CQM namely the
compact closed categories impose finite dimensionality on the Hilbert Spaces \cite{Heunen16}.  Hence, the category of 
all Hilbert spaces is $\dagger$-monoidal and not $\dagger$-compact closed. The success of CQM in 
the study of finite-dimensional processes has inspired researchers to explore strategies for extending CQM to infinite dimensional processes. 

One approach was to identify the algebraic structures which can characterize the key components of quantum mechanics 
without imposing any restriction on dimensionality. For example, in CQM, dagger Frobenius algebras (See Section \ref{Sec: observables})
provide a precise algebraic characterization of quantum observables 
in the category of finite dimensional Hilbert Spaces  \cite{CPV12}  (A quantum observable is a measurable physical property of a quantum system).  

With the aim of extending this idea to quantum observables of arbitrary dimensions,  
Abramsky and Heunen, in \cite{AbH12}, showed how Ambrose's $H^*$-algebras \cite{Amb45} 
could be used to characterize orthonormal bases, hence quantum observables, in infinite dimensional Hilbert spaces. 
However, the move from $\dagger$-Frobenius to $H^*$-algebras comes at a cost -  an
$H^*$-algebra is modelled as a semi-Frobenius algebra, that is, Frobenius algebras without the units.


Gogioso and Genovese \cite{GG17} proposed an interesting approach to reinstating the units using techniques from non-standard analysis \cite{Rob96}.  
They considered $~^\star${\bf Hilb}, the category of non-standard separable Hilbert Spaces and linear maps.  This they claimed is a 
$\dagger$-compact closed category, in which, among other things, the semi-Frobenius algebras of Abramsky and Heunen can be modelled.  
Furthermore, the counit can be reinstated because formal infinite sums are permitted.  

In \cite{CoH16}, Coecke and Heunen, in order to model infinite dimensional  quantum processes, 
took the simple step of dropping the requirement that the category is compact closed and worked in dagger 
symmetric monoidal categories ($\dagger$-SMCs). The category of all Hilbert 
spaces, $\Hilb$, is the prototypical example of a $\dagger$-SMC.  
In CQM, quantum processes are modelled as completely positive maps in the 
category of finite-dimensional Hilbert Spaces. Coecke and Heunen  
showed how to build a category of completely positive maps for  
an arbitrary $\dagger$-SMC. A downside of this approach is that by moving to $\dagger$-SMCs, 
one loses the structural richness provided by duals.

In \cite{HeR18}, Heunen and Reyes considered a different $\dagger$-SMC, namely, the category of Hilbert $\C^*$-modules.  Its objects can be equivalently viewed as 
bundles of Hilbert spaces over a locally compact Hausdorff space.  They characterized the special commutative $\dagger$-Frobenius algebras in this
category as bundles of finite dimensional Hilbert spaces (with dimensions that are uniformly bounded).  These objects, while being far from finite, do retain a 
(uniform) locally finite nature.  This example,  by using vector bundles and ideas from differential geometry, enters the domain of traditional 
theoretical physics, and serves as a reminder that $\dagger$-Frobenius algebras are not only of interest for Hilbert spaces. 


 In \cite{Vic08}, Vicary attempted to model a quantum harmonic oscillator, 
 which is inherently an infinite-dimensional 
quantum system in $\dagger$-monoidal categories using free exponential modalities. 
He proposed a notion of $\dagger$-exponentials in 
$\dagger$-symmetric monoidal categories with $\dagger$-biproducts and 
used it to derive an abstract Fock Space and the ladder operators of Fock Spaces. 
He conjectured that the category of countable dimensional inner product spaces and 
everywhere-defined linear maps is a model for this categorical setting. 





\subsection{Non-compact multiplicative $\dagger$-linear logic}

While various strategies have been tried for modelling systems of arbitrary dimensions, 
it is a general consensus in the CQM community, that the 
dimensionality constraint of the structures used in CQM is yet to be satisfactorily addressed. 
In this thesis, a different approach is taken by shifting focus beyond Hilbert spaces and other models, 
and by moving to non-compact $\dagger$-linear logic.  

Rather than insisting that the infinite dimensional structures are concretely related to Hilbert spaces, 
it is considered that there may be a system of formal types which extend the existing compact logic of CQM.  
An example of such a setting is the embedding of the category of finite dimensional Hilbert spaces  (which is $\dagger$-compact)
into the ($\dagger$-)$*$-autonomous category of Chu Spaces over complex vector spaces, $\Chus_{\sf Vec(\C)}(I)$, 
where $I$, the tensor unit of ${\sf Vec}(\C)$, is the dualizing object (See Section \ref{Section: Chu}). Yet another example is the embedding 
of finite dimensional complex matrices into the  ($\dagger$-)$*$-autonomous 
category of finiteness spaces and finiteness matrices over complex numbers \cite{Ehr05} (see Section \ref{Sec: Finiteness matrices}).

We begin our explorations with linearly distributive categories (LDCs) and $*$-autonomous categories rather 
than monoidal categories to first obtain a semantics for non-compact $\dagger$-linear logic. 
This idea is not new.  Models for quantum mechanics in $*$-autonomous categories are often described 
as ``toy models'' \cite{Abr12} and were, in particular discussed by Pavlovic \cite{Pav11} where some very similar 
directions were advocated.   Indeed, Egger \cite{Egg11}, in initiating the development of ``involutive'' categories, 
implicitly suggested that a dagger functor is not stationary ($A \neq A^\dag$) on objects in an LDC setting. 


\subsection{Linearly distributive categories}
\label{Sec: ldc-intro}

Linearly distributive categories (LDCs) \cite{CS97} provide the categorical semantics (so they are the proof theory)
for  the non-compact multiplicative fragment of linear logic (MLL) containing the 
$\ox$, $1$, $\parr$, and $\bot$. The sequent rules for the multiplicatives are shown in Figure \ref{Fig: multiplicative rules}.
In the linear logic community the multiplicative disjunction 
is often denoted by $\parr$, however this thesis follows the convention in \cite{CS97} and shall use $\oplus$. 
Under the same convention, we write the unit of  $\ox$ as top $\top$ rather than $1$ as in the linear logic community. 
The unit of $\oa$ is a bottom, $\bot$. 

\begin{figure}[h]
	\centering
		\AxiomC{$\Gamma_1, \Gamma_2  \vdash \Delta$}
		\LeftLabel{$(\top L)$}
		\UnaryInfC{$\Gamma_1, \textcolor{blue}{\top}, \Gamma_2 \vdash \Delta$}
	   \DisplayProof 
        \hspace{1.5 em}
		\AxiomC{$ $}
		\LeftLabel{$(\top R)$}
		\UnaryInfC{$\vdash \textcolor{blue}{\top}$} 
		\DisplayProof 

		\vspace{1em}
		
		\AxiomC{$\Gamma_1, \textcolor{blue}{A, B}, \Gamma_2 \vdash \Delta$}
		\LeftLabel{$(\ox  L)$}
		\UnaryInfC{$\Gamma_1, \textcolor{blue}{A \ox B}, \Gamma_2 \vdash \Delta$}
	   \DisplayProof
	   \hspace{1.5em}
		\AxiomC{$\Gamma_1 \vdash \Gamma_2, \textcolor{blue}{A}, \Gamma_3 ~~~~ \Delta_1 \vdash \Delta_2, \textcolor{blue}{B}, \Delta_3$}
		\LeftLabel{$(\ox  R)$}
		\UnaryInfC{$\Gamma_1, \Delta_1 \vdash \Gamma_2, \Delta_2, \textcolor{blue}{A \ox B}, \Gamma_2, \Delta_3$}
		\DisplayProof

		\vspace{1em}

		\AxiomC{$ $}
		\LeftLabel{$(\bot  L)$}
		\UnaryInfC{$\textcolor{blue}{\bot} \vdash $}
	    \DisplayProof 
        \hspace{1.5 em}
		\AxiomC{$\Gamma \vdash \Delta_1, \Delta_2$}
		\LeftLabel{$(\bot  R)$}
		\UnaryInfC{$\Gamma \vdash \Delta_1, \textcolor{blue}{\bot}, \Delta_2$}
	    \DisplayProof 
        
		\vspace{1 em}		
		
		\AxiomC{$\Gamma_1, \textcolor{blue}{A}, \Gamma_2 \vdash \Gamma_3 ~~~~ \Delta_1, \textcolor{blue}{B}, \Delta_2 \vdash \Delta_3$}
		\LeftLabel{$(\oa  L)$}
		\UnaryInfC{$\Gamma_1, \Delta_1,  \textcolor{blue}{A \oa B},  \Gamma_2, \Delta_2 \vdash \Gamma_3, \Delta_3$}
	    \DisplayProof 
        \hspace{1.5 em}		
		\AxiomC{$\Gamma \vdash \Delta_1, \textcolor{blue}{A, B}, \Delta_2$}
		\LeftLabel{$(\oa  R)$}
		\UnaryInfC{$\Gamma \vdash  \Delta_1, \textcolor{blue}{A \oa B}, \Delta_2$}
	    \DisplayProof 
\caption{Sequent rules for multiplicatives}
\label{Fig: multiplicative rules}
\end{figure}

\FloatBarrier

We assume that the sequent calculus is commutative (the order of premises and antecedents does not matter). 
This is accommodated by the exchange rules which allows the neighboring premises and antecedents to be swapped: 
\[ 
\AxiomC{$\Gamma_1, \textcolor{blue}{A, B}, \Gamma_2 \vdash \Delta$} 
\LeftLabel{$(exch.L)$}
\UnaryInfC{$\Gamma_1, \textcolor{blue}{B, A}, \Gamma_2 \vdash \Delta$}
\DisplayProof 
\hspace{1.5em}
\AxiomC{$\Gamma \vdash \Delta_1, \textcolor{blue}{C,D}, \Delta_2$} 
\LeftLabel{$(exch.R)$}
\UnaryInfC{$\Gamma \vdash \Delta_1, \textcolor{blue}{D, C}, \Delta_2$} 
\DisplayProof 
\]

Since LDCs are the proof theory of MLL, they come equipped with two distinct tensor products 
called the ``tensor'', $\ox$, and the ``par'', $\oplus$ 
corresponding to the multiplicative conjunction and disjunction of linear logic respectively. 
The tensor and the par are related by two natural transformations called linear distributors: 
\[ \partial^L: A \ox (B \oa C) \to (A \ox B) \oa C~~~~~~~~ \partial^R: (A \oa B) \ox C \to A \oa (B \ox C) \] 
The distributors are not isomorphisms. The following is the sequent derivation of $\partial^L$:

\vspace{0.5em}

\begin{center}
    \AxiomC{$ $}
    \RightLabel{$id$, $\oa L$, cut}
    \UnaryInfC{$B \oa C \vdash B, C$}
    
    \AxiomC{$ $}
    \RightLabel{$id$, $\ox R$, cut}
    \UnaryInfC{$A,B \vdash A \ox B$}
    
    \RightLabel{cut}
    \BinaryInfC{$A, B \oa C \vdash A \ox B, C$}
    \RightLabel{$(\ox L), (\oa R)$}
    \UnaryInfC{$A \ox (B \oa C) \vdash (A \ox B) \oa C$}
    \DisplayProof
\end{center} 

The multiplicative fragment with negation ($A^\perp$) has its categorical semantics in $*$-autonomous categories.
The sequent rules for negation allow premises to be flipped to the opposite side of the entailment:

\begin{figure}[h]
\centering
\AxiomC{$\Gamma, \textcolor{blue}{B^\perp} \vdash \Delta$}
\LeftLabel{$(Neg.L)$}
\UnaryInfC{$\Gamma \vdash \Delta, \textcolor{blue}{B}$}
\DisplayProof 
\hspace{1.5 em}
\AxiomC{$\Gamma \vdash \Delta, \textcolor{blue}{B}$}
\LeftLabel{$(Neg.R)$}
\UnaryInfC{$\Gamma, \textcolor{blue}{B^\perp}  \vdash \Delta$}
\DisplayProof 
\caption{Sequent rules for negation}
\label{Fig: negation rules}
\end{figure}

In an LDC, negation of an object is given by its categorical dual. An object $A$ has a dual $A^\perp$ (negation of $A$) 
if there exists two maps, $\eta: \top \to A \oa A^\perp$, 
and $A^\perp \ox A \to \bot$ satisfying the `snake' equations, see Definition \ref{defn: duals}.
 An LDC in which every object has a chosen dual is a $*$-autonomous category \cite{Bar91}, 
 see Section \ref{Sec: *-autonomous}. 
 The following proofs show the derivations of the $\eta$ and the $\epsilon$ maps respectively
 from the negation rules:

 \vspace{0.5em}

 \begin{center}
	\AxiomC{$ $}
	 \RightLabel{(id)}
	 \UnaryInfC{$A \vdash A$}
	 \RightLabel{Neg.R}
	 \UnaryInfC{$\vdash A^\perp, A$}
	 \RightLabel{$(\oa R)$}
	 \UnaryInfC{$ \vdash A \oa A^\perp $}
	 \RightLabel{$(\top L)$}
	 \UnaryInfC{$\top \vdash A \oa A^\perp$}	
	 \DisplayProof
	 \qquad
	  \AxiomC{$ $}
	 \RightLabel{(id)}
	 \UnaryInfC{$A \vdash A$}
	 \RightLabel{Neg.L}
	 \UnaryInfC{$A^\perp, A \vdash $}
	 \RightLabel{$(\ox L)$}
	 \UnaryInfC{$A^\perp \ox A \vdash $}
	 \RightLabel{$(\bot R)$}
	 \UnaryInfC{$A^\perp \ox A \vdash \bot$}
	 \DisplayProof
 
	 \vspace{0.75 em}
 
	 Proof for $\eta$ ~~~~~~~~~~~~~~~~ Proof for $\epsilon$
 
  \end{center}

An {\bf isomix category} is an LDC with an isomorphism, 
called the mix map, $\m: \bot \to \top$, satisfying the `mix law' \cite[Definition 6.2]{BCS00}, 
see Section \ref{Sec: mix, isomix, compact LDC}. The mix law provides a natural transformation, $\mx: A \ox B \to A \oa B$, 
from the multiplicative conjunction to the multiplicative disjunction, called the {\bf mixor}, see Section \ref{Sec: mix, isomix, compact LDC}. 
The mix law corresponds to the sequent rules shown in Figure \ref{Fig: mix map} (a)-(b). 

\begin{figure}[h]
	\centering
	$(a)~~~$ 
	\AxiomC{$\Gamma \vdash \Delta$}
	\AxiomC{$\Gamma' \vdash \Delta'$}
	\BinaryInfC{$\Gamma, \Gamma' \vdash \Delta, \Delta'$}
	\DisplayProof
	$~~~~~~(b)~~$ \AxiomC{$ $} 
	\UnaryInfC{$\top \vdash \bot$}
	\DisplayProof
	\caption{(a) Binary mix axiom (b) Nullary mix axiom}
	\label{Fig: mix map}
\end{figure}

	\FloatBarrier

In the presence of the cut rule, the rule $(a)$ is equivalent to 
the axiom $\bot \vdash \top$ \cite[Lemma 6.1]{CS97a}, see below. 

\begin{center}
    \AxiomC{$\Gamma \vdash \Delta$}
    \RightLabel{$(\bot L)$}
    \UnaryInfC{$\Gamma \vdash \Delta, \bot$}

    \AxiomC{$ $}
    \RightLabel{$id$}
    \UnaryInfC{$\bot \vdash \top$}

    \RightLabel{cut}
    \BinaryInfC{$\Gamma \vdash \Delta, \top$}

    \AxiomC{$\Gamma' \vdash \Delta'$}
    \RightLabel{ $(\top L)$}
    \UnaryInfC{$\top, \Gamma' \vdash \Delta'$}

    \RightLabel{cut}
    \BinaryInfC{$\Gamma, \Gamma' \vdash \Delta, \Delta'$}

    \DisplayProof

	\vspace{0.75em}

    Derivation of the binary mix axiom in the presence of the cut rule
\end{center}

In an isomix category, an object $U$ for which the natural transformations, $\mx_{(-, U)}$ and $\mx_{(U, -)}$ are isomorphisms
(that is for any object $A$, $A \ox U \simeq A \oa U$ and $U \ox A \simeq U \oa A$) is said to be in the {\bf core} 
of the category. The core of an isomix category determined by such objects forms a full subcategory.  

An isomix category in which the mixor is a natural isomorphism ($\ox \simeq \oa$) is called a {\bf compact} LDC. 
 The core of an isomix category is always a compact LDC. Compact LDCs are linearly equivalent to a monoidal categories.  
 Conversely, {\bf monoidal} categories \cite{Mac13} can be viewed as being degenerate  
 compact LDCs in which the mixor is the identity map ($\ox = \oa$). 
 A monoidal category in which every object has a dual, is a compact closed category \cite{Kel80}. 
 From this perspective, a compact closed category can be viewed as a compact $*$-autonomous category.  



\section{Thesis outline}

The aim of this thesis is to lay the categorical foundations for non-compact $\dagger$-linear logic and to apply it 
to categorical quantum mechanics. With this goal in mind, the thesis is divided into 
two parts: 

\begin{description} 
	\item[Part 1:] composed of chapters \ref{chap: LDC} - \ref{chap: part 1 summary} 
is focused on developing the {\bf categorical semantics} of non-compact $\dagger$-linear logic. 
     \item[Part 2:] composed of chapters \ref{Chap: CQM} - \ref{Chap: part 2 summary}, is focused on {\bf developing structures 
for CQM} in the categorical setting developed in Part 1.
\end{description}

The rest of the section outlines the contents of this thesis and indicates the contributions. 

\subsection{Part I: Dagger linear logic}

 This first part begins with Chapter  \ref{chap: LDC} which provides an introduction to 
LDCs and its variants, linear functors and transformations. This chapter also discusses the 
Ehrhard's Finiteness spaces in detail, which is used as an example throughout this thesis. 
Chapters \ref{Chap: dagger-LDC} and \ref{Chap: MUCs} 
chapters are derived from the article titled `Dagger linear logic for categorical quantum mechanics' \cite{CCS18}:
it was presented at the Symposium on Compositional Structures (SYCO I) in Birmingham (U.K.) and as a poster
 at the 15th international conference in Quantum Physics and Logic (QPL), Halifax, Canada. 
The contributions of chapters \ref{Chap: dagger-LDC} and \ref{Chap: MUCs} are outlined below. 
Chapter \ref{chap: part 1 summary} summarizes the first part of this thesis. 

 \vspace{1em}

\noindent{\bf Dagger linearly distributive categories}

\nopagebreak

\vspace{0.5em}

\nopagebreak

It is standard in CQM to interpret the dagger as a contravariant functor which is 
stationary on objects ($A = A^\dagger$) and an involution for maps ($ f^{\dag \dag} = f$). 
However, in an LDC with two tensor products - 
a tensor $\ox$, and a par $\oa$ - the dagger has to flip the tensor and par so that 
$(A \oa B)^\dagger = A^\dagger \ox B^\dagger$. Without such a flip, daggering the linear distributor would 
produce a non-permissible map. 

This implies that the dagger can no longer be viewed as being stationary on objects in an LDC. 
A non-stationary dagger implies that one has to address the coherence issues  determining how the 
dagger interacts with the structures of an LDC.  Moreover, one can replace the equality 
above by a natural isomorphism $ \lambda_\ox: A^\dag \ox B^\dag \to (A \oa B)^\dag$, and 
the involution by a natural isomorphism $\iota_A: A \to A^{\dag \dag}$. 
We deal with these natural isomorphisms and define $\dagger$-LDC, 
$\dagger$-mix and $\dagger$-isomix categories in Section \ref{Section: dagger LDC}. The sequent 
calculus for $\dagger$-linear logic is discussed in Section \ref{Sec: dagger sequent rules}.

\vspace{1em}

\noindent{\bf Unitary isomorphisms for $\dagger$-isomix categories}  

\nopagebreak

\vspace{0.5em}

\nopagebreak

In CQM, the dagger functor determines the notion of a unitary isomorphism which 
is an isomorphism $f: A \to B$ such that $f^\dag = f^{-1}$. These isomorphisms are particularly 
important since they model the unitary evolution of a quantum system. Applying this idea 
directly to $\dagger$-LDCs is not feasible, because the maps, $f^{-1}: B \to A$ and 
$f^\dag: B^\dag \to A^\dag$ now have different types, hence cannot be directly equal. 
However, minimally, if there exist isomorphisms, $\varphi_A: A \simeq A^\dag$ and $\varphi_B: B \simeq B^\dag$, 
then one can define a unitary isomorphism to be a map $f$ satisfying the following commuting diagram: 
\[  \xymatrix{
	A \ar[r]^{f} \ar[d]_{\varphi_A} & B \ar[d]^{\varphi_B} \\
	A^\dag \ar@{<-}[r]_{f^\dag} & B^\dag } \] 
Note that, in a $\dagger$-monoidal category, the isomorphisms $\varphi_A$ 
and $\varphi_B$ are simply the identity maps. 

In a $\dagger$-isomix category, the isomorphism $\varphi_A$ is referred to as a
{\em unitary structure} map  if $A$ resides in the core, and behaves coherently 
with the natural isomorphisms of the category. In this case, $A$ is referred to as unitary objects. 
The unitary objects of $\dagger$-isomix category forms a sub compact-$\dagger$-isomix category. 
Such a compact $\dagger$-isomix category in which every object is unitary is called a {\em unitary category}. 
Unitary categories are $\dagger$-linearly equivalent to a $\dagger$-monoidal categories. Moreover, 
when unitary categories have unitary duals, the category is $\dagger$-linearly equivalent to 
$\dagger$-compact closed categories. These ideas are developed in Section \ref{Sec: unitary}.

 \vspace{1em}

 \noindent {\bf Unitary construction} 

 \nopagebreak

 \vspace{0.5em}
 
 \nopagebreak
 
 Having defined the $\dagger$-LDCs and the unitary structure for $\dagger$-isomix categories, our next objective is to 
extract a unitary category, that is the traditional CQM setting, from any $\dagger$-isomix category.
For this, one collects the ``pre-unitary'' objects  of a $\dagger$-isomix category. An object $A$ 
is  {\bf pre-unitary} if $A$ is in the core and there exists an isomorphism $\varphi_A: A \to A^\dag$ and the isomorphism 
behaves coherently with the involutor $\iota$ as shown in the following commuting diagram: 
\[ \xymatrix{
	A \ar[d]_{\varphi_A} \ar[dr]^{\iota} & \\ 
	A^\dag \ar[r]_{\varphi_A^{-1 \dag}} & A^{\dag \dag}} \]
The unitary category thus extracted using the `unitary contruction' is referred to as the {\bf canonical unitary core} 
of the $\dagger$-isomix category. See Section \ref{Sec: unitary construction} 
for details. 

\vspace{1em}

\noindent{\bf Mixed unitary categories} 

\nopagebreak

\vspace{0.5em}

\nopagebreak

Collecting the structures we have developed so far, namely the $\dagger$-isomix categories and 
the unitary categories, we arrive at a general framework of `Mixed Unitary categories' (MUCs). 
MUCs embed a unitary category in a larger $\dagger$-isomix category via a $\dagger$-isomix-functor 
from the unitary category to the core of the $\dagger$-isomix category: 
\[ M: \U \to \Core(\C) \hookrightarrow \C \]  
where $\U$ is the unitary category and $\C$ is the $\dagger$-isomix category. 
Thus, a MUC is a general framework which encompasses the traditional CQM framework within its core. 
The unitary category of a MUC acting on the larger $\dagger$-isomix category is analogous to a field $K$ 
acting on a $K$-algebra as scalars. Sections \ref{Sec: MUCs} and \ref{Sec: MUC examples} are dedicated 
towards developing these ideas.

\subsection{Part II: Application to dagger linear logic to CQM}

The second part of this thesis is concerned with applying the MUC framework to quantum mechanics. 
To this end, we generalize the key algebraic structures of CQM to the MUC setting 
and study the implications of this generalization.  

This part begins with Chapter \ref{Chap: CQM} which provides an introduction to $\dagger$-monoidal categories and 
the other key algebraic structures used by CQM. Chapter \ref{Chap: positivity} is derived from 
my coauthored article titled `Complete positivity for Mixed Unitary Categories' \cite{CS19}, which was presented as a talk at the QPL 2019 held in California, USA, and 
at SYCO II held in Glasgow, UK.  Chapters \ref{Chapter: compaction} and \ref{Chapter: complementarity} are 
derived from my coauthored paper titled `Exponential modalities and complementarity' \cite{CS21}. This work 
was presented at the 4th International Conference on Applied Category Theory held in July 2021. 
Chapter \ref{Chap: free exp examples} presents examples for the structures introduced in Chapters 
\ref{Chapter: compaction} and \ref{Chapter: complementarity}. Chapter \ref{Chap: part 2 summary} summarizes the second part of the thesis. 

The contributions of chapters \ref{Chap: positivity} - \ref{Chapter: complementarity} are outlined below. 

\vspace{1em}

\noindent{\bf Completely positive maps in MUCs}

\nopagebreak

\vspace{0.5em}

\nopagebreak

In a $\dagger$-monoidal settings, completely positive maps abstract the notion of quantum processes.
Coecke and Heunen \cite{CoH16} developed the $\CP^\infty$ construction which produces the 
category of completely positive maps when applied to any $\dagger$-symmetric monoidal category.
We define the completely positive maps and the $\CP^\infty$ construction for a MUC based on the 
Coecke and Heunen's construction \cite{CoH16}. We also characterize the $\CP^\infty$ construction on MUCs using the 
notion of environment structures and purification. The characterization of $\CP^\infty$ construction for MUCs led to 
an elegant observation that it suffices for the unitary core to have environment maps to characterize 
the $\CP^\infty$ construction for MUCs. Chapter \ref{Chap: positivity} covers this discussion. 

\vspace{1em}

\noindent{\bf Measurement in MUCs} 

\nopagebreak

\vspace{0.5em}

\nopagebreak

Coecke and Pavlovic provided an algebraic description of quantum measurements in $\dagger$-monoidal categories  \cite{CoP07}.
In a MUC, measurement happens in two steps: first the system to be measured must be {\em compacted} into the 
unitary core (i.e., the traditional CQM core), and then the usual measurement process \cite{CoP07} 
as described by Coecke and Pavlocic must be applied within the core. In order to characterize the compaction process, 
we introduce {\em binary idempotents} and {\em $\dagger$}-binary idempotents. 
Compaction of an object in a $\dagger$-isomix category precisely corresponds to the splitting of certain 
$\dagger$-binary idempotents on the object, 
See Section \ref{Sec: compaction}. The $\dagger$-binary idempotents generalize Selinger's \cite{Sel08}  
$\dagger$-idempotents in $\dagger$-monoidal categories which $\dagger$-splits to produce classical types.  

\vspace{1em}

\noindent{\bf Complementary systems in isomix categories}

\nopagebreak

\vspace{0.5em}

\nopagebreak

The notion of complementary observables is central to quantum mechanics. An observable is a measurable property 
of a quantum system. Two observables are complementary if measuring (knowing) the value of one observable 
increases the uncertainty of the value of the other. A classic example of complementary observables is the 
position and momentum of an electron. 

In CQM, quantum observables are algebraically presented as special commutative $\dagger$-Frobenius algebras \cite{CPV12} in $\dagger$-monoidal categories. Moreover, 
complementary observables are two such $\dagger$-Frobenius algebras on the same object interacting bialgebraically to produce 
two Hopf algebras\footnote{Thus is often referred to as {\em strong} complementarity in CQM.} \cite{CoD11}. 

In LDCs, linear monoids with a $\ox$-monoid and a $\oa$-monoid provide a general version of Frobenius algebras. In fact, Frobenius 
algebras in monoidal categories can be viewed as linear monoids satifying an extra property.  
Linear monoids lead to the new notion of a ``linear comonoid", which can interact bialgebraically 
with a linear monoid to give a ``linear bialgebra".  Using these structures in a ($\dagger$-)isomix category, one can 
define a complementary system as a ($\dagger$-)linear bialgebra satisfying a few extra equations. 
Chapter \ref{Chapter: compaction} is dedicated to developing linear comonoids and linear bialgebras. 
Complementary systems in isomix categories are described in Section \ref{Sec: complementary systems}.

\vspace{1em}

\noindent{\bf Relating complementarity and exponential modalities} 

\nopagebreak

\vspace{0.5em}

\nopagebreak

A final but a significant contribution of this thesis is to establish the connection between exponential modalities 
of linear logic and complementarity of quantum mechanics using our MUC framework. For this, we define 
$(!,?)$-$\dagger$-LDCs, i.e., a $\dagger$-LDC with exponential modalities, and also provide the sequent rules 
for the corresponding logic, see Sections \ref{Sec: dagger exp modalities}, and \ref{Sec: sequent dagger exp}. 
We prove that every complementary 
system in a ($\dagger$-)isomix category arises as a compaction of a ($\dagger$-)linear bialgebra induced on the {\bf free} 
exponential modalities, see Section \ref{Sec: exp modalities}. 

\vspace{1em} 

Chapter \ref{Chap: conclusion} concludes this thesis and dicusses future directions.

\section{Prerequisites}
 We assume that the reader has a basic understanding of category theory including the definition of 
categories, functors, natural transformations, duality, isomorphisms, and adjunctions. We refer the reader 
to a few excellent sources \cite{Awo10, BaC90, Bor94, Coe08,  Bae10, Mac13}, particularly, references \cite{Coe08,  Bae10} 
are well-suited for someone with a Physics background while \cite{Bor94} is geared towards audience in computer science.

\section{Notation}

Composition is written in diagrammatic  order: $fg$ means apply $f$ followed by $g$.
The string diagrams are to be read from top to bottom (following the direction of gravity) and left to right.

%% file: chapter1.tex

\chapter{Categorical semantics for linear logic}
\label{chap: LDC}

Linear logic \cite{Gir87} was introduced by Girard in 1987 as a resource sensitive logic in which 
logical statements were treated as resources. Hence, proving statements in linear logic involves manipulating these 
resources, most of which cannot be duplicated or destroyed. Linear logic has been considered to be the 
logic of quantum information theory due to its resource sensitivity \cite{BaM11, Dun06}. 
In this chapter, a categorical semantics of linear logic is reviewed using linearly distributive categories. We also provide an 
interpretation of the linear logic structures in finiteness spaces \cite{Ehr05} which is used a running example in the 
rest of the thesis. 

\section{Linearly Distributive Categories (LDCs)}
\label{Sec: LDC}

Linearly distributive categories, introduced by Cockett and Seely \cite{CS97} in 1997, provide a categorical semantics of Multiplicative Linear Logic (MLL). 
These categories were originally referred to as Weakly distributive categories and were later renamed 
to linearly distributive categories (LDCs). 

\subsection{Linearly distributive categories}
\label{Subsec: LDC}

Informally, a linearly distributive category is a category having two monoidal structures 
linked by a linear distributor. We first recall the definition of monoidal categories before 
moving on to LDCs:

\begin{definition}
	\label{defn: monoidal}
A {\bf monoidal category} $(\X, \otimes, I, a_\otimes, u_\ox^l, u_\ox^r)$ is a category, $\X$, consisting of:
\begin{itemize}
\item a bifunctor, $\otimes: \X \times \X \rightarrow \X$, called tensor product;
\item a designated object, $I \in \X$, called the unit object;
\item a natural isomorphism, $(a_\otimes)_{A,B,C}: A \otimes (B \otimes C) \xrightarrow{\simeq} (A \otimes B) \otimes C$, called the associator;
\item a natural isomorphism, $ (u_\ox^l)_A: I \otimes A \xrightarrow{\simeq} A$, called the left unitor;
\item a natural isomorphism, $(u_\ox^r)_A: A \otimes I \xrightarrow{\simeq} A$, called the right unitor;
\end{itemize}
such that the following coherence diagrams commute \cite{Kel64}:
\begin{itemize}
\item Maclane's pentagon diagram:
\[ \xymatrixrowsep{8mm} \xymatrix{
 & A\otimes(B \otimes(C \otimes D)) \ar[dl]_{1 \otimes a_\otimes} \ar[dr]^{a_\otimes} &  \\
A \otimes ( (B \otimes C) \otimes D) \ar[d]_{a_\otimes} & &  (A \otimes B) \otimes (C \otimes D) \ar[d]^{a_\otimes} \\
 (A \otimes (B \otimes C)) \otimes D \ar[rr]_{a_\otimes \otimes 1} & &   ((A \otimes B) \otimes C) \otimes D }  \]

\item Kelly's unit diagram:
 \[ \xymatrix{
(A \otimes I) \otimes B \ar[d]_{a_\otimes} \ar[drr]^{u_\otimes^r \otimes 1_B}  \\
A \otimes (I \otimes B) \ar[rr]_{1_A \otimes u_\otimes^l} & & (A \otimes B)} \]
\end{itemize}
\end{definition}

A monoidal category in which the associator, the left unitor and the right unitor are identity arrows is called a {\bf strict monoidal category}. 

A {\bf symmetric monoidal category} (SMC) is a monoidal category with a natural isomorphism:
\[ (c_\otimes)_{A,B}: A \otimes B \xrightarrow{\simeq} B \otimes A \]
such that the following equations hold:
\begin{itemize}
\item (Hexagon law) $a_\ox c_\ox a_\ox = (1 \ox c_\ox) a_\ox (c_\ox \ox 1)$ 
\[ \xymatrix{
	A \otimes (B \otimes C) \ar[r]_{a_\otimes} \ar[d]^{1 \otimes c_\otimes} 
	& (A \otimes B) \otimes C  \ar[r]_{c_\otimes} 
	& C \otimes (A \otimes B) \ar[d]^{a_\ox}\\ 
	\otimes (C \otimes B) \ar[r]_{ a_\ox} 
	& (A \otimes C) \otimes B \ar[r]_{a_\ox (c_\ox \ox 1)} 
	& (C \otimes A) \otimes B 
} \]
\item (Inverse law) $(c_\ox)_{A, B} = (c_\ox)_{B,A}^{-1}$
 \[ \xymatrix{
	A\otimes B \ar@{=}[dr] \ar[d]_{(c_\otimes)_{A,B}}  \\
	B \otimes A  \ar[r]_{(c_\otimes)_{A,B}} & A \otimes  B}\]
\item (Unit law) $c_\ox u_\ox^l = u_\ox^r$
\[ \xymatrix{
	A \otimes I \ar[r]^{c_\otimes} \ar[dr]_{u^r} & I \otimes A \ar[d]^{u^l} \\
	& A & } \]
\end{itemize}

In an LDC, there exists two tensor products which are related by natural transformations called linear distributors:
\begin{definition} \cite{CS97}
A {\bf linearly distributive category}, $(\X, \ox, \top, \oa, \bot)$ is a category $\X$ 
consisting of:

\begin{itemize}

\item a monoidal structure, $(\ox, \top, a_\ox, u_\ox^L, u_\ox^R)$

($\ox$ is referred to as `the tensor' and its unit, $\top$, `the top')

\item a monoidal structure, $(\oa, \bot, a_\oa, u_\oa^L, u_\oa^R)$

($\oa$ is referred to as `the par' and its unit, $\bot$, `the bottom')

\item The tensor and the par are related by the following natural transformations 
which are called the left and the right linear {\bf distributors} respectively:
\begin{align*}
& \partial^l: A \ox (B \oa  C) \to (A \ox B) \oa C \\
& \partial^r: (A \oa B) \ox C \to A \oa (B \ox C) 
\end{align*}
\end{itemize}

satisfying the following coherence conditions:

\begin{itemize}

\item The assosicators and the unitors for the $\ox$ and the $\oa$ satisfy Maclane's pentagon diagram 
and Kelly's unit diagram, See Definition \ref{defn: monoidal}.

\item Coherence conditions for unit natural transformations and linear distributors:

\begin{enumerate}[LDC. 1]
\item \begin{enumerate}[(a)]
\item $\partial^l (u_\ox^l \oa 1_B) = u_\ox^l$

\[ \xymatrixcolsep{4pc} \xymatrix{
\top \ox (A \oa B) \ar[d]_{\partial^l} \ar[dr]^{u_\ox^l} \\
(\top \ox A) \oa B \ar[r]_{u_\ox^l \oa 1_B} & A \oa B } \]

\item $u_\ox^r = \partial^r; 1 \oa u_\ox^l$
\item $ \partial^r u_\oa^l = u_\oa^l \ox 1_B $
\item $1 \ox u_\oa^r = \partial^l; u_\oa^l$

\[ \xymatrixcolsep{4pc}
\xymatrix{ 	(\bot \oa A) \ox B \ar[r]^{\partial^r} \ar[dr]_{u_\oa^l \ox 1_B} & \bot \oa 
(A \ox B) \ar[d]^{u_\oa^l} \\
& A \ox B } \] 

\end{enumerate}
\end{enumerate}

\item Coherences for associativity natural transformations and the distributors:

\begin{enumerate}[LDC. 2]
\item \begin{enumerate}[(a)]
\item $a_\ox (1_A \ox \partial^l) \delta^l = \delta^l (a_\ox \oa 1_D)$

\[ \xymatrixcolsep{4pc} \xymatrix{
(A \ox B) \ox (C \oa D) \ar[r]^{a_\ox} \ar[dd]_{\partial^l} & A \ox (B \ox ( C \oa D )) \ar[d]^{1_A \ox \partial^l} \\
 & A \ox ((B \ox C) \oa D) \ar[d]^{\partial^l} \\
 ((A \ox B) \ox C) \oa D \ar[r]_{a_\ox \oa 1_D} & (A \ox (B \ox C)) \oa D } \]

\item $\partial^l (a_\ox \oa 1) = a_\ox (1 \ox \partial^l) \partial^l$

\item  $\partial^r a_\oa = (a_\oa \ox 1_D) \partial^r (1_A \ox \partial^r)$

\item $(1 \ox a_\oa) \partial^l = \partial^l (1 \oa \partial^l) a_\oa$
\end{enumerate}
\end{enumerate}

\item Coherences between the left and the right linear distributors: 

\begin{enumerate}[LDC. 3]
\item \begin{enumerate}[(a)]
\item $\partial^r (1_A \oa \partial^l) = \partial^l (\partial^r \oa 1_D) a_\ox$
\[ \xymatrix{
 & (A \oa B) \ox (C \oa D) \ar[ld]_{\partial^l} \ar[dr]^{\partial^r} & \\
((A \oa B) \ox C) \oa D) \ar[d]_{\partial^r \oa 1_D} &  & A \oa ( B \ox (C \oa D)) \ar[d]^{1_A \oa \partial^l} \\
(A \oa (B \ox C)) \oa D \ar[rr]_{a_\ox} &  & A \oa ((B \ox C) \oa D) } \]
\item $a_\ox (1 \ox \partial^r) \partial^l = (\partial^l \ox 1) \partial^r$
\end{enumerate}
\end{enumerate}
\end{itemize}

\end{definition}
A {\bf symmetric LDC} is an LDC in which both the tensor products are symmetric,  
with symmetry maps $c_\ox$ and $c_\oa$, such that
$\partial^R = c_\ox (1 \ox c_\oa) \partial^L (c_\ox \oa 1) c_\oa$. 
For a symmetric LDC, the left linear distributor determines the right linear distributor and vice versa.

Linearly distributive categories (LDCs) provide a categorical semantics for Multiplicative Linear 
Logic (MLL). The two tensor products, $\ox$ and $\oa$, in an LDC corresponds to the multiplicative 
conjunction and multiplicative disjunction in linear logic respectively.
 
The following are some examples of LDCs:
\begin{itemize}
\item Every monoidal category is also an LDC where the tensor and the par coincide, and the distributor is 
the associator natural isomorphism. An LDC with the linear distributors being isomorphisms is not 
necessarily monoidal. The next two examples elucidate the point. 

\item A {\bf bounded distributive lattice} is a lattice $(L, \leq, \wedge, \top, \vee, \bot)$ with a greatest element, 
$\top$, and a least element, $\bot$, such that for all $a \in L$, $\bot \leq a \leq \top$, and the join 
$(\wedge)$ and meet $(\vee)$ operations distribute over one another:
\[ a \wedge (b \vee c) = (a \wedge b) \vee (a \wedge c) ~~~~~~~~~~~
   a \vee (b \wedge c) = (a \vee b) \wedge (a \vee c) \]
A distributive lattice regarded as a category whose objects are lattice elements and the maps 
given by the preorder, is an LDC. The tensor is given by $\wedge$ with unit object $\top$,
 and the par is given by $\vee$ with unit object $\bot$. Both the tensor products are 
symmetric. The right linear distributor is given as follows:
\[ (a \vee b) \wedge c = (a \wedge c) \vee (b \wedge c) \leq a \vee (b \wedge c) \]
Any Boolean algebra is an example of a distributive lattice. 

\item Let $M$ be any set. A {\bf shift monoid} is a commutative monoid, $(M, +, 0)$ with a designated 
element $s$ such that there exists an inverse to $s$ i.e, $s - s = 0$. 
A second multiplication can be defined on  the set $M$ as follows: for all $x,y \in M$, $x \circ y = (x + y) - s$. The 
unit of the second multiplication is $s$. A shift monoid considered as a discrete category 
(the elements of the monoid are the objects and the maps are identity maps) is an LDC with 
$\ox := +$, and $\oa := \circ$. The unit objects are given by the units of the respective multiplications. 
The linear distributor is given by the following equality:
\[ x \circ (y + z) := (x + (y + z)) - s = ((x + y ) -s) + z = (x \circ y) + z \]
Note that the distributors are identity maps but $\ox$ and $\oa$ are distinct.

\item $*$-autonomous categories are LDCs with a dualizing object. See Section \ref{Sec: *-autonomous}.

\item The category of bialgebra modules and module homomorphisms of a $*$-autonomous category is an LDC. 
The tensor products are inherited from the base category \cite{CS97}. The category of Hopf modules and 
module homomorphisms from a $*$-autonomous category is $*$-autonomous \cite{PaS09}. We discuss the category of 
Hopf Modules in section \ref{Sec: HModx}. 

\item Girard's Coherence spaces \cite{Gir87}, Ehrhard's finiteness spaces \cite{Ehr05}, 
and Chu Spaces \cite{Bar06} are linearly distributive categories that are also $*$-autonomous. 
 Indeed, $\FRel$, the category of finiteness spaces and finiteness relations, and $\FMat(\C)$, 
 the category of finiteness spaces and finiteness matrices are used as primary examples for 
 structures introduced in this thesis. See Section \ref{Sec: motivating examples} for discussion of these categories. 

\item  Bicompleting a monoidal category (adding arbitrary limits and colimits) gives an LDC. 
A procedure for bicompletion of monoidal categories has been described by Joyal in \cite{Joy95}.  
Joyal also proved in the same article that if the base category is compact closed, then the resulting category is $*$-autonomous.
\end{itemize}

\subsection{Graphical calculus}
\label{Sec: graphical}
LDCs come equipped with a graphical calculus \cite{BCST96} that contains 
the calculus for monoidal categories. Every sequent rule and derivation in MLL corresponds to a 
circuit in the graphical calculus of LDCs, and vice versa. In this section, we review the fundamentals of 
the graphical calculus for LDCs. For detailed exposition, see \cite{BCST96, CS97}.  
The following are the generators of LDC circuits: wires represent objects and circles represent maps. The input wires of a map are tensored (with $\ox$), and the output wires are ``par''ed (with $\oa$). The following diagram represents a map $f: A \ox B \to C \oa D$. 
 \[ \begin{tikzpicture}
	\begin{pgfonlayer}{nodelayer}
		\node [style=circle, scale=2] (0) at (-0.25, -0.5) {};
		\node [style=none] (1) at (-0.25, -0.5) {$f$};
		\node [style=none] (2) at (-1, 0.5) {$A$};
		\node [style=none] (3) at (0.5, 0.5) {$B$};
		\node [style=none] (4) at (-1, -1.5) {$C$};
		\node [style=none] (5) at (0.5, -1.5) {$D$};
		\node [style=ox] (6) at (-0.25, 1) {};
		\node [style=oa] (7) at (-0.25, -2) {};
		\node [style=none] (8) at (-0.25, -2.75) {};
		\node [style=none] (9) at (-0.25, 1.75) {};
		\node [style=none] (10) at (-0.5, -1.65) {};
		\node [style=none] (11) at (0, -1.65) {};
		\node [style=none] (12) at (0, 0.65) {};
		\node [style=none] (13) at (-0.5, 0.65) {};
		\node [style=none] (14) at (-0.25, -3) {$f: A \ox B \to C \oa D$};
	\end{pgfonlayer}
	\begin{pgfonlayer}{edgelayer}
		\draw [bend right=45, looseness=1.25] (6) to (0);
		\draw [bend left=45, looseness=1.25] (6) to (0);
		\draw [bend left=45, looseness=1.25] (7) to (0);
		\draw [bend right=45, looseness=1.25] (7) to (0);
		\draw (9.center) to (6);
		\draw (7) to (8.center);
	\end{pgfonlayer}
\end{tikzpicture}  \]

The $\ox$-associator, the $\oa$-associator, the left linear distributor, and the right linear distributors are, respectively, drawn as follows:
\[(a) ~~~ \begin{tikzpicture}
	\begin{pgfonlayer}{nodelayer}
		\node [style=ox] (0) at (1, -2.5) {};
		\node [style=ox] (1) at (0.5, -0) {};
		\node [style=ox] (2) at (1.75, -1.5) {};
		\node [style=ox] (3) at (1.5, 1.25) {};
		\node [style=none] (4) at (1.5, 2) {};
		\node [style=none] (5) at (1, -3.25) {};
		\node [style=none] (6) at (0.75, -0.5) {};
		\node [style=none] (7) at (0.25, -0.5) {};
		\node [style=none] (8) at (1.25, 1) {};
		\node [style=none] (9) at (1.75, 1) {};
		\node [style=none] (10) at (2.75, 1.75) {$A \ox (B \ox C)$};
		\node [style=none] (11) at (2.25, -3) {$(A \ox B) \ox C$};
		\node [style=none] (12) at (0.5, 1) {$A \ox B$};
		\node [style=none] (13) at (2.25, 0.75) {$C$};
		\node [style=none] (14) at (2.25, -2.25) {$B \ox C$};
		\node [style=none] (15) at (0, -2) {$A$};
	\end{pgfonlayer}
	\begin{pgfonlayer}{edgelayer}
		\draw (5.center) to (0);
		\draw [bend right, looseness=1.00] (0) to (2);
		\draw [in=-120, out=160, looseness=1.00] (0) to (1);
		\draw (3) to (4.center);
		\draw [in=-150, out=75, looseness=0.75] (1) to (3);
		\draw [in=-60, out=120, looseness=1.25] (2) to (1);
		\draw [in=-30, out=75, looseness=1.00] (2) to (3);
	\end{pgfonlayer}
\end{tikzpicture} ~~~~~~~~ (b)~~~~ \begin{tikzpicture}
	\begin{pgfonlayer}{nodelayer}
		\node [style=oa] (0) at (1.5, 0.75) {};
		\node [style=oa] (1) at (2, -1.75) {};
		\node [style=oa] (2) at (0.75, -0.25) {};
		\node [style=oa] (3) at (1, -3) {};
		\node [style=none] (4) at (1, -3.75) {};
		\node [style=none] (5) at (1.5, 1.5) {};
		\node [style=none] (6) at (1.8, -1.35) {};
		\node [style=none] (7) at (2.15, -1.35) {};
		\node [style=none] (8) at (1.35, -2.65) {};
		\node [style=none] (9) at (0.75, -2.65) {};
		\node [style=none] (10) at (2.25, -3.5) {$A \oa (B \oa C)$};
		\node [style=none] (11) at (2.75, 1.25) {$(A \oa B) \oa C$};
		\node [style=none] (12) at (0, 0.5) {$A \oa B$};
		\node [style=none] (13) at (2.25, -2.5) {$B \oa C$};
		\node [style=none] (14) at (1.25, -1.25) {$B$};
		\node [style=none] (15) at (2.5, -0) {$C$};
		\node [style=none] (16) at (0.25, -2.5) {$A$};
	\end{pgfonlayer}
	\begin{pgfonlayer}{edgelayer}
		\draw (5.center) to (0);
		\draw [bend right, looseness=1.00] (0) to (2);
		\draw [in=60, out=-20, looseness=1.00] (0) to (1);
		\draw (3) to (4.center);
		\draw [in=30, out=-105, looseness=0.75] (1) to (3);
		\draw [in=120, out=-60, looseness=1.25] (2) to (1);
		\draw [in=135, out=-130, looseness=1.00] (2) to (3);
	\end{pgfonlayer}
\end{tikzpicture} ~~~~~~~~ (c) ~~~~ \begin{tikzpicture}
	\begin{pgfonlayer}{nodelayer}
		\node [style=ox] (0) at (1.25, 0.75) {};
		\node [style=ox] (1) at (0.5, -1.25) {};
		\node [style=oa] (2) at (2, -0.25) {};
		\node [style=oa] (3) at (1.25, -2.5) {};
		\node [style=none] (4) at (1.25, -3.25) {};
		\node [style=none] (5) at (1.25, 1.5) {};
		\node [style=none] (6) at (1.5, -2.25) {};
		\node [style=none] (7) at (1, -2.25) {};
		\node [style=none] (8) at (1, 0.5) {};
		\node [style=none] (9) at (1.5, 0.5) {};
		\node [style=none] (10) at (2.5, 1.25) {$A \ox (B \oa C)$};
		\node [style=none] (11) at (0.25, -0) {$A$};
		\node [style=none] (12) at (2.55, 0.5) {$B \oa C$};
		\node [style=none] (13) at (1.25, -0.5) {$B$};
		\node [style=none] (14) at (2.5, -1) {$C$};
		\node [style=none] (15) at (0, -2) {$A \ox B$};
		\node [style=none] (16) at (2.55, -3.25) {$(A \ox B) \oa C$};
	\end{pgfonlayer}
	\begin{pgfonlayer}{edgelayer}
		\draw (5.center) to (0);
		\draw [bend left, looseness=1.00] (0) to (2);
		\draw [in=90, out=-150, looseness=1.00] (0) to (1);
		\draw (3) to (4.center);
		\draw [in=150, out=-75, looseness=1.00] (1) to (3);
		\draw [in=60, out=-120, looseness=1.25] (2) to (1);
		\draw [in=30, out=-75, looseness=1.00] (2) to (3);
	\end{pgfonlayer}
\end{tikzpicture} ~~~~~~~~ (d) ~~~~\begin{tikzpicture}
	\begin{pgfonlayer}{nodelayer}
		\node [style=ox] (0) at (1.3, 0.75) {};
		\node [style=ox] (1) at (2.05, -1.25) {};
		\node [style=oa] (2) at (0.55, -0.25) {};
		\node [style=oa] (3) at (1.3, -2.5) {};
		\node [style=none] (4) at (1.3, -3.25) {};
		\node [style=none] (5) at (1.3, 1.5) {};
		\node [style=none] (6) at (1.05, -2.25) {};
		\node [style=none] (7) at (1.55, -2.25) {};
		\node [style=none] (8) at (1.55, 0.5) {};
		\node [style=none] (9) at (1.05, 0.5) {};
		\node [style=none] (10) at (0, 1.25) {$(A \oa B) \ox C$};
		\node [style=none] (11) at (2.3, -0) {$C$};
		\node [style=none] (12) at (0, 0.5) {$A \oa B$};
		\node [style=none] (13) at (1.3, -0.5) {$B$};
		\node [style=none] (14) at (0.04999995, -1) {$A$};
		\node [style=none] (15) at (2.55, -2) {$B \ox C$};
		\node [style=none] (16) at (0, -3.25) {$A \oa (B \ox C)$};
	\end{pgfonlayer}
	\begin{pgfonlayer}{edgelayer}
		\draw (5.center) to (0);
		\draw [bend right, looseness=1.00] (0) to (2);
		\draw [in=90, out=-30, looseness=1.00] (0) to (1);
		\draw (3) to (4.center);
		\draw [in=30, out=-105, looseness=1.00] (1) to (3);
		\draw [in=120, out=-60, looseness=1.25] (2) to (1);
		\draw [in=150, out=-105, looseness=1.00] (2) to (3);
	\end{pgfonlayer}
\end{tikzpicture}  \]

$\begin{tikzpicture}
	\begin{pgfonlayer}{nodelayer}
		\node [style=oa] (0) at (1, -3) {};
		\node [style=none] (1) at (1, -3.75) {};
		\node [style=none] (2) at (1.25, -2.75) {};
		\node [style=none] (3) at (0.75, -2.75) {};
		\node [style=none] (7) at (0.5, -2) {};
		\node [style=none] (8) at (1.5, -2) {};
	\end{pgfonlayer}
	\begin{pgfonlayer}{edgelayer}
		\draw (0) to (1.center);
		\draw [in=-90, out=150, looseness=1.00] (0) to (7.center);
		\draw [in=-90, out=30, looseness=1.00] (0) to (8.center);
	\end{pgfonlayer}
\end{tikzpicture}$ is the $\oa$-introduction ($\oa I$) rule, $\begin{tikzpicture}
	\begin{pgfonlayer}{nodelayer}
		\node [style=ox] (0) at (1, -3) {};
		\node [style=none] (1) at (1, -3.75) {};
		\node [style=none] (2) at (1.25, -2.75) {};
		\node [style=none] (3) at (0.75, -2.75) {};
		\node [style=none] (7) at (0.5, -2) {};
		\node [style=none] (8) at (1.5, -2) {};
	\end{pgfonlayer}
	\begin{pgfonlayer}{edgelayer}
		\draw (0) to (1.center);
		\draw [in=-90, out=150, looseness=1.00] (0) to (7.center);
		\draw [in=-90, out=30, looseness=1.00] (0) to (8.center);
	\end{pgfonlayer}
\end{tikzpicture}$ is $\ox$-introduction ($\ox I$) rule, $\begin{tikzpicture}
	\begin{pgfonlayer}{nodelayer}
		\node [style=ox] (0) at (1, -2.75) {};
		\node [style=none] (1) at (1, -2) {};
		\node [style=none] (2) at (1.25, -3) {};
		\node [style=none] (3) at (0.75, -3) {};
		\node [style=none] (7) at (0.5, -3.75) {};
		\node [style=none] (8) at (1.5, -3.75) {};
	\end{pgfonlayer}
	\begin{pgfonlayer}{edgelayer}
		\draw (0) to (1.center);
		\draw [in=90, out=-150, looseness=1.00] (0) to (7.center);
		\draw [in=90, out=-30, looseness=1.00] (0) to (8.center);
	\end{pgfonlayer}
\end{tikzpicture}$ is the $\ox$-elimination ($\ox E$) rule, $\begin{tikzpicture}
	\begin{pgfonlayer}{nodelayer}
		\node [style=oa] (0) at (1, -2.75) {};
		\node [style=none] (1) at (1, -2) {};
		\node [style=none] (2) at (0.5, -3.75) {};
		\node [style=none] (3) at (1.5, -3.75) {};
	\end{pgfonlayer}
	\begin{pgfonlayer}{edgelayer}
		\draw (0) to (1.center);
		\draw [in=90, out=-150, looseness=1.00] (0) to (2.center);
		\draw [in=90, out=-30, looseness=1.00] (0) to (3.center);
	\end{pgfonlayer}
\end{tikzpicture}$ is the $\oa$-elimination ($\oa E$) rule. As shown below, the rules, $(\ox I)$ and $(\ox E)$, 
correspond to the sequent rules for tensor introduction, $\ox L$ and $\ox R$ in Figure \ref{Fig: multiplicative rules}: 
\[ \infer[(\ox L)]{\Gamma, A \ox B, \Gamma' \vdash \Delta }{\Gamma, A, B, \Gamma' \vdash \Delta} 
~~~~~~~~
 \infer[(\ox R)]{\Gamma_1, \Delta_1 \vdash \Gamma_2, \Delta_2, A \ox B, \Gamma_3, \Delta_3 }
{\Gamma_1  \vdash \Gamma_2, A, \Gamma_3 & \Delta_1  \vdash \Delta_2, B, \Delta_3}\]
Similarly, $\oa$I and $\oa$E correspond to $\oa$L and $\oa$R respectively, see Figure \ref{Fig: multiplicative rules}.

The unitors are drawn as follows:
\[ (a) ~~~~ \begin{tikzpicture}
	\begin{pgfonlayer}{nodelayer}
		\node [style=circle] (0) at (0, -0) {$\top$};
		\node [style=none] (1) at (0, -2) {};
		\node [style=none] (2) at (0.75, -0) {};
		\node [style=none] (3) at (0.75, -2) {};
		\node [style=none] (4) at (1, -1) {$A$};
		\node [style=none] (5) at (0, -2.5) {$(u_\ox^L)^{-1}: A \to \top \ox A$};
	\end{pgfonlayer}
	\begin{pgfonlayer}{edgelayer}
		\draw (0) to (1.center);
	\end{pgfonlayer}
\end{tikzpicture} ~~~~~~~~ (b) ~~~~ \begin{tikzpicture}
	\begin{pgfonlayer}{nodelayer}
		\node [style=circle] (0) at (0, -0) {$\top$};
		\node [style=none] (1) at (0.75, 0.75) {};
		\node [style=none] (2) at (0.75, -2) {};
		\node [style=none] (3) at (1, -0.25) {$A$};
		\node [style=none] (4) at (0.25, -2.5) {$u_\ox^L: \top \ox A \to A$};
		\node [style=circle, scale=0.6] (5) at (0.75, -1.25) {};
		\node [style=none] (6) at (0, 0.75) {};
	\end{pgfonlayer}
	\begin{pgfonlayer}{edgelayer}
		\draw (2.center) to (1.center);
		\draw [dotted, bend right, looseness=1.25] (0) to (5);
		\draw (6.center) to (0);
	\end{pgfonlayer}
\end{tikzpicture} ~~~~~~~~ (c) ~~~ \begin{tikzpicture}
	\begin{pgfonlayer}{nodelayer}
		\node [style=circle] (0) at (0, -2.5) {$\bot$};
		\node [style=none] (1) at (0.75, -3.5) {};
		\node [style=none] (2) at (0.75, -0.5) {};
		\node [style=none] (3) at (1, -0.75) {$A$};
		\node [style=none] (4) at (0.5, -3.75) {$(u_\oa^L)^{-1}:  A \to \bot \oa A$};
		\node [style=circle, scale=0.6] (5) at (0.75, -1.5) {};
		\node [style=none] (6) at (0, -3.5) {};
	\end{pgfonlayer}
	\begin{pgfonlayer}{edgelayer}
		\draw (1.center) to (2.center);
		\draw [dotted, bend left, looseness=1.25, dotted] (0) to (5);
		\draw (6.center) to (0);
	\end{pgfonlayer}
\end{tikzpicture} ~~~~~~~ (d) ~~~ \begin{tikzpicture}
	\begin{pgfonlayer}{nodelayer}
		\node [style=circle] (0) at (0, -2.75) {$\bot$};
		\node [style=none] (1) at (0.75, -0.75) {};
		\node [style=none] (2) at (0.75, -3) {};
		\node [style=none] (3) at (1, -2.75) {$A$};
		\node [style=none] (4) at (0.25, -3.5) {$u_\oa^L:  \bot \oa A \to A$};
		\node [style=none] (5) at (0, -0.75) {};
	\end{pgfonlayer}
	\begin{pgfonlayer}{edgelayer}
		\draw (5.center) to (0);
	\end{pgfonlayer}
\end{tikzpicture}  \]

Diagram $(a)$ is called the left $\top$-introduction, $(b)$ is called the left $\top$-elimination, 
$(c)$ is the left $\bot$-introduction, and $(d)$ is the left $\bot$-elimination which correspond 
to the sequent rules $(\top R)$,  $(\top L)$,  $(\bot R)$, $(\bot L)$ in Figure \ref{Fig: multiplicative rules} respectively. 
The unit $\top$ is introduced, and the counit $\bot$ is eliminated using the thinning links 
which are shown using dotted wires in the diagrams. See \cite[Section 2.3]{BCST96} for details 
on the thinning links.

The following are a set of circuit equalities (which when oriented become reduction rewrite rules):
\[ [Reduction]: ~~~ \begin{tikzpicture}
	\begin{pgfonlayer}{nodelayer} 
		\node [style=circle] (0) at (0, -1) {$\top$};
		\node [style=none] (1) at (0.75, 0.25) {};
		\node [style=none] (2) at (0.75, -3) {};
		\node [style=none] (3) at (1, -1.25) {$A$};
		\node [style=circle, scale = 0.4] (4) at (0.75, -2.25) {};
		\node [style=circle] (5) at (0, -0) {$\top$};
	\end{pgfonlayer}
	\begin{pgfonlayer}{edgelayer}
		\draw (1.center) to (2.center);
		\draw [dotted, bend right, looseness=1.25] (0) to (4);
		\draw (5) to (0);
	\end{pgfonlayer}
\end{tikzpicture} =  \begin{tikzpicture}
	\begin{pgfonlayer}{nodelayer}
		\node [style=none] (0) at (2.25, -1.5) {$A$};
		\node [style=none] (1) at (2, 0.25) {};
		\node [style=none] (2) at (2, -3) {};
	\end{pgfonlayer}
	\begin{pgfonlayer}{edgelayer}
		\draw (1.center) to (2.center);
	\end{pgfonlayer}
\end{tikzpicture}
~~~~~~~~
\begin{tikzpicture}
	\begin{pgfonlayer}{nodelayer}
		\node [style=circle] (0) at (0, -1.75) {$\bot$};
		\node [style=none] (1) at (0.75, -3) {};
		\node [style=none] (2) at (0.75, 0.25) {};
		\node [style=none] (3) at (1, -1.5) {$A$};
		\node [style=circle, scale=0.4] (4) at (0.75, -0.5) {};
		\node [style=circle] (5) at (0, -2.75) {$\bot$};
	\end{pgfonlayer}
	\begin{pgfonlayer}{edgelayer}
		\draw (1.center) to (2.center);
		\draw [dotted, bend left, looseness=1.25] (0) to (4);
		\draw (5) to (0);
	\end{pgfonlayer}
\end{tikzpicture} = \begin{tikzpicture}
	\begin{pgfonlayer}{nodelayer}
		\node [style=none] (0) at (2.25, -1.5) {$A$};
		\node [style=none] (1) at (2, 0.25) {};
		\node [style=none] (2) at (2, -3) {};
	\end{pgfonlayer}
	\begin{pgfonlayer}{edgelayer}
		\draw (1.center) to (2.center);
	\end{pgfonlayer}
\end{tikzpicture}
~~~~~~~~
\begin{tikzpicture}
	\begin{pgfonlayer}{nodelayer}
		\node [style=ox] (0) at (2, -0) {};
		\node [style=ox] (1) at (2, -1) {};
		\node [style=none] (2) at (1.5, -2) {};
		\node [style=none] (3) at (2.5, -2) {};
		\node [style=none] (4) at (1.5, 1) {};
		\node [style=none] (5) at (2.5, 1) {};
		\node [style=none] (6) at (1.25, 0.75) {$A$};
		\node [style=none] (7) at (2.75, 0.75) {$B$};
		\node [style=none] (8) at (1.25, -1.75) {$A$};
		\node [style=none] (9) at (2.75, -1.75) {$B$};
	\end{pgfonlayer}
	\begin{pgfonlayer}{edgelayer}
		\draw (0) to (1);
		\draw [in=90, out=-135, looseness=1.00] (1) to (2.center);
		\draw [in=90, out=-45, looseness=1.00] (1) to (3.center);
		\draw [in=-90, out=135, looseness=1.00] (0) to (4.center);
		\draw [in=-90, out=45, looseness=1.00] (0) to (5.center);
	\end{pgfonlayer}
\end{tikzpicture} = \begin{tikzpicture}
	\begin{pgfonlayer}{nodelayer}
		\node [style=none] (0) at (1.5, -2) {};
		\node [style=none] (1) at (2, -2) {};
		\node [style=none] (2) at (1.5, 1) {};
		\node [style=none] (3) at (2, 1) {};
		\node [style=none] (4) at (1.25, 0.75) {$A$};
		\node [style=none] (5) at (2.25, 0.75) {$B$};
		\node [style=none] (6) at (1.25, -1.75) {$A$};
		\node [style=none] (7) at (2.25, -1.75) {$B$};
	\end{pgfonlayer}
	\begin{pgfonlayer}{edgelayer}
		\draw (2.center) to (0.center);
		\draw (3.center) to (1.center);
	\end{pgfonlayer}
\end{tikzpicture}
~~~~~~~~
\begin{tikzpicture}
	\begin{pgfonlayer}{nodelayer}
		\node [style=oa] (0) at (2, -0) {};
		\node [style=oa] (1) at (2, -1) {};
		\node [style=none] (2) at (1.5, -2) {};
		\node [style=none] (3) at (2.5, -2) {};
		\node [style=none] (4) at (1.5, 1) {};
		\node [style=none] (5) at (2.5, 1) {};
		\node [style=none] (6) at (1.25, 0.75) {$A$};
		\node [style=none] (7) at (2.75, 0.75) {$B$};
		\node [style=none] (8) at (1.25, -1.75) {$A$};
		\node [style=none] (9) at (2.75, -1.75) {$B$};
	\end{pgfonlayer}
	\begin{pgfonlayer}{edgelayer}
		\draw (0) to (1);
		\draw [in=90, out=-135, looseness=1.00] (1) to (2.center);
		\draw [in=90, out=-45, looseness=1.00] (1) to (3.center);
		\draw [in=-90, out=135, looseness=1.00] (0) to (4.center);
		\draw [in=-90, out=45, looseness=1.00] (0) to (5.center);
	\end{pgfonlayer}
\end{tikzpicture} = \begin{tikzpicture}
	\begin{pgfonlayer}{nodelayer}
		\node [style=none] (0) at (1.5, -2) {};
		\node [style=none] (1) at (2, -2) {};
		\node [style=none] (2) at (1.5, 1) {};
		\node [style=none] (3) at (2, 1) {};
		\node [style=none] (4) at (1.25, 0.75) {$A$};
		\node [style=none] (5) at (2.25, 0.75) {$B$};
		\node [style=none] (6) at (1.25, -1.75) {$A$};
		\node [style=none] (7) at (2.25, -1.75) {$B$};
	\end{pgfonlayer}
	\begin{pgfonlayer}{edgelayer}
		\draw (2.center) to (0.center);
		\draw (3.center) to (1.center);
	\end{pgfonlayer}
\end{tikzpicture} \]
The following are also circuit equalities (and when oriented become expansion rules:)
\[ [Expansion]: ~~~ \begin{tikzpicture}
	\begin{pgfonlayer}{nodelayer}
		\node [style=ox] (0) at (0, -0.75) {};
		\node [style=ox] (1) at (0, -2) {};
		\node [style=none] (2) at (0, -3) {};
		\node [style=none] (3) at (0, -0) {};
	\end{pgfonlayer}
	\begin{pgfonlayer}{edgelayer}
		\draw [bend right=60, looseness=1.50] (0) to (1);
		\draw [bend right=60, looseness=1.50] (1) to (0);
		\draw (3.center) to (0);
		\draw (1) to (2.center);
	\end{pgfonlayer}
\end{tikzpicture} = \begin{tikzpicture}
	\begin{pgfonlayer}{nodelayer}
		\node [style=none] (0) at (1, -0) {};
		\node [style=none] (1) at (1, -3) {};
		\node [style=none] (2) at (1.65, -1.75) {$A \ox B$};
	\end{pgfonlayer}
	\begin{pgfonlayer}{edgelayer}
		\draw (0.center) to (1.center);
	\end{pgfonlayer}
\end{tikzpicture}
~~~~~~~~
\begin{tikzpicture}
	\begin{pgfonlayer}{nodelayer}
		\node [style=oa] (0) at (0, -0.75) {};
		\node [style=oa] (1) at (0, -2) {};
		\node [style=none] (2) at (0, -3) {};
		\node [style=none] (3) at (0, -0) {};
	\end{pgfonlayer}
	\begin{pgfonlayer}{edgelayer}
		\draw [bend right=60, looseness=1.50] (0) to (1);
		\draw [bend right=60, looseness=1.50] (1) to (0);
		\draw (3.center) to (0);
		\draw (1) to (2.center);
	\end{pgfonlayer}
\end{tikzpicture} = \begin{tikzpicture}
	\begin{pgfonlayer}{nodelayer}
		\node [style=none] (0) at (1, -0) {};
		\node [style=none] (1) at (1, -3) {};
		\node [style=none] (2) at (1.65, -1.75) {$A \oa B$};
	\end{pgfonlayer}
	\begin{pgfonlayer}{edgelayer}
		\draw (0.center) to (1.center);
	\end{pgfonlayer}
\end{tikzpicture}
~~~~~~~~~
\begin{tikzpicture}
	\begin{pgfonlayer}{nodelayer}
		\node [style=circle] (0) at (0, -0.75) {$\top$};
		\node [style=circle] (1) at (-1, -1.5) {$\top$};
		\node [style=none] (2) at (-1, -3.25) {};
		\node [style=none] (3) at (0, 0) {};
		\node [style=circle, scale=0.4] (4) at (-1, -2.5) {};
	\end{pgfonlayer}
	\begin{pgfonlayer}{edgelayer}
		\draw (3.center) to (0);
		\draw (1) to (2.center);
		\draw [dotted, in=-90, out=30, looseness=1.25] (4) to (0);
	\end{pgfonlayer}
\end{tikzpicture} = \begin{tikzpicture}
	\begin{pgfonlayer}{nodelayer}
		\node [style=none] (0) at (1, -0) {};
		\node [style=none] (1) at (1, -3.25) {};
		\node [style=none] (2) at (1.25, -1.75) {$\top$};
	\end{pgfonlayer}
	\begin{pgfonlayer}{edgelayer}
		\draw (0.center) to (1.center);
	\end{pgfonlayer}
\end{tikzpicture}
~~~~~~~~~
\begin{tikzpicture}
	\begin{pgfonlayer}{nodelayer}
		\node [style=circle] (0) at (0, -2.5) {$\bot$};
		\node [style=circle] (1) at (-1, -1.75) {$\bot$};
		\node [style=none] (2) at (-1, 0) {};
		\node [style=none] (3) at (0, -3.25) {};
		\node [style=circle, scale=0.4] (4) at (-1, -0.75) {};
	\end{pgfonlayer}
	\begin{pgfonlayer}{edgelayer}
		\draw (3.center) to (0);
		\draw (1) to (2.center);
		\draw [dotted, in=90, out=-30, looseness=1.25] (4) to (0);
	\end{pgfonlayer}
\end{tikzpicture} = \begin{tikzpicture}
	\begin{pgfonlayer}{nodelayer}
		\node [style=none] (0) at (1, -0) {};
		\node [style=none] (1) at (1, -3.25) {};
		\node [style=none] (2) at (1.25, -1.75) {$\bot$};
	\end{pgfonlayer}
	\begin{pgfonlayer}{edgelayer}
		\draw (0.center) to (1.center);
	\end{pgfonlayer}
\end{tikzpicture} \]


As in linear logic, not all circuit diagrams constructed from these basic components represent a valid LDC circuit.  In his seminal paper on linear logic, \cite{Gir87}, Girard introduced a criterion for the correctness of his representation of proofs using proof nets based on switching links.  A valid proof structure must be connected and acyclic for all the switching link choices.  Using this correctness criterion has the disadvantage of requiring exponential time in the number of switching links.  Danos and Regnier \cite{DaR89} improved this situation significantly by providing an algorithm for correctness which takes linear time (see \cite{Gue99})
on the size of the circuit.  To verify the validity of the circuit diagrams of LDCs, Blute et.al. \cite{BCST96}, provided a boxing algorithm which was based on Danos and Regnier's more efficient algorithm which we now describe.

In order to verify that an LDC circuit is valid, circuit components are 
``boxed''  using the rules below. The primitive generating maps are 
automatically boxed. 
\[ (a_1)~~~
 \]

Double lines refer to multiple number of wires. The boxes contain circuit components including maps. $\ox$-introduction and $\oa$-elimination are boxed in $(a_1)$ and $(a_2)$ respectively. 
In $(b_1)$, it is shown how a box `eats' the $\ox$-elimination: in $(b_2)$ the dual rule shows a $\oa$-introduction being eaten. 
$(c)$ shows how boxes can be amalgamated when they are connected by a single wire. $\bot$-elimination, 
$\top$-introduction, and identity maps are boxed in $(d_1)$, $(d_2)$, and $(d_3)$ respectively. In $(e_1)$-$(e_4)$, 
it is shown how the thinning links can be boxed. By progressively enclosing the components of the circuit in boxes using these rules, 
if we end up with a single box (or a wire), precisely when the circuit is valid. As an example, we verify the validity of the left linear distributor:
\[ \begin{tikzpicture}
	\begin{pgfonlayer}{nodelayer}
		\node [style=ox] (0) at (1.25, 0.75) {};
		\node [style=ox] (1) at (0.5, -1.25) {};
		\node [style=oa] (2) at (2, -0.25) {};
		\node [style=oa] (3) at (1.25, -2.5) {};
		\node [style=none] (4) at (1.25, -3.25) {};
		\node [style=none] (5) at (1.25, 1.5) {};
	\end{pgfonlayer}
	\begin{pgfonlayer}{edgelayer}
		\draw (5.center) to (0);
		\draw [bend left, looseness=1.00] (0) to (2);
		\draw [in=90, out=-150, looseness=1.00] (0) to (1);
		\draw (3) to (4.center);
		\draw [in=150, out=-75, looseness=1.00] (1) to (3);
		\draw [in=60, out=-120, looseness=1.25] (2) to (1);
		\draw [in=30, out=-75, looseness=1.00] (2) to (3);
	\end{pgfonlayer}
\end{tikzpicture} \stackrel{a_1,a_2}{\Rightarrow} \begin{tikzpicture}
	\begin{pgfonlayer}{nodelayer}
		\node [style=ox] (0) at (1.25, 0.75) {};
		\node [style=ox] (1) at (0.5, -1.25) {};
		\node [style=oa] (2) at (2, -0.25) {};
		\node [style=oa] (3) at (1.25, -2.5) {};
		\node [style=none] (4) at (1.25, -3.25) {};
		\node [style=none] (5) at (1.25, 1.5) {};
		\node [style=none] (6) at (1.5, 0.25) {};
		\node [style=none] (7) at (2.5, 0.25) {};
		\node [style=none] (8) at (2.5, -1) {};
		\node [style=none] (9) at (1.5, -1) {};
		\node [style=none] (10) at (1.25, -0.5) {};
		\node [style=none] (11) at (0, -0.5) {};
		\node [style=none] (12) at (0, -1.75) {};
		\node [style=none] (13) at (1.25, -1.75) {};
	\end{pgfonlayer}
	\begin{pgfonlayer}{edgelayer}
		\draw (5.center) to (0);
		\draw [bend left, looseness=1.00] (0) to (2);
		\draw [in=90, out=-150, looseness=1.00] (0) to (1);
		\draw (3) to (4.center);
		\draw [in=150, out=-75, looseness=1.00] (1) to (3);
		\draw [in=60, out=-120, looseness=1.25] (2) to (1);
		\draw [in=30, out=-75, looseness=1.00] (2) to (3);
		\draw (9.center) to (8.center);
		\draw (8.center) to (7.center);
		\draw (7.center) to (6.center);
		\draw (6.center) to (9.center);
		\draw (11.center) to (12.center);
		\draw (12.center) to (13.center);
		\draw (13.center) to (10.center);
		\draw (10.center) to (11.center);
	\end{pgfonlayer}
\end{tikzpicture} \stackrel{c}{\Rightarrow} \begin{tikzpicture}
	\begin{pgfonlayer}{nodelayer}
		\node [style=ox] (0) at (1.25, 0.75) {};
		\node [style=ox] (1) at (0.5, -1.25) {};
		\node [style=oa] (2) at (2, -0.25) {};
		\node [style=oa] (3) at (1.25, -2.5) {};
		\node [style=none] (4) at (1.25, -3.25) {};
		\node [style=none] (5) at (1.25, 1.5) {};
		\node [style=none] (6) at (2.5, 0.25) {};
		\node [style=none] (7) at (2.5, -1.75) {};
		\node [style=none] (8) at (0, 0.25) {};
		\node [style=none] (9) at (0, -1.75) {};
	\end{pgfonlayer}
	\begin{pgfonlayer}{edgelayer}
		\draw (5.center) to (0);
		\draw [bend left, looseness=1.00] (0) to (2);
		\draw [in=90, out=-150, looseness=1.00] (0) to (1);
		\draw (3) to (4.center);
		\draw [in=150, out=-75, looseness=1.00] (1) to (3);
		\draw [in=60, out=-120, looseness=1.25] (2) to (1);
		\draw [in=30, out=-75, looseness=1.00] (2) to (3);
		\draw (7.center) to (6.center);
		\draw (8.center) to (9.center);
		\draw (8.center) to (6.center);
		\draw (9.center) to (7.center);
	\end{pgfonlayer}
\end{tikzpicture} \stackrel{b_1,b_2}{\Rightarrow} \begin{tikzpicture}
	\begin{pgfonlayer}{nodelayer}
		\node [style=ox] (0) at (1.25, 0.75) {};
		\node [style=ox] (1) at (0.5, -1.25) {};
		\node [style=oa] (2) at (2, -0.25) {};
		\node [style=oa] (3) at (1.25, -2.5) {};
		\node [style=none] (4) at (1.25, -3.25) {};
		\node [style=none] (5) at (1.25, 1.5) {};
		\node [style=none] (6) at (2.5, 1.25) {};
		\node [style=none] (7) at (2.5, -3) {};
		\node [style=none] (8) at (0, 1.25) {};
		\node [style=none] (9) at (0, -3) {};
	\end{pgfonlayer}
	\begin{pgfonlayer}{edgelayer}
		\draw (5.center) to (0);
		\draw [bend left, looseness=1.00] (0) to (2);
		\draw [in=90, out=-150, looseness=1.00] (0) to (1);
		\draw (3) to (4.center);
		\draw [in=150, out=-75, looseness=1.00] (1) to (3);
		\draw [in=60, out=-120, looseness=1.25] (2) to (1);
		\draw [in=30, out=-75, looseness=1.00] (2) to (3);
		\draw (7.center) to (6.center);
		\draw (8.center) to (9.center);
		\draw (8.center) to (6.center);
		\draw (9.center) to (7.center);
	\end{pgfonlayer}
\end{tikzpicture} \]
In the first step the $\ox$-introduction and $\oa$-elimination are boxed. In the second step the boxes are amalgamated along the 
single wire joining them. In the third step,  the box absorbs the $\ox$-elimination and $\oa$-introduction.

In contrast, we now show that the reverse of the linear distributor is invalid as the boxing process gets stuck (there are 
no rules to box $\ox$-elimination and $\oa$-introduction):
\[ \begin{tikzpicture}
	\begin{pgfonlayer}{nodelayer}
		\node [style=ox] (0) at (1.25, -2.5) {};
		\node [style=ox] (1) at (2, -0.25) {};
		\node [style=oa] (2) at (0.5, -1.25) {};
		\node [style=oa] (3) at (1.25, 0.75) {};
		\node [style=none] (4) at (1.25, 1.5) {};
		\node [style=none] (5) at (1.25, -3.25) {};
	\end{pgfonlayer}
	\begin{pgfonlayer}{edgelayer}
		\draw (5.center) to (0);
		\draw [in=-83, out=157, looseness=1.00] (0) to (2);
		\draw [in=-90, out=30, looseness=1.00] (0) to (1);
		\draw (3) to (4.center);
		\draw [in=-15, out=90, looseness=1.00] (1) to (3);
		\draw [in=-120, out=60, looseness=1.25] (2) to (1);
		\draw [in=-150, out=90, looseness=1.00] (2) to (3);
	\end{pgfonlayer}
\end{tikzpicture} \stackrel{a_1,a_2}{\Rightarrow} \begin{tikzpicture}
	\begin{pgfonlayer}{nodelayer}
		\node [style=ox] (0) at (1.25, -2.5) {};
		\node [style=ox] (1) at (2, -0.5) {};
		\node [style=oa] (2) at (0.5, -1.5) {};
		\node [style=oa] (3) at (1.25, 0.75) {};
		\node [style=none] (4) at (1.25, 1.5) {};
		\node [style=none] (5) at (1.25, -3.25) {};
		\node [style=none] (6) at (0.5, 1.25) {};
		\node [style=none] (7) at (2, 1.25) {};
		\node [style=none] (8) at (2, 0.25) {};
		\node [style=none] (9) at (0.5, 0.25) {};
		\node [style=none] (10) at (0.5, -2) {};
		\node [style=none] (11) at (2, -2) {};
		\node [style=none] (12) at (2, -3) {};
		\node [style=none] (13) at (0.5, -3) {};
	\end{pgfonlayer}
	\begin{pgfonlayer}{edgelayer}
		\draw (5.center) to (0);
		\draw [in=-83, out=157, looseness=1.00] (0) to (2);
		\draw [in=-90, out=30, looseness=1.00] (0) to (1);
		\draw (3) to (4.center);
		\draw [in=-15, out=90, looseness=1.00] (1) to (3);
		\draw [in=-120, out=60, looseness=1.25] (2) to (1);
		\draw [in=-150, out=90, looseness=1.00] (2) to (3);
		\draw (9.center) to (6.center);
		\draw (6.center) to (7.center);
		\draw (7.center) to (8.center);
		\draw (9.center) to (8.center);
		\draw (12.center) to (13.center);
		\draw (13.center) to (10.center);
		\draw (10.center) to (11.center);
		\draw (11.center) to (12.center);
	\end{pgfonlayer}
\end{tikzpicture}  \] 


\subsection{Mix, isomix and compact LDCs}
\label{Sec: mix, isomix, compact LDC}

In this thesis, we are predominately concerned with LDCs which have a mix map:

\begin{definition} \cite{CS97a}
	\label{Defn: mix cat}
A {\bf mix category} is an LDC with a {\bf mix map} ${\sf m}:\bot\to\top$ such that:
\[
\xymatrixcolsep{4pc}
\xymatrix{
A \ox B \ar[r]^{1 \ox u_\oa^{L^{-1}}} \ar[d]_{(u_\oa^R)^{-1} \ox 1} \ar@{.>}[ddrr]^{\mx_{A,B}} & A \ox (\bot \oa B) 
\ar[r]^{1 \ox (\m \oa 1)} & A \ox ( \top \oa B) \ar[d]^{\partial^L} \\
(A \oa \bot) \ox B \ar[d]_{\partial^R} & & ( A \ox \top ) \oa B  \ar[d]^{u_\ox^R \oa 1} \\
A \oa (\bot \ox B) \ar[r]_{1 \oa (\m \ox 1)} & A \oa (\top \ox B) \ar[r]_{1 \oa u_\ox^L} &  A \oa B
}
\]
\end{definition}

The map $\mx_{A,B}$ is a natural transformation and is called the {\bf mixor}. The coherence condition for the mix map has the following form in string diagrams (where the mix map is represented by an empty box):

{ \centering $\mx_{A,B}:=
\begin{tikzpicture}
	\begin{pgfonlayer}{nodelayer}
		\node [style=ox] (0) at (0, 0.2500001) {};
		\node [style=circ] (1) at (0.5000001, -0.2500001) {};
		\node [style=circ] (2) at (0, -1) {$\bot$};
		\node [style=map] (3) at (0, -1.75) {~};
		\node [style=circ] (4) at (0, -2.5) {$\top$};
		\node [style=circ] (5) at (-0.5000001, -3.25) {};
		\node [style=oa] (6) at (0, -3.75) {};
		\node [style=nothing] (7) at (0, 0.7499999) {};
		\node [style=nothing] (8) at (0, -4.25) {};
	\end{pgfonlayer}
	\begin{pgfonlayer}{edgelayer}
		\draw (7) to (0);
		\draw (0) to (1);
		\draw [in=45, out=-60, looseness=1.00] (1) to (6);
		\draw [in=120, out=-135, looseness=1.00] (0) to (5);
		\draw (5) to (6);
		\draw (6) to (8);
		\draw [densely dotted, in=-90, out=45, looseness=1.00] (5) to (4);
		\draw (4) to (3);
		\draw (3) to (2);
		\draw [densely dotted, in=-135, out=90, looseness=1.00] (2) to (1);
	\end{pgfonlayer}
\end{tikzpicture}
=
\begin{tikzpicture}
	\begin{pgfonlayer}{nodelayer}
		\node [style=circ] (0) at (-0.5000001, -0.2500001) {};
		\node [style=circ] (1) at (0, -1) {$\bot$};
		\node [style=map] (2) at (0, -1.75) {~};
		\node [style=circ] (3) at (0, -2.5) {$\top$};
		\node [style=circ] (4) at (0.5000001, -3.25) {};
		\node [style=nothing] (5) at (0, 0.7499999) {};
		\node [style=nothing] (6) at (0, -4.25) {};
		\node [style=oa] (7) at (0, -3.75) {};
		\node [style=ox] (8) at (0, 0.2500001) {};
	\end{pgfonlayer}
	\begin{pgfonlayer}{edgelayer}
		\draw [densely dotted, in=-90, out=150, looseness=1.00] (4) to (3);
		\draw (3) to (2);
		\draw (2) to (1);
		\draw [densely dotted, in=-45, out=90, looseness=1.00] (1) to (0);
		\draw (8) to (5);
		\draw (8) to (0);
		\draw [in=135, out=-120, looseness=1.00] (0) to (7);
		\draw (7) to (6);
		\draw (7) to (4);
		\draw [in=-45, out=60, looseness=1.00] (4) to (8);
	\end{pgfonlayer}
\end{tikzpicture} $ \par}

In a mix category, the associator, the distributor and the mix maps interact as follows. See Lemma 2, and proposition 3 in \cite{BCS00} for a proof.
\begin{equation*}
\mbox{\bf [mix.]}~~~~~~~ \xymatrix{
(A \oa B) \ox C \ar[r]^{\delta^R} \ar[d]_{\mx} \ar@{}[dr]|{(a)} &  A \oa (B \ox C) \ar[d]^{1 \oa \mx} \\
(A \oa B) \oa C \ar[r]_{a_\oa} & A \oa (B \oa C)
} ~~~~~~~~~~~~~
\xymatrix{
(A \ox B) \ox C \ar[r]^{\mx} \ar[d]_{a_\ox} \ar@{}[dr]|{(b)} & A \oa ( B \ox C ) \ar@{<-}[d]^{\delta^L} \\
A \ox (B \ox C) \ar[r]_{1 \ox \mx} & A \ox (B \oa C)
}
\end{equation*}
\[  ~~~~~~~~~~~~ \xymatrix{
C \ox (A \oa B) \ar[r]^{\delta^L} \ar[d]_{\mx} \ar@{}[dr]|{(c)} &  (C \ox A) \oa B \ar[d]^{ \mx \oa 1} \\
C \oa (A \oa B) \ar[r]_{a_\oa^{-1}} & (C \oa A) \oa B
} ~~~~~~~~~~~~~
\xymatrix{
A \ox (B \ox C) \ar[r]^{\mx} \ar[d]_{a_\ox^{-1}} \ar@{}[dr]|{(d)} & A \oa ( B \ox C ) \ar@{<-}[d]^{\delta^R} \\
(A \ox B) \ox C \ar[r]_{\mx \ox 1} & (A \oa B) \ox C
} \]

There are many examples of mix categories including coherence spaces \cite{Gir87}, 
and finiteness spaces \cite{Ehr05}.

\begin{definition}
An LDC with a mix map, ${\sf m: \bot \to \top}$ which is an isomorphism is said to be an 
{\bf isomix category}. 
\end{definition}

When ${\sf m}$ is an isomorphism, the coherence requirement for the mixor is automatically 
satisfied (see \cite[Lemma 6.6]{CS97a}). Moreover, the mix map, $\m$, being an isomorphism does
 not imply that the mixor, $\mx$, is an isomorphism. Finiteness spaces \cite{Ehr05} and Chu spaces 
 with the tensor unit as the dualizing object  \cite{Bar06}  provide examples of isomix categories.

 \begin{definition}
 A {\bf compact LDC} is an isomix category in which each mixor, $\mx_{A,B}$ is an isomorphism.  
 \end{definition}

 An important way in which compact LDCs arise is from the ``core" of an isomix category

 \begin{definition} \cite{BCS00} An object $U$ is in the {\bf core} of a mix category if and only if 
	the following natural transformations are isomorphisms: 
    \[ U \ox (\_) \to^{\mx_{U,(\_)}} U \oa (\_) ~~~~\mbox{and}~~~~ 
    (\_) \ox  U \to^{\mx_{(\_),U}} (\_) \oa U \] 
    \end{definition}
    Therefore, the core of a mix category, $\Core(\X) \subseteq \X$, is the full subcategory of $\X$ 
	with the mixor being an isomorphism.  It follows that $\Core(\X)$ is a compact LDC. 

    \begin{proposition} \cite[Proposition 3]{BCS00} 
        If $\X$ is a mix-LDC and $A,B \in \Core(\X)$ then $A \oa B$ and 
        $A \ox B \in \Core(\X)$ (and $A \oa B \simeq A \ox B$).  If $\X$ is an isomix-LDC, 
        then $\top, \bot \in \Core(\X)$.  
        \end{proposition}
        
		A monoidal category is a compact LDC with the mix, and the mixor maps coinciding with the 
		identity. Hence, the tensor and the par coincide in a monoidal category. In fact, any compact LDC is linearly equivalent 
		to a monoidal category. A detailed description of this linear equivalence is given in the section on
		linear functors and transformations. 
		
		The following schematic diagram summarizes different properties of LDCs:
		\begin{center}
		\begin{figure}[h]
			\centering
		\begin{tikzpicture}[scale=1.8]
			\begin{pgfonlayer}{nodelayer}
				\node [style=circle, scale=2, color=black, fill=red] (0) at (-5.75, 2.75) {};
				\node [style=circle, scale=2, color=black, fill=red!70] (1) at (-3.5, 2.75) {};
				\node [style=circle, scale=2, color=black, fill=red!60] (2) at (-1, 2.75) {};
				\node [style=circle, scale=2, color=black, fill=red!40] (3) at (1.75, 2.75) {};
				\node [style=circle, scale=2, color=black, fill=red!20] (5) at (4, 2.75) {};
				\node [style=none] (4) at (-7.75, 2.75) {};
				\node [style=none] (6) at (6, 2.75) {};
				\node [style=none] (7) at (-5.75, 2) {LDC};
				\node [style=none] (8) at (-3.5, 4) {Mix category};
				\node [style=none] (9) at (-3.5, 3.5) {$\m: \bot \to \top$};
				\node [style=none] (10) at (-1, 2) {Isomix category};
				\node [style=none] (11) at (1.75, 4.25) {Compact LDC};
				\node [style=none] (12) at (1.75, 3.65) {$A \ox B \to^{\mx}_{\simeq} A \oa B$};
				\node [style=none] (13) at (4, 2) {Monoidal category};
				\node [style=none] (14) at (-1, 1.4) {$\bot \to^{\m}_{\simeq} \top$};
				\node [style=none] (15) at (4, 1.5) {$\m = 1$, $\mx=1$};
				\node [style=none] (16) at (-5.75, 1.5) {$(\X, \ox, \top)$};
				\node [style=none] (17) at (-5.75, 1) {$(\X, \oa, \bot)$};
			\end{pgfonlayer}
			\begin{pgfonlayer}{edgelayer}
				\draw [dotted] (4.center) to (0);
				\draw (0) to (1);
				\draw (1) to (2);
				\draw (2) to (3);
				\draw (3) to (5);
				\draw [dotted] (5) to (6.center);
			\end{pgfonlayer}
		\end{tikzpicture}
		\caption{Schematic diagram of LDC properties}
		\label{Fig: LDCs}
	\end{figure}
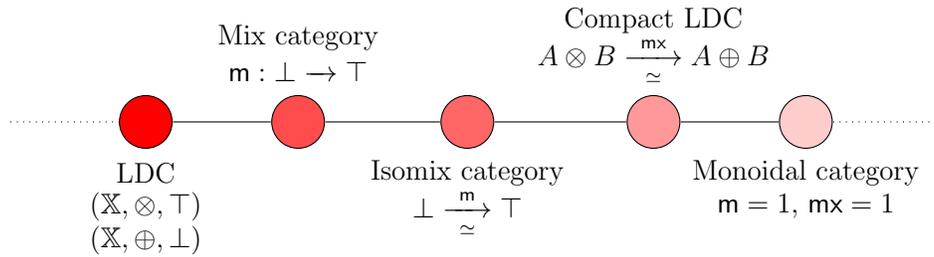
\end{center}

		\subsection{$*$-autonomous categories}
		\label{Sec: *-autonomous}

		A key notion in the theory of LDCs is the notion of a linear adjoint \cite{CKS00}.  
		Here we shall refer to linear adjoints as ``duals'' in order to avoid any confusion 
		with an adjunction of linear functors.   
		
		\begin{definition} 
			\label{defn: duals}
			Suppose $\mathbb{X}$ is a LDC and $A,B \in\X$, then $B$ is {\bf left dual}  
			(or left linear adjoint) to $A$ -- or $A$ is {\bf right dual} (right linear adjoint) 
			to $B$ -- written $(\eta, \epsilon): B \dashv \!\!\!\!\! \dashv  A$, if there 
			exists a unit map, $\eta: \top \rightarrow B \oa A$ and 
			a counit map, $\epsilon: A \ox B \rightarrow \bot$ such that the following diagrams commute:
		\[
		\xymatrix{
		B \ar[r]^{(u_\ox^L)^{-1}} \ar@{=}[d] 
		& \top \ox B \ar[r]^{\eta \ox 1} 
		& (B \oa A) \ox B \ar[d]^{\partial_R} \\
		B 
		& B \oa \bot \ar[l]^{u_\oa^R} 
		& B \oa (A \ox B) \ar[l]^{1 \oa \epsilon}
		}
		~~~~~
		\xymatrix{
		A \ar[r]^{(u_\ox^R)^{-1}} \ar@{=}[d] 
		& A \ox \top  \ar[r]^{1 \ox \eta} 
		& A  \ox  (B \oa A)\ar[d]^{\partial_L} \\
		A
		& \bot \oa A \ar[l]^{u_\oa^L} 
		& (A \ox B) \oa A   \ar[l]^{ \epsilon \oa 1} }
		\]
		\end{definition}
		
	The unit map unit and the counit maps of a dual are drawn in string diagrams as a cap and a cup: 
  \[ \eta := \begin{tikzpicture}
	\begin{pgfonlayer}{nodelayer}
		\node [style=none] (0) at (-1, 1.25) {};
		\node [style=none] (1) at (1, 1.25) {};
		\node [style=none] (2) at (-1, 2.5) {};
		\node [style=none] (3) at (1, 2.5) {};
		\node [style=none] (6) at (-1.25, 1.5) {$A$};
		\node [style=none] (7) at (1.25, 1.5) {$B$};
		\node [style=none] (8) at (0, 3.5) {$\eta$};
	\end{pgfonlayer}
	\begin{pgfonlayer}{edgelayer}
		\draw (1.center) to (3.center);
		\draw (2.center) to (0.center);
		\draw [bend left=90, looseness=1.25] (2.center) to (3.center);
	\end{pgfonlayer}
\end{tikzpicture}  ~~~~~~~\text{ and }~~~~~~~ 
\epsilon := \begin{tikzpicture}
	\begin{pgfonlayer}{nodelayer}
		\node [style=none] (0) at (-1, 3.5) {};
		\node [style=none] (1) at (1, 3.5) {};
		\node [style=none] (2) at (-1, 2.25) {};
		\node [style=none] (3) at (1, 2.25) {};
		\node [style=none] (6) at (-1.25, 3.25) {$B$};
		\node [style=none] (7) at (1.25, 3.25) {$A$};
		\node [style=none, scale=1.5] (8) at (0, 1) {$\epsilon$};
	\end{pgfonlayer}
	\begin{pgfonlayer}{edgelayer}
		\draw (1.center) to (3.center);
		\draw (2.center) to (0.center);
		\draw [bend right=90, looseness=1.50] (2.center) to (3.center);
	\end{pgfonlayer}
\end{tikzpicture} \]
		
	The commuting diagrams are called often referred to as ``snake diagrams'' because of their 
	shape in graphical calculus:
		\[	\begin{tikzpicture}
			\begin{pgfonlayer}{nodelayer}
				\node [style=none] (6) at (1, 0) {};
				\node [style=none] (7) at (1, 1) {};
				\node [style=none] (8) at (2, 1) {};
				\node [style=none] (9) at (2, 0.75) {};
				\node [style=none] (10) at (3, 0.75) {};
				\node [style=none] (11) at (3, 2) {};
				\node [style=none] (12) at (1.5, 1.75) {$\eta$};
				\node [style=none] (13) at (2.5, 0) {$\epsilon$};
				\node [style=none] (14) at (0.75, 0.25) {$A$};
				\node [style=none] (15) at (3.25, 1.75) {$A$};
			\end{pgfonlayer}
			\begin{pgfonlayer}{edgelayer}
				\draw (6.center) to (7.center);
				\draw [bend left=90, looseness=1.50] (7.center) to (8.center);
				\draw (8.center) to (9.center);
				\draw [bend right=90, looseness=1.50] (9.center) to (10.center);
				\draw (10.center) to (11.center);
			\end{pgfonlayer}
		\end{tikzpicture} = \begin{tikzpicture}
		  \draw (0,2.5) -- (0,0);
		\end{tikzpicture} ~~~~~~~~~~
		\begin{tikzpicture}
			\begin{pgfonlayer}{nodelayer}
				\node [style=none] (6) at (3, 0) {};
				\node [style=none] (7) at (3, 1) {};
				\node [style=none] (8) at (2, 1) {};
				\node [style=none] (9) at (2, 0.75) {};
				\node [style=none] (10) at (1, 0.75) {};
				\node [style=none] (11) at (1, 2) {};
				\node [style=none] (12) at (2.5, 1.75) {$\eta$};
				\node [style=none] (13) at (1.5, 0) {$\epsilon$};
				\node [style=none] (14) at (3.25, 0.25) {$B$};
				\node [style=none] (15) at (0.75, 1.75) {$B$};
			\end{pgfonlayer}
			\begin{pgfonlayer}{edgelayer}
				\draw (6.center) to (7.center);
				\draw [bend right=90, looseness=1.50] (7.center) to (8.center);
				\draw (8.center) to (9.center);
				\draw [bend left=90, looseness=1.50] (9.center) to (10.center);
				\draw (10.center) to (11.center);
			\end{pgfonlayer}
		\end{tikzpicture} = \begin{tikzpicture}
		  \draw (0,2.5) -- (0,0);
		\end{tikzpicture} \]
		The linear distributor is hidden in the above circuit diagrams due to the use to circuit expansion, 
		and circuit reduction rules, see Section \ref{Sec: graphical}.
		
		\begin{lemma} \cite{BCS00}
		\begin{enumerate}[(i)]
		\item In an LDC if $(\eta,\epsilon): B \dashvv A$ and $(\eta',\epsilon'): C \dashvv A$, then $B$ and $C$ are isomorphic;
		\item In a symmetric LDC $(\eta, \epsilon): B \dashvv A$ if and only if $(\eta c_\oa, c_\ox \epsilon): A \dashvv B$;
		\item In a mix category if $B \in \Core(\X)$ and $B \dashvv A$, then $A \in \Core(\X)$.
		\end{enumerate}
		\end{lemma}
		
    	In a monoidal category, duals coincide with the usual notion of duals. 
		Next we define the homomorphism of duals:
		
		\begin{definition}
			A homomorphism of duals, $(f,f'): (\eta,  \epsilon) \to (\tau, \gamma)$, is given by a pair of maps
			\[ \xymatrix{ A \ar[d]^f \ar@{-||}[r]^{(\eta,  \epsilon)}  &  B \\ A' \ar@{-||}[r]_{(\tau, \gamma)} &  B' \ar[u]_{f'}} \] 
		such that the following equations hold:  \[ (a) ~~~  \begin{tikzpicture}
			\begin{pgfonlayer}{nodelayer}
				\node [style=none] (0) at (-1, 2) {};
				\node [style=none] (1) at (0.5, 3) {};
				\node [style=none] (2) at (-1, 3) {};
				\node [style=none] (3) at (-0.25, 4) {$\tau$};
				\node [style=none] (4) at (0.5, 2) {};
				\node [style=none] (5) at (-1.25, 2.25) {$A'$};
				\node [style=none] (6) at (1, 3.5) {$B'$};
				\node [style=circle, scale=1.5] (7) at (0.5, 2.75) {};
				\node [style=none] (8) at (0.5, 2.75) {$f'$};
				\node [style=none] (9) at (1, 2.25) {$B$};
			\end{pgfonlayer}
			\begin{pgfonlayer}{edgelayer}
				\draw (4.center) to (7);
				\draw (1.center) to (7);
				\draw (2.center) to (0.center);
				\draw [bend left=90, looseness=1.75] (2.center) to (1.center);
			\end{pgfonlayer}
		\end{tikzpicture}
		= \begin{tikzpicture}
			\begin{pgfonlayer}{nodelayer}
				\node [style=none] (0) at (0.5, 2) {};
				\node [style=none] (1) at (-1, 3) {};
				\node [style=none] (2) at (0.5, 3) {};
				\node [style=none] (3) at (-0.25, 4) {$\eta$};
				\node [style=none] (4) at (-1, 2) {};
				\node [style=none] (5) at (0.75, 2.25) {$B$};
				\node [style=none] (6) at (-1.5, 3.5) {$A$};
				\node [style=circle, scale=1.5] (7) at (-1, 2.75) {};
				\node [style=none] (8) at (-1, 2.75) {$f$};
				\node [style=none] (9) at (-1.5, 2.25) {$A'$};
			\end{pgfonlayer}
			\begin{pgfonlayer}{edgelayer}
				\draw (4.center) to (7);
				\draw (1.center) to (7);
				\draw (2.center) to (0.center);
				\draw [bend right=90, looseness=1.75] (2.center) to (1.center);
			\end{pgfonlayer}
		\end{tikzpicture}
		  ~~~~~~(b)~~~ \begin{tikzpicture}
			\begin{pgfonlayer}{nodelayer}
				\node [style=none] (0) at (0.5, 3.75) {};
				\node [style=none] (1) at (-1, 2.75) {};
				\node [style=none] (2) at (0.5, 2.75) {};
				\node [style=none] (3) at (-0.25, 1.75) {$\gamma$};
				\node [style=none] (4) at (-1, 3.75) {};
				\node [style=none] (5) at (0.75, 3.5) {$B'$};
				\node [style=none] (6) at (-1.5, 2.25) {$A'$};
				\node [style=circle, scale=1.5] (7) at (-1, 3) {};
				\node [style=none] (8) at (-1, 3) {$f$};
				\node [style=none] (9) at (-1.5, 3.5) {$A$};
			\end{pgfonlayer}
			\begin{pgfonlayer}{edgelayer}
				\draw (4.center) to (7);
				\draw (1.center) to (7);
				\draw (2.center) to (0.center);
				\draw [bend right=90, looseness=1.75] (1.center) to (2.center);
			\end{pgfonlayer}
		\end{tikzpicture}
		 =  \begin{tikzpicture}
			\begin{pgfonlayer}{nodelayer}
				\node [style=none] (0) at (-1, 3.75) {};
				\node [style=none] (1) at (0.5, 2.75) {};
				\node [style=none] (2) at (-1, 2.75) {};
				\node [style=none] (3) at (-0.25, 1.75) {$\epsilon$};
				\node [style=none] (4) at (0.5, 3.75) {};
				\node [style=none] (5) at (-1.25, 3.5) {$A$};
				\node [style=none] (6) at (1, 2.25) {$B$};
				\node [style=circle, scale=1.5] (7) at (0.5, 3) {};
				\node [style=none] (8) at (0.5, 3) {$f'$};
				\node [style=none] (9) at (1, 3.5) {$B'$};
			\end{pgfonlayer}
			\begin{pgfonlayer}{edgelayer}
				\draw (4.center) to (7);
				\draw (1.center) to (7);
				\draw (2.center) to (0.center);
				\draw [bend right=90, looseness=1.75] (2.center) to (1.center);
			\end{pgfonlayer}
		\end{tikzpicture} \]
		\end{definition}
		
		Notice that a morphism of duals is determined completely by either of the pair of maps ($f$ 
		completely determines $f'$, and vice versa), and are referred to as Australian mates \cite{CKS00}. :
        \[ f' := \begin{tikzpicture}
			\begin{pgfonlayer}{nodelayer}
				\node [style=circle, scale=2] (0) at (0, 2) {};
				\node [style=none] (1) at (0, 2) {$f$};
				\node [style=none] (2) at (0, 2.5) {};
				\node [style=none] (3) at (1, 2.5) {};
				\node [style=none] (4) at (1, 0.75) {};
				\node [style=none] (5) at (0, 1.5) {};
				\node [style=none] (6) at (-1, 1.5) {};
				\node [style=none] (7) at (-1, 3.5) {};
				\node [style=none] (8) at (-1.25, 3) {$B'$};
				\node [style=none] (9) at (1.25, 1.25) {$B$};
			\end{pgfonlayer}
			\begin{pgfonlayer}{edgelayer}
				\draw (2.center) to (0);
				\draw (0) to (5.center);
				\draw [bend left=90, looseness=1.75] (5.center) to (6.center);
				\draw (6.center) to (7.center);
				\draw [bend left=90, looseness=2.00] (2.center) to (3.center);
				\draw (3.center) to (4.center);
			\end{pgfonlayer}
		\end{tikzpicture}
		\]
		
		The map $f$ is an isomorphism if and only if $f'$ is an isomorphism. Also, notice that 
		if $f$ or $f'$ is an isomorphism, equations $(a)$ and $(b)$ imply one another in the 
		definition of homomorphism of duals. 
		
		A dual, $(\eta, \epsilon): A \dashvv B$ such that $A$ is both left and right dual of $B$, 
		is called a {\bf cyclic dual}. 
		In a symmetric LDC, every dual $(\eta, \epsilon): A \dashvv B$ gives another 
		dual $(\eta c_\oa, c_\ox \epsilon): B \dashvv A$, which is obtained by twisting 
		the wires using the symmetry map. Thus, in a symmetric LDC, every dual is a cyclic dual. 

		\begin{definition}
		An LDC in which every object has a chosen left and right dual, 
		respectively $(\eta{*},\epsilon{*}): A^{*}  \dashvv A$ and 
		$({*}\eta,{*}\epsilon): A  \dashvv \!~^{*}A$, is a {\bf $*$-autonomous category}.   		
	    \end{definition}

		In the symmetric case a left dual gives a right dual using the symmetry: 
		thus, it is standard to assume the existence of just the left dual with the 
		right being the same object with the unit and counit given by symmetry 
		(as above). 

		In a $*$-autonomous category, taking the left (or right) linear dual of an object 
		extends to a linear functor, see Section \ref{Sec: dduals}. 

		Just as compact LDCs are linearly equivalent to monoidal categories so 
		compact $*$-autonomous categories are linearly equivalent to a compact 
		closed categories. The equivalence is given by ${\sf Mx}_\uparrow$ which 
		spreads the par onto  two tensor structures (or, indeed, by 
		${\sf Mx}_{\downarrow}$ which shows how to spread out a compact closed 
		structure on the tensor), see Section \ref{Sec: mix-functor}.
		
		
		In a symmetric $*$-autonomous category the left dual of an object is always 
		canonically isomorphic to the right dual.  Moreover, even in non-symmetric 
		$*$-autonomous categories, it is often the case that the two duals are 
		coherently isomorphic:
		
		\begin{definition}\cite{EggMcCurd12}
		A {\bf cyclor} in a $*$-autonomous category $(\X, \ox, \top, \oa, \bot, ~^{*}(\_), (\_)^*)$ is a natural isomorphism $A^* \to^{\psi} \!~^{*}A$ satisfying the following coherence conditions:
		\[
		\xymatrix{
		\bot^* \ar[rr]^{\psi} \ar[dr] & \ar@{}[d]|{\mbox{\tiny \bf [C.1]}} & ~^{*}\bot \ar[ld]^{} \\
		& \top &
		} ~~~~~~~~ \xymatrix{
		 & A \ar[ld] \ar[dr] \ar@{}[d]|{\mbox{\tiny \bf [C.2]}} & \\
		( ~^{*}A)^* \ar[r]_{\psi^*} & (A^*)^* \ar[r]_{\psi_{A^*}} & ~^{*}(A^*)
		} ~~~~~~~~ \xymatrix{
		\top^* \ar[rr]^{\psi} \ar[dr] & \ar@{}[d]|{\mbox{\tiny \bf [C.3]}} & ~^{*}\top \ar[ld]^{} \\
		& \bot &
		}
		\]
		
		\[
		\xymatrixcolsep{4pc}
		\xymatrix{
		(A \ox B)^* \ar[r]^{\psi} \ar[d]_{} \ar@{}[dr]|{\mbox{\tiny \bf [C.4]}} & ~^{*}(A \ox B) \ar[d]^{} \\
		(B^* \oa A^*) \ar[r]_{\psi \oa \psi} & ~^{*}B \oa ~^{*}A
		} ~~~~~~~~~ \xymatrix{
		(A \oa B)^* \ar[r]^{\psi} \ar[d]_{} \ar@{}[dr]|{\mbox{\tiny \bf [C.5]}}  & ~^{*}(A \oa B) \ar[d]^{} \\
		(B^* \ox A^*) \ar[r]_{\psi \ox \psi} & ~^{*}B \ox ~^{*}A
		}
		\]
		A $*$-autonomous category with a cyclor is said to be {\bf cyclic}.
		\end{definition}
		
		The coherence conditions are not independent of each other: 
		being cyclic is equivalent to any one of the following four pairs of coherences: 
		({\bf [C.1]}, {\bf [C.5]}), ({\bf [C.2]}, {\bf [C.5]}),  ({\bf [C.4]}, {\bf [C.2]}) and 
		({\bf [C.4]}, {\bf [C.3]}) \cite{EggMcCurd12}.
		
		Condition {\bf [C.2]} which is used extensively in Section \ref{daggers-duals-conjugation} is represented graphically as follows:
		\[ \begin{tikzpicture}
			\begin{pgfonlayer}{nodelayer}
				\node [style=circle] (0) at (-2, 1) {$\psi_{A^*}$};
				\node [style=circle] (1) at (-0.75, 1) {$\psi_A$};
				\node [style=none] (2) at (-2, 1.75) {};
				\node [style=none] (3) at (-0.75, 1.75) {};
				\node [style=none] (4) at (-2, -1) {};
				\node [style=none] (5) at (-0.75, -0) {};
				\node [style=none] (6) at (0.5, -0) {};
				\node [style=none] (7) at (0.5, 3) {};
				\node [style=none] (8) at (-1.5, 2.75) {$\eta*$};
				\node [style=none] (9) at (0, -0.75) {$*\epsilon$};
				\node [style=none] (10) at (-2.5, 1.75) {$A^{**}$};
				\node [style=none] (11) at (-2.7, -0.7) {$~^*(A^*)$};
				\node [style=none] (12) at (-0.5, 2) {$A^*$};
				\node [style=none] (13) at (-0.45, 0.25) {$~^*A$};
				\node [style=none] (14) at (0.75, 2.5) {$A$};
			\end{pgfonlayer}
			\begin{pgfonlayer}{edgelayer}
				\draw (0) to (4.center);
				\draw (2.center) to (0);
				\draw [bend left=90, looseness=2.00] (2.center) to (3.center);
				\draw (3.center) to (1);
				\draw (1) to (5.center);
				\draw [bend right=90, looseness=1.50] (5.center) to (6.center);
				\draw (6.center) to (7.center);
			\end{pgfonlayer}
		\end{tikzpicture}  = \begin{tikzpicture}
			\begin{pgfonlayer}{nodelayer}
				\node [style=none] (0) at (-2, 0.25) {};
				\node [style=none] (1) at (-0.75, 0.25) {};
				\node [style=none] (2) at (-2, 3) {};
				\node [style=none] (3) at (-0.75, 2) {};
				\node [style=none] (4) at (0.5, 2) {};
				\node [style=none] (5) at (0.5, -1) {};
				\node [style=none] (6) at (-0.25, 2.75) {$*\eta$};
				\node [style=none] (7) at (-1.5, -0.75) {$\epsilon*$};
				\node [style=none] (8) at (1, -0.7) {$~^*(A^*)$};
				\node [style=none] (9) at (-0.5, 1.15) {$A^*$};
				\node [style=none] (10) at (-2.25, 2.75) {$A$};
			\end{pgfonlayer}
			\begin{pgfonlayer}{edgelayer}
				\draw [bend right=90, looseness=2.00] (0.center) to (1.center);
				\draw [bend left=90, looseness=1.50] (3.center) to (4.center);
				\draw (4.center) to (5.center);
				\draw (3.center) to (1.center);
				\draw (2.center) to (0.center);
			\end{pgfonlayer}
		\end{tikzpicture}
		\]
		
	    The maps in {\bf [C.2]} are invertible:
		\[  \text{\bf [C.2]$^{-1}$} 
		  \begin{tikzpicture}
			\begin{pgfonlayer}{nodelayer}
				\node [style=circle] (0) at (0, 1) {$\psi^{-1}_A$};
				\node [style=circle] (1) at (2, 1) {$\psi^{-1}_{A^*}$};
				\node [style=none] (2) at (0, -0) {};
				\node [style=none] (3) at (2, -0) {};
				\node [style=none] (4) at (2, 3) {};
				\node [style=none] (5) at (0, 2) {};
				\node [style=none] (6) at (-1, 2) {};
				\node [style=none] (7) at (-1, -1.25) {};
				\node [style=none] (8) at (1, -1.1) {$\epsilon*$};
				\node [style=none] (9) at (-0.5, 2.75) {$*\eta$};
				\node [style=none] (10) at (2.5, 0.25) {$A^{**}$};
				\node [style=none] (11) at (-0.5, 0.25) {$A^*$};
				\node [style=none] (12) at (2.5, 2.75) {$~^*(A^*)$};
				\node [style=none] (13) at (-0.5, 1.75) {$~^*A$};
				\node [style=none] (14) at (-1.25, -0.75) {$A$};
			\end{pgfonlayer}
			\begin{pgfonlayer}{edgelayer}
				\draw (0) to (2.center);
				\draw [bend right=90, looseness=1.50] (2.center) to (3.center);
				\draw (3.center) to (1);
				\draw (4.center) to (1);
				\draw (5.center) to (0);
				\draw [bend right=90, looseness=2.00] (5.center) to (6.center);
				\draw (6.center) to (7.center);
			\end{pgfonlayer}
		\end{tikzpicture} =  \left(
		\begin{tikzpicture}
			\begin{pgfonlayer}{nodelayer}
				\node [style=circle] (0) at (-2, 1) {$\psi_{A^*}$};
				\node [style=circle] (1) at (-0.75, 1) {$\psi_A$};
				\node [style=none] (2) at (-2, 1.75) {};
				\node [style=none] (3) at (-0.75, 1.75) {};
				\node [style=none] (4) at (-2, -1) {};
				\node [style=none] (5) at (-0.75, -0) {};
				\node [style=none] (6) at (0.5, -0) {};
				\node [style=none] (7) at (0.5, 3) {};
				\node [style=none] (8) at (-1.5, 2.75) {$\eta*$};
				\node [style=none] (9) at (0, -0.75) {$*\epsilon$};
				\node [style=none] (10) at (-2.5, 1.75) {$A^{**}$};
				\node [style=none] (11) at (-2.7, -0.7) {$~^*(A^*)$};
				\node [style=none] (12) at (-0.5, 2) {$A^*$};
				\node [style=none] (13) at (-0.45, 0.25) {$~^*A$};
				\node [style=none] (14) at (0.75, 2.5) {$A$};
			\end{pgfonlayer}
			\begin{pgfonlayer}{edgelayer}
				\draw (0) to (4.center);
				\draw (2.center) to (0);
				\draw [bend left=90, looseness=2.00] (2.center) to (3.center);
				\draw (3.center) to (1);
				\draw (1) to (5.center);
				\draw [bend right=90, looseness=1.50] (5.center) to (6.center);
				\draw (6.center) to (7.center);
			\end{pgfonlayer}
		\end{tikzpicture} \right)^{-1} = \left( \begin{tikzpicture}
			\begin{pgfonlayer}{nodelayer}
				\node [style=none] (0) at (-2, 0.25) {};
				\node [style=none] (1) at (-0.75, 0.25) {};
				\node [style=none] (2) at (-2, 3) {};
				\node [style=none] (3) at (-0.75, 2) {};
				\node [style=none] (4) at (0.5, 2) {};
				\node [style=none] (5) at (0.5, -1) {};
				\node [style=none] (6) at (-0.25, 2.75) {$*\eta$};
				\node [style=none] (7) at (-1.5, -0.75) {$\epsilon*$};
				\node [style=none] (8) at (1, -0.7) {$~^*(A^*)$};
				\node [style=none] (9) at (-0.5, 1.15) {$A^*$};
				\node [style=none] (10) at (-2.25, 2.75) {$A$};
			\end{pgfonlayer}
			\begin{pgfonlayer}{edgelayer}
				\draw [bend right=90, looseness=2.00] (0.center) to (1.center);
				\draw [bend left=90, looseness=1.50] (3.center) to (4.center);
				\draw (4.center) to (5.center);
				\draw (3.center) to (1.center);
				\draw (2.center) to (0.center);
			\end{pgfonlayer}
		\end{tikzpicture} \right)^{-1}
		= \begin{tikzpicture}
			\begin{pgfonlayer}{nodelayer}
				\node [style=none] (0) at (-2, 0.25) {};
				\node [style=none] (1) at (-0.75, 0.25) {};
				\node [style=none] (2) at (-2, 3) {};
				\node [style=none] (3) at (-0.75, 2) {};
				\node [style=none] (4) at (0.5, 2) {};
				\node [style=none] (5) at (0.5, -1) {};
				\node [style=none] (6) at (-0.25, 2.75) {$\eta*$};
				\node [style=none] (7) at (-1.5, -0.75) {$*\epsilon$};
				\node [style=none] (8) at (-2.75, 2.75) {$~^*(A^*)$};
				\node [style=none] (9) at (-0.5, 1) {$A^*$};
				\node [style=none] (10) at (0.75, -0.75) {$A$};
			\end{pgfonlayer}
			\begin{pgfonlayer}{edgelayer}
				\draw [bend right=90, looseness=2.00] (0.center) to (1.center);
				\draw [bend left=90, looseness=1.50] (3.center) to (4.center);
				\draw (4.center) to (5.center);
				\draw (3.center) to (1.center);
				\draw (2.center) to (0.center);
			\end{pgfonlayer}
		\end{tikzpicture}
		\]
		
		Symmetric $*$-autonomous categories always have a canonical cyclor:
		
		\[
		\begin{tikzpicture}
			\begin{pgfonlayer}{nodelayer}
				\node [style=none] (0) at (-0.25, 2.25) {};
				\node [style=none] (1) at (-0.25, 0.5) {};
				\node [style=none] (12) at (-0.85, 1.7) {$*\eta$};
				\node [style=none] (2) at (-1.5, 0.5) {};
				\node [style=none] (3) at (-1.5, 1) {};
				\node [style=none] (4) at (-0.5, 1) {};
				\node [style=none] (34) at (-0.85, -0.4) {$\epsilon*$};
				\node [style=none] (5) at (-0.5, -1) {};
				\node [style=none] (6) at (0, 2) {$A^{*}$};
				\node [style=none] (7) at (-0.15, -0.75) {${~^*A}$};
			\end{pgfonlayer}
			\begin{pgfonlayer}{edgelayer}
				\draw (0.center) to (1.center);
				\draw [bend left=90, looseness=1.75] (1.center) to (2.center);
				\draw (2.center) to (3.center);
				\draw [bend left=90, looseness=1.75] (3.center) to (4.center);
				\draw (4.center) to (5.center);
			\end{pgfonlayer}
		\end{tikzpicture} 
		\]
		
		We shall use the cyclor in Section \ref{daggers-duals-conjugation} to show how conjugation 
		and dagger are related in the presence of dualization.
		
\section{Linear functors and transformations}
\label{Sec: linear functor}

Functors between LDCs are referred to as linear functors \cite{CS99}. Following the same 
pattern that LDCs generalize monoidal categories, linear functors generalize monoidal functors. 
Moreover, these general functors also account for the structures 
in linear logic such as exponential modalities $({!}, {?})$ and additive connectives $(+, \with)$ 
which are linear functors satisfying additional properties. 

\subsection{Linear functors and transformations}

In order to introduce linear functors, we first recall the definition of monoidal functors.
A functor $F: \X \to \Y$ between monoidal categories is a monoidal functor if its equipped with natural transformations $m_\ox: F(A) \ox F(B) \to F(A \ox B)$ and $m_I: I \to F(I)$ such that the following diagrams commute:

\[ \xymatrixcolsep{15mm}
\xymatrix{
(F(A) \ox F(B)) \ox F(C) \ar[d]_{a_\ox} \ar[r]^{m_\ox \ox 1} & F(A \ox B) \ox F(C) \ar[r]^{m_\ox} & F((A \ox B) \ox C) \ar[d]^{F(a_\ox)} \\
F(A) \ox (F(B) \ox F(C)) \ar[r]_{1 \ox m_\ox} & F(A) \ox F(B \ox C) \ar[r]_{m_\ox} & F(A \ox (B \ox C)) 
} \]

\medskip

\[ \xymatrix{
F(A) \ox I \ar[drr]^{u_\ox^R} \ar[d]_{1 \ox m_I} \\
F(A) \ox F(I) \ar[r]_{m_\ox} & F(A \ox I) \ar[r]_{F(u_\ox^L)} & F(A)}  ~~~~~~~~~~~~~~
\xymatrix{
I \ox F(A) \ar[drr]^{u_\ox^L} \ar[d]_{m_I \ox 1} \\
F(I) \ox F(A) \ar[r]_{m_\ox} & F(I \ox A) \ar[r]_{F(u_\ox^L)} & F(A)
}
\]

The first diagram for the monoidal functor is the associative law, 
and the other two diagrams are the right and the left unit laws respectively. 
Dually, a functor $(F, n_\ox, n_I): \X \to Y$ between monoidal categories is 
comonoidal if  $(F, n_\ox, n_I): \X^{\op} \to \Y^{\op}$ is monoidal. 
Monoidal functors preserve monoids and the comonoidal functors perserve comonoids.

\begin{definition} \cite[Definition 1]{CS99}
Given linearly distributive categories $\X$ and $\Y$, a {\bf linear functor} $F: \X \to \Y$ consists of 
\begin{enumerate}[(i)]
\item a pair of functors $F = (F_\ox, F_\oa)$; $(F_\ox, m_\ox, m_\top): \X \to \Y$ which is 
monoidal with respect to $\ox$ and $(F_\oa, n_\oa, n_\bot): \X \to \Y$ which is comonoidal 
with respect to $\oa$. We refer to $m_\ox$ and $n_\oa$ as {\bf tensor laxors}, 
and $m_\top$ and $n_\bot$ as {\bf unit laxors}. 

\item natural transformations:
\begin{align*}
\nu_\ox^R &: F_\ox(A \oa B) \to F_\oa(A) \oa F_\ox(B) \\
\nu_\ox^L &: F_\ox(A \oa B) \to F_\ox(A) \oa F_\oa(B) \\
\nu_\oa^R &: F_\ox(A) \ox F_\oa(B) \to F_\oa(A \ox B) \\
\nu_\oa^L &: F_\oa(A) \ox F_\ox(B) \to F_\oa( A \ox B)
\end{align*}
\end{enumerate}
such that the following coherence conditions hold:
\begin{enumerate}[{\bf \small [LF.1]}]
\item 
\begin{enumerate}[(a)]
\item $F_\ox(u^L_\oa) = \nu^R_\ox (n_\bot \oa 1) u^L_\oa$
\[ \xymatrix{
F_\ox( \bot \oa A) \ar[r]^{F_\ox(u_\oa^L)} \ar[d]_{\nu_\ox^R} & F_\ox(A) \\
F_\oa(\bot) \oa F_\ox(A) \ar[r]_{n_\bot \oa 1} & \bot \oa F_\ox(A) \ar[u]_{u_\oa^L}
} \]
\item $\nu_\ox^L ( 1 \oa n_\bot ) u_\oa^R = F_\ox( u_\oa^R) $
\item $(u_\ox^L)^{-1} (m_\top \ox 1) \nu_\oa^R = F_\oa((u_\ox^L)^{-1})$
\item $(u_\ox^R)^{-1} (m_\top \ox 1) \nu_\oa^L = F_\oa((u_\ox^R)^{-1}) $
\end{enumerate} 
\item 
\begin{enumerate}[(a)]
\item $F_\ox(a_\oa) \nu^R_\ox (1 \oa \nu^R_\ox) = \nu^R_\ox (n_\oa \oa 1) a_\oa$
\[ \xymatrix{
F_\ox((A \oa B) \oa C) \ar[r]^{F_\ox(a_\oa)} \ar[r]^{F_\ox(a_\oa)} \ar[d]_{\nu_\ox^R} & 
F_\ox(A \oa (B \oa C)) \ar[d]^{\nu_\ox^R} \\
F_\oa(A \oa B) \oa F_\ox(C) \ar[d]_{n_\oa \oa 1} &
F_\oa(A) \oa F_\ox (B \oa C) \ar[d]^{1 \oa \nu_\ox^R} \\
(F_\oa (A) \oa F_\oa(B)) \oa F_\ox(C) \ar[r]_{a_\oa} &
F_\oa(A) \oa (F_\oa(B) \oa F_\oa(C))
} \]
\item $F_\ox(a_\oa) \nu_\ox^L (1 \oa n_\oa ) = \nu^L_\oa (\nu^L \oa 1 ) a_\oa$
\item $(m_\ox \ox 1) \nu_\oa^R F_\oa(a_\ox) = a_\ox (1 \ox \nu_\oa^R) \nu_\oa^R$
\item $(\nu^R_\oa \ox 1) \nu_\oa^L F_\oa(a_\ox) = a_\ox (1 \ox m_\ox) \nu_\oa^L$
\end{enumerate}
\item 
\begin{enumerate}[(a)]
\item $F_\ox(a_\oa)\nu^R_\ox(1 \oa \nu^L_\ox) = \nu_\ox^L (\nu^R_\ox \oa 1) a_\oa$
\[ \xymatrix{
F_\ox((A \oa B) \oa C) \ar[r]^{F_\ox(a_\oa)} \ar[d]_{\nu_\ox^L} & 
F_\ox(A \oa (B \oa C))  \ar[d]^{\nu_\ox^R} \\
F_\ox(A \oa B) \oa F_\oa(C) \ar[d]_{\nu_\ox^R \oa 1} &
F_\oa(A) \oa F_\ox (B \oa C) \ar[d]^{1 \oa \nu_\ox^L} \\
(F_\oa (A) \oa F_\ox(B)) \oa F_\oa(C) \ar[r]_{a_\oa}&
F_\oa(A) \oa (F_\ox(B) \oa F_\oa(C))}
\]
\item $(\nu^R_\oa \ox 1) \nu^L_\oa F_\oa(a_\ox) = a_\ox (1 \ox \nu_\oa^L) \nu_\oa^R$
\end{enumerate}
\item 
\begin{enumerate}[(a)]
\item $(1 \ox \nu^R_\ox) \partial^L (\nu^R_\oa \oa 1) = m_\ox F_\ox(\partial^L) \nu^R_\ox$
\[ \xymatrix{
F_\ox(A) \ox F_\ox(B \oa C) \ar[r]^{1 \ox \nu_\ox^R} \ar[d]_{m_\ox} &
F_\ox(A) \ox (F_\oa(B) \oa F_\ox(C)) \ar[d]^{\partial^L} \\
F_\ox(A \ox (B \oa C)) \ar[d]_{F_\ox(\partial^L)} &
(F_\ox(A) \ox F_\oa(B)) \oa F_\ox(C) \ar[d]^{\nu_\oa^R \oa 1} \\
F_\ox((A \ox B) \oa C) \ar[r]_{\nu_\ox^R} &
F_\oa(A \oa B) \oa F_\ox(C)
} \]
\item $(\nu_\ox^L \ox 1) \partial^R (1 \oa \nu_\oa^L) = m_\ox F_\ox(\partial^R) \nu_\ox^L$
\item $(1 \ox \nu_\ox^L) \partial^L (\nu_\oa^L \oa 1) = \nu_\oa^L F_\oa(\partial^L) n_\oa$
\item $(\nu_\ox^R \ox 1) \partial^R (1 \oa \nu_\oa^R) = \nu_\oa^R F_\oa(\partial^R) n_\oa $
\end{enumerate}
\item 
\begin{enumerate}[(a)]
\item $(1 \ox \nu^L_\ox) \partial^L(m_\ox \oa 1) = m_\ox F_\ox(\partial^L) \nu^L_\ox$ 
\[ \xymatrix{
F_\ox(A) \ox F_\ox(B \oa C) \ar[r]^{1 \ox \nu_\ox^L} \ar[d]_{m_\ox} &
F_\ox(A) \ox (F_\ox(B) \oa F_\oa(C)) \ar[d]^{\partial^L} \\
F_\ox(A \ox (B \oa C)) \ar[d]_{F_\ox(\partial^L)} &
(F_\ox(A) \ox F_\ox(B)) \oa F_\oa(C) \ar[d]^{m_\ox \oa 1} \\
F_\ox((A \ox B) \oa C) \ar[r]_{\nu_\ox^L} &
F_\ox(A \ox B) \oa F_\ox(C)
}\]
\item $(\nu_\ox^R \ox 1) \partial^R (1 \oa m_\ox) = m_\ox F_\ox(\partial^R) \nu_\ox^R$
\item $(1 \ox n_\oa) \partial^L (\nu_\oa^R \oa 1) = \nu_\oa^R F_\oa(\partial^L) n_\oa $
\item $(n_\oa \ox 1) \partial^R (1 \oa \nu_\oa^L) = \nu_\oa^L F_\oa(\partial^R) n_ \oa $
\end{enumerate}
\end{enumerate}
\end{definition}

In the graphical calculus, functors are represented by linear functor boxes \cite{CS99}. 
A linear functor box can either be monoidal or comonoidal. When the functor box is monoidal $(F_\ox)$, 
it has one principal output wire (of $F_\ox$ type) represented by a port where the wire exits the box and 
the other wires (of $F_\oa$ type) are auxiliary. When the box is comonoidal ($F_\oa$), 
it has one principal input wire with a port (of $F_\oa$ type) and the other wires (of $F_\ox$ type) are auxiliary. 
The functor boxes are subject to a very natural ``box eats box'' calculus described in \cite{CS99}.  
A box can eat another box only when a ported wire meets an auxiliary wire. The linear strengths are drawn in the graphical calculus as follows:
\[\nu_\oa^L = \begin{tikzpicture}
	\begin{pgfonlayer}{nodelayer}
		\node [style=none] (0) at (-2, 2) {};
		\node [style=none] (1) at (-1.5, 2) {};
		\node [style=none] (2) at (-2.25, 2) {};
		\node [style=none] (3) at (-2.25, 1) {};
		\node [style=none] (4) at (-0.75, 1) {};
		\node [style=none] (5) at (-0.75, 2) {};
		\node [style=none] (6) at (-2, 2.75) {};
		\node [style=none] (7) at (-1, 2.75) {};
		\node [style=none] (61) at (-2.75, 2.75) {$F_\oa(A)$};
		\node [style=none] (71) at (-0.25, 2.75) {$F_\ox(B)$};
		\node [style=none] (8) at (-1.5, 0.25) {};
		\node [style=none] (81) at (-2.25, 0) {$F_\oa(A \ox B) $};
		\node [style=ox] (9) at (-1.5, 1.5) {};
		\node [style=none] (10) at (-2, 1.25) {$F_\oa$};
	\end{pgfonlayer}
	\begin{pgfonlayer}{edgelayer}
		\draw [in=-90, out=-90, looseness=1.25] (0.center) to (1.center);
		\draw [bend right=15, looseness=1.00] (6.center) to (9);
		\draw [bend right=15, looseness=0.75] (9) to (7.center);
		\draw (2.center) to (3.center);
		\draw (3.center) to (4.center);
		\draw (4.center) to (5.center);
		\draw (5.center) to (2.center);
		\draw (9) to (8.center);
	\end{pgfonlayer}
\end{tikzpicture} ~~~~~~~~ \nu_\oa^R = \begin{tikzpicture}
	\begin{pgfonlayer}{nodelayer}
		\node [style=none] (0) at (-1.5, 2) {};
		\node [style=none] (1) at (-1, 2) {};
		\node [style=none] (2) at (-2.25, 2) {};
		\node [style=none] (3) at (-2.25, 1) {};
		\node [style=none] (4) at (-0.75, 1) {};
		\node [style=none] (5) at (-0.75, 2) {};
		\node [style=none] (6) at (-2, 2.75) {};
		\node [style=none] (7) at (-1, 2.75) {};
		\node [style=none] (61) at (-2.75, 2.75) {$F_\ox(A)$};
		\node [style=none] (71) at (-0.25, 2.75) {$F_\oa(B)$};
		\node [style=none] (8) at (-1.5, 0.25) {};
		\node [style=none] (81) at (-2.25, 0) {$F_\oa(A \ox B) $};
		\node [style=ox] (9) at (-1.5, 1.5) {};
		\node [style=none] (10) at (-2, 1.25) {$F_\oa$};
	\end{pgfonlayer}
	\begin{pgfonlayer}{edgelayer}
		\draw [in=-90, out=-90, looseness=1.25] (0.center) to (1.center);
		\draw [bend right=15, looseness=1.00] (6.center) to (9);
		\draw [bend right=15, looseness=0.75] (9) to (7.center);
		\draw (2.center) to (3.center);
		\draw (3.center) to (4.center);
		\draw (4.center) to (5.center);
		\draw (5.center) to (2.center);
		\draw (9) to (8.center);
	\end{pgfonlayer}
\end{tikzpicture}  ~~~~~~~~ \nu_\ox^L = \begin{tikzpicture}
	\begin{pgfonlayer}{nodelayer}
		\node [style=none] (0) at (-2.25, 1) {};
		\node [style=none] (1) at (-2.25, 2) {};
		\node [style=none] (2) at (-0.75, 2) {};
		\node [style=none] (3) at (-0.75, 1) {};
		\node [style=none] (4) at (-2, 0.25) {};
		\node [style=none] (5) at (-1, 0.25) {};
		\node [style=none] (6) at (-1.5, 2.75) {};
		\node [style=none] (41) at (-2.75, 0.25) {$F_\ox(A)$};
		\node [style=none] (51) at (-0.25, 0.25) {$F_\oa(B)$};
		\node [style=none] (61) at (-2, 3) {$F_\ox(A \oa B)$};
		\node [style=oa] (7) at (-1.5, 1.5) {};
		\node [style=none] (8) at (-2, 1.75) {$F_\ox$};
		\node [style=none] (9) at (-2, 1) {};
		\node [style=none] (10) at (-1.5, 1) {};
	\end{pgfonlayer}
	\begin{pgfonlayer}{edgelayer}
		\draw [bend left=15, looseness=1.00] (4.center) to (7);
		\draw [bend left=15, looseness=0.75] (7) to (5.center);
		\draw (0.center) to (1.center);
		\draw (1.center) to (2.center);
		\draw (2.center) to (3.center);
		\draw (3.center) to (0.center);
		\draw (7) to (6.center);
		\draw [in=90, out=90, looseness=1.25] (9.center) to (10.center);
	\end{pgfonlayer}
\end{tikzpicture} ~~~~~~~~~ \nu_\ox^R = \begin{tikzpicture}
	\begin{pgfonlayer}{nodelayer}
		\node [style=none] (0) at (-2.25, 1) {};
		\node [style=none] (1) at (-2.25, 2) {};
		\node [style=none] (2) at (-0.75, 2) {};
		\node [style=none] (3) at (-0.75, 1) {};
		\node [style=none] (4) at (-2, 0.25) {};
		\node [style=none] (5) at (-1, 0.25) {};
		\node [style=none] (41) at (-2.75, 0.25) {$F_\oa(A)$};
		\node [style=none] (51) at (-0.25, 0.25) {$F_\ox(B)$};
		\node [style=none] (61) at (-2, 3) {$F_\ox(A \oa B)$};
		\node [style=none] (6) at (-1.5, 2.75) {};
		\node [style=oa] (7) at (-1.5, 1.5) {};
		\node [style=none] (8) at (-2, 1.75) {$F_\ox$};
		\node [style=none] (9) at (-1.5, 1) {};
		\node [style=none] (10) at (-1, 1) {};
	\end{pgfonlayer}
	\begin{pgfonlayer}{edgelayer}
		\draw [bend left=15, looseness=1.00] (4.center) to (7);
		\draw [bend left=15, looseness=0.75] (7) to (5.center);
		\draw (0.center) to (1.center);
		\draw (1.center) to (2.center);
		\draw (2.center) to (3.center);
		\draw (3.center) to (0.center);
		\draw (7) to (6.center);
		\draw [in=90, out=90, looseness=1.25] (9.center) to (10.center);
	\end{pgfonlayer}
\end{tikzpicture} \]
\[ m_\ox = \begin{tikzpicture}
	\begin{pgfonlayer}{nodelayer}
		\node [style=none] (0) at (-2.25, 2) {};
		\node [style=none] (1) at (-2.25, 1) {};
		\node [style=none] (2) at (-0.75, 1) {};
		\node [style=none] (3) at (-0.75, 2) {};
		\node [style=none] (4) at (-2, 2.75) {};
		\node [style=none] (5) at (-1, 2.75) {};
		\node [style=none] (41) at (-2.75, 2.75) {$F_\ox(A)$};
		\node [style=none] (51) at (-0.25, 2.75) {$F_\ox(B)$};
		\node [style=none] (6) at (-1.5, 0.25) {};
		\node [style=none] (61) at (-2, 0) {$F_\ox(A \ox B)$};
		\node [style=ox] (7) at (-1.5, 1.5) {};
		\node [style=none] (8) at (-2, 1.25) {$F_\ox$};
		\node [style=none] (9) at (-1.75, 1) {};
		\node [style=none] (10) at (-1.25, 1) {};
	\end{pgfonlayer}
	\begin{pgfonlayer}{edgelayer}
		\draw [bend right=15, looseness=1.00] (4.center) to (7);
		\draw [bend right=15, looseness=0.75] (7) to (5.center);
		\draw (0.center) to (1.center);
		\draw (1.center) to (2.center);
		\draw (2.center) to (3.center);
		\draw (3.center) to (0.center);
		\draw (7) to (6.center);
		\draw [in=90, out=90, looseness=1.25] (9.center) to (10.center);
	\end{pgfonlayer}
\end{tikzpicture} ~~~~~ m_\top =  \begin{tikzpicture}
	\begin{pgfonlayer}{nodelayer}
	      \node [style=none] (8) at (0, 1) {};
		\node [style=circle] (0) at (0, -0) {$\top$};
		\node [style=none] (1) at (0.75, -0.5) {};
		\node [style=none] (2) at (2.75, -0.5) {};
		\node [style=none] (3) at (0.75, -1.5) {};
		\node [style=none] (4) at (2.75, -1.5) {};
		\node [style=none] (5) at (1.75, -2.5) {};
		\node [style=circle] (6) at (1.75, -1) {$\top$};
		\node [style=circle, scale=0.5] (7) at (1.75, -2) {};
	\end{pgfonlayer}
	\begin{pgfonlayer}{edgelayer}
		\draw (1.center) to (2.center);
		\draw (2.center) to (4.center);
		\draw (4.center) to (3.center);
		\draw (3.center) to (1.center);
		\draw (6) to (5.center);
		\draw [dotted, in=-165, out=-90, looseness=1.25] (0) to (7);
		\draw (0) to (8.center);
	\end{pgfonlayer}
\end{tikzpicture} = \begin{tikzpicture}
	\begin{pgfonlayer}{nodelayer}
	      \node [style=none] (8) at (2.75, 1) {};
		\node [style=circle] (0) at (2.75, -0) {$\top$};
		\node [style=none] (1) at (2, -0.5) {};
		\node [style=none] (2) at (0, -0.5) {};
		\node [style=none] (3) at (2, -1.5) {};
		\node [style=none] (4) at (0, -1.5) {};
		\node [style=none] (5) at (1, -2.5) {};
		\node [style=circle] (6) at (1, -1) {$\top$};
		\node [style=circle, scale=0.5] (7) at (1, -2) {};
	\end{pgfonlayer}
	\begin{pgfonlayer}{edgelayer}
		\draw (1.center) to (2.center);
		\draw (2.center) to (4.center);
		\draw (4.center) to (3.center);
		\draw (3.center) to (1.center);
		\draw (6) to (5.center);
		\draw (0) to (8.center);
		\draw [dotted, in=-15, out=-90, looseness=1.25] (0) to (7);
	\end{pgfonlayer}
\end{tikzpicture} ~~~~~ n_\oa =  \begin{tikzpicture}
	\begin{pgfonlayer}{nodelayer}
		\node [style=none] (0) at (-2.25, 1) {};
		\node [style=none] (1) at (-2.25, 2) {};
		\node [style=none] (2) at (-0.75, 2) {};
		\node [style=none] (3) at (-0.75, 1) {};
		\node [style=none] (4) at (-2, 0.25) {};
		\node [style=none] (5) at (-1, 0.25) {};
		\node [style=none] (41) at (-2.75, 0.25) {$F_\oa(A)$};
		\node [style=none] (51) at (-0.25, 0.25) {$F_\oa(B)$};
		\node [style=none] (61) at (-2, 3) {$F_\oa(A \oa B)$};
		\node [style=none] (6) at (-1.5, 2.75) {};
		\node [style=oa] (7) at (-1.5, 1.5) {};
		\node [style=none] (8) at (-2, 1.75) {$F_\oa$};
		\node [style=none] (9) at (-1.25, 2) {};
		\node [style=none] (10) at (-1.75, 2) {};
	\end{pgfonlayer}
	\begin{pgfonlayer}{edgelayer}
		\draw [bend left=15, looseness=1.00] (4.center) to (7);
		\draw [bend left=15, looseness=0.75] (7) to (5.center);
		\draw (0.center) to (1.center);
		\draw (1.center) to (2.center);
		\draw (2.center) to (3.center);
		\draw (3.center) to (0.center);
		\draw (7) to (6.center);
		\draw [in=-90, out=-90, looseness=1.25] (9.center) to (10.center);
	\end{pgfonlayer}
\end{tikzpicture}
~~~~~ n_\bot = 
\begin{tikzpicture}
	\begin{pgfonlayer}{nodelayer}
		\node [style=none] (0) at (-3, 2) {};
		\node [style=none] (1) at (-3, 1) {};
		\node [style=none] (2) at (-1, 2) {};
		\node [style=none] (3) at (-1, 1) {};
		\node [style=circle] (4) at (-2, 1.5) {$\bot$};
		\node [style=circle] (5) at (0, 1) {$\bot$};
		\node [style=none] (6) at (-2, 3) {};
		\node [style=circle, scale=0.5] (7) at (-2, 2.5) {};
		\node [style=none] (8) at (0, 0) {};
	\end{pgfonlayer}
	\begin{pgfonlayer}{edgelayer}
		\draw (0.center) to (1.center);
		\draw (1.center) to (3.center);
		\draw (3.center) to (2.center);
		\draw (2.center) to (0.center);
		\draw (6.center) to (4);
		\draw (8.center) to (5);
		\draw [dotted, bend left=45, looseness=1.25] (7) to (5);
	\end{pgfonlayer}
\end{tikzpicture} = \begin{tikzpicture}
	\begin{pgfonlayer}{nodelayer}
		\node [style=none] (0) at (0, 2) {};
		\node [style=none] (1) at (0, 1) {};
		\node [style=none] (2) at (-2, 2) {};
		\node [style=none] (3) at (-2, 1) {};
		\node [style=circle] (4) at (-1, 1.5) {$\bot$};
		\node [style=circle] (5) at (-3, 1) {$\bot$};
		\node [style=none] (6) at (-1, 3) {};
		\node [style=circle, scale=0.5] (7) at (-1, 2.5) {};
		\node [style=none] (8) at (-3, 0) {};
	\end{pgfonlayer}
	\begin{pgfonlayer}{edgelayer}
		\draw (0.center) to (1.center);
		\draw (1.center) to (3.center);
		\draw (3.center) to (2.center);
		\draw (2.center) to (0.center);
		\draw (6.center) to (4);
		\draw (8.center) to (5);
		\draw [dotted, bend right=45, looseness=1.25] (7) to (5);
	\end{pgfonlayer}
\end{tikzpicture} \]

When working in the categorical doctrine of {\em symmetric} LDCs we will expect the linear functors to preserve the symmetry.   
Thus, a {\bf symmetric linear functor} is a linear functor $F= (F_\ox,F_\oa)$ 
which satisfies in addition:
\[ \xymatrix{F_\ox(A) \ox F_\ox(B) \ar[d]_{c_\ox} \ar[rr]^{m_\ox} & & F_\ox(A \ox B) \ar[d]^{F_\ox(c_\ox)} \\
                   F_\ox(B) \ox F_\ox(A) \ar[rr]_{m_\ox} & & F_\ox(B \ox A) }
    ~~~~~~
   \xymatrix{F_\oa(A \oa B)  \ar[d]_{F_\oa(c_\oa)} \ar[rr]^{n_\ox} & & F_\oa(A) \oa F_\oa(B)\ar[d]^{c_\oa} \\
                   F_\oa(B \ox A) \ar[rr]_{n_\oa} & &  F_\oa(B) \oa F_\oa(A) } \]

Linear functors preserve duals:

\begin{lemma} \cite{CKS00}
\label{Lemma: linear adjoints}
Linear functors preserve duals: when $F: \X \to \Y$ is a linear functor and $(\eta, \epsilon): A \dashvv B \in \X$, then $F_\ox(A) \dashvv F_\oa(B)$ and $F_\oa(A) \dashvv F_\ox(B)$.
\end{lemma}
\begin{proof}
The unit and counit of the adjunction $(\eta', \epsilon'): F_\ox(A) \dashvv F_\oa(B)$ is given as follows:
\[ \eta' := \top \xrightarrow{m_\top} F_\ox(\top) \xrightarrow{F_\ox(\eta)} F_\ox( A \oa B) \xrightarrow{\nu_\ox^L} F_\ox(A) \oa F_\oa(B) =
\begin{tikzpicture}
	\begin{pgfonlayer}{nodelayer}
		\node [style=none] (0) at (0.5, -0.75) {};
		\node [style=none] (1) at (-0.5, -2) {};
		\node [style=none] (2) at (0.5, -2) {};
		\node [style=none] (7) at (-1, 0) {};
		\node [style=none] (8) at (1, 0) {};
		\node [style=none] (9) at (1, -1.25) {};
		\node [style=none] (10) at (-1, -1.25) {};
		\node [style=none] (11) at (-0.75, -1.25) {};
		\node [style=none] (12) at (-0.25, -1.25) {};
		\node [style=none] (13) at (0.75, -0.25) {$F$};
		\node [style=none] (14) at (-0.5, -0.75) {};
		\node [style=none] (15) at (0, -0.25) {$\eta$};
	\end{pgfonlayer}
	\begin{pgfonlayer}{edgelayer}
		\draw (0.center) to (2.center);
		\draw (7.center) to (8.center);
		\draw (8.center) to (9.center);
		\draw (9.center) to (10.center);
		\draw (10.center) to (7.center);
		\draw [bend left=90, looseness=1.25] (11.center) to (12.center);
		\draw [bend left=90, looseness=1.25] (14.center) to (0.center);
		\draw [style=filled] (14.center) to (1.center);
	\end{pgfonlayer}
\end{tikzpicture} \]
\[ \epsilon' := F_\oa(B) \ox F_\ox(A) \xrightarrow{\nu_\oa^L} F_\oa(B \ox A) \xrightarrow{F_\oa(\epsilon)} F_\oa(\bot) \xrightarrow{n_\bot} \bot =
\begin{tikzpicture}
	\begin{pgfonlayer}{nodelayer}
		\node [style=none] (0) at (0.5, -1) {};
		\node [style=none] (1) at (-0.5, 0.25) {};
		\node [style=none] (2) at (0.5, 0.25) {};
		\node [style=none] (7) at (-1, -1.75) {};
		\node [style=none] (8) at (1, -1.75) {};
		\node [style=none] (9) at (1, -0.5) {};
		\node [style=none] (10) at (-1, -0.5) {};
		\node [style=none] (11) at (-0.75, -0.5) {};
		\node [style=none] (12) at (-0.25, -0.5) {};
		\node [style=none] (13) at (0.75, -1.5) {$F$};
		\node [style=none] (14) at (-0.5, -1) {};
		\node [style=none] (15) at (0, -1.5) {$\epsilon$};
	\end{pgfonlayer}
	\begin{pgfonlayer}{edgelayer}
		\draw (0.center) to (2.center);
		\draw (7.center) to (8.center);
		\draw (8.center) to (9.center);
		\draw (9.center) to (10.center);
		\draw (10.center) to (7.center);
		\draw [bend right=90, looseness=1.25] (11.center) to (12.center);
		\draw [bend right=90, looseness=1.25] (14.center) to (0.center);
		\draw (14.center) to (1.center);
	\end{pgfonlayer}
\end{tikzpicture}	 \]

The unit and counit of the other adjunction is given similarly.
\end{proof}			

Natural transformations between linear functors also break into two components linking 
respectively the tensor functors by a monoidal transformation and, in the opposite  direction, 
the par functors by a comonoidal transformation.

A monoidal transformation $\alpha: F \Rightarrow G$ between two monoidal functors is a natural 
transformation $\alpha: F \Rightarrow G$ such that the following diagrams commute: 
\[ \xymatrixcolsep{12mm} \xymatrix{ 
	F(A) \ox F(B) \ar[r]^{\alpha_A \ox \alpha_B} \ar[d]_{m_\ox^F} & G(A) \ox G(B) \ar[d]^{m_\ox^G} \\
	F(A \ox B) \ar[r]_{\alpha_{A \ox B}} & G(A \ox B) } ~~~~~~~~
\xymatrix{
	I \ar[d]_{m_I^F} \ar[dr]^{m_I^F} \\
	F(I) \ar[r]_{\alpha_I} & G(I) } \] 
The coherences for a comonoidal transformation are precisely the mirror images of the above coherences.

\begin{definition} \cite[Definition 3]{CS99}
A {\bf linear (natural) transformation}\footnote{In this thesis, a linear transformation 
is a natural transformation between linear functors, and is different 
from the linear transformations of linear algebra. We drop the word ``natural" for brevity.},  $\alpha: F \to G$,  between parallel linear functors $F,G: \X \to \Y$  consists of a pair of natural transformations $\alpha = (\alpha_\ox,\alpha_\oa)$ such that $\alpha_\ox: F_\ox \to G_\ox$ is a monoidal transformation and $\alpha_\oa: G_\oa \to F_\oa$ is a comonoidal transformation satisfying the following coherence conditions:
\begin{enumerate}[{\bf \small [LT.1]}]
\item $a_\ox \nu^R_\ox (a_\oa \oa 1) = \nu^R_\ox (1 \oa a_\ox)$
\[ \xymatrix{
F_\ox (A \oa B) \ar[rr]^{\alpha_\ox} \ar[d]_{\nu_\ox^R} & & G_\ox(A \oa B) \ar[d]^{\nu_\ox^R} \\
F_\oa(A) \oa F_\ox(B) \ar[dr]_{1 \oa \alpha_\ox} & & G_\oa(A) \oa G_\ox(B) \ar[ld]^{\alpha_\oa \oa 1} \\
& F_\oa(A) \oa G_\ox(B)&
} \]
\item $\alpha_\ox \nu_\ox^L (1 \oa \alpha_\oa) = \nu_\ox^L (\alpha_\ox \oa 1)$
\item $(1 \ox \alpha_\ox) \nu_\oa^L (\alpha_\oa) = (\alpha_\oa \ox 1) \nu_\oa^L$
\item $(\alpha_\ox \ox 1) \nu_\oa^R \alpha_\oa = (1 \ox \alpha_\oa) \nu_\oa^R$
\end{enumerate}
\end{definition}

Conditions {\bf [LT.1]} - {\bf [LT.4]} are represented graphically as follows:
\[
\mbox{\small\bf [LT.1]}~ \begin{tikzpicture} 
	\begin{pgfonlayer}{nodelayer}
		\node [style=none] (0) at (-0.5, -0.75) {};
		\node [style=none] (1) at (-0.5, 0.25) {};
		\node [style=none] (2) at (-2.5, 0.25) {};
		\node [style=none] (3) at (-2.5, -0.75) {};
		\node [style=oa] (4) at (-1.5, -0.25) {};
		\node [style=none] (5) at (-1, -2) {};
		\node [style=none] (6) at (-2, -2) {};
		\node [style=none] (7) at (-1.5, 1.25) {};
		\node [style=circle, scale=2] (8) at (-1, -1.25) {};
		\node [style=none] (9) at (-1, -1.25) {$\alpha_\ox$};
		\node [style=none] (10) at (-0.75, -0) {$F$};
		\node [style=none] (11) at (-1.25, -0.75) {};
		\node [style=none] (12) at (-0.75, -0.75) {};
	\end{pgfonlayer}
	\begin{pgfonlayer}{edgelayer}
		\draw (0.center) to (1.center);
		\draw (1.center) to (2.center);
		\draw (2.center) to (3.center);
		\draw (0.center) to (3.center);
		\draw [in=90, out=-150, looseness=1.00] (4) to (6.center);
		\draw (4) to (7.center);
		\draw (5.center) to (8);
		\draw [in=-30, out=90, looseness=1.00] (8) to (4);
		\draw [bend left=90, looseness=1.25] (11.center) to (12.center);
	\end{pgfonlayer}
\end{tikzpicture} 
 = \begin{tikzpicture} 
	\begin{pgfonlayer}{nodelayer}
		\node [style=none] (0) at (-0.5, -0.75) {};
		\node [style=none] (1) at (-0.5, 0.25) {};
		\node [style=none] (2) at (-2.5, 0.25) {};
		\node [style=none] (3) at (-2.5, -0.75) {};
		\node [style=oa] (4) at (-1.5, -0.25) {};
		\node [style=none] (5) at (-1, -2) {};
		\node [style=none] (6) at (-2, -2) {};
		\node [style=none] (7) at (-2, -1.25) {$\alpha_\oa$};
		\node [style=circle, scale=2] (8) at (-1.5, 1) {};
		\node [style=none] (9) at (-1.5, 2) {};
		\node [style=none] (10) at (-1.5, 1) {$\alpha_\ox$};
		\node [style=none] (11) at (-0.75, -0) {$G$};
		\node [style=circle, scale=2] (12) at (-2, -1.25) {};
		\node [style=none] (13) at (-1.5, -0.75) {};
		\node [style=none] (14) at (-1, -0.75) {};
	\end{pgfonlayer}
	\begin{pgfonlayer}{edgelayer}
		\draw (0.center) to (1.center);
		\draw (1.center) to (2.center);
		\draw (2.center) to (3.center);
		\draw (0.center) to (3.center);
		\draw (9.center) to (8);
		\draw (8) to (4);
		\draw [in=90, out=-44, looseness=0.75] (4) to (5.center);
		\draw [bend left=90, looseness=1.25] (13.center) to (14.center);
		\draw (6.center) to (12);
		\draw [bend left, looseness=1.00] (12) to (4);
	\end{pgfonlayer}
\end{tikzpicture}  ~~~\mbox{\small\bf [LT.2]}~ \begin{tikzpicture}
	\begin{pgfonlayer}{nodelayer}
		\node [style=none] (0) at (-2.5, -0.75) {};
		\node [style=none] (1) at (-2.5, 0.25) {};
		\node [style=none] (2) at (-0.5, 0.25) {};
		\node [style=none] (3) at (-0.5, -0.75) {};
		\node [style=oa] (4) at (-1.5, -0.25) {};
		\node [style=none] (5) at (-2, -2) {};
		\node [style=none] (6) at (-1, -2) {};
		\node [style=none] (7) at (-1.5, 1.25) {};
		\node [style=circle, scale=2] (8) at (-2, -1.25) {};
		\node [style=none] (9) at (-2, -1.25) {$\alpha_\ox$};
		\node [style=none] (10) at (-0.75, 0) {$F$};
		\node [style=none] (11) at (-1.75, -0.75) {};
		\node [style=none] (12) at (-2.25, -0.75) {};
	\end{pgfonlayer}
	\begin{pgfonlayer}{edgelayer}
		\draw (0.center) to (1.center);
		\draw (1.center) to (2.center);
		\draw (2.center) to (3.center);
		\draw (0.center) to (3.center);
		\draw [in=90, out=-15, looseness=1.00] (4) to (6.center);
		\draw (4) to (7.center);
		\draw (5.center) to (8);
		\draw [in=-165, out=90, looseness=1.25] (8) to (4);
		\draw [bend right=90, looseness=1.25] (11.center) to (12.center);
	\end{pgfonlayer}
\end{tikzpicture} = \begin{tikzpicture}
	\begin{pgfonlayer}{nodelayer}
		\node [style=none] (0) at (-2.5, -0.75) {};
		\node [style=none] (1) at (-2.5, 0.25) {};
		\node [style=none] (2) at (-0.5, 0.25) {};
		\node [style=none] (3) at (-0.5, -0.75) {};
		\node [style=oa] (4) at (-1.5, -0.25) {};
		\node [style=none] (5) at (-2, -2) {};
		\node [style=none] (6) at (-1, -2) {};
		\node [style=none] (7) at (-1, -1.25) {$\alpha_\oa$};
		\node [style=circle, scale=2] (8) at (-1.5, 1) {};
		\node [style=none] (9) at (-1.5, 2) {};
		\node [style=none] (10) at (-1.5, 1) {$\alpha_\ox$};
		\node [style=none] (11) at (-0.75, -0) {$G$};
		\node [style=circle, scale=2] (12) at (-1, -1.25) {};
		\node [style=none] (13) at (-1.5, -0.75) {};
		\node [style=none] (14) at (-2, -0.75) {};
	\end{pgfonlayer}
	\begin{pgfonlayer}{edgelayer}
		\draw (0.center) to (1.center);
		\draw (1.center) to (2.center);
		\draw (2.center) to (3.center);
		\draw (0.center) to (3.center);
		\draw (9.center) to (8);
		\draw (8) to (4);
		\draw [in=90, out=-135, looseness=0.75] (4) to (5.center);
		\draw [bend right=90, looseness=1.25] (13.center) to (14.center);
		\draw (6.center) to (12);
		\draw [bend right, looseness=1.00] (12) to (4);
	\end{pgfonlayer}
\end{tikzpicture} ~~~ \mbox{\small\bf [LT.3]}~ \begin{tikzpicture}
	\begin{pgfonlayer}{nodelayer}
		\node [style=none] (0) at (-0.5, 0) {};
		\node [style=none] (1) at (-0.5, -1) {};
		\node [style=none] (2) at (-2.5, -1) {};
		\node [style=none] (3) at (-2.5, -0) {};
		\node [style=ox] (4) at (-1.5, -0.5) {};
		\node [style=none] (5) at (-1, 1.25) {};
		\node [style=none] (6) at (-2, 1.25) {};
		\node [style=none] (7) at (-1.5, -2) {};
		\node [style=circle, scale=2] (8) at (-1, 0.5) {};
		\node [style=none] (9) at (-1, 0.5) {$\alpha_\oa$};
		\node [style=none] (10) at (-0.75, -0.5) {$F$};
		\node [style=none] (11) at (-1.25, 0) {};
		\node [style=none] (12) at (-0.75, 0) {};
	\end{pgfonlayer}
	\begin{pgfonlayer}{edgelayer}
		\draw (0.center) to (1.center);
		\draw (1.center) to (2.center);
		\draw (2.center) to (3.center);
		\draw (0.center) to (3.center);
		\draw [in=-90, out=165, looseness=1.00] (4) to (6.center);
		\draw (4) to (7.center);
		\draw (5.center) to (8);
		\draw [in=15, out=-90, looseness=1.25] (8) to (4);
		\draw [bend right=90, looseness=1.25] (11.center) to (12.center);
	\end{pgfonlayer}
\end{tikzpicture} = \begin{tikzpicture}
	\begin{pgfonlayer}{nodelayer}
		\node [style=none] (0) at (-0.5, 0.75) {};
		\node [style=none] (1) at (-0.5, -0.25) {};
		\node [style=none] (2) at (-2.5, -0.25) {};
		\node [style=none] (3) at (-2.5, 0.75) {};
		\node [style=ox] (4) at (-1.5, 0.25) {};
		\node [style=none] (5) at (-1, 2) {};
		\node [style=none] (6) at (-2, 2) {};
		\node [style=none] (7) at (-2, 1.25) {$\alpha_\ox$};
		\node [style=circle, scale=2] (8) at (-1.5, -1) {};
		\node [style=none] (9) at (-1.5, -2) {};
		\node [style=none] (10) at (-1.5, -1) {$\alpha_\oa$};
		\node [style=none] (11) at (-0.75, 0.5) {$F$};
		\node [style=circle, scale=2] (12) at (-2, 1.25) {};
		\node [style=none] (13) at (-1.5, 0.75) {};
		\node [style=none] (14) at (-1, 0.75) {};
	\end{pgfonlayer}
	\begin{pgfonlayer}{edgelayer}
		\draw (0.center) to (1.center);
		\draw (1.center) to (2.center);
		\draw (2.center) to (3.center);
		\draw (0.center) to (3.center);
		\draw (9.center) to (8);
		\draw (8) to (4);
		\draw [in=-90, out=44, looseness=0.75] (4) to (5.center);
		\draw [bend right=90, looseness=1.25] (13.center) to (14.center);
		\draw (6.center) to (12);
		\draw [bend right, looseness=1.00] (12) to (4);
	\end{pgfonlayer}
\end{tikzpicture}~~~ \mbox{\small\bf [LT.4]}~ \begin{tikzpicture}
	\begin{pgfonlayer}{nodelayer}
		\node [style=none] (0) at (-2.5, 0) {};
		\node [style=none] (1) at (-2.5, -1) {};
		\node [style=none] (2) at (-0.5, -1) {};
		\node [style=none] (3) at (-0.5, 0) {};
		\node [style=ox] (4) at (-1.5, -0.5) {};
		\node [style=none] (5) at (-2, 1.25) {};
		\node [style=none] (6) at (-1, 1.25) {};
		\node [style=none] (7) at (-1.5, -2) {};
		\node [style=circle, scale=2] (8) at (-2, 0.5) {};
		\node [style=none] (9) at (-2, 0.5) {$\alpha_\oa$};
		\node [style=none] (10) at (-0.75, -0.25) {$G$};
		\node [style=none] (11) at (-1.75, 0) {};
		\node [style=none] (12) at (-2.25, 0) {};
	\end{pgfonlayer}
	\begin{pgfonlayer}{edgelayer}
		\draw (0.center) to (1.center);
		\draw (1.center) to (2.center);
		\draw (2.center) to (3.center);
		\draw (0.center) to (3.center);
		\draw [in=-90, out=15, looseness=1.00] (4) to (6.center);
		\draw (4) to (7.center);
		\draw (5.center) to (8);
		\draw [in=165, out=-90, looseness=1.25] (8) to (4);
		\draw [bend left=90, looseness=1.25] (11.center) to (12.center);
	\end{pgfonlayer}
\end{tikzpicture} = \begin{tikzpicture}
	\begin{pgfonlayer}{nodelayer}
		\node [style=none] (0) at (-2.5, 0.75) {};
		\node [style=none] (1) at (-2.5, -0.25) {};
		\node [style=none] (2) at (-0.5, -0.25) {};
		\node [style=none] (3) at (-0.5, 0.75) {};
		\node [style=ox] (4) at (-1.5, 0.25) {};
		\node [style=none] (5) at (-2, 2) {};
		\node [style=none] (6) at (-1, 2) {};
		\node [style=none] (7) at (-1, 1.25) {$\alpha_\ox$};
		\node [style=circle, scale=2] (8) at (-1.5, -1) {};
		\node [style=none] (9) at (-1.5, -2) {};
		\node [style=none] (10) at (-1.5, -1) {$\alpha_\oa$};
		\node [style=none] (11) at (-0.75, 0.5) {$F$};
		\node [style=circle, scale=2] (12) at (-1, 1.25) {};
		\node [style=none] (13) at (-1.5, 0.75) {};
		\node [style=none] (14) at (-2, 0.75) {};
	\end{pgfonlayer}
	\begin{pgfonlayer}{edgelayer}
		\draw (0.center) to (1.center);
		\draw (1.center) to (2.center);
		\draw (2.center) to (3.center);
		\draw (0.center) to (3.center);
		\draw (9.center) to (8);
		\draw (8) to (4);
		\draw [in=-90, out=136, looseness=0.75] (4) to (5.center);
		\draw [bend left=90, looseness=1.00] (13.center) to (14.center);
		\draw (6.center) to (12);
		\draw [bend left, looseness=1.00] (12) to (4);
	\end{pgfonlayer}
\end{tikzpicture} 
\]

An adjunction of linear functors,  $(\eta, \epsilon): F \dashv G$ is an adjunction in the usual sense 
(i.e. satisfying the triangle equalities) in the 2-category of LDCs with linear functors and 
linear natural transformations.   In particular, such an adjunction yields a pair of 
adjunctions: $(\eta_\ox, \epsilon_\ox): F_\ox \dashv G_\ox$ which is a monoidal adjunction, 
and $(\epsilon_\oa, \eta_\oa): G_\oa \dashv F_\oa$ which is a comonoidal adjunction.  
By Kelly's results \cite{Kel97}, a functor with a right adjoint is comonoidal if and only if its 
right adjoint is monoidal. This leads to the observation that:

\begin{lemma}
\label{Lemma: strong adjunction}
If $(\eta, \epsilon): F \dashvv G$ is an adjunction of linear functors, then $F_\ox$ is iso-monoidal 
(or strong) with respect to $\ox$ and $F_\oa$ is iso-comonoidal making the 
linear functor $F$ strong.
\end{lemma}
\begin{proof}
Since $(\eta_\ox, \epsilon_\ox): F_\ox \dashv G_\ox$ is a monoidal adjunction, 
the left adjoint $(F_\ox, m_\ox, m_\top)$ is a strong monoidal functor. Similarly, 
since $(\epsilon_\oa, \eta_\oa): G_\oa \dashv F_\oa$ is a comonoidal adjunction, the 
right adjoint $(F_\oa, n_\oa, n_\bot)$ is a strong comonoidal functor.
\end{proof}

A {\bf linear equivalence} is a linear adjunction in which the unit and counit 
are linear natural isomorphisms.

\subsection{Linear functors for isomix categories}
\label{Sec: mix-functor}

Any isomix category, $(\X,\ox, \oa)$ always has two linear functors ${\sf Mx}_{\downarrow}: 
(\X,\ox, \ox) \to (\X,\ox,\oa)$  and ${\sf Mx}_\uparrow: (\X,\oa, \oa) \to (\X,\ox,\oa)$  
given by the  identity functor, that is $({\sf Mx}_\uparrow)_\ox = ({\sf Mx}_\uparrow)_\oa = 
{\sf Id} = ({\sf Mx}_{\downarrow})_\ox = ({\sf Mx}_{\downarrow})_\oa$.   
The linear strengths and monoidal maps are given by the inverse of the mix map and  the mixor.  
These mix functors take the degenerate linear structure on the tensor (respectively the par) 
and spread it out over both the tensor structures. 

\begin{lemma} \label{mix-functor}
For any isomix category $\X$ the functors ${\sf Mx}_\downarrow: (\X,\ox,\ox) \to (\X,\ox, \oa)$ 
and ${\sf Mx}_\uparrow: (\X,\oa,\oa) \to (\X,\ox, \oa)$ are linear functors.
\end{lemma}

\begin{proof}
We show that ${\sf Mx}_{{\downarrow}}: (\X,\ox, \ox) \to (\X,\ox,\oa)$ is a linear functor: 
the monoidal and comonoidal components of the functor are given by  
$(1, 1, 1)$ and $( 1, \mx, \m^{-1})$ respectively. The linear strengths are $\nu_\ox^L = \nu_\ox^R: A \ox B \to 
A \oa B := \mx$ and $\nu_\oa^L = \nu_\oa^R: A \oa B \to A \oa B := 1$.

First we show $( 1, \mx, \m^{-1}): (\X,\ox,\ox) \to (\X,\ox, \oa)$ is a monoidal functor:

\begin{itemize}
\item The associative law for monoidal functors, $(\mx \ox 1)~\mx~a_\oa = a_\ox~(1 \ox \mx)~\mx$, 
is satisfied:
\[
\begin{tikzpicture} 
\begin{pgfonlayer}{nodelayer}
\node [style=circle, scale=0.5] (0) at (-0.5, 2) {};
\node [style=circle, scale=0.5] (1) at (0.5, 1.25) {};
\node [style=map] (2) at (0, 1.75) {};
\node [style=ox] (3) at (0, 2.5) {};
\node [style=oa] (4) at (0, 0.75) {};
\node [style=ox] (5) at (1, 3) {};
\node [style=oa] (6) at (0, -1) {};
\node [style=oa] (7) at (1, -1.75) {};
\node [style=none] (8) at (1, -3) {};
\node [style=none] (9) at (1, 3.75) {};
\node [style=none] (10) at (-0.75, -3) {};
\node [style=circle, scale=0.5] (11) at (1.5, -1.25) {};
\node [style=map] (12) at (0.75, -0.5) {};
\node [style=circle, scale=0.5] (13) at (0, -0) {};
\end{pgfonlayer}
\begin{pgfonlayer}{edgelayer}
\draw [dotted, in=90, out=-15, looseness=1.25] (0) to (2);
\draw [dotted, in=165, out=-90, looseness=1.25] (2) to (1);
\draw (9.center) to (5);
\draw [bend left=45, looseness=1.00] (5) to (7);
\draw [bend right=15, looseness=1.00] (5) to (3);
\draw [bend right=15, looseness=1.00] (6) to (7);
\draw (7) to (8.center);
\draw [in=90, out=-126, looseness=1.00] (6) to (10.center);
\draw [bend right=60, looseness=1.25] (3) to (4);
\draw (4) to (6);
\draw [bend left=60, looseness=1.25] (3) to (4);
\draw [dotted, in=90, out=0, looseness=1.50] (13) to (12);
\draw [dotted, in=165, out=-90, looseness=1.25] (12) to (11);
\end{pgfonlayer}
\end{tikzpicture} = \begin{tikzpicture} 
\begin{pgfonlayer}{nodelayer}
\node [style=circle, scale=0.5] (0) at (0.25, 1) {};
\node [style=circle, scale=0.5] (1) at (1.75, -0.5) {};
\node [style=map] (2) at (1, 0.5) {};
\node [style=ox] (3) at (-0.25, 1.75) {};
\node [style=ox] (4) at (1, 3) {};
\node [style=oa] (5) at (1, -1.75) {};
\node [style=none] (6) at (1, -3) {};
\node [style=none] (7) at (1, 3.75) {};
\node [style=none] (8) at (-1, -3) {};
\node [style=circle, scale=0.5] (9) at (0.5, -0.75) {};
\node [style=map] (10) at (-0.25, -0) {};
\node [style=circle, scale=0.5] (11) at (-0.75, 1) {};
\end{pgfonlayer}
\begin{pgfonlayer}{edgelayer}
\draw [dotted, in=90, out=-15, looseness=1.25] (0) to (2);
\draw [dotted, in=165, out=-90, looseness=1.25] (2) to (1);
\draw (7.center) to (4);
\draw [bend left=45, looseness=1.00] (4) to (5);
\draw [bend right=15, looseness=1.00] (4) to (3);
\draw (5) to (6.center);
\draw [dotted, in=90, out=0, looseness=1.50] (11) to (10);
\draw [dotted, in=165, out=-90, looseness=1.25] (10) to (9);
\draw [in=90, out=-165, looseness=0.50] (3) to (8.center);
\draw [in=135, out=-45, looseness=1.00] (3) to (5);
\end{pgfonlayer}
\end{tikzpicture} 
  = \begin{tikzpicture} 
\begin{pgfonlayer}{nodelayer}
\node [style=circle, scale=0.5] (0) at (0.25, 2.5) {};
\node [style=circle, scale=0.5] (1) at (1.5, -0.5) {};
\node [style=map] (2) at (1, 0.5) {};
\node [style=ox] (3) at (-0.25, 1.75) {};
\node [style=ox] (4) at (1, 3) {};
\node [style=oa] (5) at (1, -1.25) {};
\node [style=none] (6) at (1, -3) {};
\node [style=none] (7) at (1, 3.75) {};
\node [style=none] (8) at (-1, -3) {};
\node [style=circle, scale=0.5] (9) at (1, -2.25) {};
\node [style=map] (10) at (-0.25, -0) {};
\node [style=circle, scale=0.5] (11) at (-0.75, 1) {};
\end{pgfonlayer}
\begin{pgfonlayer}{edgelayer}
\draw [dotted, in=90, out=-15, looseness=1.25] (0) to (2);
\draw [dotted, in=165, out=-90, looseness=1.25] (2) to (1);
\draw (7.center) to (4);
\draw [bend left=45, looseness=1.00] (4) to (5);
\draw [bend right=15, looseness=1.00] (4) to (3);
\draw (5) to (6.center);
\draw [dotted, in=90, out=0, looseness=1.50] (11) to (10);
\draw [dotted, in=165, out=-90, looseness=1.25] (10) to (9);
\draw [in=90, out=-165, looseness=0.50] (3) to (8.center);
\draw [in=135, out=-45, looseness=1.00] (3) to (5);
\end{pgfonlayer}
\end{tikzpicture}
\]

\item The unit laws for monoidal functors hold.   Here is the pictorial proof of $(1 \ox \m^{-1}) \mx = u_\ox^L(u_\oa^L)^{-1}$, where
the filled rectangles represent $\m^{-1}$:
\[
\begin{tikzpicture}
\begin{pgfonlayer}{nodelayer}
\node [style=map, fill=black] (0) at (0, 2) {};
\node [style=circle, scale=0.5] (1) at (0, 1.5) {};
\node [style=circle, scale=0.5] (2) at (-1.75, -0) {};
\node [style=map] (3) at (-1, 0.75) {};
\node [style=none] (4) at (0, 3.5) {};
\node [style=none] (5) at (0, -1) {};
\node [style=none] (6) at (-1.75, -1) {};
\node [style=none] (7) at (-1.75, 3.5) {};
\end{pgfonlayer}
\begin{pgfonlayer}{edgelayer}
\draw [dotted, in=90, out=-165, looseness=1.00] (1) to (3);
\draw [dotted, in=30, out=-90, looseness=1.25] (3) to (2);
\draw (4.center) to (0);
\draw (7.center) to (6.center);
\draw (0) to (5.center);
\end{pgfonlayer}
\end{tikzpicture} = \begin{tikzpicture}
\begin{pgfonlayer}{nodelayer}
\node [style=map, fill=black] (0) at (0, 2.25) {};
\node [style=circle, scale=0.5] (1) at (0, 1.5) {};
\node [style=circle, scale=0.5] (2) at (-1.75, -0.5) {};
\node [style=map] (3) at (-1, 0.25) {};
\node [style=none] (4) at (0, 3.75) {};
\node [style=none] (5) at (0.5, -1) {};
\node [style=none] (6) at (-1.75, -1) {};
\node [style=none] (7) at (-1.75, 3.75) {};
\node [style=circle] (8) at (0, 3) {$\top$};
\node [style=circle] (9) at (0, 0.25) {$\bot$};
\node [style=circle] (10) at (0.5, -0.5) {$\bot$};
\node [style=circle, scale=0.5] (11) at (0, 1) {};
\end{pgfonlayer}
\begin{pgfonlayer}{edgelayer}
\draw [dotted, dotted, in=90, out=-165, looseness=1.00] (1) to (3);
\draw [dotted, dotted, in=30, out=-90, looseness=1.25] (3) to (2);
\draw (7.center) to (6.center);
\draw (8) to (0);
\draw (0) to (1);
\draw (1) to (9);
\draw (10) to (5.center);
\draw [bend left, looseness=1.00, dotted] (11) to (10);
\draw (4.center) to (8);
\end{pgfonlayer}
\end{tikzpicture} 
  = \begin{tikzpicture}
\begin{pgfonlayer}{nodelayer}
\node [style=map, fill=black] (0) at (0, 2.25) {};
\node [style=circle, scale=0.5] (1) at (0, 1.5) {};
\node [style=circle, scale=0.5] (2) at (-1.75, 0.5) {};
\node [style=map] (3) at (-1, 1) {};
\node [style=none] (4) at (0, 3.75) {};
\node [style=none] (5) at (0, -1) {};
\node [style=none] (6) at (-1.75, -1) {};
\node [style=none] (7) at (-1.75, 3.75) {};
\node [style=circle] (8) at (0, 3) {$\top$};
\node [style=circle] (9) at (0, 0.75) {$\bot$};
\node [style=circle] (10) at (0, -0.25) {$\bot$};
\node [style=circle, scale=0.5] (11) at (-1.75, -0) {};
\end{pgfonlayer}
\begin{pgfonlayer}{edgelayer}
\draw [dotted, dotted, in=90, out=-165, looseness=1.00] (1) to (3);
\draw [dotted, dotted, in=30, out=-90, looseness=1.25] (3) to (2);
\draw (7.center) to (6.center);
\draw (8) to (0);
\draw (0) to (1);
\draw (1) to (9);
\draw (10) to (5.center);
\draw [bend left, looseness=1.00, dotted] (11) to (10);
\draw (4.center) to (8);
\end{pgfonlayer}
\end{tikzpicture} = \begin{tikzpicture}
\begin{pgfonlayer}{nodelayer}
\node [style=circle, scale=0.5] (0) at (-1.75, 1.75) {};
\node [style=none] (1) at (0, 3.75) {};
\node [style=none] (2) at (0, -1) {};
\node [style=none] (3) at (-1.75, -1) {};
\node [style=none] (4) at (-1.75, 3.75) {};
\node [style=circle] (5) at (0, 3) {$\top$};
\node [style=circle] (6) at (0, -0.25) {$\bot$};
\node [style=circle, scale=0.5] (7) at (-1.75, 0.75) {};
\end{pgfonlayer}
\begin{pgfonlayer}{edgelayer}
\draw (4.center) to (3.center);
\draw (6) to (2.center);
\draw [dotted, bend left, looseness=1.00] (7) to (6);
\draw (1.center) to (5);
\draw [dotted, bend left=45, looseness=1.00] (5) to (0);
\end{pgfonlayer}
\end{tikzpicture}
\]
The other unit law holds similarly.
\end{itemize}

${\sf Mx}_\downarrow: (\X,\ox, \ox) \to (\X,\ox,\oa)$ satisfies all the coherence requirements of a linear functor:
{\bf [LF.1]}, {\bf [LF.2]}, and {\bf [LF.3]} hold because $({\sf Mx}_{\downarrow})_\ox$ and $({\sf Mx}_\downarrow)_\oa$ are monoidal and comonoidal respectively, {\bf [LF.4]}(a) becomes 
$\mx a_\oa^{-1} = \partial^L (\mx \oa 1)$ and holds because:
\[
\begin{tikzpicture} 
\begin{pgfonlayer}{nodelayer}
\node [style=oa] (0) at (0, 2.25) {};
\node [style=oa] (1) at (0, 1) {};
\node [style=oa] (2) at (0, -1) {};
\node [style=oa] (3) at (-1, -2) {};
\node [style=none] (4) at (-1, -3) {};
\node [style=none] (5) at (0.5, -3) {};
\node [style=none] (6) at (-2, 3) {};
\node [style=none] (7) at (0, 3) {};
\node [style=map] (8) at (-1.25, 0.5) {};
\node [style=circle, scale=0.5] (9) at (0, -0.25) {};
\node [style=circle, scale=0.5] (10) at (-2, 1.5) {};
\end{pgfonlayer}
\begin{pgfonlayer}{edgelayer}
\draw (6.center) to (10);
\draw (7.center) to (0);
\draw [bend left=60, looseness=1.25] (0) to (1);
\draw [bend right=60, looseness=1.25] (0) to (1);
\draw (1) to (2);
\draw [in=90, out=-45, looseness=1.00] (2) to (5.center);
\draw [in=30, out=-150, looseness=1.50] (2) to (3);
\draw [in=-89, out=135, looseness=1.00] (3) to (10);
\draw (3) to (4.center);
\draw [in=90, out=-15, looseness=1.25,dotted] (10) to (8);
\draw [in=180, out=-90, looseness=1.25, dotted] (8) to (9);
\end{pgfonlayer}
\end{tikzpicture} 
  = \begin{tikzpicture}
\begin{pgfonlayer}{nodelayer}
\node [style=oa] (0) at (0, 2.25) {};
\node [style=oa] (1) at (0, -0) {};
\node [style=oa] (2) at (0, -1) {};
\node [style=oa] (3) at (-1, -2) {};
\node [style=none] (4) at (-1, -3) {};
\node [style=none] (5) at (0.5, -3) {};
\node [style=none] (6) at (-2, 3) {};
\node [style=none] (7) at (0, 3) {};
\node [style=map] (8) at (-1.25, 1.5) {};
\node [style=circle,scale=0.5] (9) at (-0.75, 0.75) {};
\node [style=circle, scale=0.5] (10) at (-2, 2.25) {};
\end{pgfonlayer}
\begin{pgfonlayer}{edgelayer}
\draw (6.center) to (10);
\draw (7.center) to (0);
\draw [bend left=60, looseness=1.25] (0) to (1);
\draw [bend right=60, looseness=1.25] (0) to (1);
\draw (1) to (2);
\draw [in=90, out=-45, looseness=1.00] (2) to (5.center);
\draw [in=30, out=-150, looseness=1.50] (2) to (3);
\draw [in=-89, out=135, looseness=1.00] (3) to (10);
\draw (3) to (4.center);
\draw [dotted, in=90, out=-15, looseness=1.25] (10) to (8);
\draw [dotted, in=180, out=-90, looseness=1.25] (8) to (9);
\end{pgfonlayer}
\end{tikzpicture} 
  = \begin{tikzpicture}
\begin{pgfonlayer}{nodelayer}
\node [style=oa] (0) at (0, 2.25) {};
\node [style=oa] (1) at (-1, -2) {};
\node [style=none] (2) at (-1, -3) {};
\node [style=none] (3) at (0.5, -3) {};
\node [style=none] (4) at (-2, 3) {};
\node [style=none] (5) at (0, 3) {};
\node [style=map] (6) at (-1.25, 1) {};
\node [style=circle, scale=0.5] (7) at (-0.75, -0.5) {};
\node [style=circle, scale=0.5] (8) at (-2, 2.25) {};
\end{pgfonlayer}
\begin{pgfonlayer}{edgelayer}
\draw (4.center) to (8);
\draw (5.center) to (0);
\draw [in=-89, out=135, looseness=1.00] (1) to (8);
\draw (1) to (2.center);
\draw [dotted, in=90, out=-15, looseness=1.25] (8) to (6);
\draw [dotted, in=150, out=-90, looseness=1.00] (6) to (7);
\draw [in=90, out=-45, looseness=0.50] (0) to (3.center);
\draw [in=60, out=-150, looseness=0.75] (0) to (1);
\end{pgfonlayer}
\end{tikzpicture}
\]
{\bf [LF.4]} (b) - (d) and {\bf [LF.5]} (a) - (d) are satisfied similarly. 

Thus, ${\sf Mx}_{\downarrow}$ is a linear functor. 

The proof that ${\sf Mx}_\uparrow$ is a linear functor is (linearly) dual
\end{proof}

\begin{corollary} \label{compact-mix-functor}
	When $\X$ is a compact LDC,  the mix functors, ${\sf Mx}_{\downarrow}$ and ${\sf Mx}_\uparrow$, are linear isomorphisms. 
	Consequently, compact LDCs are linearly equivalent to monoidal categories.
\end{corollary}
	
We shall denote the inverse of ${\sf Mx}_{\downarrow}$ by 
${\sf Mx}^{*}_\downarrow: (\X,\ox,\oa) \to (\X,\oa,\oa)$: this is the identity functor as a mere functor, 
strict on the par structure, and on the tensor structure having as the unit laxor ${\sf m}$ and as the tensor laxor ${\sf mx}^{-1}$.   
Similarly, we shall denote the inverse of ${\sf Mx}_\uparrow$ by ${\sf Mx}^{*}_\uparrow$. 

The linear functors $\Mx_{\downarrow}$ and $\Mx_\uparrow$ are examples of {\em  isomix Frobenius functors\/}, 
which we shall introduce formally in the next section.


\subsection{Frobenius functors}
\label{Sec: frobenius functors}

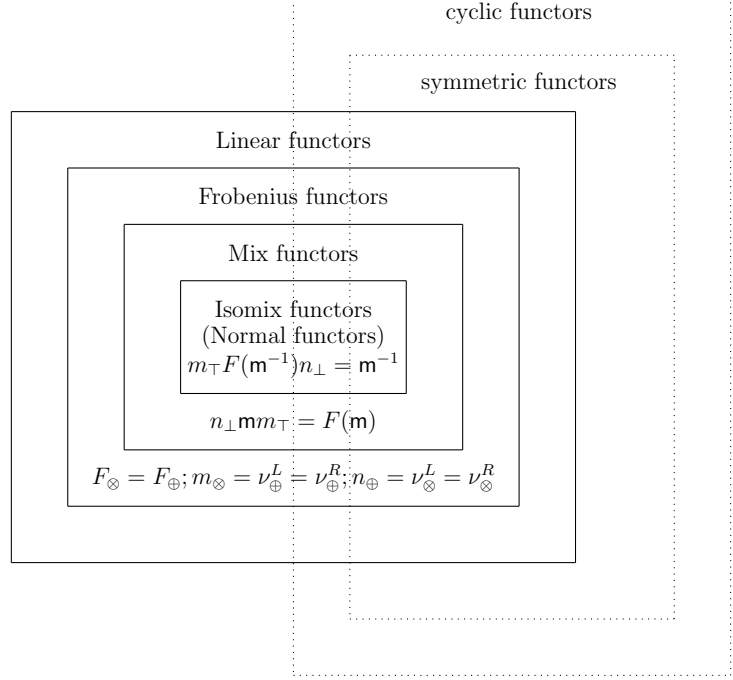
\begin{figure}[ht]
\begin{center}
\begin{tikzpicture} [scale=1.5]
	\begin{pgfonlayer}{nodelayer}
		\node [style=none] (0) at (-3, 2) {};
		\node [style=none] (1) at (-3, -0) {};
		\node [style=none] (2) at (1, -0) {};
		\node [style=none] (3) at (1, 2) {};
		\node [style=none] (4) at (-4, 3) {};
		\node [style=none] (5) at (-4, -1) {};
		\node [style=none] (6) at (2, -1) {};
		\node [style=none] (7) at (2, 3) {};
		\node [style=none] (8) at (-5, -2) {};
		\node [style=none] (9) at (3, -2) {};
		\node [style=none] (10) at (3, 4) {};
		\node [style=none] (11) at (-5, 4) {};
		\node [style=none] (12) at (-6, 5) {};
		\node [style=none] (13) at (4, 5) {};
		\node [style=none] (14) at (4, -3) {};
		\node [style=none] (15) at (-6, -3) {};
		\node [style=none] (16) at (-1, 4.5) {Linear functors};
		\node [style=none] (17) at (-1, 3.5) {Frobenius functors};
		\node [style=none] (18) at (-1, 2.5) {Mix functors};
		\node [style=none] (19) at (-1, 1.5) {Isomix functors};
		\node [style=none] (34) at (-1, 1) {(Normal functors)};
		\node [style=none] (20) at (-1, 0.5) {$m_\top F(\m^{-1})n_\bot = \m^{-1}$};
		\node [style=none] (21) at (-1, -0.5) {$n_\bot\m m_\top = F(\m)$};
		\node [style=none] (22) at (-1, -1.5) {$F_\ox = F_\oa; m_\ox = \nu_\oa^L = \nu_\oa^R; n_\oa = \nu_\ox^L = \nu_\ox^R$};
		\node [style=none] (23) at (-1, -2.75) {};
		\node [style=none] (24) at (0, 6) {};
		\node [style=none] (25) at (0, -4) {};
		\node [style=none] (26) at (5.75, 6) {};
		\node [style=none] (27) at (5.75, -4) {};
		\node [style=none] (28) at (-0.9999999, 7) {};
		\node [style=none] (29) at (-0.9999999, -5) {};
		\node [style=none] (30) at (6.75, -5) {};
		\node [style=none] (31) at (6.75, 7) {};
		\node [style=none] (32) at (3, 6.75) {cyclic functors};
		\node [style=none] (33) at (3, 5.5) {symmetric functors};
	\end{pgfonlayer}
	\begin{pgfonlayer}{edgelayer}
		\draw (0.center) to (1.center);
		\draw (1.center) to (2.center);
		\draw (2.center) to (3.center);
		\draw (3.center) to (0.center);
		\draw (4.center) to (5.center);
		\draw (5.center) to (6.center);
		\draw (6.center) to (7.center);
		\draw (7.center) to (4.center);
		\draw (11.center) to (10.center);
		\draw (10.center) to (9.center);
		\draw (9.center) to (8.center);
		\draw (8.center) to (11.center);
		\draw (12.center) to (13.center);
		\draw (13.center) to (14.center);
		\draw (12.center) to (15.center);
		\draw (15.center) to (14.center);
		\draw[dotted]  (28.center) to (29.center);
		\draw[dotted] (29.center) to (30.center);
		\draw[dotted]  (30.center) to (31.center);
		\draw[dotted]  (31.center) to (28.center);
		\draw[dotted]  (24.center) to (25.center);
		\draw[dotted]  (25.center) to (27.center);
		\draw[dotted]  (27.center) to (26.center);
		\draw[dotted] (26.center) to (24.center);
	\end{pgfonlayer}
\end{tikzpicture}
\end{center}
\caption{Linear functor family}
\label{linear-functor-family}
\end{figure}

In this thesis, we will be interested in linear functors between LDCs called the Frobenius 
functors which come in various flavours, including mix functors and isomix functors, 
as illustrated in Figure \ref{linear-functor-family}.  These functors are directly related 
to the Frobenius monoidal functors of \cite{DP08} and they are referred to as degenerate 
linear functors in \cite{BPS12}.  Furthermore, we have already seen two rather 
basic examples, namely, ${\sf Mx}_\uparrow$ and ${\sf Mx}_{\downarrow}$.

Frobenius functors preserve duals and with an additional coherence 
condition they preserve the mix map.  The coherence requirements for a dagger 
on an LDC are implied by requiring that the dagger functor be a Frobenius 
involutive equivalence.  Once the dagger  is understood we can consider $\dagger$-mix 
categories and their functors which we shall take to be mix Frobenius functors 
with a further requirement concerning the preservation of the dagger.

\begin{definition}
Let $\X$ and $\Y$ to LDCs. A {\bf Frobenius functor} is a linear functor $F: \X \to \Y$ such that:
\begin{enumerate}[{\bf \small [FLF.1]}]
\item $F_\ox = F_\oa $
\item $m_\ox = \nu_\oa^R = \nu_\oa^L $ 
\item $n_\oa = \nu_\ox^L = \nu_\ox^R$
\end{enumerate}
\end{definition}

The left and right linear strengths of $\ox$ and $\oa$ coinciding with the 
$m_\ox$ and $n_\oa$ respectively means that in the diagrammatic 
calculus, ports can be moved around freely:
\[
\begin{tikzpicture}
	\begin{pgfonlayer}{nodelayer}
		\node [style=none] (0) at (-2, 2) {};
		\node [style=none] (1) at (-1.5, 2) {};
		\node [style=none] (2) at (-2.25, 2) {};
		\node [style=none] (3) at (-2.25, 1) {};
		\node [style=none] (4) at (-0.75, 1) {};
		\node [style=none] (5) at (-0.75, 2) {};
		\node [style=none] (6) at (-2, 2.75) {};
		\node [style=none] (7) at (-1, 2.75) {};
		\node [style=none] (61) at (-2.75, 2.75) {$F_\oa(A)$};
		\node [style=none] (71) at (-0.25, 2.75) {$F_\ox(B)$};
		\node [style=none] (8) at (-1.5, 0.25) {};
		\node [style=none] (81) at (-2.25, 0) {$F_\oa(A \ox B) $};
		\node [style=ox] (9) at (-1.5, 1.5) {};
		\node [style=none] (10) at (-2, 1.25) {$F$};
	\end{pgfonlayer}
	\begin{pgfonlayer}{edgelayer}
		\draw [in=-90, out=-90, looseness=1.25] (0.center) to (1.center);
		\draw [bend right=15, looseness=1.00] (6.center) to (9);
		\draw [bend right=15, looseness=0.75] (9) to (7.center);
		\draw (2.center) to (3.center);
		\draw (3.center) to (4.center);
		\draw (4.center) to (5.center);
		\draw (5.center) to (2.center);
		\draw (9) to (8.center);
	\end{pgfonlayer}
\end{tikzpicture} = \begin{tikzpicture}
	\begin{pgfonlayer}{nodelayer}
		\node [style=none] (0) at (-1.5, 2) {};
		\node [style=none] (1) at (-1, 2) {};
		\node [style=none] (2) at (-2.25, 2) {};
		\node [style=none] (3) at (-2.25, 1) {};
		\node [style=none] (4) at (-0.75, 1) {};
		\node [style=none] (5) at (-0.75, 2) {};
		\node [style=none] (6) at (-2, 2.75) {};
		\node [style=none] (7) at (-1, 2.75) {};
		\node [style=none] (61) at (-2.75, 2.75) {$F_\ox(A)$};
		\node [style=none] (71) at (-0.25, 2.75) {$F_\oa(B)$};
		\node [style=none] (8) at (-1.5, 0.25) {};
		\node [style=none] (81) at (-2.25, 0) {$F_\oa(A \ox B) $};
		\node [style=ox] (9) at (-1.5, 1.5) {};
		\node [style=none] (10) at (-2, 1.25) {$F$};
	\end{pgfonlayer}
	\begin{pgfonlayer}{edgelayer}
		\draw [in=-90, out=-90, looseness=1.25] (0.center) to (1.center);
		\draw [bend right=15, looseness=1.00] (6.center) to (9);
		\draw [bend right=15, looseness=0.75] (9) to (7.center);
		\draw (2.center) to (3.center);
		\draw (3.center) to (4.center);
		\draw (4.center) to (5.center);
		\draw (5.center) to (2.center);
		\draw (9) to (8.center);
	\end{pgfonlayer}
\end{tikzpicture} = \begin{tikzpicture}
	\begin{pgfonlayer}{nodelayer}
		\node [style=none] (0) at (-2.25, 2) {};
		\node [style=none] (1) at (-2.25, 1) {};
		\node [style=none] (2) at (-0.75, 1) {};
		\node [style=none] (3) at (-0.75, 2) {};
		\node [style=none] (4) at (-2, 2.75) {};
		\node [style=none] (5) at (-1, 2.75) {};
		\node [style=none] (41) at (-2.75, 2.75) {$F_\ox(A)$};
		\node [style=none] (51) at (-0.25, 2.75) {$F_\ox(B)$};
		\node [style=none] (6) at (-1.5, 0.25) {};
		\node [style=none] (61) at (-2, 0) {$F_\ox(A \ox B)$};
		\node [style=ox] (7) at (-1.5, 1.5) {};
		\node [style=none] (8) at (-2, 1.25) {$F$};
		\node [style=none] (9) at (-1.75, 1) {};
		\node [style=none] (10) at (-1.25, 1) {};
	\end{pgfonlayer}
	\begin{pgfonlayer}{edgelayer}
		\draw [bend right=15, looseness=1.00] (4.center) to (7);
		\draw [bend right=15, looseness=0.75] (7) to (5.center);
		\draw (0.center) to (1.center);
		\draw (1.center) to (2.center);
		\draw (2.center) to (3.center);
		\draw (3.center) to (0.center);
		\draw (7) to (6.center);
		\draw [in=90, out=90, looseness=1.25] (9.center) to (10.center);
	\end{pgfonlayer}
\end{tikzpicture}  ~~~~~~~~ \begin{tikzpicture}
	\begin{pgfonlayer}{nodelayer}
		\node [style=none] (0) at (-2.25, 1) {};
		\node [style=none] (1) at (-2.25, 2) {};
		\node [style=none] (2) at (-0.75, 2) {};
		\node [style=none] (3) at (-0.75, 1) {};
		\node [style=none] (4) at (-2, 0.25) {};
		\node [style=none] (5) at (-1, 0.25) {};
		\node [style=none] (6) at (-1.5, 2.75) {};
		\node [style=none] (41) at (-2.75, 0.25) {$F_\ox(A)$};
		\node [style=none] (51) at (-0.25, 0.25) {$F_\oa(B)$};
		\node [style=none] (61) at (-2, 3) {$F_\ox(A \oa B)$};
		\node [style=oa] (7) at (-1.5, 1.5) {};
		\node [style=none] (8) at (-2, 1.75) {$F$};
		\node [style=none] (9) at (-2, 1) {};
		\node [style=none] (10) at (-1.5, 1) {};
	\end{pgfonlayer}
	\begin{pgfonlayer}{edgelayer}
		\draw [bend left=15, looseness=1.00] (4.center) to (7);
		\draw [bend left=15, looseness=0.75] (7) to (5.center);
		\draw (0.center) to (1.center);
		\draw (1.center) to (2.center);
		\draw (2.center) to (3.center);
		\draw (3.center) to (0.center);
		\draw (7) to (6.center);
		\draw [in=90, out=90, looseness=1.25] (9.center) to (10.center);
	\end{pgfonlayer}
\end{tikzpicture} = \begin{tikzpicture}
	\begin{pgfonlayer}{nodelayer}
		\node [style=none] (0) at (-2.25, 1) {};
		\node [style=none] (1) at (-2.25, 2) {};
		\node [style=none] (2) at (-0.75, 2) {};
		\node [style=none] (3) at (-0.75, 1) {};
		\node [style=none] (4) at (-2, 0.25) {};
		\node [style=none] (5) at (-1, 0.25) {};
		\node [style=none] (41) at (-2.75, 0.25) {$F_\oa(A)$};
		\node [style=none] (51) at (-0.25, 0.25) {$F_\ox(B)$};
		\node [style=none] (61) at (-2, 3) {$F_\ox(A \oa B)$};
		\node [style=none] (6) at (-1.5, 2.75) {};
		\node [style=oa] (7) at (-1.5, 1.5) {};
		\node [style=none] (8) at (-2, 1.75) {$F$};
		\node [style=none] (9) at (-1.5, 1) {};
		\node [style=none] (10) at (-1, 1) {};
	\end{pgfonlayer}
	\begin{pgfonlayer}{edgelayer}
		\draw [bend left=15, looseness=1.00] (4.center) to (7);
		\draw [bend left=15, looseness=0.75] (7) to (5.center);
		\draw (0.center) to (1.center);
		\draw (1.center) to (2.center);
		\draw (2.center) to (3.center);
		\draw (3.center) to (0.center);
		\draw (7) to (6.center);
		\draw [in=90, out=90, looseness=1.25] (9.center) to (10.center);
	\end{pgfonlayer}
\end{tikzpicture} = \begin{tikzpicture}
	\begin{pgfonlayer}{nodelayer}
		\node [style=none] (0) at (-2.25, 1) {};
		\node [style=none] (1) at (-2.25, 2) {};
		\node [style=none] (2) at (-0.75, 2) {};
		\node [style=none] (3) at (-0.75, 1) {};
		\node [style=none] (4) at (-2, 0.25) {};
		\node [style=none] (5) at (-1, 0.25) {};
		\node [style=none] (41) at (-2.75, 0.25) {$F_\oa(A)$};
		\node [style=none] (51) at (-0.25, 0.25) {$F_\oa(B)$};
		\node [style=none] (61) at (-2, 3) {$F_\oa(A \oa B)$};
		\node [style=none] (6) at (-1.5, 2.75) {};
		\node [style=oa] (7) at (-1.5, 1.5) {};
		\node [style=none] (8) at (-2, 1.75) {$F$};
		\node [style=none] (9) at (-1.25, 2) {};
		\node [style=none] (10) at (-1.75, 2) {};
	\end{pgfonlayer}
	\begin{pgfonlayer}{edgelayer}
		\draw [bend left=15, looseness=1.00] (4.center) to (7);
		\draw [bend left=15, looseness=0.75] (7) to (5.center);
		\draw (0.center) to (1.center);
		\draw (1.center) to (2.center);
		\draw (2.center) to (3.center);
		\draw (3.center) to (0.center);
		\draw (7) to (6.center);
		\draw [in=-90, out=-90, looseness=1.25] (9.center) to (10.center);
	\end{pgfonlayer}
\end{tikzpicture}
\]
\[ \nu_\oa^L = \nu_\oa^R = m_\ox  ~~~~~~~~~~~~~~~~~~~~~~~~~~~~~~~~~~~~  \nu_\ox^L = \nu_\ox^R = n_\oa \]
This implies that the ports can be omitted in the circuits.

A Frobenius functor is {\bf symmetric} if as a linear functor it preserves 
the symmetries of the tensor and par.  

\begin{lemma}
\label{Lemma: Frobenius}
Suppose $\X$ and $\Y$ are LDCs. The following are equivalent:
\begin{enumerate}[(a)]
\item $F: \X \to \Y$ is a Frobenius linear functor.
\item $F$ is $\ox$-monoidal and $\oa$-comonoidal such that 
\[ \xymatrixcolsep{2.5pc}
\xymatrix{
F(A) \ox F(B \oa C) \ar[r]^{1 \ox n_\oa} \ar[d]_{m_\ox} \ar@{}[dr]|{\tiny{\bf [F.1]}} & 
F(A) \ox (F(B) \oa F(C)) \ar[d]^{\delta^L} \\
F(A \ox (B \oa C)) \ar[d]_{F(\delta^L)} & (F(A) \ox F(B)) \oa F(C) \ar[d]^{m_\ox \oa 1} \\
F((A \oa B) \oa C) \ar[r]_{n_\oa} & F(A \oa B) \oa F(C)
}\] \[ \xymatrixcolsep{2.5pc} \xymatrix{
F( A \oa B) \ox F(C) \ar[r]^{n_\oa \ox 1}  \ar[d]_{m_\ox} \ar@{}[dr]|{\tiny{\bf [F.2]}} & 
(F(A) \oa F(B)) \ox F(C) \ar[d]^{\delta^R} \\
F( (A \oa B) \ox C) \ar[d]_{F(\delta^R)} & F(A) \oa (F(B) \ox F(C)) \ar[d]^{1 \oa m_\ox} \\
F(A \oa (B \ox C)) \ar[r]_{n_\oa} & F(A) \oa F(B \ox C) 
} \]
\end{enumerate}
\end{lemma}
\begin{proof}
For (a)$\Rightarrow$(b), fix $F := F_\ox = F_\oa$,  then $F$ is $\ox$-monoidal and $\oa$-comonoidal. 
Conditions {\bf \small [F.1]} and {\bf \small [F.2]} are given by {\bf  \small [LF.5]-(a)} and 
{\bf \small [LF.5]-(b)}. For the other direction, define $F_\ox = F_\oa := F$. Then it is 
straightforward to check that all the axioms of Frobenius linear functors are 
satisfied by $(F_\ox, F_\oa)$.
\end{proof}

Conditions {\bf \small  [F.1]} and {\bf \small [F.2]} in Lemma \ref{Lemma: Frobenius} 
are diagrammatically represented as follows:
\[
{\bf [F.1]}~~~~~~~~
\begin{tikzpicture}
	\begin{pgfonlayer}{nodelayer}
		\node [style=oa] (0) at (-3, 2) {};
		\node [style=ox] (1) at (-4, 1) {};
		\node [style=none] (2) at (-4.5, 3) {};
		\node [style=none] (3) at (-3, 3) {};
		\node [style=none] (4) at (-2.5, -0) {};
		\node [style=none] (5) at (-4, -0) {};
		\node [style=none] (6) at (-4.75, 2.5) {};
		\node [style=none] (7) at (-4.75, 0.5) {};
		\node [style=none] (8) at (-2.25, 0.5) {};
		\node [style=none] (9) at (-2.25, 2.5) {};
		\node [style=none] (10) at (-2.5, 2.25) {$F$};
	\end{pgfonlayer}
	\begin{pgfonlayer}{edgelayer}
		\draw [in=15, out=-150, looseness=1.00] (0) to (1);
		\draw [in=-90, out=135, looseness=1.00] (1) to (2.center);
		\draw (0) to (3.center);
		\draw [in=90, out=-45, looseness=1.00] (0) to (4.center);
		\draw (1) to (5.center);
		\draw (6.center) to (9.center);
		\draw (9.center) to (8.center);
		\draw (8.center) to (7.center);
		\draw (7.center) to (6.center);
	\end{pgfonlayer}
\end{tikzpicture} = \begin{tikzpicture}
	\begin{pgfonlayer}{nodelayer}
		\node [style=oa] (0) at (-3, 2.5) {};
		\node [style=ox] (1) at (-4, 1) {};
		\node [style=none] (2) at (-4.5, 3.5) {};
		\node [style=none] (3) at (-3, 3.5) {};
		\node [style=none] (4) at (-2.5, -0) {};
		\node [style=none] (5) at (-4, -0) {};
		\node [style=none] (6) at (-4.75, 1.5) {};
		\node [style=none] (7) at (-4.75, 0.5) {};
		\node [style=none] (8) at (-3.25, 0.5) {};
		\node [style=none] (9) at (-3.25, 1.5) {};
		\node [style=none] (10) at (-3.75, 2) {};
		\node [style=none] (11) at (-2.25, 2) {};
		\node [style=none] (12) at (-2.25, 3) {};
		\node [style=none] (13) at (-3.75, 3) {};
		\node [style=none] (14) at (-2.5, 2.75) {$F$};
		\node [style=none] (15) at (-3.5, 1.25) {$F$};
	\end{pgfonlayer}
	\begin{pgfonlayer}{edgelayer}
		\draw [in=60, out=-150, looseness=1.00] (0) to (1);
		\draw [in=-90, out=135, looseness=1.00] (1) to (2.center);
		\draw (0) to (3.center);
		\draw [in=90, out=-45, looseness=1.00] (0) to (4.center);
		\draw (1) to (5.center);
		\draw (6.center) to (9.center);
		\draw (9.center) to (8.center);
		\draw (8.center) to (7.center);
		\draw (7.center) to (6.center);
		\draw (13.center) to (10.center);
		\draw (10.center) to (11.center);
		\draw (11.center) to (12.center);
		\draw (12.center) to (13.center);
	\end{pgfonlayer}
\end{tikzpicture} 
~~~~~~~~~~~~~~~~~~ {\bf [F.2]}~~~~~~~~ \begin{tikzpicture}
	\begin{pgfonlayer}{nodelayer}
		\node [style=oa] (0) at (-3.25, 1) {};
		\node [style=ox] (1) at (-4, 1.75) {};
		\node [style=none] (2) at (-4, 3) {};
		\node [style=none] (3) at (-2.75, 3) {};
		\node [style=none] (4) at (-3.25, -0) {};
		\node [style=none] (5) at (-4.5, -0) {};
		\node [style=none] (6) at (-4.75, 2.5) {};
		\node [style=none] (7) at (-4.75, 0.5) {};
		\node [style=none] (8) at (-2.25, 0.5) {};
		\node [style=none] (9) at (-2.25, 2.5) {};
		\node [style=none] (10) at (-2.5, 2.25) {$F$};
	\end{pgfonlayer}
	\begin{pgfonlayer}{edgelayer}
		\draw [in=-45, out=165, looseness=1.25] (0) to (1);
		\draw [in=-90, out=90, looseness=1.00] (1) to (2.center);
		\draw [in=-90, out=45, looseness=1.00] (0) to (3.center);
		\draw [in=90, out=-90, looseness=1.00] (0) to (4.center);
		\draw [in=90, out=-150, looseness=0.75] (1) to (5.center);
		\draw (6.center) to (9.center);
		\draw (9.center) to (8.center);
		\draw (8.center) to (7.center);
		\draw (7.center) to (6.center);
	\end{pgfonlayer}
\end{tikzpicture} = \begin{tikzpicture}
	\begin{pgfonlayer}{nodelayer}
		\node [style=oa] (0) at (-3.25, 0.5) {};
		\node [style=ox] (1) at (-4, 1.75) {};
		\node [style=none] (2) at (-4, 3) {};
		\node [style=none] (3) at (-2.75, 3) {};
		\node [style=none] (4) at (-3.25, -0.5) {};
		\node [style=none] (5) at (-4.5, -0.5) {};
		\node [style=none] (6) at (-2.5, 1) {};
		\node [style=none] (7) at (-2.5, -0) {};
		\node [style=none] (8) at (-3.75, -0) {};
		\node [style=none] (9) at (-3.75, 1) {};
		\node [style=none] (10) at (-2.5, 2.25) {};
		\node [style=none] (11) at (-3.25, 2.25) {};
		\node [style=none] (12) at (-3.25, 1.25) {};
		\node [style=none] (13) at (-4.5, 2.25) {};
		\node [style=none] (14) at (-4.5, 1.25) {};
		\node [style=none] (15) at (-3.4, 1.5) {$F$};
		\node [style=none] (16) at (-2.75, 0.25) {$F$};
	\end{pgfonlayer}
	\begin{pgfonlayer}{edgelayer}
		\draw [in=-45, out=135, looseness=1.25] (0) to (1);
		\draw [in=-90, out=90, looseness=1.00] (1) to (2.center);
		\draw [in=-90, out=45, looseness=1.00] (0) to (3.center);
		\draw [in=90, out=-90, looseness=1.00] (0) to (4.center);
		\draw [in=90, out=-150, looseness=0.75] (1) to (5.center);
		\draw (6.center) to (9.center);
		\draw (9.center) to (8.center);
		\draw (8.center) to (7.center);
		\draw (7.center) to (6.center);
		\draw (11.center) to (13.center);
		\draw (13.center) to (14.center);
		\draw (14.center) to (12.center);
		\draw (12.center) to (11.center);
	\end{pgfonlayer}
\end{tikzpicture}
\]

Frobenius functors compose: the composition is defined as the usual composition 
of linear functors \cite{CS97}. 

It is immediate from Lemma \ref{Lemma: linear adjoints} that Frobenius functors preserve linear duals.
 In fact if $F: \X \to \Y$ is a Frobenius functor and $A \dashvv B$ is a linear dual, 
 as the duals $F_\ox(A) \dashvv F_\oa(B)$ and $F_\oa(A) \dashvv F_\ox(B)$ now coincide,  
 we just obtain the one dual $F(A) \dashvv F(B)$.   In the case when the Frobenius functor 
 is between cyclic $*$-autonomous categories we expect the 
functor to be  {\bf cyclor-preserving} in the following sense:
\[ \mbox{\bf [CFF]} ~~~~~\xymatrix{ F(X^{*}) \ar[d]_{\cong} \ar[rr]^{F(\psi)} & & 
 F(\!\!~^{*}X) \ar[d]^{\cong} \\
        F(X)^{*} \ar[rr]_{\psi} & & \!\!~^{*}F(X) } \]
where the left and right vertical arrows are respectively the maps:
\[ (u^R_\ox)^{-1} (\eta* \ox 1) \delta^R (1 \oa (m^F_\ox F(\epsilon*) n^F_\bot) u^R_\oa 
       ~~~\mbox{and}~~~(u^R)^{-1} (1 \ox *\eta) \delta^L (m^F_\ox \oa 1)((F(*\epsilon)n_\bot^F) \oa 1) u^L_\oa \]
The cyclor preserving condition maybe pictorially represented as follows:
\[ 

\]

\begin{lemma}
\label{Lemma: cyclic Frob}
Suppose $F$ is a cyclor preserving Frobenius linear functor, then
\[%
$
\end{proof}

\begin{definition}
Suppose $\X$ and $\Y$ are mix categories. $F: \X \to \Y$ is a {\bf mix functor} if it is a Frobenius functor such that 
\[
\mbox{\bf{[mix-FF]}}~~~~~
\xymatrix{
F(\bot) \ar@/_2pc/[rrr]_{F(\m)} \ar[r]^{n_\bot} & \bot \ar[r]^{\m} & \top \ar[r]^{m_\top} & F(\top) \\
}
\]
\end{definition}

The equation \mbox{\bf{[mix-FF]}} is diagrammatically given as follows:
\[%
 \]

\begin{lemma}
\label{Lemma: Mix Frobenius linear functor}
Mix functors preserve the mix map:
\[
\xymatrix{
F(A) \ox F(B) \ar[r]^{\mx} \ar[d]_{m_\ox} \ar[r]^{\mx} & F(A) \oa F(B) \\
F(A \ox B) \ar[r]_{F(\mx)} \ar[r]_{F(\mx)} & F(A \oa B) \ar[u]_{n_\oa}
}
\]
\end{lemma}
\begin{proof}
\[
%
\]
\end{proof}

Linear natural isomorphisms between Frobenius functors $(\alpha_\ox, \alpha_\oa): F \to G$ often take a special form with $\alpha_\ox = \alpha_\oa^{-1}$: this allows the coherence 
requirements to be simplified.  The next results describe some basic circumstances in which this happens:
 
\begin{lemma}
\label{Lemma: Frobenius linear transformation}
Suppose $F: \X \to \Y$ are Frobenius linear functors and $\alpha := (\alpha_\ox, \alpha_\oa): F \Rightarrow G$ is a linear natural transformation.  Then, the following are equivalent:
\begin{enumerate}[(i)]
\item One of {\bf [nat.1](a)} or {\bf [nat.1](b)} holds, and one of  $\alpha_\ox$ or $\alpha_\oa$ is an isomorphism. \[
\mbox{\bf \small [nat.1]} ~~~
\xymatrixcolsep{5pc}
\xymatrix{
\top \ar[r]^{m_\top} \ar[dr]_{m_\top} & G(\top) \ar[d]^{\alpha_\oa} \ar@{}[dl]|(.35){\tiny{(a)}} \\
& F(\top)
} ~~~~~ \text{ or }~~~~~  \xymatrix{
F(\bot) \ar[r]^{\alpha_\ox} \ar[dr]_{n_\bot} & G(\bot) \ar[d]^{n_\bot} \ar@{}[dl]|(.35){\tiny{(b)}} \\
& \bot
}
\]

\item One of {\bf [nat.1](a)} or {\bf [nat.1](b)} holds and one of the following commuting diagrams holds.
\[ \mbox{\bf \small [nat.2]}~~~~ 
\xymatrix{
G(A) \ox F(B) \ar[r]^{1 \ox \alpha_\ox} \ar[d]_{\alpha_\oa \ox 1} \ar@{}[ddr]|{\tiny{(a)}} & G(A) \ox G(B) \ar[dd]^{m_\ox^G} \\
F(A) \ox F(B) \ar[d]_{m_\ox^F} & \\
F(A \ox B) \ar[r]_{\alpha_\ox} & G(A \ox B)
} ~~~ \text{ or }~~~ \xymatrix{
F(A) \ox G(B) \ar[r]^{\alpha_\ox \ox 1} \ar[d]_{1 \ox \alpha_\oa} \ar@{}[ddr]|{\tiny{(b)}}  & G(A) \ox G(B) \ar[dd]^{m_\ox^G} \\
F(A) \ox F(B) \ar[d]_{m_\ox^F} & \\
F(A \ox B) \ar[r]_{\alpha_\ox} & G(A \ox B)
}
\]
\[
\text{or}~~~~~ 
\xymatrix{
G(A \oa B) \ar[r]^{n_\oa^G} \ar[d]_{\alpha_\oa} \ar@{}[ddr]|{\tiny{(c)}}  & G(A) \oa G(B) \ar[dd]^{1 \oa \alpha_\oa} \\
F(A \oa B) \ar[d]_{n_\oa^F} & \\
F(A) \oa F(B) \ar[r]_{\alpha_\ox \oa 1} & G(A) \oa F(B)
} ~~~\text{or}~~~
\xymatrix{
G(A \oa B) \ar[r]^{n_\oa^G} \ar[d]_{\alpha_\oa} \ar@{}[ddr]|{\tiny{(d)}} & G(A) \oa G(B) \ar[dd]^{\alpha_\oa \oa 1} \\
F(A \oa B) \ar[d]_{n_\oa} & \\
F(A) \oa F(B) \ar[r]_{1 \ox \alpha_\ox} & F(A) \oa G(B)
}
\] 
\item $\alpha_\ox^{-1} = \alpha_\oa$
\item $\alpha' := (\alpha_\oa, \alpha_\ox): G \Rightarrow F$ is a linear transformation.
\end{enumerate}
\end{lemma}

Conditions {\bf [nat.2]} are as follows in the graphical calculus:
\[ \mbox{\small (a)}~~

\]

\begin{proof}
\begin{description}
\item{$(i) \Rightarrow (iii)$:} Here is the proof assuming {\bf \small [nat.1](a)} that $\alpha_\ox \alpha_\oa = 1$:
\[ \hspace{-0.75cm}
%
\]
if either $\alpha_\ox$ or $\alpha_\oa$ are isomorphisms this implies $\alpha_\oa \alpha_\ox =1$.
\item{$(ii) \Rightarrow (iii)$:}
The assumption of {\bf [nat.1](a)} or {\bf (b)} yields, as above, that $\alpha_\ox \alpha_\oa =1$. 
Using {\bf [nat.2](c)} for example gives $\alpha_\oa \alpha_\ox =1$:
\[
%
\]
Since, $\alpha_\ox \alpha_\oa = 1$ and $\alpha_\ox \alpha_\oa = 1$ we have $\alpha_\ox = \alpha_\oa^{-1}$. The other combinations of rules are used in similar fashion.
\item{$(iii) \Rightarrow (iv)$:}
If $\alpha_\ox = \alpha_\oa^{-1}$, then 

$(\alpha_\oa \ox \alpha_\oa) m_\ox^F = m_\ox^G \alpha_\ox: G(A) \ox G(B) \to F(A \ox B)$
\[
%
\]
Thus, $\alpha_\ox$ is comonoidal. Similarly, it can be proven that $\alpha_\oa$ is monoidal. The axioms {\bf \small [LT.4] (a)}-{\bf \small (d)} 
for a linear transformation are satisfied for $(\alpha_\oa, \alpha_\ox)$ because $\alpha_\oa = \alpha_\ox^{-1}$. 

\item{$(iv) \Rightarrow (i)$ and $(ii)$:} The axioms {\bf \small [nat.1]} and {\bf \small [nat.2]} are given by the fact 
that $(\alpha_\oa, \alpha_\ox)$ is a linear transformation. 
\end{description}
\end{proof}

Frobenius functors between isomix categories are especially important in the development of dagger 
linearly distributive categories and they often satisfy an additional property:

\begin{definition}
A Frobenius functor between isomix categories is an {\bf isomix functor} in case it is a mix functor which satisfies, in addition, the following diagram:
\[ \mbox{\bf{[isomix-FF]}}~~~~~\xymatrix{ \top \ar@/^/[rrr]^{{\sf m}^{-1}} \ar[dr]_{m_\top} & & & \bot \\ & F(\top) \ar[r]_{F({\sf m}^{-1})} & F(\bot) \ar[ur]_{n_\top} } \]
\end{definition}

Recall that a linear functor is {\bf normal} in case both $m_\top$ and $n_\bot$ are isomorphisms.  We observe:

\begin{lemma} 
\label{Lemma: isomix functor}
For a mix Frobenius functor, $F: \X \to \Y$, between isomix categories the following are equivalent:
\begin{enumerate}[(i)]
\item $n_\bot: F(\bot) \to \bot$  or $m_\top: \top \to F(\top)$ is an isomorphism;
\item $F$ is a normal functor;
\item $F$ is an isomix functor.
\end{enumerate}
\end{lemma}

\begin{proof}~
\begin{description}
\item{$(i) \Rightarrow (ii)$:}  Note that, as $F$ is a mix functor $F({\sf m}) = n_\bot {\sf m}~m_\top$.  As the mix map ${\sf m}$ is an isomorphism so is $F({\sf m})$ which implies that if  $n_\bot$ is an isomorphism then $m_\top$ must be an isomorphism and vice versa.  Thus, $F$ will be a normal functor. 
\item{$(ii) \Rightarrow (iii)$:} If $F$ is normal then $n_\bot$ and $m_\top$ are isomorphisms and so 
\[ \infer={m_\top F({\sf m}^{-1}) n_\bot = {\sf m}^{-1}}{\infer={F({\sf m}^{-1}) = m_\top^{-1} {\sf m}^{-1} n_\bot^{-1}}{F({\sf m}) = n_\bot {\sf m}~m_\top}} \]
\item{$(iii) \Rightarrow (i)$:} The mix-preservation for $F$ makes  $n_\bot$ a section (and $m_\top$ a retraction) while the isomix-preservation makes $m_\bot$ a retraction (and $m_\top$ a section. 
This means $n_\bot$ is an isomorphism ($m_\top$ is an isomorphism).                                                                                                                                                                                                                                                                                                                                                                                                                                                                         
\end{description}
\end{proof}

\begin{corollary} \label{Corollary: normal-nat-iso}
$\alpha := (\alpha_\ox, \alpha_\oa)$ is a linear natural isomorphism between isomix Frobenius linear functors if and only if $\alpha_\ox = \alpha_\oa^{-1}$.
\end{corollary}

\begin{proof}
Note  that if we can establish {\bf [nat.1](a)} or {\bf (b)} then we can prove that $\alpha_\ox\alpha_\oa =1$ and, as $\alpha_\ox$ is an isomorphism it follows that $\alpha_\oa\alpha_\ox =1$.
Thus, it suffices to show that {\bf [nat.1](a)} holds:
\[ m_\top \alpha_\oa G({\sf m}^{-1}) n_\bot = m_\top F({\sf m}^{-1})\alpha_\oa  n_\bot = m_\top F({\sf m}^{-1})  n_\bot = {\sf m}^{-1} = m_\top G({\sf m}^{-1}) n_\bot  \]
However, as $G({\sf m}^{-1}) n_\bot$ is an isomorphism, it follows that $m_\top \alpha_\oa = m_\top$.
\end{proof}


Lemma \ref{Lemma: Frobenius linear transformation} and Corollary \ref{Corollary: normal-nat-iso} are generalizations of \cite[Proposition 7]{DP08}. \cite[Proposition 7]{DP08} states the following:

Let $\X$ and $\Y$ be monoidal categories and $(\eta, \epsilon): A \dashv B \in \X$. If $F,G: \X \to \Y$ are Frobenius monoidal functors with a natural transformation $\alpha: F \Rightarrow G$ which is both monoidal and comonoidal, then $\alpha_A$ is invertible. 

In Lemma \ref{Lemma: Frobenius linear transformation}, when $A \dashvv B \in \X$, then $\alpha_\oa$ is defined as follows:
\[
\alpha_\oa: G(A) \to F(A) = \begin{tikzpicture}
	\begin{pgfonlayer}{nodelayer}
		\node [style=circle] (0) at (-0.25, -0) {$\alpha_\ox$};
		\node [style=none] (1) at (1.25, -2.75) {};
		\node [style=none] (2) at (-0.25, -1) {};
		\node [style=none] (3) at (-1.75, -1) {};
		\node [style=none] (4) at (-1.75, 2.75) {};
		\node [style=circle] (5) at (-1, -1.75) {$\epsilon$};
		\node [style=none] (6) at (-0.25, 1) {};
		\node [style=none] (7) at (1.25, 1) {};
		\node [style=circle] (8) at (0.5, 1.75) {$\eta$};
		\node [style=none] (9) at (1.5, 2.5) {};
		\node [style=none] (10) at (1.5, 1) {};
		\node [style=none] (11) at (-0.5, 1) {};
		\node [style=none] (12) at (-0.5, 2.5) {};
		\node [style=none] (13) at (-2, -0.75) {};
		\node [style=none] (14) at (0, -0.75) {};
		\node [style=none] (15) at (0, -2.5) {};
		\node [style=none] (16) at (-2, -2.5) {};
		\node [style=none] (17) at (-2.25, 2.5) {$G(A)$};
		\node [style=none] (18) at (1.75, -2.5) {$F(A)$};
		\node [style=none] (19) at (0.25, 0.65) {$F(B)$};
		\node [style=none] (20) at (0.5, -0.5) {$G(B)$};
		\node [style=none] (21) at (1.25, 2.25) {$F$};
		\node [style=none] (22) at (-0.25, -2.25) {$G$};
	\end{pgfonlayer}
	\begin{pgfonlayer}{edgelayer}
		\draw (0) to (2.center);
		\draw (3.center) to (4.center);
		\draw [in=0, out=-90, looseness=1.25] (2.center) to (5);
		\draw [in=180, out=-90, looseness=1.25] (3.center) to (5);
		\draw (11.center) to (10.center);
		\draw (12.center) to (11.center);
		\draw (10.center) to (9.center);
		\draw (9.center) to (12.center);
		\draw (7.center) to (1.center);
		\draw [in=90, out=0, looseness=1.25] (8) to (7.center);
		\draw (6.center) to (0);
		\draw [in=180, out=90, looseness=1.25] (6.center) to (8);
		\draw (13.center) to (16.center);
		\draw (16.center) to (15.center);
		\draw (15.center) to (14.center);
		\draw (14.center) to (13.center);
	\end{pgfonlayer}
\end{tikzpicture}
\]

For these special linear isomorphisms with $\alpha_\ox = \alpha_\oa^{-1}$ we can simplify the coherence requirements:

\begin{lemma} \label{simplifying-coherences}
Suppose $F$ and  $G$ are Frobenius functors and $\alpha: F \to G$ is a natural isomorphism then:
\begin{enumerate}[(i)]
\item If $\alpha: F \to G$ is $\ox$-monoidal and $\oa$-comonoidal then $(\alpha,\alpha^{-1})$ is a linear transformation;
\item If $F$ and $G$ are strong Frobenius functors and $\alpha$ is $\ox$-monoidal and $\oa$-monoidal then $(\alpha,\alpha^{-1})$ is a linear transformation.
\end{enumerate}
\end{lemma}

\begin{proof}~
\begin{enumerate}[{\em (i)}]
\item If $\alpha$ is $\ox$-monoidal and $\oa$-comonoidal then so is $\alpha^{-1}$ supporting the possibility that it is a component of a linear transformation.
Considering  {\bf [LT.1]} we show that $(\alpha,\alpha^{-1})$ satisfies this requirement as:
\[ \xymatrix{ F(A \oa B) \ar[rr]^{\alpha_\ox}  \ar[d]_{n_\oa = \nu^R_\ox} & & G(A \oa B) \ar[d]^{n_\oa = \nu^R_\ox} \\
                    F(A) \oa F(B) \ar[dr]_{1 \oa \alpha} \ar[rr]^{\alpha \oa \alpha} & & G(A) \oa G(B)  \ar[dl]^{\alpha^{-1} \oa 1} \\
                    & F(A) \oa G(B) } \]
The remaining requirements follow similarly.   
\item  When the laxors for the functors are isomorphisms then being monoidal implies being comonoidal.   \qedhere
\end{enumerate}
\end{proof}


\section{Motivating examples}
\label{Sec: motivating examples}

In Section \ref{Sec: LDC}, we listed a few examples of linearly distributive categories.
In this section, we present the categories of finiteness relations, ${\sf FRel}$ and finiteness matrices
over a commutative rig $R$, $\FMat(R)$, and discuss their properties. These isomix categories are closely related 
to categorical quantum mechanics because the core of these categories are used as the standard 
categorical models to study quantum processes in finite dimensions: the core of $\FMat(R)$ is equivalent 
to the category finite dimensional matrices over a commutative rig $R$, and the core of $\FRel$ is the category of finite sets and 
relations.  

The categories $\FRel$ and $\FMat(R)$ are used as running examples 
throughout this thesis. Consequently, the purpose of this section 
is to revisit the results of \cite{Ehr05} in which these categories were introduced.

\subsection{Finiteness relations ${\sf FRel}$ and finiteness matrices ${\sf FMat}(R)$}

Finiteness spaces were introduced by Ehrhard \cite{Ehr05} as a categorical model of Girard's linear logic \cite{Gir87}.  
A finiteness space $(X,{\cal A},{\cal B})$ may be regarded as a set, $X$, called the {\bf web} of $X$, equipped 
with two sets of subsets, ${\cal A} \subseteq {\cal P}(X)$ the {\bf finitary} sets of $X$, and ${\cal B} \subseteq {\cal P}(X)$ 
the {\bf cofinitary} sets of $X$.  The finitary and cofinitary sets must be finiteness complements of each other 
(as described below) and so determine each other.  This means that a finiteness space is often presented 
with just the finitary sets.  Morphisms of a finiteness spaces are relations between the webs that preserve 
the finitary sets and reflect the cofinitary sets  The category of finiteness spaces with finiteness relations is 
an isomix $*$-autonomous category.  In this section, we aim to give a reasonably self-contained account of $\FRel$, 
the category of finitary relations, and of $\FMat(R)$ the category of finiteness matrices with entries in a commutative rig $R$, 
and to establish that these categories are symmetric $\dagger$-$*$-isomix categories. 

Let $X$ be any set and $A, B \subseteq X$ then $A$ is {\bf finiteness orthogonal} \cite{Ehr05}
\cite[Section 1]{Ehr05} to $B$, written $A \perp_f B$, in case $A \cap B$ is finite.   
Given a set of subsets of $X$, ${\cal A} \subseteq {\cal P}(X)$ set,
\[ {\cal A}^\perp = \{ B \subseteq X \mid \forall A \in {\cal A}, B \perp_f A\} \]
Thus ${\cal A}^\perp$ is the set of all subsets of $X$ which intersects with all the sets in ${\cal A}$ finitely.  Observe that:
 
 \begin{lemma} For ${\cal A}, {\cal B} \subset {\cal P}(X)$:
\begin{enumerate}[(i)]
\item ${\cal A} \subseteq {\cal B} \Rightarrow {\cal B}^{\perp} \subseteq {\cal A}^{\perp}$;
\item ${\cal A} \subseteq {\cal A}^{\perp \perp}$;
\item ${\cal A}^{\perp \perp \perp} = {\cal A}^{\perp}$;
\item ${\cal A}^\perp$ is downset closed and closed under finite unions;
\item if ${\cal A} = {\cal A}^{\perp\perp}$ then ${\cal P}_f(X) \subseteq {\cal A}$ (that is ${\cal A}$ contains all finite subsets);
\item If ${\cal A}_i= {\cal A}_i^{\perp\perp}$ for $i \in I$ then $\bigcap_{i \in I} {\cal A}_i = (\bigcup_{i \in I} {\cal A}^\perp_i)^{\perp}$ \\
(so $\bigcap_{i \in I} {\cal A}_i =  (\bigcap_{i \in I} {\cal A}_i)^{\perp\perp}$).
 \end{enumerate}
 \end{lemma}
 
The first two observations establish a Galois connection on the subsets of ${\cal P}(X)$ from which the next observation is standard.  The last two are easy consequences of the form of the finiteness orthogonality.
 
\begin{definition}
A {\bf finiteness space} is a triple $(X, {\cal A},{\cal B})$ where $X$ is a set, 
called the {\bf web} of the space, and ${\cal A}, {\cal B}  \subseteq {\cal P}(X)$, 
such that the two subsets are finiteness complements, that is 
${\cal A}^\perp = {\cal B}$ and ${\cal B}^\perp = {\cal A}$.   
Elements of ${\cal A}$ are the {\bf finitary} sets of the finiteness space, 
while elements of ${\cal B}$ are called the {\bf cofinitary} sets. 
\end{definition}

Because the cofinitary sets are completely determined by the finitary sets it is often convenient 
to write a finiteness space as $(X,F(X),F(X)^\perp)$, 
making clear that the structure is completely determined by the web and the finitary sets.  
The finitary and cofinitary sets are downward closed, that is $F(X) = {\downarrow}F(X)$, 
and are closed under arbitrary intersections and finite unions.  As finite sets always intersect 
finitely both with the finitary and cofinitary sets, the finitary and cofinitary sets of every 
finiteness space must always include all the finite subsets.  In particular this means that a 
finiteness set with a finite web must have every subset both finitary and cofinitary.  Furthermore, 
every set, $X$, always carries two ``trivial'' finiteness structures (which coincide for a finite set): 
$(X,{\cal P}_f(X),{\cal P}(X))$ and  $(X,{\cal P}(X),{\cal P}_f(X))$, where ${\cal P}_f(X)$ 
is the set of finite subsets of $X$. On the other hand, as we now 
observe, an infinite set always has infinitely many non-trivial finiteness structures: 

\begin{lemma}~
\label{Lemma: finite finite}
\begin{enumerate}[(i)]
\item  A finiteness space, $(X, F(X), F(X)^\perp)$, with a finite web, 
has its finitary and cofinitary sets completely determined to be 
$F(X) = F(X)^\perp = {\cal P}(X) = {\cal P}_f(X)$;
\item  Every infinite set carries infinitely many non-trivial finiteness structures;
\item The set of finiteness structures carried by $X$ form a complete (non-distributive) lattice with $(X,{\cal A},{\cal B}) \leq (X,{\cal A}',{\cal B}')$ if and only if ${\cal A} \subseteq {\cal A}'$ 
(equivalently ${\cal B}' \subseteq {\cal B}$).  The lattice structure is given by:
\begin{align*}
\top & := (X,{\cal P}(X),{\cal P}_f(X)) &
\bot & := (X,{\cal P}_f(X),{\cal P}(X)) \\
\bigwedge_i A_i & := \left(X, \bigcap_i {\cal A}_i, \left(\bigcup_i {\cal A}_i \right)^{\perp\perp} \right) &
\bigvee_{i \in I} A_i & := \left(X, \left(\bigcup_i {\cal A}_i \right)^{\perp\perp}, \bigcap_i {\cal B}_i \right) 
\end{align*}
\end{enumerate}
\end{lemma}

\begin{proof}~
\begin{enumerate}[{\em (i)}]
\item This is immediate.
\item There are, of course, many different possible finiteness structures for an infinite set: to see this consider the finiteness space $(X,\{ A\}^{\perp\perp},\{ A \}^\perp)$, generated by insisting that an infinite subset, $A \subseteq X$, whose complement, $\neg A$, is also infinite, is finitary.  Notice that, in this case:
\[ \{ A \}^\perp = \{ Y \mid Y \cap A~\mbox{is finite}\} = \{ Y \mid Y \subseteq \neg A \cup F~\mbox{where $F$ is finite}\} \]
\[ \{ A\}^{\perp\perp} = \{ X \mid X \subset A \cup F~\mbox{where $F$ is finite}\} \]
It is then sufficient to show that there are infinitely many of these sets which are distinct by more that a finite set:  for this note that having chosen an infinite set one can split it again 
recursively in the same manner.
\item It suffices to show that the two sets $\bigwedge_{i \in I} {\cal A}_i$ and $\bigvee_{i \in I}{\cal B}_i$ are complementary.  Note that each set $B_i \in {\cal B}_i$ is in $\bigvee_{i \in I}{\cal B}_i$ so any set in the complement must be in the complement of each ${\cal B}_i$ and so must be in $\bigwedge_{i \in I} {\cal A}_i$.  Conversely, a set in $\bigwedge_{i \in I} {\cal A}_i$ is necessarily orthogonal to each finite union $\bigcup_{i \in F} B_i$ as it intersects with each $B_i$ finitely.  Thus, $(\bigvee_{i \in I}{\cal B}_i)^\perp =\bigwedge_{i \in I} {\cal A}_i$.
\end{enumerate}
\end{proof}

To complete the description of the category ${\sf FRel}$, we describe the morphisms between these spaces:
\begin{definition}
A {\bf finiteness relation} $R: (X, F(X),F(X)^\perp) \to (Y, F(Y),F(Y)^\perp)$ between two finiteness spaces, is a relation $X \to^{R} Y$ such that:
\begin{enumerate}[(i)]
\item for all $A \in F(X), ~~ A \triangleright R := \{ y \mid \exists x \in A \text{ such that } xRy \} \in F(Y)$
\item for all $B \in {F(Y)}^\perp, ~~ R \triangleleft B := \{ x \mid \exists y \in B \text{ such that } xRy \} \in F(X)^\perp$ 
\end{enumerate}
\end{definition}

Here we are using a notation which can be translated into allegorical notation, 
see \cite{FrS89}, where it is possible to do most of these proofs more generally 
and so can be made formal for general categories of relations: we, of course, 
will have in mind sets and relations.  In this notation a finitary set in $F(X)$ 
is a subrelation of the identity relation $A \subseteq 1_{X}$ 
while a cofinitary set, $B \in  F(X)^\perp$, is a subrelation $B \subseteq 1_{Y}$.  
The orthogonality relation $A \perp_f A'$ is then the requirement that 
$A \cap A' = AA'$ is finite.  A finiteness relation is a map $R: X \to Y$ with 
$A \triangleright R:= {\sf cod}(AR) \in F(Y)$, for every $A \in F(X)$, and 
$R \triangleleft B:= {\sf dom}(RB) \in F(X)^\perp$ for every $B \in F(Y)^\perp$ 
(recall that ${\sf dom}(R) = 1_{X} \cap RR^\circ$ and ${\sf cod}(R) = 
1_{Y} \cap R^\circ R$, where $R^\circ$ is the converse of $R$).  

Finite relations are determined by having finite domains and codomains: 
so that $R: X \to Y$ is finite if and only if ${\sf cod}(R) \subseteq_f 1_{Y}$ 
and ${\sf dom}(R) \subseteq_f 1_{X}$.  Thus $R$ is a finiteness relation if and 
only if $ARB$ is a finite relation for every $A \in F(X)$ and $B \in F(Y)^\perp$.

Notice that once we require that $A \triangleright R \in F(Y)$ for every 
$A \in F(X)$ then each $B \in F(Y)^\perp$ has a finite intersection with 
$A \triangleright R$ so that 
there are finitely many $b \in (A\triangleright R)\cap B$ and whence
\begin{align*}
 ARB & = AR((A \triangleright R)\cap B)  =  AR\bigcup_{b \in A \triangleright R \cap B} \{ b\} \\
 & =  \bigcup_{b \in A \triangleright R \cap B} AR\{b\}  =  
 \bigcup_{b \in A \triangleright R \cap B} (A \cap R \triangleleft \{b\}) R \{b\} 
 \end{align*}
This means that so long as $R \triangleleft \{b\} \in F(X)^\perp$ each $A \cap R \triangleleft \{b\}$ will be finite, so, as the union over $b \in A \triangleright R \cap B$ is a finite union, this ensures that $ARB$ is a finite relation and, thus, that $R$ is a finiteness relation.

This discussion gives the following important simplifications (see Ehrhard \cite{Ehr05}) of the conditions of being a finiteness relation:

\begin{lemma} \label{character_finiteness_relation}
The following  are equivalent:
\begin{enumerate}[{(i)}]
\item $R: X \to Y$ is a finiteness relation;
\item $ARB$ is a finite relation for all $A \in F(X)$ and $B \in F(Y)^\perp$;
\item $R: X \to Y$ has $A \triangleright R \in F(Y)$ for all $A \in F(X)$ and 
for every $y \in Y$, $R \triangleleft \{ y \} \in F(X)^\perp$;
\item $R: X \to Y$ has $R \triangleleft B \in F(X)^\perp$ for all 
$B \in F(Y)^\perp$ and for every $x \in X$, $\{ x \} \triangleright R \in F(Y)$.
\end{enumerate}
\end{lemma}

\begin{proof}
The above discussion proves that {\em (i)\/} is equivalent to {\em (ii)\/}. That {\em (ii)\/} implies {\em (iii)\/} is now immediate (by converse symmetry this also gives {\em (ii)\/} implies {\em (iv)\/}).  Finally we have shown above that {\em (iii)} implies {\em (ii)\/} (and by converse symmetry that {\em (iv)\/} implies {\em (ii)\/}).
\end{proof}

The category of finiteness relations, ${\sf FRel}$, has objects finiteness spaces and maps finiteness relations with composition as in the category of relations:

\begin{lemma}
Finiteness spaces with finiteness relations, ${\sf FRel}$,  forms a category with composition 
the usual composition of relations and identities the diagonal relations.
\end{lemma}

\begin{proof}
Suppose $X \to^{R} Y$, and $Y \to^{S} Z$ are finiteness relations.  Then, for $A \in F(X)$ we know $A \triangleright R \in F(Y)$ 
and whence $A\triangleright (RS) = (A \triangleright R) \triangleright S  \in F(Z)$.  Similarly for $B \in F(Z)^\perp$ we have $S \triangleleft B \in F(Y)^\perp$ and, 
therefore, $(RS) \triangleleft B = (R \triangleleft (S \triangleleft B)) \in F(X)^\perp$. 

Clearly the identity relation is always a finiteness relation.
\end{proof}

We record an obvious yet important fact: there is a faithful functor $|\_|: {\sf FRel} \to {\sf Rel}$ 
which takes a finiteness space to its web and a finiteness relation to its underlying relation.  
This functor, we shall see, preserves the structure of ${\sf FRel}$ as an isomix *-autonomous category 
on the nose so that we can use the known coherence properties of ${\sf Rel}$ to 
obtain the corresponding coherence properties of ${\sf FRel}$.

With finiteness relations in hand one can build the category of finiteness matrices, 
${\sf FMat}(R)$, for an arbitrary commutative rig, $R$. While $R$ can be any commutative rig, 
we are most interested in the case when $R$ is the complex numbers, and, thus, in the 
category ${\sf FMat}(\C)$. 

The category ${\sf FMat}(R)$, where $R$ is any commutative rig is defined as follows:
\begin{description}
\item[Objects:]  Finiteness spaces;
\item[Maps:]  Finiteness matrices $M: X \to Y$ where 
$M = [M_{i ,j}]_{i \in X, j \in Y}$ where each $M_{i,j} \in R$ and the 
support, $|M|$, is a finiteness relation, that is $|M| =\{ (i,j) \mid M_{i,j} \neq 0\}$ is a finiteness relation;
\item[Composition:] Matrix multiplication -- where we note all the non-zero sums are finite (see Lemma \ref{Lemma: finiteness matrix} below),  $\sum_j r_{i,j}s_{j,k} = \sum_{j \in \{i\}|r| \cap |s|\{ k \}}r_{i,j}s_{j,k} .$
\item[Identities:]  Identity matrices.
\end{description}

Observe that, although finiteness matrices may be infinite dimensional, to form the composition of two finiteness matrices only requires the ability to multiply a finite number of entries and, thus, to sum only a finite number of elements.  This is because the support of each matrix is always a finiteness relation, thus, for composition one only need compute the products of non-zero elements and these lie in the intersection of a finite and cofinite set, so is finite.  This explicitly is the observation:

\begin{lemma}
\label{Lemma: finiteness matrix}
If $X \to^{R} Y \to^S Z$ is a finiteness relation,  then 
\begin{enumerate}[(i)]
\item for all $x \in X$, $\{ x \} \triangleright R \in F(Y)$, and 
\item for all $y \in Y$, $S \tl \{ z \} \in F(Y)^\perp$
\item for all $x \in X$ and $z \in Z$, $\{ x \} \tr R \perp_f S \tl \{ z \}$, that is $\{ x \} \tr R \cap S \tl \{ z \}$ is finite.
\end{enumerate}
\end{lemma}

In a finiteness matrix $X \to^{M} Y$,  each row of the matrix has support a finitary set of $Y$, and each column has a cofinitary set of $X$ as shown in Figure \ref{Diagram: finiteness matrix}. 

\begin{figure}
\centering
\begin{tikzpicture}[scale=1.5]
	\begin{pgfonlayer}{nodelayer}
		\node [style=none] (0) at (-1.5, 2) {};
		\node [style=none] (1) at (-1.5, -1) {};
		\node [style=none] (2) at (1.25, -1) {};
		\node [style=none] (3) at (1.25, 2) {};
		\node [style=none] (4) at (-0.75, 2) {};
		\node [style=none] (5) at (-0.25, 2) {};
		\node [style=none] (6) at (-0.75, -1) {};
		\node [style=none] (7) at (-0.25, -1) {};
		\node [style=none] (8) at (-1.5, 1) {};
		\node [style=none] (9) at (-1.5, 0.5) {};
		\node [style=none] (10) at (1.25, 1) {};
		\node [style=none] (11) at (1.25, 0.5) {};
		\node [style=none] (12) at (2, 0.75) {$\in F(Y)$};
		\node [style=none] (13) at (-0.5, 2.75) {$\in F(X)^\perp$};
		\node [style=none] (14) at (0.5, 1.5) {$\cdots$};
		\node [style=none] (15) at (-1, 1.5) {$\cdots$};
		\node [style=none] (16) at (0.5, -0) {$\cdots$};
		\node [style=none] (17) at (-1, -0) {$\cdots$};
		\node [style=none] (18) at (-2, 0.5) {$|X|$};
		\node [style=none] (19) at (0, -1.5) {$|Y|$};
		\node [style=none] (20) at (0.5, -1.25) {};
		\node [style=none] (21) at (1.25, -1.25) {};
		\node [style=none] (22) at (-0.5, -1.25) {};
		\node [style=none] (23) at (-1.5, -1.25) {};
		\node [style=none] (24) at (-2, 1) {};
		\node [style=none] (25) at (-2, 2) {};
		\node [style=none] (26) at (-2, -0) {};
		\node [style=none] (27) at (-2, -1) {};
		\node [style=none] (28) at (-0.5, 2.15) {};
		\node [style=none] (29) at (-0.5, 2.45) {};
	\end{pgfonlayer}
	\begin{pgfonlayer}{edgelayer}
		\draw (8.center) to (10.center);
		\draw (9.center) to (11.center);
		\draw (0.center) to (1.center);
		\draw (1.center) to (2.center);
		\draw (3.center) to (2.center);
		\draw (3.center) to (0.center);
		\draw (4.center) to (6.center);
		\draw (7.center) to (5.center);
		\draw [->] (24.center) to (25.center);
		\draw [->] (26.center) to (27.center);
		\draw [->] (22.center) to (23.center);
		\draw [->] (20.center) to (21.center);
		\draw[->] (28.center) to (29.center);
	\end{pgfonlayer}
\end{tikzpicture}
\caption{Finiteness matrix}
\label{Diagram: finiteness matrix}
\end{figure}
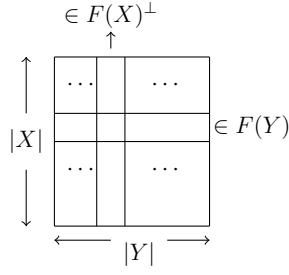

\begin{lemma}
Finiteness spaces with finiteness matrices over any commutative rig, $R$, form a category, ${\sf FMat}(R)$.
\end{lemma}

\begin{proof}
We must show that finiteness matrices compose.   We already know that, even though finiteness matrices can be infinite-dimensional, two such matrices can be multiplied using only finite sums of elements because their supports are finiteness relations. However, we still have to show that the product of two finiteness matrices is a finiteness matrix.

 Suppose $X \to^{P} Y$ and $Y \to^{Q} Z$, $[PQ]_{x,z} = \sum_{y \in Y} [P]_{x,y} [Q]_{x,y}$  can only be non-zero if $\exists ~ y \in Y$, such that $[P]_{x,y}$, and $[Q]_{y,z}$ are non-zero, although, even if there exists $y \in Y$, such that $[P]_{x,y}$ and $[Q]_{y,z}$ are non-zero, $[PQ]_{x,z}$ may still be zero. This means the support $|PQ| \subseteq |P||Q|$. Since, any subset of a finitary relation is finitary, $PQ$ is a finiteness matrix.
\end{proof}

The proof shows that taking the support of matrices gives, in general, a colax 2-functor $|\_|: {\sf FMat}(R) \to {\sf FRel}$  as $|MN| \leq |M||N|$.  In general if $f: R \to S$ is a morphism of rigs then we obtain a functor ${\sf FMat}(f): {\sf FMat}(R) \to {\sf FMat}(S)$: if $S$ is an ordered rig and $f$ is suitably colax then ${\sf FMat}(f)$ is a colax functor.  Recall that the two element lattice $\mathbbm{2}$ with join as addition and meet as multiplication has ${\sf FMat}(\mathbbm{2}) = {\sf FRel}$.

\subsection{Products, coproducts and biproducts in ${\sf FRel}$ and ${\sf FMat}(R)$}

We now begin to explore the properties of ${\sf FRel}$, and $\FMat(R)$.  First we observe that there is an obvious involution on ${\sf FRel}$ given by swapping the finitary and cofinitary sets:
\[ (\_)^{*}: {\sf FRel}^{\rm op} \to {\sf FRel}; \begin{matrix} \xymatrix{ (X,{\cal A},{\cal B}) \ar@{|->}[rr] \ar[dd]_R  & & X^{*} = (X,{\cal B},{\cal A})  \\ 
                                                                                                                ~~~~~~~~~~~~ \ar@{|->}[rr] & &~~~~~~~~~~~ \\
                                                                                                                (X',{\cal A}',{\cal B}' ) \ar@{|->}[rr] & & {X'}^{*} = (X',{\cal B}',{\cal A}' ) \ar[uu]_{R^\circ} } \end{matrix} \] 
For ${\sf FMat}(R)$ this involution is given by transposing the matrix:
\[ (\_)^{*}: {\sf FMat}(R)^{\rm op} \to {\sf Mat}(R); \begin{matrix} \xymatrix{ (X,{\cal A},{\cal B}) \ar@{|->}[rr] \ar[dd]_{\left(r_{i,j}\right)_{i \in X,j \in X'}} & & X^{*} = (X,{\cal B},{\cal A})  \\ 
                                                                                                                ~~~~~~~~~~~~ \ar@{|->}[rr] & &~~~~~~~~~~~ \\
                                                                                                                (X',{\cal A}',{\cal B}' ) \ar@{|->}[rr] & & {X'}^{*} = (X',{\cal B}',{\cal A}' ) \ar[uu]_{\left(r_{j,i}\right)_{j \in X', i \in X}} } \end{matrix} \] 
This involution clearly has the property that   $(X^*)^* = X$.  This symmetry means that  we can take a ``one-sided'' approach to establishing structures using the involution to automatically produce the dual structure.  Shortly we will see that this involution is actually the duality of the $*$-autonomous structure of both categories.                                                                                                   

The first observation is:

\begin{lemma}
\label{Lemma: prod and coprod}
 ${\sf FRel}$ and ${\sf FMat}(R)$, for any commutative rig $R$, have arbitrary products and coproducts.
\end{lemma}

\begin{proof}
Let $(X_i, F(X_i), F(X_i)^\perp)$ be a set of finiteness spaces. 
We shall work for the coproduct as the argument for the product is dual.  The coproduct  $\coprod_{i \in I} X_i$ of finiteness spaces 
has its web given by the disjoint union of the webs with injections 
$\sigma_i: X_i \to \sqcup_{i \in I} X_i$.  
Its finitary sets are given by ${\cal A} = 
\{ \bigcup_{i \in I'} A_i\sigma_i \mid A_i \in F(X_i) ~\&~ I' \subseteq_f I \}$ 
while the cofinitary sets are given by ${\cal B} = \{ \bigcup_{i \in I} B_i \sigma_i \mid B_i \in F(X_i)^\perp \}$. 
First note that the injections $\sigma_i: X_i \to \sqcup_{i \in I} X_i$ are certainly finiteness relations.  As a relation the comparison map for a family of relations, $R_i: X_i \to Y$, is unique 
and given by $\bigcup_{i \in I} \sigma_i^\circ R_i : \sqcup_{i \in I} X_i \to Y$: we must check this is a finiteness relation assuming each $R_i$ is.  Note that if $B \in F(Y)^\perp$ then 
\[ B \left( \bigcup_{i \in I} \sigma_i^\circ R_i \right)^\circ = \bigcup_{i \in I} B R_i^\circ \sigma_i \]
where $B R_i^\circ \in F(X_i)^\perp$ so this is in ${\cal B}$.  Conversely, given $\bigcup_{i \in I'} A_i\sigma_i \in {\cal A}$ mapping this forward gives:
\[ \left(\bigcup_{i \in I'} A_i\sigma_i \right) \left( \bigcup_{i \in I} \sigma_i^\circ R_i \right) = \bigcup_{i \in I'} A_i R_i \]
which is a finite union of finitary sets of $Y$ and so a finitary set.

Finally, we must check that, as we have defined it, $\coprod_{i \in I} X_i$ is a finiteness set, that is ${\cal A}^\perp = {\cal B}$ and ${\cal B}^\perp = {\cal A}$.  Toward this end 
suppose $H \subseteq \bigsqcup_{i \in I} X_i$ is orthogonal to ${\cal B}$ then, 
restricting to each $X_i$, $H \sigma_i^\circ \cap A_i$ is finite for each $A_i \in F(X_i)^\perp$ and 
there can only be finitely $i$ for which this is non-empty.  
But this means $H \in {\cal A}$.  Conversely, suppose $K \subseteq \bigsqcup_{i \in I} X_i$ is orthogonal to ${\cal A}$ then 
restricting to $X_i$ means that each $K \sigma_i^\circ \in F(X_i)^\perp$ which means 
$K = \bigcup_{i \in I} B_i \sigma$ and so is in ${\cal B}$.

For ${\sf FMat}(R)$ the coproduct has the same underlying finiteness space: it is then straightforward to check that all the required maps have the appropriate support.
\end{proof}

Note that in both ${\sf FRel}$ and ${\sf FMat}(R)$ finite products are the same as finite coproducts. 
Thus, both categories have biproducts and so are additively enriched (that is, they are enriched in commutative monoids). 

\begin{corollary}
Both ${\sf FRel}$ and ${\sf FMat}(R)$ have (finite) biproducts.
\end{corollary}

In ${\sf FRel}$ the additive enrichment is given by the union of the underlying relations.  In ${\sf FMat}(R)$ the additive enrichment is given by the pointwise addition of matrices.


\subsection{The tensor structure in ${\sf FRel}$ and ${\sf FMat}(R)$}

${\sf FRel}$ has a symmetric tensor product and this can be used to obtain a corresponding tensor product on ${\sf FMat}(R)$.   We shall start by focusing on ${\sf FRel}$ and define 
$X \ox Y := (X \times Y, F(X \ox Y), F(X \ox Y)^\perp)$, where
\[
F(X \ox Y) := {\downarrow} \{ A \x B \mid  A \in F(X), B \in F(Y) \}
\]
where ${\downarrow} {\cal A}:= \{ A' \mid A' \subseteq A, A \in {\cal A} \}$ is the downward closure of the set of subsets ${\cal A}$.  To show this is well-defined we must prove  that $F(X \ox Y) = F(X \ox Y)^{\perp\perp}$ (this is essentially \cite[Lemma 2]{Ehr05}): 

\begin{lemma} \label{Lemma: tensor characterization}
For any finiteness spaces $X$ and $Y$,  $X \ox Y = (X,F(X \ox Y),F(X \ox Y)^\perp)$, as defined above, is a well-defined finiteness space.
\end{lemma}

\begin{proof}
We must show that any $Q \subseteq X \x Y$ which intersects finitely with any $P \subseteq X \x Y$ 
which, in turn, intersects  finitely with all $R \in F(X \ox Y)$ must already be in 
$F(X \ox Y)$.   Suppose for contradiction that $Q$ is not in $F(X \ox Y)$, as defined;  
this means there is no $A \in F(X)$ and $B \in F(Y)$ such that $Q \subseteq A \x B$.  
This, in turn, means that ${\sf dom}(Q) \not\in F(X)$ or ${\sf cod}(Q) \not\in F(Y)$. 
Without the loss of generality, we may assume the former: 
this means that there is a $C \in F(X)^\perp$ with $C \cap {\sf dom}(Q)$ infinite.  
We may ``thicken'' $C \cap {\sf dom}(Q)$ to $C': X \to Y $ by choosing for each 
$c \in C \cap {\sf dom}(Q)$ a $y_c \in Y$ so that $(c,y_c) \in Q$.  
Then $C' \cap Q$ is certainly infinite, however, $C' \in F(X \ox Y)$ as 
intersecting with any $A \x B$ with $A \in F(X)$ and $B \in F(Y)$ is finite as $C \cap A$ is finite.   
But this means $Q$ cannot be in $F(X \ox Y)^{\perp\perp}$.
\end{proof}

By the discussion above Lemma \ref{character_finiteness_relation}, as $Q \in F(X \cap Y)^\perp$ precisely when $Q \cap A \x B$ is finite for all $A \in F(X)$ and $B \in F(Y)$, $Q$ may be  characterized as a finiteness relation $Q: X^{*} \to Y$ or $Q^\circ: Y^{*} \to X$.  This immediately gives:

\begin{lemma}
\label{Lemma: tensor perp characterization}
Let $X$ and $Y$ be a finiteness spaces then $P \subseteq X \x Y$ then the following are equivalent:
\begin{enumerate}[(i)]
\item $P \in F(X \ox Y)^\perp$;
\item $P \cap A \x B$ is finite for all $A \in F(X)$ and $B \in F(Y)$;
\item For all $A \in F(X)$,  $A \triangleright P  \in F(Y)^\perp$ and,  for all $B \in {F(Y)}$, $P \triangleleft B \in F(X)^\perp$;
\item For all $A \in F(X)$,  $A \triangleright P  \in F(Y)^\perp$ and, for every $y \in Y$, $P \triangleleft \{ y\} \in F(X)^\perp$;
\item For every $x \in X$, $\{ x\} \triangleright P \in F(Y)^\perp$ and, for all $B \in {F(Y)}$, $P \triangleleft B \in F(X)^\perp$.
\end{enumerate}
\end{lemma}

The unit of tensor is given by the one element set with only finiteness structure possible:  $\top :=  (\{ \star \}, {\cal P}(\{ \star \}), {\cal P}(\{ \star \}) )$.

Given finiteness relations, $X \to^{R} X'$, and $Y \to^{S} Y'$, then: 
\[ R \ox S := \{ ((x,y), (x',y')) \mid xRx', ySy' \}\]  
as in ${\sf Rel}$, is clearly a finiteness relation.  The associator and unitors for this tensor product are the same as in $\Rel$.  
Thus, the associator is the following relation: 
\[ a_\ox := {(((x,y),z), (x,(y,z))) \mid x \in X, y \in Y, z \in Z}\] 
and this clearly gives a bijection between 
$F((X \ox Y) \ox Z) = {\downarrow}\{ ((a,b),c) \mid a \in A \in F(X), b \in B \in F(Y), c \in C \in F(Z) \}$, and $F(X \ox (Y \ox Z)) = \downarrow\{ ((a,(b,c)) \mid a \in A \in F(X), b \in B \in F(Y), c \in C \in F(Z) \}$. 

The tensor product is also symmetric, with the symmetry relation as in ${\Rel}$.

\medskip

 If $M$ and $N$ are finiteness matrices, $M \ox N$ is given by the outer product,  
 while the finiteness matrices for the associators, and the unitors are given by the 
 characteristic matrices (where the entries are either 0 or 1) of the corresponding 
 coherence relations in ${\sf FRel}$. The commutativity of the rig 
 makes the tensor products symmetric. 

\begin{lemma}
\label{Lemma:tensors_FRel}
$\FRel$ and ${\sf FMat}(R)$ have a symmetric tensor product $\_ \ox \_$ and by de Morgan duality a further symmetric tensor product, called the {\bf par}, $\_ \oa \_$ where $X \oa Y := (X^* \ox Y^*)^*$.
\end{lemma}

The tensor and par do not coincide in general, however, note that, as
\begin{align*} F(X \oa Y) & := F((X^* \ox Y^*)^*) = F(X^* \ox Y^*)^\perp = \{ A' \x B' \mid A' \in F(X)^\perp, B' \in F(Y)^\perp \}^\perp \\
& = \{ G \mid \forall A' \in F(X)^\perp, B' \in F(Y)^\perp. G \cap A' \x B'~\mbox{is finite} \} 
\end{align*}
for $A \in F(X)$ and $B \in F(Y)$ we therefore have $A \x B \in F(X \oa Y)$ as $A \x B \cap A' \x B'$ is finite.  This means the map 
${\sf mx}: A \ox Y \to A \oa Y$, called the mix map,  which is the identity map on the webs is a finiteness relation.  It is, as we shall see,  not generally an isomorphism.  
The unit for the par is $\bot := \top^{*} =  (\{ \star \}, {\cal P}(\{ \star \}),{\cal P}(\{ \star \})) = \top$.  Thus, the units for the two tensors coincide.  This suggests the structure of an isomix category.  Next, we show that there is a linear distributor in both $\FRel$ and ${\sf FMat}(R)$. Preliminary to this we observe:

\begin{lemma}
\label{Lemma: perp perp} 
For any finiteness spaces $X$, $Y$, and $Z$, the following hold:
\begin{enumerate}[(i)]
\item $F(X \ox Y) \subseteq F(X \oa Y)$;
\item $F(X^* \ox Y^*) \subseteq  F(X \ox Y)^{\perp}$;
\item  $F((X^* \ox Y^*) \ox Z) \subseteq F((X \ox Y)^* \ox Z)$.
\end{enumerate}
\end{lemma}
\begin{proof}~
\begin{enumerate}[{\em (i)}]
\item As discussed above.
\item $F(X \ox Y)^{\perp} = F((X \ox Y)^*) = F((X^{**} \ox Y^{**})^*) = F(X^* \oa Y^*) \supseteq F(X^* \ox Y^*)$. 
\item $F((X^* \ox Y^*) \ox Z) = F((X^* \ox Y^*)^{**} \ox Z) = F((X \oa Y)^* \ox Z) \subseteq F((X \ox Y)^* \ox Z)$.
\end{enumerate}
\end{proof}

\begin{proposition}
\label{Lemma: FRel is a LDC}
{\sf FRel}, and $\FMat(R)$ are symmetric LDCs.
\end{proposition}
\begin{proof}
The linear distributor for $\FRel$, $\partial^L : X \ox (Y \oa Z) \to (X \ox Y) \oa Z$ on the webs is the associator in ${\sf Rel}$.  To show this is a finiteness relation (although it is not an isomorphism) it is sufficient to prove that if $A \in F( X \ox (Y \oa Z))$ then ${\sf cod}(A a_\x) \in F( (X \ox Y) \oa Z)$ as $a_\x$ certainly has the preimage of singleton sets in $F( X \ox (Y \oa Z))^\perp$ as they are singleton sets.

Observe that
\begin{align*}
F( X \ox (Y \oa Z))  &= F(X \ox (Y^* \ox Z^*)^*) \\
&= F((X)^{**} \ox  (Y^* \ox Z^*)^*)   \tag{as $X^{**} = X$ } \\
&\subseteq  F(X^* \ox (Y^* \ox Z^*))^\perp  \tag{ Lemma  \ref{Lemma: perp perp}{\em (ii)\/}} \\
&\simeq   F( (X^* \ox Y^*) \ox Z^*)^\perp  \tag{ Associativity of tensor } \\
&\subseteq  F( (X \ox Y)^*  \ox Z^*)^\perp  \tag{ Lemma \ref{Lemma: perp perp}{\em (iii)\/}}\\
&= F( (X \ox Y) \oa Z) \tag{ as $(X, F(X), F(X)^\perp)^* = (X^*, F(X^*), F(X^*)^\perp) := (X, F(X)^\perp, F(X))$ }
\end{align*}

For $\FMat(R)$, the left distributor is the characteristic matrix of $\partial^L$. 
\end{proof}

\begin{proposition}
$\FRel$, and $\FMat(R)$ are symmetric isomix $*$-autonomous categories with 
respect to the involution $(\_)^{*}$.
\end{proposition}
\begin{proof}
By Lemma \ref{Lemma: FRel is a LDC}, $\FRel$ is an LDC. It is an isomix category 
because $\m : \top \to \bot = 1_{\{ \star \}}$ is an  isomorphism.

To show this we must demonstrate that we have a ``cap''  $\top \to A \oa A^*$ and a ``cup'' $A^* \ox A \to \bot$: on the webs these are the 
corresponding maps in ${\sf Rel}$ given by the diagonal subset seen as a relation in two ways.  It suffices to show that these maps (which are dual) 
are finiteness relations.  As the required duality identities are satisfied by the underlying relations, they are also satisfied in ${\sf FRel}$.

Consider $\eta: \top \to A \ox A^*$ where $\eta:= \{ (\star,(a,a)) \mid a \in A\}$  as the finitary sets of $\top$ are just $\{\}$ and $\{ \star \}$ and $\eta^\circ$ certainly has preimages of cofinitary sets cofinitary it remain only to show that $\Delta_A := \{ (a,a) \mid a \in A\}$ is finitary in $A \ox A^{*}$.  However, this is so if and only if $\Delta_A$ is a finiteness relation from $A \to A$ (by Lemma \ref{Lemma: tensor perp characterization} suitably dualized) ... which it certainly is as it is the identity relation!

This argument translates into ${\sf FMat}(R)$ using the characteristic matrices of these relations and the fact that all the coherence maps are given by the corresponding characteristic matrices.
\end{proof}

As $\FRel$ is *-autonomous, it is self-enriched. The internal hom is given as usual by 
$X \lollipop Y =  X^* \ox Y = (X \ox Y^*)^*$. Thus,  $X \lollipop Y = 
(X \times Y, F(X \otimes Y^\perp)^\perp, F(X \otimes Y^\perp))$. 
The same argument applies to $\FMat(R)$.

\begin{lemma}
 Let $X$,  and $Y$ be  finiteness spaces, then $A \in F( X \lollipop Y)$ 
 if and only if  $A: X \to Y$ is a finiteness relation.
\end{lemma}

\subsection{The core of ${\sf FRel}$ and ${\sf FMat}(R)$}
\label{Sec: core of FMat}

Recall that, for any isomix category $\X$,  the {\bf core} of $\X$, ${\sf Core}(\X) 
\subseteq \X$, is the full subcategory consisting of objects $U$ such that for 
any object $X \in \X$, $\mx : U \ox X \to U \oa X$ (and $\mx: X \ox U \to X \oa U$) 
is an isomorphism.  The core is always closed to tensor and par and the units. Thus, 
the core of an isomix category is an isomix category with 
${\sf mx}: U \ox V \to U  \ox V$ is a natural isomorphism. 

In ${\sf FRel}$, the mix map is given by the identity map on the webs and in ${\sf FMat}(R)$ it is given by the characteristic function of this identity map.
Our objective in this subsection is to give a complete description of the core of  ${\sf FRel}$ and ${\sf FMat}(R)$:

\begin{proposition}
\label{Lemma: Core is finite}
$U \in \Core(\FRel)$  if and only if $U$ is finite.  Similarly, 
$U \in \Core(\FMat(R))$ if and only if $U$ is finite.
\end{proposition}

\begin{proof}
Suppose $U$ is a finite set then $U = U^{*}$ as $F(U)=F(U)^\perp = 
{\cal P}(U)$.  We shall prove that $F(U \ox Y) = F(U \ox Y)$ for all $Y$. 
We first prove this for $\FRel$:

\begin{description}
\item[$U$ finite $\Rightarrow$ $U \in \Core(\FRel)$:]
It suffices to show that $R \in  F(U \oa Y)$ then $R \in F(U \ox Y)$.  
However, $R \in F(U \oa Y) = F(U^{*} \oa Y)$ if and only if $R: U \to Y$ is a 
finiteness relation.  As $R$ is a finiteness relation, because 
$U \in F(U)$, ${\sf cod}(R) \in F(Y)$ but then 
$R \subseteq U \x {\sf cod}(R)$ showing $R \in  F(U \ox Y)$. 
\item[$U \in \Core(\FRel)$ $\Rightarrow$  $U$ is finite:] 
We prove the converse statement: if $U$ is infinite set, 
that there is a finiteness space, $Y$, such that $F(X \ox Y) 
\subseteq F(X \oa Y)$. Choose $Y = (Y, \mathcal{P}(Y), \mathcal{P}_f(Y))$ 
such that there is an injective function  $\alpha: U \to Y$: this forces $Y$ to be infinite.  
Because $|U|$ is infinite, there is an infinite $Q \subseteq U$ such that $Q \in F(U)^\perp$ or $Q \in F(U)$.  
Without loss of generality, we may assume that $Q \in F(U)^\perp$.

Define $R := \{ (x,\alpha(x)) \mid x \in Q \}$, then $R \in F(U \ox Y)$ because the following two conditions hold:
\begin{enumerate}[(i)]
\item For all $A \in F(U)^\perp$, $A \tr R \in F(Y)$ as $F(Y) \in \mathcal{P}(Y)$.
\item For all $y \in Y$, $R \tl \{y\} \in F(X)^\perp$, as $\alpha$ is a monic, either $R\tl \{y\} = \varnothing$, 
or $R\tr \{y\} = \{ x \}$ where $\alpha(x) = y$. Since, $R \tr \{y\}$ is a finite, $R \tr \{y\} \in F(X)^\perp$. 
\end{enumerate} 

Next we show, to complete the proof, that $R \notin F(X \ox Y)$. $R \in F(X \ox Y)$ if and only $R \subseteq A \x B$, where $A \in F(U)$ and $B \in F(Y)$.   However, ${\sf dom}(R)= Q \in F(U)^\perp$. Since, $Q$ is infinite, $Q$ cannot also be a member of $F(U)$.  Thus, $Q \notin F(X \ox Y)$.
\end{description}

The same proof applies to $\FMat(R)$ by transposing to the characteristic functions of the finiteness relationships.
\end{proof}

\begin{corollary} ~
\label{Corollary: core}
\begin{enumerate}[(i)]
\item $\Core({\FRel})$ is the category of finite sets and relations.
\item $\Core({\FMat(R)})$ is (equivalent to) the category of finite-dimensional matrices over the rig $R$.
\end{enumerate}
\end{corollary}

%% file: chapter3.tex

\chapter{Dagger linearly distributive categories}
\label{Chap: dagger-LDC}

In this chapter we define dagger linearly distributive categories, dagger linear functors and 
transformations, and provide examples. We also explore the relationship among 
dual, dagger and conjugation functors for LDCs. 

\section{Dagger for LDCs}
\label{Section: dagger LDC}

Conventionally, in categorical quantum mechanics a dagger is defined as a contravariant endofunctor which is 
stationary on objects $(A^\dagger = A)$ and an involution ($f^{\dag \dag} = f$). However, 
in an LDC, the dagger must minimally flip the tensor products to maintain the directionality of the distributor maps. 
Recall that, $\partial^\dagger: (A \ox B) \oa C \to A \ox (B \oa C)$ is not a valid map in an LDC.  Hence, 
for LDCs we cannot expect the dagger to be stationary on objects, 
however, it is still possible for it to be an involution. This section deals with the coherences 
of $\dagger$-functor for LDC and its variants, that is, mix and isomix categories. 

\subsection{Dagger linearly distributive categories}

Before proceeding to define the dagger functor for LDCs, the notion of the opposite 
LDC and the notion of a contravariant linear functors have to be developed.  

If $(\X, \ox, \top, \oa, \bot)$ is a linearly distributive category, the {\bf opposite linear distributive category} 
is $(\X, \ox, \top, \oa, \bot)^{\op} := (\X^{\op}, \oa, \bot, \ox, \top)$ where $\X^{\op}$ is the usual opposite 
category with the monoidal structures flipped as follows: 
\[\ox^{\op} := \oa ~~~~~~~ \top^{\op} := \bot ~~~~~~~ \oa^{\op} := \ox ~~~~~~~ \bot^{\op} := \top\]
$(\_)^\op$ is an endo functor for the category of LDCs and linear functors. It is also an involution:  

$(\X, \ox, \top, \oa, \bot)^{{\op}~{\op}} = (\X, \ox, \top, \oa, \bot) $. 

Let $(F_\ox, F_\oa): (\X, \ox, \top, \oa, \bot)^{\op} \to (\X, \ox, \top, \oa, \bot)$ be a linear functor. The opposite linear functor  $(F_\ox, F_\oa)^{\op}: (\X, \ox, \top, \oa, \bot) \to  (\X, \ox, \top, \oa, \bot)^{\op}$ given by the pair of opposite functors $(F_\oa^{\op}, F_\ox^{\op})$. Observe that $F^{\op}$ is  a mix Frobenius linear functor if and only if $F$ is.

\medskip

\begin{definition}
\label{Definition: daggerLDC succinct}
A {\bf dagger linearly distributive category} ($\dagger$-LDC), is an LDC, $\X$, with a contravariant Frobenius linear functor $(\_)^\dagger: \X^{\op} \to \X$ which is a linear involutive equivalence   $(\_)^\dagger ~\dashvv ~ (\_)^{\dagger^{\op}}: \X^{\op} \to \X$.
\end{definition}

We unfold this definition in Proposition \ref{Definition: daggerLDC elaborate}. However, first note that saying the dagger is an {\bf involutive} equivalence asserts that the unit and 
counit of the equivalence are the same (although one is in the opposite category).  
Thus, the adjunction expands to take the form $(\imath,\imath): (\_)^\dagger ~\dashvv ~ (\_)^{\dagger^{\op}}: \X^{\op} \to \X$.  
However, the unit and counit are linear natural transformations so $\imath$ expands to $\imath = (\imath_\ox,\imath_\oa)$.  
As the dagger functor is a left adjoint, it is strong and, thus, is normal.  
Furthermore, as the unit of an equivalence, $\imath$ is a linear natural isomorphism,   
this means  $\imath = (\imath_\ox,\imath_\oa)$ satisfies the requirements of 
Lemma \ref{Lemma: Frobenius linear transformation}, implying that $\imath_\ox^{-1}  = \imath_\oa$.  
Simplifying notation we shall set $\iota:= \imath_\oa$  so the unit linear transformation is 
$\imath := (\iota^{-1},\iota)$. We then can simplify the requirements of $\imath$ to the map  
$\iota: A \to (A^\dagger)^\dagger$ which we refer to as the {\bf involutor}.

A {\bf symmetric  $\dagger$-LDC} is a $\dagger$-LDC which is a symmetric LDC 
for which the dagger is a symmetric linear functor. A {\bf cyclic $\dagger$-$*$-autonomous category} is a 
$\dagger$-LDC with chosen left and right duals, and a cyclor which is preserved 
by the dagger.  A $\dagger$-{\bf mix} category is a $\dagger$-LDC for which $(\_)^\dagger: \X^{\op} \to \X$ 
is a mix functor.  As the dagger functor is strong (and so normal) if the category is an isomix category then being {\bf $\dagger$-mix} already implies that 
the dagger is an isomix functor.  Thus, a {\bf $\dagger$-\bf isomix} category is a $\dagger$-mix category which happens to be an isomix category.

In the remainder of the section, we unfold the definition of a $\dagger$-isomix category and give the coherence requirements explicitly.

\begin{proposition}
\label{Definition: daggerLDC elaborate}
A dagger linearly distributive category is an LDC with a functor $(\_)^\dag:\X^\op\to \X$ and 
natural isomorphisms 
\begin{align*}
\text{ \bf laxors: }  A^\dag \ox B^\dag &\xrightarrow{ \lambda_\ox} (A\oa B)^\dag ~~~~~ A^\dag \oa B^\dag \xrightarrow{ \lambda_\oa} (A\ox B)^\dag \\
\top &\xrightarrow{\lambda_\top} \bot^\dag ~~~~~~~~~~~~~~~~~~~~~ \bot \xrightarrow{\lambda_\bot} \top^\dag \\
\text{ \bf involutor: }  A &\xrightarrow{\iota} (A^\dag)^\dag 
\end{align*}
such that the following coherences hold:
\begin{enumerate}[{\bf [$\dagger$-ldc.1]}]
\item Interaction of $\lambda_\ox, \lambda_\oa$  with associators:
\[
\begin{tabular}{cc}
\xymatrix{
A^\dag \ox (B^\dag \ox C^\dag)                \ar@{->}[r]^{{a_\ox}^{-1}}    \ar@{->}[d]_{1 \ox  \lambda_\ox}   
  & (A^\dag \ox B^\dag) \ox C^\dag              \ar@{->}[d]^{\lambda_\ox \ox 1}  \\
A^\dag \ox ( B \oa C)^\dag                     \ar@{->}[d]_{\lambda_\ox}  
  & (A \oa B)^\dag \ox C^\dag                   \ar@{->}[d]^{\lambda_\ox}    \\
(A \oa (B \oa C))^\dag                          \ar@{->}[r]_{a_\oa^\dag}
  & ( (A\oa B) \oa C)^\dag
} &\xymatrix{
A^\dag \oa (B^\dag \oa C^\dag)                \ar@{->}[r]^{a_\oa^{-1}}    \ar@{->}[d]_{1 \oa \lambda_\oa }  
  & (A^\dag \oa B^\dag) \oa C^\dag              \ar@{->}[d]^{\lambda_\oa \oa 1}  \\
   A^\dag \oa (B \ox C)^\dag                 \ar@{->}[d]_{\lambda_\oa}  
  & (A \ox B)^\dag \oa C^\dag                    \ar@{->}[d]^{\lambda_\oa}    \\
(A \ox (B\ox C))^\dag                          \ar@{->}[r]_{a_\ox^\dag}
  & ((A\ox B) \ox C)^\dag
}
\end{tabular}
\]

\item Interaction of $\lambda_\top, \lambda_\bot$ with unitors: 
\begin{center}
\begin{tabular}{cc}
\xymatrix{
\top \ox A^\dag                            \ar@{->}[rr]^{\lambda_\top\ox 1} \ar@{->}[d]_{u_\ox^R}  
&  & \bot^\dag\ox A^\dag           \ar@{->}[d]^{\lambda_\ox}\\
A^\dag                                     
&  & (\bot \oa A)^\dag                      \ar@{<-}[ll]^{(u_\oa^R)^\dag}\\
} &
\xymatrix{
\bot \oa A^\dag                            \ar@{->}[rr]^{\lambda_\bot\oa 1} \ar@{->}[d]_{u_\oa^R} 
& & \top^\dag\oa A^\dag                    \ar@{->}[d]^{\lambda_\oa}\\
A^\dag                                     
&  & (\top \ox A)^\dag                      \ar@{<-}[ll]^{(u_\ox^R)^\dag}\\
} 
\end{tabular}
\end{center}
and two symmetric diagrams for $u_\ox^L$ and $u_\oa^L$ must also be satisfied.

\item Interaction of $\lambda_\ox, \lambda_\oa$ with linear distributors:
\begin{center}
\begin{tabular}{cc}
\xymatrix{
A^\dag \ox(B^\dag\oa C^\dag)              \ar@{->}[r]^{\partial^L} \ar@{->}[d]_{1\ox\lambda_\oa}
  & (A^\dag \ox B^\dag)\oa C^\dag         \ar@{->}[d]_{\lambda_\ox\oa 1}\\
A^\dag \ox (B\ox C)^\dag                  \ar@{->}[d]_{\lambda_\ox}
  & (A\oa B)^\dag \oa C^\dag              \ar@{->}[d]^{\lambda_\oa}\\
(A\oa (B\ox C))^\dag                      \ar@{->}[r]_{(\partial^R)^\dag}
  & ((A\oa B)\ox C)^\dag
} &
\xymatrix{
(A^\dag \oa B^\dag) \ox C^\dag     \ar@{->}[r]^{\partial^R} \ar@{->}[d]_{\lambda_\oa\ox 1} 
  & A^\dag \oa (B^\dag \ox C^\dag) \ar@{->}[d]_{1\oa \lambda_\ox}\\
(A\ox B)^\dag \ox C^\dag           \ar@{->}[d]_{\lambda_\ox}
  & A^\dag \oa (B\oa C)^\dag       \ar@{->}[d]^{\lambda_\oa}\\
((A\ox B)\oa C)^\dag               \ar@{->}[r]_{(\partial^L)^\dag}
  & (A\ox (B\oa C))^\dag
} 
\end{tabular}
\end{center}

\item Interaction of $\iota: A \rightarrow A^{\dagger\dagger}$ with $\lambda_\ox$, $\lambda_\oa$:
\begin{center}
\begin{tabular}{cc}
\xymatrix{
A\oa B                          \ar@{->}[r]^{\iota} \ar@{->}[d]_{\iota \oa \iota} 
  & ((A\oa B)^\dag)^\dag        \ar@{->}[d]^{\lambda_\ox^\dag}\\
(A^\dag)^\dag   \oa (B^\dag)^\dag                \ar@{->}[r]_{\lambda_\oa}
  & (A^\dag\ox B^\dag)^\dag
} &  \xymatrix{
A\ox B                          \ar@{->}[r]^{\iota} \ar@{->}[d]_{\iota \ox \iota} 
  & ((A\ox B)^\dag)^\dag        \ar@{->}[d]^{\lambda_\oa^\dag}\\ 
(A^\dag)^\dag  \ox (B^\dag)^\dag                \ar@{->}[r]_{\lambda_\ox}
  & (A^\dag\oa B^\dag)^\dag
}
\end{tabular}
\end{center}
\item Interaction of $\iota: A \rightarrow A^{\dagger\dagger}$ with $\lambda_\top$, $\lambda_\bot$:
\begin{center}
\begin{tabular}{cc}
$\begin{matrix} \xymatrix{
&\bot                   \ar@{->}[r]^{\iota} \ar@{->}[dr]_{\lambda_\bot}  
  & (\bot^\dag)^\dag   \ar@{->}[d]^{\lambda_\top^\dag}\\
&{}
  & \top^\dag 
} \end{matrix}$ & 
 $\begin{matrix} \xymatrix{
&\top                   \ar@{->}[r]^{\iota} \ar@{->}[dr]_{\lambda_\top} 
  & (\top^\dag)^\dag   \ar@{->}[d]^{\lambda_\bot^\dag} \\
&{}
  & \bot^\dag 
} \end{matrix}$
\end{tabular}
\end{center}
\item $\iota_{A^\dagger} = (\iota_A^{-1})^\dagger: A^\dagger \to A^{\dagger\dagger\dagger}$  
\end{enumerate}
\end{proposition}

The dagger structure is obtained from the previous proposition using strong monoidal laxors:  to form a linear functor the laxor $\lambda_\oa$ 
needs to be reversed by taking its inverse. Then, we have $\nu_\ox^l = \nu_\ox^r := \lambda_\oa^{-1}$ and $\nu_\oa^l = \nu_\oa^r := \lambda_\ox$.  
Once this adjustment is made all the required coherences for $\dagger$ to be a linear functor are present.
Note that {\bf [$\dagger$-ldc.6]} equivalently expresses the triangle identities of the 
adjunction $(\iota, \iota):: \dagger^{\op} \dashv \dagger : \X^{\op} \to \X$.   
The coherences for the involutor asserts that it is a monoidal transformation for both the tensor and par: 
by Lemma \ref{simplifying-coherences} (ii) this suffices to show that it is a linear transformation.

A {\bf symmetric $\dagger$-LDC} is a $\dagger$-LDC which is a symmetric LDC and for which the following additional diagrams commute:
\begin{enumerate}[{\bf [$\dagger$-ldc.7]}]
\item Interaction of $\lambda_\ox , \lambda_\oa$ with symmetry maps: 
\[
\begin{tabular}{cc}
\xymatrix{
A^\dag \ox B^\dag                          \ar@{->}[r]^{\lambda_\ox} \ar@{->}[d]_{c_\ox}         
  & (A\oa B)^\dag                          \ar@{->}[d]^{c_\oa^\dag}\\
B^\dag \ox A^\dag                          \ar@{->}[r]_{\lambda_\ox}
  & (B\oa A)^\dag\\
} & \xymatrix{
A^\dag \oa  B^\dag                          \ar@{->}[r]^{\lambda_\oa} \ar@{->}[d]_{c_\oa}  
  & (A\ox B)^\dag                          \ar@{->}[d]^{c_\ox^\dag}\\
B^\dag \oa A^\dag                          \ar@{->}[r]_{\lambda_\oa}
  & (B\ox A)^\dag\\
}
\end{tabular}
\]
\end{enumerate}

A {\bf $\dagger$-mix category} is a $\dagger$-LDC which has a mix map and satisfies the following additional coherence:

\[ \mbox{\bf [$\dagger$-\text{mix}]}  ~~~~\begin{array}[c]{c} 
\xymatrix{
\bot                 \ar@{->}[r]^{{\sf m}} \ar@{->}[d]_{\lambda_\bot}  
  & \top             \ar@{->}[d]^{\lambda_\top}\\
\top^\dag            \ar@{->}[r]_{{\sf m}^\dag}
  & \bot^\dag
} 
\end{array} \]
If ${\sf m}$ is an isomorphism, then $\X$ is a {\bf $\dagger$-isomix category} and, 
since $(\_)^\dagger$ is normal, $(\_)^\dagger$ is an isomix Frobenius functor.

\begin{lemma}
\label{lemma: mixdagger}
Suppose $\X$ is a $\dagger$-mix category then the following diagram commutes:
\[
\xymatrix{ A^\dag \ox B^\dag  \ar[r]^{\mx} \ar[d]_{\lambda_\ox}&  A^\dag \oa B^\dag \ar[d]^{\lambda_\oa} \\
                (A \oa B)^\dag \ar[r]_{\mx^\dag}  & (A \ox B)^\dag }
\]
\end{lemma}
\begin{proof}
The proof follows directly from Lemma \ref{Lemma: Mix Frobenius linear functor}.
\end{proof}

 With respect to its applications to quantum theory, this thesis primarily focuses on $\dagger$-isomix categories. 
 As we will see in Section \ref{Sec: unitary}, the notion of unitary objects and unitary isomorphisms 
 is supported only within a $\dagger$-isomix category.

It is useful to observe that the core of a mix category is closed under taking the dagger and duals.

\begin{lemma}
\label{Lemma: mixdagger}
Suppose $\X$ is a $\dagger$-mix category and $A \in \Core(\X)$ then $A^\dagger \in$ $\Core(\X)$.
\end{lemma}
\begin{proof}
The natural transformation $A^\dagger \ox X \xrightarrow{\mx} A^\dagger \oa X$ is an isomorphism as follows:
\[
\xymatrix{
A^\dagger \ox X \ar[r]^{1 \ox \iota} \ar[d]_{\mx} \ar@{}[dr]|{\scalebox{0.95}{\tiny\bf (nat. {\sf mx})}}
& A^\dagger \ox X^{\dagger\dagger} \ar[r]^{\lambda_\ox} \ar[d]_{\mx} \ar@{}[dr]|{\scalebox{.95}{\tiny\bf {lem. \ref{lemma: mixdagger}}}}
& (A \oa X^\dagger)^\dagger \ar[d]^{\mx^\dagger} \\
A^\dagger \oa X \ar[r]_{1 \oa \iota} 
& A^\dagger \oa X^{\dagger \dagger} \ar[r]_{\lambda_\oa}
&(A \ox X^\dagger)^\dagger
}
\]

\end{proof}

\begin{lemma}
Let $\X$ be $\dagger$-LDC. If $A \dashvv B$ then $B^\dagger \dashvv A^\dagger$.
\end{lemma}
\begin{proof}
The statement  follows from Lemma \ref{Lemma: linear adjoints}: Frobenius functors preserve linear adjoints. 
Explicitly, if $(\eta,\epsilon): A \dashvv B$ then $(\lambda_\top\epsilon^\dag\lambda_\oa^{-1},\lambda_\ox\eta^\dagger \lambda_\bot^{-1}): B^\dagger \dashvv A^\dagger$. 
\end{proof}

Suppose $\X$ is a $\dagger$-$*$-autonomous category and $(\eta*, \epsilon*): A^* \dashvv A$, then $((\epsilon*)^\dagger, (\eta*)^\dagger): A^\dagger \dashvv (A^*)^\dagger$, where $((\epsilon*)^\dagger, (\eta*)^\dagger) :=   (\lambda_\top\epsilon*^\dag\lambda_\oa^{-1},\lambda_\ox\eta*^\dagger \lambda_\bot^{-1})$. We draw $(\epsilon*)^\dagger$ and $(*\epsilon)^\dagger$ as dagger cups, and $(\eta*)^\dagger$ and $(*\eta)^\dagger$ as dagger caps which are pictorially represented as follows:
 \[
 \begin{tikzpicture}
	\begin{pgfonlayer}{nodelayer}
		\node [style=none] (0) at (-1, 2) {};
		\node [style=none] (1) at (1, 2) {};
		\node [style=none] (2) at (-1, 0.5) {};
		\node [style=none] (3) at (1, 0.5) {};
		\node [style=none] (4) at (-1.5, 1.5) {$X^\dagger$};
		\node [style=none] (5) at (1.5, 1.5) {$(^*X)^\dagger$};
		\node [style=none] (6) at (0, -1) {$(*\eta)^\dagger$};
	\end{pgfonlayer}
	\begin{pgfonlayer}{edgelayer}
		\draw (2.center) to (0.center);
		\draw [bend right=90, looseness=2.00] (2.center) to (3.center);
		\draw (3.center) to (1.center);
	\end{pgfonlayer}
\end{tikzpicture} ~~~~~~~~~ \begin{tikzpicture}
	\begin{pgfonlayer}{nodelayer}
		\node [style=none] (0) at (-1, 2) {};
		\node [style=none] (1) at (1, 2) {};
		\node [style=none] (2) at (-1, 0.5) {};
		\node [style=none] (3) at (1, 0.5) {};
		\node [style=none] (4) at (-1.5, 1.5) {$X^{*\dagger}$};
		\node [style=none] (5) at (1.5, 1.5) {$X^\dagger$};
		\node [style=none] (6) at (0, -1) {$(\eta*)^\dagger$};
	\end{pgfonlayer}
	\begin{pgfonlayer}{edgelayer}
		\draw (2.center) to (0.center);
		\draw [bend right=90, looseness=2.00] (2.center) to (3.center);
		\draw (3.center) to (1.center);
	\end{pgfonlayer}
\end{tikzpicture} ~~~~~~~~~ \begin{tikzpicture}
	\begin{pgfonlayer}{nodelayer}
		\node [style=none] (0) at (-1, -1) {};
		\node [style=none] (1) at (1, -1) {};
		\node [style=none] (2) at (-1, 0.5) {};
		\node [style=none] (3) at (1, 0.5) {};
		\node [style=none] (4) at (-1.5, -0.5) {$X^\dagger$};
		\node [style=none] (5) at (1.5, -0.5) {$(^*X)^\dagger$};
		\node [style=none] (6) at (0, 2) {$(*\epsilon)^\dagger$ };
	\end{pgfonlayer}
	\begin{pgfonlayer}{edgelayer}
		\draw (2.center) to (0.center);
		\draw [bend left=90, looseness=2.00] (2.center) to (3.center);
		\draw (3.center) to (1.center);
	\end{pgfonlayer}
\end{tikzpicture} ~~~~~~~~~ \begin{tikzpicture}
	\begin{pgfonlayer}{nodelayer}
		\node [style=none] (0) at (-1, -1) {};
		\node [style=none] (1) at (1, -1) {};
		\node [style=none] (2) at (-1, 0.5) {};
		\node [style=none] (3) at (1, 0.5) {};
		\node [style=none] (4) at (-1.5, -0.5) {$X^{*^\dagger}$};
		\node [style=none] (5) at (1.5, -0.5) {$X^\dagger$};
		\node [style=none] (6) at (0, 2) {$(\epsilon*)^\dagger$};
	\end{pgfonlayer}
	\begin{pgfonlayer}{edgelayer}
		\draw (2.center) to (0.center);
		\draw [bend left=90, looseness=2.00] (2.center) to (3.center);
		\draw (3.center) to (1.center);
	\end{pgfonlayer}
\end{tikzpicture} 
 \]

A $\dagger$-$*$-autonomous category is a {\bf cyclic} $\dagger$-$*$-autonomous category when the dagger preserves the cyclor in the following sense. 
\[ 
\begin{tikzpicture}
	\begin{pgfonlayer}{nodelayer}
		\node [style=none] (0) at (5.25, 2) {};
		\node [style=none] (1) at (5.25, -1) {};
		\node [style=none] (2) at (4, -0.9999999) {};
		\node [style=none] (3) at (4, 0.9999999) {};
		\node [style=none] (4) at (2.75, 0.9999999) {};
		\node [style=circle] (5) at (2.75, -0.5) {$\psi^\dagger$};
		\node [style=circle, scale=2.5] (6) at (2.75, -1.75) {};
		\node [style=none] (13) at (2.75, -1.75) {$\psi^{-1 \dagger}$};
		\node [style=none] (7) at (2.75, -2.75) {};
		\node [style=none] (8) at (5.75, 1.5) {$(A^\dagger)^*$};
		\node [style=none] (9) at (4.75, -1.75) {$\epsilon*$};
		\node [style=none] (10) at (3.25, 1.75) {$(*\epsilon)^\dagger$};
		\node [style=none] (11) at (2.25, -2.5) {$(^*A)^\dagger$};
		\node [style=none] (12) at (4.25, -0) {$A^\dagger$};
	\end{pgfonlayer}
	\begin{pgfonlayer}{edgelayer}
		\draw (0.center) to (1.center);
		\draw [bend left=90, looseness=1.25] (1.center) to (2.center);
		\draw (2.center) to (3.center);
		\draw [bend right=90, looseness=1.25] (3.center) to (4.center);
		\draw (4.center) to (5);
		\draw (5) to (6);
		\draw (6) to (7.center);
	\end{pgfonlayer}
\end{tikzpicture} = 
\begin{tikzpicture}
	\begin{pgfonlayer}{nodelayer}
		\node [style=none] (0) at (5.5, -0.9999999) {};
		\node [style=none] (1) at (4, -1) {};
		\node [style=none] (2) at (4, -3) {};
		\node [style=none] (3) at (2.5, -3) {};
		\node [style=circle, scale=1.3] (4) at (2.5, -0.9999999) {$\psi$};
		\node [style=circle, scale=2.5] (5) at (5.5, -4) {};
		\node [style=none] (13) at (5.5, -4) {$\psi^{-1 \dagger}$};
		\node [style=none] (6) at (2.5, 0.9999999) {};
		\node [style=none] (7) at (1.25, -0.25) {$(A^\dagger)^*$};
		\node [style=none] (8) at (4.75, -0.25) {$(\epsilon*)^\dagger$};
		\node [style=none] (9) at (3.25, -3.75) {$*\epsilon$};
		\node [style=none] (10) at (6.999999, -5) {$(^*A)^\dagger$};
		\node [style=none] (11) at (4.5, -2) {$A^\dagger$};
		\node [style=none] (12) at (5.5, -5.25) {};
	\end{pgfonlayer}
	\begin{pgfonlayer}{edgelayer}
		\draw [bend right=90, looseness=1.25] (0.center) to (1.center);
		\draw (1.center) to (2.center);
		\draw [bend left=90, looseness=1.25] (2.center) to (3.center);
		\draw (3.center) to (4);
		\draw (0.center) to (5);
		\draw (5) to (12.center);
		\draw (4) to (6.center);
	\end{pgfonlayer}
\end{tikzpicture} \]

\begin{lemma} 
\label{Lemma: cyclic dagger}
In a cyclic, $\dagger$-$*$-autonomous category,
\[\begin{tikzpicture}
	\begin{pgfonlayer}{nodelayer}
		\node [style=none] (0) at (6, 1) {};
		\node [style=none] (1) at (4, 1) {};
		\node [style=none] (2) at (4, -1.25) {};
		\node [style=none] (3) at (2, -1.25) {};
		\node [style=circle, scale=1.25] (4) at (2, -0.5) {$\psi$};
		\node [style=none] (5) at (2, 2.25) {};
		\node [style=none] (6) at (5, 2) {$(*\epsilon)^\dagger$};
		\node [style=none] (7) at (3, -2.25) {$*\epsilon$};
		\node [style=none] (8) at (6, -2.5) {};
		\node [style=circle] (9) at (4, -0.25) {$\psi^\dagger$};
		\node [style=none] (10) at (1.25, 2) {$A^{*\dagger*}$};
		\node [style=none] (11) at (6.5, -2.25) {$A^\dagger$};
		\node [style=none] (12) at (1.25, -1.25) {$~^*(A^{\dagger*})$};
		\node [style=none] (13) at (4.5, -1) {$A^{\dagger*}$};
		\node [style=none] (14) at (4.5, 0.5) {$~^*(A^\dagger)$};
	\end{pgfonlayer}
	\begin{pgfonlayer}{edgelayer}
		\draw [bend right=90, looseness=1.25] (0.center) to (1.center);
		\draw [bend left=90, looseness=1.25] (2.center) to (3.center);
		\draw (3.center) to (4);
		\draw (4) to (5.center);
		\draw (0.center) to (8.center);
		\draw (1.center) to (9);
		\draw (9) to (2.center);
	\end{pgfonlayer}
\end{tikzpicture} =\begin{tikzpicture}
	\begin{pgfonlayer}{nodelayer}
		\node [style=none] (15) at (13, -0.75) {};
		\node [style=none] (16) at (11, -0.75) {};
		\node [style=none] (17) at (11, 0.25) {};
		\node [style=none] (18) at (9, 0.25) {};
		\node [style=none] (19) at (10, 1.25) {$(\epsilon*)^\dagger$};
		\node [style=none] (20) at (12, -1.75) {$\epsilon*$};
		\node [style=none] (21) at (13, 2) {};
		\node [style=none] (22) at (13.75, 1.5) {$A^{*\dagger*}$};
		\node [style=none] (23) at (8.25, -1) {$A^\dagger$};
		\node [style=none] (24) at (9, -2.5) {};
		\node [style=none] (25) at (11.5, -0.5) {$A^{*\dagger}$};
	\end{pgfonlayer}
	\begin{pgfonlayer}{edgelayer}
		\draw [bend left=90, looseness=1.25] (15.center) to (16.center);
		\draw [bend right=90, looseness=1.25] (17.center) to (18.center);
		\draw (15.center) to (21.center);
		\draw (17.center) to (16.center);
		\draw (18.center) to (24.center);
	\end{pgfonlayer}
\end{tikzpicture} \]
\end{lemma}
\begin{proof}
Proved by direct application of Lemma \ref{Lemma: cyclic Frob}.
\end{proof}

\subsection{Sequent calculus for $\dagger$-linear logic}
\label{Sec: dagger sequent rules}

$\dagger$-LDCs provide a categorical semantics for the proof theory for multiplicative linear logic with the dagger ($\dagger$-MLL).
Along with the sequent rules of MLL, $\dagger$-MLL includes the additional rules $(\dagger)$, and ($\iota$), as shown in 
Figure \ref{Fig: dagger sequent}.
\begin{figure}[h]
	\centering
	\AxiomC{$\Gamma \vdash \Delta$}
	\LeftLabel{($\dagger$)}
	\UnaryInfC{$\Delta^\dag \vdash \Gamma^\dag$}
	\DisplayProof
	\hspace{1.5em}
	\AxiomC{$\Gamma_1,A, \Gamma_2 \vdash \Delta$}
	\LeftLabel{($\iota$L)}
	\UnaryInfC{$\Gamma_1, A^{\dag \dag}, \Gamma_2 \vdash \Delta$}
	\DisplayProof 
	\hspace{1.5em}
	\AxiomC{$\Gamma \vdash \Delta_1, A, \Delta_2$}
	\LeftLabel{($\iota$R)}
	\UnaryInfC{$\Gamma \vdash \Delta_1, A^{\dag \dag}, \Delta_2$}
	\DisplayProof	

  	\vspace{1em}

	where, if $\Gamma = A_1, A_2, A_3, \cdots, A_n$, then $\Gamma^\dagger = A_1^\dag, A_2^\dag, A_3^\dag, \cdots, A_n^\dag$
	\caption{Sequent rules for $\dagger$}
	\label{Fig: dagger sequent}
\end{figure}

The rule $(\dagger)$ corresponds to the contravariance of the $\dagger$-functor, and the rules ($\iota$L) and ($\iota$R) correspond 
to the involutor natural isomorphism, $\iota: A \to A^{\dagger \dagger}$. The derivation of $(\iota L)^{(-1)}$ and $(\iota R)^{(-1)}$ are 
shown in Figure \ref{Fig: iota inv}.

\begin{figure}[h]
    \centering
	\AxiomC{$A \vdash A$}
	\RightLabel{$\iota R$}
	\UnaryInfC{$A \vdash A^{\dag \dag}$} 
    \DisplayProof
	\hspace{1.5em}
	\AxiomC{$A \vdash A$}
	\RightLabel{$\iota L$}
	\UnaryInfC{$A^{\dag \dag} \vdash A$} 
    \DisplayProof
	\caption{Sequent proof:- $\iota$ is an isomorphism}
	\label{Fig: iota inv}
\end{figure}

One can derive the unit and the tensor laxors using the sequent rules, ($\dagger$), ($\iota$L), and ($\iota$R).
For example, Figure \ref{Fig: tensor laxor rule} shows the derivation of the $\oa$-laxor, and Figure \ref{Fig: unit laxor rule} 
shows the derivation of the $\bot$-laxor. The steps labelled $(*)$ in Figure \ref{Fig: tensor laxor rule} involve a cut. This 
leads to the question of whether to system satisfies cut elimination: with the current presentation, the system 
not appear to satisfy cut elimination. 

\begin{figure}[h]
    \centering
	\AxiomC{$ $}
	\RightLabel{id}
	\UnaryInfC{$A \vdash A$}
	\AxiomC{$ $}
	\RightLabel{id}
	\UnaryInfC{$B \vdash B$}
	\RightLabel{$\ox$R}
	\BinaryInfC{$A, B \vdash A \ox B$}
	\UnaryInfC{$(A \ox B)^\dag \vdash A^\dag, B^\dag$}
	\RightLabel{$\oa$R}
    \UnaryInfC{$(A \ox B)^\dag \vdash A^\dag \oa B^\dag$}
    \DisplayProof
	\hspace{1.5em}
	\AxiomC{$ $}
	\RightLabel{id}
	\UnaryInfC{$A^\dag \vdash A^\dag$}
	\AxiomC{$ $}
	\RightLabel{id}
	\UnaryInfC{$B^\dag \vdash B^\dag$}
	\RightLabel{$\oa$L}
	\BinaryInfC{$A^\dag \oa B^\dag \vdash A^\dag, B^\dag$}
	\RightLabel{$\dagger$}
	\UnaryInfC{$A^{\dag \dag}, B^{\dag \dag} \vdash (A^\dag \oa B^\dag)^\dag$}
	\RightLabel{$(*)~ 2 \times \iota L ^{(-1)}$}
	\UnaryInfC{$A, B \vdash (A^\dag \oa B^\dag)^\dag$}
	\RightLabel{$\ox$L}
	\UnaryInfC{$A \ox B \vdash (A^\dag \oa B^\dag)^\dag$}
	\RightLabel{$\iota$L}
	\UnaryInfC{$ (A^\dag \oa B^\dag)^{\dag \dag} \vdash (A \ox B)^\dag$}
	\RightLabel{$(*)~ \iota L ^{(-1)}$}
	\UnaryInfC{$(A^\dag \oa B^\dag) \vdash (A \ox B)^\dag$}
	\DisplayProof
	\caption{Derivation of the $\oa$-laxor}
	\label{Fig: tensor laxor rule}
\end{figure}

\begin{figure}
	\centering
    \AxiomC{$ $ } 
	\RightLabel{id}
	\UnaryInfC{$ \vdash \top$} 
    \RightLabel{$\dagger$}
    \UnaryInfC{$\top^\dag \vdash $}
	\RightLabel{$\bot$R}
	\UnaryInfC{$\top^{\dag} \vdash \bot$}
	\DisplayProof
	\hspace{1.5em}    
		\AxiomC{$ $}
		\RightLabel{id}
		\UnaryInfC{$\bot \vdash \bot$} 
		\RightLabel{$\iota$R}
		\UnaryInfC{$\bot \vdash \bot^{\dag \dag}$}
						\AxiomC{$ $ } 
						\RightLabel{$\bot$L}
						\UnaryInfC{$ \bot \vdash $} 
						\RightLabel{$\dagger$}
						\UnaryInfC{$ \vdash \bot^\dagger $} 
						\RightLabel{$\top$L}
						\UnaryInfC{$ \top \vdash \bot^\dagger $} 
						\RightLabel{$\dagger$}
						\UnaryInfC{$ \bot^{\dagger \dagger} \vdash \top^\dagger $} 
	\RightLabel{Cut}  
	\BinaryInfC{$\bot \vdash \top^\dagger$}
	\DisplayProof
	\caption{Derivation of $\bot$-laxor} 
    \label{Fig: unit laxor rule}
\end{figure}

\FloatBarrier

\subsection{Dagger functor box}

Suppose $\X$ is a $\dagger$-LDC and  $f: A \rightarrow B \in \X$. Then, the map $f^\dagger: B^\dagger \rightarrow A^\dagger$ is graphically depicted as follows:
\[ \begin{tikzpicture}
	\begin{pgfonlayer}{nodelayer}
		\node [style=circle, scale=1.5] (0) at (0, 1) {};
		\node [style=none] (15) at (0, 1) {$f$};
		\node [style=none] (1) at (0.75, 2) {};
		\node [style=none] (2) at (-0.75, 0) {};
		\node [style=none] (3) at (-1, 2) {};
		\node [style=none] (4) at (-1, 0) {};
		\node [style=none] (5) at (1, 0) {};
		\node [style=none] (6) at (1, 2) {};
		\node [style=none] (7) at (0, 2) {};
		\node [style=none] (8) at (0, 2.75) {};
		\node [style=none] (9) at (0, 0) {};
		\node [style=none] (10) at (0, -0.75) {};
		\node [style=none] (11) at (0.5, 1.5) {$A$};
		\node [style=none] (12) at (-0.5, 0.5) {$B$};
		\node [style=none] (13) at (0.5, -0.75) {$A^\dagger$};
		\node [style=none] (14) at (0.5, 2.75) {$B^\dagger$};
	\end{pgfonlayer}
	\begin{pgfonlayer}{edgelayer}
		\draw [style=none, in=-90, out=90, looseness=1.00] (0) to (1.center);
		\draw [style=none, in=90, out=-90, looseness=1.00] (0) to (2.center);
		\draw [style=none] (3.center) to (4.center);
		\draw [style=none] (4.center) to (5.center);
		\draw [style=none] (5.center) to (6.center);
		\draw [style=none] (6.center) to (3.center);
		\draw [style=none] (8.center) to (7.center);
		\draw [style=none] (9.center) to (10.center);
	\end{pgfonlayer}
\end{tikzpicture} \]
The rectangle is a functor box for the $\dag$-functor.  Notice how we use vertical mirroring to express the contravariance of the $\dag$-functor. By the functoriality of $(\_)^\dagger$, we have:
\begin{tikzpicture}
	\begin{pgfonlayer}{nodelayer}
		\node [style=none] (0) at (-2, 1.25) {};
		\node [style=none] (1) at (-1.5, 1.25) {};
		\node [style=none] (2) at (-2, 0.7499999) {};
		\node [style=none] (3) at (-1.5, 0.7499999) {};
		\node [style=none] (4) at (-1.75, 1.5) {};
		\node [style=none] (5) at (-1.75, 0.5000001) {};
	\end{pgfonlayer}
	\begin{pgfonlayer}{edgelayer}
		\draw (0) to (1);
		\draw (3) to (1);
		\draw (0) to (2);
		\draw (2) to (3);
		\draw (5) to (4);
	\end{pgfonlayer}
\end{tikzpicture}
=
\begin{tikzpicture}
	\begin{pgfonlayer}{nodelayer}
		\node [style=none] (0) at (-1.75, 1.5) {};
		\node [style=none] (1) at (-1.75, 0.5) {};
	\end{pgfonlayer}
	\begin{pgfonlayer}{edgelayer}
		\draw (1) to (0);
	\end{pgfonlayer}
\end{tikzpicture}.

These contravariant functor boxes compose contravariantly. Given maps $f: A \to B$ and $g: B \to C$:
\[
\begin{tikzpicture} 
	\begin{pgfonlayer}{nodelayer}
		\node [style=none] (0) at (-1, 2) {};
		\node [style=none] (1) at (0.5, 2) {};
		\node [style=none] (2) at (-1, 0.5) {};
		\node [style=none] (3) at (0.5, 0.5) {};
		\node [style=none] (4) at (0.5, -1.5) {};
		\node [style=none] (5) at (0.5, -0) {};
		\node [style=none] (6) at (-1, -1.5) {};
		\node [style=none] (7) at (-1, -0) {};
		\node [style=circle, scale=1.5] (8) at (-0.25, 1.25) {};
		\node [style=none] (9) at (-0.5, 2) {};
		\node [style=none] (10) at (0, 0.5) {};
		\node [style=none] (11) at (-0.5, -0) {};
		\node [style=circle, scale=1.5] (12) at (-0.25, -0.75) {};
		\node [style=none] (13) at (0, -1.5) {};
		\node [style=none] (14) at (-0.25, 2.75) {};
		\node [style=none] (15) at (-0.25, -2.25) {};
		\node [style=none] (16) at (-0.25, 0.5) {};
		\node [style=none] (17) at (-0.25, -0) {};
		\node [style=none] (18) at (-0.25, 2) {};
		\node [style=none] (19) at (-0.25, -1.5) {};
		\node [style=none] (20) at (-0.25, -0.75) {$f$};
		\node [style=none] (21) at (-0.25, 1.25) {$g$};
		\node [style=none] (22) at (0, 2.5) {$C^\dagger$};
		\node [style=none] (23) at (0, 0.25) {$B^\dagger$};
		\node [style=none] (24) at (0, -2) {$A^\dagger$};
	\end{pgfonlayer}
	\begin{pgfonlayer}{edgelayer}
		\draw (0.center) to (1.center);
		\draw (1.center) to (3.center);
		\draw (3.center) to (2.center);
		\draw (2.center) to (0.center);
		\draw (7.center) to (5.center);
		\draw (5.center) to (4.center);
		\draw (4.center) to (6.center);
		\draw (6.center) to (7.center);
		\draw [in=90, out=-90, looseness=1.25] (9.center) to (8);
		\draw [in=90, out=-90, looseness=1.25] (8) to (10.center);
		\draw [in=90, out=-90, looseness=1.25] (11.center) to (12);
		\draw [in=90, out=-90, looseness=1.25] (12) to (13.center);
		\draw (16.center) to (17.center);
		\draw (19.center) to (15.center);
		\draw (14.center) to (18.center);
	\end{pgfonlayer}
\end{tikzpicture} = \begin{tikzpicture}
	\begin{pgfonlayer}{nodelayer}
		\node [style=none] (0) at (-1, 2) {};
		\node [style=none] (1) at (0.5, 2) {};
		\node [style=none] (2) at (0.5, -1.5) {};
		\node [style=none] (3) at (-1, -1.5) {};
		\node [style=circle, scale=2] (4) at (-0.25, 1.25) {};
		\node [style=none] (5) at (-0.5, 2) {};
		\node [style=circle, scale=2] (6) at (-0.25, -0.75) {};
		\node [style=none] (7) at (0, -1.5) {};
		\node [style=none] (8) at (-0.25, 2.75) {};
		\node [style=none] (9) at (-0.25, -2.25) {};
		\node [style=none] (10) at (-0.25, 2) {};
		\node [style=none] (11) at (-0.25, -1.5) {};
		\node [style=none] (12) at (-0.25, -0.75) {$g$};
		\node [style=none] (13) at (-0.25, 1.25) {$f$};
		\node [style=none] (14) at (0, 2.5) {$C^\dagger$};
		\node [style=none] (15) at (0, -2) {$A^\dagger$};
	\end{pgfonlayer}
	\begin{pgfonlayer}{edgelayer}
		\draw (0.center) to (1.center);
		\draw (2.center) to (3.center);
		\draw [in=90, out=-90, looseness=1.25] (5.center) to (4);
		\draw [in=90, out=-90, looseness=1.25] (6) to (7.center);
		\draw (11.center) to (9.center);
		\draw (8.center) to (10.center);
		\draw (4) to (6);
		\draw (0.center) to (3.center);
		\draw (1.center) to (2.center);
	\end{pgfonlayer}
\end{tikzpicture}
\]

The following are the representations of the basic natural isomorphisms of a $\dagger$-LDC:
\[ \begin{array}{l l}
\lambda_\top: \top \rightarrow \bot^\dagger = 
\begin{tikzpicture}
	\begin{pgfonlayer}{nodelayer}
		\node [style=none] (0) at (-2, 0) {};
		\node [style=none] (1) at (-2, -1) {};
		\node [style=none] (2) at (-1, -1) {};
		\node [style=none] (3) at (-1, 0) {};
		\node [style=circle] (4) at (-1.5, -0.5) {$\bot$};
		\node [style=none] (5) at (-1.25, 0) {};
		\node [style=none] (6) at (-1.75, -1) {};
		\node [style=none] (7) at (-1.5, -2) {};
		\node [style=none] (8) at (-1.5, -1) {};
		\node [style=circle, scale=0.4] (9) at (-1.5, -1.5) {};
		\node [style=circle] (10) at (-2.5, -0.5) {$\top$};
		\node [style=none] (11) at (-2.5, 0.5) {};
	\end{pgfonlayer}
	\begin{pgfonlayer}{edgelayer}
		\draw [style=none] (0.center) to (1.center);
		\draw [style=none] (1.center) to (2.center);
		\draw [style=none] (2.center) to (3.center);
		\draw [style=none] (3.center) to (0.center);
		\draw [style=none] (4) to (5.center);
		\draw [style=none] (7.center) to (8.center);
		\draw [style=none] (11.center) to (10);
		\draw [densely dotted, bend right=45, looseness=1.25] (10) to (9);
	\end{pgfonlayer}
\end{tikzpicture}  & 
\lambda_\top^{-1}: \bot^\dagger \rightarrow \top =
\begin{tikzpicture}
	\begin{pgfonlayer}{nodelayer}
		\node [style=none] (0) at (-2, 0) {};
		\node [style=none] (1) at (-2, -1) {};
		\node [style=none] (2) at (-1, -1) {};
		\node [style=none] (3) at (-1, 0) {};
		\node [style=circle] (4) at (-1.5, -0.5) {$\bot$};
		\node [style=none] (5) at (-1.25, 0) {};
		\node [style=none] (6) at (-1.75, -1) {};
		\node [style=none] (7) at (-1.5, -1) {};
		\node [style=circle] (8) at (-0.5, -0.25) {$\top$};
		\node [style=none] (9) at (-1.5, 0) {};
		\node [style=none] (10) at (-1.5, 0.5) {};
		\node [style=none] (11) at (-0.5, -2) {};
		\node [style=circle, scale=0.4] (12) at (-0.5, -1.25) {};
	\end{pgfonlayer}
	\begin{pgfonlayer}{edgelayer}
		\draw [style=none] (0.center) to (1.center);
		\draw [style=none] (1.center) to (2.center);
		\draw [style=none] (2.center) to (3.center);
		\draw [style=none] (3.center) to (0.center);
		\draw [style=none] (4) to (6.center);
		\draw [style=none] (10.center) to (9.center);
		\draw [style=none] (8) to (11.center);
		\draw [densely dotted, bend right=45, looseness=1.00] (7.center) to (12);
		\draw [densely dotted] (4) to (5.center);
	\end{pgfonlayer}
\end{tikzpicture}  \\
\lambda_\bot: \bot \rightarrow \top^\dagger = \begin{tikzpicture}
	\begin{pgfonlayer}{nodelayer}
		\node [style=none] (0) at (-2, -1.5) {};
		\node [style=none] (1) at (-2, -0.5) {};
		\node [style=none] (2) at (-1, -0.5) {};
		\node [style=none] (3) at (-1, -1.5) {};
		\node [style=circle] (4) at (-1.5, -1) {$\top$};
		\node [style=none] (5) at (-1.75, -1.5) {};
		\node [style=none] (6) at (-1.25, -0.5) {};
		\node [style=none] (7) at (-1.5, -0.5) {};
		\node [style=circle] (8) at (-0.5, -1) {$\bot$};
		\node [style=none] (9) at (-1.5, -1.5) {};
		\node [style=none] (10) at (-1.5, -2) {};
		\node [style=none] (11) at (-0.5, 0.5) {};
		\node [style=circle, scale=0.4] (12) at (-0.5, 0) {};
	\end{pgfonlayer}
	\begin{pgfonlayer}{edgelayer}
		\draw [style=none] (0.center) to (1.center);
		\draw [style=none] (1.center) to (2.center);
		\draw [style=none] (2.center) to (3.center);
		\draw [style=none] (3.center) to (0.center);
		\draw [style=none] (10.center) to (9.center);
		\draw [style=none] (8) to (11.center);
		\draw [densely dotted] (4) to (5.center);
		\draw [style=none] (4) to (6.center);
		\draw [densely dotted, bend left, looseness=1.25] (7.center) to (12);
	\end{pgfonlayer}
\end{tikzpicture} &
\lambda_\bot^{-1}: \top^\dagger \rightarrow \bot = \begin{tikzpicture}
	\begin{pgfonlayer}{nodelayer}
		\node [style=none] (0) at (-2, 0) {};
		\node [style=none] (1) at (-2, -1) {};
		\node [style=none] (2) at (-1, -1) {};
		\node [style=none] (3) at (-1, 0) {};
		\node [style=circle] (4) at (-1.5, -0.5) {$\top$};
		\node [style=none] (5) at (-1.75, 0) {};
		\node [style=none] (6) at (-1.25, -1) {};
		\node [style=none] (7) at (-1.5, -1) {};
		\node [style=circle] (8) at (-0.5, -0.5) {$\bot$};
		\node [style=none] (9) at (-1.5, 0) {};
		\node [style=none] (10) at (-1.5, 0.5) {};
		\node [style=none] (11) at (-0.5, -2) {};
		\node [style=circle,scale=0.4] (12) at (-0.5, -1.25) {};
	\end{pgfonlayer}
	\begin{pgfonlayer}{edgelayer}
		\draw [style=none] (0.center) to (1.center);
		\draw [style=none] (1.center) to (2.center);
		\draw [style=none] (2.center) to (3.center);
		\draw [style=none] (3.center) to (0.center);
		\draw [style=none] (10.center) to (9.center);
		\draw [style=none] (8) to (11.center);
		\draw [densely dotted] (4) to (5.center);
		\draw [style=none] (4) to (6.center);
		\draw [densely dotted, bend right=45, looseness=1.25] (7.center) to (12);
	\end{pgfonlayer}
\end{tikzpicture}\\
\lambda_\ox: A^\dagger \ox B^\dagger \rightarrow (A \oa B)^\dagger = \begin{tikzpicture}[rotate=180]
	\begin{pgfonlayer}{nodelayer}
		\node [style=oa] (0) at (-3, 0) {};
		\node [style=none] (1) at (-3.5, 0.75) {};
		\node [style=none] (2) at (-2.5, 0.75) {};
		\node [style=none] (3) at (-3.5, -0.75) {};
		\node [style=none] (4) at (-2.5, -0.75) {};
		\node [style=ox] (5) at (-3, -1.5) {};
		\node [style=none] (6) at (-3, -2.25) {};
		\node [style=none] (7) at (-3, 1.5) {};
		\node [style=none] (8) at (-3, -0.75) {};
		\node [style=none] (9) at (-3, 0.75) {};
		\node [style=none] (10) at (-4, 0.75) {};
		\node [style=none] (11) at (-2, 0.75) {};
		\node [style=none] (12) at (-2, -0.75) {};
		\node [style=none] (13) at (-4, -0.75) {};
	\end{pgfonlayer}
	\begin{pgfonlayer}{edgelayer}
		\draw [style=none, bend right, looseness=1.00] (1.center) to (0);
		\draw [style=none, bend right, looseness=1.00] (0) to (2.center);
		\draw [style=none] (0) to (8.center);
		\draw [style=none, bend right, looseness=1.00] (3.center) to (5);
		\draw [style=none, bend right, looseness=1.00] (5) to (4.center);
		\draw [style=none] (5) to (6.center);
		\draw [style=none] (7.center) to (9.center);
		\draw [style=none] (10.center) to (13.center);
		\draw [style=none] (13.center) to (12.center);
		\draw [style=none] (12.center) to (11.center);
		\draw [style=none] (11.center) to (10.center);
	\end{pgfonlayer}
	\end{tikzpicture}  &
\lambda_\oa: A^\dagger \ox B^\dagger \rightarrow (A \oa B)^\dagger = \begin{tikzpicture}[rotate=180]
	\begin{pgfonlayer}{nodelayer}
		\node [style=ox] (0) at (-3, 0) {};
		\node [style=none] (1) at (-3.5, 0.75) {};
		\node [style=none] (2) at (-2.5, 0.75) {};
		\node [style=none] (3) at (-3.5, -0.75) {};
		\node [style=none] (4) at (-2.5, -0.75) {};
		\node [style=oa] (5) at (-3, -1.5) {};
		\node [style=none] (6) at (-3, -2.25) {};
		\node [style=none] (7) at (-3, 1.5) {};
		\node [style=none] (8) at (-3, -0.75) {};
		\node [style=none] (9) at (-3, 0.75) {};
		\node [style=none] (10) at (-4, 0.75) {};
		\node [style=none] (11) at (-2, 0.75) {};
		\node [style=none] (12) at (-2, -0.75) {};
		\node [style=none] (13) at (-4, -0.75) {};
	\end{pgfonlayer}
	\begin{pgfonlayer}{edgelayer}
		\draw [style=none, bend right, looseness=1.00] (1.center) to (0);
		\draw [style=none, bend right, looseness=1.00] (0) to (2.center);
		\draw [style=none] (0) to (8.center);
		\draw [style=none, bend right, looseness=1.00] (3.center) to (5);
		\draw [style=none, bend right, looseness=1.00] (5) to (4.center);
		\draw [style=none] (5) to (6.center);
		\draw [style=none] (7.center) to (9.center);
		\draw [style=none] (10.center) to (13.center);
		\draw [style=none] (13.center) to (12.center);
		\draw [style=none] (12.center) to (11.center);
		\draw [style=none] (11.center) to (10.center);
	\end{pgfonlayer}
\end{tikzpicture} \\
\lambda_\otimes^{-1}: (A \oa B)^\dagger \to A^\dagger \ox B^\dagger =
\begin{tikzpicture}
	\begin{pgfonlayer}{nodelayer}
		\node [style=oa] (0) at (-3, 0) {};
		\node [style=none] (1) at (-3.5, 0.75) {};
		\node [style=none] (2) at (-2.5, 0.75) {};
		\node [style=none] (3) at (-3.5, -0.75) {};
		\node [style=none] (4) at (-2.5, -0.75) {};
		\node [style=ox] (5) at (-3, -1.5) {};
		\node [style=none] (6) at (-3, -2.25) {};
		\node [style=none] (7) at (-3, 1.5) {};
		\node [style=none] (8) at (-3, -0.75) {};
		\node [style=none] (9) at (-3, 0.75) {};
		\node [style=none] (10) at (-4, 0.75) {};
		\node [style=none] (11) at (-2, 0.75) {};
		\node [style=none] (12) at (-2, -0.75) {};
		\node [style=none] (13) at (-4, -0.75) {};
	\end{pgfonlayer}
	\begin{pgfonlayer}{edgelayer}
		\draw [style=none, bend right, looseness=1.00] (1.center) to (0);
		\draw [style=none, bend right, looseness=1.00] (0) to (2.center);
		\draw [style=none] (0) to (8.center);
		\draw [style=none, bend right, looseness=1.00] (3.center) to (5);
		\draw [style=none, bend right, looseness=1.00] (5) to (4.center);
		\draw [style=none] (5) to (6.center);
		\draw [style=none] (7.center) to (9.center);
		\draw [style=none] (10.center) to (13.center);
		\draw [style=none] (13.center) to (12.center);
		\draw [style=none] (12.center) to (11.center);
		\draw [style=none] (11.center) to (10.center);
	\end{pgfonlayer}
\end{tikzpicture} 
&
\lambda_\oa^{-1}:  (A \ox B)^\dagger \to A^\dagger \oa B^\dagger =
\begin{tikzpicture}
	\begin{pgfonlayer}{nodelayer}
		\node [style=ox] (0) at (-3, 0) {};
		\node [style=none] (1) at (-3.5, 0.75) {};
		\node [style=none] (2) at (-2.5, 0.75) {};
		\node [style=none] (3) at (-3.5, -0.75) {};
		\node [style=none] (4) at (-2.5, -0.75) {};
		\node [style=oa] (5) at (-3, -1.5) {};
		\node [style=none] (6) at (-3, -2.25) {};
		\node [style=none] (7) at (-3, 1.5) {};
		\node [style=none] (8) at (-3, -0.75) {};
		\node [style=none] (9) at (-3, 0.75) {};
		\node [style=none] (10) at (-4, 0.75) {};
		\node [style=none] (11) at (-2, 0.75) {};
		\node [style=none] (12) at (-2, -0.75) {};
		\node [style=none] (13) at (-4, -0.75) {};
	\end{pgfonlayer}
	\begin{pgfonlayer}{edgelayer}
		\draw [style=none, bend right, looseness=1.00] (1.center) to (0);
		\draw [style=none, bend right, looseness=1.00] (0) to (2.center);
		\draw [style=none] (0) to (8.center);
		\draw [style=none, bend right, looseness=1.00] (3.center) to (5);
		\draw [style=none, bend right, looseness=1.00] (5) to (4.center);
		\draw [style=none] (5) to (6.center);
		\draw [style=none] (7.center) to (9.center);
		\draw [style=none] (10.center) to (13.center);
		\draw [style=none] (13.center) to (12.center);
		\draw [style=none] (12.center) to (11.center);
		\draw [style=none] (11.center) to (10.center);
	\end{pgfonlayer}
\end{tikzpicture} 
\end{array}
\]

Dagger boxes interact with involutor $A \xrightarrow{\iota} A^{\dagger\dagger}$ as follows:
$$
\begin{tikzpicture}
\begin{pgfonlayer}{nodelayer}
\node [style=circle] (0) at (-3.5, 0) {$f$};
\node [style=none] (1) at (-4.5, 1) {};
\node [style=none] (2) at (-2.5, 1) {};
\node [style=none] (3) at (-4.5, -0.5) {};
\node [style=none] (4) at (-2.5, -0.5) {};
\node [style=none] (5) at (-4.75, -1) {};
\node [style=none] (6) at (-2.25, -1) {};
\node [style=none] (7) at (-2.25, 1.25) {};
\node [style=none] (8) at (-4.75, 1.25) {};
\node [style=circle] (9) at (-4.25, 1.75) {$\iota$};
\node [style=none] (10) at (-4.25, 2.25) {};
\node [style=none] (11) at (-4.25, 1.25) {};
\node [style=none] (12) at (-3.5, -1) {};
\node [style=none] (13) at (-3.5, -1.5) {};
\node [style=none] (14) at (-2.75, 1) {};
\node [style=none] (15) at (-3.5, -0.5) {};
\node [style=none] (16) at (-3.5, 1.25) {};
\node [style=none] (17) at (-4.25, -1) {};
\node [style=none] (18) at (-4.25, -0.5) {};
\node [style=none] (19) at (-3.5, 1) {};
\node [style=circle] (20) at (-2.75, 1.75) {$\iota$};
\node [style=none] (21) at (-2.75, 1.25) {};
\node [style=none] (22) at (-2.75, 2.25) {};
\node [style=none] (23) at (-4.25, 1) {};
\node [style=none] (24) at (-2.75, -0.5) {};
\node [style=none] (25) at (-2.75, -1) {};
\node [style=none] (26) at (-3.5, 0.75) {$\cdots$};
\node [style=none] (27) at (-3.5, 1.75) {$\cdots$};
\node [style=none] (28) at (-3.5, -0.75) {$\cdots$};
\end{pgfonlayer}
\begin{pgfonlayer}{edgelayer}
\draw [style=none] (5.center) to (6.center);
\draw [style=none] (6.center) to (7.center);
\draw [style=none] (7.center) to (8.center);
\draw [style=none] (8.center) to (5.center);
\draw [style=none] (1.center) to (3.center);
\draw [style=none] (3.center) to (4.center);
\draw [style=none] (4.center) to (2.center);
\draw [style=none] (2.center) to (1.center);
\draw [style=none, in=-90, out=60, looseness=1.00] (0) to (14.center);
\draw [style=none] (0) to (15.center);
\draw [style=none] (18.center) to (17.center);
\draw [style=none] (17.center) to (5.center);
\draw [style=none] (19.center) to (16.center);
\draw [style=none] (12.center) to (13.center);
\draw [style=none] (11.center) to (9);
\draw [style=none] (10.center) to (9);
\draw [style=none] (21.center) to (20);
\draw [style=none] (22.center) to (20);
\draw [style=none, in=120, out=-90, looseness=1.00] (23.center) to (0);
\draw [style=none] (25.center) to (24.center);
\end{pgfonlayer}
\end{tikzpicture}
=
\begin{tikzpicture}
\begin{pgfonlayer}{nodelayer}
\node [style=circle] (0) at (-2.5, 1) {$f$};
\node [style=none] (1) at (-2, 2) {};
\node [style=none] (2) at (-2.5, -1) {};
\node [style=circle] (3) at (-2.5, 0) {$i$};
\node [style=none] (4) at (-3, 2) {};
\node [style=none] (5) at (-2.5, 1.75) {$\cdots$};
\end{pgfonlayer}
\begin{pgfonlayer}{edgelayer}
\draw [style=none, in=45, out=-90, looseness=1.00] (1.center) to (0);
\draw [style=none] (0) to (3);
\draw [style=none] (3) to (2.center);
\draw [style=none, in=135, out=-90, looseness=1.00] (4.center) to (0);
\end{pgfonlayer}
\end{tikzpicture}
$$

It is worth noting that one need not have a legal proof net inside a $\dagger$-box. 
This complicates the correctness criterion. However, the required 
correctness criterion is discussed in \cite{MP05}.

\subsection{Functors for $\dagger$-linearly distributive categories}
\label{Sec: dagger linear}

Clearly the functors and transformations between $\dagger$-LDCs must ``preserve'' the dagger in some sense.  Precisely we have:

\begin{definition}
	$F: \X \to \Y$ is a {\bf $\dagger$-linear functor} between $\dagger$-LDCs when $F$ is a linear functor equipped with a linear natural isomorphism 
	$\rho^F= (\rho_\ox^F: F_\ox(A^\dagger) \to F_\oa(A)^\dagger ,\rho_\oa^F:  F_\ox(A)^\dagger \to F_\oa(A^\dagger))$ called the {\bf preservator}, 
	such that  the following diagrams commute:
	\[ 
	\xymatrix{
		F_\ox(X) \ar[r]^{\iota} \ar[d]_{F_\ox(\iota)} \ar@{}[dr]|{\mbox{\tiny {\bf [$\dagger$-LF.1]}}} & 
		F_\ox(X)^{\dagger \dagger} \ar@{<-}[d]^{(\rho^F_\oa)^\dagger} \\
		F_\ox(X^{\dagger \dagger}) \ar[r]_{\rho^F_\ox} & F_\oa(X^\dagger)^\dagger
	} ~~~~~~~~~ \xymatrix{
		F_\oa(X) \ar[r]^{\iota} \ar[d]_{F_\oa(\iota)}  \ar@{}[dr]|{\mbox{\tiny {\bf [$\dagger$-LF.2]}}} & F_\oa(X)^{\dagger \dagger} \ar[d]^{(\rho^F_\ox)^\dagger} \\
		F_\oa(X^{\dagger \dagger}) \ar@{<-}[r]_{\rho^F_\oa} & F_\ox(X^\dagger)^\dagger
	}
	\]
\end{definition}

In case that $F$ is a normal mix functor between $\dagger$-isomix categories, then by Lemma \ref{Lemma: isomix functor}, $F$ is an isomix functor and, therefore by Corollary \ref{Corollary: normal-nat-iso},  the preservators become inverses, $\rho^F_\ox = (\rho^F_\oa)^{-1}$.  
This means the squares {\bf [$\dagger$-LF.1]} and {\bf [$\dagger$-LF.2]} coincide to give a single condition for the tensor preservator:
\[ \xymatrix{
	F(X) \ar[r]^{\iota} \ar[d]_{F(\iota)} \ar@{}[dr]|{\mbox{\tiny {\bf [$\dagger$-isomix]}}} & 
	F(X)^{\dagger \dagger} \ar@{->}[d]^{(\rho^F_\ox)^\dagger} \\
	F(X^{\dagger \dagger}) \ar[r]_{\rho^F_\ox} & F(X^\dagger)^\dagger
} \]


In case when $F$ is an isomix functor, by Lemma \ref{Lemma: Frobenius linear transformation}, $\rho := \rho_\ox$ is monoidal on $\ox$ and comonoidal on $\oa$:

\[
\bf{[P.1]} ~~~~~~~
\xymatrix{
	F(A^\dagger) \ox F(B^\dagger) \ar[r]^{\rho \ox \rho} \ar[d]_{m_\ox} \ar@{}[ddr]|{(a)} & F(A)^\dagger \ox F(B)^\dagger \ar[d]^{\lambda_\ox} \\
	F(A^\dagger \ox B^\dagger) \ar[d]_{F(\lambda_\ox)} & (F(A) \ox F(B))^\dagger \ar[d]^{n_\oa^\dagger} \\
	F((A \oa B)^\dagger) \ar[r]_{\rho} & (F(A\oa B))^\dagger
}~~~~~~~~~ \xymatrix{
	\top \ar[d]_{m_\top} \ar[dr]^{\lambda_\top} \ar@{}[ddr]|{(b)} & \\
	F(\top) \ar[d]_{F(\lambda_\top)} & \bot^\dagger \ar[dr]^{n_\bot^\dagger} \\
	F(\bot^\dagger) \ar[rr]_{\rho} & & (F(\bot))^\dagger
}
\]

\[
\bf{[P.2]} ~~~~~~~
\xymatrix{
	F((A \ox B)^\dagger) \ar[r]^{\rho} \ar[d]_{F(\lambda_\oa^{-1})} \ar@{}[ddr]|{(a)} & F(A \ox B)^\dagger \ar[d]^{m_\ox^\dagger} \\
	F(A^\dagger \oa B^\dagger) \ar[d]_{n_\oa^F} & (F(A) \ox F(B))^\dagger \ar[d]^{\lambda_\oa^{-1}} \\
	F(A^\dagger) \oa F(B^\dagger) \ar[r]_{\rho \oa \rho} & F(A)^\dagger \oa F(B)^\dagger
} ~~~~~~~~~ \xymatrix{
	F(\top^\dagger) \ar[dd]_{\rho} \ar[dr]^{F(\lambda_\bot^{-1})} \ar@{}[ddr]|{(b)} & \\
	& F(\bot) \ar[dr]^{n_\bot} & \\
	F(\top)^\dagger \ar[r]_{m_\top^\dagger} & \top^\dagger \ar[r]_{\lambda_\bot^{-1}} & \bot
}
\]

Pictorial representation of {\bf[P.2]-(a)} is as follows:
\[
\begin{tikzpicture} 
\begin{pgfonlayer}{nodelayer}
\node [style=none] (0) at (1.75, 3.25) {};
\node [style=none] (1) at (2.25, -0) {};
\node [style=none] (2) at (1.25, -0) {};
\node [style=none] (3) at (0.25, -0) {};
\node [style=none] (4) at (0.75, 2.25) {};
\node [style=ox] (5) at (1.75, 1.5) {};
\node [style=none] (6) at (2.75, 2.25) {};
\node [style=none] (7) at (2.75, 1) {};
\node [style=none] (8) at (3.25, 2.75) {};
\node [style=none] (9) at (1.25, 2.25) {};
\node [style=circle] (10) at (2.25, -3) {$\rho$};
\node [style=none] (11) at (1.75, 2.25) {};
\node [style=none] (12) at (1.75, 1) {};
\node [style=none] (13) at (2.25, -3.5) {};
\node [style=none] (14) at (0.75, 1) {};
\node [style=none] (15) at (3, 0.25) {$M$};
\node [style=none] (16) at (3.25, -0) {};
\node [style=circle] (17) at (1.25, -3) {$\rho$};
\node [style=none] (18) at (1.25, -3.5) {};
\node [style=none] (19) at (2.25, 2.25) {};
\node [style=none] (20) at (0.25, 2.75) {};
\node [style=none] (21) at (1.25, -0.75) {};
\node [style=ox] (22) at (1.75, -1.25) {};
\node [style=none] (23) at (2.25, -0.75) {};
\node [style=none] (24) at (2.75, -2) {};
\node [style=none] (25) at (0.75, -0.75) {};
\node [style=none] (26) at (0.75, -2) {};
\node [style=none] (27) at (1.25, -2) {};
\node [style=none] (28) at (1.75, -0.75) {};
\node [style=none] (29) at (2.75, -0.75) {};
\node [style=none] (30) at (2.25, -2) {};
\node [style=none] (31) at (2.25, 1) {};
\node [style=none] (32) at (1.25, 1) {};
\node [style=none] (33) at (1.75, -0) {};
\node [style=ox] (34) at (1.75, 0.5) {};
\node [style=none] (35) at (3, 3) {$F((A \ox B)^\dagger)$};
\node [style=none] (36) at (3, -0.25) {$F(A^\dagger \ox B^\dagger)$};
\node [style=none] (37) at (0.25, -2.5) {$F(A^\dagger)$};
\node [style=none] (38) at (3, -2.5) {$F(B^\dagger)$};
\node [style=none] (39) at (0, -3.5) {$F(A)^\dagger$};
\node [style=none] (40) at (3, -3.5) {$F(B)^\dagger$};
\end{pgfonlayer}
\begin{pgfonlayer}{edgelayer}
\draw [bend right, looseness=1.50] (9.center) to (5);
\draw [bend right, looseness=1.50] (5) to (19.center);
\draw (5) to (12.center);
\draw (4.center) to (6.center);
\draw (6.center) to (7.center);
\draw (14.center) to (7.center);
\draw (14.center) to (4.center);
\draw (20.center) to (3.center);
\draw (3.center) to (16.center);
\draw (16.center) to (8.center);
\draw (8.center) to (20.center);
\draw (17) to (18.center);
\draw (10) to (13.center);
\draw (11.center) to (0.center);
\draw [bend left, looseness=1.50] (27.center) to (22);
\draw [bend left, looseness=1.50] (22) to (30.center);
\draw (22) to (28.center);
\draw (26.center) to (24.center);
\draw (24.center) to (29.center);
\draw (25.center) to (29.center);
\draw (25.center) to (26.center);
\draw [bend right, looseness=1.50] (32.center) to (34);
\draw [bend right, looseness=1.50] (34) to (31.center);
\draw (34) to (33.center);
\draw (33.center) to (28.center);
\draw (27.center) to (17);
\draw (30.center) to (10);
\end{pgfonlayer}
\end{tikzpicture} = \begin{tikzpicture} 
\begin{pgfonlayer}{nodelayer}
\node [style=none] (0) at (1.5, 0.5) {};
\node [style=none] (1) at (-0.25, -0.75) {};
\node [style=none] (2) at (-0.25, 1.75) {};
\node [style=none] (3) at (1.25, -0.75) {};
\node [style=ox] (4) at (1, -0) {};
\node [style=none] (5) at (0.5, 0.5) {};
\node [style=none] (6) at (2.25, 1.75) {};
\node [style=none] (7) at (1.5, 0.75) {};
\node [style=none] (8) at (1.25, -1.25) {};
\node [style=none] (9) at (1.25, 1.75) {};
\node [style=ox] (10) at (1, 1.25) {};
\node [style=none] (11) at (2, 0.75) {};
\node [style=none] (12) at (0.5, 1.75) {};
\node [style=none] (13) at (1, -0.75) {};
\node [style=none] (14) at (2, -0.5) {};
\node [style=none] (15) at (0, -0.5) {};
\node [style=circle] (16) at (1.25, 2.5) {$\rho$};
\node [style=none] (17) at (0, 0.75) {};
\node [style=none] (18) at (2.25, -0.75) {};
\node [style=none] (19) at (1, 1.75) {};
\node [style=none] (20) at (1.75, -0.25) {$M$};
\node [style=none] (21) at (0.5, 0.75) {};
\node [style=none] (22) at (1.25, 3.25) {};
\node [style=none] (23) at (0, -1.25) {};
\node [style=none] (24) at (0.5, -1.25) {};
\node [style=ox] (25) at (1, -2) {};
\node [style=none] (26) at (1, -2.5) {};
\node [style=none] (27) at (2, -1.25) {};
\node [style=none] (28) at (0, -2.5) {};
\node [style=none] (29) at (1.5, -1.25) {};
\node [style=none] (30) at (2, -2.5) {};
\node [style=none] (31) at (0.5, -2.5) {};
\node [style=none] (32) at (0.5, -3) {};
\node [style=none] (33) at (1.5, -3) {};
\node [style=none] (34) at (1.5, -2.5) {};
\node [style=none] (35) at (2.25, 3) {$F((A \ox B)^\dagger)$};
\node [style=none] (36) at (3.25, -1) {$F((A) \ox F(B))^\dagger$};
\node [style=none] (37) at (-0.5, -2.75) {$F(A)^\dagger$};
\node [style=none] (38) at (2.75, -2.75) {$F(B)^\dagger$};
\node [style=none] (39) at (2.25, 2.25) {};
\end{pgfonlayer}
\begin{pgfonlayer}{edgelayer}
\draw [bend left, looseness=1.50] (5.center) to (10);
\draw [bend left, looseness=1.50] (10) to (0.center);
\draw (1.center) to (18.center);
\draw (18.center) to (6.center);
\draw (2.center) to (6.center);
\draw (2.center) to (1.center);
\draw [bend right, looseness=1.50] (21.center) to (4);
\draw [bend right, looseness=1.50] (4) to (7.center);
\draw (4) to (13.center);
\draw (17.center) to (11.center);
\draw (11.center) to (14.center);
\draw (15.center) to (14.center);
\draw (15.center) to (17.center);
\draw (16) to (9.center);
\draw (19.center) to (10);
\draw (3.center) to (8.center);
\draw (22.center) to (16);
\draw [bend right, looseness=1.50] (24.center) to (25);
\draw [bend right, looseness=1.50] (25) to (29.center);
\draw (25) to (26.center);
\draw (23.center) to (27.center);
\draw (27.center) to (30.center);
\draw (28.center) to (30.center);
\draw (28.center) to (23.center);
\draw (31.center) to (32.center);
\draw (34.center) to (33.center);
\end{pgfonlayer}
\end{tikzpicture}
\]


For linear natural transformations $(\beta_\ox, \beta_\oa): F \to G$ between $\dagger$-linear functors, we demand that $\beta_\ox$ and $\beta_\oa$ are related by:

\[ \xymatrix{F_\ox(A^\dagger) \ar[d]_{\rho^F_\ox} \ar[rr]^{\beta_\ox} && G_\ox(A^\dagger)  \ar[d]^{\rho^G_\ox} \\
	(F_\oa(X))^\dagger \ar[rr]_{\beta_\oa^\dagger} & &  (G_\oa(X))^\dagger}
~~~~~~
\xymatrix{(G_\ox(X))^\dagger \ar[d]_{\rho^G_\oa} \ar[rr]^{\beta_\ox^\dagger} &&  (F_\ox(X))^\dagger \ar[d]^{\rho^F_\oa} \\
	G_\oa(A^\dagger) \ar[rr]_{\beta_\oa} & &  F_\oa(A^\dagger)} \]
Notice that this means that $\beta_\ox$ is completely determined by $\beta_\oa$ in the following sense:
\[ \xymatrix{ F_\ox(A) \ar[d]_{F_\ox(\iota)}  \ar[rr]^{\beta_\ox} && G_\ox(A) \ar[d]^{G_\ox(\iota)} \\
	F_\ox(A^{\dagger\dagger})  \ar[d]_{\rho^F_\ox} \ar[rr]^{\beta_\ox} & & G_\ox(A^{\dagger\dagger}) \ar[d]^{\rho^G_\ox}  \\
	F_\oa(A^\dagger)^\dagger \ar[rr]_{\beta_\oa^\dagger} && G_\oa(A^\dagger)^\dagger } \]
Because the vertical maps are isomorphisms, this diagram can be used to express $\beta_\ox$ in terms of $\beta_\oa$.  Similarly $\beta_\oa$ can be expressed in terms 
of $\beta_\ox$.  Thus, it is possible to express the coherences in terms of just one of these transformations.


\section{Examples: $\dagger$-LDCs}
\label{subsection: Dagger LDC examples}

In this section, we discuss a few basic examples of $\dagger$-isomix categories. 
The first two are compact LDCs. All these examples have a non-stationary dagger functor.  
More examples of $\dagger$-isomix categories can be found in the next chapter when we discuss the dagger 
and the conjugation functor.


\subsection{Every $\dagger$-monoidal category is a $\dagger$-LDC}

A $\dagger$-monoidal category \cite{Sel07} is defined as a symmetric monoidal category, $\X$, with a  
contravariant functor $\dagger: \X^\op \to \X$ which is stationary on objects ($A = A^\dag$) such that:
\begin{enumerate}[(i)]
	\item for all $f$, $f^{\dag \dag} = f$
	\item for all $f$ and $g$, $(f \ox g)^\dag = f^\dag \ox g^\dag$
	\item $a_\ox^\dag = a_\ox^{-1}$
	\item $(u_\ox^l)^\dag = (u_\ox^l)^{-1}$
	\item $(u_\ox^r)^\dag = (u_\ox^r)^{-1}$
\end{enumerate}

Note that every $\dagger$-monoidal category is a compact $\dagger$-LDC in which the 
the laxors and the involutor are identity transformations. 

\subsection{Finite dimensional framed vector spaces}
\label{subsection:fdfv}

In this section we describe the category of ``framed'' finite dimensional vector spaces, where a frame in this context is just a 
choice of basis.  Thus, the objects in this category are vector spaces with a chosen basis while the maps, ignoring the basis, are simply homomorphisms of 
the vector spaces.

The category of finite dimensional framed vectors spaces, ${\sf FFVec}_K$, is a monoidal category defined as follows:
\begin{description}
\item[Objects:]  The objects are pairs $(V,{\cal V})$ where $V$ is a finite dimensional $K$-vector space and ${\cal V} = \{ v_1,...,v_n \}$ is a basis. 
\item[Maps:]  A map $(V,{\cal V}) \xrightarrow{f} (W,{\cal W})$ is a linear map  $V \xrightarrow{f} W$ in ${\sf FdVec}_K$.
\item[Tensor:] $(V,{\cal V}) \ox (W,{\cal W}) = (V \ox W,\{ v \ox w | v \in {\cal V}, w \in {\cal W} \})$ where $V \ox W$ is the usual tensor product.  The unit is 
$(K,\{ e \})$ where $e$ is the unit of the field $K$.
\end{description}

To define the ``dagger'' we  must first choose a conjugation $\overline{(\_)}: K \to K$ (see more details in Section \ref{Sec: conjugation}), that is a field homomorphism with $k = \overline{(\overline{k})}$. The 
canonical example of conjugation is conjugation of the complex numbers, however, the conjugation can be arbitrarily chosen -- so could also, for example, be the identity.  
This conjugation then can be extended to a (covariant) functor:
\[ \overline{(\_)}: {\sf FFVec}_K \to {\sf FFVec}_K; \begin{array}[c]{c} \xymatrix{(V,{\cal V}) \ar[d]^{f} \\ (W,{\cal W})} \end{array} 
                 \mapsto \begin{array}[c]{c} \xymatrix{\overline{(V,{\cal V}) } \ar[d]^{\overline{f}} \\ \overline{(W,{\cal W}) }} \end{array} \]
where $\overline{(V,{\cal V})}$ is the vector space with the same basis but with the conjugate action $c~\overline{\cdot}~v = \overline{c} \cdot  v$.  The conjugate 
homomorphism, $\overline{f}$, is then the same underlying map which is homomorphism between the conjugate spaces.

${\sf FFVec}_K$ is also a compact closed category with $(V,{\cal B})^{*} = (V^{*}, \{ \widetilde{b_i} | b_i \in {\cal B} \})$ where 
\[ V^{*} = V \multimap K~~~~\mbox{and}~~~~\widetilde{b_i}: V \to K; \left(\sum_j \beta_j \cdot b_j \right) \mapsto  \beta_i \]
This makes $(\_)^{*}: {\sf FFVec}_K^{\rm op} \to {\sf FFVec}_K$ a contravariant functor whose action is determined by precomposition.   Finally, we 
define the ``dagger'' to be the composite $(V,{\cal B})^\dagger = \overline{(V,{\cal B})^{*}}$.

This is a monoidal category with tensor and par being identified (so the linear distribution is the associator) and is isomix.  We must show that it is a $\dagger$-LDC.  
Towards this aim we define the required natural transformations on the basis:

\[ \lambda_\ox = \lambda_\oa: (V,{\cal V})^\dag \ox (W,{\cal W})^\dag \to ((V,{\cal V}) \ox (W,{\cal W}))^\dag;  \widetilde{v_i} \ox \widetilde{w_j} \mapsto \widetilde{v_i \ox w_j} \]

\[ \lambda_\top = \lambda_\bot: (K,\{ e\}) \to (K,\{ e\})^\dag; k \mapsto \overline{k} \]
\[ \iota: (V,{\cal V}) \to ((V,{\cal V})^\dag)^\dag; v \mapsto \lambda f. f(v) \]
Note that the last transformation is given in a basis independent manner. Importantly, it may also be given in a basis dependent manner as
$\iota(v_i) = \widetilde{\widetilde{v_i}}$ as the behaviour of these two maps is the  same when applied to the basis of $(V,{\cal V})^\dag$ namely 
the elements $\widetilde{v_j}$:
\[ \iota(v_i)(\widetilde{v_j}) = (\lambda f. f(v_i)) \widetilde{v_j} = \widetilde{v_j} v_i = \partial_{i,j} = \widetilde{\widetilde{v_i}}(\widetilde{v_j})  \]
Also note that $\widetilde{v_i \ox w_j} = (\widetilde{v_i} \ox \widetilde{w_j}) u_\ox$, where $u_\ox: K \ox K \to K$ is the multiplication of the field.
With these definitions in hand it is straightforward to check that this gives a $\dagger$-LDC by checking the required coherences on basis elements.
To demonstrate the technique consider the coherence {\bf [$\dagger$-ldc.4]}:
\[ \xymatrix{A \oa B \ar[d]_{\iota \oa \iota} \ar[rr]^\iota & & ((A \oa B)^\dag)^\dag \ar[d]^{\lambda_\ox^\dag} \\
                   (A^\dag)^\dag \oa (B^\dag)^\dag  \ar[rr]_{\lambda_\oa} & & (A^\dag \ox B^\dag)^\dag} \]
 We must show (identifying tensor and par) that $\lambda_\ox^\dag (\iota(a_i \ox b_j)) = \lambda_\ox(\iota \ox \iota(a_i \ox b_j))$.  Now the result  is a 
 higher-order term so it suffices to show the evaluations on basis elements are the same.  This means we need to show:
 $\lambda_\ox^\dag (\iota(a_i \ox b_j))(\widetilde{a_p} \ox \widetilde{b_q}) = \lambda_\ox(\iota \ox \iota(a_i \ox b_j))(\widetilde{a_p} \ox \widetilde{b_q})$
 \begin{eqnarray*}
(\lambda_\ox(\iota \ox \iota(a_i \ox b_j)))(\widetilde{a_p} \ox \widetilde{b_q}) & = & (\lambda_\ox(\widetilde{\widetilde{a_i}} \ox \widetilde{\widetilde{b_j}}))(\widetilde{a_p} \ox \widetilde{b_q}) \\
&  = & (\widetilde{\widetilde{a_i} \ox \widetilde{b_j}}) (\widetilde{a_p} \ox \widetilde{b_q}) \\
& = & (\widetilde{a_p} \ox \widetilde{b_q})(\widetilde{\widetilde{a_i}} \ox \widetilde{\widetilde{b_j}}) u_\ox ~~~\mbox{(diagrammatic order)}\\
& = & \partial_{p,i} \partial_{q,j} \\
(\lambda_\ox^\dag (\iota(a_i \ox b_j)))(\widetilde{a_p} \ox \widetilde{b_q}) 
& = & (\lambda_\ox^\dag (\widetilde{\widetilde{a_i \ox b_j}}))(\widetilde{a_p} \ox \widetilde{b_q})  \\ 
& = & (\widetilde{a_p} \ox \widetilde{b_q}) \lambda_\ox \widetilde{\widetilde{a_i \ox b_j}} ~~~\mbox{(diagrammatic order)}\\
& = & \widetilde{a_p \ox b_q} \widetilde{\widetilde{a_i \ox b_j}} \\
& = & \partial_{p,i} \partial_{q,j} 
\end{eqnarray*}

Thus, $ {\sf FFVec}_K$ is a compact $\dagger$-isomix category where the $\dagger$ functor shifts objects i.e., $A \neq A^\dagger$.


\subsection{Category of abstract state spaces}
\label{Sec: Asp}

This source of examples for $\dagger$-isomix categories is inspired by the category of convex operational models \cite{BaW11}. 
The following is a way to construct a new $\dagger$-isomix category, the category of abstract state spaces, from an existing one.

\begin{definition}
	Let $\X$ be a $\dagger$-isomix category. Define $\Asp(\X)$ as follows:
	\begin{description}
		\item[Objects:] $(A, e_A:A \to \bot, u_A: \top \to A)$
		\item[Arrows:]  $f: A \to B \in \X$ such that the following diagram commutes:
		
		$
		\xymatrix{
			& \top \ar[dl]_{u_A}  \ar[dr]^{u_B} & \\
			A \ar[rr]^{f} \ar[dr]_{e_A}  & & B \ar[ld]^{e_B} \\
			& \bot &
		}
		$
	\end{description}
	Identity arrow and composition are inherited directly from $\X$. $\Asp(\X)$ is a LDC:
	\begin{description}
		\item[$\ox$ on objects:] $(A, e_A, u_A) \ox (B, e_B, u_B) := (A \ox B, e', u')$ where, $e' := \mx (e_A \oa e_B) u_\oa$ and $u' := u_\ox^{-1} (u_A \ox u_B)$. The unit of $\ox$ is given by $(\top, \m^{-1}: \top \to \bot, 1_\top)$.
		\item[$\oa$ on objects:] $(A, e_A, u_A) \oa (B, e_B, u_B) := (A \oa B, e', u')$ where, $e' := (e_A \oa e_B) u_\oa$ and $u' := u_\ox^{-1} (u_A \ox u_B) \mx$. The unit of $\oa$ is $(\bot, 1_\bot, \m^{-1}: \top \to \bot)$
	\end{description}
\end{definition}

$\Asp(\X)$ is also $\dagger$-isomix category with \[ (A, e, u)^\dagger := 
(A^\dagger, u^\dagger \lambda_\bot^{-1}, \lambda_\top e^\dagger)\] 
All the basic natural isomorphisms are inherited from $\X$. Hence, $\Asp(\X)$ is a $\dagger$-isomix category.

%% file: chapter4.tex

\section{Dagger, duals, and conjugation}
\label{daggers-duals-conjugation}


The goal of this section is to review the interaction of the dualizing, conjugation and dagger functors. 
In dagger compact closed categories, the dagger functor $(\_)^\dagger$, and the dualizing functor $(\_)^*$ 
commute with each other and their composite gives the conjugate functor $(\_)_*$. Similary, $(\_)_*$ and $(\_)^*$ 
when composed gives the dagger functor.  Our aim is to  generalize these interactions to $\dagger$-LDCs and to achieve 
this at a reasonable level of abstraction.   To achieve this we shall need the notion which we here refer to as ``conjugation'' 
but was investigated by Egger in \cite{Egg11} under the moniker of ``involution'' (which clashes with our usage).  


\subsection{Duals}
\label{Sec: dduals}
The reverse of an LDC, $\X$, written $\X^\rev := (\X, \ox, \top, \oa, \bot)^{\sf rev} = (\X, \ox^{\sf rev}, 
\top, \oa^{\sf rev}, \bot)$ where,
\[ A \ox^{\sf rev} B := B \ox A ~~~~~~~~~~~ A \oa^{\sf rev} B := B \oa A \]
and the associators and distributors are adjusted accordingly.  Similar to the opposite of an LDC, we have 
$(\X^{\sf rev})^{\sf rev} = \X$.  

In a $*$-autonomous category, taking the left (or right) linear dual of an object extends to a Frobenius linear 
functor as follows:
\begin{align*}
(\_)^*: (\X^{\sf op})^{\sf rev} &\to \X  ; ~~ A \mapsto A^* ; ~~
\begin{tikzpicture}
	\begin{pgfonlayer}{nodelayer}
		\node [style=none] (0) at (0, 2) {};
		\node [style=none] (1) at (-0.5, 1) {};
		\node [style=none] (2) at (0.5, 1) {};
		\node [style=none] (3) at (-0.5, -0) {};
		\node [style=none] (4) at (0.5, -0) {};
		\node [style=none] (5) at (0, 1) {};
		\node [style=none] (6) at (0, -1) {};
		\node [style=none] (7) at (0, -0) {};
		\node [style=none] (8) at (0, 0.5) {$f$};
		\node [style=none] (9) at (0.25, 1.75) {$A$};
		\node [style=none] (10) at (0.25, -0.75) {$B$};
	\end{pgfonlayer}
	\begin{pgfonlayer}{edgelayer}
		\draw (1.center) to (3.center);
		\draw (3.center) to (4.center);
		\draw (4.center) to (2.center);
		\draw (2.center) to (1.center);
		\draw (0.center) to (5.center);
		\draw (7.center) to (6.center);
	\end{pgfonlayer}
\end{tikzpicture} \mapsto  \begin{tikzpicture}
	\begin{pgfonlayer}{nodelayer}
		\node [style=none] (0) at (0, 1.5) {};
		\node [style=none] (1) at (-0.5, 1) {};
		\node [style=none] (2) at (0.5, 1) {};
		\node [style=none] (3) at (-0.5, -0) {};
		\node [style=none] (4) at (0.5, -0) {};
		\node [style=none] (5) at (0, 1) {};
		\node [style=none] (6) at (0, -0.5) {};
		\node [style=none] (7) at (0, -0) {};
		\node [style=none] (8) at (-1, 1.5) {};
		\node [style=none] (9) at (-1, -1) {};
		\node [style=none] (10) at (1, -0.5) {};
		\node [style=none] (11) at (1, 2) {};
		\node [style=none] (12) at (0, 0.5) {$f$};
		\node [style=none] (13) at (1.25, 1.75) {$B^*$};
		\node [style=none] (14) at (-1.25, -0.75) {$A^*$};
	\end{pgfonlayer}
	\begin{pgfonlayer}{edgelayer}
		\draw (1.center) to (3.center);
		\draw (3.center) to (4.center);
		\draw (4.center) to (2.center);
		\draw (2.center) to (1.center);
		\draw (0.center) to (5.center);
		\draw (7.center) to (6.center);
		\draw [bend right=90, looseness=1.25] (0.center) to (8.center);
		\draw (8.center) to (9.center);
		\draw [bend right=90, looseness=1.25] (6.center) to (10.center);
		\draw (10.center) to (11.center);
	\end{pgfonlayer}
\end{tikzpicture}
\end{align*}

The $(\_)^*$ functor is both contravariant and, op-monoidal and op-comonoidal:

\[ m_{\ox}:  A^* \ox B^* \to (B \oa A)^* :=  \begin{tikzpicture}
\begin{pgfonlayer}{nodelayer}
\node [style=none] (0) at (-0.25, 2) {};
\node [style=ox] (1) at (-0.25, 1) {};
\node [style=none] (2) at (-0.5, 0.25) {};
\node [style=none] (3) at (0.25, -0) {};
\node [style=none] (4) at (-1.25, 0.25) {};
\node [style=none] (5) at (-2, -0) {};
\node [style=oa] (6) at (-1.5, 1) {};
\node [style=none] (7) at (-1.5, 1.5) {};
\node [style=none] (8) at (-2.25, 1.5) {};
\node [style=none] (9) at (-2.25, -1) {};
\node [style=none] (10) at (0.5, 1.75) {$B^* \ox A^*$};
\node [style=none] (11) at (-0.8, 0.5) {$B^*$};
\node [style=none] (12) at (0.5, 0.5) {$A^*$};
\node [style=none] (13) at (-3.2, -0.75) {$(A \oa B)^* $};
\end{pgfonlayer}
\begin{pgfonlayer}{edgelayer}
\draw [in=90, out=-45, looseness=1.00] (1) to (3.center);
\draw [in=90, out=-135, looseness=1.00] (1) to (2.center);
\draw [bend left=90, looseness=1.25] (2.center) to (4.center);
\draw [in=-60, out=90, looseness=1.00] (4.center) to (6);
\draw (6) to (7.center);
\draw [bend right=90, looseness=1.25] (7.center) to (8.center);
\draw (8.center) to (9.center);
\draw [in=90, out=-135, looseness=0.75] (6) to (5.center);
\draw [in=-90, out=-90, looseness=1.00] (5.center) to (3.center);
\draw (1) to (0.center);
\end{pgfonlayer}
\end{tikzpicture} 
~~~~~~~~~~~~~~~ m_\top: \top \to \bot^* :=  \begin{tikzpicture} 
\begin{pgfonlayer}{nodelayer}
\node [style=circle] (0) at (0, 1.75) {$\top$};
\node [style=none] (1) at (0, 3) {};
\node [style=circle, scale=0.5] (2) at (-1, 1) {};
\node [style=circle] (3) at (-1, -0.25) {$\bot$};
\node [style=none] (4) at (-2, -1) {};
\node [style=none] (5) at (-2, 1.5) {};
\node [style=none] (6) at (-1, 1.5) {};
\node [style=none] (7) at (0.25, 2.75) {$\top$};
\node [style=none] (8) at (-2.5, -0.75) {$\bot^*$};
\end{pgfonlayer}
\begin{pgfonlayer}{edgelayer}
\draw (1.center) to (0);
\draw (3) to (6.center);
\draw [bend left=90, looseness=2.00] (5.center) to (6.center);
\draw (5.center) to (4.center);
\draw [in=0, out=-90, looseness=1.50, dotted] (0) to (2);
\end{pgfonlayer}
\end{tikzpicture} \]
\[ n_{\oa}:  (A \ox B)^* \to B^* \oa A^* :=  \begin{tikzpicture} 
\begin{pgfonlayer}{nodelayer}
\node [style=none] (0) at (2, 2) {};
\node [style=ox] (1) at (1, 0.25) {};
\node [style=none] (2) at (2, -0.5) {};
\node [style=none] (3) at (1, -0.5) {};
\node [style=none] (4) at (0.5, 0.75) {};
\node [style=none] (5) at (1.5, 0.75) {};
\node [style=oa] (6) at (-0.75, 0.25) {};
\node [style=none] (7) at (-0.75, -1.75) {};
\node [style=none] (8) at (-0.25, 0.75) {};
\node [style=none] (9) at (-1.25, 0.75) {};
\node [style=none] (10) at (2, 2.25) {$(A \ox B)^*$};
\node [style=none] (11) at (-0.75, -2) {$B^* \oa A^*$};
\end{pgfonlayer}
\begin{pgfonlayer}{edgelayer}
\draw (0.center) to (2.center);
\draw [bend left=90, looseness=1.50] (2.center) to (3.center);
\draw (3.center) to (1);
\draw (6) to (7.center);
\draw [in=150, out=-90, looseness=1.00] (9.center) to (6);
\draw [in=-90, out=11, looseness=1.00] (6) to (8.center);
\draw [in=165, out=-90, looseness=1.25] (4.center) to (1);
\draw [in=-90, out=15, looseness=1.25] (1) to (5.center);
\draw [bend left=90, looseness=2.00] (8.center) to (4.center);
\draw [bend left=90, looseness=1.25] (9.center) to (5.center);
\end{pgfonlayer}
\end{tikzpicture} ~~~~~~~~~~~~~~~ n_\bot: \top^* \to \bot :=  \begin{tikzpicture}
\begin{pgfonlayer}{nodelayer}
\node [style=circle] (0) at (-2, 1) {$\top$};
\node [style=none] (1) at (-1, 2) {};
\node [style=none] (2) at (-2, -0) {};
\node [style=none] (3) at (-1, -0) {};
\node [style=circle] (4) at (0, -1) {$\bot$};
\node [style=circle, scale=0.5] (5) at (-1, 1) {};
\node [style=none] (6) at (0, -2) {};
\node [style=none] (7) at (-0.6, 1.75) {$\top^*$};
\end{pgfonlayer}
\begin{pgfonlayer}{edgelayer}
\draw (0) to (2.center);
\draw [bend right=90, looseness=1.75] (2.center) to (3.center);
\draw (3.center) to (1.center);
\draw [in=90, out=-15, looseness=1.25, dotted] (5) to (4);
\draw (6.center) to (4);
\end{pgfonlayer}
\end{tikzpicture} \]

These maps are op-monoidal and op-comonoidal laxors, hence are isomorphisms, which satisfy the obvious coherences. 
Thus, $(\_)^*$ is a strong Frobenius linear functor. 

In the rest of the section, we will write $(\X^\op)^\rev$ as $\X^{\op\rev}$.

\begin{lemma} If $\X$ is an isomix category, then $(\_)^*: \X^{\op\rev} \to \X$ is an isomix functor. \end{lemma}
\begin{proof}
Because  $(\_)^*$ is a strong Frobenius functor, by Lemma \ref{Lemma: isomix functor}, it suffices  to prove that 
$(\_)^*$ preserves mix, i.e.,  $(\_)^*$  is a mix functor i.e., we need  to show that $n_\bot \!~\m~ m_\top : \top^* \to \bot^* =  
\m^* $. The proof is as follows:
\[ n_\bot \m m_\top = \begin{tikzpicture}
	\begin{pgfonlayer}{nodelayer}
		\node [style=circle] (0) at (-3, 2.75) {$\top$};
		\node [style=none] (1) at (-2, 2) {};
		\node [style=none] (2) at (-2, 3) {};
		\node [style=circle] (3) at (-1, 1.25) {$\bot$};
		\node [style=map] (4) at (-1, 0.5) {};
		\node [style=circle] (5) at (-1, -0.25) {$\top$};
		\node [style=circle] (6) at (-2, -2) {$\bot$};
		\node [style=none] (7) at (-3, -1.25) {};
		\node [style=none] (8) at (-3, -2.25) {};
		\node [style=circle, scale=0.5] (9) at (-2, 1.75) {};
		\node [style=circle, scale=0.5] (10) at (-2, -0.75) {};
		\node [style=none] (11) at (-3, 2) {};
		\node [style=none] (12) at (-2, -1.25) {};
	\end{pgfonlayer}
	\begin{pgfonlayer}{edgelayer}
		\draw (1.center) to (2.center);
		\draw (7.center) to (8.center);
		\draw [bend left, looseness=1.00, dotted] (9) to (3);
		\draw [bend left=45, looseness=1.00, dotted] (5) to (10);
		\draw (3) to (4);
		\draw (4) to (5);
		\draw [bend right=90, looseness=3.25] (11.center) to (1.center);
		\draw (0) to (11.center);
		\draw [bend left=90, looseness=3.75] (7.center) to (12.center);
		\draw (12.center) to (6);
	\end{pgfonlayer}
\end{tikzpicture} =\begin{tikzpicture}
	\begin{pgfonlayer}{nodelayer}
		\node [style=circle] (0) at (1.25, 3) {$\top$};
		\node [style=none] (1) at (0.25, 3) {};
		\node [style=circle] (2) at (-1, 1.25) {$\bot$};
		\node [style=map] (3) at (-1, 0.5) {};
		\node [style=circle] (4) at (-1, -0.25) {$\top$};
		\node [style=circle] (5) at (-2, -2.25) {$\bot$};
		\node [style=none] (6) at (-3, -2.25) {};
		\node [style=circle, scale=0.5] (7) at (0.25, 2) {};
		\node [style=circle, scale=0.5] (8) at (-2, -0.75) {};
		\node [style=none] (9) at (0.25, 1.5) {};
		\node [style=none] (10) at (1.25, 1.5) {};
		\node [style=none] (11) at (-2, -0.25) {};
		\node [style=none] (12) at (-3, -0.25) {};
	\end{pgfonlayer}
	\begin{pgfonlayer}{edgelayer}
		\draw [dotted, bend right, looseness=1.25] (7) to (2);
		\draw [dotted, bend left=45, looseness=1.00] (4) to (8);
		\draw (2) to (3);
		\draw (3) to (4);
		\draw (1.center) to (9.center);
		\draw [bend right=90, looseness=2.00] (9.center) to (10.center);
		\draw (10.center) to (0);
		\draw (5) to (11.center);
		\draw (12.center) to (6.center);
		\draw [bend left=90, looseness=1.75] (12.center) to (11.center);
	\end{pgfonlayer}
\end{tikzpicture} = \begin{tikzpicture} 
	\begin{pgfonlayer}{nodelayer}
		\node [style=circle] (0) at (1.25, 3) {$\top$};
		\node [style=none] (1) at (0.25, 2.75) {};
		\node [style=circle] (2) at (-1, 1.25) {$\bot$};
		\node [style=map] (3) at (-1, 0.5) {};
		\node [style=circle] (4) at (-1, -0.25) {$\top$};
		\node [style=circle] (5) at (-2, -2.25) {$\bot$};
		\node [style=none] (6) at (-3, -2) {};
		\node [style=circle, scale=0.5] (7) at (0.25, 2) {};
		\node [style=circle, scale=0.5] (8) at (-2, -0.75) {};
		\node [style=none] (9) at (0.25, -1.75) {};
		\node [style=none] (10) at (1.25, -1.75) {};
		\node [style=none] (11) at (-2, 2.5) {};
		\node [style=none] (12) at (-3, 2.5) {};
	\end{pgfonlayer}
	\begin{pgfonlayer}{edgelayer}
		\draw [dotted, bend right, looseness=1.25] (7) to (2);
		\draw [dotted, bend left=45, looseness=1.00] (4) to (8);
		\draw (2) to (3);
		\draw (3) to (4);
		\draw (1.center) to (9.center);
		\draw [bend right=90, looseness=2.00] (9.center) to (10.center);
		\draw (10.center) to (0);
		\draw (5) to (11.center);
		\draw (12.center) to (6.center);
		\draw [bend left=90, looseness=1.75] (12.center) to (11.center);
	\end{pgfonlayer}
\end{tikzpicture} \stackrel{\mx}{=} \begin{tikzpicture} 
	\begin{pgfonlayer}{nodelayer}
		\node [style=circle] (0) at (1.25, 3) {$\top$};
		\node [style=none] (1) at (0.25, 2.75) {};
		\node [style=circle] (2) at (-0.75, 1.25) {$\bot$};
		\node [style=map] (3) at (-0.75, 0.5) {};
		\node [style=circle] (4) at (-0.75, -0.25) {$\top$};
		\node [style=circle] (5) at (-2, -2.25) {$\bot$};
		\node [style=none] (6) at (-3, -2) {};
		\node [style=circle, scale=0.5] (7) at (-2, 2) {};
		\node [style=circle, scale=0.5] (8) at (0.25, -0.75) {};
		\node [style=none] (9) at (0.25, -1.75) {};
		\node [style=none] (10) at (1.25, -1.75) {};
		\node [style=none] (11) at (-2, 2.5) {};
		\node [style=none] (12) at (-3, 2.5) {};
	\end{pgfonlayer}
	\begin{pgfonlayer}{edgelayer}
		\draw [dotted, bend left, looseness=1.25] (7) to (2);
		\draw [dotted, bend right=45, looseness=1.00] (4) to (8);
		\draw (2) to (3);
		\draw (3) to (4);
		\draw (1.center) to (9.center);
		\draw [bend right=90, looseness=2.00] (9.center) to (10.center);
		\draw (10.center) to (0);
		\draw (5) to (11.center);
		\draw (12.center) to (6.center);
		\draw [bend left=90, looseness=1.75] (12.center) to (11.center);
	\end{pgfonlayer}
\end{tikzpicture} = \begin{tikzpicture}
	\begin{pgfonlayer}{nodelayer}
		\node [style=map] (0) at (-0.75, 1) {};
		\node [style=none] (1) at (-1.75, -1.75) {};
		\node [style=none] (2) at (-0.75, -0.25) {};
		\node [style=none] (3) at (0.25, -0.25) {};
		\node [style=none] (4) at (-0.75, 2.25) {};
		\node [style=none] (5) at (-1.75, 2.25) {};
		\node [style=none] (6) at (0.25, 3.25) {};
	\end{pgfonlayer}
	\begin{pgfonlayer}{edgelayer}
		\draw [bend right=90, looseness=2.00] (2.center) to (3.center);
		\draw (5.center) to (1.center);
		\draw [bend left=90, looseness=1.75] (5.center) to (4.center);
		\draw (4.center) to (0);
		\draw (0) to (2.center);
		\draw (6.center) to (3.center);
	\end{pgfonlayer}
\end{tikzpicture} = \m^*
\]
\end{proof}

\begin{lemma}
$(\eta, \epsilon) :: (\_)^* \dashvv ~ {^*(\_)^{\op\rev}} : \X^{\op \rev} \to \X$  
\[ 
\eta_\ox: X \to ~^*(X^*) := \begin{tikzpicture} 
	\begin{pgfonlayer}{nodelayer}
		\node [style=none] (0) at (1, -2) {};
		\node [style=none] (1) at (1, 0.75) {};
		\node [style=none] (2) at (0, 0.75) {};
		\node [style=none] (12) at (0.5, 1.75) {$*\eta$};		
		\node [style=none] (3) at (0, -0.25) {};
		\node [style=none] (4) at (-1, -0.25) {};
		\node [style=none] (12) at (-0.5, -1.5) {$\epsilon*$};		
		\node [style=none] (5) at (-1, 2.5) {};
		\node [style=none] (6) at (-1.25, 2) {$X$};
		\node [style=none] (7) at (-0.35, 0.5) {$X^*$};
		\node [style=none] (8) at (1.6, -1.75) {$^*(X^*)$};
	\end{pgfonlayer}
	\begin{pgfonlayer}{edgelayer}
		\draw (4.center) to (5.center);
		\draw [bend right=90, looseness=3.25] (4.center) to (3.center);
		\draw (3.center) to (2.center);
		\draw [bend left=90, looseness=2.50] (2.center) to (1.center);
		\draw (1.center) to (0.center);
	\end{pgfonlayer}
\end{tikzpicture} \in \X
~~~~~~~~~~~~ \eta_\oa := \eta_\ox^{-1} \]
\[
\epsilon_\oa: X \to  (^*X)^*:= \begin{tikzpicture}
	\begin{pgfonlayer}{nodelayer}
		\node [style=none] (0) at (-0.75, -2) {};
		\node [style=none] (1) at (-0.75, 0.75) {};
		\node [style=none] (2) at (0.25, 0.75) {};
		\node [style=none] (12) at (-0.5, 1.75) {$\eta*$};
		\node [style=none] (3) at (0.25, -0.25) {};
		\node [style=none] (4) at (1.25, -0.25) {};
		\node [style=none] (34) at (0.65, -1.5) {$*\epsilon$};
		\node [style=none] (5) at (1.25, 2.5) {};
		\node [style=none] (6) at (1.65, 2) {$X$};
		\node [style=none] (7) at (0.6, 0.5) {$^*X$};
		\node [style=none] (8) at (-1.25, -1.75) {$(^*X)^*$};
	\end{pgfonlayer}
	\begin{pgfonlayer}{edgelayer}
		\draw (4.center) to (5.center);
		\draw [bend left=90, looseness=3.25] (4.center) to (3.center);
		\draw (3.center) to (2.center);
		\draw [bend right=90, looseness=2.50] (2.center) to (1.center);
		\draw (1.center) to (0.center);
	\end{pgfonlayer}
\end{tikzpicture} \in \X ~~~~~~~~~~~ \epsilon_\ox := \epsilon_\oa^{-1}
\]
is a linear equivalence of Frobenius linear functors.
\end{lemma}
\begin{proof}
The proof is straightforward in the graphical calculus.
\end{proof}

For a cyclic $*$-autonomous category, we can straighten out this equivalence to be a dualizing involutive equivalence 
(i.e. so that the unit and counit are equal):

\begin{lemma}
$(\eta', \epsilon') :: (\_)^* \dashvv  ((\_)^{*})^{\op\rev}: \X^{\op\rev} \to \X$  where $\eta'_\ox = {\eta'_\oa}^{-1} := 
\eta_\ox \psi^{-1}$, $\epsilon'_\ox = \epsilon'_\oa:= \epsilon \psi^*$ and $\eta' = \epsilon'$.
\end{lemma}
\begin{proof}

The unit and counit are drawn as follows:

\[
\eta_\ox' =  \begin{tikzpicture}
	\begin{pgfonlayer}{nodelayer}
		\node [style=none] (0) at (-3, 3) {};
		\node [style=none] (1) at (-3, 1) {};
		\node [style=none] (2) at (-2, 1) {};
		\node [style=none] (3) at (-2, 2) {};
		\node [style=none] (4) at (-1, 2) {};
		\node [style=circle] (5) at (-1, -0) {$\psi^{-1}$};
		\node [style=none] (6) at (-1, -1) {};
		\node [style=none] (7) at (-2.5, 0.25) {$\epsilon*$};
		\node [style=none] (8) at (-1.5, 2.75) {$~^*\eta$};
		\node [style=none] (9) at (-3.25, 2.75) {$X$};
		\node [style=none] (10) at (-0.5, 1.75) {$~^*(X^*)$};
		\node [style=none] (11) at (-0.5, -0.7) {$X^{**}$};
		\node [style=none] (12) at (-2.35, 1.5) {$X^*$};
	\end{pgfonlayer}
	\begin{pgfonlayer}{edgelayer}
		\draw (0.center) to (1.center);
		\draw [bend right=90, looseness=1.75] (1.center) to (2.center);
		\draw (2.center) to (3.center);
		\draw [bend left=90, looseness=2.00] (3.center) to (4.center);
		\draw (4.center) to (5);
		\draw (5) to (6.center);
	\end{pgfonlayer}
\end{tikzpicture} \in \X ~~~~~~~~~~ \epsilon_\ox' = \begin{tikzpicture}
	\begin{pgfonlayer}{nodelayer}
		\node [style=none] (0) at (-0.75, -0.5) {};
		\node [style=none] (1) at (-0.75, 2.25) {};
		\node [style=none] (2) at (0.25, 2.25) {};
		\node [style=none] (3) at (0.25, 1) {};
		\node [style=none] (4) at (1.25, 1) {};
		\node [style=none] (5) at (1.25, 3) {};
		\node [style=none] (6) at (0.75, 0.25) {$*\epsilon$};
		\node [style=none] (7) at (-0.25, 3) {$\eta*$};
		\node [style=none] (8) at (1.5, 2.5) {$X^*$};
		\node [style=none] (9) at (-1.75, -0.5) {};
		\node [style=circle, scale=2] (10) at (-1.75, 0.5) {};
		\node [style=none] (11) at (-1.75, 1.25) {};
		\node [style=none] (12) at (-2.75, 1.25) {};
		\node [style=none] (13) at (-2.75, -1.25) {};
		\node [style=none] (14) at (-1.25, -1.25) {$\epsilon*$};
		\node [style=none] (15) at (-2.25, 2) {$\eta*$};
		\node [style=none] (16) at (-1.75, 0.5) {$\psi$};
		\node [style=none] (17) at (-3.25, -0.75) {$X^{**}$};
		\node [style=none] (18) at (-0.25, -0) {$~^*(X^*)$};
	\end{pgfonlayer}
	\begin{pgfonlayer}{edgelayer}
		\draw [bend left=90, looseness=2.00] (1.center) to (2.center);
		\draw (2.center) to (3.center);
		\draw [bend right=90, looseness=2.00] (3.center) to (4.center);
		\draw (4.center) to (5.center);
		\draw (1.center) to (0.center);
		\draw [bend right=90, looseness=1.75] (9.center) to (0.center);
		\draw (11.center) to (10);
		\draw (10) to (9.center);
		\draw [bend right=90, looseness=1.75] (11.center) to (12.center);
		\draw (12.center) to (13.center);
	\end{pgfonlayer}
\end{tikzpicture} = 
 \begin{tikzpicture}
	\begin{pgfonlayer}{nodelayer}
		\node [style=none] (0) at (-0.75, 3) {};
		\node [style=none] (1) at (-0.75, -0) {};
		\node [style=none] (2) at (-1.75, -0) {};
		\node [style=none] (3) at (-1.75, 2) {};
		\node [style=none] (4) at (-2.75, 2) {};
		\node [style=none] (5) at (-2.75, -1) {};
		\node [style=none] (6) at (-1.25, -0.75) {$*\epsilon$};
		\node [style=none] (7) at (-2.25, 2.75) {$\eta^*$};
		\node [style=none] (8) at (-0.5, 2.75) {$X$};
		\node [style=none] (9) at (-3.25, -0.5) {$X^{**}$};
		\node [style=circle] (10) at (-1.75, 1) {$\psi$};
		\node [style=none] (11) at (-2.25, 0.25) {$~^*X$};
		\node [style=none] (12) at (-2.1, 1.7) {$X^*$};
	\end{pgfonlayer}
	\begin{pgfonlayer}{edgelayer}
		\draw (0.center) to (1.center);
		\draw [bend left=90, looseness=1.75] (1.center) to (2.center);
		\draw [bend right=90, looseness=2.00] (3.center) to (4.center);
		\draw (4.center) to (5.center);
		\draw (3.center) to (10);
		\draw (10) to (2.center);
	\end{pgfonlayer}
\end{tikzpicture} \stackrel{{\bf \tiny [C.2]}}{=} 
 \begin{tikzpicture}
	\begin{pgfonlayer}{nodelayer}
		\node [style=none] (0) at (-3, 3) {};
		\node [style=none] (1) at (-3, 1) {};
		\node [style=none] (2) at (-2, 1) {};
		\node [style=none] (3) at (-2, 2) {};
		\node [style=none] (4) at (-1, 2) {};
		\node [style=circle, scale=2.3] (5) at (-1, -0) {};
		\node [style=none] (13) at (-1, -0) {$\psi^{-1}$};
		\node [style=none] (6) at (-1, -1) {};
		\node [style=none] (7) at (-2.5, 0.25) {$\epsilon*$};
		\node [style=none] (8) at (-1.5, 2.75) {$~^*\eta$};
		\node [style=none] (9) at (-3.25, 2.75) {$X$};
		\node [style=none] (10) at (-0.5, 1.75) {$~^*(X^*)$};
		\node [style=none] (11) at (-0.5, -0.7) {$X^{**}$};
		\node [style=none] (12) at (-2.35, 1.5) {$X^*$};
	\end{pgfonlayer}
	\begin{pgfonlayer}{edgelayer}
		\draw (0.center) to (1.center);
		\draw [bend right=90, looseness=1.75] (1.center) to (2.center);
		\draw (2.center) to (3.center);
		\draw [bend left=90, looseness=2.00] (3.center) to (4.center);
		\draw (4.center) to (5);
		\draw (5) to (6.center);
	\end{pgfonlayer}
\end{tikzpicture} \in \X
\]

The cyclor is a linear transformation which is an isomorphism as it is monoidal with respect to both tensor 
and par and adjoints are determined only upto isomorphism. It remains to check that the triangle identities hold:
\[
\eta_{\ox'_{X^*}}(\epsilon_\ox')^* = 
\begin{tikzpicture} 
	\begin{pgfonlayer}{nodelayer}
		\node [style=circle] (0) at (0.75, 0.25) {$\psi^{-1}$};
		\node [style=none] (1) at (0.75, -2.5) {};
		\node [style=none] (2) at (0.75, 2.25) {};
		\node [style=none] (3) at (-0.25, 2.25) {};
		\node [style=none] (4) at (-0.25, 1) {};
		\node [style=none] (5) at (-1.25, 1) {};
		\node [style=none] (6) at (-1.25, 3) {};
		\node [style=none] (7) at (1.5, -0.5) {$X^{***}$};
		\node [style=none] (8) at (-0.75, 0.25) {$\epsilon*$};
		\node [style=none] (9) at (0.25, 3) {$*\eta$};
		\node [style=none] (10) at (-1.5, 2.5) {$X^*$};
		\node [style=none] (11) at (-1, -2.5) {};
		\node [style=circle] (12) at (-1, -1.5) {$\psi^{-1}$};
		\node [style=none] (13) at (-1, -0.75) {};
		\node [style=none] (14) at (-2, -0.75) {};
		\node [style=none] (15) at (-2, -3.25) {};
		\node [style=none] (16) at (-1.5, -0.1) {$*\eta$};
		\node [style=none] (17) at (0, -3.35) {$\epsilon*$};
		\node [style=none] (18) at (-2.5, -3) {};
	\end{pgfonlayer}
	\begin{pgfonlayer}{edgelayer}
		\draw (2.center) to (0);
		\draw (0) to (1.center);
		\draw [bend right=90, looseness=2.00] (2.center) to (3.center);
		\draw (3.center) to (4.center);
		\draw [bend left=90, looseness=2.00] (4.center) to (5.center);
		\draw (5.center) to (6.center);
		\draw [bend right=75, looseness=1.25] (11.center) to (1.center);
		\draw (11.center) to (12);
		\draw (12) to (13.center);
		\draw [bend right=90, looseness=1.50] (13.center) to (14.center);
		\draw (14.center) to (15.center);
	\end{pgfonlayer}
\end{tikzpicture} = \begin{tikzpicture} 
	\begin{pgfonlayer}{nodelayer}
		\node [style=none] (0) at (0, -2.5) {};
		\node [style=none] (1) at (0, 2.5) {};
		\node [style=none] (2) at (1, 2.5) {};
		\node [style=none] (3) at (1, 1.25) {};
		\node [style=none] (4) at (2, 1.25) {};
		\node [style=none] (5) at (2, 3.25) {};
		\node [style=none] (6) at (-0.5, 1.25) {$X^{***}$};
		\node [style=none] (7) at (1.5, 0.5) {$*\epsilon$};
		\node [style=none] (8) at (0.5, 3.25) {$\eta*$};
		\node [style=none] (9) at (2.25, 2.75) {$X^*$};
		\node [style=none] (10) at (-1, -2.5) {};
		\node [style=circle] (11) at (-1, -1.5) {$\psi^{-1}$};
		\node [style=none] (12) at (-1, -0.75) {};
		\node [style=none] (13) at (-2, -0.75) {};
		\node [style=none] (14) at (-2, -3.25) {};
		\node [style=none] (15) at (-1.5, -0) {$*\eta$};
		\node [style=none] (16) at (-0.5, -3) {$\epsilon*$};
		\node [style=circle] (17) at (1, 2) {$\psi$};
		\node [style=none] (18) at (-2.25, -2.75) {$X^*$};
	\end{pgfonlayer}
	\begin{pgfonlayer}{edgelayer}
		\draw [bend left=90, looseness=2.00] (1.center) to (2.center);
		\draw [bend right=90, looseness=2.00] (3.center) to (4.center);
		\draw (4.center) to (5.center);
		\draw [bend right=75, looseness=1.25] (10.center) to (0.center);
		\draw (10.center) to (11);
		\draw (11) to (12.center);
		\draw [bend right=90, looseness=1.50] (12.center) to (13.center);
		\draw (13.center) to (14.center);
		\draw (1.center) to (0.center);
		\draw (2.center) to (17);
		\draw (17) to (3.center);
	\end{pgfonlayer}
\end{tikzpicture} = \begin{tikzpicture}
	\begin{pgfonlayer}{nodelayer}
		\node [style=none] (0) at (-1, -1) {};
		\node [style=none] (1) at (0, -1) {};
		\node [style=none] (2) at (0, 3.25) {};
		\node [style=none] (3) at (-0.5, -1.75) {$*\epsilon$};
		\node [style=none] (4) at (0.25, 2.75) {$X^*$};
		\node [style=circle] (5) at (-1, 1.5) {$\psi^{-1}$};
		\node [style=none] (6) at (-1, 2.25) {};
		\node [style=none] (7) at (-2, 2.25) {};
		\node [style=none] (8) at (-2, -3.25) {};
		\node [style=none] (9) at (-1.5, 3) {$*\eta$};
		\node [style=circle] (10) at (-1, -0.25) {$\psi$};
		\node [style=none] (11) at (-2.25, -2.75) {$X^*$};
	\end{pgfonlayer}
	\begin{pgfonlayer}{edgelayer}
		\draw [bend right=90, looseness=2.00] (0.center) to (1.center);
		\draw (1.center) to (2.center);
		\draw (5) to (6.center);
		\draw [bend right=90, looseness=1.50] (6.center) to (7.center);
		\draw (7.center) to (8.center);
		\draw (10) to (0.center);
		\draw (5) to (10);
	\end{pgfonlayer}
\end{tikzpicture} = ~ \begin{tikzpicture} \draw (2, 3.25) -- (2,-3.5); \end{tikzpicture} ~ = ~ 1
\] The other triangle identity holds similarly.
\end{proof}

The equality of $\eta'$ and $\epsilon'$ is immediate from {\bf [C.2]} for cyclors with the 
map $\eta'=\epsilon'$ being the {\bf dualizor}.   In the symmetric case, the dualizor of this 
equivalence may be drawn as:
 \[\begin{tikzpicture}
	\begin{pgfonlayer}{nodelayer}
		\node [style=none] (0) at (-3, 3) {};
		\node [style=none] (1) at (-3, 1) {};
		\node [style=none] (2) at (-2, 1) {};
		\node [style=none] (3) at (-2, 2) {};
		\node [style=none] (4) at (-1, 2) {};
		\node [style=circle, scale=2] (5) at (-1, -0) {};
		\node [style=none] (5) at (-1, -0) {$\psi^{-1}$};
		\node [style=none] (6) at (-1, -1) {};
		\node [style=none] (7) at (-2.5, 0.25) {$\epsilon*$};
		\node [style=none] (8) at (-1.5, 2.75) {$*\eta$};
		\node [style=none] (9) at (-3.25, 2.75) {$A$};
		\node [style=none] (10) at (-0.5, 1.75) {$~^*(A^*)$};
		\node [style=none] (11) at (-0.5, -0.7) {$A^{**}$};
		\node [style=none] (12) at (-2.25, 1.5) {$A^*$};
	\end{pgfonlayer}
	\begin{pgfonlayer}{edgelayer}
		\draw (0.center) to (1.center);
		\draw [bend right=90, looseness=1.75] (1.center) to (2.center);
		\draw (2.center) to (3.center);
		\draw [bend left=90, looseness=2.00] (3.center) to (4.center);
		\draw (4.center) to (5);
		\draw (5) to (6.center);
	\end{pgfonlayer}
\end{tikzpicture} =
\begin{tikzpicture}
\begin{pgfonlayer}{nodelayer}
\node [style=none] (0) at (-1.75, 2) {};
\node [style=none] (1) at (-1.75, -0) {};
\node [style=none] (2) at (0, -0) {};
\node [style=none] (12) at (-0.75, -1) {$\epsilon*$};
\node [style=none] (3) at (0, 1) {};
\node [style=none] (4) at (-1, 1) {};
\node [style=none] (34) at (-0.5, 1.75) {$\eta*$};
\node [style=none] (5) at (-1, -2.75) {};
\node [style=none] (6) at (-2, 1.75) {$A$};
\node [style=none] (8) at (-1.5, -2.25) {$A^{**}$};
\end{pgfonlayer}
\begin{pgfonlayer}{edgelayer}
\draw (0.center) to (1.center);
\draw [bend right=90, looseness=1.50] (1.center) to (2.center);
\draw (2.center) to (3.center);
\draw [bend right=90, looseness=1.75] (3.center) to (4.center);
\draw (4.center) to (5.center);
\end{pgfonlayer}
\end{tikzpicture} =  \begin{tikzpicture}
	\begin{pgfonlayer}{nodelayer}
		\node [style=none] (0) at (-0.75, 3) {};
		\node [style=none] (1) at (-0.75, -0) {};
		\node [style=none] (2) at (-1.75, -0) {};
		\node [style=none] (3) at (-1.75, 2) {};
		\node [style=none] (4) at (-2.75, 2) {};
		\node [style=none] (5) at (-2.75, -1) {};
		\node [style=none] (6) at (-1.25, -0.75) {$*\epsilon$};
		\node [style=none] (7) at (-2.25, 2.75) {$\eta*$};
		\node [style=none] (8) at (-0.5, 2.75) {$A$};
		\node [style=none] (9) at (-3.25, -0.5) {$A^{**}$};
		\node [style=circle] (10) at (-1.75, 1) {$\psi$};
		\node [style=none] (11) at (-2.25, 0.25) {$~^*A$};
		\node [style=none] (12) at (-2, 1.7) {$A^*$};
	\end{pgfonlayer}
	\begin{pgfonlayer}{edgelayer}
		\draw (0.center) to (1.center);
		\draw [bend left=90, looseness=1.75] (1.center) to (2.center);
		\draw [bend right=90, looseness=2.00] (3.center) to (4.center);
		\draw (4.center) to (5.center);
		\draw (3.center) to (10);
		\draw (10) to (2.center);
	\end{pgfonlayer}
\end{tikzpicture} \]

\subsection{Conjugation}
\label{Sec: conjugation}

Recall the following structure from Egger \cite{Egg11}:

\begin{definition}
A {\bf conjugation} for a monoidal category $(X, \otimes, I)$ consists of a functor $\overline{(\_)}: \X^\rev \to \X$ 
with natural isomorphisms:
\[ \bar{A} \otimes \bar{B} \to^{\chi} \bar{B \otimes A} ~~~~~~~~~~~~~ \bar{\bar{A}} \to^{\varepsilon} A  \]
called respectively the (tensor reversing) {\bf conjugating laxor} and the {\bf conjugator}
such that 
\[\bar{\bar{\bar{A}}} \to^{\bar{\varepsilon_A} = \varepsilon_{\bar{A}} } \bar{A} \] 
and  
\[
\xymatrix{
(\bar{A} \ox \bar{B}) \ox \bar{C} \ar[rr]^{a_\ox} \ar[d]_{\chi \ox 1}  \ar@{}[ddrr]|{\bf [CF.1]_\ox}  
& & \bar{A} \ox (\bar{B} \ox \bar{C}) \ar[d]^{1 \ox \chi} \\
\bar{(B \ox A)} \ox \bar{C} \ar[d]_{\chi} & & \bar{A} \ox \bar{(C \ox B)} \ar[d]^{\chi} \\
\bar{C \ox (B \ox A)} \ar[rr]_{\bar{a_\otimes^{-1}}} & & \bar{(C \ox B) \ox A}
} ~~~~~~~~ \xymatrix{
\bar{\bar{A}} \ox \bar{\bar{B}} \ar[rr]^{\chi} \ar[dd]_{\varepsilon \ox \varepsilon}  \ar@{}[ddrr]|
{\bf [CF.2]_\ox}  & & \bar{\bar{B} \ox \bar{A}} \ar[dd]^{\bar{\chi}} \\ 
& & \\
A \ox B & & \bar{\bar{A \ox B}} \ar[ll]^{\varepsilon}
}
\]
\end{definition}

A monoidal category is {\bf conjugative} when it has a conjugation functor.

A symmetric monoidal category, which is conjugative, is {\bf symmetric conjugative} in case it 
satisfies the additional coherence:
\[ 
\xymatrix{
\bar{A} \otimes \bar{B} \ar[d]_{c_\otimes} \ar[rr]^{\chi}  \ar@{}[drr]|{\bf [CF.3]_\ox}  &  & 
\bar{B \otimes A} \ar[d]^{\bar{c_\otimes}} \\
\bar{B} \otimes \bar{A} \ar[rr]_{\chi} &  & \bar{A \otimes B} 
}
\]

No coherences have been specified for the unit $I$ because the expected coherences are automatic:

\begin{lemma}
\label{Lemma: unit conjugate} \cite[Lemma 2.3]{Egg11}
For every conjugative monoidal category, there exists a unique isomorphism 
$I \xrightarrow{\chi^{\!\!\!\circ}} \bar{I}$ such that 
\[
\xymatrix{
I \otimes \bar{A} \ar[rr]^{\chi^{\!\!\!\circ} \otimes 1} \ar[d]_{u_\otimes}  \ar@{}[drr]|{\bf [CF.4]_\top}  
& & \bar{I} \otimes \bar{A} \ar[d]^{\chi} \\
\bar{A} \ar[rr]_{\bar{u_\otimes^{-1}}} & & \bar{A \otimes I } } ~~~~~ 
\xymatrix{ \bar{A} \ox I \ar[rr]^{1 \otimes \chi^{\!\!\!\circ}} \ar[d]_{u_\otimes}  \ar@{}[drr]|{\bf [CF.5]_\top}  & & \bar{A} 
\ox \bar{I} \ar[d]^{\chi} \\
\bar{A} \ar[rr]_{\bar{u_\otimes^{-1}}} & & \bar{I \otimes A} } ~~~~~
\xymatrix{
I \ar[rr]^{\chi^{\!\!\!\circ}} \ar@{=}[d]  \ar@{}[drr]|{\bf [CF.6]_\top}  & & \bar{I} \ar[d]^{\bar{\chi^{\!\!\!\circ}}} \\
I \ar[rr]_{\varepsilon^{-1}} & & \bar{\bar{I}}
}
\]
\end{lemma}

\begin{definition} \cite{Egg11}
A {\bf conjugative LDC} is a linearly distributive category $(\X, \ox, \top, \oa, \bot)$ together with a 
conjugating functor $\bar{(\_)}: \X \to \X$ and natural isomorphisms:
\[ \bar{A} \ox \bar{B} \xrightarrow{\chi_\ox} \bar{B \ox A} ~~~~~~~~~~~~ \bar{A \oa B} 
\xrightarrow{\chi_\oa} \bar{B} \oa \bar{A} ~~~~~~~~~~~~ \bar{\bar{A}} \xrightarrow{\varepsilon} A \] 
\end{definition}

such that $(\X, \ox, \top, \chi_\ox, \varepsilon)$ and $(\X, \oa, \bot, \chi_\oa^{-1}, \varepsilon)$ are 
conjugative (symmetric) monoidal categories with respect to the conjugating functor and  the following diagrams commute:

\[ 
\xymatrix{
\bar{B \oa C} \ox \bar{A} \ar[rr]^{\chi_\oa \ox 1} \ar[d]_{\chi_\ox}  \ar@{}[ddrr]|{\bf [CF.7]}  & & 
(\bar{C} \oa \bar{B}) \ox \bar{A} \ar[d]^{\partial} \\
\bar{(A \ox (B \oa C))} \ar[d]_{\bar{\partial}} & & \bar{C} \oa ( \bar{B} \ox \bar{A}) \ar[d]^{1 \oa \chi_\ox} \\
\bar{((A \ox B) \oa C)} \ar[rr]_{\chi_\oa} & & \bar{C} \oa \bar{A \ox B}
} ~~~~~~~~~~ 
\xymatrix{
\bar{A} \ox \bar{C \oa B} \ar[rr]^{\chi_\ox} \ar[d]_{1 \ox \chi_\oa} \ar@{}[ddrr]|{\bf [CF.8]}   & & 
\bar{(C \oa B) \ox A} \ar[d]^{\bar{\partial}} \\
\bar{A} \ox (\bar{B} \oa \bar{C}) \ar[d]_{\partial} & & \bar{C \oa (B \ox A)} \ar[d]^{\chi_\oa} \\
(\bar{A} \ox \bar{B}) \oa \bar{C} \ar[rr]_{\chi_\ox \oa 1} & & \bar{(B \ox A)} \oa \bar{C}
} 
\]

Note, by Lemma \ref{Lemma: unit conjugate}, there exists canonical isomorphisms 
$\top \xrightarrow{\chi^{\!\!\!\circ}_\top} \bar{\top}$ and $\bot \xrightarrow{\chi^{\!\!\!\circ}_\bot} \bar{\bot}$, 
hence conjugation is a normal functor.  However, the conjugation is not necessarily a mix functor when $\X$ is a mix category.  
 For conjugation to be a mix functor, the following extra condition must be satisfied:
\[ {\bf \small [CF.9]}~~~~~ \xymatrix{\overline{\bot} \ar[rd]_{(\chi^{\!\!\!\circ}_\bot)^{-1}} \ar@/^/[rrr]^{\overline{{\sf m}}} 
& & & \overline{\top} \\ & \bot \ar[r]_{{\sf m}} & \top \ar[ur]_{\chi^{\!\!\!\circ}_\top} } \]

\begin{proposition}
A conjugative LDC is precisely a LDC, $\X$, with a Frobenius adjoint $(\epsilon^{-1}, \epsilon): 
\overline{(\_)} \dashv \overline{(\_)}^\rev: \X^\rev \to \X$ where $\epsilon := (\varepsilon,\varepsilon^{-1})$.  
Furthermore, if $\X$ is an isomix category and conjugation is a mix functor then conjugation is an isomix equivalence.
\end{proposition}

\begin{proof}
It is clear that $\overline{(\_)}$ is a strong Frobenius functor so being mix implies isomix.   
Also, $\varepsilon$ is clearly monoidal for tensor and par.  The triangle equalities give $ \overline{\varepsilon^{-1}} 
\varepsilon = 1: \overline{A} \to \overline{A}$ thus $\varepsilon = \overline{\varepsilon}$.
\end{proof}

Clearly conjugation should flip left duals into right duals:

\begin{lemma}
\label{Lemma: involutive linear adjoint}
If $B \dashvv A$ is a linear dual then $\bar{A} \dashvv \bar{B}$ is a linear dual.
\end{lemma}

\begin{proof}
Suppose $(\eta, \varepsilon): B \dashvv A$. Then, $(\chi^{\!\!\!\circ}_\top \bar{\eta} \chi_\oa,\chi_\ox \bar{\varepsilon} 
\chi^{\!\!\!\circ}_\bot): \bar{A} \dashvv \bar{B}$.
\end{proof}

When a $*$-autonomous category is cyclic one expects that conjugation should interact with the cyclor in a coherent fashion:

\begin{definition}\cite{EggMcCurd12}
\label{Defn: conjugative cyclic}
A {\bf conjugative cyclic $*$-autonomous category} is a conjugative $*$-autonomous category together 
with a cyclor $A^* \to^{\psi} \!\!~^{*}\!A$ such that 
\[
\xymatrix{
(\bar{A})^* \ar[rr]^{\psi} \ar[d]_{\simeq}  & & ^{*}(\bar{A}) \ar[d]^{\simeq} \\
\bar{(^{*}A)} \ar[rr]_{\bar{\psi^{-1}}} & & \bar{(A^*)}
}
\]
which gives a map $\sigma: (\overline{A})^{*} \to \overline{(A^{*})}$.
\end{definition}

The above condition is drawn as follows:
\[ \sigma =
\begin{tikzpicture} 
	\begin{pgfonlayer}{nodelayer}
		\node [style=none] (0) at (-1, 0.75) {$\overline{\psi^{-1}}$};
		\node [style=circle, scale=2.5] (0) at (-1, 0.75) {};
		\node [style=none] (1) at (-1, -0.5) {};
		\node [style=none] (2) at (-1, 2.25) {};
		\node [style=none] (3) at (0, 2.25) {};
		\node [style=none] (4) at (0, 1) {};
		\node [style=none] (5) at (1, 1) {};
		\node [style=none] (6) at (1, 3) {};
		\node [style=none] (7) at (1.5, 2.75) {$(\overline{A})^*$};
		\node [style=none] (8) at (0.5, 0.2) {$\epsilon^*$};
		\node [style=none] (9) at (-0.5, 3.1) {$\overline{~^*\eta}$};
		\node [style=none] (10) at (-1.5, -0.25) {$\overline{A^*}$};
		\node [style=none] (12) at (-0.25, 1.5) {$\overline{A}$};
		\node [style=none] (11) at (-1.5, 1.75) {$\overline{~^*A}$};
	\end{pgfonlayer}
	\begin{pgfonlayer}{edgelayer}
		\draw (2.center) to (0);
		\draw (0) to (1.center);
		\draw [bend left=90, looseness=2.00] (2.center) to (3.center);
		\draw (3.center) to (4.center);
		\draw [bend right=90, looseness=2.00] (4.center) to (5.center);
		\draw (5.center) to (6.center);
	\end{pgfonlayer}
\end{tikzpicture} = \begin{tikzpicture} 
	\begin{pgfonlayer}{nodelayer}
		\node [style=circle] (0) at (-1, 2.5) {$\psi$};
		\node [style=none] (1) at (-1, 3.25) {};
		\node [style=none] (2) at (-1, 0.5) {};
		\node [style=none] (3) at (0, 0.5) {};
		\node [style=none] (4) at (0, 1.75) {};
		\node [style=none] (5) at (1, 1.75) {};
		\node [style=none] (6) at (1, -0.25) {};
		\node [style=none] (7) at (-1.75, 3.25) {$(\overline{A})^*$};
		\node [style=none] (8) at (0.5, 2.6) {$\overline{\eta^*}$};
		\node [style=none] (9) at (-0.5, -0.35) {$~^*\epsilon$};
		\node [style=none] (10) at (1.5, 0.25) {$\overline{A^*}$};
		\node [style=none] (11) at (-1.75, 1.5) {$~^*(\overline{A})$};
		\node [style=none] (12) at (-0.25, 1.25) {$\overline{A}$};
	\end{pgfonlayer}
	\begin{pgfonlayer}{edgelayer}
		\draw (2.center) to (0);
		\draw (0) to (1.center);
		\draw [bend right=90, looseness=2.00] (2.center) to (3.center);
		\draw (3.center) to (4.center);
		\draw [bend left=90, looseness=2.00] (4.center) to (5.center);
		\draw (5.center) to (6.center);
	\end{pgfonlayer}
\end{tikzpicture}
\]

When the $*$-autonomous category is symmetric, conjugation automatically preserves the canonical cyclor.

\begin{lemma}
\label{Lemma: varepsi monoidal}
In a conjugative $*$-autonomous category, 
\[ 
\begin{tikzpicture} 
	\begin{pgfonlayer}{nodelayer}
		\node [style=circle] (0) at (0, 1) {$\varepsilon$};
		\node [style=none] (1) at (0, -0) {};
		\node [style=none] (2) at (-1.75, -0) {};
		\node [style=none] (3) at (-1.75, 2) {};
		\node [style=none] (4) at (0, 2) {};
		\node [style=none] (5) at (-1, 3) {$\overline{\overline{\eta*}}$};
		\node [style=none] (6) at (-2.25, 0.3) {$\overline{\overline{X^*}}$};
		\node [style=none] (7) at (0.5, 1.7) {$\overline{\overline{X}}$};
		\node [style=none] (8) at (0.5, 0.25) {$X$};
	\end{pgfonlayer}
	\begin{pgfonlayer}{edgelayer}
		\draw (4.center) to (0);
		\draw (0) to (1.center);
		\draw [bend left=90, looseness=1.50] (3.center) to (4.center);
		\draw (3.center) to (2.center);
	\end{pgfonlayer}
\end{tikzpicture} = 
\begin{tikzpicture}
	\begin{pgfonlayer}{nodelayer}
		\node [style=circle] (0) at (-1.75, 1) {$\varepsilon^{-1}$};
		\node [style=none] (1) at (-1.75, -0) {};
		\node [style=none] (2) at (0, -0) {};
		\node [style=none] (3) at (0, 2) {};
		\node [style=none] (4) at (-1.75, 2) {};
		\node [style=none] (5) at (-1, 3) {$\eta*$};
		\node [style=none] (6) at (-2.5, 0.3) {$\overline{\overline{X^*}}$};
		\node [style=none] (7) at (0.5, 0.5) {$X$};
		\node [style=none] (8) at (-2.25, 1.7) {$X^*$};
	\end{pgfonlayer}
	\begin{pgfonlayer}{edgelayer}
		\draw (4.center) to (0);
		\draw (0) to (1.center);
		\draw [bend right=90, looseness=1.50] (3.center) to (4.center);
		\draw (3.center) to (2.center);
	\end{pgfonlayer}
\end{tikzpicture}
\]
\[
\chi_\top^{\!\!\!\circ}~ \overline{\chi_\top^{\!\!\!\circ}} ~ \overline{\overline{\eta*}} ~ \overline{\chi_\oa}^{-1} \chi_\oa^{-1} (1 \oa \varepsilon) = \eta* (\varepsilon^{-1} \oa 1)  : \top \to \overline{\overline{X^*}} \oa X 
\]
\end{lemma}
\begin{proof}
\begin{align*}
\chi_\top^{\!\!\!\circ} ~ \overline{\chi_\top^{\!\!\!\circ}} ~ \overline{\overline{\eta*}} ~ \overline{\chi}^{-1} \chi^{-1} (1 \oa \varepsilon) &=  \chi_\top^{\!\!\!\circ} ~ \overline{\chi_\top^{\!\!\!\circ}} ~ \overline{\overline{\eta*}} ~ \overline{\chi}^{-1} \chi^{-1} (\varepsilon \varepsilon^{-1} \oa \varepsilon) \\
&=\chi_\top^{\!\!\!\circ} ~ \overline{\chi_\top^{\!\!\!\circ}}  ~ \overline{\overline{\eta*}} ~ \overline{\chi}^{-1} \chi^{-1} (\varepsilon  \oa \varepsilon) (\varepsilon^{-1} \oa 1) \\
&\stackrel{{\bf \small [CF.2]_\oa}}{=} \chi^{\!\!\!\circ} ~ \overline{\chi_\top^{\!\!\!\circ}} ~ \overline{\overline{\eta*}} \varepsilon (\varepsilon^{-1} \oa 1) \\
&\stackrel{{\bf \small nat.}}{=} \chi_\top^{\!\!\!\circ} ~ \overline{\chi_\top^{\!\!\!\circ}} ~  \varepsilon \eta\!* (\varepsilon^{-1} \oa 1) \\
&\stackrel{{\bf \small [CF.6]_\top}}{=} \eta\!* (\varepsilon^{-1} \oa 1)
\end{align*}
\end{proof}

\subsection{Dagger and conjugation}
The interaction of the dagger and conjugation for cyclic $*$-autonomous categories in the presence of the dualizing functor is 
illustrated by the following diagram:
\[
\xymatrixcolsep{5pc}
\xymatrixrowsep{5pc}
\xymatrix{
\X^{\op} \ar@/^1pc/[rr]^{(\_)^\dagger} \ar@/^1pc/[dr]|{((\_)^*)^\rev} \ar@{}[rr]|{\bot} &~ \ar@{}[d]|{\cong} & 
\X  \ar@/^1pc/[ll]|{((\_)^{\dagger})^{\op}} \ar@/_1pc/[dl]|{\overline{(\_)}^\rev} \ar@{}[dl]|{\dashv} \\
 & \X^\rev \ar@/^1pc/[ul]^{(\_)^{*^{\op}}} \ar@/_1pc/[ur]_{\overline{(\_)}} \ar@{}[ul]|{\dashv} &
 }
\]

Specifically we have: 

\begin{proposition}
	\label{Prop: dagger+dualizing}
Every cyclic $\dagger$-$*$-autonomous category is a conjugative $*$-autonomous category.
\end{proposition}
\begin{proof}
Let $\X$ be a cyclic, $\dagger$-$*$-autonomous category. Then composing adjoints gives the equivalence  $(\_)^{\dagger^*} 
\dashv (\_)^{*^\dagger}$. To build a conjugation, however,
we need an equivalence between the same functors: to obtain such an equivalence we use the natural equivalence $\omega: 
(\_)^{\dagger*} \to (\_)^{*\dagger} $ from the cyclor preserving condition for Frobenius linear functors.  
A conjugative equivalence, in addition, requires that the unit and counit of the equivalence be inverses of each of other.  
The unit and counit of the equivalence are given by {\em (a)} and {\em (b)} respectively;
\[ \mbox{\em (a)}~~~~~~~~~  \xymatrix{
\X^{\sf rev} \ar@{=}[rrrr] \ar[dr]_{(\_)^{*^\op}} & 
&  
\ar@{}[d]|{\Downarrow ~ \eta_\ox'}&  
& 
\X^{\sf rev} \ar@{<-}[ld]^{(\_)^{*^\op}}  
& \\
&
 \X^{\op} \ar@{=}[rr] \ar[dr]_{\dagger} &  
 \ar@{}[d]|{\Downarrow ~ \iota^{-1}} & 
 \X^{\sf op} \ar@{<-}[ld]^{\dagger^\op} \ar@{}[rr]_{\omega}  \ar@{}[rr]^{\Longrightarrow}  & 
 & 
 \X^{\sf oprev} \ar@{<-}@/^1pc/[llld]^{(\_)^{*^{\sf oprev}}} \ar@/_1pc/[ul]_{\dagger^\rev}  \\
& & \X & &
}
\]

\[
\mbox{\em (b)}~~~~~~~~~~~~~~ \xymatrix{
& & &
\X^{\rev}  \ar@{<-}@/_1pc/[llld]_{\dagger^\rev} \ar@{<-}[dl]_{(\_)^{*^{\rev}}} \ar[dr]^{(\_)^{*^{\sf op}}} \ar@{}[d]|{\Downarrow ~ \epsilon'_\ox} & &\\
\X^{\sf op rev} \ar@{<-}[rd]_{(\_)^{*^{\sf oprev}}} \ar@{}[rr]_{\omega^{-1}}  \ar@{}[rr]^{\Longrightarrow}  &  & 
\X^{\sf op} \ar@{=}[rr] \ar@{<-}[dl]^{\dagger^\op} & \ar@{}[d]|{\Downarrow ~ \iota^{-1}} &
\X^{\sf op} \ar[dr]^{\dagger} \\
& \X \ar@{=}[rrrr]&  & & & \X }
\] where the isomorphism $\omega: (\_)^{\dagger*} \to (\_)^{*\dagger}$ is from the cyclor preserving condition, {\bf [CFF]}, for Frobenius linear functors:
\[
\omega := 

\]

It remains to show that the unit and the counit maps are inverses of each other in $\X$:

\[ (a) ~~~~~ 
 %
 \stackrel{\ref{Lemma:  cyclic dagger}}{=} 
%
\]
 $(*)$ holds because $\dagger$ preserves the cyclor. Thus, $(a)$ and $(b)$ are inverses of each other.
\end{proof}

Next, we show that a conjugation functor together with a dualizing functors gives a $\dagger$:

\begin{proposition}
\label{Theorem: conjugation+dualizing}
Every cyclic, conjugative $*$-autonomous category  is also a $\dagger$-$*$-autonomous category.
\end{proposition}
\begin{proof}
Let $\X$ be a cyclic, conjugative $*$-autonomous category then  $\bar{(\_)^*} \dashv \bar{(\_)}^*$ 
is an equivalence. To build a dagger we need an equivalence on the same 
functor: we obtain this by using the natural equivalence $\sigma: \bar{(\_)^*} \to \bar{(\_)}^*$ from Definition
 \ref{Defn: conjugative cyclic}.  An involutive equivalence, in addition, 
requires the unit and counit of the (contravariant) equivalence to be the same map (which we called the involutor, $\iota$). 
We show that this is the case: 

The unit and counit of the equivalence is given by {\em (a)} and {\em (b)} respectively;
\[ \mbox{\em (a)}~~~~~~~~~  \xymatrix{
\X^{\op} \ar@{=}[rrrr] \ar[dr]_{(\_)^{*^{\sf rev}}} & 
&  
\ar@{}[d]|{\Downarrow ~ \eta_\ox'}&  
& 
\X^{\op} \ar@{<-}[ld]^{(\_)^{*^{\sf op}}}  
& \\
&
 \X^{\sf rev} \ar@{=}[rr] \ar[dr]_{\bar{(\_)}} &  
 \ar@{}[d]|{\Downarrow ~ \varepsilon^{-1}} & 
 \X^{\sf rev} \ar@{<-}[ld]^{\bar{(\_)}^{\sf rev}} \ar@{}[rr]_{\sigma}  \ar@{}[rr]^{\Longrightarrow}  & 
 & 
 \X^{\sf oprev} \ar@{<-}@/^1pc/[llld]^{(\_)^{*^{\sf oprev}}} \ar@/_1pc/[ul]_{\bar{(\_)}^\op}  \\
& & \X & &
}
\]

\[
\mbox{\em (b)}~~~~~~~~~~~~~~ \xymatrix{
& & &
\X^{\op}  \ar@{<-}@/_1pc/[llld]_{\bar{(\_)}^\op} \ar@{<-}[dl]_{(\_)^{*^{\op}}} \ar[dr]^{(\_)^{*^{\sf rev}}} 
\ar@{}[d]|{\Downarrow ~ \epsilon'_\ox} & &\\
\X^{\sf op rev} \ar@{<-}[rd]_{(\_)^{*^{\sf oprev}}} \ar@{}[rr]_{\sigma^{-1}}  \ar@{}[rr]^{\Longrightarrow}  &  & 
\X^{\sf rev} \ar@{=}[rr] \ar@{<-}[dl]^{\bar{(\_)}^{\sf rev}} & \ar@{}[d]|{\Downarrow ~ \varepsilon} &
\X^{\sf rev} \ar[dr]^{\bar{(\_)}} \\
& \X \ar@{=}[rrrr]&  & & & \X }
\] where $\sigma: \bar{A}^* \to \bar{A^*}$ is given in Definition \ref{Defn: conjugative cyclic}.  Below we 
show that the unit and counit coincide in $\X$.

\[
(a)~~~~~ 
 \stackrel{\ref{Lemma: varepsi monoidal}}{=}
%
 =: \iota^{-1}
\] 
\end{proof}

Observe that for composition of the dualizing functor and the conjugation functor to yield a dagger, and 
vice versa, a $*$-autonomous category is required to be cyclic with the cyclor being preserved by the 
conjugation (see Definition \ref{Defn: conjugative cyclic}) and the dagger (see just before Lemma \ref{Lemma: cyclic dagger}). 

Combining Propositions \ref{Prop: dagger+dualizing} and \ref{Theorem: conjugation+dualizing}, we get:

\begin{theorem} 
	Every cyclic $*$-autonomous category is conjugative $*$-autonomous if and only if it is 
	$\dagger$-$*$-autonomous.
\end{theorem}

\section{Examples: Dagger and conjugation}
\label{Sec: Examples Dagger and conjugation}

Let us now look at some examples of $\dagger$-isomix categories in which the dagger is given 
by dualizing and conjugation functors. 

\subsection{A group with conjugation considered as a category}

\begin{definition}
A {\bf group with conjugation} is a group $(G, ., e)$  together with a function $\overline{(\_)}: G \to G$ such that, 
for all $g \in G$, $\overline{\overline{g}} = g$, and for all $g, h \in G$, $\overline{g . h} = \overline{h} \overline{g}$, and 
$\overline e = e$.
\end{definition}

Let $(G,., e)$ be a group with conjugation. The discrete category $\D{ (G,., e})$ whose objects are the elements 
of the group is a monoidal category with the tensor product given by $g \ox h := g.h$, and the monoidal unit $e$. 
Moreover, $\D{(G,.,e)}$ is a compact closed category where $g^* := g^{-1}$, and it has a trivial conjugative cyclor 
(see Definition \ref{Defn: conjugative cyclic}). Thus, $\D{(G,.,e})$ which is a compact $\dagger$-isomix-$*$-autonomous 
category with $g^\dagger := \overline{g^*}$ is an example of how the conjugation gives rise to a dagger.

Here are some examples of groups with conjugation and the discrete categories given by them:
\begin{itemize}
\item Suppose we fix the group to be $(\C, +, 0)$ where the objects are complex numbers and the tensor product is addition. 
The dual and conjugation of complex numbers are given as follows: $(a+ib)^* = -a - ib$ and $\overline{a + ib} := a - ib$. 
Hence, \[(a + ib)^\dagger := \overline{(a+ib)^*} = \overline{(-a-ib)} = -a + ib\] 
\item Consider the multiplicative group $(\C^*, ., 1)$ where the objects are non-zero complex numbers and the 
tensor product is given by multiplication. The dualizing and the conjugation functors are given as follows: 
\[ (a+ib)^* = c+id, \text{ where } ac-bd=1 \text { and }  ad+bc=0 \]  \[ \overline{a + ib} := a-ib \] $(a+ib)^\dagger$ 
is given by $\overline{(a+ib)^*}$.
\item Suppose the group is fixed to be $\D{(P(x), +, 0})$ where $P(x)$ is a polynomial ring. $\D{(P(x), +, 0})$ is a 
conjugative compact closed category: $P(x)^* = -P(x)$ and $\overline{P(x)} = P(-x)$. Then, $P(x)^\dagger = -P(-x)$.
\item Consider the general linear group of degree 2, $(\mathbb{M}_2, . ,I_2)$ over complex numbers. Then, 
the discrete category $\D(\mathbb{M}_2, ., I_2)$ has a dualizing functor given by matrix inverse and conjugation 
is given by conjugate transpose: $\overline{\left(
\begin{matrix}
a+ib & m+in \\
c+id & p+iq
\end{matrix}
\right)} :=  \left(
\begin{matrix}
a-ib & c-id \\
m-in & p-iq
\end{matrix}
\right)
$. Then,  $\D(\mathbb{M}_2, .,I_2)$ is a $\dagger$-isomix $*$-autonomous category with:
 \[ \left(
  \begin{matrix}
 a+ib & m+in \\
 c+id & p+iq 
 \end{matrix}
 \right)^\dagger := \overline{\left(
 \begin{matrix}
a+ib & m+in \\
c+id & p+iq 
 \end{matrix}
 \right)^*} =  \left(
 \begin{matrix}
a-ib & c-id \\
m-in & p-iq 
\end{matrix}
 \right)^{-1} \] 
\end{itemize}

\subsection{Finiteness matrices and finiteness relations}
\label{Sec: Finiteness matrices}

We now describe the conjugation, $\overline{(\_)}$ and, thus, the dagger functor, $(\_)^\dagger$, for ${\sf FRel}$, and ${\sf FMat}(R)$.  Recall that the dagger functor is obtained by composing the conjugacy dualizing functors:  $X^\dagger = \overline{X}^*$.

\begin{lemma}
	\label{Lemma: conjugative cat}
$\FRel$, $\FMat(R)$, where $R$ is a conjugative rig, are conjugative isomix $*$-autonomous categories.
\end{lemma}

A conjugative rig is a rig with a conjugation $\overline{(\_)}: R \to R$ such that $r = \overline{\overline{r}}$ and the conjugation preserves addition, $\overline{0} = 0$ and $\overline{r + s} =\overline{r}+\overline{s}$, and preserves the multiplication $\overline{1} = 1$ and $\overline{r_1 \cdot r_2} = \overline{r_1} \cdot \overline{r_2}$.  

\begin{proof}
For $\FRel$ the conjugation functor is identity on objects and arrows and thus the dagger is the duality, $X^\dagger = X^{*}$.  $\FRel$ is a conjugative $*$-autonomous category with the required natural isomorphisms defined as follows:
\begin{itemize}
\item $\overline{X} \ox \overline{Y} \to^{\chi_\ox} \overline{Y \ox X}$ and $\chi_\oa$ are the symmetry maps. 
\item $\overline{ \overline{ A}} \to^{\varepsilon} A$ is the identity map.
\end{itemize} 

For $\FMat(R)$ the conjugation functor for is determined by the conjugation on the rig. The functor is still identity on the objects and the symmetry map is used to provide $\chi_\ox$ and $\chi_\oa$. For a finiteness matrix, $M$, $\overline{M}$ is given by conjugating all the elements of $M$. $\chi_\ox$, $\chi_\oa$, $\varepsilon$ are the characteristic function of their corresponding finiteness relations.

Both are a conjugative isomix category because $\top = \bot$.
\end{proof}

\subsection{Chu spaces}
\label{Section: Chu}

Applications of Chu Spaces to represent quantum systems have been studied in \cite{Abr12}, \cite{Abr13}. 
In this section we show that the Chu construction over a closed conjugative monoidal category,
 which  has pullbacks, produces a $\dagger$-isomix LDC, ${\sf Chu}_\X(I)$.  To get 
 the $*$-autonomous category and $\dagger$-structure on ${\sf Chu}_\X(I)$ we shall start 
 by explaining how one can produce  conjugative structure on the ${\sf Chu}$ category.  
 To achieve this we develop the structure of this category, starting with a 
 conjugative closed monoidal category, $\X$, which is not necessarily  symmetric.  
 Note that the fact that it is conjugative means that it is both left and right closed which 
 allows us to consider the non-commutative ${\sf Chu}$ construction: in this regard  
 we shall follow J\"urgen Koslowski's construction \cite{Jur06} using simplified ``Chu-cells'' on the 
 same dualizing object to obtain not a ${*}$-linear bicategory but a cyclic  $*$-autonomous category.  
 Furthermore, we shall choose a dualizing object which is conjugative in order  to obtain a conjugative 
 cyclic $*$-autonomous  category.

A conjugative object is an object $D$ of $\X$ with an isomorphism $d: \overline{D} \to D$ such that 
$\overline{d}d = \varepsilon: D \to \overline{\overline{D}}$.  
We can then define ${\sf Chu}_\X(D)$ as follows:

\begin{description}
\item[Objects:] $(A, B, \psi_0, \psi_1)$ where $\psi_0: A \ox B \to D$ and $\psi_1: B \ox A \to D$  in $\X$ 
(these are the simplified Chu cells). 
\item[Arrows:] $(f,g): (A, B, \psi_0,\psi_1) \to (A', B', \psi_0',\psi_1')$ where $f: A \to A'$ and $g: B' \to B$ 
and the following diagrams commutes:
\[ \xymatrix{
& A \ox B' \ar[dl]_{1 \ox g} \ar[dr]^{f \ox 1} & \\
A \ox B \ar[dr]_{\psi_0} & & A' \ox B' \ar[ld]^{\psi'_0} \\
& D &}
~~~~~\xymatrix{
& B' \ox A \ar[dl]_{g \ox 1} \ar[dr]^{1 \ox f} & \\
B \ox A \ar[dr]_{\psi_1} & & B' \ox A' \ar[ld]^{\psi'_1} \\
& D &}
\]
\item[Compositon:] $(f,g)(f',g') := (ff', g'g)$. Composition is well-defined as:
\[
\xymatrix{
&& A \ox B'' \ar[dl]_{1 \ox g'} \ar[dr]^{f \ox 1}   && \\
& A \ox B' \ar[dr]_{f \ox 1}  \ar[dl]_{1 \ox g}  & & A' \ox B'' \ar[dl]^{1 \ox g'}  \ar[dr]^{f' \ox 1} \\
A \ox B \ar[drr]^{\psi_0} & & A' \ox B' \ar[d]^{\psi'_0} && A'' \ox B'' \ar[lld]_{\psi_0''} \\
& & D & & 
}
\]
and similarly for the reverse Chu-maps: $\psi_1$, $\psi'_1$ and $\psi_1''$.
The {\bf identity maps} are $(1_A,1_B): (A, B, \psi_0,\psi_1) \to (A, B, \psi_0,\psi_1)$ as expected.

\item[Tensor product $\ox$:] $(A, B, \psi_0, \psi_1) \ox (A', B', \psi_0', \psi_1') := (A \ox A', E, \gamma_0, \gamma_1)$, 
where $E$ is the pullback in the following diagram:
\[
\xymatrix{
& E \ar[rd]^{\pi_1} \ar[ld]_{\pi_0} & \\ 
A' \multimap B \ar[rd]_-{1 \multimap \tilde{\psi_1}} & & B' \poppilol A \ar[ld]^-{\tilde{\psi_1'} \poppilol A}  \\
& A' \multimap (D \poppilol A) \to^{\simeq} (A' \multimap D) \poppilol A &
}
\]
with
\[
\infer{B \ox A \to^{\psi_1} B}{B \to^{\tilde{\psi_1}} D \poppilol A} ~~~~~~~~~~ 
\infer{A' \ox B' \to D}{B' \to^{\tilde{\psi_1'}} A' \multimap D}
\]
and, 
\[
\gamma_0 := (A \ox A') \ox E \to^{1 \ox \pi_0} (A \ox A') \ox (A' \multimap B) \to^{a_\ox} A \ox ( A' \ox A' \multimap B) 
\to^{1 \ox eval_{\multimap}} A \ox B \to^{\psi_0'} D
\]
\[
\gamma_1 := E \ox (A \ox A')  \to^{\pi_1 \ox 1} (B' \poppilol A) \ox (A \ox A')  \to^{a_\ox^{-1}} 
(B' \poppilol A \ox A) \ox A' \to^{eval_{\poppilol} \ox 1} B' \ox A' \to^{\psi_1'} D
\]
The tensor unit is $(I, D, u_\ox^l, u_\ox^r)$.
\end{description}

It is standard that ${\sf Chu}_\X(D)$  is a (non-commutative) $*$-autonomous category. Furthermore, it is cyclic because 
$$~^{*}(A,B,\psi_0,\psi_1) = (A,B,\psi_0,\psi_1)^{*} = (B,A,\psi_1,\psi_0).$$ 
In addition, ${\sf Chu}_\X(D)$  is conjugative with 
$$\overline{(A,B,\psi_0,\psi_1)} := (\overline{A},\overline{B},\chi \overline{\psi_1}d,\chi \overline{\psi_0}d)$$
and $\overline{(f,g)} = (\overline{f},\overline{g})$.
Finally being conjugative cyclic $*$-autonomous implies that one has a dagger!

In the case that $\X$ is a symmetric monoidal closed category we may recapture the usual Chu construction \cite{Bar06}, 
which we denote ${\sf Chus}_\X(D)$.  Consider the full subcategory 
of Chu-objects with special Chu-cells of the form $(A,B, \psi,c_\otimes \psi)$ in which the symmetry map is used to 
obtain the second cell, this gives an inclusion 
${\sf Chus}_\X(D) \to {\sf Chu}_\X(D)$.  

We observe that $\X$ is symmetric conjugative when this subcategory is closed under the conjugation:

\begin{lemma}
	\label{Lemma: Chus lemma}
If $\X$  is an conjugative symmetric monoidal closed category  and $d: \overline{D} \to D$ is an involutive object, 
then ${\sf Chus}_\X(D)$ is a conjugative symmetric $*$-autonomous  category.
\end{lemma}
\begin{proof}
It suffices to observe that the Chu-cells of $\overline{(A,B,\psi,c_\otimes\psi)}$ have the right form.  Using the coherence of the 
involution with symmetry, the first Chu-cell of this object has  $\chi \overline{c_\otimes \psi}d = c_\ox \chi \overline{\psi}d$
which is exactly the symmetry map applied to the second Chu-cell of the object as desired.
\end{proof}

To obtain an isomix category one can choose D = I. ${\sf Chus}_\X(I)$ is an isomix category because the unit for tensor 
and par are the same (namely $\top = \bot = (I,I, u_\ox^l = u_\ox^r)$).  The tensor unit is always a conjugative object 
since $(\chi^{\!\!\!\circ})^{-1}: \overline{I} \to I$; therefore, this is immediately a conjugative symmetric 
$*$-autonomous category.  Composing the conjugation with the dualizing functor gives us a dagger.   

\subsection{Category of Hopf modules in a $*$-autonomous category}

\label{Sec: HModx}
In this example\footnote{We thank J-S. P. Lemay for bringing our attention to this example.}, 
we start with any symmetric $*$-autonomous category, $\X$, and build the category of modules over a 
Hopf Algebra which is in turn a $\dagger$-$*$-autonomous category.  

First of all, it has been already proven in \cite{PaS09}, that the category of Hopf modules over a $\ox$-Hopf 
algebra in any symmetric $*$-autonomous category is also a $*$-autonomous category. Then we note that, whenever the Hopf Algebra is cocommutative, the resulting $*$-autonomous category has a conjugation functor. One can construct the dagger functor by composing the conjugation functor and dualizing functor as in Theorem \ref{Theorem: conjugation+dualizing}. We establish some basic definitions before describing the category of modules over a Hopf Algerba, ${\mbox{\bf H-Mod}}_\X$.

\begin{definition}
A {\bf bialgebra}  in a symmetric monoidal category is a 4-tuple 
$$(\nabla: B \ox B \to B, e:  I \to B, \Delta : B \to B \ox B, u: B \to I)$$
such that $(A, \nabla,e)$ is a monoid and $(A, \Delta, u)$  is a comonoid and $\nabla$ and $e$ are coalgebra 
homomorphisms with respect to the comultiplication and the counit.
\end{definition}

Note that instead of requiring that  $\nabla$ and $e$ are coalgebra homomorphisms, one could equivalently 
require $\Delta$ and $u$ are algebra homomorphims with respect to the multiplication and the unit.

 The components of a bialgebra are graphically depicted as follows:
 
 \[
 \begin{tikzpicture}
	\begin{pgfonlayer}{nodelayer}
		\node [style=none] (0) at (-3, 1) {};
		\node [style=none] (1) at (-2.5, 1) {};
		\node [style=none] (2) at (-2.75, 0.75) {};
		\node [style=none] (3) at (-3.25, 1.5) {};
		\node [style=none] (4) at (-2.25, 1.5) {};
		\node [style=none] (5) at (-2.75, 0.25) {};
	\end{pgfonlayer}
	\begin{pgfonlayer}{edgelayer}
		\draw (0.center) to (1.center);
		\draw (1.center) to (2.center);
		\draw (0.center) to (2.center);
		\draw [in=-90, out=30, looseness=1.00] (1.center) to (4.center);
		\draw [in=-90, out=150, looseness=1.00] (0.center) to (3.center);
		\draw (2.center) to (5.center);
	\end{pgfonlayer}
\end{tikzpicture} : A \ox A \to A ~~~~~~~ \begin{tikzpicture}
	\begin{pgfonlayer}{nodelayer}
		\node [style=none] (0) at (-3, 0.75) {};
		\node [style=none] (1) at (-2.5, 0.75) {};
		\node [style=none] (2) at (-2.75, 1) {};
		\node [style=none] (3) at (-3.25, 0.25) {};
		\node [style=none] (4) at (-2.25, 0.25) {};
		\node [style=none] (5) at (-2.75, 1.5) {};
	\end{pgfonlayer}
	\begin{pgfonlayer}{edgelayer}
		\draw (0.center) to (1.center);
		\draw (1.center) to (2.center);
		\draw (0.center) to (2.center);
		\draw [in=90, out=-30, looseness=1.00] (1.center) to (4.center);
		\draw [in=90, out=-150, looseness=1.00] (0.center) to (3.center);
		\draw (2.center) to (5.center);
	\end{pgfonlayer}
\end{tikzpicture} : A \to A \ox A ~~~~~~~~ \begin{tikzpicture}
	\begin{pgfonlayer}{nodelayer}
		\node [style=none] (0) at (-3, -0) {};
		\node [style=none] (1) at (-2.5, -0) {};
		\node [style=none] (2) at (-2.75, 0.25) {};
		\node [style=none] (3) at (-2.75, 1.5) {};
	\end{pgfonlayer}
	\begin{pgfonlayer}{edgelayer}
		\draw (0.center) to (1.center);
		\draw (1.center) to (2.center);
		\draw (0.center) to (2.center);
		\draw (2.center) to (3.center);
	\end{pgfonlayer}
\end{tikzpicture}: A \to I ~~~~~~~~~~ \begin{tikzpicture}
	\begin{pgfonlayer}{nodelayer}
		\node [style=none] (0) at (-3, 1.5) {};
		\node [style=none] (1) at (-2.5, 1.5) {};
		\node [style=none] (2) at (-2.75, 1.25) {};
		\node [style=none] (3) at (-2.75, 0) {};
	\end{pgfonlayer}
	\begin{pgfonlayer}{edgelayer}
		\draw (0.center) to (1.center);
		\draw (1.center) to (2.center);
		\draw (0.center) to (2.center);
		\draw (2.center) to (3.center);
	\end{pgfonlayer}
\end{tikzpicture} : I \to A
 \]
This gives a succinct graphical depiction of the coalgebra homomorphism laws; namely:
\[
\begin{tikzpicture}
	\begin{pgfonlayer}{nodelayer}
		\node [style=none] (0) at (-2, 2) {};
		\node [style=none] (1) at (-2.25, 1.75) {};
		\node [style=none] (2) at (-1.75, 1.75) {};
		\node [style=none] (3) at (-2, 2.75) {};
		\node [style=none] (4) at (-0.5, 2.75) {};
		\node [style=none] (5) at (-0.75, 1.75) {};
		\node [style=none] (6) at (-0.5, 2) {};
		\node [style=none] (7) at (-0.25, 1.75) {};
		\node [style=none] (8) at (-1.75, 0.25) {};
		\node [style=none] (9) at (-0.5, -0.75) {};
		\node [style=none] (10) at (-2.25, 0.25) {};
		\node [style=none] (11) at (-2, -0) {};
		\node [style=none] (12) at (-2, -0.75) {};
		\node [style=none] (13) at (-0.5, -0) {};
		\node [style=none] (14) at (-0.25, 0.25) {};
		\node [style=none] (15) at (-0.75, 0.25) {};
	\end{pgfonlayer}
	\begin{pgfonlayer}{edgelayer}
		\draw (0.center) to (1.center);
		\draw (1.center) to (2.center);
		\draw (2.center) to (0.center);
		\draw (3.center) to (0.center);
		\draw (6.center) to (5.center);
		\draw (5.center) to (7.center);
		\draw (7.center) to (6.center);
		\draw (4.center) to (6.center);
		\draw (11.center) to (10.center);
		\draw (10.center) to (8.center);
		\draw (8.center) to (11.center);
		\draw (12.center) to (11.center);
		\draw (13.center) to (15.center);
		\draw (15.center) to (14.center);
		\draw (14.center) to (13.center);
		\draw (9.center) to (13.center);
		\draw [in=15, out=-165, looseness=1.00] (5.center) to (8.center);
		\draw [bend right, looseness=1.25] (1.center) to (10.center);
		\draw [in=150, out=-30, looseness=1.00] (2.center) to (15.center);
		\draw [bend left, looseness=1.00] (7.center) to (14.center);
	\end{pgfonlayer}
\end{tikzpicture} = \begin{tikzpicture}
	\begin{pgfonlayer}{nodelayer}
		\node [style=none] (0) at (-2, -1.5) {};
		\node [style=none] (1) at (-2.25, -1.75) {};
		\node [style=none] (2) at (-1.75, -1.75) {};
		\node [style=none] (3) at (-2, -0.75) {};
		\node [style=none] (4) at (-1.75, 0.25) {};
		\node [style=none] (5) at (-2.25, 0.25) {};
		\node [style=none] (6) at (-2, -0) {};
		\node [style=none] (7) at (-2, -0.75) {};
		\node [style=none] (8) at (-2.75, 1) {};
		\node [style=none] (9) at (-1.25, 1) {};
		\node [style=none] (10) at (-2.75, -2.5) {};
		\node [style=none] (11) at (-1.25, -2.5) {};
	\end{pgfonlayer}
	\begin{pgfonlayer}{edgelayer}
		\draw (0.center) to (1.center);
		\draw (1.center) to (2.center);
		\draw (2.center) to (0.center);
		\draw (3.center) to (0.center);
		\draw (6.center) to (5.center);
		\draw (5.center) to (4.center);
		\draw (4.center) to (6.center);
		\draw (7.center) to (6.center);
		\draw [bend right=45, looseness=0.75] (4.center) to (9.center);
		\draw [bend right, looseness=1.25] (8.center) to (5.center);
		\draw [bend right, looseness=0.75] (1.center) to (10.center);
		\draw [bend left, looseness=1.00] (2.center) to (11.center);
	\end{pgfonlayer}
\end{tikzpicture} ~~~~~~~ \begin{tikzpicture} 
	\begin{pgfonlayer}{nodelayer}
		\node [style=none] (0) at (-1, 0.75) {};
		\node [style=none] (1) at (-1.25, 1) {};
		\node [style=none] (2) at (-0.75, 1) {};
		\node [style=none] (3) at (-1, -0) {};
		\node [style=none] (4) at (-1.25, -0.25) {};
		\node [style=none] (5) at (-0.75, -0.25) {};
		\node [style=none] (6) at (-1.75, -1) {};
		\node [style=none] (7) at (-0.25, -1) {};
	\end{pgfonlayer}
	\begin{pgfonlayer}{edgelayer}
		\draw (1.center) to (2.center);
		\draw (2.center) to (0.center);
		\draw (0.center) to (1.center);
		\draw (4.center) to (5.center);
		\draw (5.center) to (3.center);
		\draw (3.center) to (4.center);
		\draw (0.center) to (3.center);
		\draw [in=71, out=-135, looseness=1.00] (4.center) to (6.center);
		\draw [bend left, looseness=0.75] (5.center) to (7.center);
	\end{pgfonlayer}
\end{tikzpicture} =   \begin{tikzpicture}
	\begin{pgfonlayer}{nodelayer}
		\node [style=none] (0) at (-0.75, 0.75) {};
		\node [style=none] (1) at (-1, 1) {};
		\node [style=none] (2) at (-0.5, 1) {};
		\node [style=none] (3) at (-0.75, -1) {};
		\node [style=none] (4) at (0, 0.75) {};
		\node [style=none] (5) at (0.25, 1) {};
		\node [style=none] (6) at (0, -1) {};
		\node [style=none] (7) at (-0.25, 1) {};
	\end{pgfonlayer}
	\begin{pgfonlayer}{edgelayer}
		\draw (1.center) to (2.center);
		\draw (2.center) to (0.center);
		\draw (0.center) to (1.center);
		\draw (0.center) to (3.center);
		\draw (7.center) to (5.center);
		\draw (5.center) to (4.center);
		\draw (4.center) to (7.center);
		\draw (4.center) to (6.center);
	\end{pgfonlayer}
\end{tikzpicture} ~~~~~~~~ \begin{tikzpicture} 
	\begin{pgfonlayer}{nodelayer}
		\node [style=none] (0) at (-1, -0.75) {};
		\node [style=none] (1) at (-1.25, -1) {};
		\node [style=none] (2) at (-0.75, -1) {};
		\node [style=none] (3) at (-1, 0) {};
		\node [style=none] (4) at (-1.25, 0.25) {};
		\node [style=none] (5) at (-0.75, 0.25) {};
		\node [style=none] (6) at (-1.75, 1) {};
		\node [style=none] (7) at (-0.25, 1) {};
	\end{pgfonlayer}
	\begin{pgfonlayer}{edgelayer}
		\draw (1.center) to (2.center);
		\draw (2.center) to (0.center);
		\draw (0.center) to (1.center);
		\draw (4.center) to (5.center);
		\draw (5.center) to (3.center);
		\draw (3.center) to (4.center);
		\draw (0.center) to (3.center);
		\draw [in=-71, out=135, looseness=1.00] (4.center) to (6.center);
		\draw [bend right, looseness=0.75] (5.center) to (7.center);
	\end{pgfonlayer}
\end{tikzpicture} = \begin{tikzpicture}
	\begin{pgfonlayer}{nodelayer}
		\node [style=none] (0) at (-1, -0.75) {};
		\node [style=none] (1) at (-1.25, -1) {};
		\node [style=none] (2) at (-0.75, -1) {};
		\node [style=none] (3) at (-1, 1) {};
		\node [style=none] (4) at (0, -0.75) {};
		\node [style=none] (5) at (0.25, -1) {};
		\node [style=none] (6) at (0, 1) {};
		\node [style=none] (7) at (-0.25, -1) {};
	\end{pgfonlayer}
	\begin{pgfonlayer}{edgelayer}
		\draw (1.center) to (2.center);
		\draw (2.center) to (0.center);
		\draw (0.center) to (1.center);
		\draw (0.center) to (3.center);
		\draw (7.center) to (5.center);
		\draw (5.center) to (4.center);
		\draw (4.center) to (7.center);
		\draw (4.center) to (6.center);
	\end{pgfonlayer}
\end{tikzpicture} ~~~~~~~~ \begin{tikzpicture}
	\begin{pgfonlayer}{nodelayer}
		\node [style=none] (0) at (-1, -0.75) {};
		\node [style=none] (1) at (-1.25, -1) {};
		\node [style=none] (2) at (-0.75, -1) {};
		\node [style=none] (3) at (-1, 0.75) {};
		\node [style=none] (4) at (-1.25, 1) {};
		\node [style=none] (5) at (-0.75, 1) {};
	\end{pgfonlayer}
	\begin{pgfonlayer}{edgelayer}
		\draw (1.center) to (2.center);
		\draw (2.center) to (0.center);
		\draw (0.center) to (1.center);
		\draw (0.center) to (3.center);
		\draw (4.center) to (3.center);
		\draw (3.center) to (5.center);
		\draw (5.center) to (4.center);
	\end{pgfonlayer}
\end{tikzpicture} = I
\]

\begin{definition} An {\bf antipode} for a bialgebra $(B, \nabla, \bialgunitmap{0.8}, \Delta, 
    \bialgcounitmap{0.8})$ is an endomorphism $s: B \to B$ such that 
\[
\begin{tikzpicture} 
	\begin{pgfonlayer}{nodelayer}
		\node [style=none] (0) at (-2, 2.25) {};
		\node [style=none] (1) at (-2.25, 2) {};
		\node [style=none] (2) at (-1.75, 2) {};
		\node [style=none] (3) at (-2, 2.75) {};
		\node [style=none] (4) at (-1.75, 0.5) {};
		\node [style=none] (5) at (-2.25, 0.5) {};
		\node [style=none] (6) at (-2, 0.25) {};
		\node [style=none] (7) at (-2, -0.5) {};
		\node [style=circle, scale=1.5] (8) at (-1.5, 1.25) {};
		\node [style=none] (9) at (-1.5, 1.25) {$s$};
	\end{pgfonlayer}
	\begin{pgfonlayer}{edgelayer}
		\draw (0.center) to (1.center);
		\draw (1.center) to (2.center);
		\draw (2.center) to (0.center);
		\draw (3.center) to (0.center);
		\draw (6.center) to (5.center);
		\draw (5.center) to (4.center);
		\draw (4.center) to (6.center);
		\draw (7.center) to (6.center);
		\draw [bend left, looseness=1.00] (2.center) to (8);
		\draw [bend left, looseness=0.75] (8) to (4.center);
		\draw [bend right=45, looseness=1.00] (1.center) to (5.center);
	\end{pgfonlayer}
\end{tikzpicture} = \begin{tikzpicture}
	\begin{pgfonlayer}{nodelayer}
		\node [style=none] (0) at (-1.75, 2.25) {};
		\node [style=none] (1) at (-1.5, 2) {};
		\node [style=none] (2) at (-2, 2) {};
		\node [style=none] (3) at (-1.75, 2.75) {};
		\node [style=none] (4) at (-2, 0.5) {};
		\node [style=none] (5) at (-1.5, 0.5) {};
		\node [style=none] (6) at (-1.75, 0.25) {};
		\node [style=none] (7) at (-1.75, -0.5) {};
		\node [style=circle, scale=1.5] (8) at (-2.25, 1.25) {};
		\node [style=none] (9) at (-2.25, 1.25) {$s$};
	\end{pgfonlayer}
	\begin{pgfonlayer}{edgelayer}
		\draw (0.center) to (1.center);
		\draw (1.center) to (2.center);
		\draw (2.center) to (0.center);
		\draw (3.center) to (0.center);
		\draw (6.center) to (5.center);
		\draw (5.center) to (4.center);
		\draw (4.center) to (6.center);
		\draw (7.center) to (6.center);
		\draw [bend right, looseness=1.00] (2.center) to (8);
		\draw [bend right, looseness=0.75] (8) to (4.center);
		\draw [bend left=45, looseness=1.00] (1.center) to (5.center);
	\end{pgfonlayer}
\end{tikzpicture} = \begin{tikzpicture}
	\begin{pgfonlayer}{nodelayer}
		\node [style=none] (0) at (1, 1.5) {};
		\node [style=none] (1) at (1, 0.5) {};
		\node [style=none] (2) at (1, -0.5) {};
		\node [style=none] (3) at (1, 2.75) {};
		\node [style=none] (4) at (0.75, 0.75) {};
		\node [style=none] (5) at (1.25, 0.75) {};
		\node [style=none] (6) at (1, 0.5) {};
		\node [style=none] (7) at (1.25, 1.25) {};
		\node [style=none] (8) at (0.75, 1.25) {};
		\node [style=none] (9) at (1, 1.5) {};
	\end{pgfonlayer}
	\begin{pgfonlayer}{edgelayer}
		\draw (3.center) to (0.center);
		\draw (1.center) to (2.center);
		\draw (4.center) to (5.center);
		\draw (5.center) to (6.center);
		\draw (6.center) to (4.center);
		\draw (8.center) to (7.center);
		\draw (7.center) to (9.center);
		\draw (9.center) to (8.center);
	\end{pgfonlayer}
\end{tikzpicture}
\]

A {\bf Hopf algebra} is a bialgebra with an antipode. An {\bf involutive Hopf algebra} is a hopf algebra where 
the antipode is self-inverse.
\end{definition}

A standard example of a Hopf algebra is a group algebra over a field:  for all group elements $g$, $\nabla : 
g \mapsto g \ox g$, $\bialgcounitmap{0.8}: g \mapsto 1$, $\Delta: g \ox h \mapsto gh$ and $s: g \mapsto g^{-1}$.

\begin{lemma}
\label{Lemma: involutive s}
Suppose $\X$ is a symmetric monoidal category, then:

\begin{enumerate}[(i)]

\item \cite[Theroem 3.5]{Blu96} If $H$ is a commutative or a cocommutative Hopf Algebra in $\X$, then $s^2 = 1$ where $s$ is the antipode: 
so it is an involutive Hopf algebra.

\item \cite[Lemma 2.11]{Lem19} If $H$ is a commutative Hopf Algebra, then $s$ is a monoid homomorphism. If $H$ is a cocommutative 
Hopf Algebra, then $s$ is  a comonoid homomorphism.

\end{enumerate}

\end{lemma}

\begin{definition}
A {\bf left module} for a bialgebra $(B, \nabla, u, \Delta, e)$ is a tuple $(M, a_M^l:B \ox M \to M)$ 
such that $a_M^l$ is a $B$-action i.e., the following diagram commutes:

\[
\xymatrix{
M \ar[r]^{u_\ox^l} \ar@{=}[dr] & \top \ox M \ar[d]^{\bialgunitmap{0.8}} \\ 
& M
}
\]
\end{definition}

We graphically depict $a_m^l$ as follows:
\[
\begin{tikzpicture}[xscale=-1]
	\begin{pgfonlayer}{nodelayer}
		\node [style=none] (0) at (-2.75, 1.25) {};
		\node [style=none] (1) at (-2.5, 1.25) {};
		\node [style=none] (2) at (-2.75, 1) {};
		\node [style=none] (3) at (-2.75, 0.25) {};
		\node [style=none] (4) at (-2, 2) {};
		\node [style=none] (5) at (-2.75, 2) {};
	\end{pgfonlayer}
	\begin{pgfonlayer}{edgelayer}
		\draw (0.center) to (1.center);
		\draw (1.center) to (2.center);
		\draw (0.center) to (2.center);
		\draw (2.center) to (3.center);
		\draw (5.center) to (0.center);
		\draw [in=-105, out=30, looseness=1.25] (1.center) to (4.center);
	\end{pgfonlayer}
\end{tikzpicture}: B \ox M \to M
\]

giving the graphical presentation of the module laws:

\[
\begin{tikzpicture}
	\begin{pgfonlayer}{nodelayer}
		\node [style=none] (0) at (-1.5, 2.25) {};
		\node [style=none] (1) at (-0.5, 2.25) {};
		\node [style=none] (2) at (0.25, 2.25) {};
		\node [style=none] (3) at (-1, 0.5) {};
		\node [style=none] (4) at (0, -0.25) {};
		\node [style=none] (5) at (0.25, -0.25) {};
		\node [style=none] (6) at (0.25, -0.5) {};
		\node [style=none] (7) at (0.25, -1) {};
		\node [style=none] (8) at (-1.25, 0.75) {};
		\node [style=none] (9) at (-0.75, 0.75) {};
		\node [style=none] (10) at (-1, 0.5) {};
	\end{pgfonlayer}
	\begin{pgfonlayer}{edgelayer}
		\draw (2.center) to (5.center);
		\draw (6.center) to (5.center);
		\draw (5.center) to (4.center);
		\draw (4.center) to (6.center);
		\draw [in=-90, out=150, looseness=1.25] (4.center) to (3.center);
		\draw (6.center) to (7.center);
		\draw (8.center) to (9.center);
		\draw (9.center) to (10.center);
		\draw (10.center) to (8.center);
		\draw [bend right=15, looseness=1.00] (0.center) to (8.center);
		\draw [bend right=15, looseness=1.00] (9.center) to (1.center);
	\end{pgfonlayer}
\end{tikzpicture} = \begin{tikzpicture}
	\begin{pgfonlayer}{nodelayer}
		\node [style=none] (0) at (-1.75, 2) {};
		\node [style=none] (1) at (-1.5, 2) {};
		\node [style=none] (2) at (-1.5, 1.75) {};
		\node [style=none] (3) at (-1.5, 3) {};
		\node [style=none] (4) at (-2.25, 3) {};
		\node [style=none] (5) at (-0.75, -0) {};
		\node [style=none] (6) at (-0.75, 3) {};
		\node [style=none] (7) at (-0.75, 0.25) {};
		\node [style=none] (8) at (-1.5, 1.25) {};
		\node [style=none] (9) at (-1, 0.25) {};
		\node [style=none] (10) at (-0.75, -0.5) {};
	\end{pgfonlayer}
	\begin{pgfonlayer}{edgelayer}
		\draw (3.center) to (1.center);
		\draw (2.center) to (1.center);
		\draw (1.center) to (0.center);
		\draw (0.center) to (2.center);
		\draw [in=-90, out=150, looseness=1.25] (0.center) to (4.center);
		\draw (6.center) to (7.center);
		\draw (5.center) to (7.center);
		\draw (7.center) to (9.center);
		\draw (9.center) to (5.center);
		\draw [in=-90, out=150, looseness=1.25] (9.center) to (8.center);
		\draw (2.center) to (8.center);
		\draw (5.center) to (10.center);
	\end{pgfonlayer}
\end{tikzpicture} ~~~~\text{ and }~~ \begin{tikzpicture} 
	\begin{pgfonlayer}{nodelayer}
		\node [style=none] (0) at (-0.75, 1) {};
		\node [style=none] (1) at (-0.75, 3) {};
		\node [style=none] (2) at (-0.75, 1.25) {};
		\node [style=none] (3) at (-1, 1.25) {};
		\node [style=none] (4) at (-0.75, -0) {};
		\node [style=none] (5) at (-1.75, 3) {};
		\node [style=none] (6) at (-1.25, 3) {};
		\node [style=none] (7) at (-1.5, 2.75) {};
	\end{pgfonlayer}
	\begin{pgfonlayer}{edgelayer}
		\draw (1.center) to (2.center);
		\draw (0.center) to (2.center);
		\draw (2.center) to (3.center);
		\draw (3.center) to (0.center);
		\draw (0.center) to (4.center);
		\draw (5.center) to (6.center);
		\draw (6.center) to (7.center);
		\draw (7.center) to (5.center);
		\draw [bend right, looseness=1.00] (7.center) to (3.center);
	\end{pgfonlayer}
\end{tikzpicture} = \begin{tikzpicture}
	\begin{pgfonlayer}{nodelayer}
		\node [style=none] (0) at (-0.75, 1.75) {};
		\node [style=none] (1) at (-0.75, 3.25) {};
		\node [style=none] (2) at (-0.75, 1.75) {};
		\node [style=none] (3) at (-0.75, 0.25) {};
	\end{pgfonlayer}
	\begin{pgfonlayer}{edgelayer}
		\draw (1.center) to (2.center);
		\draw (0.center) to (2.center);
		\draw (0.center) to (3.center);
	\end{pgfonlayer}
\end{tikzpicture}
\] 

\begin{definition}
Let $\X$ be a $*$-autonomous category and $H$ be a Hopf $\ox$-algebra in $\X$. The category of left H-modules in $\X$, ${\mbox{\bf H-Mod}}_\X$ has:
\begin{description}
\item[Objects:]  Left $H$-modules $(A, a_A^l:H \ox A\to A)$:
\item[Arrows:] A module homomorphism $(A, a_A^L:H \ox A\to A) \to^{f}(B, a_B^L:H \ox B\to B)$ is a map $A \to^{f} B$ such that the following diagram commutes:
\[
\xymatrix{
H \ox A \ar[r]^{ a_A^L} \ar[d]_{1 \ox f} & A \ar[d]^{f} \\
H \ox B \ar[r]_{ a_B^L} & B
}
\]

This is graphically depicted as follows:
\[
 \begin{tikzpicture} 
	\begin{pgfonlayer}{nodelayer}
		\node [style=none] (0) at (0, 1) {};
		\node [style=none] (1) at (0, 0.75) {};
		\node [style=none] (2) at (-0.25, 1) {};
		\node [style=none] (3) at (0, 2) {};
		\node [style=none] (4) at (0, -0.5) {};
		\node [style=none] (5) at (-0.75, 2) {};
		\node [style=none] (6) at (-0.75, 1.5) {};
		\node [style=circle, scale=1.5] (7) at (0, 0.25) {};
		\node [style=none] (8) at (0, 0.25) {$f$};
	\end{pgfonlayer}
	\begin{pgfonlayer}{edgelayer}
		\draw (2.center) to (0.center);
		\draw (0.center) to (1.center);
		\draw (1.center) to (2.center);
		\draw [bend left, looseness=1.00] (2.center) to (6.center);
		\draw (6.center) to (5.center);
		\draw (3.center) to (0.center);
		\draw (1.center) to (7);
		\draw (7) to (4.center);
	\end{pgfonlayer}
\end{tikzpicture} =
\begin{tikzpicture} 
	\begin{pgfonlayer}{nodelayer}
		\node [style=none] (0) at (0, 0.5) {};
		\node [style=none] (1) at (0, 0.25) {};
		\node [style=none] (2) at (-0.25, 0.5) {};
		\node [style=none] (3) at (0, 2) {};
		\node [style=none] (4) at (0, -0.25) {};
		\node [style=none] (5) at (-0.75, 2) {};
		\node [style=none] (6) at (-0.75, 1.5) {};
		\node [style=circle, scale=1.5] (7) at (0, 1.25) {};
		\node [style=none] (8) at (0, 1.25) {$f$};
	\end{pgfonlayer}
	\begin{pgfonlayer}{edgelayer}
		\draw (2.center) to (0.center);
		\draw (0.center) to (1.center);
		\draw (1.center) to (2.center);
		\draw (1.center) to (4.center);
		\draw (3.center) to (7);
		\draw (7) to (0.center);
		\draw [bend left, looseness=1.00] (2.center) to (6.center);
		\draw (6.center) to (5.center);
		\draw[fill=black] (0.center) -- (1.center) -- (2.center) -- (0.center);
	\end{pgfonlayer}
\end{tikzpicture} \]

\end{description}
\end{definition}

Observe that any left action is indeed a module homomorphism.

\begin{theorem} \cite{PaS09}
Let $\X$ be symmetric $*$-autonomous category and $H$ be a $\ox$-Hopf Algebra in $\X$ with bijective antipode ($s^2 = 1$). Then, ${\mbox{\bf H-Mod}}_\X$ is a $*$-autonomous category. If the Hopf Algebra, $H$, is cocommutative, then ${\mbox{\bf H-Mod}}_\X$ is a symmetric $*$-autonomous category.
\end{theorem}
\begin{proof} (Sketch)
The monoidal product $\ox$ for ${\mbox{\bf H-Mod}}_\X$ is defined as follows:
\[ (A, \leftaction{0.4}{white}) \ox (B, \leftaction{0.4}{black}) := (A \ox B, \leftaction{0.4}{gray} ) 
\text{ where, } \leftaction{0.8}{gray} := \begin{tikzpicture}
	\begin{pgfonlayer}{nodelayer}
		\node [style=none] (0) at (-0.25, -0.5) {};
		\node [style=none] (1) at (-0.25, -0.75) {};
		\node [style=none] (2) at (-0.5, -0.5) {};
		\node [style=none] (3) at (-0.25, 0.25) {};
		\node [style=none] (4) at (-1, 0.75) {};
		\node [style=none] (5) at (-1, -0) {};
		\node [style=none] (6) at (-0.25, -1.25) {};
		\node [style=none] (7) at (1, -0.75) {};
		\node [style=none] (8) at (1, 0.25) {};
		\node [style=none] (9) at (1, -1.25) {};
		\node [style=none] (10) at (0.25, -0) {};
		\node [style=none] (11) at (1, -0.5) {};
		\node [style=none] (12) at (-0.5, 0.75) {};
		\node [style=none] (13) at (0.75, -0.5) {};
		\node [style=none] (14) at (0.25, 1.75) {};
		\node [style=none] (15) at (1, 1.75) {};
		\node [style=none] (16) at (-0.75, 1.75) {};
		\node [style=none] (17) at (-0.75, 1) {};
		\node [style=none] (18) at (-1, 0.75) {};
		\node [style=none] (19) at (-0.5, 0.75) {};
	\end{pgfonlayer}
	\begin{pgfonlayer}{edgelayer}
		\draw (2.center) to (0.center);
		\draw (0.center) to (1.center);
		\draw (1.center) to (2.center);
		\draw [bend left, looseness=1.00] (2.center) to (5.center);
		\draw (5.center) to (4.center);
		\draw (3.center) to (0.center);
		\draw (1.center) to (6.center);
		\draw (13.center) to (11.center);
		\draw (11.center) to (7.center);
		\draw (7.center) to (13.center);
		\draw [bend left, looseness=1.00] (13.center) to (10.center);
		\draw [in=-45, out=90, looseness=1.00] (10.center) to (12.center);
		\draw (8.center) to (11.center);
		\draw (7.center) to (9.center);
		\draw [in=-90, out=90, looseness=1.00] (3.center) to (14.center);
		\draw (15.center) to (8.center);
		\draw (18.center) to (19.center);
		\draw (19.center) to (17.center);
		\draw (17.center) to (18.center);
		\draw (16.center) to (17.center);
		\draw[fill=black] (0.center) --  (1.center) --  (2.center) --  (0.center);
	\end{pgfonlayer}
\end{tikzpicture} \]

 The unit of $\ox$ is given by $(\top, H \ox \top \to^{u_\ox^R} H \to^{e} \top)$, the left action is drawn as 
 $\begin{tikzpicture} 
	\begin{pgfonlayer}{nodelayer}
		\node [style=circle, scale=0.5] (0) at (-1, 0.25) {};
		\node [style=circle, scale=1.5] (1) at (0.25, 1) {};
		\node [style=none] (2) at (0.25, 1) {$\top$};
		\node [style=none] (3) at (-1, 2) {};
		\node [style=none] (4) at (0.25, 2) {};
		\node [style=none] (5) at (-0.25, -0.75) {$\top$};
		\node [style=none] (6) at (-1.5, 1.75) {$H$};
		\node [style=none] (7) at (0.5, 1.75) {$\top$};
		\node [style=none] (8) at (-1.25, -1) {};
		\node [style=none] (9) at (-0.75, -1) {};
		\node [style=none] (10) at (-1, -0.75) {};
		\node [style=none] (11) at (-1, -1) {};
	\end{pgfonlayer}
	\begin{pgfonlayer}{edgelayer}
		\draw (3.center) to (0);
		\draw [bend left, looseness=1.25, dotted] (1) to (0);
		\draw (4.center) to (1);
		\draw (8.center) to (9.center);
		\draw (9.center) to (10.center);
		\draw (8.center) to (10.center);
		\draw (0) to (10.center);
	\end{pgfonlayer}
\end{tikzpicture}$

 The par product is defined as: $(A, \leftaction{0.4}{white}) \oa (B, \leftaction{0.4}{black}) := (A \oa B, \leftaction{0.4}{gray} ) $
where,
\begin{align*}
 \leftaction{0.8}{gray} &:= H \ox (A \oa B) \to^{\Delta \ox 1} (H \ox H) \ox (A \oa B) \to^{a_\ox} H \ox (H \ox (A \oa B)) \\
  &\to^{c_\ox} (H \ox (A \oa B)) \ox H \to^{\partial^L \ox 1} ((H \ox A) \oa B) \ox H \to^{\partial^R} (H \ox A) \oa (B \ox H) \\
  &\to^{1 \oa c_\ox} (H \ox A) \oa (H \ox B) \to^{\leftaction{0.4}{white} \oa \leftaction{0.4}{black}} A \oa B
\end{align*}

and  the unit of $\oa$ is
\[
\bot := ( \bot, H \ox \bot \to^{u \ox \bot} \top \ox \bot \to^{u_\ox} \bot)
\]

All the basic natural isomorphisms are inherited directly from $\X$, and they are module homomorphisms. Thus, {\bf HMod}$_\X$ is a LDC.

The dualizing functor $(\_)^*$ is given as follows:$(A, \leftaction{0.4}{white}: H \ox A \to A)^* := (A^*, \leftaction{0.4}{white}*: H \ox A^* \to A^*) \ \text{ where, }$ \[ \leftaction{0.4}{white}* := 
\begin{tikzpicture} 
	\begin{pgfonlayer}{nodelayer}
		\node [style=none] (0) at (1.5, -0.25) {};
		\node [style=none] (1) at (1.5, -1) {};
		\node [style=none] (2) at (1.5, -0.5) {};
		\node [style=none] (3) at (2.25, -1) {};
		\node [style=none] (4) at (0.75, 0.5) {};
		\node [style=none] (5) at (0.75, 0.25) {};
		\node [style=none] (6) at (1.25, -0.25) {};
		\node [style=none] (7) at (1.5, 0.25) {};
		\node [style=none] (8) at (-0.25, 0.25) {};
		\node [style=none] (9) at (2.25, 3) {};
		\node [style=none] (10) at (-0.25, -1.25) {};
		\node [style=none] (11) at (0.75, 3) {};
		\node [style=none] (12) at (1, 2.75) {$H$};
		\node [style=none] (13) at (2.5, 2.75) {$A^*$};
		\node [style=none] (14) at (-0.25, -3) {};
		\node [style=none] (15) at (0, -2.75) {$A^*$};
		\node [style=circle, scale=1.5] (16) at (0.75, 2) {};
		\node [style=none] (17) at (0.75, 2) {$s$};
	\end{pgfonlayer}
	\begin{pgfonlayer}{edgelayer}
		\draw (6.center) to (0.center);
		\draw (0.center) to (2.center);
		\draw (2.center) to (6.center);
		\draw [bend left, looseness=1.00] (6.center) to (5.center);
		\draw (5.center) to (4.center);
		\draw (2.center) to (1.center);
		\draw [bend left=90, looseness=2.00] (8.center) to (7.center);
		\draw [bend right=90, looseness=2.00] (1.center) to (3.center);
		\draw (9.center) to (3.center);
		\draw (8.center) to (10.center);
		\draw (11.center) to (16);
		\draw (16) to (4.center);
		\draw (14.center) to (10.center);
		\draw (7.center) to (0.center);
	\end{pgfonlayer}
\end{tikzpicture}\]
 Equationally, 

\begin{align*}
\leftaction{0.4}{white}* &:= H \ox A^* \to^{s \ox 1} H \ox A^* \to^{u_\ox^{-1} \ox 1} (H \ox \top) \ox A^* \to^{1 \ox \eta \ox 1} (H \ox (A^* \oa A)) \ox A^* \\
&\to^{c_\ox \ox 1}  ((A^* \oa A) \ox H) \ox A^* \to^{\partial \ox 1} (A^* \oa (A \ox H)) \ox A^* \to^{1 \ox c_\ox \ox 1} (A^* \oa (H \ox A)) \ox A^* \\
& \to^{(1 \oa \leftaction{0.3}{white}) \ox 1} (A^* \oa A) \ox A^* \to^{\partial} A^* \oa (A \ox A^*) \to^{1 \oa \epsilon} A \oa \bot \to^{u_\oa^R} A^*
\end{align*}

The cups and caps are inherited directly from $\X$, hence the snake diagrams hold. 
The antipode in the definition of $\leftaction{0.4}{white}*: H \ox A^* \to A^*$ makes the cup and cap module morphisms.

Suppose $(A, \leftaction{0.4}{white}) \to^{f} (B, \leftaction{0.4}{black})$ is as a module morphism, 
then $f^* := B^* \to^{f^*} A^* \in \X$ which is also a module morphism. Thus, {\bf H-Mod}$_\X$ is a 
monoidal category with a dualizing functor, hence a $*$-autonomous category.

If $H$ is cocommutative, then $(A, \leftaction{0.4}{white}) \ox (B, \leftaction{0.4}{black}) \to^{c_\otimes} 
(B, \leftaction{0.4}{black}) \ox (A, \leftaction{0.4}{white})$ is a module homomorphism.

In that case, {\bf H-Mod}$_\X$ is a symmetric $*$-autonomous category.
\end{proof}

Futhermore, we can show that the category of Hopf modules is conjugative.

\begin{lemma}
Let $\X$ be a symmetric $*$-autonomous category. ${\mbox{\bf H-Mod}}_\X$, the category of modules over a 
cocommutative Hopf Algebra H is a conjugative symmetric $*$-autonomous category .
\end{lemma}
\begin{proof}
We already know that ${\mbox{\bf H-Mod}}_\X$ is a symmetric $*$-autonomous category. 
We define the conjugation functor $\bar{(\_)}: {\mbox{\bf H-Mod}}_\X \to {\mbox{\bf H-Mod}}_\X$ as follows:

\begin{itemize}
\item 
$\overline{(A, \leftaction{0.4}{white})} := (A, \overline{\leftaction{0.4}{white}})$ where,
$ \overline{\leftaction{0.5}{white}} := 
\begin{tikzpicture}[scale=1]
	\begin{pgfonlayer}{nodelayer}
		\node [style=none] (0) at (0, 1) {};
		\node [style=none] (1) at (0, 0.75) {};
		\node [style=none] (2) at (-0.25, 1) {};
		\node [style=none] (3) at (0, 2) {};
		\node [style=none] (4) at (0, 0.5) {};
		\node [style=circle] (5) at (-0.75, 1.5) {};
		\node [style=none] (6) at (-0.75, 2) {};
		\node [style=none] (7) at (-0.75, 1.5) {$s$};
	\end{pgfonlayer}
	\begin{pgfonlayer}{edgelayer}
		\draw (2.center) to (0.center);
		\draw (0.center) to (1.center);
		\draw (1.center) to (2.center);
		\draw (1.center) to (4.center);
		\draw (3.center) to (0.center);
		\draw [bend right, looseness=1.00] (5) to (2.center);
		\draw (5) to (6.center);
	\end{pgfonlayer}
\end{tikzpicture}
$
\item Suppose $f: (A, \leftaction{0.4}{white}) \to (B, \leftaction{0.4}{black})$, then $\overline{f} := f$
\end{itemize}

The  basic natural isomorphisms are given by: \[ \overline{(B, \leftaction{0.4}{black})} \ox 
\overline{(A, \leftaction{0.4}{white})}  \to^{\chi}  \overline{(A, \leftaction{0.4}{white}) \ox (B, \leftaction{0.4}{black})} := 
B \ox A \to^{(c_\ox)_{B,A}} A \ox B\]  \[(A, \overline{\overline{\leftaction{0.4}{white}}}) \to^{\varepsilon} 
(A, \leftaction{0.4}{white}) := 1\]

 The natural isormorphisms satisfy all the coherences of conjugative symmetric $*$-autonomous category.
\end{proof}

\begin{lemma}
\label{Lemma conjugative}
Suppose $\X$ is a symmetric (iso)mix $*$-autonomous category, then {\bf H-Mod}$_\X$, the category of Hopf 
modules over a cocommutative Hopf Algebra H is a (iso)mix conjugative symmetric $*$-autonomous category.
\end{lemma}
\begin{proof}~
The mix map $\m: \bot \to \top$ is inherited directly from $\X$.
\end{proof}

\begin{corollary}
Suppose $\X$ is a symmetric (iso)mix $*$-autonomous, then {\bf H-Mod}$_\X$, the category of modules over a 
cocommutative Hopf Algebra H is a symmetric $\dagger$ (iso)mix $*$-autonomous category.
\end{corollary}
\begin{proof}
From Lemma \ref{Lemma conjugative},  ${\mbox{\bf H-Mod}}_\X$ is an  (iso)mix conjugative symmetric 
$*$-autonomous category. Then, by Theorem \ref{Theorem: conjugation+dualizing} one can construct a 
dagger functor by composing the conjugation and the dualizing functor as follows:  $(\_)^\dagger := \overline{(\_)^*}:
{\mbox{\bf H-Mod}}_\X^{\op} \to {\mbox{\bf H-Mod}}_\X$. Therefore,

$(A, \leftaction{0.4}{white})^\dagger := (A^*, \overline{ \leftaction{0.4}{white}^*})$ where,
\[
(A, \overline{\leftaction{0.4}{white}}^*) := 
\begin{tikzpicture} 
	\begin{pgfonlayer}{nodelayer}
		\node [style=none] (0) at (1.5, -1) {};
		\node [style=none] (1) at (2.25, -1) {};
		\node [style=none] (2) at (1.5, 0.5) {};
		\node [style=none] (3) at (1.5, 0.75) {};
		\node [style=none] (4) at (0.25, 0.75) {};
		\node [style=none] (5) at (2.25, 3.5) {};
		\node [style=none] (6) at (0.25, -1.25) {};
		\node [style=none] (7) at (0.75, 3.5) {};
		\node [style=none] (8) at (1, 3.25) {$H$};
		\node [style=none] (9) at (2.5, 3.25) {$A^*$};
		\node [style=none] (10) at (0.25, -1.75) {};
		\node [style=none] (11) at (0.5, -1.5) {$A^*$};
		\node [style=circle, scale=1.5] (12) at (0.75, 2) {};
		\node [style=none] (13) at (0.75, 2.25) {};
		\node [style=none] (14) at (0.75, 2) {$s$};
		\node [style=circle, scale=1.5] (15) at (0.75, 2.75) {};
		\node [style=none] (16) at (0.75, 2.75) {$s$};
		\node [style=none] (17) at (1.5, 0.5) {};
		\node [style=none] (18) at (1.5, 0.25) {};
		\node [style=none] (19) at (1.25, 0.5) {};
	\end{pgfonlayer}
	\begin{pgfonlayer}{edgelayer}
		\draw [in=-90, out=90, looseness=1.00] (2.center) to (3.center);
		\draw [bend left=90, looseness=1.75] (4.center) to (3.center);
		\draw [bend right=90, looseness=2.00] (0.center) to (1.center);
		\draw (5.center) to (1.center);
		\draw (4.center) to (6.center);
		\draw (13.center) to (12);
		\draw (7.center) to (15);
		\draw (15) to (13.center);
		\draw (6.center) to (10.center);
		\draw (19.center) to (17.center);
		\draw (17.center) to (18.center);
		\draw (18.center) to (19.center);
		\draw [bend right=15, looseness=1.25] (12) to (19.center);
		\draw (18.center) to (0.center);
	\end{pgfonlayer}
\end{tikzpicture} =
\begin{tikzpicture} 
	\begin{pgfonlayer}{nodelayer}
		\node [style=none] (0) at (1.5, -1) {};
		\node [style=none] (1) at (2.25, -1) {};
		\node [style=none] (2) at (1.5, 0.5) {};
		\node [style=none] (3) at (1.5, 1) {};
		\node [style=none] (4) at (-0.25, 1) {};
		\node [style=none] (5) at (2.25, 3) {};
		\node [style=none] (6) at (0.5, 2.75) {$H$};
		\node [style=none] (7) at (2.5, 2.75) {$A^*$};
		\node [style=none] (8) at (-0.25, -1.75) {};
		\node [style=none] (9) at (0, -1.5) {$A^*$};
		\node [style=none] (10) at (0.75, 3) {};
		\node [style=none] (11) at (1.25, 0.5) {};
		\node [style=none] (12) at (1.5, 0.5) {};
		\node [style=none] (13) at (1.5, 0.25) {};
		\node [style=none] (14) at (2.75, 3.75) {};
	\end{pgfonlayer}
	\begin{pgfonlayer}{edgelayer}
		\draw [in=-90, out=90, looseness=1.00] (2.center) to (3.center);
		\draw [bend left=90, looseness=2.00] (4.center) to (3.center);
		\draw [bend right=90, looseness=2.00] (0.center) to (1.center);
		\draw (5.center) to (1.center);
		\draw (11.center) to (12.center);
		\draw (11.center) to (13.center);
		\draw (13.center) to (12.center);
		\draw (0.center) to (13.center);
		\draw [in=-90, out=135, looseness=1.00] (11.center) to (10.center);
		\draw (4.center) to (8.center);
	\end{pgfonlayer}
\end{tikzpicture}
\]
\end{proof}

Thus, one can generate a $\dagger$-isomix category from a symmetric isomix $*$-autonomous 
category by choosing the Hopf modules over any cocommutative $\ox$- Hopf Algebra. 

%% file: chapter5.tex

\chapter{Mixed unitary categories (MUCs)}
\label{Chap: MUCs}


The notion of unitary isomorphism is important in categorical quantum mechanics since these 
isomorphisms model the unitary evolution of a quantum system. An isomorphism 
in a $\dagger$-monoidal category is said to be unitary when the inverse of the map 
coincides with its dagger. This idea cannot be directly applied to define 
unitary isomorphisms in $\dagger$-LDCs due to the non-stationary dagger functor ($A \neq A^\dagger$). 
This arises the following question: {\em what are unitary isomorphisms in $\dagger$-LDCs?}
The objective of this chapter is resolve this question and to introduce mixed unitary 
categories (MUCs).   



\section{Unitary categories}
\label{Sec: unitary}


\subsection{Unitary structure}

 Categorically, within a $\dagger$-monoidal category, 
 a unitary map is an isomorphism $f: A \to B$ such that $f^{-1} =  f^\dagger$. 
This definition of unitary isomorphism cannot be used directly within the framework of $\dagger$-LDCs 
since the types of $f^{-1}: B \to A$ and  $f^\dagger: B^\dagger \to A^\dagger$ are different. 
However, one can define such a unitary isomorphism if, minimally, 
$A \simeq A^\dagger$ and $B \simeq B^\dagger$, and the isomorphisms behave 
coherently with the $\dagger$-linear structure. We call such isomorphisms 
{\em unitary structure} maps and the objects equipped with such isomorphisms 
as {\em unitary objects}:

\begin{definition}
	\label{defn: unitary structure}
A  $\dagger$-isomix category, $\X$ has {\bf unitary structure} in case there is an essentially small class of objects $\mathcal{U}$, called the {\bf unitary objects} of $\X$ such that
\begin{enumerate}[{\bf [U.1]}]
\item for all $A \in \mathcal{U}$, $A \in  \Core(\X)$, and $A$ is equipped with an isomorphism, $\varphi_A: A \to A^\dag$, called the {\bf unitary structure map} of $A$
\item $\mathcal{U}$ is closed under $(\_)^\dag$ so that for all $A \in \mathcal{U}$, $\varphi_{A^\dag} = ((\varphi_A)^{-1})^\dag$ 
\item for all $A \in \mathcal{U}$, the following diagram commutes:
 \[   \xymatrix{  A   \ar[d]_{\varphi_A} \ar[drrr]^{\iota}  & \\ A^\dag \ar[rrr]_{\varphi_{A^\dag}}  & & & (A^\dag)^\dag  } \]
\item $\bot, \top \in \mathcal{U}$ satisfy:
\[ \xymatrixcolsep{2pc}
\xymatrix{
\bot \ar[r]^{\varphi_\bot} \ar[d]_{\lambda_\bot} \ar[dr]^{\m} & \bot^\dagger \ar[d]^{\lambda_\top^{-1}}  \\
\top^\dagger \ar[r]_{\varphi_\top^{-1}} & \top
}
\]
\item If $A , B \in \mathcal{U}$, then $A \ox B$ and $A \oa B \in \mathcal{U}$ satisfy:
\[ (a) ~~~~~ \xymatrixcolsep{3pc}
\xymatrix{
A \ox B \ar[r]^{\varphi_A \ox \varphi_B}_{\simeq} \ar@/_2pc/[rrr]_{\mx}&
 A^\dagger \ox B^\dagger \ar[r]^{\lambda_\oa}_{\simeq} & 
 (A \oa B) ^\dagger \ar[r]^{\varphi_{A \oa B}^{-1}} _{\simeq} &
A \oa B
}
\]
\[ (b) ~~~~~ \xymatrixcolsep{3pc}
\xymatrix{
A \ox B \ar[r]^{\varphi_{A \ox B}}_{\simeq} \ar@/_2pc/[rrr]_{\mx}&
 (A \ox B)^\dagger \ar[r]^{\lambda_\ox^{-1}}_{\simeq} & 
 A^\dagger \oa B^\dagger \ar[r]^{\varphi_A^{-1} \oa \varphi_B^{-1}} _{\simeq} &
A \oa B
}
\]
 \end{enumerate}
\end{definition}


\begin{lemma}
\label{Lemma: square root tensor unitary}
When $A$ and $B$ are unitary objects in a $\dagger$-isomix category then, $\varphi_{A^{\dagger\dagger}} = (\varphi_A)^{\dagger \dagger}: A^{\dagger\dagger} \to A^{\dagger \dagger \dagger}$.
\end{lemma}
\begin{proof}~
\[ \varphi_{(A^\dagger)^{\dagger}} = ((\varphi_{A^\dagger})^{-1})^{\dagger} = ((((\varphi_A)^{-1})^\dagger)^{-1})^\dagger = ((((\varphi_A)^{-1})^{-1})^\dagger)^\dagger = ((\varphi_A)^\dagger)^\dagger \]
\end{proof}


Often we shall want the unitary objects to have linear adjoints (or duals) but we shall need the analogue of $\dagger$-duals $(\eta^\dagger = c_\ox \epsilon$ and $\epsilon^\dagger = \eta c_\ox)$ from categorical quantum mechanics:

\begin{definition} \label{defn: unitary dual}
A {\bf unitary linear duality} $(\eta, \epsilon): A \dashvv_{~u} B$ between unitary objects  $A$ and $B$ is a linear duality satisfying in addition:
\[
\begin{matrix}
\xymatrix{ \\
{\bf [Udual.]} \\
}~~~
\xymatrix{
\top \ar@{}[ddrr]|{(a)} \ar[rr]^{\eta} \ar[d]_{\lambda_\top}  & & A \oa B \ar[d]^{\varphi_A \oa \varphi_B} \\
\bot^\dagger \ar[d]_{\epsilon^\dag} & & A^\dagger \oa B^\dagger \ar[d]^{c_\oa} \\ 
(B \ox A)^\dag \ar[rr]_{\lambda_\oa^{-1}} & & B^\dagger \oa A^\dagger} 
~~~~~ & \text{(or)} & ~~~~~~
\xymatrix{
A \ox B \ar@{}[ddrr]|{(b)} \ar[rr]^{\varphi_A \ox \varphi_B} \ar[d]_{c_\ox} & & A^\dag \ox B^\dag \ar[d]^{\lambda_\ox} \\
B \ox A \ar[d]_{\epsilon} & & (A \oa B)^\dagger \ar[d]^{\eta^\dagger} \\
\bot \ar[rr]_{\lambda_\bot} & & \top^\dagger } 
\end{matrix}
\]
\end{definition}

Observe that ${\bf [Udual.]} (a) \Leftrightarrow (b)$. In a compact $\dagger$-LDC, $\top \dashvv_{~u} \bot$. {\bf [Udual] (a)} is shown diagrammatically as follows:
\[\begin{tikzpicture}
	\begin{pgfonlayer}{nodelayer}
		\node [style=none] (0) at (-2, 4) {};
		\node [style=none] (1) at (1, 4) {};
		\node [style=none] (2) at (-2, 2) {};
		\node [style=none] (3) at (1, 2) {};
		\node [style=circle] (4) at (-0.5, 2.5) {$\epsilon$};
		\node [style=none] (5) at (-1.25, 4) {};
		\node [style=none] (6) at (0.25, 4) {};
		\node [style=none] (7) at (-1.25, 2) {};
		\node [style=none] (8) at (0.25, 2) {};
		\node [style=none] (9) at (-1.25, 1) {};
		\node [style=none] (10) at (0.25, 1) {};
	\end{pgfonlayer}
	\begin{pgfonlayer}{edgelayer}
		\draw [in=150, out=-90, looseness=1.25] (5.center) to (4);
		\draw [in=30, out=-90, looseness=1.25] (6.center) to (4);
		\draw (0.center) to (1.center);
		\draw (1.center) to (3.center);
		\draw (3.center) to (2.center);
		\draw (2.center) to (0.center);
		\draw (7.center) to (9.center);
		\draw (8.center) to (10.center);
	\end{pgfonlayer}
\end{tikzpicture}  = \begin{tikzpicture}
	\begin{pgfonlayer}{nodelayer}
		\node [style=circle] (0) at (0.5, 5.75) {$\eta$};
		\node [style=none] (1) at (-0.5, 4.75) {};
		\node [style=none] (2) at (0, 4.75) {};
		\node [style=none] (3) at (-0.25, 4.5) {};
		\node [style=none] (4) at (1, 4.75) {};
		\node [style=none] (5) at (1.5, 4.75) {};
		\node [style=none] (6) at (1.25, 4.5) {};
		\node [style=none] (7) at (-0.25, 3) {};
		\node [style=none] (8) at (1.25, 3) {};
		\node [style=none] (9) at (-0.25, 4.75) {};
		\node [style=none] (10) at (1.25, 4.75) {};
	\end{pgfonlayer}
	\begin{pgfonlayer}{edgelayer}
		\draw (1.center) to (2.center);
		\draw (2.center) to (3.center);
		\draw (3.center) to (1.center);
		\draw (4.center) to (5.center);
		\draw (5.center) to (6.center);
		\draw (6.center) to (4.center);
		\draw [in=90, out=-150, looseness=1.00] (0) to (9.center);
		\draw [in=90, out=-30, looseness=1.00] (0) to (10.center);
		\draw [in=90, out=-75, looseness=0.75] (6.center) to (7.center);
		\draw [in=90, out=-90, looseness=1.00] (3.center) to (8.center);
	\end{pgfonlayer}
\end{tikzpicture} \]

\begin{lemma}
Suppose $(\eta_1, \epsilon_1): V_1 \dashvv_{~u} U_1$ and $(\eta_2, \epsilon_2): V_2 \dashvv_{~u} U_2$. Then, $(V_1 \otimes V_2) \dashvv_{~u} (U_1 \oa U_2)$.
\end{lemma}
\begin{proof}
Define $(\eta', \epsilon'): (V_1 \otimes V_2) \dashvv_{~u} (U_1 \oa U_2)$ so that  
$\eta' = \begin{tikzpicture} 
	\begin{pgfonlayer}{nodelayer}
		\node [style=circle] (0) at (-4, 3) {$\eta_1$};
		\node [style=circle] (1) at (-2, 3) {$\eta_2$};
		\node [style=ox] (2) at (-4, 1.75) {};
		\node [style=oa] (3) at (-2, 1.75) {};
		\node [style=none] (4) at (-4, 1) {};
		\node [style=none] (5) at (-2, 1) {};
	\end{pgfonlayer}
	\begin{pgfonlayer}{edgelayer}
		\draw [style=none, in=15, out=-165, looseness=1.00] (1) to (2);
		\draw [style=none, bend left, looseness=1.25] (1) to (3);
		\draw [style=none, in=180, out=-15, looseness=1.00] (0) to (3);
		\draw [style=none, bend left=45, looseness=1.25] (2) to (0);
		\draw [style=none] (2) to (4.center);
		\draw [style=none] (3) to (5.center);
	\end{pgfonlayer}
\end{tikzpicture} ~~~~~~~ 
\epsilon' = \begin{tikzpicture} 
	\begin{pgfonlayer}{nodelayer}
		\node [style=circle] (0) at (-4, 1) {$\epsilon_1$};
		\node [style=circle] (1) at (-2, 1) {$\epsilon_2$};
		\node [style=oa] (2) at (-4, 2.25) {};
		\node [style=ox] (3) at (-2, 2.25) {};
		\node [style=none] (4) at (-4, 3) {};
		\node [style=none] (5) at (-2, 3) {};
	\end{pgfonlayer}
	\begin{pgfonlayer}{edgelayer}
		\draw [style=none, in=-15, out=165, looseness=1.00] (1) to (2);
		\draw [style=none, bend right, looseness=1.25] (1) to (3);
		\draw [style=none, in=180, out=15, looseness=1.00] (0) to (3);
		\draw [style=none, bend right=45, looseness=1.25] (2) to (0);
		\draw [style=none] (2) to (4.center);
		\draw [style=none] (3) to (5.center);
	\end{pgfonlayer}
\end{tikzpicture}
$. This is easily checked to be a unitary linear adjoint.
\end{proof}

We can now define what it means for an isomorphism to be unitary:

\begin{definition}
Suppose $A$ and $B$ are unitary objects. An isomorphism $A\xrightarrow{f} B$ is said to be a {\bf unitary isomorphism} if the following diagram commutes:
\[  \xymatrix{A   \ar[r]^{\varphi_A}    \ar[d]_{f} \ar[r]^{\varphi_A} & A^\dag \\ B  \ar[r]_{\varphi_B} & B^\dag  \ar[u]_{f^\dag}  }  \]
\end{definition}

Observe that $\varphi$ is ``twisted'' natural for all unitary isomorphisms, thus, unitary isomorphisms compose and contain the identity maps. In a category in which the unitary structure maps are identity morphisms, one recovers the usual notion of unitary isomorphisms.

Our next objective is to show that all the coherence isomorphisms between unitary objects are unitary maps. First a warmup: 

\begin{lemma}
\label{lemma:MUCProperties}
In a $\dagger$-isomix category with unitary structure:
\begin{enumerate}[(i)]
\item If $f$ is a unitary isomorphism, then so is $f^\dagger$;
\item If $f$ and $g$ are unitary, then so are $f \ox g$ and $f \oa g$;
\item Unitary isomorphisms are closed under composition.
\end{enumerate}
\end{lemma}

\begin{proof}~
\begin{enumerate}[{\em (i)}]
\item Recall that $\varphi_{A^\dag} = (\varphi_A^{-1})^\dag$,  then $f^\dagger$ is unitary because 
\[ \xymatrix{B^\dag \ar[d]_{(\varphi_B^{-1})^\dag = \varphi_{B^\dag}} \ar[rr]^{f^\dag} & & A^\dag \ar[d]^{(\varphi_A^{-1})^\dag = \varphi_{A^\dag}} \\
   B^{\dag\dag}  & & A^{\dag\dag} \ar[ll]^{f^{\dag\dag}}} \]
is just the dagger functor applied to the unitary diagram of $f$.
\item Suppose $f$ and $g$ are unitary morphisms, then:
\[
\xymatrix{
A \ox B \ar@{->}[rrr]^{\varphi_{A \ox B}} \ar[ddd]_{f \ox g} \ar[dr]_{\mx}  \ar@{}[dddr]|{\mbox{\tiny \bf (nat. $\mx$)}~~~}
&   \ar@{}[dr]|{\mbox{\tiny {\bf [U.5(b)]}}} &  &  (A \ox B)^\dagger \ar@{}[lddd]|{~~~~~\mbox{\tiny \bf (nat. $\lambda_\oa)$}} \\
& A \oa B \ar[r]^{\varphi_A \oa \varphi_B} \ar[d]_{f \oa g} 
& A^\dagger \oa B^\dagger \ar[ur]_{\lambda_\oa}  & 
\\ & A' \oa B' \ar[r]_{\varphi_{A'} \ox \varphi_{B'}} \ar@{}[dr]|{\mbox{\tiny { \bf [U.5(b)]}}}
& A'^\dagger \oa B'^\dagger \ar[dr]^{\lambda_\oa} \ar[u]_{f^\dagger \oa g^\dagger}
& \\ A' \ox B' \ar[rrr]_{\varphi_{A' \ox B'}} \ar[ur]_{\mx}
& & & (A' \ox B')^\dagger \ar[uuu]_{(f \ox g)^\dagger}
}
\]
The inner square commutes because $f$ and $g$ are unitary maps.
Similarly, using {\bf [U.5(b)]}, one can show that if $f$ and $g$ are unitary, then $f \oa g$ is unitary. 

\item
The proof is trivial.
\end{enumerate}
\end{proof}

The following lemma will be used to prove that the natural isomorphisms in a $\dagger$-isomix category are 
unitary for unitary objects.
\begin{lemma}
\label{lemma: auxiliary}
The following diagram commutes:
\[ \xymatrix{
(A \ox B) \ox C \ar[r]^{\mx} \ar[d]_{a_\ox}  & (A \ox B) \oa C  \ar[r]^{\mx \oa 1}  & (A \oa B) \oa C \ar[d]^{a_\oa} \\
A \ox (B \ox C) \ar[r]_{\mx} & A \oa ( B \ox C)  \ar[r]_{1 \oa \mx} & A \oa (B \oa C) }
\]
\end{lemma}
\begin{proof} The given diagram commutes due to the naturality of the mixor, and due to the
	rules governing the interaction of mixor, associator and distributor, see Section \ref{Sec: mix, isomix, compact LDC}. 
\[ \xymatrix{
(A \ox B) \ox C \ar[r]^{\mx} \ar[d]_{a_\ox} \ar@{}[dr]|{{\sf \bf mix}~(b)} & (A \ox B) \oa C \ar@{<-}[d]^{\partial^L}  \ar[r]^{\mx \oa 1} \ar@{}[dr]|{{\sf \bf mix}~(a)}  & (A \oa B) \oa C \ar[d]^{a_\oa} \\
A \ox (B \ox C) \ar[r]_{1 \ox \mx}  \ar@/{_1pc}/[dr]_{\mx}  & A \ox (B \oa C) \ar[r]_{\mx} \ar@{}[d]|{nat. \mx} & A \oa (B \oa C) \\
& A \oa ( B \ox C)  \ar@/{_1pc}/[ur]_{1 \oa \mx}  & } \]
\end{proof}

\begin{lemma}
\label{lemma:cohUnitary}
	
Suppose $\X$ is a $\dagger$-isomix category with unitary structure and 
$A$, $B$, and $C$ are unitary objects. Then the following are unitary isomorphisms:

\begin{multicols}{2}
\begin{enumerate}[(i)]
\item $\lambda_\ox: A^\dagger \ox B^\dagger \rightarrow (A \oa B)^\dagger$
\item $\lambda_\oa:  A^\dagger \oa B^\dagger \rightarrow (A \ox B)^\dagger$
\item $\lambda_\top: \top \rightarrow \bot^\dagger$
\item $\lambda_\bot: \bot \to \top^\dagger$
\item $\varphi_A: A \rightarrow A^\dagger$ 
\item $m: \top \rightarrow \bot$
\item $\mx_{A,B}: A \ox B \rightarrow A \oa B$
\item $\iota : A \rightarrow (A^{\dagger})^\dagger$
\item $a_\ox: (A \ox B) \ox C \rightarrow A \ox (B \ox C)$
\item $a_\oa: (A \oa  B) \oa C \rightarrow A \oa (B \oa C)$
\item $c_\ox: A \ox B \rightarrow B \ox A$
\item $c_\oa: A \oa B \rightarrow B \oa A$
\item $\partial_L: A \ox (B \oa C) \rightarrow (A \ox B) \oa C$
\item $\partial_R: (A \oa B) \ox C \rightarrow A \oa (B \ox C)$
\end{enumerate}
\end{multicols}
\end{lemma}

\begin{proof}~
\begin{enumerate}[(i)]
\item $\lambda_\ox: A^\dagger \ox B^\dagger \rightarrow (A \oa B)^\dagger$ is a unitary map because:

\[\xymatrixcolsep{4pc}\xymatrix{
	{}&&&&\\
	A^\dag\ox B^\dag \ar[r]^{\phi_A^{-1}\ox\phi_B^{-1}} \ar[d]^{\lambda_\ox} \ar@{=}@/^3pc/[rr]     \ar@{}[dr]|{\mbox{\tiny {\bf [U.5(a)] }}}
	 &  A\ox B \ar[r]^{\phi_A\ox \phi_B}  \ar[d]^{\mx} 	                     \ar@{}[dr]|{\mbox{\tiny {\bf nat.}}}
	 &  A^\dag \ox B^\dag \ar[r]^{\phi_{A^\dag \ox B^\dag}} \ar[d]^{\mx}									 \ar@{}[dr]|{\mbox{\tiny {\bf [U.5(a)]}}}
	 &  (A^\dag \ox B^\dag)^\dag \ar[d]_{\lambda_\pr^{-1}} \ar@{=}@/^4pc/[ddd]\\
	(A\pr B)^\dag \ar[r]^{\phi_{A\pr B}^{-1}}  \ar[ddr]_{\phi_{(A\pr B)^\dag}} 								  \ar@{}[dr]|{\mbox{\tiny { \bf [U.3]}}}
	 & A\pr B  \ar[r]^{\phi_A\pr\phi_B} \ar@{=}[d]
	 & A^\dag \pr B^\dag \ar[r]^{\phi_{A^\dag}\pr\phi_{B^\dag}}										 \ar@{}[d]|{\mbox{\tiny { \bf [U.3]}}\ \pr\mbox{\tiny { \bf [U.3]}}}
	 & (A^\dag)^\dag \pr (B^\dag)^\dag  \ar@{=}[d]\\
   {}
     & A \pr B  \ar[rr]^{\iota \pr \iota} \ar[d]^{\iota}
     & {}																					\ar@{}[d]|{\mbox{\tiny {\bf [$\dagger$-ldc.5(a)]}}}
     & (A^\dag)^\dag \pr (B^\dag)^\dag  \ar[d]_{\lambda_\pr}\\
   {}
     & ((A \pr B)^\dag)^\dag  \ar[rr]_{\lambda_\ox^\dag}
     & {}
     & (A^\dag \ox B^\dag)^\dag
}\]

\item $\lambda_\oa$ is unitary because:

\[
\xymatrix{
A^\dag\pr B^\dag                            \ar[rrr]^{\phi_{A^\dag\pr B^\dag}} \ar[dr]^{\mx^{-1}}   \ar[ddd]_{\lambda_\pr} \ar@{}[dddr]|{\mbox{\tiny {\bf Lem. \ref{lemma: mixdagger}}}} 
  &
  &
  &
  (A^\dag\pr B^\dag)^\dag               \ar@{}[dddl]|{\mbox{\tiny {\bf (Lem. \ref{lemma: mixdagger})}}^\dag} \\
{}
  & A^\dag\ox B^\dag                       \ar[r]^{\phi_{A^\dag\ox B^\dag}} \ar[d]_{\lambda_\ox}   \ar@{}[ur]|{\mbox{\tiny {\bf Lem. \ref{lemma:cohUnitary} (vi)}}} \ar@{}[dr]|{\mbox{\tiny {\bf Lem. \ref{lemma:cohUnitary} (i)}}}
  & (A^\dag\ox B^\dag)^\dag           \ar[ur]^{(\mx^{-1})^\dag}
  &\\
{}
  & (A\pr B)^\dag                       \ar[r]^{\phi_{(A\pr B)^\dag}}  \ar@{}[dr]|{\mbox{\tiny {\bf Lems. \ref{lemma:cohUnitary} (vi), \ref{lemma:MUCProperties} (i)}}}
  & ((A\pr B)^\dag)^\dag                \ar[u]_{\lambda_\ox^\dag} \ar[dr]^{((\mx^{-1})^\dag)^\dag}
  &\\
(A\ox B)^\dag                            \ar[rrr]^{\phi_{(A\pr B)^\dag}} \ar[ur]^{(\mx^{-1})^\dag} 
  &
  &
  &
  ((A\ox B)^\dag)^\dag              \ar[uuu]_{\lambda_\pr^\dag}
}
\]

\item $\lambda_\bot: \bot \rightarrow \top^\dagger$ is unitary because:
\[ 
\xymatrix{
\bot \ar[d]_{\lambda_\bot} \ar[rr]^{\varphi_\bot} & & \bot^\dagger \ar[d]^{(\lambda_\bot^{-1})^{\dagger}} \\
\top^\dagger \ar[urr]_{m^\dagger} \ar[rr]_{\varphi_{\top^\dagger} = (\varphi_{\top}^{-1})^\dagger} &  &\top^{\dagger \dagger}
}
\]
The left triangle commutes by {\bf [U.4]} and  {\bf [$\dagger$-mix]}.  The right triangle commutes by {\bf [U.4]} and the functoriality of $\dag$.

\item $\lambda_\top: \top \rightarrow \bot^\dagger$ is unitary because:

\[ 
\xymatrix{
\top \ar[d]_{\lambda_\top} \ar[rr]^{\varphi_\top} & & \top^\dagger \ar[d]^{(\lambda_\top^{-1})^{\dagger}} \\
\bot^\dagger \ar[urr]_{(m^{-1})^\dagger} \ar[rr]_{\varphi_{\bot^\dagger} = (\varphi_{\bot}^{-1})^\dagger} &  &\bot^{\dagger \dagger}
}
\]

The left triangle commutes by {\bf [U.4]} and  {\bf [$\dagger$-mix]}.  The right triangle commutes by {\bf [U.4]} and the functoriality of $\dag$.

\item $\varphi_A$ is unitary because the following square commutes by {\bf [U.3]} and {\bf [U.4]}.
\[
\xymatrix{
A \ar[r]^{\varphi_A} \ar[d]_{\varphi_A} & A^\dagger \ar[d]^{(\varphi^{-1})^\dagger} \\
A^\dagger \ar[r]^{\varphi_{A^\dagger}} & A^{\dagger \dagger}
}
\]

\item $m: \bot \rightarrow \top$ is unitary because:
\[
\xymatrix{
\bot \ar[r]^{\varphi_\bot} \ar[d]_{\m} & \bot^\dagger \ar[d]^{(\m^{-1})^\dagger} \\
\top \ar[r]_{\varphi_\top} \ar[ur]^{\lambda_\top} & \top^\dagger
}
\]
The left and right triangles commute by {\bf [U.4]} and {\bf [$\dagger$-mix]} respectively. Hence, the outer squares commutes.

\item $\mx_{A,B}: A \ox B \rightarrow A \oa B$ is unitary as:
\[ \xymatrix{
{}
  &
  &
  &\\
A \ox B \ar[dd]_\mx\ar@/^2.5pc/[rrr]^{\varphi_{A \ox B}} \ar[r]^\mx \ar@/_1.5pc/[drr]_{\varphi_A \ox \varphi_B}  \ar@{}[drr]|{\mbox{\tiny {\bf nat.}}}   \ar@{}[urrr]|{\mbox{\tiny {\bf [U.5(a)]}}}  
  & A \oa B \ar[r]^{\varphi_A \oa \varphi_B} 
  & A^\dag \oa B^\dag \ar[r]^{\lambda_\oa}  
  & (A \ox B)^\dag \\
{}
  &
  & A^\dag \ox B^\dag \ar[u]^\mx \ar[dr]^{\lambda_\ox}   \ar@{}[ur]|{\mbox{\tiny {\bf Lem. \ref{lemma: mixdagger}}}} \\
A \oa B \ar[rrr] _{\varphi_{A\oa B}}    \ar@{}[urr]|{\mbox{\tiny {\bf [U.3]}}} 
  &
  &
  & (A \oa B)^\dag \ar[uu]_{\mx^\dag}
}
\]
                    
\item $\iota: A \rightarrow A^{\dagger \dagger}$ is unitary as in
\[
\xymatrix{
A \ar[d]_{\iota} \ar[r]^{\varphi_A} & A^\dagger \ar[d]^{(\iota^{-1})^\dagger} \ar[ld]_{\varphi_{A^\dagger}} \\ 
A^{\dagger \dagger} \ar[r]_{\varphi_{A^{\dagger \dagger}}} & A^{\dagger \dagger \dagger}
}
\] 
the left triangle commutes by {\bf [U.3]} and the right triangle commutes by: 
\begin{eqnarray*}
(\iota^{-1})^\dagger &= & ((\varphi_{A^\dagger})^{-1} \varphi_A^{-1})^\dagger =  (((\varphi_{A}^{-1})^\dagger)^{-1} \varphi_A^{-1})^\dagger \\
                         & =  & ((\varphi_A^\dagger)(\varphi_A^{-1}) )^\dagger = ( \varphi_A^{-1} )^{\dagger} ( \varphi_A )^{\dagger \dagger} \\
                         & = &  \varphi_{A^\dagger} (\varphi_A)^{\dagger \dagger} = \varphi_{A^\dagger} (\varphi_{A^{\dagger \dagger}})
\end{eqnarray*}

\item $a_\ox$ is unitary as:
%
%

\[
\xymatrix{
(A \ox B) \ox C                                           \ar[rrr]^{\varphi_{(A \ox B) \ox C}} \ar[ddddd]_{a_\ox} \ar[dr]_{\mx}   \ar@{}[drrr]|{\mbox{\tiny {\bf [U.5(b)]}}} \ar@{}[dddddr]|{\mbox{\tiny \bf Lem.~ \ref{lemma: auxiliary}}}
 & {} 
 & {}
 & ( (A \ox B) \ox C )^\dagger                     \\ 
 & {}
(A \ox B) \pr C                                             \ar[r]^{\varphi_{A\ox B} \pr \varphi_{C}} \ar[d]_{\mx\pr 1}   \ar@{}[dr]|{\mbox{\tiny {\bf  [U.5(b)]$\pr$(id)}}} 
 & (A \ox B)^\dag \pr C^\dag                               \ar[d]^{ \lambda_\oa^{-1}\pr 1}   \ar[ur]_{\lambda_\oa}   \ar@{}[dddr]|{\mbox{\tiny {\bf [\dag-ldc.1]}}} 
 & {}  \\ 
{}
 & (A \pr B) \pr C                                            \ar[r]_{(\varphi_A \oa \varphi_B) \oa \varphi_C} \ar[d]_{a_\oa}  \ar@{}[dr]|{\mbox{\tiny {\bf nat.}}} 
 & (A^\dag \pr B^\dag) \pr C^\dag                         \ar[d]^{a_\oa} 
 & {} \\ 
{}
 & A \pr (B \pr C)                                           \ar[r]_{\phi_A \pr (\phi_B \pr \phi_C)}   \ar@{}[dr]|{\mbox{\tiny {\bf (id)$\pr$[U.5(b)] }}} 
 & A^\dag \pr (B^\dag \pr C^\dag)       \ar[d]^{1\oa\lambda_\oa}
 & {} \\ 
{}
 & A \pr (B \ox C)                                          \ar[r]_{\varphi_{A} \pr \varphi_{B\ox C}} \ar[u]^{1\pr \mx}
 & A^\dag \pr (B \ox C)^\dag               \ar[dr]^{\lambda_\oa} 
 & {} \\ 
A \ox (B \ox C)                                          \ar[rrr]^{\varphi_{A \ox (B \ox C)}} \ar[ur]^{\mx}   \ar@{}[urrr]|{\mbox{\tiny {\bf [U.5(b)]}}} 
 & {}
 & {}
 & ( A \ox (B \ox C)  )^\dagger \ar[uuuuu]_{a_\ox^\dag} 
}
\]

\item $a_\oa$ is unitary because:

\[
\xymatrix{
(A \pr B) \pr C                                           \ar[rrr]^{\varphi_{(A \pr B) \pr C}} \ar[ddddd]_{a_\pr} \ar[dr]_{\mx^{-1}}   \ar@{}[drrr]|{\mbox{\tiny {\bf [U.5(a)]}}}  \ar@{}[dddddr]|{\mbox{\tiny \bf Lem.~ \ref{lemma: auxiliary}}} 
 & {} 
 & {}
 & ( (A \pr B) \pr C )^\dagger                     \\ 
 & {}
(A \pr B) \ox C                                             \ar[r]^{\varphi_{A\pr B} \ox \varphi_{C}} \ar[d]_{\mx^{-1}\ox 1}   \ar@{}[dr]|{\mbox{\tiny {\bf  [U.5(a)]$\ox$(id)}}} 
 & (A \pr B)^\dag \ox C^\dag                               \ar[d]^{ \lambda_\ox^{-1}\ox 1}   \ar[ur]_{\lambda_\pr}   \ar@{}[dddr]|{\mbox{\tiny {\bf [\dag-ldc.1]}}} 
 & {}  \\ 
{}
 & (A \ox B) \ox C                                            \ar[r]_{(\varphi_A \pr \varphi_B) \pr \varphi_C} \ar[d]_{a_\ox}  \ar@{}[dr]|{\mbox{\tiny {\bf nat.}}} 
 & (A^\dag \ox B^\dag) \ox C^\dag                         \ar[d]^{a_\ox} 
 & {} \\ 
{}
 & A \ox (B \ox C)                                           \ar[r]_{\phi_A \ox (\phi_B \ox \phi_C)}   \ar@{}[dr]|{\mbox{\tiny {\bf (id)$\ox$[U.5(a)] }}} 
 & A^\dag \ox (B^\dag \ox C^\dag)       \ar[d]^{1\ox\lambda_\ox}
 & {} \\ 
{}
 & A \ox (B \pr C)                                          \ar[r]_{\varphi_{A} \ox \varphi_{B\pr C}} \ar[u]^{1\ox \mx^{-1}}
 & A^\dag \ox (B \oa C)^\dag               \ar[dr]^{\lambda_\ox} 
 & {} \\ 
A \oa (B \oa C)                                          \ar[rrr]^{\varphi_{A \oa (B \oa C)}} \ar[ur]^{\mx^{-1}}   \ar@{}[urrr]|{\mbox{\tiny {\bf [U.5(a)]}}} 
 & {}
 & {}
 & ( A \oa (B \oa C)  )^\dagger \ar[uuuuu]_{a_\oa^\dag} 
}
\]

\item $c_\ox$ is unitary because:

\[
\xymatrix{
A\ox B                            \ar[rrr]^{\phi_{A\ox B}} \ar[dr]^{\mx} \ar[ddd]_{c_\ox}
 & {}                                  \ar@{}[dr]|{\mbox{\tiny {\bf [U.5(b)]}}}
 & {}
 & (A\ox B)^\dag              \\
{}                                    
 & A \pr B                        \ar[r]^{\phi_A \pr \phi_B} \ar[d]_{c_\pr}   \ar@{}[dr]|{\mbox{\tiny {\bf nat.}}}
 & A^\dag \pr B^\dag       \ar[d]^{c_\pr} \ar[ur]^{\lambda_\pr}           \ar@{}[dr]|{\mbox{\tiny {\bf [$\dagger$-ldc.2(b)]}}}
 & {} \\
{}
 & B \pr A                       \ar[r]^{\phi_B \pr \phi_A}
 & B^\dag \pr A^\dag      \ar[dr]^{\lambda_\pr}
 & {} \\
B\ox A                           \ar[ur]^{\mx} \ar[rrr]^{\phi_{B\ox A}}
& {}                                 \ar@{}[ur]|{\mbox{\tiny {\bf [U.5(b)]}}}
& {}
& (B\ox A)^\dag                 \ar[uuu]_{(c_\ox^{-1})^\dag = c_\ox^\dag}\\
}
\]

where the left square commutes because

$$
\begin{tikzpicture}
	\begin{pgfonlayer}{nodelayer}
		\node [style=circ] (0) at (0.5, -0.25) {};
		\node [style=circ] (1) at (0, -1) {$\top$};
		\node [style=map] (2) at (0, -1.75) {};
		\node [style=circ] (3) at (0, -2.5) {$\bot$};
		\node [style=circ] (4) at (-0.5, -3.25) {};
		\node [style=none] (5) at (0.5, -3.25) {};
		\node [style=none] (6) at (-0.5, -4.25) {};
		\node [style=none] (7) at (0.5, -4.25) {};
		\node [style=none] (8) at (-0.5, 0.5) {};
		\node [style=none] (9) at (0.5, 0.5) {};
	\end{pgfonlayer}
	\begin{pgfonlayer}{edgelayer}
		\draw [densely dotted, in=-90, out=45, looseness=1.00] (4) to (3);
		\draw (3) to (2);
		\draw (2) to (1);
		\draw [densely dotted, in=-135, out=90, looseness=1.00] (1) to (0);
		\draw [style=none] (0) to (9);
		\draw [style=none] (8) to (4);
		\draw [style=none] (0) to (5);
		\draw [style=none, in=90, out=-90, looseness=1.00] (5) to (6);
		\draw [style=none, in=90, out=-90, looseness=1.00] (4) to (7);
	\end{pgfonlayer}
\end{tikzpicture}
=
\begin{tikzpicture}
	\begin{pgfonlayer}{nodelayer}
		\node [style=circ] (0) at (0, -1.25) {$\top$};
		\node [style=map] (1) at (0, -2) {};
		\node [style=circ] (2) at (0.5, -3.5) {};
		\node [style=circ] (3) at (0, -2.75) {$\bot$};
		\node [style=circ] (4) at (-0.5, -0.5) {};
		\node [style=none] (5) at (0.5, -0.5) {};
		\node [style=none] (6) at (-0.5, 0.5) {};
		\node [style=none] (7) at (0.5, 0.5) {};
		\node [style=none] (8) at (0.5, -4.25) {};
		\node [style=none] (9) at (-0.5, -4.25) {};
	\end{pgfonlayer}
	\begin{pgfonlayer}{edgelayer}
		\draw [densely dotted, in=-90, out=150, looseness=1.25] (2) to (3);
		\draw (3) to (1);
		\draw (1) to (0);
		\draw [densely dotted, in=-45, out=90, looseness=1.00] (0) to (4);
		\draw [style=none] (8) to (2);
		\draw [style=none] (9) to (4);
		\draw [style=none, in=-90, out=90, looseness=1.00] (4) to (7);
		\draw [style=none, in=-90, out=90, looseness=1.00] (5) to (6);
		\draw [style=none] (5) to (2);
	\end{pgfonlayer}
\end{tikzpicture}
=
\begin{tikzpicture}
	\begin{pgfonlayer}{nodelayer}
		\node [style=circ] (0) at (0.5, -0.25) {};
		\node [style=circ] (1) at (0, -1) {$\top$};
		\node [style=map] (2) at (0, -1.75) {};
		\node [style=circ] (3) at (0, -2.5) {$\bot$};
		\node [style=circ] (4) at (-0.5, -3.25) {};
		\node [style=none] (5) at (-0.5, -4) {};
		\node [style=none] (6) at (0.5, -4) {};
		\node [style=none] (7) at (-0.5, 0.75) {};
		\node [style=none] (8) at (0.5, 0.75) {};
		\node [style=none] (9) at (-0.5, -0.25) {};
	\end{pgfonlayer}
	\begin{pgfonlayer}{edgelayer}
		\draw [densely dotted, in=-90, out=45, looseness=1.00] (4) to (3);
		\draw (3) to (2);
		\draw (2) to (1);
		\draw [densely dotted, in=-135, out=90, looseness=1.00] (1) to (0);
		\draw [style=none] (6) to (0);
		\draw [style=none] (5) to (4);
		\draw [style=none] (4) to (9);
		\draw [style=none, in=-90, out=90, looseness=1.00] (9) to (8);
		\draw [style=none, in=-90, out=90, looseness=1.00] (0) to (7);
	\end{pgfonlayer}
\end{tikzpicture}
$$

\item $c_\pr$ is unitary because:

\[
\xymatrix{
A\pr B                            \ar[rrr]^{\phi_{A\pr B}}  \ar[dr]^{\mx^{-1}}   \ar[ddd]_{c_\pr}  \ar@{}[drrr]|{\mbox{\tiny {\bf Lem. \ref{lemma:cohUnitary} (vii)}}}
  &
  &
  & (A\pr B)^\dag \\
{}
  & A \ox B                      \ar[r]^{\phi_{A\ox B}} \ar[d]_{c_\ox}    \ar@{}[dr]|{\mbox{\tiny {\bf Lem. \ref{lemma:cohUnitary} (xi)}}}
  & (A\ox B)^\dag     \ar[ur]^{(\mx^{-1})^\dag} &\\
{}
  & B \ox A                      \ar[r]^{\phi_{B\ox A}} 
  & (B \ox A)^\dag    \ar[dr]^{(\mx^{-1})^\dag}  \ar[u]_{c_\ox^\dag} &\\
B\pr A                            \ar[rrr]^{\phi_{B\pr A}}   \ar[ur]^{\mx^{-1}}   \ar@{}[urrr]|{\mbox{\tiny {\bf Lem. \ref{lemma:cohUnitary} (vii)}}}
  &
  &
  & (B\pr A)^\dag              \ar[uuu]_{c_\pr^\dag} 
  }
\]

where the left square commutes for the same reason and the right square is the dagger of the left square.

\item $\partial_L$ is unitary see Figure \ref{Fig: linear dist. unitary}.

\begin{figure}
\newpage

\begin{sideways}
\scalebox{.84}{
$
\xymatrix{
{}
& {}
& {}
& {}
& {}
& {}
& {}
& {}
& {}\\
{}
& {}
& {}
& {}
& {}
& {}
& {}
& {}
& {}\\
{}
 & {}
 & A\ox (B \pr C)                      \ar[rr]^{\phi_{A}\ox \phi_{B\pr C}}    \ar[d]_{\mx}       \ar@/^4pc/[rrrrrr]^{\phi_{A\ox (B\pr C)}} \ar@/_6pc/[dddddd]_{\partial_L}
 & {}                                                                                                                                          \ar@{}[d]|{\mbox{\tiny {\bf [U.6(b)]}}}
 & A^\dag \ox (B \pr C)^\dag            \ar@=[r]   \ar@=[d]                                                                \ar@{}[dr]|{\mbox{\tiny {\bf id}}}
 & A^\dag \ox (B \pr C)^\dag            \ar[r]^{\mx}    \ar[d]^{1 \ox \lambda_\ox^{-1} }     \ar@{}[dr]|{\mbox{\tiny {\bf nat.}}}    \ar@{}[u]|{\mbox{\tiny {\bf [U6.(b)]}}}  
 & A^\dag \pr (B \pr C)^\dag            \ar[rr]^{\lambda_\pr}    \ar[d]^{1 \pr \lambda_\ox^{-1}}                                              \ar@{}[ddrr]|{\mbox{\tiny {\bf [$\dagger$-ldc.4(b)]}}}  
 & {}
 & (A\ox (B\pr C))^\dag                         \\
{}
 & {}
 & A\pr (B \pr C)                      \ar[r]^{\phi_{A\pr (B \pr C)}}     \ar[d]_{a_\pr^{-1}}                                \ar@{}[dr]|{\mbox{\tiny {\bf Lem. \ref{lemma:cohUnitary} (ix)}}}
 & (A\pr (B\pr C))^\dag               \ar[r]^{\lambda_\ox^{-1}}      \ar[d]^{(a_\pr)^\dag}
 & A^\dag \ox (B\pr C)^\dag           \ar[r]^{1\ox \lambda_\ox^{-1}}                                                 \ar@{}[d]|{\mbox{\tiny {\bf [$\dagger$-ldc.1(b)]}}}
 & A^\dag \ox (B^\dag \ox C^\dag)     \ar[r]^{\mx }    \ar[d]^{a_\ox^{-1}}                                                 \ar@{}[dr]|{\mbox{\tiny {\bf mx.~cat.}}}
 & A^\dag \pr (B^\dag \ox C^\dag)             
 & {}
 & {}\\
{}
 & {}
 & (A\pr B)\pr C                       \ar[r]^{\phi_{(A\pr B)\pr C}}     \ar@=[d]
 & ((A\pr B)\pr C)^\dag                \ar[r]^{\lambda_\ox^{-1}}                                                            \ar@{}[d]|{\mbox{\tiny {\bf [U.6(a)]}}}
 & (A\pr B)^\dag \ox C^\dag           \ar[r]^{\lambda_\ox^{-1}\ox 1}    \ar[d]_{\mx}                             \ar@{}[dr]|{\mbox{\tiny {\bf nat.}}}
 & (A^\dag \ox B^\dag) \ox C^\dag      \ar[r]^{ \mx\ox 1}     \ar[d]^{\mx}                                            \ar@{}[dr]|{\mbox{\tiny {\bf nat.}}}
 & (A^\dag \pr B^\dag) \ox C^\dag                 \ar[d]^{\mx}\ar@=[r]          \ar[u]_{\partial_R}                                          \ar@{}[ddddr]|{\mbox{\tiny {\bf nat.}}} 
 & (A^\dag \pr B^\dag) \ox C^\dag                 \ar[dddd]^{\lambda_\pr\ox 1}
 & {}\\
{}
 & {}
 & (A\pr B)\pr C                      \ar[rr]^{\phi_{A\pr B}\pr \phi_{C}}     \ar[d]_{\mx^{-1} \oa 1}                     \ar@{}[ll]|{\mbox{\tiny {\bf mx ~ cat.}}}     \ar@{}[drr]|{\mbox{\tiny{\bf (id)$\ox$(Lem. \ref{lemma:cohUnitary}   (vii))} } }
 & {}                                                                                                                                              
 & (A\pr B)^\dag \pr C^\dag           \ar[r]^{\lambda_\ox^{-1}\pr 1}     \ar[d]_{\mx^\dag\pr 1}         \ar@{}[dr]|{\mbox{\tiny {\bf 1$\ox$(Lem. \ref{lemma: mixdagger})  }}}
 & (A^\dag \ox B^\dag) \pr C^\dag                \ar[d]^{\mx\pr 1}  \ar[r]^{\mx\pr 1}                          \ar@{}[ddr]|{\mbox{\tiny {\bf id}}}
 & (A^\dag \pr B^\dag) \pr C^\dag                \ar@=[dd]
 & {}
 & {}\\
{}
 & {}
 & (A\ox B)\oa C                       \ar[rr]^{\phi_{A\pr B}\oa \phi_{C}}    \ar@=[d]                                    
 & {}                                                                                                                                              \ar@{}[d]|{\mbox{\tiny {\bf [U.6(a)] }}}
 & (A\ox B)^\dag \oa C^\dag           \ar[r]^{\lambda_\ox^{-1}\pr 1}     \ar[d]_{\lambda_\ox}            \ar@{}[dr]|{\mbox{\tiny {\bf id}}}
 & (A^\dag\oa B^\dag) \oa C^\dag                \ar[d]^{\lambda_\oa \oa 1}
 & {}
 & {}
 & {}\\
{}
 & {}
 & (A\ox B)\oa C                       \ar[r]^{\mx^{-1}}     \ar@=[d]
 & (A\pr B)\ox C                       \ar[r]^{\phi_{(A\pr B)\pr C}}                                                            \ar@{}[dr]|{\mbox{\tiny {\bf [U.6(a)]}}}
 & ((A\ox B) \ox C)^\dag                \ar[r]^{\lambda_\ox^{-1}}
 & (A\ox B)^\dag \oa C^\dag           \ar[r]^{ \lambda_\oa^{-1} \ox 1}     \ar@=[d]                               \ar@{}[dr]|{\mbox{\tiny {\bf id}}}
 & (A^\dag \oa B^\dag) \oa C^\dag       \ar[d]^{\lambda_\oa \ox 1}
 & {}
 & {}\\
{}
 & {}
 & (A\ox B) \oa C                      \ar[rrr]^{\phi_{A\ox B} \oa \phi_{C}}                           \ar@/_4pc/[rrrrrr]_{\phi_{(A\ox B) \oa C}}
 & {}
 & {}
 & (A \ox B)^\dag \oa C^\dag          \ar@=[r]                                                                                   \ar@{}[d]|{\mbox{\tiny {\bf [U.6(a)]}}}
 & (A \ox B)^\dag \oa C^\dag          \ar[r]^{\mx^{-1}}
 & (A \ox B)^\dag \ox C^\dag          \ar[r]^{\lambda_\ox}
 &  ((A \ox B) \oa C)^\dag               \ar[uuuuuu]^{\partial_L^\dag}\\
{}
 & {}
 & {}
 & {}
 & {}
 & {}
 & {}
 & {}
 & {}\\
}
$
}

\end{sideways}
\caption{$\partial_L$ is a unitary isomorphism}
\label{Fig: linear dist. unitary}
\newpage
\end{figure}

\item $\partial_R$ is unitary because:

\[ \hspace{-1.25cm} 
\xymatrixrowsep{4pc}
\xymatrix{
{}
  & {}
  & {}
  & {}
  & {}
  & {}\\
(\!A\!\pr\! B\!) \!\ox\! C                           \ar[r]^{\mx} \ar[d]_{\partial_R}    \ar@/^4pc/[rrrrr]^{\phi_{(A\pr B) \ox C }}
  & (\!A\!\pr\! B\!) \!\pr\! C                      \ar[r]^{\mx^{-1}\pr 1}                                                        \ar@{}[d]|{\mbox{\tiny {\bf }}}
  & (\!A\!\ox\! B\!) \!\pr\! C                      \ar[r]^{\phi_{(A\ox B) \pr C}}                                              \ar@{}[dr]|{\mbox{\tiny {\bf Lem. \ref{lemma:cohUnitary}  (xiii)}}}  \ar@{}[ur]|{\mbox{\tiny {\bf Lem. \ref{lemma:cohUnitary} (vii), \ref{lemma:MUCProperties}}}}
  & (\!(\!A\!\ox\! B\!) \!\pr\! C \!)^\dag         \ar[r]^{(\mx^{-1}\pr 1 )^\dag} \ar[d]_{\partial_L^\dag}   
  & (\!(\!A\!\pr\! B\!) \!\pr\! C \!)^\dag         \ar[r]^{\mx^\dag}                                                               \ar@{}[d]|{\mbox{\tiny {\bf }}}
  & (\!(\!A\!\pr\! B\!) \!\ox\! C \!)^\dag \\
A\!\pr\! (\!B \!\ox\! C\!)                           \ar[r]^{\mx^{-1}}   \ar@/_4pc/[rrrrr]_{\phi_{A\pr (B \ox C) }}
  & A\!\pr\! (\!B \!\pr\! C\!)                      \ar[r]^{1 \pr \mx} 
  & A\!\ox\! (\!B \!\pr\! C\!)                      \ar[r]^{\phi_{A\ox (B \pr C)}} \ar[u]_{\partial_L}  \ar@{}[dr]|{\mbox{\tiny {\bf Lem. \ref{lemma:cohUnitary}  (vii), \ref{lemma:MUCProperties}}}}
  & (\!A\!\ox\! (\!B\!\pr\! C\!) \!)^\dag         \ar[r]^{(\mx \pr 1 )^\dag}
  & (\!A\!\pr\! (\!B \!\pr\! C\!) \!)^\dag         \ar[r]^{(\mx^{-1})^\dag}
  & (\!A\!\pr\! (\!B \!\ox\! C\!) \!)^\dag        \ar[u]_{\partial_R^\dag}\\
{}
  & {}
  & {}
  & {}
  & {}
  & {}
}
\]
\end{enumerate}
\end{proof}


\subsection{Unitary categories}
\label{subsection: unitary categories}

With the notion of unitary objects in place, one can consider $\dagger$-isomix categories in which all the objects are unitary: 
these are called {\em unitary categories\/}. This section develops the theory of unitary categories.

\begin{definition}
A {\bf (symmetric) unitary category} is a (symmetric) $\dagger$-isomix category with a unitary structure for which every object is unitary.
\end{definition}

Clearly, a unitary category must be a compact $\dagger$-LDC, since the mixor is a unitary isomorphism, see Lemma \ref{lemma:cohUnitary}-(iii).

 A $\dagger$-monoidal category is a strict unitary category in which the unitary structure map and the mix map are identity morphisms. Similarily, a $\dagger$-compact closed category is a strict unitary category in which all objects have unitary duals.

 In the rest of this subsection, we show that any unitary category is $\dagger$-linearly equivalent to a conventional dagger monoidal category. 
 A unitary category being a compact LDC is linearly equivalent, using ${\sf Mx}^*_\uparrow: (\X, \ox,\oa) \to 
 (\X,\oa,\oa)$ (see Corollary \ref{compact-mix-functor}) to the underlying monoidal category based on the 
 par (and the tensor). We now show that for a unitary category one can induce a stationary on objects dagger 
 on $(\X,\oa,\oa)$. We denote this dagger by $(\_)^\ddagger$ and define it by $f^\ddagger := 
 \varphi_Bf^\dagger\varphi_A^{-1}$ as illustrated by the left diagram below:
 
 \[ \xymatrix{ B \ar[d]_{\varphi_B}\ar[rr]^{f^\ddagger} \ar@{}[rrd]|{:=}& & A \ar[d]^{\varphi_A} \\
 	B^\dagger \ar[rr]_{f^\dagger} && A^\dagger} ~~~~~~~~~~~~~
 \xymatrix{ A \ar@/_1pc/[dd]_{\iota} \ar[d]^{\varphi_A}\ar[rr]^{f^{\ddagger\ddagger}} & & B \ar[d]_{\varphi_B} \ar@/^1pc/[dd]^{\iota} \\	
 	A^\dagger \ar[d]^{(\varphi_A^{-1})^\dagger} \ar[rr]_{(f^\ddagger)^\dagger} && B^\dagger \ar[d]_{(\varphi_B^{-1})^\dagger} \\ 	
 	A^{\dagger\dagger} \ar[rr]_{f^{\dagger\dagger}} & & B^{\dagger\dagger} } \]
 
 This new dagger clearly preserves composition and is also a stationary on objects involution as proven by the second diagram:
 the lower square of this diagram is the dagger of the inverted definition and the resulting outer square is the naturality of $\iota$ forcing $f^{\ddagger\ddagger} = f$.
 
 Next, we observe that $u: X \to Y$ is a unitary isomorphism in $\X$ if and only if $u^{-1}= u^\ddagger$.  This makes unitary isomorphisms in the traditional sense of categorical quantum mechanics coincide 
 with the notion introduced here.   Thus, $u$ is unitary in the sense here if and only if the diagram below commutes
 \[ \xymatrix{ B \ar[d]_{\varphi_B} \ar[rr]^{u^{-1}} && A  \ar[d]^{\varphi_A} \\ B^\dagger \ar[rr]_{u^\dagger} && A} \]
 but this diagram commutes if and only if $u^{-1} = u^\ddagger$. 
 
 \begin{definition} \label{preserving-unitary-structure}
 	A $\dagger$-Frobenius mix functor, $F: \X \to \Y$, between compact $\dagger$-isomix categories with unitary structure {\bf preserves unitary structure} if
 	\begin{enumerate}[(i)]
 		\item for all unitary objects $A \in \X$, $F(A)$ is a unitary object such that $\varphi_{F(A)} = F(\varphi_A) \rho^F$ 
 		\item Either $n_\bot^F$ or $m_\top^F$ are unitary isomorphisms i.e.,
 		\[
 		\xymatrix{
 			F(\bot) \ar[r]^{F(\varphi_\bot)} \ar[d]_{n_\bot} & F(\bot^\dagger) \ar[r]^{\rho} & F(\bot)^\dagger \\
 			\bot \ar[rr]_{\varphi_\bot} & & \bot^\dagger \ar[u]_{n_\bot^\dagger}
 		} (or)  \xymatrix{
 			\top \ar[rr]^{\varphi_\top} \ar[d]_{m_\top} & & \top^\dagger \\
 			F(\top) \ar[r]_{F(\varphi_\top)} & F(\top^\dagger) \ar[r]_{\rho} & F(\top)^\dagger \ar[u]_{m_\top^\dagger}
 		}
 		\]
 	\end{enumerate}
 \end{definition}
 
Notice that if $F$ preserves unitary structure, it must be an isomix functor by Lemma \ref{Lemma: isomix functor}. Also, when $A \in \X$ is a unitary object,  then $F(A)$ must be a unitary object, and so $F(A)$ is in the core.
 
 We now show that ${\sf Mx}_\uparrow: (\X,\oa,\oa) \to (\X,\ox,\oa)$ provides a unitary structure preserving equivalence of a dagger monoidal category into a unitary category:
 
 \begin{proposition}  \label{unitary-2-dagger}
 	Unitary categories are $\dagger$-linearly equivalent via the mix functor ${\sf Mx}_\uparrow: (\X,\oa,\oa) \to (\X,\ox,\oa)$ to the underlying dagger monoidal category on the par.  
 	Furthermore, closed unitary categories under this equivalence become dagger compact closed categories.
 \end{proposition}
 
 \begin{proof}
 	We must exhibit a preservator, that is a natural transformation showing that the involution is preserved:
 	\[ \infer={A \to_{\varphi_A} A^\dagger}{{\sf Mx}_\uparrow(A^\ddagger) \to^{\varphi_A} {\sf Mx}_\uparrow(A)^\dagger} \]
 	Note that $\varphi$ is a natural transformation by the definition of $(\_)^\ddagger$ and its coherence requirements 
 	make it a linear natural equivalence.  Making this the preservator immediately means that unitary structure is preserved.
 	
 	Finally, we must show that unitary linear duals under ${\sf Mx}^{*}_\uparrow$ become $\ddagger$-duals.  Given $(\eta,\epsilon): A \dashvv_u B$  we must 
 	show that under ${\sf Mx}^{*}_\uparrow$ this produces a dagger dual.  ${\sf Mx}^{*}_\uparrow(\eta) = {\sf m} ~\eta: \bot \to A \oa B$ and 
 	${\sf Mx}^{*}_\uparrow(\epsilon) = {\sf mx}^{-1} \epsilon: B \oa A \to \bot$
 	We then require that $c_\oa {\sf Mx}^{*}_\uparrow(\epsilon) = {\sf Mx}^{*}_\uparrow(\eta)^\ddagger$. This is provided by:
 	\[ \xymatrix{A \oa B  \ar[d]^{{\sf mx}^{-1}}  \ar@/_2pc/[ddd]_{\varphi_{A \oa B}} \ar[r]^{c_\oa} & B \oa A \ar@/^1pc/[rr]^{{\sf Mx}^{*}_\uparrow(\epsilon)} \ar[r]_{{\sf mx}^{-1}} 
 		& B \ox A \ar[r]_{\epsilon} & \bot \ar[dd]_{{\sf m}} \ar@/^2pc/[ddd]^{\varphi_\bot} \ar@/_/[dddl]_{\lambda_\bot}\\
 		A \ox B \ar[d]^{\varphi_A \ox \varphi_B} \ar@/_/[rru]^{c_\ox} \\
 		A^\dagger \ox B^\dagger \ar[d]^{\lambda_\ox} & & & \top \ar[d]_{\lambda_\top} \\
 		(A \oa B)^\dagger \ar@{}[rrruuu]|{{\rm Defn.} ~\ref{defn: unitary dual}~(b)} \ar@/_2pc/[rrr]_{{\sf Mx}^{*}_\uparrow(\eta)^\dagger} \ar[rr]_{\eta^\dagger} & & \top^\dagger \ar[r]_{{\sf m}^\dagger} & \bot^\dagger } \]
 \end{proof}

\subsection{The unitary construction}
\label{Sec: unitary construction}
A $\dagger$-isomix category can have many different unitary structures, as 
we shall describe in this section, 
thus it is {\em structure\/}, and not a property.   The requirements, however, 
do mean that for a $\dagger$-isomix category, $\X$, there is always a 
smallest unitary structure, referred to as the ``trivial'' unitary structure, 
that produces a full unitary subcategory in $\X$.  In this subsection, we 
provide a construction called the unitary construction which produces this unitary category from any 
$\dagger$-isomix category.   The construction is based on identifying 
objects with pre-unitary structure: the tensor units always have a 
canonical ``pre-unitary" structure so the construction always produces 
a non-empty category.   However, to ensure that an application of 
the construction yields a unitary category in which there are objects 
which are not isomorphic to the units, one must exhibit concretely 
such objects.  Fortunately this is often not difficult to do, making 
the construction quite applicable.

\begin{definition} ~
\label{defn: pre-unitary}
\begin{enumerate}[(i)] 
\item In a $\dagger$-isomix category, a {\bf pre-unitary object} is an object $U \in \Core(\X)$, together with an isomorphism $\alpha: U \to U^\dagger$ such that  $\alpha (\alpha^{-1})^\dagger = \iota$. 

\item Suppose $\X$ is a $\dagger$-isomix category, then define ${\sf Unitary}(\X)$, the {\bf canonical unitary core} of $\X$, as follows:
\begin{description}
\item[Objects:] Pre-unitary objects $(U, \alpha)$,
\item[Maps:] $(U, \alpha) \to^f (V, \beta)$ where $U \to^f V $ is any map of $\X$.
\end{description}
\end{enumerate}
\end{definition}

We note that any object which is isomorphic to a preunitary object is also pre-unitary:
\begin{lemma}
In a $\dagger$-isomix category, if $U$ is a pre-unitary object 
and there exists an isomorphism $f: U \to U'$, then $U'$ is pre-unitary. 
\end{lemma}

Our objective is to show that Unitary($\X$) is endowed with all the structure of a unitary category.

\begin{lemma}
For any $\dagger$-isomix category, its canonical unitary core is a compact $\dagger$-LDC with tensor and par defined by
\[ (\top,{\sf m}^{-1}\lambda_\bot: \top \to \top^\dagger) ~~~~~(A, \alpha) \ox (B, \beta) := (A \ox B, \mx(\alpha \oa \beta) \lambda_\oa: A \ox B \to (A \ox B)^\dagger)\]
\[ (\bot, {\sf m} ~\lambda_\top: \bot \to \bot^\dagger)  ~~~~~(A, \alpha) \oa (B, \beta) := (A \oa B, \mx^{-1}(\alpha \ox \beta) \lambda_\ox: A \oa B \to (A \oa B)^\dagger) \]
and $(U,\alpha)^\dagger := (U^\dagger, (\alpha^{-1})^\dagger)$.
\end{lemma}

\begin{proof}
The proof uses the techniques of Lemma \ref{Lemma: square root tensor unitary}. 

 Note that, as the map and tensor structure is inherited from $\X$, it suffices to show that these objects are all pre-unitary objects.  Starting with $(U \alpha)^\dagger$ 
we have:
\[ (\alpha^{-1})^\dagger (((\alpha^{-1})^\dagger)^{-1})^\dagger = (\alpha^{-1})^\dagger (\alpha^\dagger)^\dagger = (\alpha^\dagger \alpha^{-1})^\dagger = (\iota^{-1})^\dagger = \iota \]
For the tensor and par we have:
\begin{eqnarray*}
{\sf m}^{-1}\lambda_\bot (({\sf m}^{-1}\lambda_\bot)^{-1})^\dagger & = & {\sf m}^{-1}\lambda_\bot {\sf m}^\dagger \lambda_\bot^\dagger \\
& \stackrel{\text{\tiny {\bf [$\dagger$-mix]}}}{=} &  {\sf m}^{-1} {\sf m} \lambda_\top  \lambda_\bot^\dagger  = \iota \\
\mx^{-1}(\alpha \oa \beta) \lambda_\oa ((\mx^{-1}(\alpha \oa \beta) \lambda_\oa)^{-1})^\dagger 
& = & \mx^{-1}(\alpha \oa \beta) \lambda_\oa (\mx^\dagger) (\alpha^{-1} \oa \beta^{-1})^\dagger (\lambda_\oa^{-1})^\dagger \\
& = & \mx^{-1}(\alpha \oa \beta) \mx \lambda_\ox (\alpha^{-1} \oa \beta^{-1})^\dagger (\lambda_\oa^{-1})^\dagger \\
& = & (\alpha \ox \beta) \lambda_\ox (\alpha^{-1} \oa \beta^{-1})^\dagger (\lambda_\oa^{-1})^\dagger \\
& = & (\alpha \ox \beta) ((\alpha^{-1})^\dagger \ox (\beta^{-1})^\dagger) \lambda_\ox (\lambda_\oa^{-1})^\dagger \\
&\stackrel{\text{\tiny {\bf Defn \ref{defn: pre-unitary}-(i)}}}{=}& (\iota \ox \iota) \lambda_\ox (\lambda_\oa^{-1})^\dagger \\
&\stackrel{\text{\tiny {\bf [$\dagger$-ldc.4]}}}{=}& \iota 
\end{eqnarray*}
\begin{eqnarray*}
{\sf m} \lambda_\top (({\sf m} \lambda_\top)^{-1})^\dagger & = & {\sf m} \lambda_\top ({\sf m}^{-1})^\dagger  (\lambda_\top^{-1})^\dagger  \\
& = & {\sf m} ~{\sf m}^{-1} \lambda_\bot (\lambda_\top^{-1})^\dagger = \iota \\
\mx(\alpha \ox \beta) \lambda_\ox ((\mx(\alpha \ox \beta) \lambda_\ox)^{-1})^\dagger 
& = & \mx(\alpha \ox \beta) \lambda_\ox (\mx^{-1})^\dagger (\alpha^{-1} \ox \beta^{-1})^\dagger (\lambda_\ox^{-1})^\dagger \\
& = & (\alpha \oa \beta) \mx~ \mx^{-1} \lambda_\oa (\alpha^{-1} \ox \beta^{-1})^\dagger (\lambda_\ox^{-1})^\dagger \\
&= & (\alpha \oa \beta) ((\alpha^{-1})^\dagger \oa (\beta^{-1})^\dagger) \lambda_\oa (\lambda_\ox^{-1})^\dagger \\
& = & (\iota \oa \iota) \lambda_\oa (\lambda_\ox^{-1})^\dagger  = \iota.
\end{eqnarray*}
\end{proof}

This makes $\Unitary(\X)$ into a compact $\dagger$-LDC with all the structure inherited directly from $\X$. 
However, more is true: each object now has an obvious unitary structure.  This gives:

\begin{proposition}
For any $\dagger$-isomix category, $\X$, $\Unitary(\X)$ is a unitary category with a full and faithful underlying $\dagger$-isomix functor $U: {\sf Unitary}(\X) \to \X$.
\end{proposition}

\begin{proof}
The laxors are all identity maps so that the underlying functors is immediately a $\dagger$-mix functor.

It remains to show that every object is unitary:  we set the unitary structure of an object to be $\alpha: (X,\alpha) \to (X,\alpha)^\dagger$.   However, {\bf [U.1]} -- {\bf [U.5]} are immediately satisfied by construction implying this provides unitary structure for every object.
\end{proof}

Next, we prove the couniversal property of the unitary construction. Define ${\sf UCat}$ to be the category of unitary categories and $\dagger$-isomix functors that preserve   unitary structure in the sense of Definition \ref{preserving-unitary-structure}, thus, whenever ${\varphi_A}$ is the unitary structure  then $F'(\varphi_A) \rho^{F'}$ is unitary structure. Define {\sf Kompact} to be the category of compact $\dagger$-LDCs and $\dagger$-isomix functors.

We now show that the unitary construction produces a right adjoint to the underlying functor $U: {\sf UCat} \to {\sf Kompact}$ which is the identity functor. Preliminary to this result we prove that Frobenius functors preserve preunitary objects:

\begin{lemma}
\label{Lemma: Frobenius preunitary}
If $F: \X \to \Y$ is a $\dagger$-isomix functor between compact $\dagger$-LDCs and $(A,\varphi)$ is a preunitary object of $\X$, then $(F(A),F(\varphi)\rho)$ is a preunitary object of $\Y$.
\end{lemma}
\begin{proof}
To prove that $(F(A),F(\varphi)\rho)$ is a preunitary object, one has the following computation:
\begin{eqnarray*}
F(\varphi)\rho ((F(\varphi) \rho)^{-1})^\dagger 
& = & F(\varphi)\rho F(\varphi^{-1})^\dagger (\rho^{-1})^\dagger \\
& = & F(\varphi (\varphi^{-1})^\dagger) \rho (\rho^{-1})^\dagger \\
& = & F(\iota) \rho (\rho^{-1})^\dagger \stackrel{{\bf [\dagger-isomix]}}{=} \iota.
\end{eqnarray*}
\end{proof}

\begin{proposition}
\label{Prop: Couniversal}
$U: {\sf UCat} \to {\sf Kompact}$ has a right adjoint ${\sf Unitary}: {\sf Kompact} \to {\sf UCat}; \C \mapsto {\sf Unitary}(\C)$.
\end{proposition}
\begin{proof}
The couniversal diagram is as follows:
\[ \xymatrix{ U(\U) \ar[rr]^{F} \ar@{.>}[d]_{U(F^\flat)} && \C \\
                    U({\sf Unitary}(\C)) \ar[urr]_{\epsilon}} \]

Since $F$ is a $\dagger$-isomix functor it preserves preunitary structure (see Lemma \ref{Lemma: Frobenius preunitary}).  This means that each $(U,\varphi_U)$ in $\U$ is carried by $F$ onto a preunitary object in $\C$, $(F(U),F(\varphi)\rho^F)$.  But a preunitary object in $\C$ is an object of ${\sf Unitary}(\C)$ and this determines $F^\flat$.   The functor $F^\flat$ is uniquely determined as it must preserve the unitary structure.
\end{proof}

\section{Examples:  The unitary construction}
\label{Sec: Examples The unitary construction}

In Section \ref{daggers-duals-conjugation}, we discussed examples of $\dagger$-isomix categories in which the $\dagger$ 
is given by composing the conjugation functor and the dualizing functor. In the rest of the section, we apply the unitary 
construction to each of those examples to construct a unitary category:


\subsection{Category of abstract state spaces}
In Section \ref{Sec: Asp}, we discussed a construction on a $\dagger$-isomix category, $\X$, that produces a 
category of abstract state spaces, $\Asp(\X)$, which is a $\dagger$-isomix category. In this section, we examine 
the preunitary objects of $\Asp(\X)$. Since all the basic natural isomorphisms are inherited from $\X$, $\Core(\X)$ 
determines $\Core(\Asp(\X))$. If $(A ,\alpha)$ is a preunitary object for $\X$, and $(A, e_A, u_A) \in \Asp(\X)$ then, 
$((A, e_A, u_A), \alpha)$ is a preunitary object for $\Asp(\X)$ if $u_A \alpha = \lambda_\top e_A^\dagger$.


\subsection{Category of a group with involution}
We discussed a source of examples of compact $\dagger$-LDCs which are given by groups with conjugation. Applying 
unitary construction to each of the example categories results in the following unitary categories. It could be noticed 
that the preunitary objects in each of these categories includes those group elements such that $\overline{g^{-1}} = g$. 
More explicitly, the preunitary objects are $(g,1)$ such that $\overline{g^{-1}} = g$.

\begin{itemize}
    \item In the discrete category of complex numbers, $\D(\C, +, 0)$, \[(a + ib)^\dagger := \overline{(a+ib)^*} = \overline{(-a-ib)} = 
    -a + ib\] The preunitary objects in this category are given by all complex numbers, i.e., $(ib, 1)$. 
    
    \item In the discrete category of non-zero complex numbers, $\D(\C, ., 1)$, the preunitary objects are given by complex 
    numbers on a unit circle.
    
    \item In the discrete category, $\D(P(x), +, 0)$, where $P(x)$ is a polynomial ring, $P(x)^\dagger = -P(-x)$ and the preunitary 
    objects are polynomials $ P(x) = \sum_n a_n x^n$ such that n is odd. 
    
    \item In $\D(\mathbb{M}_2, \cdot, I_2)$ where $\mathbb{M}_2$ is the group of $2 \times 2$ invertible matrices over 
    $\mathbb{C}$. The $\dagger$ structure is as follows:
     \[\left(
      \begin{matrix}
     a+ib & m+in \\
     c+id & p+iq 
     \end{matrix}
     \right)^\dagger := \overline{\left(
     \begin{matrix}
     a+ib & m+in \\
     c+id & p+iq 
     \end{matrix}
     \right)^*} = \left(
     \begin{matrix}
    a-ib & c-id \\
    m-in & p-iq
     \end{matrix}
     \right)^{-1} \]
     The preunitary objects in this category are the unitary matrices.
    \end{itemize}
\subsection{Category of Hopf modules in a $*$-automonous category}
In Section \ref{Sec: HModx}, we described a construction of $\dagger$-isomix categories  from any  symmetric isomix 
$*$-autonomous category, $\X$, by choosing the Hopf Modules over a  cocommutative $\ox$-Hopf Algebra. 
We referred to the resulting category as ${\mbox{\bf H-Mod}}_\X$. Now we shall look at the preunitary objects in 
${\mbox{\bf H-Mod}}_\X$ in order to apply the unitary construction to this category. We begin by identifying the 
objects in the core of ${\mbox{\bf H-Mod}}_\X$:


\begin{lemma}
Suppose $\X$ is a mix $*$-autonomous category and $H$ is a cocommutative Hopf Algebra in $\X$. If $(A, \leftaction{0.4}{white})$ is a H-Module and $A \in \Core(\X)$, then $ (A, \leftaction{0.4}{white}) \in \Core({\mbox{\bf H-Mod}}_\X)$.
\end{lemma}
\begin{proof}
The mixor $\mx: A \ox B \to A \oa B$ is inherited directly from $\X$. Henc,e  $ (A, \leftaction{0.4}{white}) \in \Core(${\bf H-Mod$_\X)$}.
\end{proof}

Now that we identified the objects in the core, we prove a lemma that will be used later to identify the preunitary objects from the core:

\begin{lemma}
\label{Lemma: aux} The following equality holds for a Frobenius algebra: 
\[ \begin{tikzpicture}
	\begin{pgfonlayer}{nodelayer}
		\node [style=circle] (0) at (-2.5, 0.5) {};
		\node [style=none] (1) at (-2, 1) {};
		\node [style=none] (3) at (-3, 1) {};
		\node [style=none] (4) at (-2.5, 0) {};
		\node [style=none] (5) at (-1.5, 0) {};
		\node [style=circle] (6) at (-2, -1) {};
		\node [style=circle] (7) at (-2, -0.5) {};
		\node [style=none] (8) at (-3, 1) {};
		\node [style=none] (9) at (-4, 1) {};
		\node [style=circle] (10) at (-3.5, 2.25) {};
		\node [style=circle] (11) at (-3.5, 1.75) {};
		\node [style=circle] (12) at (-3.5, 2.75) {};
		\node [style=none] (13) at (-2, 1) {};
		\node [style=none] (14) at (-4.75, 1) {};
		\node [style=circle] (15) at (-3.5, 3.25) {};
		\node [style=none] (16) at (-1.5, 3.5) {};
		\node [style=none] (17) at (-4, -1) {};
		\node [style=none] (18) at (-4.75, -1) {};
	\end{pgfonlayer}
	\begin{pgfonlayer}{edgelayer}
		\draw [in=-90, out=30, looseness=1.25] (0) to (1.center);
		\draw [in=150, out=-90] (3.center) to (0);
		\draw [in=-90, out=45] (7) to (5.center);
		\draw (7) to (6);
		\draw [in=127, out=-90, looseness=0.75] (4.center) to (7);
		\draw [in=90, out=-30] (11) to (8.center);
		\draw (11) to (10);
		\draw [in=-150, out=90] (9.center) to (11);
		\draw [in=90, out=-30, looseness=0.75] (12) to (13.center);
		\draw (12) to (15);
		\draw [in=-150, out=90, looseness=0.75] (14.center) to (12);
		\draw (0) to (4.center);
		\draw (5.center) to (16.center);
		\draw (17.center) to (9.center);
		\draw (18.center) to (14.center);
	\end{pgfonlayer}
\end{tikzpicture}
 = \begin{tikzpicture}
	\begin{pgfonlayer}{nodelayer}
		\node [style=none] (19) at (1.75, -1) {};
		\node [style=circle] (20) at (1, 1.25) {};
		\node [style=none] (21) at (0.25, -1) {};
		\node [style=none] (22) at (1, 3.5) {};
	\end{pgfonlayer}
	\begin{pgfonlayer}{edgelayer}
		\draw (20) to (22.center);
		\draw [in=-150, out=90] (21.center) to (20);
		\draw [in=90, out=-30] (20) to (19.center);
	\end{pgfonlayer}
\end{tikzpicture} \]
\end{lemma}
\begin{proof}
\[\begin{tikzpicture}
	\begin{pgfonlayer}{nodelayer}
		\node [style=circle] (0) at (-2.5, 0.75) {};
		\node [style=none] (1) at (-2, 1.25) {};
		\node [style=none] (2) at (-2.5, 0.25) {};
		\node [style=none] (3) at (-3, 1.25) {};
		\node [style=none] (4) at (-2.5, 0.25) {};
		\node [style=none] (5) at (-1.5, 0.25) {};
		\node [style=circle] (6) at (-2, -1) {};
		\node [style=circle] (7) at (-2, -0.25) {};
		\node [style=none] (8) at (-3, 1.25) {};
		\node [style=none] (9) at (-4, 1.25) {};
		\node [style=circle] (10) at (-3.5, 2.5) {};
		\node [style=circle] (11) at (-3.5, 1.75) {};
		\node [style=circle] (12) at (-3.5, 3.25) {};
		\node [style=none] (13) at (-2, 1.25) {};
		\node [style=none] (14) at (-4.75, 1.25) {};
		\node [style=circle] (15) at (-3.5, 4) {};
		\node [style=none] (16) at (-1.5, 4) {};
		\node [style=none] (17) at (-4, -1) {};
		\node [style=none] (18) at (-4.75, -1) {};
	\end{pgfonlayer}
	\begin{pgfonlayer}{edgelayer}
		\draw [in=-90, out=30] (0) to (1.center);
		\draw (0) to (2.center);
		\draw [in=150, out=-90] (3.center) to (0);
		\draw [in=-90, out=60, looseness=0.75] (7) to (5.center);
		\draw (7) to (6);
		\draw [in=127, out=-90, looseness=0.75] (4.center) to (7);
		\draw [in=90, out=-15] (11) to (8.center);
		\draw (11) to (10);
		\draw [in=-165, out=90] (9.center) to (11);
		\draw [in=90, out=-30, looseness=0.75] (12) to (13.center);
		\draw (12) to (15);
		\draw [in=-150, out=90, looseness=0.75] (14.center) to (12);
		\draw (16.center) to (5.center);
		\draw (9.center) to (17.center);
		\draw (14.center) to (18.center);
	\end{pgfonlayer}
\end{tikzpicture}  = \begin{tikzpicture}
	\begin{pgfonlayer}{nodelayer}
		\node [style=none] (19) at (3.25, 1) {};
		\node [style=circle] (20) at (2.75, -1) {};
		\node [style=circle] (21) at (2.75, -0.25) {};
		\node [style=circle] (22) at (1, 3.25) {};
		\node [style=none] (23) at (0, -1) {};
		\node [style=circle] (24) at (1, 4) {};
		\node [style=circle] (25) at (1, 2.5) {};
		\node [style=none] (26) at (2, 2.5) {};
		\node [style=circle] (27) at (1.5, 1.75) {};
		\node [style=none] (28) at (2, 2.5) {};
		\node [style=none] (30) at (1, -1) {};
		\node [style=circle] (32) at (1.5, 1) {};
		\node [style=none] (33) at (3.25, 4) {};
		\node [style=none] (34) at (3.25, 4) {};
	\end{pgfonlayer}
	\begin{pgfonlayer}{edgelayer}
		\draw [in=-90, out=45] (21) to (19.center);
		\draw (21) to (20);
		\draw (22) to (24);
		\draw [in=-150, out=90, looseness=0.75] (23.center) to (22);
		\draw [in=-90, out=30] (27) to (28.center);
		\draw [in=150, out=-90, looseness=1.25] (25) to (27);
		\draw [in=-127, out=90, looseness=0.75] (30.center) to (32);
		\draw [in=90, out=0, looseness=1.25] (22) to (26.center);
		\draw (27) to (32);
		\draw (32) to (21);
		\draw (34.center) to (19.center);
	\end{pgfonlayer}
\end{tikzpicture}  =
\begin{tikzpicture}
	\begin{pgfonlayer}{nodelayer}
		\node [style=none] (35) at (7.5, 2) {};
		\node [style=circle] (36) at (7, -1) {};
		\node [style=circle] (37) at (7, 1) {};
		\node [style=circle] (38) at (5.75, 3.25) {};
		\node [style=none] (39) at (5, 2) {};
		\node [style=circle] (40) at (5.75, 4) {};
		\node [style=none] (41) at (5.75, 1) {};
		\node [style=circle] (43) at (6.25, 2) {};
		\node [style=none] (44) at (5, -1) {};
		\node [style=none] (45) at (5.75, -1) {};
		\node [style=none] (46) at (7.5, 4) {};
	\end{pgfonlayer}
	\begin{pgfonlayer}{edgelayer}
		\draw [in=-90, out=30] (37) to (35.center);
		\draw (37) to (36);
		\draw (38) to (40);
		\draw [in=-165, out=90] (39.center) to (38);
		\draw [in=-127, out=90] (41.center) to (43);
		\draw [in=90, out=-30, looseness=1.25] (38) to (43);
		\draw (43) to (37);
		\draw (44.center) to (39.center);
		\draw (45.center) to (41.center);
		\draw (35.center) to (46.center);
	\end{pgfonlayer}
\end{tikzpicture} = \begin{tikzpicture}
	\begin{pgfonlayer}{nodelayer}
		\node [style=none] (47) at (11.5, 4) {};
		\node [style=none] (48) at (10.75, -1) {};
		\node [style=circle] (49) at (10.75, 1) {};
		\node [style=circle] (50) at (10, 2) {};
		\node [style=none] (51) at (9.25, -1) {};
		\node [style=circle] (52) at (10, 4) {};
	\end{pgfonlayer}
	\begin{pgfonlayer}{edgelayer}
		\draw [in=-90, out=30, looseness=0.75] (49) to (47.center);
		\draw (49) to (48.center);
		\draw (50) to (52);
		\draw [in=-150, out=90, looseness=0.75] (51.center) to (50);
		\draw (50) to (49);
	\end{pgfonlayer}
\end{tikzpicture} = \begin{tikzpicture}
	\begin{pgfonlayer}{nodelayer}
		\node [style=none] (0) at (3, -1) {};
		\node [style=circle] (1) at (2, 2) {};
		\node [style=none] (2) at (1, -1) {};
		\node [style=none] (3) at (2, 4) {};
	\end{pgfonlayer}
	\begin{pgfonlayer}{edgelayer}
		\draw (1) to (3.center);
		\draw [in=-150, out=90, looseness=0.75] (2.center) to (1);
		\draw [in=90, out=-30, looseness=0.75] (1) to (0.center);
	\end{pgfonlayer}
\end{tikzpicture} \]
\end{proof}

In the following Proposition we identify the preunitary objects in the core:

\begin{proposition} 
Suppose $\X$ is a symmetric mix $*$-autonomous category and $H$ is a cocommutative Hopf Algebra in $\X$. If $A \in 
\Core(\X)$ and $(A, \mulmap{1.2}{white}, \unitmap{1.2}{white}, \comulmap{1.2}{white}, \counitmap{1.2}{white})$ 
is a cocommutative Frobenius Algebra with an algebra homomorphism $H \to^{h} A$ then, 
\begin{enumerate}[(a)]
\item $(A, \leftaction{0.4}{white})$ is a H-Module where, $\leftaction{0.4}{white}: H \ox A \to A := \begin{tikzpicture}
	\begin{pgfonlayer}{nodelayer}
		\node [style=circle] (0) at (0, -0) {};
		\node [style=none] (1) at (-0.5, 1.75) {};
		\node [style=none] (2) at (0.75, 1.75) {};
		\node [style=none] (3) at (0, -0.75) {};
		\node [style=circle, scale=2] (4) at (-0.5, 1) {};
		\node [style=none] (5) at (-0.5, 1) {$h$};
	\end{pgfonlayer}
	\begin{pgfonlayer}{edgelayer}
		\draw (1.center) to (4);
		\draw [bend left, looseness=1.00] (0) to (4);
		\draw [in=-90, out=15, looseness=1.00] (0) to (2.center);
		\draw (0) to (3.center);
	\end{pgfonlayer}
\end{tikzpicture}$ 
 \item $\overline{(A, \leftaction{0.4}{white})^*} = (A, \leftaction{0.4}{white})$ where $A^*$ is the self-dual Frobenius 
 Algebra with cups and caps defined as
$
\begin{tikzpicture}
	\begin{pgfonlayer}{nodelayer}
		\node [style=circle] (0) at (0, -0) {};
		\node [style=none] (1) at (-0.5, 1) {};
		\node [style=none] (2) at (0.75, 1) {};
		\node [style=circle] (3) at (0, -1) {};
		\node [style=circle] (4) at (0, -1) {};
	\end{pgfonlayer}
	\begin{pgfonlayer}{edgelayer}
		\draw [in=-90, out=15, looseness=1.00] (0) to (2.center);
		\draw [in=150, out=-90, looseness=1.00] (1.center) to (0);
		\draw (0) to (3);
	\end{pgfonlayer}
\end{tikzpicture}  and  \begin{tikzpicture}
	\begin{pgfonlayer}{nodelayer}
		\node [style=circle] (0) at (0, 0) {};
		\node [style=none] (1) at (-0.5, -1) {};
		\node [style=none] (2) at (0.75, -1) {};
		\node [style=circle] (3) at (0, 1) {};
		\node [style=circle] (4) at (0, 1) {};
	\end{pgfonlayer}
	\begin{pgfonlayer}{edgelayer}
		\draw [in=90, out=-15, looseness=1.00] (0) to (2.center);
		\draw [in=-150, out=90, looseness=1.00] (1.center) to (0);
		\draw (0) to (3);
	\end{pgfonlayer}
\end{tikzpicture}
$ respectively. Hence, $A^* = A$ and $(A, \leftaction{0.4}{white})^\dagger =  (A, \leftaction{0.4}{white})$.
\end{enumerate}
\end{proposition}
\begin{proof}~
\begin{enumerate}[(a)]
\item $\begin{tikzpicture} 
	\begin{pgfonlayer}{nodelayer}
		\node [style=circle] (0) at (0, -0) {};
		\node [style=none] (1) at (0.75, 1.25) {};
		\node [style=none] (2) at (0, -0.75) {};
		\node [style=circle, scale=2] (3) at (-0.5, 0.5) {};
		\node [style=none] (4) at (-0.5, 0.5) {$h$};
		\node [style=none] (5) at (-0.5, 1.25) {};
	\end{pgfonlayer}
	\begin{pgfonlayer}{edgelayer}
		\draw [bend left, looseness=1.00] (0) to (3);
		\draw [in=-90, out=15, looseness=1.00] (0) to (1.center);
		\draw (0) to (2.center);
		\draw (5.center) to (3);
	\end{pgfonlayer}
\end{tikzpicture}: H \ox A \to A$ is a left action because $h: H \to A$ is an algebra homomorphism.
\item $
\begin{tikzpicture} 
	\begin{pgfonlayer}{nodelayer}
		\node [style=none] (0) at (1.25, -1) {};
		\node [style=none] (1) at (2.25, -1) {};
		\node [style=none] (2) at (1.5, 0.5) {};
		\node [style=none] (3) at (1.5, 1) {};
		\node [style=none] (4) at (0.25, 1) {};
		\node [style=none] (5) at (2.25, 3) {};
		\node [style=none] (6) at (0.25, -1.25) {};
		\node [style=none] (7) at (1, 2.75) {$H$};
		\node [style=none] (8) at (2.5, 2.75) {$A^*$};
		\node [style=none] (9) at (0.25, -3) {};
		\node [style=none] (10) at (0.5, -2.75) {$A^*$};
		\node [style=circle, scale=1.5] (11) at (0.75, 0.25) {};
		\node [style=none] (12) at (0.75, 3) {};
		\node [style=circle] (13) at (1.25, -0.5) {};
		\node [style=none] (14) at (0.75, 0.25) {$h$};
	\end{pgfonlayer}
	\begin{pgfonlayer}{edgelayer}
		\draw [in=-90, out=90, looseness=1.00] (2.center) to (3.center);
		\draw [bend left=90, looseness=2.75] (4.center) to (3.center);
		\draw [bend right=90, looseness=2.00] (0.center) to (1.center);
		\draw (5.center) to (1.center);
		\draw (4.center) to (6.center);
		\draw [bend right, looseness=1.00] (11) to (13);
		\draw [bend right=15, looseness=1.00] (13) to (2.center);
		\draw (13) to (0.center);
		\draw (12.center) to (11);
		\draw (6.center) to (9.center);
	\end{pgfonlayer}
\end{tikzpicture} =
\begin{tikzpicture} 
	\begin{pgfonlayer}{nodelayer}
		\node [style=none] (0) at (1.25, -1) {};
		\node [style=none] (1) at (2.25, -1) {};
		\node [style=none] (2) at (1.5, 0.5) {};
		\node [style=none] (3) at (0.25, 0.25) {};
		\node [style=none] (4) at (1.5, 1) {};
		\node [style=none] (5) at (-0.25, 1) {};
		\node [style=none] (6) at (2.25, 3) {};
		\node [style=none] (7) at (0.25, -1.25) {};
		\node [style=none] (8) at (-0.25, -1.25) {};
		\node [style=none] (9) at (-0.5, -2) {};
		\node [style=none] (10) at (0.25, -2) {};
		\node [style=none] (11) at (-1.5, -2) {};
		\node [style=none] (12) at (-1.5, 3) {};
		\node [style=none] (13) at (0.25, -3) {};
		\node [style=none] (14) at (0.75, 0.25) {};
		\node [style=circle] (15) at (1.25, -0.5) {};
		\node [style=circle, scale=1.5] (16) at (-1.5, -0.75) {};
		\node [style=none] (17) at (-0.25, 1) {};
		\node [style=none] (18) at (1.5, 1) {};
		\node [style=circle] (19) at (0.5, 2.75) {};
		\node [style=circle] (20) at (0.5, 2) {};
		\node [style=none] (21) at (0.25, 0.25) {};
		\node [style=circle] (22) at (0.5, 1.5) {};
		\node [style=none] (23) at (0.75, 0.25) {};
		\node [style=circle] (24) at (0.5, 0.75) {};
		\node [style=none] (25) at (1.25, -1) {};
		\node [style=circle] (26) at (1.75, -2.75) {};
		\node [style=none] (27) at (2.25, -1) {};
		\node [style=circle] (28) at (1.75, -2) {};
		\node [style=none] (29) at (-1.5, -0.75) {$h$};
	\end{pgfonlayer}
	\begin{pgfonlayer}{edgelayer}
		\draw [in=-90, out=90, looseness=1.00] (2.center) to (4.center);
		\draw (6.center) to (1.center);
		\draw (3.center) to (7.center);
		\draw (5.center) to (8.center);
		\draw [in=105, out=-90, looseness=1.25] (8.center) to (10.center);
		\draw [in=90, out=-75, looseness=0.75] (7.center) to (9.center);
		\draw [bend right=90, looseness=1.25] (11.center) to (9.center);
		\draw (13.center) to (10.center);
		\draw [bend right=15, looseness=1.00] (15) to (2.center);
		\draw (15) to (0.center);
		\draw (16) to (11.center);
		\draw [in=-90, out=135, looseness=1.00] (15) to (14.center);
		\draw [bend right, looseness=1.00] (20) to (17.center);
		\draw [bend left, looseness=1.00] (20) to (18.center);
		\draw (19) to (20);
		\draw [bend right, looseness=1.00] (24) to (21.center);
		\draw [bend left, looseness=1.00] (24) to (23.center);
		\draw (22) to (24);
		\draw [bend left, looseness=1.00] (28) to (25.center);
		\draw [bend right, looseness=1.00] (28) to (27.center);
		\draw (26) to (28);
		\draw (12.center) to (16);
	\end{pgfonlayer}
\end{tikzpicture} \stackrel{Lemma ~ \ref{Lemma: aux}}{=} 
\begin{tikzpicture} 
	\begin{pgfonlayer}{nodelayer}
		\node [style=none] (0) at (0.75, 0.25) {};
		\node [style=none] (1) at (-0.25, 0.25) {};
		\node [style=none] (2) at (0.25, -2) {};
		\node [style=none] (3) at (-1.5, -2) {};
		\node [style=none] (4) at (-1.5, 1.75) {};
		\node [style=none] (5) at (0.25, -3.5) {};
		\node [style=circle, scale=1.5] (6) at (-1.5, 0.25) {};
		\node [style=none] (7) at (-1.5, 0.25) {$h$};
		\node [style=none] (8) at (0.75, 0.25) {};
		\node [style=none] (9) at (0.25, 1.75) {};
		\node [style=circle] (10) at (0.25, 1) {};
		\node [style=none] (11) at (-0.25, 0.25) {};
		\node [style=none] (12) at (-1.5, -2) {};
		\node [style=none] (13) at (-0.5, -2) {};
		\node [style=circle] (14) at (-1, -3.25) {};
		\node [style=circle] (15) at (-1, -2.5) {};
		\node [style=none] (16) at (-1.5, -2) {};
	\end{pgfonlayer}
	\begin{pgfonlayer}{edgelayer}
		\draw (5.center) to (2.center);
		\draw (6) to (3.center);
		\draw (4.center) to (6);
		\draw (7.center) to (6);
		\draw [bend right, looseness=1.00] (10) to (11.center);
		\draw [bend left, looseness=1.00] (10) to (8.center);
		\draw (9.center) to (10);
		\draw [bend left, looseness=1.00] (15) to (16.center);
		\draw [bend right, looseness=1.00] (15) to (13.center);
		\draw (14) to (15);
		\draw [in=90, out=-90, looseness=1.00] (1.center) to (2.center);
		\draw [in=90, out=-90, looseness=1.00] (0.center) to (13.center);
	\end{pgfonlayer}
\end{tikzpicture} \stackrel{\text{cocomm.}}{=} 
\begin{tikzpicture} 
	\begin{pgfonlayer}{nodelayer}
		\node [style=none] (0) at (0.75, 0.25) {};
		\node [style=none] (1) at (-0.25, 0.25) {};
		\node [style=none] (2) at (0.75, -2) {};
		\node [style=none] (3) at (-1.5, -2) {};
		\node [style=none] (4) at (-1.5, 1.75) {};
		\node [style=none] (5) at (0.75, -3.5) {};
		\node [style=circle, scale=1.5] (6) at (-1.5, 0.25) {};
		\node [style=none] (7) at (-1.5, 0.25) {$h$};
		\node [style=none] (8) at (0.75, 0.25) {};
		\node [style=none] (9) at (0.25, 1.75) {};
		\node [style=circle] (10) at (0.25, 1) {};
		\node [style=none] (11) at (-0.25, 0.25) {};
		\node [style=none] (12) at (-1.5, -2) {};
		\node [style=none] (13) at (-0.5, -2) {};
		\node [style=circle] (14) at (-1, -3.25) {};
		\node [style=circle] (15) at (-1, -2.5) {};
		\node [style=none] (16) at (-1.5, -2) {};
	\end{pgfonlayer}
	\begin{pgfonlayer}{edgelayer}
		\draw (5.center) to (2.center);
		\draw (6) to (3.center);
		\draw (4.center) to (6);
		\draw (7.center) to (6);
		\draw [bend right, looseness=1.00] (10) to (11.center);
		\draw [bend left, looseness=1.00] (10) to (8.center);
		\draw (9.center) to (10);
		\draw [bend left, looseness=1.00] (15) to (16.center);
		\draw [bend right, looseness=1.00] (15) to (13.center);
		\draw (14) to (15);
		\draw (1.center) to (13.center);
		\draw (0.center) to (2.center);
	\end{pgfonlayer}
\end{tikzpicture} =
\begin{tikzpicture}
	\begin{pgfonlayer}{nodelayer}
		\node [style=circle] (17) at (3.75, -1.5) {};
		\node [style=none] (18) at (4.25, 0.5) {};
		\node [style=none] (19) at (3.75, -2.25) {};
		\node [style=circle, scale=2] (20) at (3.25, -0.5) {};
		\node [style=none] (21) at (3.25, -0.5) {$h$};
		\node [style=none] (22) at (3.25, 0.5) {};
		\node [style=none] (23) at (3.25, 1.75) {};
		\node [style=none] (24) at (4.25, 1.75) {};
		\node [style=none] (25) at (3.75, -3.5) {};
	\end{pgfonlayer}
	\begin{pgfonlayer}{edgelayer}
		\draw [bend left] (17) to (20);
		\draw [in=-90, out=30, looseness=0.75] (17) to (18.center);
		\draw (17) to (19.center);
		\draw (20) to (22.center);
		\draw (23.center) to (22.center);
		\draw (24.center) to (18.center);
		\draw (19.center) to (25.center);
	\end{pgfonlayer}
\end{tikzpicture}$
\end{enumerate}
\end{proof}

\begin{corollary}
$(((A, \mulmap{1.2}{white}, \unitmap{1.2}{white}, \comulmap{1.2}{white}, \counitmap{1.2}{white}), \leftaction{0.4}{white}), 1)$ is a preunitary object.
\end{corollary}

Thus, we have a source of non-trivial preunitary objects so that we can form a non-trivial unitary category.


\section{Mixed unitary categories}
\label{Sec: MUCs}

A mixed unitary category has a unitary core which is a model of classical categorical quantum mechanics 
extended by a larger setting in which possibly infinite dimensional objects can 
be modelled. We are now ready for the definition of mixed unitary categories, which is the key structure developed in 
the first part of this thesis.

\subsection{Mixed unitary category}

\begin{definition}
A {\bf mixed unitary category} (MUC) is a $\dagger$-isomix category, $\C$, equipped with a strong $\dagger$-isomix functor 
$M: \U \to \C$ from a unitary category $\U$ to $\C$ such that there exists the following natural transformations:
\[ \mx': M(U) \oa X \to M(U) \ox X  \text{ with } \mx  ~\mx' = 1 \text{ and }\mx' ~ \mx = 1 \]
\[ \mx'': X \oa M(U) \to X \ox M(U) \text{ with } \mx  ~\mx'' = 1 \text{ and }\mx'' ~ \mx = 1 \]
A mixed unitary category, $M: \U \to \C$ is {\bf symmetric} if the functor $M$, the 
unitary category $\U$, and the $\dagger$-isomix category $\C$ are symmetric.
\end{definition}

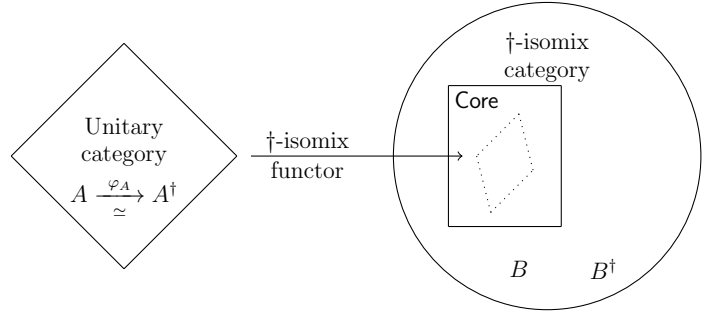
\begin{figure}[h]
	\centering
  \begin{tikzpicture} [scale=1.5]
	\begin{pgfonlayer}{nodelayer}
		\node [style=circle, scale=14] (0) at (4.5, 0) {};
		\node [style=none] (1) at (-3, 2) {};
		\node [style=none] (2) at (-5, 0) {};
		\node [style=none] (3) at (-3, -2) {};
		\node [style=none] (4) at (-1, 0) {};
		\node [style=none] (5) at (-3, -0.75) {$A \to^{\varphi_A}_{\simeq} A^\dagger$};
		\node [style=none] (6) at (-3, 0.5) {Unitary};
		\node [style=none] (7) at (-3, 0) {category};
		\node [style=none] (8) at (-0.75, 0) {};
		\node [style=none] (9) at (3, 0) {};
		\node [style=none] (10) at (0.25, 0.25) {$\dagger$-isomix};
		\node [style=none] (11) at (0.25, -0.25) {functor};
		\node [style=none] (12) at (4.5, 2) {$\dagger$-isomix};
		\node [style=none] (13) at (4.5, 1.5) {category};
		\node [style=none] (14) at (4, -2) {$B$};
		\node [style=none] (15) at (5.5, -2) {$B^\dagger$};
		\node [style=none] (16) at (4, 0.75) {};
		\node [style=none] (17) at (3.25, 0) {};
		\node [style=none] (18) at (3.5, -1) {};
		\node [style=none] (19) at (4.25, -0.25) {};
		\node [style=none] (20) at (2.75, 1.25) {};
		\node [style=none] (21) at (4.75, 1.25) {};
		\node [style=none] (22) at (4.75, -1.25) {};
		\node [style=none] (23) at (2.75, -1.25) {};
		\node [style=none] (24) at (3.25, 1) {$\Core$};
	\end{pgfonlayer}
	\begin{pgfonlayer}{edgelayer}
		\draw (1.center) to (2.center);
		\draw (2.center) to (3.center);
		\draw (3.center) to (4.center);
		\draw (4.center) to (1.center);
		\draw [->] (8.center) to (9.center);
		\draw [dotted] (16.center) to (17.center);
		\draw [dotted] (17.center) to (18.center);
		\draw [dotted] (18.center) to (19.center);
		\draw [dotted] (19.center) to (16.center);
		\draw (20.center) to (21.center);
		\draw (21.center) to (22.center);
		\draw (22.center) to (23.center);
		\draw (23.center) to (20.center);
	\end{pgfonlayer}
\end{tikzpicture}
\caption{Schematic diagram for MUC}
\label{Fig: MUC}
\end{figure}

In the definition of a MUC, the requirement of a transformation $\mx'$ which is inverse to $\mx$ ensures that the functor 
$M: \U \to \C$ factors through $\Core(\C)$. Figure \ref{Fig: MUC} is a schematic diagram of a MUC.

The $\dagger$-isomix category of a MUC is to be thought of as a larger space inside which a (small) unitary category embeds. 
Within the unitary category, $A \simeq A^\dagger$ by the means of the unitary structure map, 
however, outside the unitary core,  an object is not in general isomorphic to its dagger. 
	
In the rest of the thesis, MUCs are exclusively assumed to be symmetric unless stated otherwise.

Mix unitary categories organize themselves into a 2-category ${\sf MUC}$ (although we shall not discuss the 2-cell structure):
\begin{description}
\item[0-cells:]  Are mix unitary categories $M: \U \to \X$;
\item[1-cells:]  Are MUC morphisms: these are squares of $\dagger$-isomix functors $(F',F,\gamma): M \to N$ commuting up to a
 $\dagger$-linear natural isomorphism $\gamma: MF \Rightarrow F'N$:
 \[ \xymatrix{ \U \ar[d]_{F'} \ar@{}[drr]|{\Downarrow~\gamma} \ar[rr]^M & & \X \ar[d]^F \\ \V \ar[rr]_{N} & & \Y} \]
 The functor $F': \U \to \V$ is between unitary categories, and we demand of it that it preserves unitary structure in the 
 sense of Definition \ref{preserving-unitary-structure}, thus, whenever ${\varphi_A}$ is the unitary structure  then $F'(\varphi_A) 
 \rho^F$ is unitary structure.
 \item[2-cells:] These are ``pillows''  of natural transformations. $(\beta, \beta') : (F, F', \gamma_F) \Rightarrow (G, G', \gamma_G)$ 
is a 2-cell if and only if it satisfies the following equality:
\[ \xymatrix{ \U \ar@/_1pc/[dd]_{G'} \ar@/^1pc/[dd]^{F'}   \ar@{}[ddrr]|
{\Downarrow~\gamma_F} \ar[rr]^M & & \X \ar@/^1pc/[dd]^F \\ 
{\xLeftarrow{\beta'}} & & \\ 
\V \ar[rr]_{N} & & \Y} ~ \xymatrix{ \\ = \\ } ~ \xymatrix{ \U \ar@/_1pc/[dd]_{G'} \ar@{}[ddrr]|
{\Downarrow~\gamma_G} \ar[rr]^M & & \X \ar@/_1pc/[dd]_G \ar@/^1pc/[dd]^F \\ 
 & & {\xLeftarrow{\beta}} \\ 
\V \ar[rr]_{N} & & \Y} \]
 \end{description}
 
 We remark that we have observed that any MUC can be ``simplified'' to a dagger monoidal category with a strong $\dagger$-mix 
 Frobenius functor into a $\dagger$--isomix category: this is achieved by precomposing with ${\sf Mx}_\downarrow$.   This may 
 seem a worthwhile simplification, but it should be recognized that it simply transfers complexity from the unitary category itself 
 onto the preservator which must now ``create'' unitary structure:
 \[ \xymatrix{\U  \ar[d]_{{\sf Mx}^{*}_{\downarrow}} \ar[rr]^M & & \C \ar@{=}[d] \\
                     \U_\downarrow \ar[rr]_{{\sf Mx}_\downarrow;M} & & \C} \]
 Here $\U_\downarrow = (\U,\oa,\oa)$ is viewed as a dagger monoidal category and ${\sf Mx}_\downarrow^{*}$ is the inverse 
 of ${\sf Mx}_{\downarrow}$.  
 The point is that the preservator of the lower arrow ${\sf Mx}_\downarrow;M$ is non-trivial as it must encode the unitary 
 structure of $\U$.
 
\subsection{Canonical mixed unitary categories}

Our objective is now to show that the unitary construction of the previous section gives rise to a right adjoint to the underlying 
2-functor $U: {\sf MUC} \to {\sf MCC}$  where 
the 2-category ${\sf MCC}$ is defined as:

\begin{description}
\item[0-cells:] Its objects are  {\bf mixed $\dagger$-compact categories} (MCC), that is strong $\dagger$-Frobenius functors 
$V: \C \to \Y$ where $\C$ is a compact $\dagger$-LDC, $\Y$ is a $\dagger$-isomix category, 
and $V$ factors through the core of $\Y$ i.e, for all $\forall$ objects $C \in \C$, $Y \in \Y$,   $\exists$ $\mx': V(C) \oa Y \to V(C) \ox Y$ 
such that $\mx~\mx' = 1$ and $\mx' ~ \mx = 1$.
\item[1-cells:]   The 1-cells are squares of mix Frobenius functors which commute up to a linear natural isomorphism;
\item[2-cells:]    Are pillows of natural transformations (which we shall ignore).
\end{description}

An example of a mix $\dagger$-compact category is, of course, the inclusion of the core into a $\dagger$-isomix category 
$C:{\sf Core}(\X) \hookrightarrow \X$;

\begin{proposition}
$U: {\sf MUC} \to {\sf MCC}$ has a right adjoint ${\sf Unitary}: {\sf MCC} \to {\sf MUC}; (\C \to^V \X) \mapsto ({\sf Unitary}(\C) 
\to^{U;V} \X)$.
\end{proposition}

\begin{proof}

The couniversal diagram is as follows:
\[ \xymatrix{ {\U \to^M \X} \ar[rr]^{(F,G,\gamma)} \ar[d]_{(F^\flat,G,\gamma^\flat)} && {\C \to^V \Y} \\
                    {{\sf Unitary}(\C) \to_{U;V} \Y} \ar[urr]_{\epsilon} } \]
where $\epsilon$ is the square on the left and $(F^\flat,G,\gamma^\flat)$ is the square on the right:
 \[ \xymatrix{{\sf Unitary}(\C) \ar[d]_U \ar[r]^{~~~~U} & \C \ar[r]^V & \Y \ar@{=}[d] \\
                         \C \ar[rr]_V && \Y} 
     ~~~~~~~~
     \xymatrix{\U \ar[dr]_F \ar[d]^{F^\flat} \ar@{}[drr]|{~~~~~~~~~\uparrow~\gamma} \ar[rr]^M & & \X \ar[d]^{G} \\
                      {\sf Unitary}(\C) \ar[r]_{~~~~U} & \C \ar[r]_V & \Y} \]
                      
It follows from Proposition \ref{Prop: Couniversal} that the couniversal diagram commute.

\end{proof}

This proposition means that in building a non-trivial MUC from a mixed $\dagger$-compact categories it suffices to show that the 
compact $\dagger$-LDC contains non-trivial pre-unitary objects. 


\section{Examples: Mixed unitary categories}
\label{Sec: MUC examples}

We have already noted that dagger monoidal categories are automatically 
unitary categories in which the unitary structure is given by identity maps.  
The identity functors then give a rather trivial MUC. 
One can construct a MUC from any $\dagger$-isomix category using the unitary construction: 
for any $\dagger$-isomix category, $\X$, ${\sf Unitary}(\Core(\X)) \to^{U} \Core(\X) \hookrightarrow \X$ is a MUC. 
More excitingly one can take the bicompletion \cite{Joy95} of the $\dagger$-monoidal category: this is a non-trivial $\dagger$-isomix $*$-autonomous category 
extension of the original $\dagger$-monoidal category and provides, thus, an interesting example of how MUCs arise.  

Our purpose in this section is to exhibit some non-trivial manifestations of the various structural components of a MUC.  
To this end we discuss in some detail three basic examples.


\subsection{Finite dimensional framed vector spaces}
\label{Sec: FFVec}
In this section we show that the example ${\sf FFVec}_K$, the category of 
finite dimensional framed vector spaces defined in 
Section \ref{subsection:fdfv} is a unitary category (hence is immediately a mixed unitary category). The unitary structure map on each object $(V, {\cal V})$ is defined as follows:
\[ \varphi_{(V,{\cal  V})}: (V,{\cal  V}) \to (V,{\cal  V})^\dag; v_i \mapsto \widetilde{v_i} \]
and it remains to check the coherences {\bf [U.3]}--{\bf [U.6]}.  First note that {\bf [U.4]} holds immediately by the observation above that 
$\iota(v_i) = \widetilde{\widetilde{v_i}}$.  For {\bf [U.3]} we require that $\varphi_{A^\dag}(\widetilde{a_i}) = (\varphi_A^{-1})^\dag (\widetilde{a_i})$ 
the result is a higher-order term so, we may check that the evaluations are the same on basis elements:
\begin{eqnarray*}
	(\varphi_{A^\dag}(\widetilde{a_i}) ) (\widetilde{a_j}) & = & \widetilde{\widetilde{a_i}}(\widetilde{a_j}) = \partial_{i,j} \\
	((\varphi_{A}^{-1})^\dag(\widetilde{a_i}))(\widetilde{a_j}) & = & \widetilde{a_i} ((\varphi_{A}^{-1})^\dag(\widetilde{a_j})) = \widetilde{a_i}(a_j) =  \partial_{i,j}
\end{eqnarray*}
Note that {\bf [U.5]}(a) and {\bf [U.5]}(b), in this example, require $\lambda_\top = \varphi_\top$ which can easily be verified as each reduces to conjugation.
{\bf [U.6]}(a) and {\bf [U.6]}(b), in this example, are the same requirement which is verified by:
\[ \lambda_\ox(\varphi_A \ox \varphi_B(a_i \ox b_j) ) = \lambda_\ox (\widetilde{a_i} \ox \widetilde{b_j}) = \widetilde{a_i \ox b_j} = \varphi_{A \ox B} (a_i \ox b_j) \]

This gives:

\begin{proposition}
	${\sf FFVec}_K$ with the unitary structure above is a MUC.
\end{proposition}

This raises the question of what precisely the unitary maps of this example are.  To elucidate this we note that  a functor can easily be constructed
$U:{\sf FFVec}_K \to {\sf Mat}(K)$ where, for each object in ${\sf FFVec}_K$ we choose a total order on the elements of the basis and note that 
any map is then given by a matrix acting on the bases: thus a matrix in ${\sf Mat}(K)$ with the appropriate dimensions.  We now observe:

\begin{lemma} 
	An isomorphism $u: (A,{\cal A}) \to (B,{\cal B})$ in ${\sf FFVec}_K$ is unitary 
	if and only if $U(f)$ is unitary in ${\sf Mat}(K)$.
\end{lemma}

\proof While $U$ does not preserve $(\_)^\dag$ on the nose it does so up to the natural equivalence determined by $U(\varphi_A)$ which being a basis 
permutation is a unitary equivalence.   Thus, it is not hard to see that the following diagram commutes:
\[ \xymatrix{ U(B,{\cal B})  \ar[d]_{U(f)^\dag} \ar[rr]^{U(\varphi_B)} & &U((B,{\cal B})^\dag) \ar[d]^{U(f^\dag)} \\
	U(A,{\cal A}) \ar[rr]_{U(\varphi_A)} & & U((A,{\cal A})^\dag) } \]
Recall that in the category of matrices, the dagger is stationary on objects so $U(B,{\cal B}) = U(B,{\cal B})^\dag$.  

Now suppose $u$ is unitary in ${\sf FFVec}_K$  then $u^{-1} = \varphi_B u^\dagger \varphi_A^{-1}$ so that 
\[ U(u)^{-1} = U(u^{-1}) = U(\varphi_B u^\dagger \varphi_A^{-1}) = U(\varphi_B) U(u^\dagger) U(\varphi_A^{-1}) = U(u)^\dagger \]
so that its underlying map is unitary. Conversely, if $U(u)$ is unitary then 
\[ U(u^{-1}) = U(u)^{-1} = U(u)^\dag = U(\varphi_B u^\dagger \varphi_A^{-1}) \]
which immediately implies, as $U$ is faithful, that $u$ is unitary in ${\sf FFVec}_K$.
\endproof

One might reasonably regard this as a rather roundabout way to describe the standard notion of a unitary map.  However, two things of importance have been 
achieved.  First an example of a unitary category with a non-stationary dagger and, thus, a non-identity unitary structure, has been exhibited.  Second we have 
shown how the standard unitary structure may be re-expressed in this formalism using non-stationary constructs. 


\subsection{Finiteness matrices}
In Section \ref{Sec: Finiteness matrices}, we described the category of finiteness matrices, ${\sf FMat}(\C)$. 
The core of ${\sf FMat}(\C)$ is the subcategory determined by objects whose webs are finite sets, that is the objects are $(X, P(X))$ 
where $X$ is a finite set. 
Clearly, $\Core({\sf FMat}(\C))$ is then equivalent to the category of finite dimensional matrices, ${\sf Mat}(\C)$.  
This is a well-known $\dagger$-compact closed category, which is a unitary category with unitary structure given by 
identity maps (as $(\_)^\dagger$ is stationary on objects).  The inclusion ${\cal I}: {\sf Mat}(\C) \to {\sf FMat}(\C)$ provides an important example of a MUC.   


\subsection{The embedding of finite dimensional Hilbert Spaces into Chu Spaces}
\label{Sec: CHU MUC}
In Section \ref{Section: Chu}, we showed that the Chu construction applied to a symmetric conjugative closed monoidal 
category, $\X$, with pullbacks gives a $\dagger$-isomix category. Recall that the dagger in the resulting category of Chu 
spaces is given by composing the conjugation with the dualizing functor.  In this section, we start by discussing, in  general, 
the construction of a mixed unitary category from a Chu category ${\sf Chus}_\X(I)$.  A crucial step in this is to identify objects 
which are in the core of this category.

Recall that a  symmetric monoidal closed category, $\X$, is (degenerately) a compact linearly distributive category and, thus, 
there may be objects which have linear adjoints: these are called {\bf nuclear} objects \cite{HiR89}.  Explicitly a nuclear 
object $A$ in a symmetric monoidal closed  category is an object with $A \multimap B \cong A^{*} \ox B$, where 
$A^* := A \multimap I$.   The nuclear objects form a compact closed subcategory of $\X$ which is conjugative when $\X$ is 
conjugative.  In ${\sf Vec}_\mathbb{C}$ the nucleus consist precisely of the finite dimensional vector spaces.  
If $(\eta, \epsilon): A \dashv\!\!\!\dashv B$ is witness that $A$ (and $B$) are nuclear in $\X$ then the object 
$(A,B,\epsilon,c_\otimes\epsilon)$ is in the core of ${\sf Chus}_\X(I)$ because in the second component of the tensor product 
with any other object $(X,Y,\nu,c_\ox \nu)$ one has the degenerate pullback:
\[ \xymatrix{
& Y \ox B \ar[rd]^{\simeq} \ar[ld] & \\ 
X \multimap B \ar[rd]^{\simeq}  & & Y \ox A^{*} \ar[ld]  \\
& X \multimap A^* \to^{\simeq} (X \multimap I) \ox A &
} \]
where we use the isomorphism $B \to^{\simeq} A^{*}$.

In this manner the nuclear objects of ${\sf Nuclear}(\X)$, which form a compact closed category with a dagger, 
may be embedded into the core of ${\sf Chus}_\X(I)$.  To obtain a unitary category it suffices then to use the unitary 
construction for which, to obtain a non-trivial result, we need to show that there are non-trivial examples of pre-unitary objects. 
To achieve this we consider an object $H$ for which $(e, n): H \dashv\!\!\!\dashv \overline{H}$ and such that $e$ satisfies:
\[ \xymatrix{\overline{\overline{H}} \otimes \overline{H} \ar[d]_{\varepsilon \otimes 1} \ar[r]^{\chi} & \overline{\overline{H} \otimes 
\overline{H}} \ar[dd]^{\overline{e}}\\
                   H \otimes \overline{H} \ar[d]_{e} \\
                   I \ar[r]_{\chi^{\!\!\!\circ}} & \overline{I} } \]
               
 For such an object we note:
 \begin{align*}
{ \overline{(H,\overline{H},e,c_\otimes e)}}^* & =  (\overline{H},\overline{\overline{H}},\chi\overline{c_\otimes e} 
(\chi^{\!\!\!\circ})^{-1},\chi\overline{e} (\chi^{\!\!\!\circ})^{-1})^* \\
 & =  (\overline{\overline{H}},\overline{H},\chi\overline{e} (\chi^{\!\!\!\circ})^{-1},\chi\overline{c_\otimes e} (\chi^{\!\!\!\circ})^{-1})
 \end{align*}
 
This makes
 $$(\varepsilon^{-1},1) : (H,\overline{H},e,c_\otimes e) \to (\overline{\overline{H}},\overline{H},\chi\overline{e} 
 (\chi^{\!\!\!\circ})^{-1},\chi\overline{c_\otimes e} (\chi^{\!\!\!\circ})^{-1})$$
a preunitay map.  Note that it is a Chu map by the commuting diagram above and as $\overline{\varepsilon} =\varepsilon$ we have 
$$(\varepsilon^{-1},1) (1,\overline{\varepsilon}) = (\varepsilon^{-1},1) (1,\varepsilon) = (\varepsilon^{-1},\varepsilon)$$
where $(\varepsilon^{-1},\varepsilon)$ is the involutor.

In ${\sf Vec}_{\mathbb{C}}$ a map $e: H \otimes \overline{H} \to \mathbb{C}$ is a ``sesquilinear form'' and the diagram above 
asserts that it is in addition a symmetric form.  Any  Hilbert space with its inner product, thus, satisfies the above conditions.   
Thus,  it is clear that the embedding of the category of finite dimensional Hilbert Spaces into Chu spaces, ${\sf FHilb} 
\hookrightarrow {\sf Chus}_{{\sf Vec}_\C} (\C)$ is a mixed unitary category.  The embedding is in fact a full and faithful 
embedding which extends to {\em all\/} Hilbert spaces (although only the finite dimensional ones land in the core).  

Explicitly the embedding is defined as follows: suppose $H$ is a (finite dimensional) Hilbert Space, then the corresponding 
Chu Space is given by $(H, \overline{H}, \langle - | - \rangle_H)$, where $\langle - | - \rangle_H: H \ox \overline{H} \to \C$ is 
the inner product. For any linear map $H \to^{f} K$ between Hilbert Spaces, the corresponding Chu map is given by $(f, f^\dagger): 
(H, \overline{H}, \langle - | - \rangle_H) \to (K, \overline{K}, \langle - | - \rangle_K)$, where $f^\dagger$ is the Hermitian adjoint of 
$f$ so, $\langle f(a) | b \rangle = \langle a | f^\dagger(b) \rangle$.

Furthermore, observe that $(H, \overline{H}, \langle - | - \rangle_H)^\dagger :=  \overline{(H, \overline{H}, \langle - | - \rangle_H)^*} = 
(H, \overline{H}, \langle - | - \rangle_H)$. Hence, this embedding preserves the (stationary) dagger for all Hilbert spaces.  However, 
the par of two infinite dimensional Hillbert spaces in this Chu category is not a Hilbert space so that the duality cannot be seen within 
the category of Hilbert spaces.






\chapter{Summary}
\label{chap: part 1 summary}

Chapters \ref{Chap: dagger-LDC}, and \ref{Chap: MUCs} extend the theory of $\dagger$-monoidal categories and $\dagger$-compact closed categories to linearly distributive 
and $*$-autonomous settings to obtain the categorical semantics of (multiplicative) $\dagger$-linear logic.  In these linear settings, the 
two different tensor products (tensor and par) must be flipped by the dagger.  Thus, one cannot have a stationary (identity on objects) 
dagger, and hence one is forced to replace the conventional dagger by a contravariant structure-preserving involution.  This has coherence  
consequences: section \ref{Section: dagger LDC} is dedicated to understanding the details of these coherences.

If multiplicative  $\dagger$-linear logic is to provide a semantics for a generalized categorical quantum mechanics (CQM), then notions 
such as isometry and unitary isomorphism, which are central to CQM, should have an expression in this logic.  In section \ref{Sec: unitary} 
we showed that with additional ``unitary structure'' one can recapture classical CQM as a ``unitary core''  of multiplicative  $\dagger$-linear logic.
Furthermore, we showed how, from any $\dagger$-isomix category, it is always possible to extract  
a ``unitary core'' which is, up to equivalence, a $\dagger$-monoidal category (i.e a classical semantic setting for CQM).   

This led to the notion of a mixed unitary category (MUC) given by a $\dagger$-isomix category with a chosen unitary 
core as our proposal for an extension of CQM.   A MUC can be viewed as an extension much as a $K$-algebra 
extends a field $K$ and permits the expression of properties which are difficult to express within $K$ itself.  
In the extended setting of a MUC -- finiteness matrices with its core for example -- provides an  
extension of the classical CQM setting in which infinite dimensional types, such as those given by the exponential modalities,  
 are present.   Furthermore, in the extended setting one can bend, and yank wires without the category being compact.

This concludes the first part of this thesis. The second part discusses the application of MUCs to CQM. 

%% file: chapter2.tex

\chapter{Categorical Quantum Mechanics}
\label{Chap: CQM}

The program of Categorical Quantum Mechanics (CQM) was started by Coecke and Abramsky \cite{AC04} 
in 2004 with the aim of developing a high-level formal language for quantum mechanics while moving 
away from the standard formalism based on Hilbert spaces. CQM derived ideas from logic and 
computer science, and used the diagrammatic language of compact closed categories for an intuitive 
but mathematically rigorous presentation of the fundamental axioms of quantum theory. 

The purpose this chapter is to review the fundamentals of Categorical Quantum Mechanics 
and to describe its essential features used to study quantum mechanics. 

\section{Dagger monoidal categories}
\label{Sec: dag mon cat}

Categorical quantum mechanics (CQM) models physical systems as objects within a monoidal category and 
the physical processes as the maps in the category. The identity maps equate to the do-nothing processes.
 Sequential composition of processes is given by the composition of maps while parallel composition 
 is modelled using the tensor product.  Categorifying the notion of inner product  in traditional quantum 
 theory produces  a $\dagger$-functor for monoidal categories giving rise to the theory of 
 $\dagger$-monoidal categories. 

 \subsection{Graphical calculus for monoidal categories}
 Monoidal categories come equipped with a graphical calculus \cite{Sel10},
 which allows for diagrammatic representation and manipulation of systems and processes. The 
 availability of this graphical calculus is perhaps the most attractive aspect of using monoidal 
 categories to study quantum mechanics.

 In the graphical calculus of monoidal categories, the objects are represented by wires, and arrows by circles
\footnote{In the CQM community, maps are often drawn as boxes. A circle and a box are topologically same.}. 
 An identity  arrow is given by a wire without by circle.
 \[ \begin{tikzpicture}
	\begin{pgfonlayer}{nodelayer}
		\node [style=none] (0) at (-2, 3) {};
		\node [style=none] (1) at (-2, -0.25) {};
		\node [style=none] (2) at (-2.5, 1.5) {$A$};
		\node [style=none, scale=2] (3) at (-2, -1) {An object $A$};
	\end{pgfonlayer}
	\begin{pgfonlayer}{edgelayer}
		\draw (0.center) to (1.center);
	\end{pgfonlayer}
\end{tikzpicture} ~~~~~~~~
\begin{tikzpicture}
	\begin{pgfonlayer}{nodelayer}
		\node [style=none] (0) at (-2, 3) {};
		\node [style=none] (1) at (-2, -0.25) {};
		\node [style=none] (2) at (-2.5, 2.5) {$A$};
		\node [style=none, scale=2] (3) at (-2, -1) {$f: A \to B$};
		\node [style=none] (4) at (-2.5, 0.5) {$B$};
		\node [style=circle, scale=1.5] (5) at (-2, 1.5) {};
		\node [style=none] (6) at (-2, 1.5) {$f$};
	\end{pgfonlayer}
	\begin{pgfonlayer}{edgelayer}
		\draw (0.center) to (5);
		\draw (1.center) to (5);
	\end{pgfonlayer}
\end{tikzpicture} ~~~~~~~~
\begin{tikzpicture}
	\begin{pgfonlayer}{nodelayer}
		\node [style=none] (0) at (-2, 3) {};
		\node [style=none] (1) at (-2, -0.25) {};
		\node [style=none] (2) at (-2.25, 1.5) {$A$};
		\node [style=none, scale=2] (3) at (-2, -1) {Identity arrow $1_A$};
	\end{pgfonlayer}
	\begin{pgfonlayer}{edgelayer}
		\draw (1.center) to (0.center);
	\end{pgfonlayer}
\end{tikzpicture}\]

Note that the diagram for the object $A$ is same as the identity arrow of $A$ which is to be interpreted as 
follows: a physical system which does nothing is same as the system itself. 

Composition of maps is given by connecting the wires sequentially. 
Tensor product of object is given by parallel juxtaposition of wires. The unit object is given by an 
empty circle representing the empty system. Recall that in LDCs, when the tensor and the par units are different, we used 
labelled wires to represent units. 
  \[ \begin{tikzpicture}
	 \begin{pgfonlayer}{nodelayer}
		 \node [style=none] (0) at (-2, 3) {};
		 \node [style=none] (1) at (-2, -0.25) {};
		 \node [style=none] (2) at (-2.5, 2.75) {$A$};
		 \node [style=none, scale=2] (3) at (-2, -1) {Composing maps};
		 \node [style=none] (4) at (-2.5, 1.5) {$B$};
		 \node [style=circle, scale=1.5] (5) at (-2, 2) {};
		 \node [style=none] (6) at (-2, 2) {$f$};
		 \node [style=circle, scale=1.5] (7) at (-2, 0.75) {};
		 \node [style=none] (8) at (-2, 0.75) {$g$};
		 \node [style=none] (9) at (-2.5, -0) {$C$};
	 \end{pgfonlayer}
	 \begin{pgfonlayer}{edgelayer}
		 \draw (0.center) to (5);
		 \draw (1.center) to (7);
		 \draw (7) to (5);
	 \end{pgfonlayer}
 \end{tikzpicture} ~~~~~~~~
  \begin{tikzpicture}
	 \begin{pgfonlayer}{nodelayer}
		 \node [style=none] (0) at (-2.25, 3) {};
		 \node [style=none] (1) at (-2.25, -0.25) {};
		 \node [style=none] (2) at (-2.5, 2.5) {$A$};
		 \node [style=none, scale=2] (3) at (-1.75, -1) {$f : A \ox C \to B \ox D$};
		 \node [style=none] (4) at (-2.5, 0.25) {$B$};
		 \node [style=none] (5) at (-1, 2.5) {$C$};
		 \node [style=none] (6) at (-1.25, 3) {};
		 \node [style=none] (7) at (-1, 0.25) {$D$};
		 \node [style=none] (8) at (-1.25, -0.25) {};
		 \node [style=circle, scale=1.5] (9) at (-1.75, 1.5) {};
		 \node [style=none] (10) at (-1.75, 1.5) {$f$};
	 \end{pgfonlayer}
	 \begin{pgfonlayer}{edgelayer}
		 \draw [in=-90, out=120, looseness=1.00] (9) to (0.center);
		 \draw [in=60, out=-90, looseness=1.00] (6.center) to (9);
		 \draw [in=-120, out=90, looseness=1.00] (1.center) to (9);
		 \draw [in=-60, out=90, looseness=1.00] (8.center) to (9);
	 \end{pgfonlayer}
 \end{tikzpicture}  ~~~~~~~~
 \begin{tikzpicture}
	 \begin{pgfonlayer}{nodelayer}
		 \node [style=none] (0) at (-2, 3) {};
		 \node [style=none] (1) at (-2, -0.25) {};
		 \node [style=none] (2) at (-2.5, 2.5) {$A$};
		 \node [style=none, scale=2] (3) at (-1.75, -1) {$f \ox g : A \ox C \to B \ox D$};
		 \node [style=circle, scale=1.5] (4) at (-2, 1.25) {};
		 \node [style=none] (5) at (-2.5, 0.25) {$B$};
		 \node [style=none] (6) at (-2, 1.25) {$f$};
		 \node [style=none] (7) at (-0.8, 2.5) {$C$};
		 \node [style=circle, scale=1.5] (8) at (-1.25, 1.25) {};
		 \node [style=none] (9) at (-1.25, 3) {};
		 \node [style=none] (10) at (-0.8, 0.25) {$D$};
		 \node [style=none] (11) at (-1.25, 1.25) {$g$};
		 \node [style=none] (12) at (-1.25, -0.25) {};
	 \end{pgfonlayer}
	 \begin{pgfonlayer}{edgelayer}
		 \draw (1.center) to (4);
		 \draw (4) to (0.center);
		 \draw (12.center) to (8);
		 \draw (8) to (9.center);
	 \end{pgfonlayer}
 \end{tikzpicture} ~~~~~~~~ I := 
 \begin{tikzpicture}
	 \begin{pgfonlayer}{nodelayer}
		 \node [style=none] (0) at (0, 2) {};
	 \end{pgfonlayer}
 \end{tikzpicture} \]

 In CQM, maps from the tensor unit to any other object are referred to as {\bf states}, maps 
 into the tensor unit from any other object are referred to as {\bf effects}, and maps that 
 start and end in the tensor unit are referred to as {\bf scalars}. States and effects are represented 
 using triangles: 
 \[ \begin{tikzpicture}
	\begin{pgfonlayer}{nodelayer}
		\node [style=none] (0) at (0, 2.75) {$\phi$};
		\node [style=none] (1) at (0, 3.25) {};
		\node [style=none] (2) at (-0.5, 2.5) {};
		\node [style=none] (3) at (0.5, 2.5) {};
		\node [style=none] (4) at (0, 2.5) {};
		\node [style=none] (5) at (0, 1) {};
		\node [style=none, scale=2] (7) at (0, 0.5) {State};
	\end{pgfonlayer}
	\begin{pgfonlayer}{edgelayer}
		\draw (1.center) to (2.center);
		\draw (2.center) to (3.center);
		\draw (3.center) to (1.center);
		\draw (4.center) to (5.center);
	\end{pgfonlayer}
\end{tikzpicture} ~~~~~~~~~~~~ \begin{tikzpicture}
	\begin{pgfonlayer}{nodelayer}
		\node [style=none] (0) at (0, 1.5) {$\phi$};
		\node [style=none] (1) at (0, 1) {};
		\node [style=none] (2) at (-0.5, 1.75) {};
		\node [style=none] (3) at (0.5, 1.75) {};
		\node [style=none] (4) at (0, 1.75) {};
		\node [style=none] (5) at (0, 3.25) {};
		\node [style=none, scale=2] (7) at (0, 0.5) {Effect};
	\end{pgfonlayer}
	\begin{pgfonlayer}{edgelayer}
		\draw (1.center) to (2.center);
		\draw (2.center) to (3.center);
		\draw (3.center) to (1.center);
		\draw (4.center) to (5.center);
	\end{pgfonlayer}
\end{tikzpicture} \]

 The assosciativity, left unitor, right unitor and the symmetry natural isomorphisms are given
  as follows:
 \[ \begin{tikzpicture}
	\begin{pgfonlayer}{nodelayer}
		\node [style=none] (0) at (-2.5, 3) {};
		\node [style=none] (1) at (-2.5, -0.25) {};
		\node [style=none] (2) at (-2.5, 3.25) {$A$};
		\node [style=none, scale=2] (3) at (-1.75, -1) {$a_\ox: (A \ox B) \ox C \to A \ox (B \ox C)$};
		\node [style=none] (4) at (-1.75, 3.25) {$B$};
		\node [style=none] (5) at (-1, 3) {};
		\node [style=none] (6) at (-1, -0.25) {};
		\node [style=none] (7) at (-1.75, 3) {};
		\node [style=none] (8) at (-1.75, -0.25) {};
		\node [style=none] (9) at (-1, 3.25) {$C$};
	\end{pgfonlayer}
	\begin{pgfonlayer}{edgelayer}
		\draw (0.center) to (1.center);
		\draw (7.center) to (8.center);
		\draw (5.center) to (6.center);
	\end{pgfonlayer}
\end{tikzpicture}
~~~~~~~~
\begin{tikzpicture}
	\begin{pgfonlayer}{nodelayer}
		\node [style=none] (0) at (-2.25, 3) {};
		\node [style=none] (1) at (-2.25, -0.25) {};
		\node [style=none] (2) at (-2.25, 3.25) {$A$};
		\node [style=none, scale=2] (3) at (-2.25, -1) {$u_\ox^r: A \ox I \to A$};
	\end{pgfonlayer}
	\begin{pgfonlayer}{edgelayer}
		\draw (0.center) to (1.center);
	\end{pgfonlayer}
\end{tikzpicture} ~~~~~~~~~~~~
\begin{tikzpicture}
	\begin{pgfonlayer}{nodelayer}
		\node [style=none] (0) at (-2.25, 3) {};
		\node [style=none] (1) at (-2.25, -0.25) {};
		\node [style=none] (2) at (-2.25, 3.25) {$A$};
		\node [style=none, scale=2] (3) at (-2.25, -1) {$u_\ox^l: I \ox A \to  A$};
	\end{pgfonlayer}
	\begin{pgfonlayer}{edgelayer}
		\draw (0.center) to (1.center);
	\end{pgfonlayer}
\end{tikzpicture} 
 \]
 
Note that the associativity isomorphisms, the left and the right unitors are identity maps. In fact, this 
is the graphical caluclus for a strict monoidal category. Since every monoidal category is equivalent to a strict 
monoidal category \cite{Mac13} allows one to use this calculus on any monoidal category. 

The symmetry map is represented using crossed wires as shown in the left. The inverse law for a 
symmetric monoidal category is given by the diagrammatic equation in the right: 
\[  \begin{tikzpicture}
	\begin{pgfonlayer}{nodelayer}
		\node [style=none] (0) at (-1.5, 3.25) {};
		\node [style=none] (1) at (0, -0.5) {};
		\node [style=none] (2) at (-1.75, 2.5) {$A$};
		\node [style=none, scale=2] (3) at (-0.75, -1.25) {$c_\ox: A \ox B \to B \ox A$};
		\node [style=none] (4) at (0, 3.25) {};
		\node [style=none] (5) at (-1.5, -0.5) {};
		\node [style=none] (6) at (0.25, 2.5) {$B$};
		\node [style=none] (7) at (-1.75, 0) {$B$};
		\node [style=none] (8) at (0.25, 0) {$A$};
	\end{pgfonlayer}
	\begin{pgfonlayer}{edgelayer}
		\draw [in=90, out=-90, looseness=1.25] (0.center) to (1.center);
		\draw [in=-90, out=90, looseness=1.25] (5.center) to (4.center);
	\end{pgfonlayer}
\end{tikzpicture}
~~~~~~~~
\begin{tikzpicture}
	\begin{pgfonlayer}{nodelayer}
		\node [style=none] (0) at (-1.75, 3) {};
		\node [style=none] (1) at (-0.75, 3) {};
		\node [style=none] (2) at (-1.75, -0.5) {};
		\node [style=none] (3) at (-0.75, -0.5) {};
		\node [style=none] (4) at (0.75, 3) {};
		\node [style=none] (5) at (1.75, 3) {};
		\node [style=none] (6) at (1.75, -0.5) {};
		\node [style=none] (7) at (0.75, -0.5) {};
		\node [style=none] (8) at (0, 1.25) {$=$};
		\node [style=none, scale=2] (9) at (0, -1.25) {$(c_\ox)_{B,A} = (c_\ox)_{A,B}^{-1}$};
		\node [style=none] (10) at (-1.75, 1.25) {};
		\node [style=none] (11) at (-0.75, 1.25) {};
	\end{pgfonlayer}
	\begin{pgfonlayer}{edgelayer}
		\draw [in=90, out=-90, looseness=1.25] (1.center) to (10.center);
		\draw [in=90, out=-90, looseness=1.25] (0.center) to (11.center);
		\draw [in=90, out=-90, looseness=1.25] (11.center) to (2.center);
		\draw [in=90, out=-90, looseness=1.25] (10.center) to (3.center);
		\draw (4.center) to (7.center);
		\draw (6.center) to (5.center);
	\end{pgfonlayer}
\end{tikzpicture} \]

An equation holds in a SMC if and only if it holds in the graphical calculus. 
Proving equations using graphical calculus is simpler than proving equations by 
reasoning with mathematical symbols because human brain excels at processing 
visual information. 

Two diagrams are the same in the graphical calculus if one can be transformed to another 
up to planar isotopy or by using an axiom. For details on graphical calculus for monoidal 
categories see \cite{Sel10}. 
 
\subsection{Dagger monoidal categories}

Complex Hilbert spaces are used as the de facto framework for describing quantum 
processes. The inner product structure on these spaces allows the notion of 
adjoint which is in turn needed to define quantum observables: 
self-adjoint operators on the space.  
In quantum computing, every quantum logic gate is represented by a unitary (in other words a self-adjoint) matrix.
CQM abstracts the notion of inner products as dagger \cite{Sel07, CoK17} functors for $\dagger$-monoidal categories. 

\begin{definition}
A {\bf dagger} category is a category, $\X$, equipped with a involutive 
$(f^{\dag \dag} = f)$ contravariant functor 
$\dagger: \X^\op \to \X$ which is the identity on objects $(A = A^\dag)$. 
\end{definition}

Note that, a given category can have more than one $\dagger$-functor, for example for 
the category of complex matrices (see Section \ref{Sec: Mat(R)}) conjugate transpose and transpose 
operations on matrices gives two different dagger functors. Hence, $\dagger$ is a 
structure rather than a property for a category. Because $\dagger$ is a contravariant functor, $(fg)^\dag = g^\dag f^\dag$.

Given any map $f: A \to B$, the map $f^\dagger: B \to A$ is referred to as its {\bf adjoint}. 
A {\bf unitary} isomorphism in a $\dagger$-category is an isomorphism $f$ such that $f^\dag = f^{-1}$. 
A map $f: A \to B$ is an isometry if $f f^\dag = 1_A$. 

\begin{definition}
A {\bf $\dagger$-symmetric monoidal} category is a symmetric monoidal category 
which is also a $\dagger$-category such that the
$\dagger$ behaves coherently with the monoidal structure:
\begin{enumerate}
	\item for all maps $f$, $g$, $(f \ox g)^\dag = f^\dag \ox g^\dag$
	\item the assosciator, the unitors and the symmetry map are unitary natural isomorphisms
\end{enumerate}
\end{definition}

Next, we discuss a key correspondence often used in 
 quantum information theory called the operator-state duality, also referred to as the Choi-Jamiolkowski isomorphism: 
every linear map between (finite dimensional) Hilbert spaces $H$ and $K$ corresponds precisely to a state 
in the tensor product space $H \ox K$. This correspondence is abstracted as follows using the compact structure \cite{Kel97}
for monoidal categories. 

\begin{definition}
	\label{defn: right dual}
In a monoidal category, an object $B$ is {\bf right dual} to an object $A$  if there exists maps:
\[ \eta: I \to A \ox B ~~~~~~~~~~	\epsilon: B \ox A \to I \]
such that:
\[ (1 \ox \eta) (\epsilon \ox 1) = 1_B ~~~~~~~~~~~~
(\eta \ox 1)(1 \ox \epsilon) = 1_A \]
\end{definition}

The maps $\eta$ and $\epsilon$ are represented in graphical calculus by a cap and a cup respectively:
\[ \eta: I \to A \ox B = \begin{tikzpicture}
		\begin{pgfonlayer}{nodelayer}
			\node [style=none] (0) at (0, 2) {};
			\node [style=none] (1) at (0, 1) {};
			\node [style=none] (2) at (1, 1) {};
			\node [style=none] (3) at (1, 2) {};
		\end{pgfonlayer}
		\begin{pgfonlayer}{edgelayer}
			\draw (1.center) to (0.center);
			\draw [bend left=90, looseness=2.25] (0.center) to (3.center);
			\draw (3.center) to (2.center);
		\end{pgfonlayer}
	\end{tikzpicture} ~~~~~~~~~~~ 
	\epsilon: B \ox A \to I = \begin{tikzpicture}
		\begin{pgfonlayer}{nodelayer}
			\node [style=none] (0) at (0, 1.75) {};
			\node [style=none] (1) at (0, 2.75) {};
			\node [style=none] (2) at (1, 2.75) {};
			\node [style=none] (3) at (1, 1.75) {};
		\end{pgfonlayer}
		\begin{pgfonlayer}{edgelayer}
			\draw (1.center) to (0.center);
			\draw [bend right=90, looseness=2.25] (0.center) to (3.center);
			\draw (3.center) to (2.center);
		\end{pgfonlayer}
	\end{tikzpicture}	
\]
The equations satisfied by a right dual 
are also referred to as the snake equations owing to their shape 
in the graphical calculus:
	\[ 	\begin{tikzpicture}
		\begin{pgfonlayer}{nodelayer}
			\node [style=none] (0) at (0, 1.75) {};
			\node [style=none] (1) at (0, 3.5) {};
			\node [style=none] (2) at (1, 2.75) {};
			\node [style=none] (3) at (1, 1.75) {};
			\node [style=none] (4) at (2, 2.75) {};
			\node [style=none] (5) at (2, 0.75) {};
			\node [style=none] (6) at (-0.25, 3.25) {$B$};
			\node [style=none] (7) at (0.65, 2.25) {$A$};
			\node [style=none] (8) at (2.25, 1) {$B$};
		\end{pgfonlayer}
		\begin{pgfonlayer}{edgelayer}
			\draw (1.center) to (0.center);
			\draw [bend right=90, looseness=2.25] (0.center) to (3.center);
			\draw (3.center) to (2.center);
			\draw (5.center) to (4.center);
			\draw [bend right=90, looseness=1.75] (4.center) to (2.center);
		\end{pgfonlayer}
	\end{tikzpicture} = \begin{tikzpicture}
		\begin{pgfonlayer}{nodelayer}
			\node [style=none] (0) at (3.25, 2.25) {$B$};
			\node [style=none] (1) at (3, 3.5) {};
			\node [style=none] (2) at (3, 0.75) {};
		\end{pgfonlayer}
		\begin{pgfonlayer}{edgelayer}
			\draw (2.center) to (1.center);
		\end{pgfonlayer}
	\end{tikzpicture}
	~~~~~~~~~~~~
	\begin{tikzpicture}
		\begin{pgfonlayer}{nodelayer}
			\node [style=none] (0) at (2, 1.75) {};
			\node [style=none] (1) at (2, 3.5) {};
			\node [style=none] (2) at (1, 2.75) {};
			\node [style=none] (3) at (1, 1.75) {};
			\node [style=none] (4) at (0, 2.75) {};
			\node [style=none] (5) at (0, 0.75) {};
			\node [style=none] (6) at (2.35, 3.25) {$A$};
			\node [style=none] (7) at (1.25, 2.25) {$B$};
			\node [style=none] (8) at (-0.35, 1) {$A$};
		\end{pgfonlayer}
		\begin{pgfonlayer}{edgelayer}
			\draw (1.center) to (0.center);
			\draw [bend left=90, looseness=2.25] (0.center) to (3.center);
			\draw (3.center) to (2.center);
			\draw (5.center) to (4.center);
			\draw [bend left=90, looseness=1.75] (4.center) to (2.center);
		\end{pgfonlayer}
	\end{tikzpicture} = \begin{tikzpicture}
		\begin{pgfonlayer}{nodelayer}
			\node [style=none] (0) at (3.15, 2.25) {$A$};
			\node [style=none] (1) at (3, 3.5) {};
			\node [style=none] (2) at (3, 0.75) {};
		\end{pgfonlayer}
		\begin{pgfonlayer}{edgelayer}
			\draw (2.center) to (1.center);
		\end{pgfonlayer}
	\end{tikzpicture} \]

Similarly, an object $B$ is {\bf left dual} to an object $A$ if there exists maps 
$\eta: I \to B \ox A$, and $\epsilon: A \ox B \to I$ such that the corresponding snake equations hold.
	
\begin{definition}
	A {\bf compact closed category} (KCC) is a monoidal category in which each 
	object $A$ is equipped with chosen right and left duals.
\end{definition}
For a KCC which is also a SMC, a right dual of an object is also its left dual, and vice versa.  
Such a dual of an object $A$ is written as $A^*$. 
	
\begin{definition} \cite{AC04}
A {\bf $\dagger$-compact closed} category ($\dagger$-KCC) is a $\dagger$-symmetric monoidal category which is also a 
compact closed category such that for each object $\eta^\dag = c_\ox \epsilon$ 
(equivalently $\epsilon^\dag = \eta c_\ox$).
\end{definition}

Dagger compact closed categories were introduced by Coecke and Abramsky \cite{AC04} 
as an axiomatic framework for quantum information theory under the title strongly compact closed categories. 
The operator-state correspondence is straightforward in a $\dagger$-compact closed category. 
Every map $f: A \to B$ precisely corresponds to a state $\ceil{f}: I \to B \ox A^*$  
and an effect $\floor{f}: B^* \ox A \to I$ defined as follows:
\[  \ceil{f} := \begin{tikzpicture}
	\begin{pgfonlayer}{nodelayer}
		\node [style=circle, scale=1.5] (0) at (0, 2) {};
		\node [style=none] (10) at (0, 2) {$f$};
		\node [style=none] (1) at (0, 2.75) {};
		\node [style=none] (2) at (0, 1) {};
		\node [style=none] (3) at (0.75, 2.75) {};
		\node [style=none] (4) at (0.75, 1) {};
		\node [style=none] (5) at (-0.5, 2.5) {$A$};
		\node [style=none] (6) at (-0.5, 1.25) {$B$};
		\node [style=none] (7) at (1.25, 1.25) {$A^*$};
	\end{pgfonlayer}
	\begin{pgfonlayer}{edgelayer}
		\draw (0) to (2.center);
		\draw (0) to (1.center);
		\draw (3.center) to (4.center);
		\draw [bend left=90, looseness=1.75] (1.center) to (3.center);
	\end{pgfonlayer}
\end{tikzpicture} 
~~~~~~~~~~~~
\floor{f} :=  \begin{tikzpicture}
	\begin{pgfonlayer}{nodelayer}
		\node [style=circle] (0) at (0.75, 1.75) {$f$};
		\node [style=none] (1) at (0.75, 1) {};
		\node [style=none] (2) at (0.75, 2.75) {};
		\node [style=none] (3) at (0, 1) {};
		\node [style=none] (4) at (0, 2.75) {};
		\node [style=none] (5) at (1.25, 1.25) {$B$};
		\node [style=none] (6) at (1.25, 2.5) {$A$};
		\node [style=none] (7) at (-0.5, 2.5) {$B^*$};
	\end{pgfonlayer}
	\begin{pgfonlayer}{edgelayer}
		\draw (0) to (2.center);
		\draw (0) to (1.center);
		\draw (3.center) to (4.center);
		\draw [bend left=90, looseness=1.75] (1.center) to (3.center);
	\end{pgfonlayer}
\end{tikzpicture} \]

In CQM, $\dagger$-KCCs are the fundamental framework for quantum information theory and 
quantum computing.  We will discuss the presentation of quantum channels, quantum observables, 
and strong complementarity within these categories in the later sections.

\section{Examples}
In this section, we recall a few standard examples of $\dagger$-KCCs which are used in CQM. 

\subsection{Categories of sets and relations, $\Rel$}
\label{Sec: relations}
In CQM, the category of sets and relations is used as a non-standard model of 
quantum  information theory. This category is often studied in comparison with the category 
of Hilbert spaces and bounded linear maps to distinguish quantum versus non-quantum 
features \cite{HeV19, HeT15}. 
The category $\Rel$ is 
defined as follows:
\begin{description}
	\item[Objects:] Sets
	\item[Arrows:] $R:X \rightarrow Y$ where $R \subseteq X \times Y$ 
	\item[Identity maps:] $I_X:= \{ (x,x) | ~ \ x \in X \}$
	\item[Composition:] Suppose $A \xrightarrow{R} B \xrightarrow{S} C$, 
	then \[ RS: A \rightarrow C := \{ (x,z) | ~ \exists y \in Y, (x,y) \in R \text{ and } (y,z) \in S \} \]
	\item[Tensor product:] $A \otimes B := A \times B$  is the cartesian product on both sets 
	(and arrows), the associativity map is $a_\ox \subseteq ((A \times B) \times C) \times (A \times (B \times C ))
	:= \{ (((a,b),c),(a,(b,c))) | a \in A,b \in B, c \in C\}$, and the unit is the one element set $I = \{ \star \}$.
    \item[Symmetry map:] $c_\otimes \subseteq (A \times B) \times (B \otimes A) := \{ ((a,b),(b,a)) | (a,b) \in A \times B \} $
	\item[Dagger:] Given $R: A \to B$, $B \xrightarrow{R^\dagger} A := \{ (b,a) | (a,b) \in R \}$ is the converse relation of $f$.
\end{description}
Thus, $\Rel$ is a symmetric $\dagger$-monoidal category. In ${\sf Rel}$, every object is self-dual: $(\eta, \epsilon): A \dashvv A$ with 
$\eta : I \to A \times A := \{ ( \star, (x,x)) | ~ x \in A \}$, and $\epsilon$ is given by the converse relation. 
So, $\epsilon: A \times A \to I :=  \{ ( (x,x), \star) | ~ x \in A \}$.  This makes $\Rel$ a 
$\dagger$-compact closed category.

\subsection{Categories of Hilbert spaces and bounded linear maps, $\Hilb$}

The category of Hilbert Spaces and bounded linear maps, $\Hilb$, is a $\dagger$-SMC. 
The subcategory $\FHilb$ of $\Hilb$ consisting of only finite-dimensional 
Hilbert spaces is a $\dagger$-KCC. $\FHilb$ is used in CQM as a standard setting 
for studying processes in quantum information theory and quantum computing.

The category $\Hilb$ is defined as follows:
\begin{description}
	\item[Objects:] Hilbert spaces
	\item[Maps:] Bounded linear maps
	\item[Composition:] Usual composition of linear maps 
	\item[Tensor product:] Standard tensor product of of the underlying vector spaces. The tensor 
	product on the maps is given by the Kronecker product. The tensor unit is $\C$.  
	\item[Dagger:] Given any map $f: A \to B$, its adjoint $f^\dag: B \to A$ is defined as 
	the map satisfying the following equation for all $a: I \to A$, and $b: I \to B$.
   \[ \left< b f^\dag, a \right> = \left< b, af \right> \]
\end{description}
$\Hilb$ is a $\dagger$-symmetric monoidal category. In finite-dimensions, Hilbert Spaces are equipped 
with a {\bf compact structure}. Given any finite-dimensional 
Hilbert Space, $H^*$ refers to the dual space of all functionals $H \rightarrow H$. 
Recall that every finite-dimensional Hilbert space has an orthonormal basis. 
Suppose $\{ e_i \}_{i = 1}^{dimH}$ is an orthonormal basis for $H$, then: 
\[ \eta: 1 \rightarrow H^* \otimes H ; 1 \mapsto \sum_i e_i^* \otimes e_i 
\text{ and }
\epsilon: H \otimes H^* \rightarrow I ; e_i \otimes e_j^* \mapsto \delta_{ij} \] 
Note that one cannot define a counit $\epsilon$ as above for a Hilbert space since it 
could possibly result in an infinite sum for an infinite vector. 

The $\eta$ and $\epsilon$ satisfy snake equations as follows. Suppose $a \in H$ and $a = \sum_i a_i e_i$.
	\begin{align*}
	(1 \otimes \eta)(\epsilon \otimes 1) (\sum_i a_i e_i) &= (\epsilon \otimes 1) (\sum_i a_i e_i \otimes \sum_j e_j^* \otimes e_j) \\
	&= \sum_i \sum_j a_i ((\epsilon \otimes 1) e_i \otimes e_j^* \otimes e_j))  \text{ (By linearlity) }\\
	&= \sum_{ij} a_i (\delta_{ij} \otimes e_j) =  \sum_j a_j e_j = a
\end{align*}
The subcategory finite-dimensional Hilbert Spaces, $\FHilb$, is $\dagger$-compact closed.

\subsection{Categories of finite matrices over a commutative rig $R$, $\Mat(R)$}
\label{Sec: Mat(R)}
Consider the category $\Mat(R)$ defined as follows:
\begin{description}
	\item[Objects]: $n \in \N$;
	\item[Arrows:] $n \to^{M} m$, where $M$ is a $n \times m$ matrix over a commutative 
	rig\footnote{Ring without negatives, hence a semiring} $R$;
	\item[Composition:] matrix multiplication;
	\item[Tensor product:] $n \ox m := n \cdot m$; $M \ox N$ is the outer product of matrices. 
	\end{description} 
	$\Mat(R)$ is a strict symmetric monoidal category. 
	Every object in ${\sf Mat}(R)$ is self-dual: 
	$(\eta_n, \epsilon_n): n \dashvv n$, where the $\eta_n$ and $\epsilon_n$ are defined as follows:
	\[
	  \eta_n: 1 \to n \ox n := \sum_{i=1}^{n} (e_i \ox e_i) ~~~~ \epsilon_n: n \ox n \to 1 := \sum_{i=1}^{n} (e_i^T \ox e_i^T)
	\]
	Here $e_i$ is a standard  basis vector for ${\sf Mat}(R)$, and $e_i^T$ is the transpose 
	of $e_i$. For example, for $n=2$, $e_1 := [1~0]$, $e_2 := [0~1]$ and, thus,  
	$\eta_2:=[1~0~0~1]$. If $M: n \to m$ then $M^*: m \to n$, the dual of $M$, 
	$M^{*}$, is just the transpose, $M^T$. 
	
	When the commutative rig, $R$, has a conjugation, $\overline{(\_)}: R \to R$ such that 
	$\overline{r}+\overline{s} = \overline{r+s}$, $\overline{0} = 0$, $\overline{r}~\overline{s} = 
	\overline{rs}$ and $\overline{1} = 1$, then ${\sf Mat}(R)$ has a dagger given by 
	conjugate transpose.   In particular, this means ${\sf Mat}(\C)$, finite-dimensional matrices over the 
	complex numbers, in addition, has a dagger given by conjugate transpose.  
	In fact, the category $\Mat(\C)$ is equivalent to $\FHilb$ \cite[Example 1.34]{HeV19}.
	
	Note that the two element ordered set, $\mathbbm{2}$, with join as addition and meet as multiplication is a 
	rig and ${\sf Mat}(\mathbbm{2})$ is then equivalent to the category of finite sets and relations.  
	This equivalence can be turned into an isomorphism if one inflates ${\sf Mat}(R)$ 
	so that it has objects finite sets, $I,J \in {\sf Set}_f$, and maps matrices given by maps $M: I \x J \to R$.  

\section{Complete positivity}
In this section we discuss quantum channels in $\dagger$-SMCs and $\dagger$-KCCs. 
In quantum mechanics, the representation of a physical state can be either {\em pure} or {\em mixed}. 
A representation is mixed when it is a statistical ensemble of possible states of the system. 
If the representation is pure, then one knows the exact quantum state of the system.  

 
While a pure state is represented as a vector in a Hilbert Space, a mixed state is represented as a 
positive self-adjoint operator on the Hilbert Space. The mixed state formalism of 
quantum mechanics is very useful in practice since in an experimental setting our knowledge about the state of 
a quantum system is often limited. In the mixed state formalism, a quantum process is a 
completely positive map (sending positive self-adjoint operators to positive self-adjoint operators)  
which preserves trace. In this section, we discuss how completely positive maps and traces are  
abstracted in a categorical setting.

\subsection{The CPM construction}

We begin our discussion with linear operators on finite-dimensional Hilbert Spaces. 
A {\bf linear operator} is a linear map from a Hilbert space to itself. 
The space of all linear operators, $\mathcal{L}(H)$, on a finite-dimensional Hilbert Space 
$H$ is a $\C^\star$-algebra with the $\star$ on $\mathcal{L}(H)$ defined to be the adjoint operation. 
Considering the linear map in matrix form, the adjoint of the matrix is its conjugate transpose.
An element $a$ of a $\C^\star$ algebra is {\bf positive} if there exists an element $b$ in 
the algebra such that $a = b^\star b$. Note that a positive element is 
self-adjoint, that is $a^\star = a$. The mixed state of quantum system is given a by positive element of 
$\mathcal{L}(H)$ of norm 1. 

A linear map between two $\C^\star$ algebras is {\bf positive} 
if  it preserves the positive elements, that is, $f(x^\star x) = y^\star y$. A linear map is 
said to be {\bf completely positive} if it is positive and for all $n > 0$, the map $1_n \ox f : \C^{n \times n} 
\ox A \to \C^{n \times n} \ox B$ is positive. If the representation of a quantum state is mixed, 
then completely positive maps that preserve trace are used to represent quantum processes. 

The following theorem by Choi characterizes the form of completely positive maps:
\begin{theorem} \cite[Theorem 1]{Cho75} Let $H$ and $K$ be finite-dimensional Hilbert spaces.  
	A linear map $\phi: \mathcal{L}(H) \to \mathcal{L}(K)$ is completely 
	positive if and only if there exists a collection of  linear maps, $\{M_i | M_i \in \mathcal{L}(H,K)$ 
	such that for all $A \in \mathcal{L}(H)$, 
	\[ \phi(A) = \sum_i M_i^\dag A M_i \] 
	where $M_i^\dag$ is a adjoint of $M_i$.
\end{theorem}

The equation in the above theorem is referred to as the {\bf Kraus decomposition} of $\phi$, and the 
collection of maps, $\{M_i | M_i \in \mathcal{L}(H,K)$, are referred to as the {\bf Kraus operators}.  

Selinger \cite{Sel07} abstracted the notion of completely positive maps to 
$\dagger$-compact closed categories using the characterization discussed above.  
In a $\dagger$-compact closed category, a map $f: A \to A$ is {\bf positive} if there exists a $g: A \to B$ 
such that $f = g g^\dag$. A map $f: A^* \ox A \to B^* \ox B$ is {\bf completely positive} if $f$ is of the 
following form: 
\begin{equation}
	\label{eqn: CPM map}
 f = \begin{tikzpicture}
	\begin{pgfonlayer}{nodelayer}
		\node [style=twocircle] (0) at (-1, 3) {};
		\node [style=twocircle] (1) at (1, 3) {};
		\node [style=none, scale=1.5] (2) at (1, 3) {$g$};
		\node [style=none, scale=1.5] (3) at (-1, 3) {$\overline{g}$};
		\node [style=none] (4) at (-1.75, 2) {};
		\node [style=none] (5) at (-0.5, 2) {};
		\node [style=none] (6) at (0.5, 2) {};
		\node [style=none] (7) at (1.75, 2) {};
		\node [style=none] (8) at (-1, 4) {};
		\node [style=none] (9) at (1, 4) {};
		\node [style=none] (10) at (-1.5, 3.75) {$A^*$};
		\node [style=none] (11) at (1.5, 3.75) {$A$};
		\node [style=none] (12) at (-2.25, 1.75) {$B^*$};
		\node [style=none] (13) at (0.25, 2.75) {$E$};
		\node [style=none] (14) at (1.75, 1.25) {};
		\node [style=none] (15) at (2.25, 1.75) {$B$};
		\node [style=none] (16) at (-1.75, 1.25) {};
	\end{pgfonlayer}
	\begin{pgfonlayer}{edgelayer}
		\draw [in=90, out=-150] (1) to (6.center);
		\draw [bend left] (1) to (7.center);
		\draw [in=90, out=-30] (0) to (5.center);
		\draw [bend right] (0) to (4.center);
		\draw (8.center) to (0);
		\draw (9.center) to (1);
		\draw [bend right=90, looseness=1.75] (5.center) to (6.center);
		\draw (7.center) to (14.center);
		\draw (16.center) to (4.center);
	\end{pgfonlayer}
\end{tikzpicture} 
\end{equation} 
where $\overline{g} := g^{\dag *} = g^{*\dag}$. The object $E$ in the above diagram 
is interpreted as the {\em environment} which is {\em discarded} after the process $f$ is complete. More on 
environment and discarding will be discussed later in Section \ref{Sec: environment}. The map 
$g: A \to E \ox B$ can be interpreted as a Kraus operator and the diagram itself as a description of the 
Kraus decomposition. 

Selinger also introduced the completely positive maps (CPM) construction which reconciles
 the pure states and mixed states formalisms of quantum theory within the theory 
of $\dagger$-compact closed categories.  
Given a $\dagger$-compact closed category $\X$, the category CPM$(\X)$ consists of objects of the form 
$A^* \ox A$ for all $A \in \C$ and the maps are chosen to be the completely positive maps. CPM$(\X)$ is 
again a $\dagger$-compact closed category \cite[Theorem 4.20]{Sel07} with the same $\dagger$ functor 
as $\C$. The CPM construction applied to the category of finite dimensional Hilbert spaces produces 
a category containing mixed states and quantum processes.  

\subsection{The $\CP^\infty$ construction}

Note that the CPM construction uses the compact structure of $\dagger$-KCCs, thus is applicable only 
to the category of finite dimensional Hilbert spaces. 
Coecke and Heunen \cite{CoH16} generalized the CPM construction to $\dagger$-monoidal categories 
thereby eliminating the restriction on dimensions. Accordingly, their 
generalized construction is referred to as the $\CP^\infty$ construction.

In $\dagger$-SMCs, one cannot bend wires, hence the representation of a 
completely positive map as shown in equation \ref{eqn: CPM map} cannot be used. 
In order to define a completely positive map within a $\dagger$-SMC, Coecke and Heunen 
straightened the wire for the environment in diagram \ref{eqn: CPM map}, 
obtaining the following form for {\bf completely positive} maps in this setting:
\begin{equation}
	\label{eqn: CP map}
	 \begin{tikzpicture}
		\begin{pgfonlayer}{nodelayer}
			\node [style=twocircle] (0) at (1, -0.25) {};
			\node [style=twocircle] (1) at (1, 3.25) {};
			\node [style=none] (2) at (1, 3.25) {$g$};
			\node [style=none] (3) at (1, -0.25) {$g^\dag$};
			\node [style=none] (4) at (1.75, 0.75) {};
			\node [style=none] (5) at (0.25, 0.75) {};
			\node [style=none] (6) at (0.25, 2.25) {};
			\node [style=none] (7) at (1.75, 2.25) {};
			\node [style=none] (8) at (1, -1.25) {};
			\node [style=none] (9) at (1, 4.25) {};
			\node [style=none] (10) at (0.75, -1) {$A$};
			\node [style=none] (11) at (0.75, 4) {$A$};
			\node [style=none] (13) at (0, 1.5) {$E$};
			\node [style=none] (15) at (2, 0.25) {$B$};
			\node [style=none] (16) at (2, 2.75) {$B$};
			\node [style=circle, scale=3, dashed] (17) at (1.75, 1.5) {};
		\end{pgfonlayer}
		\begin{pgfonlayer}{edgelayer}
			\draw [in=90, out=-150] (1) to (6.center);
			\draw [in=90, out=-30, looseness=1.25] (1) to (7.center);
			\draw [in=-90, out=150] (0) to (5.center);
			\draw [in=270, out=30, looseness=1.25] (0) to (4.center);
			\draw (8.center) to (0);
			\draw (9.center) to (1);
			\draw (6.center) to (5.center);
		\end{pgfonlayer}
	\end{tikzpicture}	
\end{equation}
Within the dashed circle, `test maps' are plugged in as shown in the equation \ref{eqn: CP equivalence}. 
The map $g: A \to E \ox B$ is called a {\bf Kraus} map and $E$ is referred to as the 
environment or the ancillary system. 

Any two Kraus maps, $f: A \to E \ox B$, and $g: A \to F \ox B$ are said to be {\bf equivalent}, 
that is $f \sim g$, if for all $h: B \ox C \to D$, they satisfy the following equation:
\begin{equation}
\label{eqn: CP equivalence}
\begin{tikzpicture}
	\begin{pgfonlayer}{nodelayer}
		\node [style=twocircle] (0) at (1, -0.25) {};
		\node [style=twocircle] (1) at (1, 3.25) {};
		\node [style=none] (2) at (1, 3.25) {$f$};
		\node [style=none] (3) at (1, -0.25) {$f^\dag$};
		\node [style=none] (5) at (0.25, 0.75) {};
		\node [style=none] (6) at (0.25, 2.25) {};
		\node [style=none] (8) at (1, -1.25) {};
		\node [style=none] (9) at (1, 4.25) {};
		\node [style=none] (10) at (0.75, -1) {$A$};
		\node [style=none] (11) at (0.75, 4) {$A$};
		\node [style=none] (13) at (0, 1.5) {$E$};
		\node [style=none] (15) at (2.75, -1) {$B$};
		\node [style=none] (16) at (2.75, 4) {$B$};
		\node [style=twocircle] (18) at (1.75, 0.75) {};
		\node [style=twocircle] (19) at (1.75, 2.25) {};
		\node [style=none] (20) at (2.5, 4.25) {};
		\node [style=none] (21) at (2.5, -1.25) {};
		\node [style=none] (22) at (2, 1.5) {$D$};
		\node [style=none] (23) at (1.75, 2.25) {$h$};
		\node [style=none] (24) at (1.75, 0.75) {$h^\dag$};
	\end{pgfonlayer}
	\begin{pgfonlayer}{edgelayer}
		\draw [in=90, out=-150] (1) to (6.center);
		\draw [in=-90, out=150] (0) to (5.center);
		\draw (8.center) to (0);
		\draw (9.center) to (1);
		\draw (6.center) to (5.center);
		\draw [in=-90, out=15, looseness=1.25] (0) to (18);
		\draw [in=-15, out=90] (19) to (1);
		\draw (19) to (18);
		\draw [in=-90, out=45] (19) to (20.center);
		\draw [in=90, out=-30] (18) to (21.center);
	\end{pgfonlayer}
\end{tikzpicture} = \begin{tikzpicture}
	\begin{pgfonlayer}{nodelayer}
		\node [style=twocircle] (0) at (1, -0.25) {};
		\node [style=twocircle] (1) at (1, 3.25) {};
		\node [style=none] (2) at (1, 3.25) {$g$};
		\node [style=none] (3) at (1, -0.25) {$g^\dag$};
		\node [style=none] (5) at (0.25, 0.75) {};
		\node [style=none] (6) at (0.25, 2.25) {};
		\node [style=none] (8) at (1, -1.25) {};
		\node [style=none] (9) at (1, 4.25) {};
		\node [style=none] (10) at (0.75, -1) {$A$};
		\node [style=none] (11) at (0.75, 4) {$A$};
		\node [style=none] (13) at (0, 1.5) {$F$};
		\node [style=none] (15) at (2.75, -1) {$B$};
		\node [style=none] (16) at (2.75, 4) {$B$};
		\node [style=twocircle] (18) at (1.75, 0.75) {};
		\node [style=twocircle] (19) at (1.75, 2.25) {};
		\node [style=none] (20) at (2.5, 4.25) {};
		\node [style=none] (21) at (2.5, -1.25) {};
		\node [style=none] (22) at (2, 1.5) {$D$};
		\node [style=none] (23) at (1.75, 2.25) {$h$};
		\node [style=none] (24) at (1.75, 0.75) {$h^\dag$};
	\end{pgfonlayer}
	\begin{pgfonlayer}{edgelayer}
		\draw [in=90, out=-150] (1) to (6.center);
		\draw [in=-90, out=150] (0) to (5.center);
		\draw (8.center) to (0);
		\draw (9.center) to (1);
		\draw (6.center) to (5.center);
		\draw [in=-90, out=15, looseness=1.25] (0) to (18);
		\draw [in=-15, out=90] (19) to (1);
		\draw (19) to (18);
		\draw [in=-90, out=45] (19) to (20.center);
		\draw [in=90, out=-30] (18) to (21.center);
	\end{pgfonlayer}
\end{tikzpicture}
\end{equation}
For any Kraus map, $f: A \to E \ox B$, we will write its equivalence class as $[f]$.

\begin{definition} \cite[Definition 3]{CoH16} The {\bf $\CP^\infty$ construction} is defined as follows. 
Let $\C$ be any $\dagger$-monoidal category.
\begin{description}
	\item[Objects:] Same as objects as $\C$
	\item[Maps:] Equivalence classes of Kraus maps in $\C$
	\item[Identity maps:] Equivalence class of Kraus maps given by the left 
unitor, $[(u_\ox^l)^{-1} : A \to I \ox A]$. 
	\item[Composition:] The composition of two maps $f: A \to B$, and $g: B \to C \in \CP^\infty(\C)$  
	is given by composing their respective Kraus maps upto equivalence:
	\[ fg := \left[ \begin{tikzpicture}
		\begin{pgfonlayer}{nodelayer}
			\node [style=twocircle] (1) at (1, 3.25) {};
			\node [style=none] (2) at (1, 3.25) {$f$};
			\node [style=none] (5) at (0.25, 0.75) {};
			\node [style=none] (6) at (0.25, 2.25) {};
			\node [style=none] (9) at (1, 4.25) {};
			\node [style=none] (10) at (0.65, 1) {$E'$};
			\node [style=none] (11) at (0.75, 4) {$A$};
			\node [style=none] (13) at (0, 1) {$E$};
			\node [style=none] (15) at (1.85, 3) {$B$};
			\node [style=twocircle] (19) at (1.75, 2.25) {};
			\node [style=none] (22) at (1, 0.75) {};
			\node [style=none] (23) at (2.5, 0.75) {};
			\node [style=none] (24) at (1.75, 2.25) {$g$};
			\node [style=none] (25) at (2.75, 1) {$C$};
		\end{pgfonlayer}
		\begin{pgfonlayer}{edgelayer}
			\draw [in=90, out=-150] (1) to (6.center);
			\draw (9.center) to (1);
			\draw (6.center) to (5.center);
			\draw [in=-15, out=90] (19) to (1);
			\draw [in=90, out=-150, looseness=1.25] (19) to (22.center);
			\draw [in=90, out=-30, looseness=1.25] (19) to (23.center);
		\end{pgfonlayer}
	\end{tikzpicture} \right] \]   
\end{description} 
The composition is well-defined because if $f \sim f'$ and $g \sim g'$, then $fg \sim f'g'$.
\end{definition}

If $\C$ is a $\dagger$-SMC, then $\CP^\infty(\C)$ is a SMC \cite[Proposition 5]{CoH16}. 
Note that, in this case, $\CP^\infty(\C)$ does not have a $\dagger$-functor.
Moreover, if $\C$ is a $\dagger$-KCC, then $\CP^\infty(\C)$ isomorphic to CPM$(\C)$ 
\cite[Proposition 6]{CoH16}. 

The set of all bounded linear maps from a Hilbert space to itself is a Von Neumann algebra. 
The category of Von Neumann algebras and normal completely positive maps is isomorphic to $\CP^\infty(\Hilb)$ 
\cite[Theorem 13]{CoH16}.

\subsection{Environment structures}
\label{Sec: environment}

In the previous sections, we discussed the constructions for transforming a category of pure states 
and processes into a category of mixed states and completely positive maps. In this section we 
discuss characterizations of these constructions using the notion of {\em environment} \cite{CoH16,Coecke10,Coe13}. 
While pure states contain information {\em purely} about the system, 
mixed states contain information about the system and its environment together.
One can thus characterize the categories of mixed states by the presence of environment 
structure and discarding maps. 

\begin{definition}
An {\bf environment structure} for a $\dagger$-SMC, $\C_{pure}$, consists of a strict monoidal functor 
$F: \C_{pure} \to \C$ where $\C$ is a SMC  with designated map, $\gamma_A: A \to I$, called the 
{\bf discarding} map for all objects $A \in \C$ such that the following equations hold:
\begin{description}
\item[Env 1:] For all objects $A, B \in \C_{pure}$,  $\gamma_{F(A \ox B)} = \gamma_{F(A)} \ox \gamma_{F(B)} (u_\ox)_I $
\item[Env 2:] $\gamma_I = id_I$ (F(I) = I)
\item[Env 3:] For all Kraus maps $f: A \to X \ox B, g: A  \to Y \ox B \in \C_{pure}$, 

$f \sim g$ if and only if $F(f) (\gamma_{F(X)} \ox 1_B) = F(g) (\gamma_{F(Y)} \ox 1_B)$.
\end{description}
\end{definition}
The axioms for an environment structure are drawn as follows:
\[ {\bf [Env~1]:}~~~ \begin{tikzpicture}
	\begin{pgfonlayer}{nodelayer}
		\node [style=none] (0) at (-1, 0) {};
		\node [style=none] (1) at (0, 0) {};
		\node [style=none] (4) at (-0.8, -0.25) {};
		\node [style=none] (5) at (-0.2, -0.25) {};
		\node [style=none] (6) at (-0.6, -0.5) {};
		\node [style=none] (7) at (-0.4, -0.5) {};
		\node [style=none] (2) at (-0.5, 0) {};
		\node [style=none] (3) at (-0.5, 1) {};
		\node [style=none] (12) at (0.5, 0.75) {$F(A \ox B)$};
	\end{pgfonlayer}
	\begin{pgfonlayer}{edgelayer}
		\draw (3.center) to (2.center);
		\draw[thick] (0.center) to (1.center);
		\draw[thick] (4.center) to (5.center);
		\draw[thick] (6.center) to (7.center);
	\end{pgfonlayer}
\end{tikzpicture} = \begin{tikzpicture}
	\begin{pgfonlayer}{nodelayer}
		\node [style=none] (0) at (-1, 0) {};
		\node [style=none] (1) at (0, 0) {};
		\node [style=none] (14) at (-0.8, -0.25) {};
		\node [style=none] (15) at (-0.2, -0.25) {};
		\node [style=none] (16) at (-0.6, -0.5) {};
		\node [style=none] (17) at (-0.4, -0.5) {};
		\node [style=none] (2) at (-0.5, 0) {};
		\node [style=none] (3) at (-0.5, 1) {};
		\node [style=none] (4) at (0.5, 0) {};
		\node [style=none] (5) at (1.5, 0) {};
		\node [style=none] (6) at (0.7, -0.25) {};
		\node [style=none] (7) at (1.3, -0.25) {};
		\node [style=none] (8) at (0.9, -0.5) {};
		\node [style=none] (9) at (1.1, -0.5) {};
		\node [style=none] (10) at (1, 0) {};
		\node [style=none] (11) at (1, 1) {};
		\node [style=none] (12) at (-1.2, 0.75) {$F(A)$};
		\node [style=none] (13) at (1.5, 0.75) {$F(B)$};
	\end{pgfonlayer}
	\begin{pgfonlayer}{edgelayer}
		\draw (3.center) to (2.center);
		\draw[thick] (0.center) to (1.center);
		\draw[thick] (14.center) to (15.center);
		\draw[thick] (16.center) to (17.center);
		\draw (11.center) to (10.center);
		\draw[thick] (8.center) to (9.center);
		\draw[thick] (4.center) to (5.center);
		\draw[thick] (6.center) to (7.center);
	\end{pgfonlayer}
\end{tikzpicture}
~~~~~~~~
{\bf [Env~2]:} ~~~\begin{tikzpicture}
	\begin{pgfonlayer}{nodelayer}
		\node [style=none] (0) at (-1, 0) {};
		\node [style=none] (1) at (0, 0) {};
		\node [style=none] (4) at (-0.8, -0.25) {};
		\node [style=none] (5) at (-0.2, -0.25) {};
		\node [style=none] (6) at (-0.6, -0.5) {};
		\node [style=none] (7) at (-0.4, -0.5) {};
		\node [style=none] (2) at (-0.5, 0) {};
		\node [style=none] (3) at (-0.5, 1) {};
		\node [style=none] (12) at (-0.25, 0.75) {$I$};
	\end{pgfonlayer}
	\begin{pgfonlayer}{edgelayer}
		\draw (3.center) to (2.center);
		\draw[thick] (0.center) to (1.center);
		\draw[thick] (4.center) to (5.center);
		\draw[thick] (6.center) to (7.center);
	\end{pgfonlayer}
\end{tikzpicture} = id_I \] \[ {\bf [Env~3]:}~~~ f \sim g ~~~\Leftrightarrow~~~
\begin{tikzpicture}
	\begin{pgfonlayer}{nodelayer}
		\node [style=none] (0) at (-1, 0) {};
		\node [style=none] (1) at (0, 0) {};
		\node [style=none] (11) at (-0.8, -0.25) {};
		\node [style=none] (12) at (-0.2, -0.25) {};
		\node [style=none] (13) at (-0.6, -0.5) {};
		\node [style=none] (14) at (-0.4, -0.5) {};
		\node [style=none] (2) at (-0.5, 0) {};
		\node [style=none] (3) at (-0.5, 0.5) {};
		\node [style=circle, scale=2.5] (4) at (0.25, 1.25) {};
		\node [style=none] (5) at (1, 0.5) {};
		\node [style=none] (6) at (0.25, 2.25) {};
		\node [style=none] (7) at (0.25, 1.25) {$F(f)$};
		\node [style=none] (8) at (1, 2) {$F(A)$};
		\node [style=none] (9) at (-1.25, 0.5) {$F(X)$};
		\node [style=none] (10) at (1.75, -0.5) {$F(B)$};
		\node [style=none] (15) at (1, -0.75) {};
	\end{pgfonlayer}
	\begin{pgfonlayer}{edgelayer}
		\draw (3.center) to (2.center);
		\draw [thick] (0.center) to (1.center);
		\draw [thick] (11.center) to (12.center);
		\draw [thick] (13.center) to (14.center);
		\draw [in=90, out=-15, looseness=1.25] (4) to (5.center);
		\draw [in=90, out=-165, looseness=1.25] (4) to (3.center);
		\draw (6.center) to (4);
		\draw (5.center) to (15.center);
	\end{pgfonlayer}
\end{tikzpicture} = \begin{tikzpicture}
	\begin{pgfonlayer}{nodelayer}
		\node [style=none] (0) at (-1, 0) {};
		\node [style=none] (1) at (0, 0) {};
		\node [style=none] (11) at (-0.8, -0.25) {};
		\node [style=none] (12) at (-0.2, -0.25) {};
		\node [style=none] (13) at (-0.6, -0.5) {};
		\node [style=none] (14) at (-0.4, -0.5) {};
		\node [style=none] (2) at (-0.5, 0) {};
		\node [style=none] (3) at (-0.5, 0.5) {};
		\node [style=circle, scale=2.5] (4) at (0.25, 1.25) {};
		\node [style=none] (5) at (1, 0.5) {};
		\node [style=none] (6) at (0.25, 2.25) {};
		\node [style=none] (7) at (0.25, 1.25) {$F(g)$};
		\node [style=none] (8) at (1, 2) {$F(A)$};
		\node [style=none] (9) at (-1.25, 0.5) {$F(Y)$};
		\node [style=none] (10) at (1.75, -0.5) {$F(B)$};
		\node [style=none] (15) at (1, -0.75) {};
	\end{pgfonlayer}
	\begin{pgfonlayer}{edgelayer}
		\draw (3.center) to (2.center);
		\draw [thick] (0.center) to (1.center);
		\draw [thick] (11.center) to (12.center);
		\draw [thick] (13.center) to (14.center);
		\draw [in=90, out=-15, looseness=1.25] (4) to (5.center);
		\draw [in=90, out=-165, looseness=1.25] (4) to (3.center);
		\draw (6.center) to (4);
		\draw (5.center) to (15.center);
	\end{pgfonlayer}
\end{tikzpicture} \]

A process $f: A \to B$ in a $\dagger$-SMC, $\C_{pure}$, with an environment structure 
$(F: \C_{pure} \to \C, \envmap)$ is said to be {\bf normalised} if $F(f) \gamma_{F(B)} = \gamma_{F(A)}$. 
The normalised processes of $\C_{pure}$ form a sub-$\dagger$-SMC of $\C$ with the same 
environment structure. In that case, for the sub-$\dagger$-SMC, the discarding map $\gamma$ is a 
monoidal transformation for from the functor $F$ to the constant endofunctor that sends all objects to $I$ 
and all maps to $id_I$. 

\begin{definition}
A $\dagger$-SMC, $\C_{pure}$, with an environment structure $(F: \C_{pure} \to \C, \envmap)$ is said to allow 
{\bf purification} if  every map in $\C$ is of the following form for some Kraus map $f: A \to X \ox B \in \C_{pure}$
\[ {\bf[Env~4]}~~~ \begin{tikzpicture}
	\begin{pgfonlayer}{nodelayer}
		\node [style=none] (0) at (-1, 0) {};
		\node [style=none] (1) at (0, 0) {};
		\node [style=none] (11) at (-0.8, -0.25) {};
		\node [style=none] (12) at (-0.2, -0.25) {};
		\node [style=none] (13) at (-0.6, -0.5) {};
		\node [style=none] (14) at (-0.4, -0.5) {};
		\node [style=none] (2) at (-0.5, 0) {};
		\node [style=none] (3) at (-0.5, 0.5) {};
		\node [style=circle, scale=2.5] (4) at (0.25, 1.25) {};
		\node [style=none] (5) at (1, 0.5) {};
		\node [style=none] (6) at (0.25, 2.25) {};
		\node [style=none] (7) at (0.25, 1.25) {$F(f)$};
		\node [style=none] (8) at (1, 2) {$F(A)$};
		\node [style=none] (9) at (-1.25, 0.5) {$F(X)$};
		\node [style=none] (10) at (1.75, -0.5) {$F(B)$};
		\node [style=none] (15) at (1, -0.75) {};
	\end{pgfonlayer}
	\begin{pgfonlayer}{edgelayer}
		\draw (3.center) to (2.center);
		\draw [thick] (0.center) to (1.center);
		\draw [thick] (11.center) to (12.center);
		\draw [thick] (13.center) to (14.center);
		\draw [in=90, out=-15, looseness=1.25] (4) to (5.center);
		\draw [in=90, out=-165, looseness=1.25] (4) to (3.center);
		\draw (6.center) to (4);
		\draw (5.center) to (15.center);
	\end{pgfonlayer}
\end{tikzpicture}  \]
\end{definition}

The idea of {\bf [Env 3]} is that every process in $\C$ is equal to a pure process followed by discarding information. 

Every $\dagger$-SMC, $\C_{pure}$,  comes with a canonical environment structure which allows purification \cite[Theorem 16]{CoH16}:
\[ F: \C_{pure} \to \CP^\infty(\C_{pure}); A \to^f B \mapsto  A \to^{[f (u_\ox^l)^{-1}]} B \]
\[ \gamma_A : A \to I \in \CP^\infty(\C_{pure}) := [(u_\ox^r)^{-1}] \] 
The environment structure also allows for purification:  for each $[f]: A \to B$ in $\CP^\infty(\C_{pure})$, 
$F(f) (\gamma_I \ox 1_B) = [f]$. 

\begin{theorem}\cite[Theorem 15]{CoH16}
	For any $\dagger$-SMC, $\C_{pure}$, $\CP^\infty(\C_{pure}) \simeq \C$ 
	is an isomorphism of monoidal categories if $\C_{pure}$ is equipped 
	with an environment structure and purification. 
\end{theorem}

The above result can be extended to $\dagger$-compact closed categories  when 
the environment structure behaves coherently with the compact structure in the 
following sense:
\begin{definition}
	An {\bf environment structure} for a $\dagger$-KCC, $\C_{pure}$,  consists of a 
	strict $\dagger$-monoidal functor $F: \C_{pure} \to \C$, where $\C$ is also a $\dagger$-KCC
	and the following conditions hold:
	\begin{enumerate}[(i)]
	\item For all $A \in \C_{pure}$, $F(A^*) = F(A)^*$ 
	\item For all $A \in \C_{pure}$, $F(A)^* \simeq F(A)$
	\item For each object $X \in \C$, there exists a designated map $\gamma_X: X \to I$ such that 
	{\bf [Env 1]}-{\bf[Env 3]} holds and:
	\[ \begin{tikzpicture}
		\begin{pgfonlayer}{nodelayer}
			\node [style=none] (0) at (-1, 0) {};
			\node [style=none] (1) at (0, 0) {};
			\node [style=none] (4) at (-0.8, -0.25) {};
			\node [style=none] (5) at (-0.2, -0.25) {};
			\node [style=none] (6) at (-0.6, -0.5) {};
			\node [style=none] (7) at (-0.4, -0.5) {};
			\node [style=none] (2) at (-0.5, 0) {};
			\node [style=none] (3) at (-0.5, 1) {};
			\node [style=none] (12) at (1.25, -0.25) {$F(A)^*$};
			\node [style=none] (13) at (0.5, 1) {};
			\node [style=none] (14) at (0.5, -0.75) {};
			\node [style=none] (15) at (-1, 0.5) {};
			\node [style=none] (16) at (1, 0.5) {};
			\node [style=none] (17) at (1, 1.75) {};
			\node [style=none] (18) at (-1, 1.75) {};
			\node [style=none] (19) at (0.75, 0.75) {$F$};
		\end{pgfonlayer}
		\begin{pgfonlayer}{edgelayer}
			\draw (3.center) to (2.center);
			\draw [thick] (0.center) to (1.center);
			\draw [thick] (4.center) to (5.center);
			\draw [thick] (6.center) to (7.center);
			\draw [bend left=90, looseness=1.50] (3.center) to (13.center);
			\draw (14.center) to (13.center);
			\draw (15.center) to (16.center);
			\draw (16.center) to (17.center);
			\draw (17.center) to (18.center);
			\draw (18.center) to (15.center);
		\end{pgfonlayer}
	\end{tikzpicture} = \left( \begin{tikzpicture}
		\begin{pgfonlayer}{nodelayer}
			\node [style=none] (0) at (-1, 0) {};
			\node [style=none] (1) at (0, 0) {};
			\node [style=none] (4) at (-0.8, -0.25) {};
			\node [style=none] (5) at (-0.2, -0.25) {};
			\node [style=none] (6) at (-0.6, -0.5) {};
			\node [style=none] (7) at (-0.4, -0.5) {};
			\node [style=none] (2) at (-0.5, 0) {};
			\node [style=none] (3) at (-0.5, 1.25) {};
			\node [style=none] (8) at (0, 1.2) {$F(A)$};
		\end{pgfonlayer}
		\begin{pgfonlayer}{edgelayer}
			\draw (3.center) to (2.center);
			\draw [thick] (0.center) to (1.center);
			\draw [thick] (4.center) to (5.center);
			\draw [thick] (6.center) to (7.center);
		\end{pgfonlayer}
	\end{tikzpicture} \right)^\dag	 \] 
\end{enumerate}
The environment structure allows purification if {\bf[Env 5]} holds.
\end{definition}

Every $\dagger$-compact closed category comes with a canonical environment structure given 
as follows:
\[ F: \C_{pure} \to \text{CPM}(\C_{pure}); A \to^f B \mapsto (A^* \ox A) \to^{f^* \ox f}  B^* \ox B\]
\[ \gamma_A :=  \begin{tikzpicture}
	\begin{pgfonlayer}{nodelayer}
		\node [style=none] (0) at (-5.25, 6.75) {};
		\node [style=none] (1) at (-3.75, 6.75) {};
		\node [style=none] (2) at (-5.65, 6.5) {$A^*$};
		\node [style=none] (3) at (-3.5, 6.5) {$A$};
	\end{pgfonlayer}
	\begin{pgfonlayer}{edgelayer}
		\draw [bend right=90, looseness=2.25] (0.center) to (1.center);
	\end{pgfonlayer}
\end{tikzpicture} \]

The discarding map in the category CPM$(\FHilb)$ is referred to as the {\bf trace}. 
A map $g: A \to B$ in CPM$(\FHilb)$ as shown below is precisely a {\bf quantum channel} (completely positive and trace 
preserving) when $g \in \FHilb$ is normalized i.e, $F(g) \gamma_{F(E_\ox B)} = \gamma_{F(A)}$. 
\[ \begin{tikzpicture}
	\begin{pgfonlayer}{nodelayer}
		\node [style=twocircle] (0) at (-1, 3) {};
		\node [style=twocircle] (1) at (1, 3) {};
		\node [style=none, scale=1.5] (2) at (1, 3) {$g$};
		\node [style=none, scale=1.5] (3) at (-1, 3) {$\overline{g}$};
		\node [style=none] (4) at (-1.75, 2) {};
		\node [style=none] (5) at (-0.5, 2) {};
		\node [style=none] (6) at (0.5, 2) {};
		\node [style=none] (7) at (1.75, 2) {};
		\node [style=none] (8) at (-1, 4) {};
		\node [style=none] (9) at (1, 4) {};
		\node [style=none] (10) at (-1.35, 3.75) {$A^*$};
		\node [style=none] (11) at (1.35, 3.75) {$A$};
		\node [style=none] (12) at (-2.15, 1.75) {$B^*$};
		\node [style=none] (13) at (0.25, 2.75) {$E$};
		\node [style=none] (15) at (2, 1.75) {$B$};
		\node [style=none] (16) at (-1.75, 1) {};
		\node [style=none] (17) at (1.75, 1) {};
	\end{pgfonlayer}
	\begin{pgfonlayer}{edgelayer}
		\draw [in=90, out=-150] (1) to (6.center);
		\draw [bend left] (1) to (7.center);
		\draw [in=90, out=-30] (0) to (5.center);
		\draw [bend right] (0) to (4.center);
		\draw (8.center) to (0);
		\draw (9.center) to (1);
		\draw [bend right=90, looseness=1.75] (5.center) to (6.center);
		\draw (4.center) to (16.center);
		\draw (7.center) to (17.center);
	\end{pgfonlayer}
\end{tikzpicture}\]

The environment structure in the above example also allows for purification because every map in CPM$(\C_{pure})$ is of the 
following form for some $g: A \to E \ox B \in \C_{pure}$:
\[ \begin{tikzpicture}
	\begin{pgfonlayer}{nodelayer}
		\node [style=twocircle] (0) at (-1, 3) {};
		\node [style=twocircle] (1) at (1, 3) {};
		\node [style=none, scale=1.5] (2) at (1, 3) {$g$};
		\node [style=none, scale=1.5] (3) at (-1, 3) {$\overline{g}$};
		\node [style=none] (4) at (-1.75, 2) {};
		\node [style=none] (5) at (-0.5, 2) {};
		\node [style=none] (6) at (0.5, 2) {};
		\node [style=none] (7) at (1.75, 2) {};
		\node [style=none] (8) at (-1, 4) {};
		\node [style=none] (9) at (1, 4) {};
		\node [style=none] (10) at (-1.5, 3.75) {$A^*$};
		\node [style=none] (11) at (1.5, 3.75) {$A$};
		\node [style=none] (12) at (-2.25, 1.75) {$B^*$};
		\node [style=none] (13) at (0.25, 2.75) {$E$};
		\node [style=none] (14) at (1.75, 1.25) {};
		\node [style=none] (15) at (2.25, 1.75) {$B$};
		\node [style=none] (16) at (-1.75, 1.25) {};
	\end{pgfonlayer}
	\begin{pgfonlayer}{edgelayer}
		\draw [in=90, out=-150] (1) to (6.center);
		\draw [bend left] (1) to (7.center);
		\draw [in=90, out=-30] (0) to (5.center);
		\draw [bend right] (0) to (4.center);
		\draw (8.center) to (0);
		\draw (9.center) to (1);
		\draw [bend right=90, looseness=1.75] (5.center) to (6.center);
		\draw (7.center) to (14.center);
		\draw (16.center) to (4.center);
	\end{pgfonlayer}
\end{tikzpicture} = \begin{tikzpicture}
	\begin{pgfonlayer}{nodelayer}
		\node [style=none] (0) at (-1, 0) {};
		\node [style=none] (1) at (0, 0) {};
		\node [style=none] (11) at (-0.8, -0.25) {};
		\node [style=none] (12) at (-0.2, -0.25) {};
		\node [style=none] (13) at (-0.6, -0.5) {};
		\node [style=none] (14) at (-0.4, -0.5) {};
		\node [style=none] (2) at (-0.5, 0) {};
		\node [style=none] (3) at (-0.5, 0.5) {};
		\node [style=circle, scale=2.5] (4) at (0.25, 1.25) {};
		\node [style=none] (5) at (1, 0.5) {};
		\node [style=none] (6) at (0.25, 2.25) {};
		\node [style=none] (7) at (0.25, 1.25) {$F(g)$};
		\node [style=none] (8) at (1, 2) {$F(A)$};
		\node [style=none] (9) at (-1.25, 0.5) {$F(E)$};
		\node [style=none] (10) at (1.75, -0.5) {$F(B)$};
		\node [style=none] (15) at (1, -0.75) {};
	\end{pgfonlayer}
	\begin{pgfonlayer}{edgelayer}
		\draw (3.center) to (2.center);
		\draw [thick] (0.center) to (1.center);
		\draw [thick] (11.center) to (12.center);
		\draw [thick] (13.center) to (14.center);
		\draw [in=90, out=-15, looseness=1.25] (4) to (5.center);
		\draw [in=90, out=-165, looseness=1.25] (4) to (3.center);
		\draw (6.center) to (4);
		\draw (5.center) to (15.center);
	\end{pgfonlayer}
\end{tikzpicture}  \]

\section{Measurement and complementarity}
We covered quantum channels in the previous section. In this section we review the 
notions of quantum observables, measurement and complementarity categorically. 

\subsection{Quantum observables}
\label{Sec: observables}

A quantum observable is a physical property which can be measured. 
Frobenius algebras are one of the pillars of CQM because they are used to 
abstract the notion of quantum observables. 
In traditional quantum mechanics, an observable is 
represented by a self-adjoint operator a.k.a Hermitian operator 
on a Hilbert space. The set of all eigenvectors 
for an observable gives an orthogonal basis for the state space of the quantum system.
After a measurement, the state of the quantum system will be one of these basis states.
The eigenvalue corresponding to an eigenvector represents its probability amplitude 
i.e., the probablity that the quantum system will end up in the particular basis state after 
measurement. Frobenius algebras provide a neat abstraction of these ideas in 
$\dagger$-monoidal categories. In this section, we review Frobenius algebras and 
their correspondence to quantum observables. 

In a SMC, a monoid $\monoid{A}$ consists of an object $A$ with a multiplication map, 
$\mulmap{1.5}{white}: A \ox A \to A$, and a unit map, $\unitmap{1.5}{white}: I \to A$, 
such that the multiplication is assosciative (see diagram (a)) and the unit law holds 
(see diagram (b)).
\[ (a)~~~	
 \]
Similarly, a comonoid $\comonoid{A}$ consists of an object $A$ with a comultiplication 
$\comulmap{1.5}{white}: A \to A \ox A$, and a counit map $e: A \to I$ such that the 
coumutiplication is coassociative (vertical reflection of equation for associativity), 
and that the counit law holds (vertical reflection of equation for unit law).

\begin{definition}
In a SMC, a {\bf Frobenius algebra}, $\Frob{A}$, consists of a monoid $\monoid{A}$ and 
a comonoid $\comonoid{A}$ such that the multiplication and the comultiplication 
interacts as follows:
\[  %
 \]
\end{definition}
The above equation is referred to as the Frobenius law. It can be proven that the Frobenius 
law holds for a monoid and a comonoid on $A$ if and only if the following equation holds:
\[ %
 \]
The diagram on the right of the equation is fondly referred to as a `spider' in the CQM community. 

In an SMC, if an object $A$ is Frobenius algbera, then $A$ is a self-dual object with the following 
cup and cap: 
\[ \eta: I \to A \ox A := %
 \]

For a Frobenius algebra, $\Frob{A}$,  the monoid $\monoid{A}$, and the comonoid $\comonoid{A}$ are 
dual to one another by the means of the self-dual cup and cap: 
\[ %
 \] 
Conversely, if $A$ is a self-dual object $A \dashvv A$ and a monoid $\monoid{A}$, and there exists a 
map $\counitmap{1.5}{white}: A \to I$ such that the self-dual cup is given by $
%
 $, then $A$ has a Frobenius structure, see \cite[Proposition 5.16]{HeV19}. The map 
$\counitmap{1.5}{white}: A \to I$ is referred to the non-degenerate form of the Frobenius structure.

\begin{definition}
	In a $\dagger$-SMC, a {\bf $\dagger$-Frobenius algebra} is a Frobenius algebra $\Frob{A}$ such that 
	\[ \left( ~ %
	 \]
\end{definition}

In $\dagger$-KCCs, {\em pants algebra} provide an important class of examples for $\dagger$-FAs. 
To define pants algebra, note that in a KCC, for each object $A$, there exists a monoid on 
$A \ox A^*$ with the following multiplication and unit respectively:
\[ %
 \]
The above algebra is referred to as a {\bf pants algebra} due to the shape of its multiplication map. 

Every monoid in a $KCC$ embeds into the pants monoid:
\begin{lemma}
\label{Lemma: monoidal pants embedding}
In a KCC, every monoid, $(A, \mulmap{1.5}{white}, \unitmap{1.5}{white})$, embeds into its 
pants monoid on $A^* \ox A$ by means of the following monoid morphism:
\[ %
  \]  
\end{lemma}

\begin{lemma}
In a $\dagger$-KCC, pants algebra are $\dagger$-Frobenius.
\end{lemma}

Applying the operator-state duality to pants $\dagger$-algebra we note the following. Let $f, g: A \to A$ 
be any two processes in a $\dagger$-KCC. Multliplying the states corresponding to these processes 
using the pants monoid gives the state corresponding to composition $fg$, see equation $(a)$. Moreover, 
applying the counit to the state of $f$ gives the trace of $f$, see diagram $(b)$.
\[ (a) ~~~ %
 \] 

Let us review a concrete example of pants algebra. Consider
the space of $n \times n$ complex matrices written as $M_n$ 
which gives the pants algebra over the space $\C^n$ in $\FHilb$
 \cite[Example 4.12]{HeV19}. We will show that $M_n$ is the pants algebra for $\C^n$.
The space of $n \times n$ complex matrices is a Hilbert space 
with the inner product given by $\langle A | B \rangle := Tr(A^\dagger B)$ ($A^\dagger$ is 
conjugate transpose of $A$) and comes with a canonical basis $\{ e_{ij} | i,j=1,\cdots,n \}$. 
The basis $e_{ij}$ is a $n \times n$ matrix with zero for all entries except the entry $(i,j)$. 

The algebra of $n \times n$ complex matrices has a $\dagger$ Frobenius structure in $\FHilb$:

{\bf Monoid structure:} The multiplication, $m$, is given by matrix multiplication 
and the unit, $u$, is the $n \times n$ identity matrix.

{\bf Comonoid structure:} Define the counit map $e := u^\dagger$. Then, by definition of $\dagger$ for 
in Hilbert spaces, $ \langle e(e_{ij}) | 1 \rangle = \langle (e_{ij}) | e^\dagger(1) \rangle = \langle e_{ij} 
| u (1) \rangle = \langle e_{ij} | I_A \rangle = \delta_{ij} = Tr(e_{ij})$. By extension of linearity, for 
all $A \in M_n$, $e(A) = Tr(A)$.

Similarly, define the comultiplication to be, $d := m^\dagger$. 
By the definition of $\dagger$ for $\mathsf{FHilb}$, $\langle m(e_{ij} \otimes e_{kl}) | 
e_{pq} \rangle = \langle e_{ij} \otimes e_{kl} | m^\dagger(e_{pq}) \rangle$. 
As before, in order to derive the definition of $m^\dagger$ we expand the equation on both sides.

\begin{align*}
 \langle m^\dagger(e_{ij}) | e_{kl} \otimes e_{pq} \rangle &= \langle e_{ij} | m(e_{kl} \otimes e_{pq}) \rangle \\
&= \langle e_{ij} | \delta_{lp} e_{kq} \rangle \\
&= \delta_{lp} Tr(e_{ij}^* e_{kq}) \\
&= \delta_{lp} Tr(e_{ij}^* e_kq) \\
&=\delta_{lp} \delta_{ik} \delta_{jq}. 
\end{align*}
By defining $m^\dagger(e_{ij})  := \sum_l e_{il} \otimes e_{lj}$, $\langle m^\dagger(e_{ij}) | e_{kl} \otimes e_{pq} \rangle$ evaluates to $\delta_{lp} \delta_{ik} \delta_{jq}$.

A $\dagger$-Frobenius algebra in a symmetric $\dagger$-SMC is said to be {\bf special} if $(a)$ holds, 
{\bf commutative} if $(b)$ holds, and {\bf symmetric} if $(c)$ holds:
\[ (a)~~~ \begin{tikzpicture}
	\begin{pgfonlayer}{nodelayer}
		\node [style=circle] (0) at (0, 4) {};
		\node [style=circle] (1) at (0, 3) {};
		\node [style=none] (2) at (0, 2.25) {};
		\node [style=none] (3) at (0, 4.75) {};
	\end{pgfonlayer}
	\begin{pgfonlayer}{edgelayer}
		\draw [bend left=60, looseness=1.25] (0) to (1);
		\draw [bend right=60, looseness=1.25] (0) to (1);
		\draw (1) to (2.center);
		\draw (0) to (3.center);
	\end{pgfonlayer}
\end{tikzpicture} = \begin{tikzpicture}
	\begin{pgfonlayer}{nodelayer}
		\node [style=none] (4) at (1, 4.75) {};
		\node [style=none] (5) at (1, 2.25) {};
	\end{pgfonlayer}
	\begin{pgfonlayer}{edgelayer}
		\draw (4.center) to (5.center);
	\end{pgfonlayer}
\end{tikzpicture}  ~~~~~~~~ (b) ~~~ \begin{tikzpicture}
	\begin{pgfonlayer}{nodelayer}
		\node [style=circle] (0) at (0.25, 4) {};
		\node [style=none] (1) at (0.25, 3.25) {};
		\node [style=none] (2) at (0.75, 5.5) {};
		\node [style=none] (3) at (-0.25, 5.5) {};
	\end{pgfonlayer}
	\begin{pgfonlayer}{edgelayer}
		\draw [in=45, out=-90, looseness=1.75] (3.center) to (0);
		\draw [in=135, out=-90, looseness=1.75] (2.center) to (0);
		\draw (0) to (1.center);
	\end{pgfonlayer}
\end{tikzpicture} = \begin{tikzpicture}
	\begin{pgfonlayer}{nodelayer}
		\node [style=circle] (0) at (0.25, 4.25) {};
		\node [style=none] (1) at (0.25, 3.25) {};
		\node [style=none] (2) at (-0.5, 5.5) {};
		\node [style=none] (3) at (1, 5.5) {};
	\end{pgfonlayer}
	\begin{pgfonlayer}{edgelayer}
		\draw [bend left] (3.center) to (0);
		\draw [bend right] (2.center) to (0);
		\draw (0) to (1.center);
	\end{pgfonlayer}
\end{tikzpicture}  ~~~~~~~~ (c) ~~~ \begin{tikzpicture}
	\begin{pgfonlayer}{nodelayer}
		\node [style=circle] (0) at (0.25, 4) {};
		\node [style=circle] (1) at (0.25, 3.25) {};
		\node [style=none] (2) at (0.75, 5.5) {};
		\node [style=none] (3) at (-0.25, 5.5) {};
	\end{pgfonlayer}
	\begin{pgfonlayer}{edgelayer}
		\draw [in=45, out=-90, looseness=1.75] (3.center) to (0);
		\draw [in=135, out=-90, looseness=1.75] (2.center) to (0);
		\draw (0) to (1);
	\end{pgfonlayer}
\end{tikzpicture} = \begin{tikzpicture}
	\begin{pgfonlayer}{nodelayer}
		\node [style=circle] (0) at (0.25, 4.25) {};
		\node [style=circle] (1) at (0.25, 3.25) {};
		\node [style=none] (2) at (-0.5, 5.5) {};
		\node [style=none] (3) at (1, 5.5) {};
	\end{pgfonlayer}
	\begin{pgfonlayer}{edgelayer}
		\draw [bend left] (3.center) to (0);
		\draw [bend right] (2.center) to (0);
		\draw (0) to (1);
	\end{pgfonlayer}
\end{tikzpicture} \] 

Commutativity is  a stronger condition than the symmetry. For example, the pants algebra 
is non-commutative but symmetric. 

The connection between Frobenius algebras and quantum observables given by the 
fact that in the $\FHilb$ category, every special commutative $\dagger$-FA ($\dagger$-SCFA) precisely 
corresponds to an orthonormal basis. This correspondence arises from the 
notion of {\bf classical states} for a $\dagger$-FA: the states, $a: I \to A$, which can be copied, and 
deleted as shown below. 
\[ \text{copy:}~~~\begin{tikzpicture}
	\begin{pgfonlayer}{nodelayer}
		\node [style=none] (0) at (0, 5.25) {};
		\node [style=none] (1) at (-0.5, 4.5) {};
		\node [style=none] (2) at (0.5, 4.5) {};
		\node [style=circle] (3) at (0, 3.75) {};
		\node [style=none] (4) at (-0.5, 2.75) {};
		\node [style=none] (5) at (0.5, 2.75) {};
		\node [style=none] (6) at (0, 4.5) {};
		\node [style=none] (7) at (0, 4.75) {$a$};
	\end{pgfonlayer}
	\begin{pgfonlayer}{edgelayer}
		\draw (1.center) to (2.center);
		\draw (2.center) to (0.center);
		\draw (0.center) to (1.center);
		\draw [bend right] (3) to (4.center);
		\draw [bend left] (3) to (5.center);
		\draw (6.center) to (3);
	\end{pgfonlayer}
\end{tikzpicture} = \begin{tikzpicture}
	\begin{pgfonlayer}{nodelayer}
		\node [style=none] (0) at (-0.5, 5.25) {};
		\node [style=none] (1) at (-1, 4.5) {};
		\node [style=none] (2) at (0, 4.5) {};
		\node [style=none] (4) at (-0.5, 2.75) {};
		\node [style=none] (6) at (-0.5, 4.5) {};
		\node [style=none] (7) at (-0.5, 4.75) {$a$};
		\node [style=none] (8) at (0.75, 5.25) {};
		\node [style=none] (9) at (0.25, 4.5) {};
		\node [style=none] (10) at (1.25, 4.5) {};
		\node [style=none] (11) at (0.75, 2.75) {};
		\node [style=none] (12) at (0.75, 4.5) {};
		\node [style=none] (13) at (0.75, 4.75) {$a$};
	\end{pgfonlayer}
	\begin{pgfonlayer}{edgelayer}
		\draw (1.center) to (2.center);
		\draw (2.center) to (0.center);
		\draw (0.center) to (1.center);
		\draw (6.center) to (4.center);
		\draw (9.center) to (10.center);
		\draw (10.center) to (8.center);
		\draw (8.center) to (9.center);
		\draw (12.center) to (11.center);
	\end{pgfonlayer}
\end{tikzpicture}
~~~~~~~~~~~ \text{delete:}~~~\begin{tikzpicture}
	\begin{pgfonlayer}{nodelayer}
		\node [style=none] (0) at (0, 5.25) {};
		\node [style=none] (1) at (-0.5, 4.5) {};
		\node [style=none] (2) at (0.5, 4.5) {};
		\node [style=circle] (3) at (0, 3) {};
		\node [style=none] (6) at (0, 4.5) {};
		\node [style=none] (7) at (0, 4.75) {$a$};
	\end{pgfonlayer}
	\begin{pgfonlayer}{edgelayer}
		\draw (1.center) to (2.center);
		\draw (2.center) to (0.center);
		\draw (0.center) to (1.center);
		\draw (6.center) to (3);
	\end{pgfonlayer}
\end{tikzpicture} = id_I\] 

\begin{theorem}
	\label{Theorem: scfa}
	\cite[Theorem 5.1]{CPV12}
In $\FHilb$, the set of classical states for a $\dagger$-SCFA on an Hilbert space $H$ 
precisely corresponds to an orthonormal basis for $H$.
\end{theorem}
Hence, there exists a bijective correspondence between $\dagger$-SCFAs and 
orthonormal bases. The basis states are the only 
states of the quantum system that can be copied and deleted, in other words, 
a classical state. Hence, $\dagger$-SCFA are also referred to as {\em classical structures}. 
In the previous theorem if we drop the keyword special, then corresponding basis is orthogonal. 

A quantum measurement is the process of extracting classical data 
(can be copied and deleted) from a quantum state. 
Categorically, Coecke and Pavlovic \cite{CoP07} described a ``demolition'' measurement in a 
$\dagger$-monoidal category as a process, $m: A \to X$, with $m^\dagger m = 1_X$, to a 
special commutative $\dagger$-Frobenius algebra, $X$. The object $A$ refers to the 
state space of a quantum system. Demolition means that we ignore the resulting state of the 
quantum system after measurement and preserve only the classical data. 

In Theorem \ref{Theorem: scfa}, we saw that $\dagger$-SCFA model classical data. 
However, Frobenius algebras which are non-commutative but symmetric model quantum information:  
\begin{theorem}\cite{Vic10}
In $\FHilb$, every special symmetric $\dagger$-FA precisely corresponds to a 
$\C^*$-algebra.
\end{theorem}
The proof of the above theorem relies on Lemma \ref{Lemma: monoidal pants embedding}.
Due to their correspondence to $\C^*$-algebras, 
special symmetric $\dagger$-FAs are also referred to as {\em quantum algebras}. Note that 
quantum algebras are non-commutative. 

\subsection{Strong complementarity}
\label{Sec: strong comp}

Bohr's complementarity \cite{Gri18} is a key feature that distinguishes quantum mechanics from classical mechanics. 
Two quantum observables are said to be complementary if measuring one observable leads to 
maximum uncertainty regarding the value of the other. An example of complementary observables is 
position and momentum of an electron. All physical properties occur in complementary pairs due to the 
wave and the particle nature of matter. 

In CQM, complementary interaction of two observables  are axiomatized using Hopf algebras \cite{Abe04}, 
which are bialgebras with an antipode: 

\begin{definition}\cite{Han08}
In a SMC, a {\bf bialgebra} $\bialg{A}$ consists of a monoid $\monoid{A}$ and a comonoid $\comonoidb{A}$ 
satisfying the following equations: 
\[ (a)~~~  \begin{tikzpicture}[yscale=-1]
	\begin{pgfonlayer}{nodelayer}
		\node [style=circle] (0) at (0, 4) {};
		\node [style=black] (1) at (0, 5) {};
		\node [style=none] (2) at (-0.5, 3) {};
		\node [style=none] (3) at (0.5, 3) {};
	\end{pgfonlayer}
	\begin{pgfonlayer}{edgelayer}
		\draw [bend right, looseness=1.25] (0) to (2.center);
		\draw [bend left, looseness=1.25] (0) to (3.center);
		\draw (1) to (0);
	\end{pgfonlayer}
\end{tikzpicture} = \begin{tikzpicture}[yscale=-1]
	\begin{pgfonlayer}{nodelayer}
		\node [style=black] (1) at (0, 5) {};
		\node [style=none] (3) at (0, 3) {};
		\node [style=black] (4) at (0.75, 5) {};
		\node [style=none] (5) at (0.75, 3) {};
	\end{pgfonlayer}
	\begin{pgfonlayer}{edgelayer}
		\draw (3.center) to (1);
		\draw (5.center) to (4);
	\end{pgfonlayer}
\end{tikzpicture} ~~~~~~~~~~~~ (b)~~~ \begin{tikzpicture}
	\begin{pgfonlayer}{nodelayer}
		\node [style=black] (0) at (0, 4) {};
		\node [style=circle] (1) at (0, 5) {};
		\node [style=none] (2) at (-0.5, 3) {};
		\node [style=none] (3) at (0.5, 3) {};
	\end{pgfonlayer}
	\begin{pgfonlayer}{edgelayer}
		\draw [bend right, looseness=1.25] (0) to (2.center);
		\draw [bend left, looseness=1.25] (0) to (3.center);
		\draw (1) to (0);
	\end{pgfonlayer}
\end{tikzpicture} = \begin{tikzpicture}
	\begin{pgfonlayer}{nodelayer}
		\node [style=circle] (1) at (0, 5) {};
		\node [style=none] (3) at (0, 3) {};
		\node [style=circle] (4) at (0.75, 5) {};
		\node [style=none] (5) at (0.75, 3) {};
	\end{pgfonlayer}
	\begin{pgfonlayer}{edgelayer}
		\draw (3.center) to (1);
		\draw (5.center) to (4);
	\end{pgfonlayer}
\end{tikzpicture} ~~~~~~~~~~~ (c)~~~ \begin{tikzpicture}
	\begin{pgfonlayer}{nodelayer}
		\node [style=circle] (1) at (0, 5) {};
		\node [style=black] (2) at (0, 3) {};
	\end{pgfonlayer}
	\begin{pgfonlayer}{edgelayer}
		\draw (1) to (2);
	\end{pgfonlayer}
\end{tikzpicture} = id_I  ~~~~~~~~~~~~ (d)~~~ \begin{tikzpicture}
	\begin{pgfonlayer}{nodelayer}
		\node [style=black] (0) at (0, 5.25) {};
		\node [style=black] (1) at (1.25, 5.25) {};
		\node [style=circle] (2) at (0, 3.75) {};
		\node [style=circle] (3) at (1.25, 3.75) {};
		\node [style=none] (4) at (0, 3) {};
		\node [style=none] (5) at (1.25, 3) {};
		\node [style=none] (6) at (0, 6) {};
		\node [style=none] (7) at (1.25, 6) {};
		\node [style=none] (8) at (1.25, 3) {};
	\end{pgfonlayer}
	\begin{pgfonlayer}{edgelayer}
		\draw (1) to (2);
		\draw (0) to (3);
		\draw [bend right=45] (0) to (2);
		\draw [bend left=45] (1) to (3);
		\draw (3) to (8.center);
		\draw (2) to (4.center);
		\draw (6.center) to (0);
		\draw (7.center) to (1);
	\end{pgfonlayer}
\end{tikzpicture} = \begin{tikzpicture}
	\begin{pgfonlayer}{nodelayer}
		\node [style=none] (4) at (0, 3) {};
		\node [style=none] (6) at (0, 6) {};
		\node [style=none] (7) at (1.5, 6) {};
		\node [style=none] (8) at (1.5, 3) {};
		\node [style=black] (9) at (0.75, 4) {};
		\node [style=circle] (10) at (0.75, 5.25) {};
	\end{pgfonlayer}
	\begin{pgfonlayer}{edgelayer}
		\draw [in=90, out=-165] (9) to (4.center);
		\draw [in=90, out=-15] (9) to (8.center);
		\draw (9) to (10);
		\draw [in=255, out=15] (10) to (7.center);
		\draw [in=285, out=165] (10) to (6.center);
	\end{pgfonlayer}
\end{tikzpicture} \] 
\end{definition}

The final equation is often referred to as the {\em bialgebra rule}. 
Observing equations $(a)$-$(d)$, we note that the multiplication map acts as a comonoid morphism for the 
tensor comonoid $A \ox A$. Equivalently, the comultiplication acts as a monoid morphism for the 
tensor monoid $A \ox A$. The following lemma gives alternate descriptions for a bialgebra:
\begin{lemma} In an SMC, the following are equivalent:
\begin{enumerate}[(i)]
\item  $\bialg{A}$ is a bialgebra.
\item $\monoid{A}$ is a monoid and $\comonoidb{A}$ is a comonoid such that $\comulmap{1.5}{black}$ is a monoid morphism 
for the tensor monoid $(A\ox A, \twinmul{1}{white}, \twinunit{1}{white})$.
\item $\monoid{A}$ is a monoid and $\comonoidb{A}$ is a comonoid such that $\mulmap{1.5}{white}$ 
is a comonoid morphism for the tensor comonoid $(A\ox A, \twincomul{1}{black}, \twincounit{1}{black})$.
\end{enumerate}
\end{lemma}

The following are a few examples of bialgebras:
\begin{itemize}
\item In a category with biproducts, every object has a bialgebra structure given by the product and the coproduct maps:
\[ A + A \to^{\begin{bmatrix} id_A \\ id_A \end{bmatrix}} A ~~~~~~~~ 0 \to^{0_{0,A}} A 
~~~~~~~~ A \to^{\left< id_A, id_A \right>} A + A~~~~~~~~ A \to^{0_A,0} 0 \]
The category of vector spaces and linear maps has biproducts: `+' is given by direct sum and $0$ is the trivial vector space.
\item Consider the category of sets and functions, ${\sf Set}$. 
The category ${\sf Set}$ does not have biproducts, however, 
it is an SMC with the tensor product given by the cartesian product and the 
unit object being a chosen singleton set. In ${\sf Set}$, every monoid $(A, \circ , u)$ (that is, 
$A$ is a set with a binary operation $\circ$ and a unit $u$),  
is a bialgebra with the comonoid structure given 
as follows:
\[ \text{ for all } a \in A, ~ \comulmap{1.5}{black}: A \to A \times A;  a \mapsto (a,a) 
~~~~~~~~  \counitmap{1.5}{black}: A \to \{*\} ; a \mapsto * \]
\end{itemize}

We are now ready to define Hopf algebras:
\begin{definition} \cite{Blu96, BlS04}
In a SMC, a {\bf Hopf algebra} is a bialgbera $\bialg{A}$ with a map $s: A \to A$ referred to as an {\bf antipode} such that:
\[ \begin{tikzpicture}
	\begin{pgfonlayer}{nodelayer}
		\node [style=black] (0) at (0, 3.75) {};
		\node [style=circle] (1) at (0, 1.75) {};
		\node [style=none] (2) at (0, 4.5) {};
		\node [style=none] (3) at (0, 1) {};
		\node [style=onehalfcircle] (4) at (-0.5, 2.75) {};
		\node [style=none] (5) at (-0.5, 2.75) {$s$};
	\end{pgfonlayer}
	\begin{pgfonlayer}{edgelayer}
		\draw (1) to (3.center);
		\draw (2.center) to (0);
		\draw [bend right] (0) to (4);
		\draw [bend right] (4) to (1);
		\draw [bend right=60] (1) to (0);
	\end{pgfonlayer}
\end{tikzpicture} = \begin{tikzpicture}
	\begin{pgfonlayer}{nodelayer}
		\node [style=black] (0) at (-0.5, 3.25) {};
		\node [style=circle] (1) at (-0.5, 2.25) {};
		\node [style=none] (2) at (-0.5, 4.5) {};
		\node [style=none] (3) at (-0.5, 1) {};
	\end{pgfonlayer}
	\begin{pgfonlayer}{edgelayer}
		\draw (1) to (3.center);
		\draw (2.center) to (0);
	\end{pgfonlayer}
\end{tikzpicture} = \begin{tikzpicture}
	\begin{pgfonlayer}{nodelayer}
		\node [style=black] (0) at (-0.5, 3.75) {};
		\node [style=circle] (1) at (-0.5, 1.75) {};
		\node [style=none] (2) at (-0.5, 4.5) {};
		\node [style=none] (3) at (-0.5, 1) {};
		\node [style=onehalfcircle] (4) at (0, 2.75) {};
		\node [style=none] (5) at (0, 2.75) {$s$};
	\end{pgfonlayer}
	\begin{pgfonlayer}{edgelayer}
		\draw (1) to (3.center);
		\draw (2.center) to (0);
		\draw [bend left] (0) to (4);
		\draw [bend left] (4) to (1);
		\draw [bend left=60] (1) to (0);
	\end{pgfonlayer}
\end{tikzpicture} \]
\end{definition}
The above equation is referred to as the {\em Hopf law}. 
The following are a few examples of Hopf algberas:

\begin{itemize}
\item  In the examples of bialgebras, we saw that every monoid $(A, \circ, e)$ in the category ${\sf Set}$ 
is a bialgbera. If this monoid is also a group, that is, each element of $A$ is equipped with an inverse, 
then the bialgebra is Hopf with anitpode $s: A \to A$ defined to be, for all $a \in A$, $s(a) := a^{-1}$.  

\item Another example of a Hopf algbera is a group $K$-algebra \cite{Maj00}. Given a  
finite group $(G, \circ, u)$,  (that is, $G$ is a set with a binary operation $\circ$, 
a unit $u$ and an inverse for each element) and a field $K$, it is 
the free  $K$-vector space, $K[G]$, over $G$ with basis $\{e_g\}_{g \in G}$ .  
$K[G]$ is a Hopf algebra with the following $K$-linear maps:  for all $g,h \in G$, 
\[ m: K[G] \ox K[G] \to K[G];  e_g \ox e_h \mapsto e_{g \circ h}
~~~~~~~~ u: K \to K[G]; 1 \mapsto e_u \] 
\[ d: K[G] \to K[G] \ox K[G]; e_g \mapsto e_g \ox e_g 
~~~~~~~~ e: K[G] \to K ; e_g \mapsto 1 \]
Finally, the antipode is:
\[ s: K[G] \to K[G]; e_g \mapsto e_{g^{-1}} \]

\end{itemize}

The bialgebras and the Frobenius algberas capture interactions between 
a monoid and comonoid on a single underlying object. However, the 
axioms for a bialgbera are characterized by disconnected diagrams while the 
axioms for a Frobenius algebra are characterized by connected diagrams. 
The Hopf law depicts the maximum disconnect and can be interpreted as 
as a process in which there is no information flow from the 
multiplication to the comultiplication. This makes Hopf algberas 
an appealing algebraic structure to model complementarity of 
two quantum observables \cite{CoD11}. 

It is worthwhile to note that an object, $A$, which is both a Frobenius algbera and 
a bialgebra is trivial i.e., $A \simeq I$ :
\begin{lemma} \cite[Theorem 6.23]{HeV19} In a $\dagger$-SMC, suppose 
a monoid $\monoid{A}$ and a comonoid $\comonoidb{A}$ is both a 
bialgbera and Frobenius algebra, then:  
\[ \begin{tikzpicture}
	\begin{pgfonlayer}{nodelayer}
		\node [style=black] (1) at (0.25, 4) {};
		\node [style=none] (2) at (0.25, 5) {};
		\node [style=circle] (3) at (0.25, 3) {};
		\node [style=none] (4) at (0.25, 2) {};
	\end{pgfonlayer}
	\begin{pgfonlayer}{edgelayer}
		\draw (2.center) to (1);
		\draw (4.center) to (3);
	\end{pgfonlayer}
\end{tikzpicture} = \begin{tikzpicture}
	\begin{pgfonlayer}{nodelayer}
		\node [style=none] (2) at (0.25, 2) {};
		\node [style=none] (4) at (0.25, 5) {};
	\end{pgfonlayer}
	\begin{pgfonlayer}{edgelayer}
		\draw (2.center) to (4.center);
	\end{pgfonlayer}
\end{tikzpicture}  \] 
\end{lemma}
\begin{proof}
\[ \begin{tikzpicture}
	\begin{pgfonlayer}{nodelayer}
		\node [style=black] (1) at (0.25, 4) {};
		\node [style=none] (2) at (0.25, 5) {};
		\node [style=circle] (3) at (0.25, 3) {};
		\node [style=none] (4) at (0.25, 2) {};
	\end{pgfonlayer}
	\begin{pgfonlayer}{edgelayer}
		\draw (2.center) to (1);
		\draw (4.center) to (3);
	\end{pgfonlayer}
\end{tikzpicture} \stackrel{(1)}{=}  \begin{tikzpicture}
	\begin{pgfonlayer}{nodelayer}
		\node [style=circle] (0) at (-0.75, 5) {};
		\node [style=none] (1) at (-0.75, 3) {};
		\node [style=circle] (2) at (0, 5) {};
		\node [style=black] (3) at (0.5, 3) {};
		\node [style=circle] (4) at (0.5, 4) {};
		\node [style=none] (5) at (1, 5) {};
	\end{pgfonlayer}
	\begin{pgfonlayer}{edgelayer}
		\draw (0) to (1.center);
		\draw [bend left] (5.center) to (4);
		\draw [bend right] (2) to (4);
		\draw (4) to (3);
	\end{pgfonlayer}
\end{tikzpicture} \stackrel{(2)}{=}\begin{tikzpicture}
	\begin{pgfonlayer}{nodelayer}
		\node [style=circle] (0) at (-0.5, 5.25) {};
		\node [style=none] (1) at (-1, 3) {};
		\node [style=black] (3) at (0.25, 3) {};
		\node [style=circle] (4) at (0.25, 4) {};
		\node [style=none] (5) at (0.75, 5.25) {};
		\node [style=black] (6) at (-0.5, 4.5) {};
	\end{pgfonlayer}
	\begin{pgfonlayer}{edgelayer}
		\draw [in=45, out=-90] (5.center) to (4);
		\draw (4) to (3);
		\draw (6) to (4);
		\draw [in=90, out=-135, looseness=1.25] (6) to (1.center);
		\draw (0) to (6);
	\end{pgfonlayer}
\end{tikzpicture} \stackrel{(3)}{=} \begin{tikzpicture}
	\begin{pgfonlayer}{nodelayer}
		\node [style=circle] (0) at (-0.25, 5.25) {};
		\node [style=none] (1) at (-0.25, 3) {};
		\node [style=black] (3) at (0.75, 3) {};
		\node [style=circle] (4) at (0.25, 4.5) {};
		\node [style=none] (5) at (0.75, 5.25) {};
		\node [style=black] (6) at (0.25, 3.75) {};
	\end{pgfonlayer}
	\begin{pgfonlayer}{edgelayer}
		\draw [in=30, out=-90] (5.center) to (4);
		\draw [bend right, looseness=1.25] (6) to (1.center);
		\draw [bend left, looseness=1.25] (6) to (3);
		\draw [in=150, out=-90] (0) to (4);
		\draw (4) to (6);
	\end{pgfonlayer}
\end{tikzpicture} = \begin{tikzpicture}
	\begin{pgfonlayer}{nodelayer}
		\node [style=none] (1) at (0.75, 3) {};
		\node [style=none] (5) at (0.75, 5.25) {};
	\end{pgfonlayer}
	\begin{pgfonlayer}{edgelayer}
		\draw (5.center) to (1.center);
	\end{pgfonlayer}
\end{tikzpicture} \]
\end{proof}
Step $(1)$ of the proof uses the unit law for monoids, step $(2)$ uses equation $(b)$ of bialgebras, 
and step $(3)$ uses the Frobenius law. Thus, for a Frobenius algebra which is also a bialgebra, 
the unit map is inverse of the counit map. 

In the previous section, we saw that quantum observables are precisely 
$\dagger$-SCFAs in $\FHilb$. We now present the conditions for any two  $\dagger$-SCFA 
to be complementary:

\begin{definition} \cite{HeV19}
In a $\dagger$-SMC, any two $\dagger$-SCFAs, say $\Frob{A}$ and $\bFrob{A}$ are 
{\bf complementary} if $\bialg{A}$ ( equivalently $\bialgb{A}$) is a Hopf algebra 
with antipode:
\[ \begin{tikzpicture}
	\begin{pgfonlayer}{nodelayer}
		\node [style=black] (0) at (-0.5, 3.5) {};
		\node [style=black] (1) at (-0.5, 4.5) {};
		\node [style=circle] (2) at (0.25, 2.25) {};
		\node [style=circle] (3) at (0.25, 1.25) {};
		\node [style=none] (4) at (-1, 1.25) {};
		\node [style=none] (5) at (0.75, 4.5) {};
	\end{pgfonlayer}
	\begin{pgfonlayer}{edgelayer}
		\draw (1) to (0);
		\draw (2) to (3);
		\draw [in=90, out=-135] (0) to (4.center);
		\draw (0) to (2);
		\draw [in=-90, out=45] (2) to (5.center);
	\end{pgfonlayer}
\end{tikzpicture} = \begin{tikzpicture}
	\begin{pgfonlayer}{nodelayer}
		\node [style=black] (0) at (0.25, 3.5) {};
		\node [style=black] (1) at (0.25, 4.5) {};
		\node [style=circle] (2) at (-0.5, 2.25) {};
		\node [style=circle] (3) at (-0.5, 1.25) {};
		\node [style=none] (4) at (0.75, 1.25) {};
		\node [style=none] (5) at (-1, 4.5) {};
	\end{pgfonlayer}
	\begin{pgfonlayer}{edgelayer}
		\draw (1) to (0);
		\draw (2) to (3);
		\draw [in=90, out=-45] (0) to (4.center);
		\draw (0) to (2);
		\draw [in=-90, out=135] (2) to (5.center);
	\end{pgfonlayer}
\end{tikzpicture} \]
\end{definition}
Note that in the CQM literature, the above definition is referred to as {\em strong complementarity}. 

Let us look at an example of complementary observables in $\FHilb$.
A simple yet significant example is given by the spin of an electron 
along the $X$, $Y$ and $Z$ axes. The spin of an electron is either 
{\em up} or {\em down} or a superposition these two states. This system 
is referred to as a qubit in quantum computation, and its state space is $\C^2$. 

Each of the $X$, $Y$, and $Z$ observables are complementary to one another. 
These observables are given by {\em Pauli matrices}  $X$, $Y$ and $Z$ defined as follows:
 \[ X =  \begin{pmatrix}  0 & 1 \\ 1 & 0  \end{pmatrix} 
 ~~~~~~~~ Y = \begin{pmatrix}  0 & -i \\ i & 0  \end{pmatrix} 
~~~~~~~~ Z = \begin{pmatrix}  1 & 0 \\ 0 & -1  \end{pmatrix} \]
The eigenbasis of $X$, $Y$, and $Z$ are $\left \{(|0 \rangle + | 1 \rangle),  (|0 \rangle - | 1 \rangle) \right \}$, 
$Y$ is $\left \{ (|0 \rangle + i | 1 \rangle), (|0 \rangle - i | 1 \rangle) \right \}$, 
and $Z$ is $\{ |0 \rangle, | 1 \rangle \}$ respectively, where:
\[ | 0 \rangle = \begin{pmatrix}  1  \\  0   \end{pmatrix} ~~~~~~~  | 1 \rangle = \begin{pmatrix}  0  \\  1   \end{pmatrix} \]
These eigenvectors provide an orthonormal basis for $\C^2$, hence they are $\dagger$-SCFAs in $\FHilb$, see \cite[Example 2.11]{CoD11}. Note that, 
when used in calculations, the $X$ and $Y$-eigenvectors are normalized so that for any of these vectors $v$, $\langle v | v \rangle  = v v^\dag = 1$. 
In CQM, normalization is handled by introducing scalars $(a: I \to I)$ in the bialgebra equations \cite{CoD11, Bac15}. 

For a pair of complementary observables, the antipode is the identity map if and only if the self-dual 
cups and caps of the Frobenius algebras coincide. This is true in the case of $Z$ and $X$ observables.
(and not for $Z$-$Y$ and $X$-$Y$ pairs). Coecke and Duncan developed a diagrammatic 
calculus for $Z$ and $X$ observables, called the ZX-calculus \cite{CoD11}. 
The bialgebraic interaction between the $Z$ and the $X$ observables is used in the calculus 
 to construct logic gates in quantum computation. The ZX calculus is considered one of the most significant 
outcomes of the CQM program and is used extensively for quantum circuits optimization.

We end this section by providing an alternate and a useful characterization of complementary systems: 

\begin{theorem}\cite[Theorem 6.4]{DuK16} 
In a $\dagger$-SMC, two $\dagger$-SCFA $\Frob{A}$, and $\bFrob{A}$ are complementary if and only if 
the following equations hold: 
\[ \begin{tikzpicture}
	\begin{pgfonlayer}{nodelayer}
		\node [style=black] (0) at (0.25, 3.5) {};
		\node [style=black] (1) at (0.25, 4.5) {};
		\node [style=circle] (2) at (-0.25, 2.25) {};
		\node [style=none] (4) at (0.75, 2.25) {};
	\end{pgfonlayer}
	\begin{pgfonlayer}{edgelayer}
		\draw (1) to (0);
		\draw [in=90, out=-30] (0) to (4.center);
		\draw [in=90, out=-150] (0) to (2);
	\end{pgfonlayer}
\end{tikzpicture} = \begin{tikzpicture}
	\begin{pgfonlayer}{nodelayer}
		\node [style=circle] (1) at (0.25, 4.5) {};
		\node [style=none] (2) at (0.25, 2.25) {};
	\end{pgfonlayer}
	\begin{pgfonlayer}{edgelayer}
		\draw (2.center) to (1);
	\end{pgfonlayer}
\end{tikzpicture}
~~~~~~~~ \begin{tikzpicture}
	\begin{pgfonlayer}{nodelayer}
		\node [style=circle] (0) at (0.25, 3.25) {};
		\node [style=circle] (1) at (0.25, 2.25) {};
		\node [style=black] (2) at (-0.25, 4.5) {};
		\node [style=none] (4) at (0.75, 4.5) {};
	\end{pgfonlayer}
	\begin{pgfonlayer}{edgelayer}
		\draw (1) to (0);
		\draw [in=-90, out=30] (0) to (4.center);
		\draw [in=-90, out=150] (0) to (2);
	\end{pgfonlayer}
\end{tikzpicture} = \begin{tikzpicture}
	\begin{pgfonlayer}{nodelayer}
		\node [style=black] (1) at (0.25, 2.25) {};
		\node [style=none] (2) at (0.25, 4.5) {};
	\end{pgfonlayer}
	\begin{pgfonlayer}{edgelayer}
		\draw (2.center) to (1);
	\end{pgfonlayer}
\end{tikzpicture}\] 
\end{theorem}

In the previous section, for every Frobenius algebra, the multiplication is dual the comultiplication, and the 
unit is to dual to the counit via the self-dual cup and cap of the Frobenius algebra. 
The above theorem implies two $\dagger$-SCFAs are complementary if and only if 
the unit and the counit of each of these algebras are 
dual to one another using the self-dual cup and the cap of its complementary algebra. 

%% file: chapter-CP.tex

\chapter{Completely positive maps}
\label{Chap: positivity}

Categorical quantum mechanics (CQM) has mostly focused on quantum information theory 
and quantum computation which are finite dimensional branches of quantum mechanics. 
In order to widen the scope of CQM, there have been various proposals \cite{CoH16, GG17, AbH12} 
to extend the structures in CQM to infinite dimensional systems. In Part \ref{Part: dll} of this thesis, we generalized the framework 
of dagger compact closed categories to what we have called mixed unitary categories (MUCs).  In this chapter 
we define completely positive maps for a MUC. 

\section{String calculus for MUCs}
To facilitate reasoning within MUCs, it is useful to employ a circuit calculus built on the circuit calculus for LDCs introduced 
in \cite{BCST96}. The extended circuit calculus for mixed unitary categories includes dagger boxes,  components for unitary 
structure maps, and inverse mixor morphisms.

\subsection{Unitary structure map}

A unitary object is an object equipped with an isomorphism $A \xrightarrow{\varphi_A} A^\dagger$, called the unitary 
structure map, which is drawn as a downward pointing triangle:
$$

$$

Diagrammatic representation of axioms {\bf [U.5]}-(a), {\bf [U.4]}-(a) in Definition \ref{defn: unitary structure} 
for unitary structure, and {\bf [Udual]}-(a) in Definition \ref{defn: unitary dual} for unitary duals are as follows. 

\begin{center}
%
\end{center}

\subsection{Inverse of the mixor}

The mixor (obtained from the mix map $\bot \to^{\m} \top$ as below) $U_1 \ox U_2 \to^{\mx} U_1 \oa U_2$ is an 
isomorphism whenever $U_1$ and $U_2$ are unitary objects and its inverse are represented as follows:

$$
\mx_{U_1,U_2}:=
%
$$

Observe that $\mx^{-1}$ maps asscocitated to different unitary objects slide past each other, and indeed by naturality, 
$\mx^{-1}$ slides over components in circuits:
$$
%
$$

\section{Quantum channels for MUCs}

\subsection{Kraus maps}
A Kraus map $(f,U): A \to B$ in a mixed unitary category, $M: \U \to \C$, is a map $f: A \to M(U) \oa B \in \C$ for some 
$U \in \U$. $U$ is called the ancillary system of $f$.  We glue the Kraus map to its dagger along its ancillary system giving 
rise to a combinator which acts on ``test maps''. Two Kraus maps are equivalent when their effects on test maps are indistinguishable.

\begin{definition} Given a MUC,  $\U \to^{M} \C$, two Kraus maps $(f,U) , (g,V): A \to B$ are {\bf equivalent}, $(f,U) \sim (g,V)$, 
    if for all  unitary objects $X$ and all maps $h: B \ox C \rightarrow V $ (called {\bf test maps}), the following equation holds:
\[ \begin{tikzpicture} 
	\begin{pgfonlayer}{nodelayer}
		\node [style=none] (0) at (-1.75, -3.5) {};
		\node [style=none] (1) at (-1.75, -3.5) {};
		\node [style=none] (2) at (-1.25, -2) {};
		\node [style=none] (3) at (-0.75, -0.25) {};
		\node [style=none] (4) at (-0.75, -3) {};
		\node [style=none] (5) at (-0.5, -3) {};
		\node [style=none] (6) at (-0.25, -2) {};
		\node [style=none] (7) at (-0.25, -5.5) {};
		\node [style=none] (8) at (-0.5, -2) {};
		\node [style=none] (9) at (-0.75, -2) {};
		\node [style=circle,scale=1.5] (10) at (-0.75, -2.5) {};
		\node [style=none] (100) at (-0.75, -2.5) {$h$};
		\node [style=none] (11) at (-1, -3) {};
		\node [style=none] (12) at (-1.25, -3) {};
		\node [style=none] (13) at (-1, -2) {};
		\node [style=none] (14) at (-0.25, -3) {};
		\node [style=none] (15) at (-1, -0) {};
		\node [style=none] (16) at (-0.75, -0.25) {};
		\node [style=none] (17) at (-0.5, -0) {};
		\node [style=none] (18) at (-0.5, 3.25) {};
		\node [style=circle, scale=1.5] (19) at (-0.75, 1.25) {};
		\node [style=none] (190) at (-0.75, 1.25) {$h$};
		\node [style=none] (20) at (-0.75, -0) {};
		\node [style=none] (21) at (-2, 3.25) {};
		\node [style=none] (22) at (-2, -5) {};
		\node [style=none] (23) at (-2, -5.5) {};
		\node [style=circle, scale=1.5] (24) at (-2, -4) {};
		\node [style=none] (240) at (-2, -4) {$f$};
		\node [style=none] (25) at (-2.5, -3.5) {};
		\node [style=none] (26) at (-2.5, -4.5) {};
		\node [style=none] (27) at (-1.5, -4.5) {};
		\node [style=none] (28) at (-1.5, -3.5) {};
		\node [style=none] (29) at (-1.75, -4.5) {};
		\node [style=none] (30) at (-2.25, -4.5) {};
		\node [style=none] (31) at (-2, -3.5) {};
		\node [style=circle, scale=1.5] (32) at (-2, 2.25) {};
		\node [style=none] (320) at (-2, 2.25) {$f$};
		\node [style=none] (33) at (-2.5, -0) {};
		\node [style=none] (34) at (-2.25, -3.5) {};
		\node [style=none] (35) at (-1.75, -3.5) {};
		\node [style=none] (36) at (-2.75, -0) {};
		\node [style=none] (37) at (-2.25, -0) {};
		\node [style=none] (38) at (-2.5, -0.25) {};
		\node [style=none] (39) at (-2.5, -0) {};
		\node [style=none] (40) at (-2, -4.5) {};
		\node [style=none] (41) at (-2.75, 0.75) {};
		\node [style=none] (42) at (-0.5, 0.75) {};
		\node [style=none] (43) at (-0.5, 0.5) {};
		\node [style=none] (44) at (-3, 0.25) {};
		\node [style=none] (45) at (-2, 0.25) {};
		\node [style=none] (46) at (-3, -0.75) {};
		\node [style=none] (47) at (-2, -0.75) {};
		\node [style=none] (48) at (-2.25, -0.5) {$M$};
		\node [style=none] (49) at (-1.25, 0.25) {};
		\node [style=none] (50) at (-0.25, 0.25) {};
		\node [style=none] (51) at (-1.25, -0.75) {};
		\node [style=none] (52) at (-0.25, -0.75) {};
		\node [style=none] (53) at (-0.5, -0.5) {$M$};
		\node [style=circle, scale=1.5] (54) at (-2.5, -1.5) {};
		\node [style=circle, scale=1.5] (55) at (-0.75, -1.25) {};
		\node [style=none] (540) at (-2.5, -1.5) {$\rho$};
		\node [style=none] (550) at (-0.75, -1.25) {$\rho$};
	\end{pgfonlayer}
	\begin{pgfonlayer}{edgelayer}
		\draw [style=none] (0.center) to (1.center);
		\draw [style=none] (13.center) to (10);
		\draw [style=none] (10) to (8.center);
		\draw [style=none] (10) to (4.center);
		\draw [style=none] (2.center) to (6.center);
		\draw [style=none] (6.center) to (14.center);
		\draw [style=none] (14.center) to (12.center);
		\draw [style=none] (12.center) to (2.center);
		\draw [style=none, in=90, out=-90, looseness=1.00] (5.center) to (7.center);
		\draw [style=none] (15.center) to (17.center);
		\draw [style=none] (17.center) to (16.center);
		\draw [style=none] (16.center) to (15.center);
		\draw [style=none, in=-90, out=60, looseness=0.75] (19) to (18.center);
		\draw [style=none] (19) to (20.center);
		\draw [style=none] (22.center) to (23.center);
		\draw (25.center) to (26.center);
		\draw (26.center) to (27.center);
		\draw (28.center) to (27.center);
		\draw (25.center) to (28.center);
		\draw (24) to (31.center);
		\draw (24) to (30.center);
		\draw (24) to (29.center);
		\draw [in=90, out=-135, looseness=0.75] (32) to (33.center);
		\draw (36.center) to (37.center);
		\draw (37.center) to (38.center);
		\draw (38.center) to (36.center);
		\draw (40.center) to (22.center);
		\draw [in=180, out=-45, looseness=1.25] (32) to (19);
		\draw (21.center) to (32);
		\draw [in=-150, out=30, looseness=1.00] (41.center) to (42.center);
		\draw [in=90, out=-90, looseness=1.00] (11.center) to (0.center);
		\draw (44.center) to (46.center);
		\draw (46.center) to (47.center);
		\draw (47.center) to (45.center);
		\draw (45.center) to (44.center);
		\draw (49.center) to (50.center);
		\draw (50.center) to (52.center);
		\draw (52.center) to (51.center);
		\draw (51.center) to (49.center);
		\draw (38.center) to (54);
		\draw (3.center) to (55);
		\draw (55) to (9.center);
		\draw [in=98, out=-90, looseness=1.00] (54) to (34.center);
	\end{pgfonlayer}
\end{tikzpicture} =   \begin{tikzpicture} 
	\begin{pgfonlayer}{nodelayer}
		\node [style=none] (0) at (-1.75, -3.5) {};
		\node [style=none] (1) at (-1.75, -3.5) {};
		\node [style=none] (2) at (-1.25, -2) {};
		\node [style=none] (3) at (-0.75, -0.25) {};
		\node [style=none] (4) at (-0.75, -3) {};
		\node [style=none] (5) at (-0.5, -3) {};
		\node [style=none] (6) at (-0.25, -2) {};
		\node [style=none] (7) at (-0.25, -5.5) {};
		\node [style=none] (8) at (-0.5, -2) {};
		\node [style=none] (9) at (-0.75, -2) {};
		\node [style=circle,scale=1.5] (10) at (-0.75, -2.5) {};
		\node [style=none] (100) at (-0.75, -2.5) {$h$};
		\node [style=none] (11) at (-1, -3) {};
		\node [style=none] (12) at (-1.25, -3) {};
		\node [style=none] (13) at (-1, -2) {};
		\node [style=none] (14) at (-0.25, -3) {};
		\node [style=none] (15) at (-1, -0) {};
		\node [style=none] (16) at (-0.75, -0.25) {};
		\node [style=none] (17) at (-0.5, -0) {};
		\node [style=none] (18) at (-0.5, 3.25) {};
		\node [style=circle, scale=1.5] (19) at (-0.75, 1.25) {};
		\node [style=none] (190) at (-0.75, 1.25) {$h$};
		\node [style=none] (20) at (-0.75, -0) {};
		\node [style=none] (21) at (-2, 3.25) {};
		\node [style=none] (22) at (-2, -5) {};
		\node [style=none] (23) at (-2, -5.5) {};
		\node [style=circle, scale=1.5] (24) at (-2, -4) {};
		\node [style=none] (240) at (-2, -4) {$g$};
		\node [style=none] (25) at (-2.5, -3.5) {};
		\node [style=none] (26) at (-2.5, -4.5) {};
		\node [style=none] (27) at (-1.5, -4.5) {};
		\node [style=none] (28) at (-1.5, -3.5) {};
		\node [style=none] (29) at (-1.75, -4.5) {};
		\node [style=none] (30) at (-2.25, -4.5) {};
		\node [style=none] (31) at (-2, -3.5) {};
		\node [style=circle, scale=1.5] (32) at (-2, 2.25) {};
		\node [style=none] (320) at (-2, 2.25) {$g$};
		\node [style=none] (33) at (-2.5, -0) {};
		\node [style=none] (34) at (-2.25, -3.5) {};
		\node [style=none] (35) at (-1.75, -3.5) {};
		\node [style=none] (36) at (-2.75, -0) {};
		\node [style=none] (37) at (-2.25, -0) {};
		\node [style=none] (38) at (-2.5, -0.25) {};
		\node [style=none] (39) at (-2.5, -0) {};
		\node [style=none] (40) at (-2, -4.5) {};
		\node [style=none] (41) at (-2.75, 0.75) {};
		\node [style=none] (42) at (-0.5, 0.75) {};
		\node [style=none] (43) at (-0.5, 0.5) {};
		\node [style=none] (44) at (-3, 0.25) {};
		\node [style=none] (45) at (-2, 0.25) {};
		\node [style=none] (46) at (-3, -0.75) {};
		\node [style=none] (47) at (-2, -0.75) {};
		\node [style=none] (48) at (-2.25, -0.5) {$M$};
		\node [style=none] (49) at (-1.25, 0.25) {};
		\node [style=none] (50) at (-0.25, 0.25) {};
		\node [style=none] (51) at (-1.25, -0.75) {};
		\node [style=none] (52) at (-0.25, -0.75) {};
		\node [style=none] (53) at (-0.5, -0.5) {$M$};
		\node [style=circle, scale=1.5] (54) at (-2.5, -1.5) {};
		\node [style=circle, scale=1.5] (55) at (-0.75, -1.25) {};
		\node [style=none] (540) at (-2.5, -1.5) {$\rho$};
		\node [style=none] (550) at (-0.75, -1.25) {$\rho$};
	\end{pgfonlayer}
	\begin{pgfonlayer}{edgelayer}
		\draw [style=none] (0.center) to (1.center);
		\draw [style=none] (13.center) to (10);
		\draw [style=none] (10) to (8.center);
		\draw [style=none] (10) to (4.center);
		\draw [style=none] (2.center) to (6.center);
		\draw [style=none] (6.center) to (14.center);
		\draw [style=none] (14.center) to (12.center);
		\draw [style=none] (12.center) to (2.center);
		\draw [style=none, in=90, out=-90, looseness=1.00] (5.center) to (7.center);
		\draw [style=none] (15.center) to (17.center);
		\draw [style=none] (17.center) to (16.center);
		\draw [style=none] (16.center) to (15.center);
		\draw [style=none, in=-90, out=60, looseness=0.75] (19) to (18.center);
		\draw [style=none] (19) to (20.center);
		\draw [style=none] (22.center) to (23.center);
		\draw (25.center) to (26.center);
		\draw (26.center) to (27.center);
		\draw (28.center) to (27.center);
		\draw (25.center) to (28.center);
		\draw (24) to (31.center);
		\draw (24) to (30.center);
		\draw (24) to (29.center);
		\draw [in=90, out=-135, looseness=0.75] (32) to (33.center);
		\draw (36.center) to (37.center);
		\draw (37.center) to (38.center);
		\draw (38.center) to (36.center);
		\draw (40.center) to (22.center);
		\draw [in=180, out=-45, looseness=1.25] (32) to (19);
		\draw (21.center) to (32);
		\draw [in=-150, out=30, looseness=1.00] (41.center) to (42.center);
		\draw [in=90, out=-90, looseness=1.00] (11.center) to (0.center);
		\draw (44.center) to (46.center);
		\draw (46.center) to (47.center);
		\draw (47.center) to (45.center);
		\draw (45.center) to (44.center);
		\draw (49.center) to (50.center);
		\draw (50.center) to (52.center);
		\draw (52.center) to (51.center);
		\draw (51.center) to (49.center);
		\draw (38.center) to (54);
		\draw (3.center) to (55);
		\draw (55) to (9.center);
		\draw [in=98, out=-90, looseness=1.00] (54) to (34.center);
	\end{pgfonlayer}
\end{tikzpicture}   \]
\end{definition}

Note that the $\mx^{-1}$ map can be slid up and down along the wires of the unitary objects $M(U)$ and $M(V)$ by 
naturality of the $\mx$ map. The diagram includes covariant functor boxes for $M$ and contravariant functor boxes for the 
dagger. The diagram on the left is given equationally as follows:
\[
A \ox C \to^{f \ox 1} (M(U) \oa B) \ox C \to^{\partial} M(U) \oa (B \ox C) \to^{1 \oa h} M(U) \oa M(V) \to^{\mx^{-1}} M(U) 
\ox M(V)\]\[ \to^{M(\varphi_{U}) \ox M(\varphi_V)} M(U^\dagger) \ox M(V^\dagger) \to^{\rho \ox \rho} M(U)^\dagger 
\ox M(V)^\dagger \to^{1 \ox (h^\dagger \lambda_\oa^{-1})}  M(U)^\dagger \ox (B^\dagger \oa C^\dagger) \] \[\to^{\delta} 
(M(U)^\dagger \ox B^\dagger) \oa C^\dagger \to^{\lambda_\ox \oa 1} (M(U) \oa B)^\dagger \oa C^\dagger \to^{f^\dagger \oa 1} 
A^\dagger \oa C^\dagger \to^{\lambda_\oa} (A \ox C)^\dagger
\]

The natural isomorphism $\rho: M(U^\dagger) \to M(U)^\dagger$ is the preservator of the $\dagger$-isomix functor, $M$, 
which ensures coherence with the $\dagger$ from $\U$ to $\C$. In the rest of this chapter, 
the covariant functor boxes represent the $\dagger$-isomix functor $M$ of the MUC unless specified 
otherwise.

By forgetting the test maps and gluing Kraus map with its dagger, one gets a notationally convenient combinator which can 
be diagrammatically represented by:
 \[  
 \]

 An equivalence class of Kraus morphisms is a {\bf quantum channel}. If $(f,U): A \to B$ and $(g, V): A \to B$ are Kraus 
 maps for which $U$ and $V$ are unitarily isomorphic, they are necessarily equivalent with respect to this  relation:

\begin{lemma}
\label{Lemma: unitary equivalence} 
Let $(f,U), (g,V) : A \to B$ be Kraus morphisms. If $~U \to^{\alpha} V$ is a unitary isomorphism and $f( \alpha \oa 1) = 
g$, then $(f, U) \sim (g, V)$.
\end{lemma}
\begin{proof} Let $h: B \ox C \to M(X)$ be any test map. Then, 
\[ %
 \]
\end{proof}

In a $*$-MUC, the equivalence relation on Kraus maps is given equationally as follows:

\begin{lemma}
Suppose $M: \U \to \C$,  is a $*$-MUC, that is every object in $\C$ has  linear adjoint, then two Kraus maps $(f, U)$ and 
$(g, V)$ are equivalent if and only if
\[
%
 \]
\end{lemma}
\begin{proof}
It is obvious that the given equation holds when the Kraus maps $(f, U)$ and $(g, V)$ are equivalent. For the 
converse, assume that the given equation holds. Then, the Kraus maps $(f, U)$ and $(g, V)$ are equivalent because, for any 
test map $h$, we have that: 
\[ %
 \]
\end{proof}

Let us now examine Kraus maps in our running examples: 
\begin{itemize}
\item 
In the MUC, $\R \subset \C$, let $c, c'$ be any two complex numbers. Kraus maps in $\R \subset \C$ are $ (=, r): c \to c' $ 
such that $c = rc'.$ If $c' \neq 0$, then there is at most one Kraus map $(=,r): c \to c'$. If $c'=0$, then $c=0$ and for all 
$r' \in \R$, $(=,r) \sim (=,r'): c \to c'$. Thus,  in the complex plane, there are only Kraus maps between those complex numbers 
that can be connected by a line that extends through the origin making it 
the projective line $P^1(\C)$, where $\C$ refers to the complex field. 
\[
\begin{tikzpicture} [scale=2]
	\begin{pgfonlayer}{nodelayer}
		\node [style=none] (0) at (0, 2.5) {};
		\node [style=none] (1) at (0, -2.25) {};
		\node [style=none] (2) at (-2.5, -0) {};
		\node [style=none] (3) at (2.5, -0) {};
		\node [style=none] (4) at (1.5, 1.5) {};
		\node [style=none] (5) at (1.5, -1.5) {};
		\node [style=none] (6) at (-1.5, -1.5) {};
		\node [style=none] (7) at (-1.5, 1.5) {};
		\node [style=none] (8) at (2.5, 0.25) {$\R$};
		\node [style=none] (9) at (0.25, 2.5) {$\iota$};
		\node [style=none] (10) at (0.75, 2) {};
		\node [style=none] (11) at (-0.75, -2) {};
	\end{pgfonlayer}
	\begin{pgfonlayer}{edgelayer}
		\draw [<-] (0.center) to (1.center);
		\draw [->](2.center) to (3.center);
		\draw [blue] (7.center) to (5.center);
		\draw [blue] (6.center) to (4.center);
		\draw [blue] (10.center) to (11.center);
	\end{pgfonlayer}
\end{tikzpicture} \]

\item In $\Mat_\C \subset \FMat_\C$, every Kraus map $(M, \C^n): (X, \mathcal{A}, \mathcal{A}^\perp) \to (Y, \mathcal{B}, 
\mathcal{B}^\perp)$ is given by the sum of {\bf pure completely positive maps} i.e., Kraus maps with $\C$ as ancillary object:
\[ \begin{tikzpicture}
	\begin{pgfonlayer}{nodelayer}
		\node [style=none] (0) at (0, 3.75) {};
		\node [style=circle] (1) at (0, 2.5) {$M$};
		\node [style=none] (2) at (-0.8, 3.5) {$(X, \mathcal{A})$};
		\node [style=circle] (4) at (0, -2.5) {$M^\dagger$};
		\node [style=none] (5) at (0, -3.75) {};
		\node [style=none] (6) at (-0.95, -3.5) {$(X, \mathcal{A}^\perp)$};
		\node [style=circle] (7) at (1, 1.25) {$N$};
		\node [style=none] (8) at (1.5, 3.75) {};
		\node [style=circle] (9) at (1, -1) {$N^\dagger$};
		\node [style=none] (10) at (1.5, -3.75) {};
		\node [style=none] (11) at (2.2, 3.5) {$(Y, \mathcal{B})$};
		\node [style=none] (12) at (1.35, -0) {$\C^m$};
		\node [style=none] (13) at (2.25, -3.5) {$(Y, \mathcal{B}^\perp)$};
		\node [style=none] (14) at (0.75, 2.25) {};
		\node [style=none] (15) at (0.4, -1.75) {};
	\end{pgfonlayer}
	\begin{pgfonlayer}{edgelayer}
		\draw (0.center) to (1);
		\draw (5.center) to (4);
		\draw [bend right=15, looseness=1.00] (7) to (8.center);
		\draw [bend left=15, looseness=1.00] (9) to (10.center);
		\draw (7) to (9);
		\draw [bend left=15, looseness=1.00] (9) to (4);
		\draw [bend left=15, looseness=1.00] (1) to (7);
	\end{pgfonlayer}
\end{tikzpicture} 
 \]

 {\bf Choi's theorem} states that every completely positive map can be written as a sum of pure completely positive maps. Analogously, every Kraus map in the category $\FMat_\C$ can be written as a sum of pure maps as above.  Given a Kraus map $(M, \C^m)$, here is the argument:
\begin{align*}
&(M \ox 1) (1 \ox NN^\dagger)(M^\dagger \oa 1) = (M \ox 1) \left( \left( \sum_i \invamalg_i \amalg_i  \right) \ox NN^\dagger \right) (M^\dagger \oa 1) \\
&=  \sum_i \left(  (M(\invamalg_i \ox 1) \ox 1) (1 \ox hh^\dagger) (((\amalg_i \oa 1)M^\dagger) \oa 1)  \right) = \sum_i (M_i \ox 1) (1 \ox NN^\dagger) (M_i^\dagger \oa 1)
\end{align*}
\end{itemize}

\subsection{$\CP^\infty$-construction for MUCs}
The CPM-construction \cite{Sel07} on $\dagger$-compact closed categories applied to the concrete category 
of finite-dimensional Hilbert Spaces and linear maps produces a category of mixed states 
and completely positive maps. Coecke and Heunen \cite{CoH16} 
generalized the CPM-construction to $\dag$-symmetric monoidal categories, and thus, to infinite dimensions. 
They call the generalized construction the $\CP^\infty$-construction.  
In this section, we generalize the $\CP^\infty$ construction to MUCs: our construction coincides with the original 
$\CP^\infty$-construction when the MUC is a $\dagger$-monoidal category.

The  $\CP^\infty$-construction on MUCs is defined as follows:

\begin{definition}
Given a MUC, $M: \U \to \C$, define $\CP^\infty(M: \U \to \C)$ to have:
\begin{description} 	
\item[Objects:] Same as $\C$
\item[Maps:] $[(f,U)]:A \to B := \{ f: A \to M(U) \oa B \in  \C) \} \ / \sim$
\item[Composition:] Composition of $[(f,U)]:A \to B$ and $[(g,V)]:B \to C$ is defined as follows: 
$$[(f,U)][(g,V)] := (A \xrightarrow{f} M(U) \oa B \xrightarrow{1 \oa g} M(U) \oa (M(V) \oa C) \xrightarrow{a_\oa} (M(U \oa V)) \oa C \in \C) \ / \sim$$
\item[Identity:]  $1_A$ is defined as $[A \xrightarrow{(u_\oa^L)^{-1}} \bot \oa A \to^{(n_\bot^M)^{-1} \oa 1} M(\bot) \oa A]  \in \X$
\end{description}
\end{definition}
In $\CP^\infty(M: \U \to \C)$, composition of $f: A \to M(U) \oa B$, and $g: B \to M(V) \oa C$ is drawn as follows:
\[  \begin{tikzpicture} 
	\begin{pgfonlayer}{nodelayer}
		\node [style=circle, scale=1.5] (0) at (0.5, 2) {};
		\node [style=circle, scale=1.5] (1) at (1, 1) {};
		\node [style=none] (100) at (0.5, 2) {$f$};
		\node [style=none] (101) at (1, 1) {$g$};
		\node [style=oa] (2) at (-0.25, -0.25) {};
		\node [style=none] (3) at (-0.75, 0.25) {};
		\node [style=none] (4) at (0.25, 0.25) {};
		\node [style=none] (5) at (0.25, -1) {};
		\node [style=none] (6) at (-0.75, -1) {};
		\node [style=none] (7) at (0, -0.75) {$M$};
		\node [style=none] (8) at (0.5, 2.75) {};
		\node [style=none] (9) at (1.5, -1.5) {};
		\node [style=none] (10) at (-0.25, -1.5) {};
	\end{pgfonlayer}
	\begin{pgfonlayer}{edgelayer}
		\draw (3.center) to (6.center);
		\draw (6.center) to (5.center);
		\draw (5.center) to (4.center);
		\draw (4.center) to (3.center);
		\draw [in=90, out=-45, looseness=1.00] (0) to (1);
		\draw [in=45, out=-135, looseness=1.00] (1) to (2);
		\draw [in=120, out=-135, looseness=1.00] (0) to (2);
		\draw [in=90, out=-45, looseness=1.00] (1) to (9.center);
		\draw (2) to (10.center);
		\draw (0) to (8.center);
	\end{pgfonlayer}
\end{tikzpicture} \]

Observe that our $\CP^\infty$-construction on $M: \U \to \C$  coincides with 
the original $\CP^\infty$-construction \cite{CoH16} when $M: \U \to \C$ is dagger monoidal category 
$\U = \C$ and $M = id$. 

Before we prove that $\CP^{\infty}(M: \U \to \C)$ is a category, we observe the following result about  unitary objects:
\begin{lemma}
\label{Lemma: rho-tensor-par}
Suppose $C$ and $D$ are unitary objects. Then, the following diagrams
 commute:
\[ \xymatrixcolsep{1in} \xymatrix{
M(C) \ox M(D) \ar[d]_{m_\ox^M}  \ar[r]^{M(\varphi) \ox M(\varphi)} \ar@{}[dddr]|{\bf (a)} & M(C^\dagger) \ox M(D^\dagger) \ar[d]^{\rho \ox \rho} \\
M(C \ox D) \ar[d]_{M(\varphi)} & M(C)^\dagger \ox M(D)^\dagger \ar[d]^{\lambda_\ox} \\
M((C \ox D)^\dagger) \ar[d]_{\rho} & (M(C) \oa M(D))^\dagger \ar[d]^{\mx^\dagger}\\
(M(C) \ox M(D))^\dagger \ar[r]_{(m_\ox^M)^\dagger} & (M(C) \ox M(D))^\dagger
} \] \[ \xymatrixcolsep{1in} \xymatrix{
 M(C) \oa M(D) \ar[r]^{M(\varphi_C) \oa M(\varphi_D)} \ar[d]_{n_\oa^{-1}} \ar@{}[dddr]|{\bf (b)} & M(C^\dagger) \oa M(D^\dagger) \ar[d]^{\rho \oa \rho} \\
M(C \oa D) \ar[d]_{M(\varphi)} & M(C)^\dagger \oa M(D)^\dagger \ar[d]^{\lambda_\oa} \\
M((C \oa D)^\dagger) \ar[d]_{\rho} & (M(C) \ox M(D))^\dagger \ar[d]^{(\mx^{-1})^\dagger}  \\
(M(C \oa D))^\dagger & (M(C) \oa M(D))^\dagger \ar[l]^{n_\oa^\dagger}
} \]
\end{lemma}
\begin{proof}
    \begin{align*}

        \end{align*}
        
        The commuting diagram (b) is proved similarly.
\end{proof}

Diagrammatic representation of Lemma \ref{Lemma: rho-tensor-par} (b):
\[
%
 \]

\newpage

\begin{proposition}
$\CP^\infty(M: \U \to \C)$ is a category.
\end{proposition}
\begin{proof}~
\begin{itemize}
\item{Composition is well-defined:}
 That is to say, if $(f,U) \sim (f',U')$ then $(f,U)(g,V) \sim (f',U')(g,V)$.  First we observe that: 
It suffices to show that:

$%
$

$(g,V) \sim (g',V')  \Rightarrow (f,U)(g,V) \sim (f,U)(g',V')$ is proved similarly.

\item{Identity laws hold: $[(u_\oa^{L})^{-1} (n_\bot^M)^{-1}][f] = [f] = [f][(u_\oa^{L})^{-1} (n_\bot^M)^{-1}]$}
It suffices to prove the following: 
$$
%
$$

The proof uses the facts that $\rho$ is a monoidal transformation (diagram on the left), and that for any unitary object $U$, the diagram on the right holds:
\[ 
\xymatrixcolsep{4pc}
\xymatrix{
\top \ar[d]_{m_\top^M} \ar[rr]^{\lambda_\top} &  & \bot^\dagger \ar[d]_{(n_\bot^M)^\dagger}\ar@{=}[dr] \\
M(\top) \ar[r]_{M(\lambda_\top)} & M(\bot^\dagger) \ar[r]_{\rho} & M(\bot)^\dagger \ar[r]_{((n_\bot^M)^{-1})^\dagger} &  \bot^\dagger
} ~~~~~~~~~~ \xymatrix{
\top \ox U \ar[r]^{u_\ox^L} \ar[d]_{\m^{-1} \ox 1} & U \ar[d]^{(u_\oa^L)^{-1}} \\
\bot \ox U \ar[r]_{mx} & \bot \oa U
}
\]

\item{Composition is associative:}
Suppose $(f, U_1): A \to B,  (g, U_2): B \to C, (h, U_3): C \to D \in \CP^\infty(M: \U \to \C)$. Then,

$$
((f,U_1)(g, U_2)) (h, U_3) := \begin{tikzpicture} 
	\begin{pgfonlayer}{nodelayer}
		\node [style=circle] (0) at (-3.5, 1.75) {$f$};
		\node [style=circle] (1) at (-3, 0.75) {$g$};
		\node [style=circle] (2) at (-2.5, -0.25) {$h$};
		\node [style=none] (3) at (-2, -2) {};
		\node [style=oa] (4) at (-4, -0.5) {};
		\node [style=oa] (5) at (-3.5, -1.25) {};
		\node [style=none] (6) at (-3.5, -2) {};
		\node [style=none] (7) at (-3.5, 2.5) {};
		\node [style=none] (8) at (-4.5, -0) {};
		\node [style=none] (9) at (-3, -0) {};
		\node [style=none] (10) at (-3, -1.75) {};
		\node [style=none] (11) at (-4.5, -1.75) {};
		\node [style=none] (12) at (-4.25, -1.5) {$M$};
	\end{pgfonlayer}
	\begin{pgfonlayer}{edgelayer}
		\draw [style=none, in=90, out=-45, looseness=1.00] (0) to (1);
		\draw [style=none, in=90, out=-15, looseness=1.00] (1) to (2);
		\draw [style=none, in=90, out=-45, looseness=1.00] (2) to (3.center);
		\draw [style=none, in=45, out=-150, looseness=1.00] (1) to (4);
		\draw [style=none, in=-120, out=120, looseness=1.25] (4) to (0);
		\draw [style=none, in=150, out=-75, looseness=1.00] (4) to (5);
		\draw [style=none, in=75, out=-135, looseness=1.00] (2) to (5);
		\draw [style=none] (5) to (6.center);
		\draw [style=none] (7.center) to (0);
		\draw (9.center) to (10.center);
		\draw (10.center) to (11.center);
		\draw (11.center) to (8.center);
		\draw (8.center) to (9.center);
	\end{pgfonlayer}
\end{tikzpicture} \sim
\begin{tikzpicture}
	\begin{pgfonlayer}{nodelayer}
		\node [style=circle] (0) at (1.25, 0.25) {$f$};
		\node [style=oa] (1) at (0, -0.5) {};
		\node [style=circle] (2) at (0.25, 1.75) {$h$};
		\node [style=none] (3) at (0.25, 2.5) {};
		\node [style=none] (4) at (-0.5, -2) {};
		\node [style=none] (5) at (1.75, -2) {};
		\node [style=circle] (6) at (0.75, 1) {$g$};
		\node [style=oa] (7) at (-0.5, -1.25) {};
		\node [style=none] (8) at (-1, -0) {};
		\node [style=none] (9) at (0.5, -0) {};
		\node [style=none] (10) at (0.5, -1.75) {};
		\node [style=none] (11) at (-1, -1.75) {};
		\node [style=none] (12) at (0.25, -1.5) {$M$};
	\end{pgfonlayer}
	\begin{pgfonlayer}{edgelayer}
		\draw [style=none, in=90, out=-45, looseness=1.00] (2) to (6);
		\draw [style=none, in=90, out=-15, looseness=1.00] (6) to (0);
		\draw [style=none, in=90, out=-45, looseness=1.00] (0) to (5.center);
		\draw [style=none, in=90, out=-150, looseness=1.00] (6) to (1);
		\draw [style=none] (3.center) to (2);
		\draw [style=none, in=-150, out=15, looseness=1.25] (1) to (0);
		\draw [style=none, in=30, out=-124, looseness=1.25] (1) to (7);
		\draw [style=none, in=111, out=-135, looseness=1.00] (2) to (7);
		\draw [style=none] (7) to (4.center);
		\draw (8.center) to (9.center);
		\draw (9.center) to (10.center);
		\draw (10.center) to (11.center);
		\draw (11.center) to (8.center);
	\end{pgfonlayer}
\end{tikzpicture}
 =: (f, U_1) ((g, U_2)(h, U_3))
$$

Since $(U_1 \oa U_2) \oa U_3 \to^{a_\oa} U_1 \oa (U_2 \oa U_3)$ is a unitary isomorphism, by Lemma \ref{Lemma: unitary equivalence},  $(fg)h \sim f(gh) \in \C \Rightarrow ((f, U_1)(g, U_2))(h, U_3) = (f, U_1)((g, U_2)(h, U_3)) \in \CP^\infty(M:\U \to \C)$.
\end{itemize}
\end{proof}

There is a canonical functor $Q: \C \to  \CP^\infty(M: \U \to \C)$ of the original category into the category of channels:   

\begin{lemma}
\label{Lemma: embedding}
Let $M: \U \to \C $ be a mixed unitary category, then there is a canonical functor: 
\[ Q: \C \to  \CP^\infty(M: \U \to \C) ; 
\begin{matrix}
\xymatrix{
A \ar[d]_{f} \\
B
} \end{matrix}   \mapsto  \begin{matrix}
\xymatrix{
A \ar[d]^{f (u_\oa^L)^{-1}((n_\bot^M)^{-1} \oa 1)} \\
M(\bot) \oa B
}
\end{matrix}
\]
\end{lemma}
\begin{proof}
$Q$ preserves identity maps and composition because $ f(u_\oa^L)^{-1} \sim f \sim (u_\oa^L)^{-1} f $. 
\end{proof}
There is no reason why this functor should be faithful and, indeed, in many cases it will 
{\em not\/} be faithful \cite{CoH16}.

\begin{theorem}
	$\CP^\infty(M: \U \to \C)$ is an isomix category.
\end{theorem}
\begin{proof} We know that $\CP^\infty(M: \U \to \C)$ is well-defined category. Indeed, it has two tensors: $\widehat{\ox}$ and $\widehat{\oa}$ given by the following Kraus maps:
    $$f \widehat{\ox} g := \begin{tikzpicture}
        \begin{pgfonlayer}{nodelayer}
            \node [style=circle, scale=2] (0) at (-2, 1) {};
            \node [style=none] (12) at (-2, 1) {$f$};
            \node [style=circle, scale=2] (1) at (0, 1) {};
            \node [style=none] (13) at (0, 1) {$g$};
            \node [style=ox] (2) at (-1, 2) {};
            \node [style=oa] (3) at (-2, -0.5) {};
            \node [style=ox] (4) at (0, -0.5) {};
            \node [style=none] (5) at (0, -1.5) {};
            \node [style=none] (6) at (-2, -1.5) {};
            \node [style=none] (7) at (-1, 3) {};
            \node [style=none] (8) at (-2.5, -0) {};
            \node [style=none] (9) at (-1.5, -0) {};
            \node [style=none] (10) at (-1.5, -1) {};
            \node [style=none] (11) at (-2.5, -1) {};
        \end{pgfonlayer}
        \begin{pgfonlayer}{edgelayer}
            \draw (7.center) to (2);
            \draw [in=90, out=-90, looseness=1.00] (1) to (3);
            \draw [in=105, out=-75, looseness=1.00] (0) to (4);
            \draw [in=90, out=-165, looseness=1.00] (2) to (0);
            \draw [in=90, out=-15, looseness=1.00] (2) to (1);
            \draw (4) to (5.center);
            \draw (3) to (6.center);
            \draw [bend right, looseness=1.25] (0) to (3);
            \draw [bend left, looseness=1.25] (1) to (4);
            \draw (8.center) to (9.center);
            \draw (9.center) to (10.center);
            \draw (10.center) to (11.center);
            \draw (11.center) to (8.center);
        \end{pgfonlayer}
    \end{tikzpicture}
    \hspace*{1cm}
    f \widehat{\oa} g := \begin{tikzpicture}
        \begin{pgfonlayer}{nodelayer}
            \node [style=circle, scale=2] (0) at (-2, 1) {};
            \node [style=none] (12) at (-2, 1) {$f$};
            \node [style=circle, scale=2] (1) at (0, 1) {};
            \node [style=none] (13) at (0, 1) {$g$};
            \node [style=oa] (2) at (-1, 2) {};
            \node [style=oa] (3) at (-2, -0.5) {};
            \node [style=oa] (4) at (0, -0.5) {};
            \node [style=none] (5) at (0, -1.5) {};
            \node [style=none] (6) at (-2, -1.5) {};
            \node [style=none] (7) at (-1, 3) {};
            \node [style=none] (8) at (-2.5, -0) {};
            \node [style=none] (9) at (-1.5, -0) {};
            \node [style=none] (10) at (-1.5, -1) {};
            \node [style=none] (11) at (-2.5, -1) {};
        \end{pgfonlayer}
        \begin{pgfonlayer}{edgelayer}
            \draw (7.center) to (2);
            \draw [in=90, out=-90, looseness=1.00] (1) to (3);
            \draw [in=105, out=-75, looseness=1.00] (0) to (4);
            \draw [in=90, out=-165, looseness=1.00] (2) to (0);
            \draw [in=90, out=-15, looseness=1.00] (2) to (1);
            \draw (4) to (5.center);
            \draw (3) to (6.center);
            \draw [bend right, looseness=1.25] (0) to (3);
            \draw [bend left, looseness=1.25] (1) to (4);
            \draw (8.center) to (9.center);
            \draw (9.center) to (10.center);
            \draw (10.center) to (11.center);
            \draw (11.center) to (8.center);
        \end{pgfonlayer}
    \end{tikzpicture}$$
    
    The units for $\widehat{\ox}$ and $\widehat{\oa}$ are $\top$ and $\bot$ respectively. 
    
    The linear distribution maps and all the basic natural isomorphisms  are inherited from $\X$ by composing each one of them with $(u_\oa^L)^{-1}$ i.e.,
    
    \[
    \hfil \infer{ A \widehat{\ox} (B \widehat{\ox} C)  \xrightarrow{a_{\widehat{\ox}} := a_\ox (u_\oa^L)^{-1}} (A \widehat{\ox} B) \widehat{\ox} C \in  \CP^\infty(M: \U \to \C)} { A \ox (B \ox C) \xrightarrow{a_\ox} (A \ox B) \ox C \xrightarrow{(u_\oa^L)^{-1}}  \C \ / \sim}
    \]

    We prove that the associators and the other maps as defined above are natural isomorphisms in $\CP^\infty(M: \U \to \C)$:  From Lemma $\ref{Lemma: embedding}$, $Q: \C \hookrightarrow \CP^\infty(M: \U \to \C)$ is functorial which means that all commuting diagrams and isomorphisms are preserved. It remains to show that Q preserves the linear structure and the mix map:
    
    \begin{itemize}
    
    \item $Q$ preserves $\ox$: Suppose $f:A \to A'$ and $g:B \to B' \in \C$. Then,  $Q(f) \widehat{\ox} Q(g) = Q(f \ox g)$: 
    
    \begin{align*}
    Q(f) \widehat{\ox} Q(g) &:= A \widehat{\ox} B 
    \xrightarrow{f (u_\oa^L)^{-1} \widehat{\ox} g (u_\oa^L)^{-1}}
    (\bot \widehat{\oa} \bot) \widehat{\oa} (A' \widehat{\ox}  B') \\
    Q(f \ox g) &:=  A \widehat{\ox} B \xrightarrow{ (f \ox g) u_\oa^{-1}} \bot \widehat{\oa} (A' \widehat{\ox} B') 
    \end{align*}
    
    Since, $\bot \oa \bot \xrightarrow{u_\oa^L} \bot$ is a unitary isomorphism and, $f (u_\oa^L)^{-1} \widehat{\ox} g (u_\oa^L)^{-1} (u_\oa^L \oa 1) = (f \ox g) (u_\oa^L)^{-1} \in \C$, by Lemma \ref{Lemma: unitary equivalence}, 
    
    \[ f (u_\oa^L)^{-1}) \widehat{\ox} (g (u_\oa^L)^{-1} \sim (f \ox g) (u_\oa^L)^{-1} \]
    
    Therefore, $Q(f) \widehat{\ox} Q(g) = Q(f \ox g)$. Similarly,  $Q(f) \widehat{\oa} Q(g) = Q(f \oa g)$.
    
    \item $Q$ preserves all basic natural isomorphisms  (associators, unitors, symmetry maps, mix map) and linear distributions:
    
    To prove that $a_{\widehat{\ox}}$ is natural in $\CP^\infty(M:\U \to \C)$, we need to prove that the following diagram commutes in $\CP^\infty( M: \U \to \C )$:
    
    \[
    \xymatrix{
    A \widehat{\ox} (B \widehat{\ox} C) \ar[r]^{a_{\widehat{\ox}}}  \ar[d]_{(f \widehat{\ox} g) \widehat{\ox} h}
    & (A \widehat{\ox} B) \widehat{\ox} C \ar[d]^{f \widehat{\ox} (g \widehat{\ox} h)} \\
    A' \widehat{\ox} (B' \widehat{\ox} C') \ar[r]^{a_{\widehat{\ox}}}
    & (A' \widehat{\ox} B') \widehat{\ox} C'
    }
    \]
    
    In other words, we need to show that the two compositions in $\C$ are equivalent as Kraus maps. This is follows from Lemma \ref{Lemma: unitary equivalence} as there is a unitary isomorphism between the ancillary objects $\bot \oa (U_1 \oa (U_2 \oa U_3))$ and $\bot \oa (U_1 \oa U_2) \oa U_3$.  Similarly, we can show that the other basic linearly distrbutive transformations as defined are natural transformations. Since $Q$ is functorial, it preserves isomorphisms and commuting diagrams so that the coherence diagrams automatically commute.
    \end{itemize}    
\end{proof}

$\CP^\infty(M: \U \to \C)$, in general, does not have a dagger even when $\C$ is a $\dagger$-isomix category. 
However, if $M: \U \to \C$ is a $*$-MUdC, that is, a mixed unitary category in which every object in $\U$ has a unitary 
dual and if $\C$ is a $\dagger$-isomix $*$-autonomous category, then $\CP^\infty(M: \U \to \C)$ has an obvious dagger 
as shown below:

\begin{lemma}
    If {\em M:} $\U \to \C$ is a $*${\em -MUdC} then $\CP^\infty($M: $\U \to \C)$ is a $\dagger$-isomix category and 
    \[ N: \U \to \CP^\infty(M: \U \to \C) ;

    \] is a $*${\em -MUdC}.
    \end{lemma}
    \begin{proof}
    (Sketch) We first oberve that $\CP^\infty(M: \U \to \C)$ is a $\dagger$-isomix category. We have already proven that $\CP^\infty(M: \U \to \C)$ is an isomix category. Suppose  $\X$ is a $*$-MUdC, then the $\dagger$ functor for $\CP^\infty(M: \U \to \C)$ is defined as follows: Suppose $(f, U): A \to B$ and $(\eta, \epsilon): V \dashvv_{~u} U$ with $A' \dashvv A$ and $B' \dashvv B \in \C$ then,
    \[
    \dagger : \CP^\infty(M: \U \to \C)^{\rm op} \to \CP^\infty(M: \U \to \C); 
	%
    \]
    
    We need to prove that $\dagger$ is well-defined:  $f \sim g \Rightarrow f^\dagger \sim g^\dagger$. 
    \begin{align*}
    %
    \end{align*}
    
    The equality is proved by using the snake diagrams and {\bf [U.2]}, {\bf [U.5](b)}, and {\bf [Udual.]}. 
    
    Suppose $f: A \rightarrow U_1 \oa B$ and $g: B \to U_2 \oa C$ with $(\eta_1, \epsilon_1): U_1 \dashvv_{~u} V_1$ and 
    $(\eta_2, \epsilon_2): U_2 \dashvv_{~u} V_2$, then $\dagger$ preserves composition, that is  $(fg)^\dagger = g^\dagger f^\dagger$:
    
    \[
    \hfil \infer{ (fg)^\dagger: C^\dagger \rightarrow (V_1 \ox V_2)^\dagger \oa A^\dag} { (fg): A \to (U_1 \oa U_2) \oa C}
    \]
    
    \[
    \hfil \infer{ (g^\dagger f^\dagger): C^\dagger \rightarrow (V_2^\dag \ox V_1^\dagger) \oa A^\dag } 
    { g^\dagger: C^\dagger \to U_2^\dagger \oa B^\dagger ~~~~~~~~~ f^\dag: B^\dag \to U_1^\dag \oa A^\dag}
    \]
    
    To prove that $(fg)^\dagger = (g^\dagger f^\dagger)$ in $\CP^\infty(M: \U \to \C)$, represent the maps in circuit calculus 
    and fuse the $\dagger$-boxes. Once the $\dagger$-boxes are fused, use 
    Lemma \ref{Lemma: unitary equivalence} to show that both Kraus operations belong to the same equivalence. 
    $\dag$ preserves identity map since $((u_\oa^R)^{-1}, u_\oa^L): \top \dashvv_{~u} \bot$. Hence, $\dag$ is a functor.
    
    All the basic natural isomorphisms associated with $\dagger$ functor - $\lambda_\oa, \lambda_\ox, 
    \lambda_\bot, \lambda_\top, \iota$ - are lifted from $\C$ using $Q: \C \hookrightarrow \CP^\infty(M: \U \to \C)$ 
    which is defined in Lemma \ref{Lemma: embedding}. The lifted morphisms are natural in $\CP^\infty(M : \U \to \C)$ 
    since their ancillaries are unitarily isomorphic. Since $\dagger$ is functorial, all commuting diagrams are preserved. 
    By the same argument, unitary structure is preserved under $Q$.
    
    Thus, $\CP^\infty(M: \U \to \C)$ is a mixed unitary category: as $Q$ preserves all unitary linear adjoints this makes  
    $\CP^\infty(M: \U \to \C)$ a $*$-MUdC. 
    \end{proof}

\begin{lemma}
The $\CP^\infty$-construction is functorial on $*$-MUdCs.
\end{lemma}
\begin{proof}
    Let $M: \U \to \C$ and $M': \U' \to \C'$ be $*$-MUdCs and the following square be a MUC morphism:
    \[
    \xymatrix{
    \U  \ar[rr]^{M} \ar[d]_{F_u} & \ar@{}[d]|{\Downarrow~\alpha} & \C \ar[d]^{F} \\
    \U' \ar[rr]_{M'} & &\C'
    }
    \]
    $F_u$ and $F$ are $\dagger$-isomix functors and $F_u$ preserves unitary structure i.e., $F_u(\varphi_A) \rho^{F_u} = 
    \varphi_{F_u(A)}$ and ($n_\bot^{F_u}$ or $m_\top^{F_u}$) is a unitary isomorphism. $\alpha$ is a $\dagger$-linear 
    natural isomorphism. 
    
    Then, the $\CP^\infty$-construction is functorial if there is a MUC morphism:
    \[
    \xymatrix{
    \U  \ar[rr]^{MQ} \ar[d]_{F_u} & \ar@{}[d]|{\Downarrow~\alpha} & \CP^\infty(M: \U \to \C) \ar[d]^{G} \\
    \U' \ar[rr]_{M'Q} & & \CP^\infty(M': \U' \to \C')
    }
    \]
    
    Recall the functor $Q: \C \hookrightarrow \CP^\infty(M: \U \to \C)$ from Lemma \ref{Lemma: embedding}. 
    Define $G: \CP^\infty(M: \U \to \C) \to \CP^\infty(M': \U' \to \C')$ as follows:
    
    \[ G: \CP^\infty(M: \U \to \C) \to \CP^\infty(M': \U' \to \C') ; ~~~~
    \begin{matrix}
    \xymatrix{
    A \ar[d]_{[(f, U)]} \\
    B
    }
    \end{matrix}
    \mapsto
    \begin{matrix}
    \xymatrix{
    F(A) \ar[d]^{[(F(f) n_\oa^F (\alpha \oa 1), F_u(U))]} \\
    F(B)
    }
    \end{matrix}
    \]
    The action of functor $G: \CP^\infty(M: \U \to \C) \to \CP^\infty(M': \U' \to \C')$ on maps is drawn as follows:
\[    \begin{tikzpicture}
	\begin{pgfonlayer}{nodelayer}
		\node [style=onehalfcircle] (15) at (2.75, 1) {};
		\node [style=none] (16) at (2.25, 0.25) {};
		\node [style=none] (17) at (3.25, 0.25) {};
		\node [style=none] (18) at (2.75, 2.25) {};
		\node [style=none] (19) at (2.25, -0.75) {};
		\node [style=none] (20) at (3.25, -0.75) {};
		\node [style=none] (13) at (2.75, 1) {$f$};
	\end{pgfonlayer}
	\begin{pgfonlayer}{edgelayer}
		\draw (18.center) to (15);
		\draw [bend right] (15) to (16.center);
		\draw [bend left] (15) to (17.center);
		\draw (19.center) to (16.center);
		\draw (20.center) to (17.center);
	\end{pgfonlayer}
\end{tikzpicture}    \mapsto     \begin{tikzpicture}
	\begin{pgfonlayer}{nodelayer}
		\node [style=onehalfcircle] (0) at (0, 1.25) {};
		\node [style=none] (1) at (0.75, -0.75) {};
		\node [style=none] (2) at (0, 2.25) {};
		\node [style=onehalfcircle] (3) at (-0.75, 0) {};
		\node [style=none] (4) at (-0.75, -0.75) {};
		\node [style=none] (5) at (-1, 0.5) {};
		\node [style=none] (6) at (1, 0.5) {};
		\node [style=none] (7) at (1, 1.75) {};
		\node [style=none] (8) at (-1, 1.75) {};
		\node [style=none] (9) at (-0.75, 1.5) {$F$};
		\node [style=none] (10) at (0.75, 2.25) {F(A)};
		\node [style=none] (11) at (1.25, -0.5) {$F(B)$};
		\node [style=none] (12) at (-2, -0.5) {$M'(F_u(U))$};
		\node [style=none] (13) at (0, 1.25) {$f$};
		\node [style=none] (14) at (-0.75, 0) {$\alpha$};
	\end{pgfonlayer}
	\begin{pgfonlayer}{edgelayer}
		\draw (2.center) to (0);
		\draw [in=90, out=-30, looseness=1.25] (0) to (1.center);
		\draw [in=90, out=-150, looseness=1.25] (0) to (3);
		\draw (3) to (4.center);
		\draw (8.center) to (7.center);
		\draw (6.center) to (7.center);
		\draw (6.center) to (5.center);
		\draw (5.center) to (8.center);
	\end{pgfonlayer}
\end{tikzpicture} 	 \]
    The natural isomorphism $\alpha$ lifts to $\CP^{\infty}(M': \U' \to \C')$ as follows: \[ \alpha':= [(\alpha_u (u_\oa^L)^{-1} 
    ((n_\bot^{M'})^{-1} \oa 1), \bot)]: F(M(U)) \to M'(F_u(U)) \]
    
    It is immediate that $\alpha'_U:  (G(Q(M(U_1) := F(M(U_1))) \to (Q(M'(F_u(U_2))) := M'(F_u(U_1))$ is an isomorphism. 
    Let $U_1 \to^{f} U_2 \in \U$. To prove that $\alpha'$ is natural, we show that the following diagram commutes in 
    $\CP^\infty(M': \U' \to \C')$.
    \[
    \xymatrixcolsep{15pc}
    \xymatrix{
    F(M(U_1)) \ar[d]_{\alpha'_{U_1}} \ar[r]^{[( (F(M(f)(u_\oa^L)^{-1}((n_\bot^M)^{-1} \oa 1)) n_\oa^F (\alpha \oa 1)),F_u(\bot))]} &  
    F(M(U_2)) \ar[d]^{\alpha'_{U_2}}\\
    M'(F_u(U_1)) \ar[r]_{[(M'(F_u(f)) (u_\oa^L)^{-1} ((n_\bot^{M'})^{-1} \oa 1), \bot)]} & M'(F_u(U_2))
    } \]
    
    The underlying Kraus maps for both the compositions are equivalent since $F_u(\bot) \to^{(n_\bot^{F_u})^{-1}} \bot$ is 
    unitarily isomorphic: 
     \[   \begin{tikzpicture}
		\begin{pgfonlayer}{nodelayer}
			\node [style=circle] (0) at (0, 2) {$\alpha$};
			\node [style=onehalfcircle] (1) at (0, 0.25) {};
			\node [style=none] (2) at (1, 0) {$F_uM'$};
			\node [style=none] (3) at (-0.5, 0.75) {};
			\node [style=none] (4) at (-0.5, -0.25) {};
			\node [style=none] (5) at (1.5, -0.25) {};
			\node [style=none] (6) at (1.5, 0.75) {};
			\node [style=onehalfcircle] (7) at (-0.75, -1.5) {};
			\node [style=none] (8) at (0, -1.75) {$M'$};
			\node [style=none] (9) at (-1.25, -1) {};
			\node [style=none] (10) at (-1.25, -2) {};
			\node [style=none] (11) at (0.25, -2) {};
			\node [style=none] (12) at (0.25, -1) {};
			\node [style=none] (13) at (0.75, -4.25) {};
			\node [style=circle, scale=0.3] (14) at (0.5, -0.75) {};
			\node [style=onehalfcircle] (15) at (-2.25, 0.25) {};
			\node [style=circle, scale=0.3] (16) at (0, 1.25) {};
			\node [style=none] (17) at (-2.75, 0.75) {};
			\node [style=none] (18) at (-1.25, 0.75) {};
			\node [style=none] (19) at (-1.25, -0.25) {};
			\node [style=none] (20) at (-2.75, -0.25) {};
			\node [style=none] (21) at (-1.5, 0) {$M'$};
			\node [style=oa] (22) at (-1.75, -3.25) {};
			\node [style=none] (23) at (-1.75, -4.25) {};
			\node [style=none] (24) at (-2.5, -2.75) {};
			\node [style=none] (25) at (-1, -2.75) {};
			\node [style=none] (26) at (-1, -3.75) {};
			\node [style=none] (27) at (-2.5, -3.75) {};
			\node [style=none] (28) at (-2.25, -3.5) {$M$};
			\node [style=none] (29) at (-3, -4.25) {$M(\bot \oa \bot)$};
			\node [style=none] (30) at (2, -4) {$M'(F_u(U_2)$};
			\node [style=none] (31) at (0, 2.5) {};
			\node [style=none] (32) at (1.25, 2.25) {$F(M(U_1)$};
			\node [style=none] (33) at (-2.25, 0.25) {$\bot$};
			\node [style=none] (34) at (0, 0.25) {$f$};
			\node [style=none] (35) at (-0.75, -1.5) {$\bot$};
		\end{pgfonlayer}
		\begin{pgfonlayer}{edgelayer}
			\draw (3.center) to (6.center);
			\draw (6.center) to (5.center);
			\draw (5.center) to (4.center);
			\draw (4.center) to (3.center);
			\draw (9.center) to (12.center);
			\draw (12.center) to (11.center);
			\draw (11.center) to (10.center);
			\draw (10.center) to (9.center);
			\draw [in=90, out=-60] (1) to (13.center);
			\draw (1) to (0);
			\draw [dotted, bend right=45] (14) to (7);
			\draw [dotted, in=75, out=174] (16) to (15);
			\draw (17.center) to (18.center);
			\draw (18.center) to (19.center);
			\draw (19.center) to (20.center);
			\draw (20.center) to (17.center);
			\draw [in=120, out=-90, looseness=0.75] (15) to (22);
			\draw [in=45, out=-90] (7) to (22);
			\draw (22) to (23.center);
			\draw (24.center) to (25.center);
			\draw (25.center) to (26.center);
			\draw (26.center) to (27.center);
			\draw (27.center) to (24.center);
			\draw (31.center) to (0);
		\end{pgfonlayer}
	\end{tikzpicture}	 \sim \begin{tikzpicture}
		\begin{pgfonlayer}{nodelayer}
			\node [style=onehalfcircle] (0) at (1.75, 1.5) {};
			\node [style=onehalfcircle] (1) at (-0.25, 0) {};
			\node [style=circle] (2) at (-0.25, -1) {$\alpha$};
			\node [style=circle] (3) at (1.75, -1) {$\alpha$};
			\node [style=onehalfcircle] (4) at (1, -2.75) {};
			\node [style=none] (5) at (1.75, 2.5) {};
			\node [style=none] (6) at (1.75, -4.75) {};
			\node [style=none] (7) at (1, -3.75) {};
			\node [style=none] (8) at (-0.25, -3.75) {};
			\node [style=circle, scale=0.3] (9) at (1.75, 0.75) {};
			\node [style=circle, scale=0.3] (10) at (1.75, -1.75) {};
			\node [style=none] (11) at (0, -2.25) {};
			\node [style=none] (12) at (0, -3.25) {};
			\node [style=none] (13) at (1.5, -3.25) {};
			\node [style=none] (14) at (1.5, -2.25) {};
			\node [style=none] (15) at (0.25, -3) {$M$};
			\node [style=none] (16) at (-0.75, 2) {};
			\node [style=none] (17) at (-0.75, -0.5) {};
			\node [style=none] (18) at (2.75, -0.5) {};
			\node [style=none] (19) at (2.75, 2) {};
			\node [style=none] (20) at (-0.25, 1.5) {$MF$};
			\node [style=oa] (21) at (0.5, -4.25) {};
			\node [style=none] (22) at (0.5, -4.75) {};
			\node [style=none] (23) at (-1, -4.75) {$M'(F_u(\bot) \oa \bot)$};
			\node [style=none] (24) at (2.75, -4.75) {$M'(F_u(U_2))$};
			\node [style=none] (25) at (0.75, 2.5) {$F(M(U_1))$};
			\node [style=none] (26) at (1.75, 1.5) {$f$};
			\node [style=none] (27) at (-0.25, 0) {$\bot$};
			\node [style=none] (28) at (1, -2.75) {$\bot$};
		\end{pgfonlayer}
		\begin{pgfonlayer}{edgelayer}
			\draw (5.center) to (0);
			\draw (0) to (3);
			\draw (3) to (6.center);
			\draw [dotted, in=75, out=172, looseness=0.75] (9) to (1);
			\draw (1) to (2);
			\draw [dotted, bend right] (10) to (4);
			\draw (4) to (7.center);
			\draw (2) to (8.center);
			\draw (11.center) to (12.center);
			\draw (12.center) to (13.center);
			\draw (13.center) to (14.center);
			\draw (14.center) to (11.center);
			\draw (16.center) to (19.center);
			\draw (19.center) to (18.center);
			\draw (18.center) to (17.center);
			\draw (17.center) to (16.center);
			\draw [in=165, out=-90] (8.center) to (21);
			\draw [in=30, out=-90] (7.center) to (21);
			\draw (21) to (22.center);
		\end{pgfonlayer}
	\end{tikzpicture}	  \]
    Hence, the diagram commutes and $\alpha'$ is natural in $\CP^\infty(M': \U' \to \C')$.    
\end{proof}

The following table summarizes the structures inherited by $\CP^\infty(M: \U \to \C)$ from $ M: \U \to \C$:
\vspace{-0.20cm}
\begin{center}
\begin{tabular}{ |c|c| }
\hline
$ M: \U \to \C $ & $ \CP^\infty (M: \U \to \C) $  \\ 
\hline \hline
mixed unitary category & isomix category  \\  \hline
 $*$-mixed unitary category with unitary duals &  $*$-mixed unitary category with unitary duals\\ \hline
 $\dagger$-symmetric monoidal category & symmetric monoidal category    \\ \hline
 $\dagger$-compact closed category & $\dagger$-compact closed category ($\CP^\infty(\X) \simeq$ CPM$(\X)$)\\ \hline
\end{tabular}
\end{center}

\section{Environment structure for MUCs}

In this section, we describe when a given isomix category is of the form $\CP^\infty( M:\U \to \C)$ by generalizing the 
notion of environment structures from $\dagger$-symmetric monoidal categories (see \cite{CoH16}) to mixed unitary categories. 
We then show that an environment structure over $M$ which has purification is isomorphic to $\CP^\infty(M: \U \to \C)$.

\subsection{Environment structure}

We first define environment structure for MUCs and give examples:
\begin{definition} 
\label{Definition: env structure}
An {\bf environment structure} for a mixed unitary category $M: \U \to \C$ is a strict isomix functor $F: \C \to \D$ where $\D$ is an isomix category and 
 a family of maps $\envmap_U: F(M(U)) \to \bot$ indexed by objects $ U \in \U$ such that the following conditions hold:
\begin{enumerate}[{\bf [Env.1]}]
\item For $U, V \in \U$, the following diagrams commute:
\[ \mbox{\bf (a)} ~~~~~~ \xymatrixcolsep{3.5pc}
\xymatrix{
F(M(U)) \ox F(M(V)) \ar[r]^{\mx} \ar[d]_{m_\ox} & F(M(U)) \oa F(M(V)) \ar[r]^{~~~~~~\envmap \oa \envmap} & \bot \oa \bot \ar[d]^{u_\oa} \\
F(M(U \ox V)) \ar[rr]_{\envmap} & & \bot
}
\]
\[ \mbox{\bf (b)} ~~~~~~~~~~~~~~~~~~~~~~ \xymatrixcolsep{7pc} \xymatrix{
F(M(U \oa V)) \ar@/^1pc/[drr]^{\envmap } \ar[d]_{n_\oa}\\
F(M(U)) \oa F(M(V)) \ar[r]_{\envmap \oa \envmap} & \bot \oa \bot \ar[r]_{u_\oa} & \bot
} \]
\item  Kraus maps $(f, U) \sim (g, V) \in \C$ if and only if  the following diagram commutes:
\[ \xymatrix{
& F(M(A)) \ar[ld]_{F(M(f))} \ar[dr]^{F(M(g))} & \\
F(M(U \oa B)) \ar[d]_{\nu_\ox} &  & F(M(U \oa B)) \ar[d]_{\nu_\ox} \\
F(M(U)) \oa MF(B) \ar[d]_{\envmap \oa 1} &  & F(M(U)) \oa F(M(B)) \ar[d]_{\envmap \oa 1} \\
\bot \oa F(M(B)) \ar@{=}[rr] &  & \bot \oa F(M(B)) 
} \]
\end{enumerate}
\end{definition}

The conditions are represented diagrammatically as follows:
\[
\mbox{ \bf{ [Env.1a]}}~~~~~ 
\begin{tikzpicture}
	\begin{pgfonlayer}{nodelayer}
		\node [style=ox] (0) at (1.75, 1.5) {};
		\node [style=none] (1) at (1.75, -0) {};
		\node [style=none] (2) at (1, 3) {};
		\node [style=none] (3) at (2.5, 3) {};
		\node [style=none] (4) at (0.5, 2.25) {};
		\node [style=none] (5) at (0.5, 0.75) {};
		\node [style=none] (6) at (2.75, 0.75) {};
		\node [style=none] (7) at (2.75, 2.25) {};
		\node [style=none] (8) at (2.25, 1) {$MF$};
		\node [style=none] (9) at (1.5, -0) {};
		\node [style=none] (10) at (2, -0) {};
		\node [style=none] (11) at (1.6, -0.1) {};
		\node [style=none] (12) at (1.9, -0.1) {};
		\node [style=none] (13) at (1.7, -0.2) {};
		\node [style=none] (14) at (1.8, -0.2) {};
		\node [style=none] (15) at (1.75, -0.5) {};
		\node [style=none] (16) at (1.75, -1) {};
		\node [style=circle, scale=2] (17) at (1.75, -0.1) {};
	\end{pgfonlayer}
	\begin{pgfonlayer}{edgelayer}
		\draw (4.center) to (7.center);
		\draw (7.center) to (6.center);
		\draw (6.center) to (5.center);
		\draw (5.center) to (4.center);
		\draw (0) to (1.center);
		\draw [bend right=15, looseness=1.00] (2.center) to (0);
		\draw [bend right=15, looseness=1.00] (0) to (3.center);
		\draw (15.center) to (16.center);
		\draw (9.center) to (10.center);
		\draw (11.center) to (12.center);
		\draw (13.center) to (14.center);
	\end{pgfonlayer}
\end{tikzpicture} = \begin{tikzpicture}
	\begin{pgfonlayer}{nodelayer}
		\node [style=none] (0) at (-3, 3) {};
		\node [style=none] (1) at (-1, 3) {};
		\node [style=circle, scale=0.2] (2) at (-3, 2.5) {};
		\node [style=circle, scale=0.2] (3) at (-1, 1.5) {};
		\node [style=map] (4) at (-2, 2) {};
		\node [style=none] (5) at (-3, 1) {};
		\node [style=none] (6) at (-1, 1) {};
		\node [style=none] (7) at (-3.25, 1) {};
		\node [style=none] (8) at (-2.75, 1) {};
		\node [style=none] (9) at (-3.15, 0.9) {};
		\node [style=none] (10) at (-2.85, 0.9) {};
		\node [style=none] (11) at (-3.05, 0.8) {};
		\node [style=none] (12) at (-2.95, 0.8) {};
		\node [style=none] (13) at (-1.25, 1) {};
		\node [style=none] (14) at (-0.75, 1) {};
		\node [style=none] (15) at (-1.15, 0.9) {};
		\node [style=none] (16) at (-0.85, 0.9) {};
		\node [style=none] (17) at (-1.05, 0.8) {};
		\node [style=none] (18) at (-0.95, 0.8) {};
		\node [style=none] (19) at (-1, 0.5) {};
		\node [style=none] (20) at (-3, 0.5) {};
		\node [style=none] (21) at (-1, -0.25) {};
		\node [style=circle] (22) at (-3, -0.25) {$\bot$};
		\node [style=circle, scale=2] (23) at (-3, 0.85) {};
		\node [style=circle, scale=2] (24) at (-1, 0.85) {};
	\end{pgfonlayer}
	\begin{pgfonlayer}{edgelayer}
		\draw (0.center) to (5.center);
		\draw (1.center) to (6.center);
		\draw [dotted, bend left=45, looseness=1.25] (2) to (4);
		\draw [dotted, bend right=45, looseness=1.25] (4) to (3);
		\draw (7.center) to (8.center);
		\draw (9.center) to (10.center);
		\draw (11.center) to (12.center);
		\draw (13.center) to (14.center);
		\draw (15.center) to (16.center);
		\draw (17.center) to (18.center);
		\draw (20.center) to (22);
		\draw (19.center) to (21.center);
	\end{pgfonlayer}
\end{tikzpicture}  ~~~~~~~~~~ 
\mbox{ \bf{ [Env.1b]}}~~~~~\begin{tikzpicture}
	\begin{pgfonlayer}{nodelayer}
		\node [style=oa] (0) at (1.75, 1.5) {};
		\node [style=none] (1) at (1.75, 3) {};
		\node [style=none] (2) at (1, 0) {};
		\node [style=none] (3) at (2.5, 0) {};
		\node [style=none] (4) at (0.5, 0.75) {};
		\node [style=none] (5) at (0.5, 2.25) {};
		\node [style=none] (6) at (2.75, 2.25) {};
		\node [style=none] (7) at (2.75, 0.75) {};
		\node [style=none] (8) at (0.85, 1) {$MF$};
		\node [style=none] (9) at (2.25, 0) {};
		\node [style=none] (10) at (2.75, 0) {};
		\node [style=none] (11) at (2.35, -0.1) {};
		\node [style=none] (12) at (2.65, -0.1) {};
		\node [style=none] (13) at (2.45, -0.2) {};
		\node [style=none] (14) at (2.55, -0.2) {};
		\node [style=none] (15) at (2.5, -0.5) {};
		\node [style=none] (16) at (2.5, -1.25) {};
		\node [style=circle, scale=2] (17) at (2.5, -0.1) {};
		\node [style=none] (18) at (0.75, 0) {};
		\node [style=none] (19) at (1.25, 0) {};
		\node [style=none] (20) at (0.85, -0.1) {};
		\node [style=none] (21) at (1.15, -0.1) {};
		\node [style=none] (22) at (0.95, -0.2) {};
		\node [style=none] (23) at (1.05, -0.2) {};
		\node [style=none] (24) at (1, -0.5) {};
		\node [style=circle, scale=2] (25) at (1, -0.1) {};
		\node [style=circle] (26) at (1, -1.25) {$\bot$};
	\end{pgfonlayer}
	\begin{pgfonlayer}{edgelayer}
		\draw (4.center) to (7.center);
		\draw (7.center) to (6.center);
		\draw (6.center) to (5.center);
		\draw (5.center) to (4.center);
		\draw (0) to (1.center);
		\draw [bend left=15] (2.center) to (0);
		\draw [bend left=15] (0) to (3.center);
		\draw (15.center) to (16.center);
		\draw (9.center) to (10.center);
		\draw (11.center) to (12.center);
		\draw (13.center) to (14.center);
		\draw (18.center) to (19.center);
		\draw (20.center) to (21.center);
		\draw (22.center) to (23.center);
		\draw (24.center) to (26);
	\end{pgfonlayer}
\end{tikzpicture} = \begin{tikzpicture}
	\begin{pgfonlayer}{nodelayer}
		\node [style=none] (0) at (-3, 3) {};
		\node [style=none] (1) at (-3, 0.5) {};
		\node [style=none] (2) at (-3.25, 0.5) {};
		\node [style=none] (3) at (-2.75, 0.5) {};
		\node [style=none] (4) at (-3.15, 0.4) {};
		\node [style=none] (5) at (-2.85, 0.4) {};
		\node [style=none] (6) at (-3.05, 0.3) {};
		\node [style=none] (7) at (-2.95, 0.3) {};
		\node [style=none] (8) at (-3, -0) {};
		\node [style=none] (9) at (-3, -1) {};
		\node [style=circle, scale=2] (10) at (-3, 0.4) {};
		\node [style=none] (11) at (-3, 3.5) {$F(M(U \oa V))$};
	\end{pgfonlayer}
	\begin{pgfonlayer}{edgelayer}
		\draw (2.center) to (3.center);
		\draw (4.center) to (5.center);
		\draw (6.center) to (7.center);
		\draw (0.center) to (1.center);
		\draw (8.center) to (9.center);
	\end{pgfonlayer}
\end{tikzpicture} \] \[ \mbox{ \bf{ [Env.2]}}~~~~~
\begin{tikzpicture}[scale=0.8]
	\begin{pgfonlayer}{nodelayer}
		\node [style=circle] (0) at (0, -0) {$f$};
		\node [style=none] (1) at (0, 1.5) {};
		\node [style=none] (2) at (-0.5, -1) {};
		\node [style=none] (3) at (0.5, -1) {};
		\node [style=none] (4) at (-0.5, -1.75) {};
		\node [style=none] (5) at (0.5, -2.75) {};
		\node [style=none] (6) at (-1, -1) {};
		\node [style=none] (7) at (1, -1) {};
		\node [style=none] (8) at (-1, 0.5) {};
		\node [style=none] (9) at (1, 0.5) {};
		\node [style=none] (10) at (-0.75, -1.75) {};
		\node [style=none] (11) at (-0.25, -1.75) {};
		\node [style=none] (12) at (-0.65, -1.85) {};
		\node [style=none] (13) at (-0.35, -1.85) {};
		\node [style=none] (14) at (-0.55, -1.95) {};
		\node [style=none] (15) at (-0.45, -1.95) {};
		\node [style=none] (16) at (-0.5, -2.25) {};
		\node [style=circle] (17) at (-0.5, -3.25) {$\bot$};
		\node [style=circle, scale=2] (18) at (-0.5, -1.85) {};
		\node [style=none] (19) at (0.75, -0.75) {$F$};
	\end{pgfonlayer}
	\begin{pgfonlayer}{edgelayer}
		\draw (1.center) to (0);
		\draw [bend right=15, looseness=1.00] (0) to (2.center);
		\draw [bend left=15, looseness=1.00] (0) to (3.center);
		\draw (2.center) to (4.center);
		\draw (3.center) to (5.center);
		\draw (8.center) to (6.center);
		\draw (6.center) to (7.center);
		\draw (7.center) to (9.center);
		\draw (9.center) to (8.center);
		\draw (16.center) to (17);
		\draw (10.center) to (11.center);
		\draw (12.center) to (13.center);
		\draw (14.center) to (15.center);
	\end{pgfonlayer}
\end{tikzpicture} = \begin{tikzpicture}[scale=0.8]
	\begin{pgfonlayer}{nodelayer}
		\node [style=circle] (0) at (0, -0) {$g$};
		\node [style=none] (1) at (0, 1.5) {};
		\node [style=none] (2) at (-0.5, -1) {};
		\node [style=none] (3) at (0.5, -1) {};
		\node [style=none] (4) at (-0.5, -1.75) {};
		\node [style=none] (5) at (0.5, -2.75) {};
		\node [style=none] (6) at (-1, -1) {};
		\node [style=none] (7) at (1, -1) {};
		\node [style=none] (8) at (-1, 0.5) {};
		\node [style=none] (9) at (1, 0.5) {};
		\node [style=none] (10) at (-0.75, -1.75) {};
		\node [style=none] (11) at (-0.25, -1.75) {};
		\node [style=none] (12) at (-0.65, -1.85) {};
		\node [style=none] (13) at (-0.35, -1.85) {};
		\node [style=none] (14) at (-0.55, -1.95) {};
		\node [style=none] (15) at (-0.45, -1.95) {};
		\node [style=none] (16) at (-0.5, -2.25) {};
		\node [style=circle] (17) at (-0.5, -3.25) {$\bot$};
		\node [style=circle, scale=2] (18) at (-0.5, -1.85) {};
		\node [style=none] (19) at (0.75, -0.75) {$F$};
	\end{pgfonlayer}
	\begin{pgfonlayer}{edgelayer}
		\draw (1.center) to (0);
		\draw [bend right=15, looseness=1.00] (0) to (2.center);
		\draw [bend left=15, looseness=1.00] (0) to (3.center);
		\draw (2.center) to (4.center);
		\draw (3.center) to (5.center);
		\draw (8.center) to (6.center);
		\draw (6.center) to (7.center);
		\draw (7.center) to (9.center);
		\draw (9.center) to (8.center);
		\draw (16.center) to (17);
		\draw (10.center) to (11.center);
		\draw (12.center) to (13.center);
		\draw (14.center) to (15.center);
	\end{pgfonlayer}
\end{tikzpicture} \Leftrightarrow \begin{tikzpicture}
	\begin{pgfonlayer}{nodelayer}
		\node [style=circle] (0) at (0, -0) {$f$};
		\node [style=none] (1) at (0, 1.5) {};
		\node [style=none] (2) at (-0.5, -1) {};
		\node [style=none] (3) at (0.5, -1) {};
		\node [style=none] (4) at (-0.5, -2.25) {};
		\node [style=none] (5) at (0.5, -2.25) {};
	\end{pgfonlayer}
	\begin{pgfonlayer}{edgelayer}
		\draw (1.center) to (0);
		\draw [bend right=15, looseness=1.00] (0) to (2.center);
		\draw [bend left=15, looseness=1.00] (0) to (3.center);
		\draw (2.center) to (4.center);
		\draw (3.center) to (5.center);
	\end{pgfonlayer}
\end{tikzpicture} \sim \begin{tikzpicture}
	\begin{pgfonlayer}{nodelayer}
		\node [style=circle] (0) at (0, -0) {$g$};
		\node [style=none] (1) at (0, 1.5) {};
		\node [style=none] (2) at (-0.5, -1) {};
		\node [style=none] (3) at (0.5, -1) {};
		\node [style=none] (4) at (-0.5, -2.25) {};
		\node [style=none] (5) at (0.5, -2.25) {};
	\end{pgfonlayer}
	\begin{pgfonlayer}{edgelayer}
		\draw (1.center) to (0);
		\draw [bend right=15, looseness=1.00] (0) to (2.center);
		\draw [bend left=15, looseness=1.00] (0) to (3.center);
		\draw (2.center) to (4.center);
		\draw (3.center) to (5.center);
	\end{pgfonlayer}
\end{tikzpicture} \]

\begin{definition}
An environment structure $F: \C \to \D$ with $\envmap$ for a mixed unitary category $M: \U \to \C$ 
has {\bf purification} if 
\begin{itemize}
\item $F$ is bijective on objects, and 
\item for all $f: F(A) \to F(B) \in \D$, there exists a Kraus map $(f', U): A \to B \in \C$ such that 
\[
\mbox{ \bf{[Env.3]} }~~~~~~
 f = \begin{tikzpicture}
	\begin{pgfonlayer}{nodelayer}
		\node [style=circle, scale=2] (0) at (0, -0) {};
		\node [style=none] (100) at (0, -0) {$f'$};
		\node [style=none] (1) at (0, 1.5) {};
		\node [style=none] (2) at (-0.5, -1) {};
		\node [style=none] (3) at (0.5, -1) {};
		\node [style=none] (4) at (-0.5, -1.75) {};
		\node [style=none] (5) at (0.5, -2.75) {};
		\node [style=none] (6) at (-1, -1) {};
		\node [style=none] (7) at (1, -1) {};
		\node [style=none] (8) at (-1, 0.5) {};
		\node [style=none] (9) at (1, 0.5) {};
		\node [style=none] (10) at (-0.75, -1.75) {};
		\node [style=none] (11) at (-0.25, -1.75) {};
		\node [style=none] (12) at (-0.65, -1.85) {};
		\node [style=none] (13) at (-0.35, -1.85) {};
		\node [style=none] (14) at (-0.55, -1.95) {};
		\node [style=none] (15) at (-0.45, -1.95) {};
		\node [style=none] (16) at (-0.5, -2.25) {};
		\node [style=circle] (17) at (-0.5, -3.25) {$\bot$};
		\node [style=circle, scale=2] (18) at (-0.5, -1.85) {};
		\node [style=none] (19) at (0.75, -0.75) {$F$};
	\end{pgfonlayer}
	\begin{pgfonlayer}{edgelayer}
		\draw (1.center) to (0);
		\draw [bend right=15, looseness=1.00] (0) to (2.center);
		\draw [bend left=15, looseness=1.00] (0) to (3.center);
		\draw (2.center) to (4.center);
		\draw (3.center) to (5.center);
		\draw (8.center) to (6.center);
		\draw (6.center) to (7.center);
		\draw (7.center) to (9.center);
		\draw (9.center) to (8.center);
		\draw (16.center) to (17);
		\draw (10.center) to (11.center);
		\draw (12.center) to (13.center);
		\draw (14.center) to (15.center);
	\end{pgfonlayer}
\end{tikzpicture}
\]
Equationally, 
\[  F(A) \to^{f} F(B) = F(A) \to^{F(f')} F(M(U) \oa B) \to^{n_\oa} M(F(U)) \oa F(B) \to^{\envmap \oa 1} \bot 
\oa F(B) \to^{u_\oa} F(B) \]
\end{itemize}
\end{definition}

\begin{lemma}
\label{Lemma: Env example}
Every mixed unitary category $M: \U \to \C$ has an environment structure given by 
\[ F: \C \to  \CP^\infty(M: \U \to \C) ; 
\begin{matrix}
\xymatrix{
A \ar[d]_{f} \\
B
} \end{matrix}   \mapsto  \begin{matrix}
\xymatrix{
A \ar[d]^{[(f(u_\oa^L)^{-1} (n_\bot^M)^{-1}, \bot)]} \\
M(\bot) \oa B
}
\end{matrix}
\]  and \[ \envmap: M(U) \to \bot \in \CP^\infty(M: \U \to \C) := [((u_\oa^R)^{-1}, U)] \] for each object $U \in \U$. 
Moreover, this environment structure has purification. 
\end{lemma}

\begin{proof}
$F$ is functorial since $f(u_\oa^L)^{-1} ((n_\bot^M)^{-1} \oa 1) \sim f \sim (u_\oa^L)^{-1} (n_\bot^M)^{-1} (1 \oa f) $. Define 
$F_\ox = F_\oa := F$. Note that, $F$ is a strict monoidal functor and an isomix functor.
In order to prove that $F: \C \to \D$ with $\envmap$ satisfy axioms for environment structures, the properties of isomix 
functor and Lemma \ref{Lemma: unitary equivalence} are used.

 To prove that this environment structure has purification, consider any map $[(f, U)]: A \to B \in \CP^\infty(M: \U \to \C)$.  
 Then, there exists a Kraus map $(f,U): A \to B \in \C$. Then, the map in equation {\bf [Env. 3] } is drawn as follows:

\[ 
\begin{tikzpicture}[scale=1.5]
\draw (0.75,1) -- (0.75,1.5);
\draw (0.25,0.5) -- (1.25, 0.5) -- (1.25, 1) -- (0.25, 1) -- (0.25, 0.5);
\draw (0.75,0.75) node{F(f)};
\draw (0.55,0.5) -- (0.55,0);
\draw (0.95,0.5) -- (0.95,0);
\draw (0.40, 0) -- (0.70, 0);
\draw (0.45, -0.1) -- (0.65, -0.1);
\draw (0.50, -0.2) -- (0.60, -0.2);
\end{tikzpicture} = \left[
\begin{tikzpicture}
	\begin{pgfonlayer}{nodelayer}
		\node [style=circle, scale=0.5] (0) at (-2.25, 0.25) {};
		\node [style=circle, scale=1.25] (1) at (-1.5, -0.5) {};
		\node [style=none] (100) at (-1.5, -0.5) {$\bot$};
		\node [style=circle, scale=0.5] (2) at (0, 0.25) {};
		\node [style=circle, scale=1.25] (3) at (-0.75, -0.5) {};
		\node [style=none] (30) at (-0.75, -0.5) {$\bot$};
		\node [style=oa] (4) at (-2.25, -1.75) {};
		\node [style=oa] (5) at (0, -1.75) {};
		\node [style=none] (6) at (0, -3) {};
		\node [style=none] (7) at (-2.5, 0.5) {$U$};
		\node [style=none] (8) at (0.25, 0.5) {$B$};
		\node [style=none] (9) at (0, -3.25) {$\bot \oa B$};
		\node [style=circle, scale=1.55] (10) at (-1, 1.5) {};
		\node [style=none] (101) at (-1, 1.5) {$f$};
		\node [style=none] (11) at (-1, 2.5) {};
		\node [style=circle, scale=0.5] (12) at (-1, 2) {};
		\node [style=circle, scale=1.25] (13) at (-3.25, -1.75) {};
		\node [style=none] (31) at (-3.25, -1.75) {$\bot$};
		\node [style=oa] (14) at (-2.75, -2.5) {};
		\node [style=none] (15) at (-2.75, -3) {};
	\end{pgfonlayer}
	\begin{pgfonlayer}{edgelayer}
		\draw [dotted, bend left, looseness=1.00] (0) to (1);
		\draw [dotted, bend right, looseness=1.25] (2) to (3);
		\draw [in=140, out=-90, looseness=1.25] (1) to (5);
		\draw (2) to (5);
		\draw [in=40, out=-90, looseness=1.00] (3) to (4);
		\draw (0) to (4);
		\draw (5) to (6.center);
		\draw [in=-150, out=90, looseness=1.00] (0) to (10);
		\draw [in=90, out=-30, looseness=1.00] (10) to (2);
		\draw (11.center) to (12);
		\draw (12) to (10);
		\draw [dotted, bend right, looseness=1.25] (12) to (13);
		\draw (14) to (15.center);
		\draw [bend right, looseness=1.00] (13) to (14);
		\draw [bend left, looseness=1.00] (4) to (14);
	\end{pgfonlayer}
\end{tikzpicture}
\right] \]  Because, $\begin{tikzpicture}
	\begin{pgfonlayer}{nodelayer}
		\node [style=circle, scale=0.5] (0) at (-2.25, 0.25) {};
		\node [style=circle, scale=1.25] (1) at (-1.5, -0.5) {};
		\node [style=none] (100) at (-1.5, -0.5) {$\bot$};
		\node [style=circle, scale=0.5] (2) at (0, 0.25) {};
		\node [style=circle, scale=1.25] (3) at (-0.75, -0.5) {};
		\node [style=none] (30) at (-0.75, -0.5) {$\bot$};
		\node [style=oa] (4) at (-2.25, -1.75) {};
		\node [style=oa] (5) at (0, -1.75) {};
		\node [style=none] (6) at (0, -3) {};
		\node [style=none] (7) at (-2.5, 0.5) {$U$};
		\node [style=none] (8) at (0.25, 0.5) {$B$};
		\node [style=none] (9) at (0, -3.25) {$\bot \oa B$};
		\node [style=circle, scale=1.55] (10) at (-1, 1.5) {};
		\node [style=none] (101) at (-1, 1.5) {$f$};
		\node [style=none] (11) at (-1, 2.5) {};
		\node [style=circle, scale=0.5] (12) at (-1, 2) {};
		\node [style=circle, scale=1.25] (13) at (-3.25, -1.75) {};
		\node [style=none] (31) at (-3.25, -1.75) {$\bot$};
		\node [style=oa] (14) at (-2.75, -2.5) {};
		\node [style=none] (15) at (-2.75, -3) {};
	\end{pgfonlayer}
	\begin{pgfonlayer}{edgelayer}
		\draw [dotted, bend left, looseness=1.00] (0) to (1);
		\draw [dotted, bend right, looseness=1.25] (2) to (3);
		\draw [in=140, out=-90, looseness=1.25] (1) to (5);
		\draw (2) to (5);
		\draw [in=40, out=-90, looseness=1.00] (3) to (4);
		\draw (0) to (4);
		\draw (5) to (6.center);
		\draw [in=-150, out=90, looseness=1.00] (0) to (10);
		\draw [in=90, out=-30, looseness=1.00] (10) to (2);
		\draw (11.center) to (12);
		\draw (12) to (10);
		\draw [dotted, bend right, looseness=1.25] (12) to (13);
		\draw (14) to (15.center);
		\draw [bend right, looseness=1.00] (13) to (14);
		\draw [bend left, looseness=1.00] (4) to (14);
	\end{pgfonlayer}
\end{tikzpicture} \sim 
\begin{tikzpicture} 
	\begin{pgfonlayer}{nodelayer}
		\node [style=circle] (0) at (-1, 0.25) {$f$};
		\node [style=none] (1) at (-1.5, -1.5) {};
		\node [style=none] (2) at (-0.5, -1.5) {};
		\node [style=none] (3) at (-1, 2.75) {};
	\end{pgfonlayer}
	\begin{pgfonlayer}{edgelayer}
		\draw [style=none] (3.center) to (0);
		\draw [style=none, in=90, out=-135, looseness=1.00] (0) to (1.center);
		\draw [style=none, in=90, out=-45, looseness=1.00] (0) to (2.center);
	\end{pgfonlayer}
\end{tikzpicture}
$, the environment structure $(Q: \C \to \CP^\infty(M: \U \to \C), \envmap)$ has purification.
\end{proof}

The following are the environment structures for our running examples:
\begin{itemize}
\item Consider the MUC, $\R^* \subset \C$. Then, \[ (\R^* \to^{Q} \CP^\infty( \R \subset \C), \envmap_r: r \to 1)\] 
is an environment structure where, $ \envmap_r := (=, 1/r) : r \to 1  $
\item Consider the MUC, ${\sf Mat}_\C \to \FMat_\C$. Then, 
\[ {\sf Mat}_\C \to^{Q} \CP^\infty({\sf Mat}_\C \subset \FMat_\C) \] is an environment structure where, 
$\envmap_{\C^n}: \C^n \to \C ; \rho \mapsto Tr(\rho) $.
\end{itemize}

\subsection{Characterizing the $\CP^\infty$ construction}

In the rest of the section, we show that any environment structure with purification is initial in the category of environment structures.  
Given Lemma \ref{Lemma: Env example}, this shows that any environment structure over a $M:\U \to \C$ with purification 
is isomorphic to $\CP^\infty(M: \U \to \C)$.  Thus, environment structures with purification captures the abstract structure 
of $\CP^\infty( M: \U \to \C)$.

\begin{definition}
Let $M: \U \to \C$ be a mixed unitary category. Define a category $\mathsf{Env}(M: \U \to \C)$ as follows:
\begin{description}
\item[Objects:] Environment structures for $M: \U \to \C$
\item[Arrows:] Suppose $D: \C \to \D$ with $\envmap$, and $D': \C \to \D'$ with $\anotherenvmap{1.8}$ are two 
environment structures. Then,  a morphism of environment structures is a strict isomix functor $F: \D \to \D'$ such that 
\begin{itemize}
 \item $DF = D'$
 \item $F(\envmap) = \anotherenvmap{1.8}$
\end{itemize}
\item[Identity arrows:] Identity functor
\item[Composition:] Linear functor composition
\end{description}
\end{definition}

\begin{lemma}
\label{Lemma: Env initial}
Let $M: \U \to \C$ be a mixed unitary category. Suppose $D: \C \to \D$ with $\envmap$ is an environment structure 
with purification, then it is initial in ${\sf Env}(M: \U \to \C)$.
\end{lemma}
\begin{proof}
Suppose $(D': \X \rightarrow \Y', \anotherenvmap{1.5}) \in \mathsf{Env(\C)}$. We show that there is a unique strict isomix functor $F: \Y \rightarrow \Y'$ such that $DF = D'$ and $F(\envmap) = \anotherenvmap{1.5}$.

Define $F: \D \rightarrow \D'$ as follows:
\begin{itemize}
\item  Since $(D, \envmap)$ has purification, $D: \C \to \D$ is bijective on objects. Then for all $A \in \D$, $A = D(X)$ for a unique $X \in \C$. Then,
\[ F(A) := D'(X) \]
\item Let $f: A \rightarrow B \in \D$. Since $(D: \C \to \D, \envmap)$ has purification,   
\begin{align*}
\begin{tikzpicture}
	\begin{pgfonlayer}{nodelayer}
		\node [style=circle] (0) at (0, 1) {$f$};
		\node [style=none] (1) at (0, 2) {};
		\node [style=none] (2) at (0, -0) {};
		\node [style=none] (3) at (0.2, 2.2) {$A$};
		\node [style=none] (4) at (0.2, -0.2) {$B$};
	\end{pgfonlayer}
	\begin{pgfonlayer}{edgelayer}
		\draw (1.center) to (0);
		\draw (0) to (2.center);
	\end{pgfonlayer}
\end{tikzpicture} = \begin{tikzpicture}
	\begin{pgfonlayer}{nodelayer}
		\node [style=circle, scale=3] (0) at (0, 1) {};
		\node [style=none] (5) at (0, 1) {$D(f')$};
		\node [style=none] (1) at (0, 2) {};
		\node [style=none] (2) at (-0.6, -0) {};
		\node [style=none] (3) at (0.6, -0) {};
		\node [style=none] (4) at (0.9, -0.3) {$D(Y)$};
		\node [style=none] (6) at (0.2, 2.2) {$D(X)$};
	\end{pgfonlayer}
	\begin{pgfonlayer}{edgelayer}
		\draw (1.center) to (0);
		\draw [bend right, looseness=1.] (0) to (2.center);
		\draw [bend left, looseness=1] (0) to (3.center);
\draw (-0.85, 0) -- (-0.35, 0);
\draw (-0.75, -0.1) -- (-0.45, -0.1);
\draw (-0.65, -0.2) -- (-0.55, -0.2);
	\end{pgfonlayer}
\end{tikzpicture}  
 \xmapsto{F} ~~~ \begin{tikzpicture}
	\begin{pgfonlayer}{nodelayer}
		\node [style=circle, scale=3] (0) at (0, 1) {};
		\node [style=none] (5) at (0, 1) {$D'({f'})$};
		\node [style=none] (1) at (0, 2) {};
		\node [style=none] (2) at (-0.6, -0) {};
		\node [style=none] (3) at (0.6, -0) {};
		\node [style=none] (4) at (0.9, -0.3) {$D'(Y)$};
		\node [style=none] (6) at (0.2, 2.2) {$D'(X)$};
	\end{pgfonlayer}
	\begin{pgfonlayer}{edgelayer}
		\draw (1.center) to (0);
		\draw [bend right, looseness=1.] (0) to (2.center);
		\draw [bend left, looseness=1] (0) to (3.center);
		\draw (-0.85, 0) -- (-0.35, 0);
		\draw  [bend right=90, looseness=1.25] (-0.85, 0) to (-0.35, 0);
	\end{pgfonlayer}
\end{tikzpicture} 
\end{align*}
where $F(\envmap) = \anotherenvmap{1.5}$.
\end{itemize}

This fixes the definition of $F$. To prove that $F$ is well-defined on arrows we need to show that 
$f=g \Rightarrow F(f) = F(g)$. Since, $(D, \envmap)$ has purification, let
\[
\begin{tikzpicture}
	\begin{pgfonlayer}{nodelayer}
		\node [style=circle] (0) at (0, 1) {$f$};
		\node [style=none] (1) at (0, 2) {};
		\node [style=none] (2) at (0, -0) {};
	\end{pgfonlayer}
	\begin{pgfonlayer}{edgelayer}
		\draw (1.center) to (0);
		\draw (0) to (2.center);
	\end{pgfonlayer}
\end{tikzpicture} = \begin{tikzpicture}
	\begin{pgfonlayer}{nodelayer}
		\node [style=circle, scale=3] (0) at (0, 1) {};
		\node [style=none] (5) at (0, 1) {$D(\overline{f})$};
		\node [style=none] (1) at (0, 2) {};
		\node [style=none] (2) at (-0.6, -0) {};
		\node [style=none] (3) at (0.6, -0) {};
	\end{pgfonlayer}
	\begin{pgfonlayer}{edgelayer}
		\draw (1.center) to (0);
		\draw [bend right, looseness=1.] (0) to (2.center);
		\draw [bend left, looseness=1] (0) to (3.center);
\draw (-0.85, 0) -- (-0.35, 0);
\draw (-0.75, -0.1) -- (-0.45, -0.1);
\draw (-0.65, -0.2) -- (-0.55, -0.2);
	\end{pgfonlayer}
\end{tikzpicture} ~~~~~~~ \begin{tikzpicture}
	\begin{pgfonlayer}{nodelayer}
		\node [style=circle] (0) at (0, 1) {$g$};
		\node [style=none] (1) at (0, 2) {};
		\node [style=none] (2) at (0, -0) {};
	\end{pgfonlayer}
	\begin{pgfonlayer}{edgelayer}
		\draw (1.center) to (0);
		\draw (0) to (2.center);
	\end{pgfonlayer}
\end{tikzpicture} = \begin{tikzpicture}
	\begin{pgfonlayer}{nodelayer}
		\node [style=circle, scale=3] (0) at (0, 1) {};
		\node [style=none] (5) at (0, 1) {$D(\overline{g})$};
		\node [style=none] (1) at (0, 2) {};
		\node [style=none] (2) at (-0.6, -0) {};
		\node [style=none] (3) at (0.6, -0) {};
	\end{pgfonlayer}
	\begin{pgfonlayer}{edgelayer}
		\draw (1.center) to (0);
		\draw [bend right, looseness=1.] (0) to (2.center);
		\draw [bend left, looseness=1] (0) to (3.center);
\draw (-0.85, 0) -- (-0.35, 0);
\draw (-0.75, -0.1) -- (-0.45, -0.1);
\draw (-0.65, -0.2) -- (-0.55, -0.2);
	\end{pgfonlayer}
\end{tikzpicture} 
\]
Then,
\[
 \begin{tikzpicture}
	\begin{pgfonlayer}{nodelayer}
		\node [style=circle, scale=3] (0) at (0, 1) {};
		\node [style=none] (5) at (0, 1) {$D(\overline{f})$};
		\node [style=none] (1) at (0, 2) {};
		\node [style=none] (2) at (-0.6, -0) {};
		\node [style=none] (3) at (0.6, -0) {};
	\end{pgfonlayer}
	\begin{pgfonlayer}{edgelayer}
		\draw (1.center) to (0);
		\draw [bend right, looseness=1.] (0) to (2.center);
		\draw [bend left, looseness=1] (0) to (3.center);
\draw (-0.85, 0) -- (-0.35, 0);
\draw (-0.75, -0.1) -- (-0.45, -0.1);
\draw (-0.65, -0.2) -- (-0.55, -0.2);
	\end{pgfonlayer}
\end{tikzpicture}  =  \begin{tikzpicture}
	\begin{pgfonlayer}{nodelayer}
		\node [style=circle, scale=3] (0) at (0, 1) {};
		\node [style=none] (5) at (0, 1) {$D(\overline{g})$};
		\node [style=none] (1) at (0, 2) {};
		\node [style=none] (2) at (-0.6, -0) {};
		\node [style=none] (3) at (0.6, -0) {};
	\end{pgfonlayer}
	\begin{pgfonlayer}{edgelayer}
		\draw (1.center) to (0);
		\draw [bend right, looseness=1.] (0) to (2.center);
		\draw [bend left, looseness=1] (0) to (3.center);
\draw (-0.85, 0) -- (-0.35, 0);
\draw (-0.75, -0.1) -- (-0.45, -0.1);
\draw (-0.65, -0.2) -- (-0.55, -0.2);
	\end{pgfonlayer}
\end{tikzpicture}  \Leftrightarrow  \begin{tikzpicture}
	\begin{pgfonlayer}{nodelayer}
		\node [style=circle, scale=2] (0) at (0, 1) {};
		\node [style=none] (5) at (0, 1) {$\overline{f}$};
		\node [style=none] (1) at (0, 2) {};
		\node [style=none] (2) at (-0.6, -0) {};
		\node [style=none] (3) at (0.6, -0) {};
	\end{pgfonlayer}
	\begin{pgfonlayer}{edgelayer}
		\draw (1.center) to (0);
		\draw [bend right, looseness=1.] (0) to (2.center);
		\draw [bend left, looseness=1] (0) to (3.center);
	\end{pgfonlayer}
\end{tikzpicture} = \begin{tikzpicture}
	\begin{pgfonlayer}{nodelayer}
		\node [style=circle, scale=2] (0) at (0, 1) {};
		\node [style=none] (5) at (0, 1) {$\overline{g}$};
		\node [style=none] (1) at (0, 2) {};
		\node [style=none] (2) at (-0.6, -0) {};
		\node [style=none] (3) at (0.6, -0) {};
	\end{pgfonlayer}
	\begin{pgfonlayer}{edgelayer}
		\draw (1.center) to (0);
		\draw [bend right, looseness=1.] (0) to (2.center);
		\draw [bend left, looseness=1] (0) to (3.center);
	\end{pgfonlayer}
\end{tikzpicture} \Leftrightarrow \begin{tikzpicture}
	\begin{pgfonlayer}{nodelayer}
		\node [style=circle, scale=3] (0) at (0, 1) {};
		\node [style=none] (5) at (0, 1) {$D'(\overline{f})$};
		\node [style=none] (1) at (0, 2) {};
		\node [style=none] (2) at (-0.6, -0) {};
		\node [style=none] (3) at (0.6, -0) {};
	\end{pgfonlayer}
	\begin{pgfonlayer}{edgelayer}
		\draw (1.center) to (0);
		\draw [bend right, looseness=1.] (0) to (2.center);
		\draw [bend left, looseness=1] (0) to (3.center);
		\draw (-0.85, 0) -- (-0.35, 0);
		\draw  [bend right=90, looseness=1.25] (-0.85, 0) to (-0.35, 0);
	\end{pgfonlayer}
\end{tikzpicture}  =  \begin{tikzpicture}
	\begin{pgfonlayer}{nodelayer}
		\node [style=circle, scale=3] (0) at (0, 1) {};
		\node [style=none] (5) at (0, 1) {$D'(\overline{g})$};
		\node [style=none] (1) at (0, 2) {};
		\node [style=none] (2) at (-0.6, -0) {};
		\node [style=none] (3) at (0.6, -0) {};
	\end{pgfonlayer}
	\begin{pgfonlayer}{edgelayer}
		\draw (1.center) to (0);
		\draw [bend right, looseness=1.] (0) to (2.center);
		\draw [bend left, looseness=1] (0) to (3.center);
		\draw (-0.85, 0) -- (-0.35, 0);
		\draw  [bend right=90, looseness=1.25] (-0.85, 0) to (-0.35, 0);
	\end{pgfonlayer}
\end{tikzpicture} 
\]

$F: \D \to \D'$ preserves identity: \[ F(1_A) = F(1_{D(X)}) = F(D(1_X)) = D'(1_X) = 1_{D'(X)} = 1_{F(D(X))} = 1_{F(A)} \]

$F: \D \to \D'$ preserves composition:
\[ F \left( 
\begin{tikzpicture}
	\begin{pgfonlayer}{nodelayer}
		\node [style=circle] (0) at (0, 2) {$f$};
		\node [style=circle] (1) at (0, 0.75) {$g$};
		\node [style=none] (2) at (0, -0.5) {};
		\node [style=none] (3) at (0, 3) {};
	\end{pgfonlayer}
	\begin{pgfonlayer}{edgelayer}
		\draw (3.center) to (0);
		\draw (0) to (1);
		\draw (1) to (2.center);
	\end{pgfonlayer}
\end{tikzpicture}
\right ) = 
F \left(\begin{tikzpicture} 
	\begin{pgfonlayer}{nodelayer}
		\node [style=circle, scale=2.5] (0) at (-1, 2) {};
		\node [style=circle, scale=2.5] (1) at (-1, 0.5) {};
		\node [style=none] (100) at (-1, 0.5) {$D(\overline{g})$};
		\node [style=none] (101) at (-1, 2) {$D(\overline{f})$};
		\node [style=none] (2) at (-1, -0.5) {};
		\node [style=none] (3) at (-1, 3) {};
		\node [style=none] (4) at (-2.75, -0.5) {};
		\node [style=none] (5) at (-1.75, -0.5) {};

		\node [style=none] (6) at (-3, -0.5) {};
		\node [style=none] (8) at (-2.9, -0.6) {};
		\node [style=none] (9) at (-2.8, -0.7) {};
		\node [style=none] (7) at (-2.5, -0.5) {};
		\node [style=none] (10) at (-2.6, -0.6) {};
		\node [style=none] (11) at (-2.7, -0.7) {};
				
		\node [style=none] (13) at (-1.5, -0.5) {};
		\node [style=none] (15) at (-1.6, -0.6) {};
		\node [style=none] (17) at (-1.7, -0.7) {};
		\node [style=none] (12) at (-2, -0.5) {};
		\node [style=none] (14) at (-1.9, -0.6) {};
		\node [style=none] (16) at (-1.8, -0.7) {};
	\end{pgfonlayer}
	\begin{pgfonlayer}{edgelayer}
		\draw (3.center) to (0);
		\draw (0) to (1);
		\draw (1) to (2.center);
		\draw [bend right=45, looseness=0.75] (1) to (5.center);
		\draw [bend right, looseness=1.00] (0) to (4.center);
		\draw (6.center) to (7.center);
		\draw (8.center) to (10.center);
		\draw (9.center) to (11.center);
		\draw (16.center) to (17.center);
		\draw (14.center) to (15.center);
		\draw (12.center) to (13.center);
	\end{pgfonlayer}
\end{tikzpicture} \right) \stackrel{\tiny{\bf{Env.1a} }}{=} F \left( 
\begin{tikzpicture} 
	\begin{pgfonlayer}{nodelayer}
		\node [style=circle, scale=2.8] (0) at (-1, 2) {};
		\node [style=circle, scale=2.8] (1) at (-1, 0.5) {};
		\node [style=none] (2) at (-1, -0.5) {};
		\node [style=none] (3) at (-1, 3) {};
		\node [style=none] (4) at (-2, -0.25) {};
		\node [style=none] (5) at (-1.75, -0.25) {};
		\node [style=none] (6) at (-1, 2) {$D'(\overline{f})$};
		\node [style=none] (7) at (-1, 0.5) {$D'(\overline{g})$};
		\node [style=none] (8) at (-2.25, -0.25) {};
		\node [style=none] (9) at (-1.5, -0.25) {};
		\node [style=none] (10) at (-2.1, -0.35) {};
		\node [style=none] (11) at (-1.65, -0.35) {};	
		\node [style=none] (12) at (-1.95, -0.45) {};
		\node [style=none] (13) at (-1.85, -0.45) {};		
	\end{pgfonlayer}
	\begin{pgfonlayer}{edgelayer}
		\draw (3.center) to (0);
		\draw (0) to (1);
		\draw (1) to (2.center);
		\draw [bend right=45, looseness=0.75] (1) to (5.center);
		\draw [bend right, looseness=1.00] (0) to (4.center);
		\draw (8.center) to (9.center);
		\draw (10.center) to (11.center);
		\draw (12.center) to (13.center);
	\end{pgfonlayer}
\end{tikzpicture}
 \right) :=
\begin{tikzpicture} 
	\begin{pgfonlayer}{nodelayer}
		\node [style=circle, scale=2.8] (0) at (-1, 2) {};
		\node [style=circle, scale=2.8] (1) at (-1, 0.5) {};
		\node [style=none] (2) at (-1, -0.5) {};
		\node [style=none] (3) at (-1, 3) {};
		\node [style=none] (4) at (-2, -0.25) {};
		\node [style=none] (5) at (-1.75, -0.25) {};
		\node [style=none] (6) at (-1, 2) {$D'(\overline{f})$};
		\node [style=none] (7) at (-1, 0.5) {$D'(\overline{g})$};
		\node [style=none] (8) at (-2.25, -0.25) {};
		\node [style=none] (9) at (-1.5, -0.25) {};
	\end{pgfonlayer}
	\begin{pgfonlayer}{edgelayer}
		\draw (3.center) to (0);
		\draw (0) to (1);
		\draw (1) to (2.center);
		\draw [bend right=45, looseness=0.75] (1) to (5.center);
		\draw [bend right, looseness=1.00] (0) to (4.center);
		\draw (8.center) to (9.center);		
		\draw [bend right = 90, looseness=1.00] (8.center) to (9.center);
	\end{pgfonlayer}
\end{tikzpicture}
\stackrel{\tiny{\bf{Env.1a} }}{=} \begin{tikzpicture} 
	\begin{pgfonlayer}{nodelayer}
		\node [style=circle, scale=2.8] (0) at (-1, 2) {};
		\node [style=circle, scale=2.8] (1) at (-1, 0.5) {};
		\node [style=none] (100) at (-1, 0.5) {$D'(\overline{g})$};
		\node [style=none] (101) at (-1, 2) {$D'(\overline{f})$};
		\node [style=none] (2) at (-1, -0.5) {};
		\node [style=none] (3) at (-1, 3) {};
		\node [style=none] (4) at (-2.75, -0.5) {};
		\node [style=none] (5) at (-1.75, -0.5) {};

		\node [style=none] (6) at (-3, -0.5) {};
		\node [style=none] (7) at (-2.5, -0.5) {};
				
		\node [style=none] (13) at (-1.5, -0.5) {};
		\node [style=none] (12) at (-2, -0.5) {};
	\end{pgfonlayer}
	\begin{pgfonlayer}{edgelayer}
		\draw (3.center) to (0);
		\draw (0) to (1);
		\draw (1) to (2.center);
		\draw [bend right=45, looseness=0.75] (1) to (5.center);
		\draw [bend right, looseness=1.00] (0) to (4.center);
		\draw (6.center) to (7.center);
		\draw (12.center) to (13.center);
		\draw [bend right=90, looseness=1.20] (6.center) to (7.center);
		\draw [bend right=90, looseness=1.20] (12.center) to (13.center);
	\end{pgfonlayer}
\end{tikzpicture} = \begin{tikzpicture}
	\begin{pgfonlayer}{nodelayer}
		\node [style=circle, scale=2.8] (0) at (0, 2) {};
		\node [style=none] (100) at (0, 2) {$F(f)$};
		\node [style=circle, scale=2.8] (1) at (0, 0.75) {};
		\node [style=none] (101) at (0, 0.75) {$F(g)$};
		\node [style=none] (2) at (0, -0.5) {};
		\node [style=none] (3) at (0, 3) {};
	\end{pgfonlayer}
	\begin{pgfonlayer}{edgelayer}
		\draw (3.center) to (0);
		\draw (0) to (1);
		\draw (1) to (2.center);
	\end{pgfonlayer}
\end{tikzpicture}
\]

$F: \D \to \D'$ is strict monoidal in $\otimes$:
\[ 
F \left(
\begin{tikzpicture}
	\begin{pgfonlayer}{nodelayer}
		\node [style=circle] (0) at (0, 1) {$f$};
		\node [style=none] (1) at (0, 2) {};
		\node [style=none] (2) at (0, -0) {};
	\end{pgfonlayer}
	\begin{pgfonlayer}{edgelayer}
		\draw (1.center) to (0);
		\draw (0) to (2.center);
	\end{pgfonlayer}
\end{tikzpicture} \ox \begin{tikzpicture}
	\begin{pgfonlayer}{nodelayer}
		\node [style=circle] (0) at (0, 1) {$g$};
		\node [style=none] (1) at (0, 2) {};
		\node [style=none] (2) at (0, -0) {};
	\end{pgfonlayer}
	\begin{pgfonlayer}{edgelayer}
		\draw (1.center) to (0);
		\draw (0) to (2.center);
	\end{pgfonlayer}
\end{tikzpicture} \right) = F \left( \begin{tikzpicture}
	\begin{pgfonlayer}{nodelayer}
		\node [style=circle, scale=3] (0) at (0, 1) {};
		\node [style=none] (5) at (0, 1) {$D(\overline{f})$};
		\node [style=none] (1) at (0, 2) {};
		\node [style=none] (2) at (-0.6, -0) {};
		\node [style=none] (3) at (0.6, -0) {};
	\end{pgfonlayer}
	\begin{pgfonlayer}{edgelayer}
		\draw (1.center) to (0);
		\draw [bend right, looseness=1.] (0) to (2.center);
		\draw [bend left, looseness=1] (0) to (3.center);
\draw (-0.85, 0) -- (-0.35, 0);
\draw (-0.75, -0.1) -- (-0.45, -0.1);
\draw (-0.65, -0.2) -- (-0.55, -0.2);
	\end{pgfonlayer}
\end{tikzpicture} \ox \begin{tikzpicture}
	\begin{pgfonlayer}{nodelayer}
		\node [style=circle, scale=3] (0) at (0, 1) {};
		\node [style=none] (5) at (0, 1) {$D(\overline{g})$};
		\node [style=none] (1) at (0, 2) {};
		\node [style=none] (2) at (-0.6, -0) {};
		\node [style=none] (3) at (0.6, -0) {};
	\end{pgfonlayer}
	\begin{pgfonlayer}{edgelayer}
		\draw (1.center) to (0);
		\draw [bend right, looseness=1.] (0) to (2.center);
		\draw [bend left, looseness=1] (0) to (3.center);
\draw (-0.85, 0) -- (-0.35, 0);
\draw (-0.75, -0.1) -- (-0.45, -0.1);
\draw (-0.65, -0.2) -- (-0.55, -0.2);
	\end{pgfonlayer}
\end{tikzpicture} \right) = F\left( 
\begin{tikzpicture}
	\begin{pgfonlayer}{nodelayer}
		\node [style=circle, scale=2.8] (0) at (-2, 1) {};
		\node [style=none] (1) at (-2, 2.25) {};
		\node [style=circle, scale=2.8] (2) at (-0.5, 1) {};
		\node [style=none] (3) at (-0.5, 2.25) {};
		\node [style=none] (4) at (-2.75, -0.75) {};
		\node [style=none] (5) at (-2.25, -0.75) {};
		\node [style=none] (6) at (-1, -1) {};
		\node [style=none] (7) at (0, -1) {};
		\node [style=none] (8) at (-2, 1) {$D(\overline{f})$};
		\node [style=none] (9) at (-0.5, 1) {$D(\overline{g})$};
		\node [style=none] (10) at (-3, -0.75) {};
		\node [style=none] (11) at (-2.85, -0.85) {};
		\node [style=none] (12) at (-2.7, -0.95) {};
		\node [style=none] (13) at (-2, -0.75) {};
		\node [style=none] (14) at (-2.15, -0.85) {};
		\node [style=none] (15) at (-2.3, -0.95) {};
	\end{pgfonlayer}
	\begin{pgfonlayer}{edgelayer}
		\draw [bend left, looseness=1.25] (0) to (6.center);
		\draw [bend left=15, looseness=1.00] (2) to (7.center);
		\draw [bend right=15, looseness=1.00] (2) to (5.center);
		\draw [bend right=15, looseness=1.00] (0) to (4.center);
		\draw (0) to (1.center);
		\draw (3.center) to (2);
		\draw (10.center) to (13.center);
		\draw (11.center) to (14.center);
		\draw (12.center) to (15.center);
	\end{pgfonlayer}
\end{tikzpicture}
\right)=  
\begin{tikzpicture}
	\begin{pgfonlayer}{nodelayer}
		\node [style=circle, scale=3] (0) at (-2, 1) {};
		\node [style=none] (1) at (-2, 2.25) {};
		\node [style=circle, scale=3] (2) at (-0.5, 1) {};
		\node [style=none] (3) at (-0.5, 2.25) {};
		\node [style=none] (4) at (-2.75, -0.75) {};
		\node [style=none] (5) at (-2.25, -0.75) {};
		\node [style=none] (6) at (-1, -1) {};
		\node [style=none] (7) at (0, -1) {};
		\node [style=none] (8) at (-2, 1) {$D'(\overline{f})$};
		\node [style=none] (9) at (-0.5, 1) {$D'(\overline{g})$};
		\node [style=none] (10) at (-3, -0.75) {};
		\node [style=none] (13) at (-2, -0.75) {};
	\end{pgfonlayer}
	\begin{pgfonlayer}{edgelayer}
		\draw [bend left, looseness=1.25] (0) to (6.center);
		\draw [bend left=15, looseness=1.00] (2) to (7.center);
		\draw [bend right=15, looseness=1.00] (2) to (5.center);
		\draw [bend right=15, looseness=1.00] (0) to (4.center);
		\draw (0) to (1.center);
		\draw (3.center) to (2);
		\draw (10.center) to (13.center);
		\draw [bend right=90, looseness=1.25] (10.center) to (13.center);
	\end{pgfonlayer}
\end{tikzpicture}
= \begin{tikzpicture}
	\begin{pgfonlayer}{nodelayer}
		\node [style=circle, scale=3] (0) at (0, 1) {};
		\node [style=none] (5) at (0, 1) {$D'(\overline{f})$};
		\node [style=none] (1) at (0, 2) {};
		\node [style=none] (2) at (-0.6, -0) {};
		\node [style=none] (3) at (0.6, -0) {};
	\end{pgfonlayer}
	\begin{pgfonlayer}{edgelayer}
		\draw (1.center) to (0);
		\draw [bend right, looseness=1.] (0) to (2.center);
		\draw [bend left, looseness=1] (0) to (3.center);
		\draw (-0.85, 0) -- (-0.35, 0);
		\draw  [bend right=90, looseness=1.25] (-0.85, 0) to (-0.35, 0);
	\end{pgfonlayer}
\end{tikzpicture}  \ox \begin{tikzpicture}
	\begin{pgfonlayer}{nodelayer}
		\node [style=circle, scale=3] (0) at (0, 1) {};
		\node [style=none] (5) at (0, 1) {$D'(\overline{g})$};
		\node [style=none] (1) at (0, 2) {};
		\node [style=none] (2) at (-0.6, -0) {};
		\node [style=none] (3) at (0.6, -0) {};
	\end{pgfonlayer}
	\begin{pgfonlayer}{edgelayer}
		\draw (1.center) to (0);
		\draw [bend right, looseness=1.] (0) to (2.center);
		\draw [bend left, looseness=1] (0) to (3.center);
		\draw (-0.85, 0) -- (-0.35, 0);
		\draw  [bend right=90, looseness=1.25] (-0.85, 0) to (-0.35, 0);
	\end{pgfonlayer}
\end{tikzpicture} 
\]
\[ \text{and, } F((u_\ox^L)_A) = F((u_\ox^L)_{D(X)}) = F(D((u_\ox^L)_X)) = D'((u_\ox^L)_X) = (u_\ox^L)_{D'(X)} = (u_\ox^L)_{F(D(X))} = (u_\ox^L)_{F(A)}\]

Similarly, it can be proved that $F$ is strict comonoidal in $\oa$. 

Define $F_\ox = F_\oa := F$ and linear strengths to be identity maps. Thus, $F$ is a unique strict Frobenius functor. $F$ is an isomix functor because D and D' preserve the mix map $\m$ on the nose.
\end{proof}

\begin{corollary}
Suppose $D: \C \to \D$ with $\envmap$ is an environment structure with purification for $M:\U \to \C$. Then, 
$\D \simeq \CP^\infty(M: \U \to \C)$.
\end{corollary}
\begin{proof}
By Lemma \ref{Lemma: Env initial}, $D: \C \to \D$ with $\envmap$ is initial in $\mathsf{Env}(M: \U \to \C)$. By Lemma 
\ref{Lemma: Env example},  $F: \C \to \CP^\infty(M: \U \to \C)$ with $[(u_\oa^R)^{-1},U)]$ is an environment structure 
for $M: \U \to \C$ which has purification, hence it is also an initial object in $\mathsf{Env}(M: \U \to \C)$. Since, initial 
objects of a category are isomorphic, there exists a strict isomix functor $\D \xrightarrow{F} \CP^\infty(M: \U \to \C)$ 
that is full and faithful.
\end{proof}

%% file: chapter-compaction.tex

\chapter{Linear monoids, comonoids, and bialgebras}
\label{Chapter: compaction}

In a $\dagger$-monoidal category, (strong) complementary observables are axiomatized 
by two special commutative $\dagger$-Frobenius algebras interacting bialgebraically 
to produce two Hopf algebras, see Section \ref{Sec: strong comp}. In an LDC, Frobenius algebras are generalized to 
linear monoids. However, bialgebraic interaction between two linear monoids is prohibited 
due to the directionality of the linear distributor. This raises the following question:   
{\em how one can model complementary observables in a $\dagger$-isomix setting?} 

In this chapter, we develop new structures called linear comonoids and linear bialgebras 
which provide the basis to describe complementary systems in $\dagger$-isomix categories. 

\section{Linear monoids}
\label{Section: linear monoids}

For a Frobenius algebra, the monoid and the comonoid occur on the same object, and the 
object is self-dual. In an LDC, the presence of distinct tensor products - the tensor and 
the par - allows a `Frobenius interaction' between a monoid and its dual comonoid, however  
the monoid and the comonoid now occur on distinct objects. A linear monoid, loosely defined, 
is a $\ox$-monoid on an object $A$ and $A$ has a dual object. 

\subsection{Duals}
\label{Sec: duals}

Adding the sequent rules for negation (see Figure \ref{Fig: negation rules}) to LDCs results in the notion of duals. 
 We briefly discussed duals \cite{CKS00} in LDCs in 
 Section \ref{Sec: *-autonomous}. We delve deeper into the notion of duals, subsequently introduce $\dagger$-duals 
in this section. 

In a compact closed category or a $*$-autonomous category there is a canonical dual associated with each object.  
In an LDC if an object has a dual, it must be specified.  Of course, any two duals of an object are isomorphic. 
A {\bf self-duality} is a dual, $(\eta, \epsilon): A \dashvv B$ in which $A$ is isomorphic to $B$ (or $A=B$). 
A dual, $(\eta, \epsilon): A \dashvv B$ such that $A$ is both left and right dual of $B$, 
is called a {\bf cyclic dual}. 
In a symmetric LDC, every dual $(\eta, \epsilon): A \dashvv B$ gives another 
dual $(\eta c_\oa, c_\ox \epsilon): B \dashvv A$, which is obtained by twisting 
the wires using the symmetry map. Thus, in a symmetric LDC, every dual is a cyclic dual. 

\begin{lemma}
In an LDC, if $(a,b): A \dashvv B$, and $(c,d): C \dashvv D$ are duals then $(x,y): A \ox C\dashvv D \oa B$, and 
$(x',y') : A \oa C \dashvv D \ox B$ are duals, where 
\[ x := \begin{tikzpicture}
	\begin{pgfonlayer}{nodelayer}
		\node [style=none] (0) at (-1.75, 1.75) {};
		\node [style=none] (1) at (-0.75, 1.75) {};
		\node [style=none] (2) at (-1.75, 2) {};
		\node [style=none] (3) at (-0.75, 2) {};
		\node [style=none] (6) at (-0.25, 2.25) {};
		\node [style=none] (7) at (-2.25, 2.25) {};
		\node [style=none] (8) at (-0.25, 1.75) {};
		\node [style=none] (9) at (-2.25, 1.75) {};
		\node [style=none] (12) at (-1.25, 2.75) {$c$};
		\node [style=none] (13) at (-2.75, 1) {$A \ox C$};
		\node [style=none] (16) at (0.25, 1) {$D \oa B$};
		\node [style=ox] (17) at (-2, 1.25) {};
		\node [style=none] (18) at (-2, 0.75) {};
		\node [style=none] (19) at (-1.75, 1.75) {};
		\node [style=none] (20) at (-2.25, 1.75) {};
		\node [style=oa] (21) at (-0.5, 1.25) {};
		\node [style=none] (22) at (-0.5, 0.75) {};
		\node [style=none] (23) at (-0.25, 1.75) {};
		\node [style=none] (24) at (-0.75, 1.75) {};
		\node [style=none] (25) at (-1.25, 3.25) {$a$};
	\end{pgfonlayer}
	\begin{pgfonlayer}{edgelayer}
		\draw (0.center) to (2.center);
		\draw [bend left=90, looseness=1.75] (2.center) to (3.center);
		\draw (3.center) to (1.center);
		\draw (6.center) to (8.center);
		\draw (18.center) to (17);
		\draw [in=-90, out=150, looseness=1.25] (17) to (20.center);
		\draw [in=-90, out=30, looseness=1.25] (17) to (19.center);
		\draw (22.center) to (21);
		\draw [in=-90, out=150, looseness=1.25] (21) to (24.center);
		\draw [in=-90, out=30, looseness=1.25] (21) to (23.center);
		\draw (7.center) to (20.center);
		\draw [bend right=90, looseness=1.25] (6.center) to (7.center);
	\end{pgfonlayer}
\end{tikzpicture}
~~~~ y := \begin{tikzpicture}
	\begin{pgfonlayer}{nodelayer}
		\node [style=none] (0) at (-0.75, 2.25) {};
		\node [style=none] (1) at (-1.75, 2.25) {};
		\node [style=none] (2) at (-0.75, 2) {};
		\node [style=none] (3) at (-1.75, 2) {};
		\node [style=none] (6) at (-2.25, 1.75) {};
		\node [style=none] (7) at (-0.25, 1.75) {};
		\node [style=none] (8) at (-2.25, 2.25) {};
		\node [style=none] (9) at (-0.25, 2.25) {};
		\node [style=none] (12) at (-1.25, 1.25) {$b$};
		\node [style=none] (13) at (0.25, 3) {$A \ox C$};
		\node [style=none] (16) at (-2.75, 3) {$D \oa B$};
		\node [style=ox] (17) at (-0.5, 2.75) {};
		\node [style=none] (18) at (-0.5, 3.25) {};
		\node [style=none] (19) at (-0.75, 2.25) {};
		\node [style=none] (20) at (-0.25, 2.25) {};
		\node [style=oa] (21) at (-2, 2.75) {};
		\node [style=none] (22) at (-2, 3.25) {};
		\node [style=none] (23) at (-2.25, 2.25) {};
		\node [style=none] (24) at (-1.75, 2.25) {};
		\node [style=none] (25) at (-1.25, 0.75) {$d$};
	\end{pgfonlayer}
	\begin{pgfonlayer}{edgelayer}
		\draw (0.center) to (2.center);
		\draw [bend left=90, looseness=1.75] (2.center) to (3.center);
		\draw (3.center) to (1.center);
		\draw (6.center) to (8.center);
		\draw (18.center) to (17);
		\draw [in=90, out=-30, looseness=1.25] (17) to (20.center);
		\draw [in=90, out=-150, looseness=1.25] (17) to (19.center);
		\draw (22.center) to (21);
		\draw [in=90, out=-30, looseness=1.25] (21) to (24.center);
		\draw [in=90, out=-150, looseness=1.25] (21) to (23.center);
		\draw (7.center) to (20.center);
		\draw [bend right=90, looseness=1.25] (6.center) to (7.center);
	\end{pgfonlayer}
\end{tikzpicture}
 ~~~~ x' := \begin{tikzpicture}
	\begin{pgfonlayer}{nodelayer}
		\node [style=none] (0) at (-1.75, 1.75) {};
		\node [style=none] (1) at (-0.75, 1.75) {};
		\node [style=none] (2) at (-1.75, 2) {};
		\node [style=none] (3) at (-0.75, 2) {};
		\node [style=none] (6) at (-0.25, 2.25) {};
		\node [style=none] (7) at (-2.25, 2.25) {};
		\node [style=none] (8) at (-0.25, 1.75) {};
		\node [style=none] (9) at (-2.25, 1.75) {};
		\node [style=none] (12) at (-1.25, 2.75) {$c$};
		\node [style=none] (13) at (-2.75, 1) {$A \oa C$};
		\node [style=none] (16) at (0.25, 1) {$D \ox B$};
		\node [style=oa] (17) at (-2, 1.25) {};
		\node [style=none] (18) at (-2, 0.75) {};
		\node [style=none] (19) at (-1.75, 1.75) {};
		\node [style=none] (20) at (-2.25, 1.75) {};
		\node [style=ox] (21) at (-0.5, 1.25) {};
		\node [style=none] (22) at (-0.5, 0.75) {};
		\node [style=none] (23) at (-0.25, 1.75) {};
		\node [style=none] (24) at (-0.75, 1.75) {};
		\node [style=none] (25) at (-1.25, 3.25) {$a$};
	\end{pgfonlayer}
	\begin{pgfonlayer}{edgelayer}
		\draw (0.center) to (2.center);
		\draw [bend left=90, looseness=1.75] (2.center) to (3.center);
		\draw (3.center) to (1.center);
		\draw (6.center) to (8.center);
		\draw (18.center) to (17);
		\draw [in=-90, out=150, looseness=1.25] (17) to (20.center);
		\draw [in=-90, out=30, looseness=1.25] (17) to (19.center);
		\draw (22.center) to (21);
		\draw [in=-90, out=150, looseness=1.25] (21) to (24.center);
		\draw [in=-90, out=30, looseness=1.25] (21) to (23.center);
		\draw (7.center) to (20.center);
		\draw [bend right=90, looseness=1.25] (6.center) to (7.center);
	\end{pgfonlayer}
\end{tikzpicture}
 ~~~~ y' := \begin{tikzpicture}
	\begin{pgfonlayer}{nodelayer}
		\node [style=none] (0) at (-0.75, 2.25) {};
		\node [style=none] (1) at (-1.75, 2.25) {};
		\node [style=none] (2) at (-0.75, 2) {};
		\node [style=none] (3) at (-1.75, 2) {};
		\node [style=none] (6) at (-2.25, 1.75) {};
		\node [style=none] (7) at (-0.25, 1.75) {};
		\node [style=none] (8) at (-2.25, 2.25) {};
		\node [style=none] (9) at (-0.25, 2.25) {};
		\node [style=none] (12) at (-1.25, 1.25) {$b$};
		\node [style=none] (13) at (0.25, 3) {$A \oa C$};
		\node [style=none] (16) at (-2.75, 3) {$D \ox B$};
		\node [style=oa] (17) at (-0.5, 2.75) {};
		\node [style=none] (18) at (-0.5, 3.25) {};
		\node [style=none] (19) at (-0.75, 2.25) {};
		\node [style=none] (20) at (-0.25, 2.25) {};
		\node [style=ox] (21) at (-2, 2.75) {};
		\node [style=none] (22) at (-2, 3.25) {};
		\node [style=none] (23) at (-2.25, 2.25) {};
		\node [style=none] (24) at (-1.75, 2.25) {};
		\node [style=none] (25) at (-1.25, 0.75) {d};
	\end{pgfonlayer}
	\begin{pgfonlayer}{edgelayer}
		\draw (0.center) to (2.center);
		\draw [bend left=90, looseness=1.75] (2.center) to (3.center);
		\draw (3.center) to (1.center);
		\draw (6.center) to (8.center);
		\draw (18.center) to (17);
		\draw [in=90, out=-30, looseness=1.25] (17) to (20.center);
		\draw [in=90, out=-150, looseness=1.25] (17) to (19.center);
		\draw (22.center) to (21);
		\draw [in=90, out=-30, looseness=1.25] (21) to (24.center);
		\draw [in=90, out=-150, looseness=1.25] (21) to (23.center);
		\draw (7.center) to (20.center);
		\draw [bend right=90, looseness=1.25] (6.center) to (7.center);
	\end{pgfonlayer}
\end{tikzpicture}\]
\end{lemma}

\begin{definition} A {\bf morphism of duals}, $(f,g): ((\eta,  \epsilon): A \dashvv B) \to ((\tau, \gamma): A' \dashvv B')$, is given by a pair of maps
$f: A \to A'$, and $g: B' \to B$ such that:  
\begin{center}
(a) ~~~  \begin{tikzpicture}
	\begin{pgfonlayer}{nodelayer}
		\node [style=none] (0) at (-1, 2) {};
		\node [style=none] (1) at (0.5, 3) {};
		\node [style=none] (2) at (-1, 3) {};
		\node [style=none] (3) at (-0.25, 4) {$\tau$};
		\node [style=none] (4) at (0.5, 2) {};
		\node [style=none] (5) at (-1.25, 2.25) {$A'$};
		\node [style=none] (6) at (1, 3.5) {$B'$};
		\node [style=circle, scale=1.5] (7) at (0.5, 2.75) {};
		\node [style=none] (8) at (0.5, 2.75) {$g$};
		\node [style=none] (9) at (1, 2.25) {$B$};
	\end{pgfonlayer}
	\begin{pgfonlayer}{edgelayer}
		\draw (4.center) to (7);
		\draw (1.center) to (7);
		\draw (2.center) to (0.center);
		\draw [bend left=90, looseness=1.75] (2.center) to (1.center);
	\end{pgfonlayer}
\end{tikzpicture}
= \begin{tikzpicture}
	\begin{pgfonlayer}{nodelayer}
		\node [style=none] (0) at (0.5, 2) {};
		\node [style=none] (1) at (-1, 3) {};
		\node [style=none] (2) at (0.5, 3) {};
		\node [style=none] (3) at (-0.25, 4) {$\eta$};
		\node [style=none] (4) at (-1, 2) {};
		\node [style=none] (5) at (0.75, 2.25) {$B$};
		\node [style=none] (6) at (-1.5, 3.5) {$A$};
		\node [style=circle, scale=1.5] (7) at (-1, 2.75) {};
		\node [style=none] (8) at (-1, 2.75) {$f$};
		\node [style=none] (9) at (-1.5, 2.25) {$A'$};
	\end{pgfonlayer}
	\begin{pgfonlayer}{edgelayer}
		\draw (4.center) to (7);
		\draw (1.center) to (7);
		\draw (2.center) to (0.center);
		\draw [bend right=90, looseness=1.75] (2.center) to (1.center);
	\end{pgfonlayer}
\end{tikzpicture}
  ~~~~~~(b)~~~ \begin{tikzpicture}
	\begin{pgfonlayer}{nodelayer}
		\node [style=none] (0) at (0.5, 3.75) {};
		\node [style=none] (1) at (-1, 2.75) {};
		\node [style=none] (2) at (0.5, 2.75) {};
		\node [style=none] (3) at (-0.25, 1.75) {$\gamma$};
		\node [style=none] (4) at (-1, 3.75) {};
		\node [style=none] (5) at (0.75, 3.5) {$B'$};
		\node [style=none] (6) at (-1.5, 2.25) {$A'$};
		\node [style=circle, scale=1.5] (7) at (-1, 3) {};
		\node [style=none] (8) at (-1, 3) {$f$};
		\node [style=none] (9) at (-1.5, 3.5) {$A$};
	\end{pgfonlayer}
	\begin{pgfonlayer}{edgelayer}
		\draw (4.center) to (7);
		\draw (1.center) to (7);
		\draw (2.center) to (0.center);
		\draw [bend right=90, looseness=1.75] (1.center) to (2.center);
	\end{pgfonlayer}
\end{tikzpicture}
 =  \begin{tikzpicture}
	\begin{pgfonlayer}{nodelayer}
		\node [style=none] (0) at (-1, 3.75) {};
		\node [style=none] (1) at (0.5, 2.75) {};
		\node [style=none] (2) at (-1, 2.75) {};
		\node [style=none] (3) at (-0.25, 1.75) {$\epsilon$};
		\node [style=none] (4) at (0.5, 3.75) {};
		\node [style=none] (5) at (-1.25, 3.5) {$A$};
		\node [style=none] (6) at (1, 2.25) {$B$};
		\node [style=circle, scale=1.5] (7) at (0.5, 3) {};
		\node [style=none] (8) at (0.5, 3) {$g$};
		\node [style=none] (9) at (1, 3.5) {$B'$};
	\end{pgfonlayer}
	\begin{pgfonlayer}{edgelayer}
		\draw (4.center) to (7);
		\draw (1.center) to (7);
		\draw (2.center) to (0.center);
		\draw [bend right=90, looseness=1.75] (2.center) to (1.center);
	\end{pgfonlayer}
\end{tikzpicture}
\end{center}
\end{definition}

A morphism of duals is determined by either of the maps, as $f$ is 
dual to $g$. They are, thus, Australian mates, see \cite{CKS00}. 

In a $\dagger$-LDC, if $A$ is dual to $B$, then $B^\dagger$ is dual to $A^\dagger$:
\begin{lemma} \label{daggering-a-dual}
	In a $\dag$-LDC, if $(\eta, \epsilon): A \dashvv B$ is a dual  the $(\epsilon\dagger, \eta\dagger): B^\dagger \dashvv A^\dagger$ is a dual where:
	\begin{align*}
	\epsilon\dagger &:= \top \to^{\lambda_\top} \bot^\dagger \to^{\epsilon^\dagger} 
	(B \ox A)^\dagger \to^{\lambda_\oa^{-1}} B^\dagger \oa A^\dagger \\
	\eta\dagger &:= A^\dagger \ox B^\dagger \to^{\lambda_\ox} (A \oa B)^\dagger 
	\to^{\eta^\dagger} \top^\dagger \to^{\lambda_\bot^{-1}} \bot
	\end{align*}
\end{lemma}

This leads to the notion of $\dagger$-duals:

\begin{definition}
	\label{defn: right dagger dual}
	In a $\dagger$-LDC, a {\bf right $\dagger$-dual}, $A \dagdual A^\dagger$ is a dual $(\eta, \epsilon):A \dashvv A^\dagger$ such that 
	
	\medskip

	{ \centering $(\iota_A, 1_{A^\dagger}): (\eta, \epsilon):A \dashvv A^\dagger \to (\eta^\dagger,\epsilon^\dagger):A^{\dagger\dagger} \dashvv A^\dagger$ \par}
       
	is an isomorphism of duals.  A {\bf left $\dagger$-dual}, $A^\dag \dagdual A$ is a dual 
	$(\eta, \epsilon): A^\dag \dashvv A$ such that 

	\medskip

	{ \centering $(1_{A^\dagger}, \iota_A): \eta^\dagger,\epsilon^\dagger):A^{\dagger\dagger} \dashvv A^\dagger \to (\eta, \epsilon): A^\dagger \dashvv A$ \par}
	
	\medskip

	is an isomorphism of duals. A {\bf self $\dagger$-dual} is a right (or left) $\dagger$-dual 
	with an isomorphism $\alpha: A \to A^\dagger$ such that $\alpha \alpha^{-1 \dagger}= \iota$.
\end{definition}
$(\iota_A, 1_{A^\dagger})$ being an isomorphism of the duals means that the following equations hold: 
\begin{equation}
	\label{eqn: right dagger dual} 
		(a)~~~ \begin{tikzpicture}
			\begin{pgfonlayer}{nodelayer}
				\node [style=none] (0) at (0.5, 0.25) {};
				\node [style=none] (1) at (-1, 1.5) {};
				\node [style=none] (2) at (0.5, 1.5) {};
				\node [style=none] (3) at (-0.25, 2.5) {$\eta$};
				\node [style=none] (4) at (-1, 0.25) {};
				\node [style=none] (5) at (-1.5, 0.25) {$A^{\dagger  \dagger}$};
				\node [style=none] (6) at (1, 0.25) {$A^{\dagger}$};
				\node [style=circle, scale=1.5] (7) at (-1, 1) {};
				\node [style=none] (8) at (-1, 1) {$\iota$};
			\end{pgfonlayer}
			\begin{pgfonlayer}{edgelayer}
				\draw (0.center) to (2.center);
				\draw (1.center) to (7);
				\draw (4.center) to (7);
				\draw [bend left=90, looseness=1.75] (1.center) to (2.center);
			\end{pgfonlayer}
		\end{tikzpicture}  =
		\begin{tikzpicture}
			\begin{pgfonlayer}{nodelayer}
				\node [style=none] (0) at (0.25, 3.25) {};
				\node [style=none] (1) at (-1, 5.5) {};
				\node [style=none] (2) at (0.4999997, 5.5) {};
				\node [style=none] (3) at (-0.25, 4.5) {$\epsilon$};
				\node [style=none] (4) at (-0.75, 3.25) {};
				\node [style=none] (5) at (-0.65, 5.25) {$A^\dagger$};
				\node [style=none] (6) at (-1.25, 3.25) {$A^{\dagger \dagger}$};
				\node [style=none] (7) at (0.5, 3.25) {$A^\dagger$};
				\node [style=none] (8) at (0.75, 4) {};
				\node [style=none] (9) at (-1.25, 4) {};
				\node [style=none] (10) at (-1.25, 5.5) {};
				\node [style=none] (11) at (0.75, 5.5) {};
				\node [style=none] (12) at (0.25, 4) {};
				\node [style=none] (13) at (-0.75, 4) {};
				\node [style=none] (14) at (0.25, 5.25) {$A$};
			\end{pgfonlayer}
			\begin{pgfonlayer}{edgelayer}
				\draw (8.center) to (9.center);
				\draw (10.center) to (9.center);
				\draw (10.center) to (11.center);
				\draw (11.center) to (8.center);
				\draw (13.center) to (4.center);
				\draw (0.center) to (12.center);
				\draw [bend right=90, looseness=1.75] (1.center) to (2.center);
			\end{pgfonlayer}
		\end{tikzpicture}
		~~~~~~~ \text{(equivalently)} ~~~~~~~ (b)~~~
		\begin{tikzpicture}
			\begin{pgfonlayer}{nodelayer}
				\node [style=none] (0) at (0.25, 2.25) {};
				\node [style=none] (1) at (-1.25, 0.75) {};
				\node [style=none] (2) at (0.25, 0.75) {};
				\node [style=none] (3) at (-0.5, -0.25) {$\epsilon$};
				\node [style=none] (4) at (-1.25, 2.25) {};
				\node [style=none] (5) at (-1.25, 1.5) {};
				\node [style=none] (6) at (-1.5, 2.25) {$A^{\dagger}$};
				\node [style=none] (7) at (0.5, 2.25) {$A$};
			\end{pgfonlayer}
			\begin{pgfonlayer}{edgelayer}
				\draw (1.center) to (4.center);
				\draw (2.center) to (0.center);
				\draw [bend right=90, looseness=1.75] (1.center) to (2.center);
			\end{pgfonlayer}
		\end{tikzpicture}
		= \begin{tikzpicture}
			\begin{pgfonlayer}{nodelayer}
				\node [style=none] (0) at (-0.75, 2.25) {};
				\node [style=none] (1) at (0.5, -0) {};
				\node [style=none] (2) at (-1, -0) {};
				\node [style=none] (3) at (-0.25, 1) {$\eta$};
				\node [style=none] (4) at (0.2499997, 2.25) {};
				\node [style=none] (5) at (0.15, 0.25) {$A^\dagger$};
				\node [style=none] (6) at (0.75, 2.25) {$A$};
				\node [style=none] (7) at (-1.25, 2.25) {$A^{\dagger}$};
				\node [style=none] (8) at (-1.25, 1.25) {};
				\node [style=none] (9) at (0.75, 1.25) {};
				\node [style=none] (10) at (0.75, 0) {};
				\node [style=none] (11) at (-1.25, 0) {};
				\node [style=none] (12) at (-0.75, 1.25) {};
				\node [style=none] (13) at (0.2499997, 1.25) {};
				\node [style=none] (14) at (-0.75, 0.25) {$A$};
				\node [style=circle, scale=1.5] (15) at (0.25, 1.75) {};
				\node [style=none] (16) at (0.25, 1.75) {$\iota$};
			\end{pgfonlayer}
			\begin{pgfonlayer}{edgelayer}
				\draw (8.center) to (9.center);
				\draw (10.center) to (9.center);
				\draw (10.center) to (11.center);
				\draw (11.center) to (8.center);
				\draw (12.center) to (0.center);
				\draw (13.center) to (15);
				\draw (4.center) to (15);
				\draw [bend left=90, looseness=1.75] (2.center) to (1.center);
			\end{pgfonlayer}
		\end{tikzpicture} \end{equation}

		In the sequent calclulus of $\dagger$-linear logic, a premise is a right $\dagger$-dual 
		if it satisfies the rules in Figure \ref{Fig: rules dag duals}. As shown in Section \ref{Sec: ldc-intro}, 
		one can derive $\eta: \top \to A \oa A^\dag$, and $A^\dag \ox A \to \bot$ for right $\dagger$-duals using these rules. 
		Note that the left $\dagger$-duals can be derived using the same rules, the $(\dagger)$ rule and the $(\iota)$ rules 
		in Section \ref{Sec: dagger sequent rules}.
		\begin{figure}[h]
			\centering
		    \AxiomC{$\Gamma, A \vdash \Delta$}
			\LeftLabel{($\dagger$-dual:R)}
            \UnaryInfC{ $\Gamma, \vdash A^\dag, \Delta$} 
			 \DisplayProof
			\qquad
			\AxiomC{$\Gamma, A^\dag \vdash \Delta$}
			\LeftLabel{($\dagger$-dual:L)}
			\UnaryInfC{$\Gamma \vdash A, \Delta$}
			\DisplayProof
			\caption{Sequent rules for right $\dagger$-duals}
			\label{Fig: rules dag duals}
	    \end{figure}
		\FloatBarrier

    \begin{lemma}
	\label{Lemma: dagger equiv}
	In a $\dagger$-LDC, the following are equivalent for a dual $(\eta, \epsilon): A \dashvv B$:
	\begin{enumerate}[(i)]
	\item There exists an isomorphism from the dual to its right $\dagger$-dual (illustrated in diagram $(a)$ below).
		\item There exists an isomorphism from the dual to its left $\dagger$-dual (illustrated in diagram $(b)$ below).
	\item There exists a morphism $(p,q)$ between the dual and its dagger such that $(p,q)$ satisfies any two of the following three conditions: 
	(1) $pq^\dagger = \iota_A$, (2) $p^\dagger q = \iota_B^{-1}$, and (3) $p$ or $q$ are isomorphisms.
	\end{enumerate}
	\[ (a) ~~ 
 \]
\end{lemma}
\begin{proof}

\begin{description}
	
\item[$(i) \Rightarrow (ii)$] 
Given that $(\eta, \epsilon): A \dashvv B \to^{(1,q)} (\tau, \gamma): A \dashvv A^\dagger$ 
is an isomorphism of duals. Then, define $p: A \to B^\dagger$ as follows:

{ \centering $ p :=%
 $ \par}

We note that $\iota_B p^\dag q  = 1_B$:

$ %
 $

$(1)$ and $(4)$ are because $(1,q)$ is a morphism of duals, $(2)$ is because of the right $\dagger$-duality of $A$, 
and $(3)$ is because $\iota_B^{-1 \dagger} = \iota_{B^\dag}$. 

Define the left $\dagger$-dual $(\tau', \gamma'): B^\dag \dashvv B$ as follows:
\[  \tau'= %
 \]
It is clear that $\tau'$ and $\gamma'$ satisfy snake diagrams. $(\tau', \gamma'): B^\dag \dashvv B$ 
can be proven to be a left $\dagger$-dual using the equation $\iota_B p^\dagger q = 1_B$. 
It is straightforward from the definition of $\tau'$ and $\gamma'$ that $(p,1)$ is a morphism of duals.

\item[$(ii) \Rightarrow (iii)$] Given that $(p,1)$ is a morphism of duals as illustrated by diagram $(b)$. 
Define, 

{\centering  $ q: A^\dagger \to B := %
 $ \par}

$(p,q): ((\eta, \epsilon): A \dashvv B) \to ((\eta\dag, \epsilon\dag): B^\dag \dashvv A^\dag)$ is a morphism of duals because:

{ \centering $%
 $ \par}

Since $p$ and $q$ are isomorphisms, it remains to prove that either $pq^\dagger = i_A$ 
or $p^\dagger q = \iota_B^{-1}$. Using the defintion of $q$, we note that, 

{ \centering $ p: A \to B^\dag =%
 $ \par}

From the proof of $(i) \Rightarrow (ii)$ we know that $p^\dagger q = \iota_B^{-1}$. Thus, statement $(iii)$ holds.

\item[$(iii) \Rightarrow (i)$]  Assume that statement $(iii)$ holds. Define a dual $(\tau, \gamma): A \dashvv B$ as follows: 
\[ %
  \]
$(\tau, \gamma): A \dashvv A^\dag$ is a right $\dagger$-dual because $pq^\dagger = \iota_A$. 
It is straightforward that $(1_A, q)$ is a morphism of duals. 
\end{description}
\end{proof}

A {\bf $\dagger$-dual} is a dual that satisfies any of the three equivalent statements of the previous Lemma. 
Notice,  however, in statement (i) of the previous lemma, that the morphism to a given right $\dagger$-dual is 
uniquely determined by demanding the map is $(1_A,q)$ as $1_A$ completely determines $q$. 
If {\em there exists} a right $\dagger$-dual for $A$, then  any mere linear dual of $A$ is a $\dagger$-dual. 
 Thus, $\dagger$-duality really is a property of the object $A$ rather than of the duality. Hence, in the rest of the 
 thesis, we will use right and left $\dagger$-duals synonymously with $\dagger$-duals.  

 In a dagger monoidal setting, conventionally, we have a dagger dual \cite{Sel07} when the symmetry map acts 
 an isomorphism between the dual and its dagger:
 \[  \begin{tikzpicture}[scale=1.5]
    \begin{pgfonlayer}{nodelayer}
        \node [style=none] (0) at (-0.75, 2) {$A$};
        \node [style=none] (1) at (2, 2) {$B$};
        \node [style=none] (2) at (-1.25, -0.5) {$B^\dagger = B$};
        \node [style=none] (3) at (2.5, -0.5) {$A^\dagger = A$};
        \node [style=none] (4) at (-0.25, 2) {};
        \node [style=none] (5) at (1.5, 2.25) {};
        \node [style=none] (6) at (1.5, 1.75) {};
        \node [style=none] (7) at (1.35, 2.25) {};
        \node [style=none] (8) at (1.35, 1.75) {};
        \node [style=none] (9) at (1.5, 2) {};
        \node [style=none] (10) at (1.5, -0.75) {};
        \node [style=none] (11) at (1.5, -0.5) {};
        \node [style=none] (12) at (1.35, -0.25) {};
        \node [style=none] (13) at (1.5, -0.25) {};
        \node [style=none] (14) at (-0.25, -0.5) {};
        \node [style=none] (15) at (1.35, -0.75) {};
        \node [style=none] (16) at (-0.75, 1.75) {};
        \node [style=none] (17) at (2, -0.25) {};
        \node [style=none] (18) at (-0.75, -0.25) {};
        \node [style=none] (19) at (2, 1.75) {};
        \node [style=none] (20) at (0.5, 2.25) {$(\eta, \epsilon)$};
        \node [style=none] (21) at (0.5, -0.25) {$(\epsilon^\dagger, \eta^\dagger)$};
        \node [style=none, scale=1.45] (22) at (0.5, -1.15) {$\epsilon^\dagger = \eta c_\ox; \eta^\dagger = c_\ox \epsilon$};
    \end{pgfonlayer}
    \begin{pgfonlayer}{edgelayer}
        \draw (7.center) to (8.center);
        \draw (5.center) to (6.center);
        \draw (4.center) to (9.center);
        \draw (12.center) to (15.center);
        \draw (13.center) to (10.center);
        \draw (14.center) to (11.center);
        \draw [in=90, out=-90, looseness=0.75] (19.center) to (18.center);
        \draw [in=90, out=-90, looseness=0.75] (16.center) to (17.center);
    \end{pgfonlayer}
\end{tikzpicture} \]

 Using the symmetry map in this manner requires the dagger functor to be stationary on objects ($A = A^\dag$) 
 which does not hold in a $\dagger$-LDC in general ($A \neq A^\dag$). The following diagram represents a $dagger$-dual in a $\dagger$-LDC:
 \[ \begin{tikzpicture}[scale=1.5]
	\begin{pgfonlayer}{nodelayer}
		\node [style=none] (0) at (-0.75, 2.5) {$A$};
		\node [style=none] (1) at (2, 2.5) {$B$};
		\node [style=none] (2) at (-0.75, 0.5) {$B^\dagger$};
		\node [style=none] (3) at (2, 0.5) {$A^\dagger$};
		\node [style=none] (4) at (-0.25, 2.5) {};
		\node [style=none] (5) at (1.35, 2.75) {};
		\node [style=none] (6) at (1.35, 2.25) {};
		\node [style=none] (7) at (1.35, 2.5) {};
		\node [style=none] (8) at (1.5, 0.25) {};
		\node [style=none] (9) at (1.5, 0.5) {};
		\node [style=none] (10) at (1.35, 0.75) {};
		\node [style=none] (11) at (1.5, 0.75) {};
		\node [style=none] (12) at (-0.25, 0.5) {};
		\node [style=none] (13) at (1.35, 0.25) {};
		\node [style=none] (14) at (-0.75, 2.25) {};
		\node [style=none] (15) at (-0.75, 0.75) {};
		\node [style=none] (16) at (2, 0.75) {};
		\node [style=none] (17) at (2, 2.25) {};
		\node [style=none] (18) at (0.5, 2.75) {$(\eta, \epsilon)$};
		\node [style=none] (19) at (0.5, 0.75) {$(\epsilon^\dagger, \eta^\dagger)$};
		\node [style=none, scale=1.25] (20) at (0.5, -0.15) {$\epsilon^\dagger = \eta (q^{-1} \ox p) \lambda_\oa$};
		\node [style=none] (21) at (1.5, 2.75) {};
		\node [style=none] (22) at (1.5, 2.25) {};
		\node [style=none] (23) at (-1, 1.5) {$p$};
		\node [style=none] (24) at (2.25, 1.5) {$q$};
	\end{pgfonlayer}
	\begin{pgfonlayer}{edgelayer}
		\draw (5.center) to (6.center);
		\draw (4.center) to (7.center);
		\draw (10.center) to (13.center);
		\draw (11.center) to (8.center);
		\draw (12.center) to (9.center);
		\draw [<-, line width=0.15mm] (17.center) to (16.center);
		\draw [->, line width=0.15mm] (14.center) to (15.center);
		\draw (21.center) to (22.center);
	\end{pgfonlayer}
\end{tikzpicture}  \]

 In the category of complex finiteness matrices, ${\sf FMat(\C)}$,
which is a $\dagger$-$*$-isomix category, every dual is a $\dagger$-dual. 
 The category of finite-dimensional complex matrices, ${\sf Mat(\C)}$ which is a 
 $\dagger$-compact closed category is isomorphic to the core of ${\sf FMat(\C)}$, 
 see Section \ref{Sec: core of FMat}.  Thus, every dual in ${\sf Mat}(\C)$ is also $\dagger$-dual. 
 In fact every dual in ${\sf FMat(\C)}$  is both a $\dagger$-dual and a unitary dual 
 ($\eta^\dag = \epsilon c_\ox$). These examples are discussed in detail in Section \ref{Sec: duals examples}. 

Observe $\dagger$-linear functors preserve $\dagger$-duals:

 \begin{lemma}
	 \label{Lemma: daglin dagdual}
	 $\dagger$-linear functors preserve $\dagger$-duals.
 \end{lemma}
 \begin{proof}
 Let $\X$, and $\Y$ be $\dagger$-LDCs. Let $(F_\ox, F_\oa): \X \to \Y$ be a $\dagger$-linear functor. 
 Then, for each object A, we have maps:
	 \[
	 \s_A: F_\ox(A^\dagger) \to^{\simeq} (F_\oa(A))^\dagger ~~~~~~~ \t_A: (F_\ox(A))^\dagger \to^{\simeq} F_\oa(A^\dagger)
	 \]
 such that $(\s, \t): (F_\ox((-)^\dagger), (F_\oa(-)^\dagger)) \Rightarrow ( (F_\oa(-))^\dagger, F_\ox((-)^\dagger))$ is a linear 
 natural isomorphism.  We must prove that if $A \dagmonw A^\dagger$ is a right $\dagger$-linear 
 dual, then $F_\ox(A)$ has a right $\dagger$-dual. 
 
 Suppose $A \dagdual A^\dagger$ is a right $\dagger$-dual. Since, linear functors preserve linear duals, 
 $F_\ox(A) \dashvv F_\oa(A^\dagger)$ is also a linear dual with the cap and the cup given by:
 \[ \top \to^{\eta'} F_\ox(A) \oa F_\ox(A)^\dagger := \begin{tikzpicture}
	 \begin{pgfonlayer}{nodelayer}
		 \node [style=none] (0) at (0.5, 3) {};
		 \node [style=circle, scale=2] (1) at (0.5, 2) {};
		 \node [style=none] (2) at (0.5, 1) {};
		 \node [style=none] (3) at (0.5, 2) {$\t^{-1}$};
		 \node [style=none] (4) at (-0.5, 1) {};
		 \node [style=none] (5) at (-1, 3) {};
		 \node [style=none] (6) at (1, 3) {};
		 \node [style=none] (7) at (-1, 4.5) {};
		 \node [style=none] (8) at (1, 4.5) {};
		 \node [style=none] (9) at (0.75, 4.25) {$F_\oa$};
		 \node [style=none] (10) at (-1.25, 1.5) {$F_\ox(A)$};
		 \node [style=none] (11) at (1.5, 2.5) {$F_\oa(A^\dagger)$};
		 \node [style=none] (12) at (1.5, 1.25) {$(F_\ox(A))^\dagger$};
		 \node [style=none] (13) at (-0.75, 3.5) {$A$};
		 \node [style=none] (14) at (0.75, 3.5) {$A^\dagger$};
		 \node [style=none] (15) at (-0.5, 3) {};
		 \node [style=none] (16) at (-0.75, 3) {};
		 \node [style=none] (17) at (-0.25, 3) {};
	 \end{pgfonlayer}
	 \begin{pgfonlayer}{edgelayer}
		 \draw (0.center) to (1);
		 \draw (1) to (2.center);
		 \draw (5.center) to (6.center);
		 \draw (6.center) to (8.center);
		 \draw (8.center) to (7.center);
		 \draw (7.center) to (5.center);
		 \draw (4.center) to (15.center);
		 \draw [bend left=90, looseness=3.50] (15.center) to (0.center);
		 \draw [bend left=60, looseness=1.25] (16.center) to (17.center);
	 \end{pgfonlayer}
 \end{tikzpicture} ~~~~~~~~ F_\ox(A)^\dagger \ox F_\ox (A) \to^{\epsilon'} \bot := \begin{tikzpicture}
	 \begin{pgfonlayer}{nodelayer}
		 \node [style=none] (0) at (0.75, 2.5) {};
		 \node [style=none] (1) at (-0.25, 2.5) {};
		 \node [style=circle, scale=1.5] (2) at (-0.25, 3.5) {};
		 \node [style=none] (3) at (-0.25, 4.5) {};
		 \node [style=none] (4) at (-0.25, 3.5) {$\t$};
		 \node [style=none] (5) at (0.75, 4.5) {};
		 \node [style=none] (6) at (1.25, 2.5) {};
		 \node [style=none] (7) at (-0.75, 2.5) {};
		 \node [style=none] (8) at (1.25, 1) {};
		 \node [style=none] (9) at (-0.75, 1) {};
		 \node [style=none] (10) at (-0.5, 1.25) {$F_\ox$};
		 \node [style=none] (11) at (-1.25, 4) {$(F_\ox(A))^\dagger$};
		 \node [style=none] (12) at (-1, 3) {$F_\oa(A^\dagger)$};
		 \node [style=none] (13) at (1.5, 4) {$F_\ox(A)$};
		 \node [style=none] (14) at (-0.5, 2) {$A^\dagger$};
		 \node [style=none] (15) at (1, 2) {$A$};
		 \node [style=none] (16) at (-0.5, 2.5) {};
		 \node [style=none] (17) at (0, 2.5) {};
	 \end{pgfonlayer}
	 \begin{pgfonlayer}{edgelayer}
		 \draw (5.center) to (0.center);
		 \draw (1.center) to (2);
		 \draw (2) to (3.center);
		 \draw (6.center) to (7.center);
		 \draw (7.center) to (9.center);
		 \draw (9.center) to (8.center);
		 \draw (8.center) to (6.center);
		 \draw [bend left=90, looseness=3.50] (0.center) to (1.center);
		 \draw [bend right=75, looseness=1.25] (16.center) to (17.center);
	 \end{pgfonlayer}
 \end{tikzpicture}\]
 
 To show that it is a right $\dagger$-dual:
 \[ \begin{tikzpicture}
	 \begin{pgfonlayer}{nodelayer}
		 \node [style=none] (0) at (-0.5, 4.5) {};
		 \node [style=none] (1) at (0.5, 4.5) {};
		 \node [style=circle, scale=2] (2) at (0.5, 2.25) {};
		 \node [style=none] (3) at (0.5, 1) {};
		 \node [style=none] (4) at (0.5, 2.25) {$\t^{-1}$};
		 \node [style=none] (5) at (-0.5, 1) {};
		 \node [style=none] (6) at (-1, 4.5) {};
		 \node [style=none] (7) at (1, 4.5) {};
		 \node [style=none] (8) at (-1, 6) {};
		 \node [style=none] (9) at (1, 6) {};
		 \node [style=none] (10) at (0.75, 5.75) {$F_\oa$};
		 \node [style=none] (11) at (-1.25, 1.5) {$F_\ox(A)$};
		 \node [style=none] (12) at (1.25, 3.75) {$F_\oa(A^\dagger)$};
		 \node [style=none] (13) at (1.5, 1.25) {$(F_\ox(A))^\dagger$};
		 \node [style=none] (14) at (-0.75, 5) {$A$};
		 \node [style=none] (15) at (0.75, 5) {$A^\dagger$};
		 \node [style=circle, scale=1.5] (16) at (-0.5, 3.25) {};
		 \node [style=none] (17) at (-0.5, 3.25) {$\iota$};
		 \node [style=none] (18) at (-0.75, 4.5) {};
		 \node [style=none] (19) at (-0.25, 4.5) {};
	 \end{pgfonlayer}
	 \begin{pgfonlayer}{edgelayer}
		 \draw (1.center) to (2);
		 \draw (2) to (3.center);
		 \draw (6.center) to (7.center);
		 \draw (7.center) to (9.center);
		 \draw (9.center) to (8.center);
		 \draw (8.center) to (6.center);
		 \draw [bend left=90, looseness=3.50] (0.center) to (1.center);
		 \draw (5.center) to (16);
		 \draw (0.center) to (16);
		 \draw [bend left=90, looseness=1.25] (18.center) to (19.center);
	 \end{pgfonlayer}
 \end{tikzpicture} = \begin{tikzpicture}
	 \begin{pgfonlayer}{nodelayer}
		 \node [style=none] (0) at (-0.5, 3.75) {};
		 \node [style=none] (1) at (0.5, 3.75) {};
		 \node [style=circle, scale=2] (2) at (0.5, 2.25) {};
		 \node [style=none] (3) at (0.5, 0.25) {};
		 \node [style=none] (4) at (0.5, 2.25) {$\t^{-1}$};
		 \node [style=none] (5) at (-0.5, 0.25) {};
		 \node [style=none] (6) at (-1, 3.75) {};
		 \node [style=none] (7) at (1, 3.75) {};
		 \node [style=none] (8) at (-1, 6) {};
		 \node [style=none] (9) at (1, 6) {};
		 \node [style=none] (10) at (0.75, 5.75) {$F_\oa$};
		 \node [style=none] (11) at (-1.25, 0.5) {$F_\ox(A)^{\dagger \dagger}$};
		 \node [style=none] (12) at (1.5, 2.75) {$F_\oa(A^\dagger)$};
		 \node [style=none] (13) at (1.5, 0.75) {$(F_\ox(A))^\dagger$};
		 \node [style=none] (14) at (-0.75, 5.5) {$A$};
		 \node [style=circle, scale=1.5] (15) at (-0.5, 2.5) {};
		 \node [style=none] (16) at (-0.5, 2.5) {$\s$};
		 \node [style=circle, scale=1.5] (17) at (-0.5, 1.25) {};
		 \node [style=none] (18) at (-0.5, 1.25) {$\t$};
		 \node [style=circle, scale=1.5] (19) at (-0.5, 4.75) {};
		 \node [style=none] (20) at (0.5, 5) {};
		 \node [style=none] (21) at (-0.5, 4.75) {$\iota$};
		 \node [style=none] (22) at (-1.25, 3.25) {$F_\ox(A^{\dagger \dagger})$};
		 \node [style=none] (23) at (-1.75, 2) {$F_\oa(A^\dagger)^\dagger$};
		 \node [style=none] (24) at (-1, 0.75) {};
		 \node [style=none] (25) at (0, 0.75) {};
		 \node [style=none] (26) at (0, 1.75) {};
		 \node [style=none] (27) at (-1, 1.75) {};
		 \node [style=none] (28) at (-0.5, 1.75) {};
		 \node [style=none] (29) at (-0.5, 2.25) {};
		 \node [style=none] (30) at (-0.5, 0.75) {};
		 \node [style=none] (31) at (-0.5, 0.25) {};
		 \node [style=none] (32) at (-0.75, 1.75) {};
		 \node [style=none] (33) at (-0.25, 0.75) {};
		 \node [style=none] (34) at (-0.75, 3.75) {};
		 \node [style=none] (35) at (-0.25, 3.75) {};
	 \end{pgfonlayer}
	 \begin{pgfonlayer}{edgelayer}
		 \draw (1.center) to (2);
		 \draw (2) to (3.center);
		 \draw (6.center) to (7.center);
		 \draw (7.center) to (9.center);
		 \draw (9.center) to (8.center);
		 \draw (8.center) to (6.center);
		 \draw (0.center) to (15);
		 \draw (0.center) to (19);
		 \draw [bend right=90, looseness=2.75] (20.center) to (19);
		 \draw (20.center) to (1.center);
		 \draw (24.center) to (25.center);
		 \draw (25.center) to (26.center);
		 \draw (26.center) to (27.center);
		 \draw (27.center) to (24.center);
		 \draw (5.center) to (30.center);
		 \draw (28.center) to (29.center);
		 \draw (32.center) to (17);
		 \draw (33.center) to (17);
		 \draw [bend left=90, looseness=1.50] (34.center) to (35.center);
	 \end{pgfonlayer}
 \end{tikzpicture} = \begin{tikzpicture}
	\begin{pgfonlayer}{nodelayer}
		\node [style=none] (0) at (-0.5, 3.75) {};
		\node [style=none] (1) at (0.5, 3.75) {};
		\node [style=circle, scale=2] (2) at (0.5, 2.25) {};
		\node [style=none] (3) at (0.5, 0.25) {};
		\node [style=none] (4) at (0.5, 2.25) {$\t^{-1}$};
		\node [style=none] (5) at (-0.5, 0.25) {};
		\node [style=none] (6) at (-1.25, 3.75) {};
		\node [style=none] (7) at (1, 3.75) {};
		\node [style=none] (8) at (-1.25, 6) {};
		\node [style=none] (9) at (1, 6) {};
		\node [style=none] (10) at (0.75, 5.75) {$F_\oa$};
		\node [style=none] (11) at (-1.25, 0.5) {$F_\ox(A)^{\dagger \dagger}$};
		\node [style=none] (12) at (1.5, 2.75) {$F_\oa(A^\dagger)$};
		\node [style=none] (13) at (1.5, 0.75) {$(F_\ox(A))^\dagger$};
		\node [style=circle, scale=1.5] (14) at (-0.5, 2.5) {};
		\node [style=none] (15) at (-0.5, 2.5) {$\s$};
		\node [style=circle, scale=1.5] (16) at (-0.5, 1.25) {};
		\node [style=none] (17) at (-0.5, 1.25) {$\t$};
		\node [style=none] (18) at (-1.25, 3.25) {$F_\ox(A^{\dagger \dagger})$};
		\node [style=none] (19) at (-1.75, 2) {$F_\oa(A^\dagger)^\dagger$};
		\node [style=none] (20) at (-1, 0.75) {};
		\node [style=none] (21) at (0, 0.75) {};
		\node [style=none] (22) at (0, 1.75) {};
		\node [style=none] (23) at (-1, 1.75) {};
		\node [style=none] (24) at (-0.5, 1.75) {};
		\node [style=none] (25) at (-0.5, 2.25) {};
		\node [style=none] (26) at (-0.5, 0.75) {};
		\node [style=none] (27) at (-0.5, 0.25) {};
		\node [style=none] (28) at (-0.75, 1.75) {};
		\node [style=none] (29) at (-0.25, 0.75) {};
		\node [style=none] (30) at (-0.5, 5.5) {};
		\node [style=none] (31) at (0.5, 5.5) {};
		\node [style=none] (32) at (-1, 5.5) {};
		\node [style=none] (33) at (-1, 4.5) {};
		\node [style=none] (34) at (0.75, 4.5) {};
		\node [style=none] (35) at (0.75, 5.5) {};
		\node [style=none] (36) at (0.5, 4.5) {};
		\node [style=none] (37) at (-0.5, 4.5) {};
		\node [style=none] (38) at (-0.25, 3.75) {};
		\node [style=none] (39) at (-0.75, 3.75) {};
	\end{pgfonlayer}
	\begin{pgfonlayer}{edgelayer}
		\draw (1.center) to (2);
		\draw (2) to (3.center);
		\draw (6.center) to (7.center);
		\draw (7.center) to (9.center);
		\draw (9.center) to (8.center);
		\draw (8.center) to (6.center);
		\draw (0.center) to (14);
		\draw (20.center) to (21.center);
		\draw (21.center) to (22.center);
		\draw (22.center) to (23.center);
		\draw (23.center) to (20.center);
		\draw (5.center) to (26.center);
		\draw (24.center) to (25.center);
		\draw (28.center) to (16);
		\draw (29.center) to (16);
		\draw [bend right=90, looseness=2.50] (30.center) to (31.center);
		\draw (33.center) to (34.center);
		\draw (34.center) to (35.center);
		\draw (35.center) to (32.center);
		\draw (32.center) to (33.center);
		\draw (37.center) to (0.center);
		\draw (36.center) to (1.center);
		\draw [bend right=90, looseness=1.25] (38.center) to (39.center);
	\end{pgfonlayer}
\end{tikzpicture} =\begin{tikzpicture}
	 \begin{pgfonlayer}{nodelayer}
		 \node [style=none] (0) at (0.5, 4.25) {};
		 \node [style=circle, scale=2] (1) at (0.5, 2.25) {};
		 \node [style=none] (2) at (0.5, 0.25) {};
		 \node [style=none] (3) at (0.5, 2.25) {$\t^{-1}$};
		 \node [style=none] (4) at (-0.5, 0.25) {};
		 \node [style=none] (5) at (-1.25, 4.25) {};
		 \node [style=none] (6) at (1, 4.25) {};
		 \node [style=none] (7) at (-1.25, 6) {};
		 \node [style=none] (8) at (1, 6) {};
		 \node [style=none] (9) at (0.5, 4.75) {$F_\ox$};
		 \node [style=none] (10) at (-1.25, 0.5) {$F_\ox(A)^{\dagger \dagger}$};
		 \node [style=none] (11) at (1.5, 2.75) {$F_\oa(A^\dagger)$};
		 \node [style=none] (12) at (1.5, 0.75) {$(F_\ox(A))^\dagger$};
		 \node [style=circle, scale=1.5] (13) at (-0.5, 2.5) {};
		 \node [style=none] (14) at (-0.5, 2.5) {$\s$};
		 \node [style=circle, scale=1.5] (15) at (-0.5, 1.25) {};
		 \node [style=none] (16) at (-0.5, 1.25) {$\t$};
		 \node [style=none] (17) at (-1.75, 2) {$F_\oa(A^\dagger)^\dagger$};
		 \node [style=none] (18) at (-1, 0.75) {};
		 \node [style=none] (19) at (0, 0.75) {};
		 \node [style=none] (20) at (0, 1.75) {};
		 \node [style=none] (21) at (-1, 1.75) {};
		 \node [style=none] (22) at (-0.5, 1.75) {};
		 \node [style=none] (23) at (-0.5, 2.25) {};
		 \node [style=none] (24) at (-0.5, 0.75) {};
		 \node [style=none] (25) at (-0.5, 0.25) {};
		 \node [style=none] (26) at (-0.75, 1.75) {};
		 \node [style=none] (27) at (-0.25, 0.75) {};
		 \node [style=none] (28) at (-0.5, 5.5) {};
		 \node [style=none] (29) at (0.5, 5.5) {};
		 \node [style=none] (30) at (-1, 5.5) {};
		 \node [style=none] (31) at (-1, 4.5) {};
		 \node [style=none] (32) at (0.75, 4.5) {};
		 \node [style=none] (33) at (0.75, 5.5) {};
		 \node [style=none] (34) at (0.5, 4.5) {};
		 \node [style=none] (35) at (-0.5, 4.5) {};
		 \node [style=none] (36) at (0.5, 6) {};
		 \node [style=none] (37) at (-0.5, 6) {};
		 \node [style=none] (38) at (0.5, 3) {};
		 \node [style=circle, scale=2] (39) at (-0.5, 3.5) {};
		 \node [style=circle, scale=1.5] (40) at (0.5, 3.25) {};
		 \node [style=none] (41) at (-0.5, 3.5) {$\s^{-1}$};
		 \node [style=none] (42) at (0.5, 3.25) {$\t$};
		 \node [style=none] (43) at (-0.5, 4.25) {};
		 \node [style=none] (44) at (1.25, 3.75) {$(F_\ox(A))^\dagger$};
		 \node [style=none] (45) at (-1.75, 3.75) {$F_\oa(A^\dagger)^\dagger$};
		 \node [style=none] (46) at (-0.75, 5.5) {};
		 \node [style=none] (47) at (-0.25, 5.5) {};
	 \end{pgfonlayer}
	 \begin{pgfonlayer}{edgelayer}
		 \draw (1) to (2.center);
		 \draw (5.center) to (6.center);
		 \draw (6.center) to (8.center);
		 \draw (8.center) to (7.center);
		 \draw (7.center) to (5.center);
		 \draw (18.center) to (19.center);
		 \draw (19.center) to (20.center);
		 \draw (20.center) to (21.center);
		 \draw (21.center) to (18.center);
		 \draw (4.center) to (24.center);
		 \draw (22.center) to (23.center);
		 \draw (26.center) to (15);
		 \draw (27.center) to (15);
		 \draw [bend right=90, looseness=2.50] (28.center) to (29.center);
		 \draw (31.center) to (32.center);
		 \draw (32.center) to (33.center);
		 \draw (33.center) to (30.center);
		 \draw (30.center) to (31.center);
		 \draw (36.center) to (29.center);
		 \draw (1) to (38.center);
		 \draw (0.center) to (40);
		 \draw (43.center) to (39);
		 \draw (13) to (39);
		 \draw (37.center) to (28.center);
		 \draw [bend right=75, looseness=1.25] (46.center) to (47.center);
	 \end{pgfonlayer}
 \end{tikzpicture} = \begin{tikzpicture}
	 \begin{pgfonlayer}{nodelayer}
		 \node [style=none] (0) at (0.5, 2.75) {};
		 \node [style=none] (1) at (0.5, 0.25) {};
		 \node [style=none] (2) at (-0.75, 0.25) {};
		 \node [style=none] (3) at (-1.25, 2.75) {};
		 \node [style=none] (4) at (1, 2.75) {};
		 \node [style=none] (5) at (-1.25, 6) {};
		 \node [style=none] (6) at (1, 6) {};
		 \node [style=none] (7) at (0.5, 3.25) {$F_\ox$};
		 \node [style=none] (8) at (-1.25, 0.5) {$F_\ox(A)^{\dagger \dagger}$};
		 \node [style=none] (9) at (1.5, 0.75) {$(F_\ox(A))^\dagger$};
		 \node [style=circle, scale=2] (10) at (-0.5, 5.25) {};
		 \node [style=none] (11) at (-0.5, 5.25) {$\t$};
		 \node [style=none] (12) at (-0.75, 0.25) {};
		 \node [style=none] (13) at (-0.5, 6) {};
		 \node [style=none] (14) at (-0.5, 4.25) {};
		 \node [style=none] (15) at (0.5, 4.25) {};
		 \node [style=none] (16) at (-1, 4.25) {};
		 \node [style=none] (17) at (-1, 3) {};
		 \node [style=none] (18) at (0.75, 3) {};
		 \node [style=none] (19) at (0.75, 4.25) {};
		 \node [style=none] (20) at (0.5, 3) {};
		 \node [style=none] (21) at (-0.5, 3) {};
		 \node [style=none] (22) at (0.5, 6) {};
		 \node [style=none] (23) at (-0.5, 6) {};
		 \node [style=none] (24) at (-0.75, 2.75) {};
		 \node [style=none] (25) at (-0.75, 3.5) {$A^\dagger$};
		 \node [style=none] (26) at (-0.75, 4.25) {};
		 \node [style=none] (27) at (-0.25, 4.25) {};
	 \end{pgfonlayer}
	 \begin{pgfonlayer}{edgelayer}
		 \draw (3.center) to (4.center);
		 \draw (4.center) to (6.center);
		 \draw (6.center) to (5.center);
		 \draw (5.center) to (3.center);
		 \draw (13.center) to (10);
		 \draw [bend right=90, looseness=2.50] (14.center) to (15.center);
		 \draw (17.center) to (18.center);
		 \draw (18.center) to (19.center);
		 \draw (19.center) to (16.center);
		 \draw (16.center) to (17.center);
		 \draw (22.center) to (15.center);
		 \draw (2.center) to (24.center);
		 \draw (1.center) to (0.center);
		 \draw (14.center) to (10);
		 \draw [bend right=75, looseness=1.25] (26.center) to (27.center);
	 \end{pgfonlayer}
 \end{tikzpicture}  \]
 \end{proof}   

In a symmetric $\dagger$-LDC, if $(a,b): A \dagdual A^\dag$ and $(c,d): B \dagdual B^\dag$ are $\dagger$-duals, 
then $(A \ox B) \dashvv (A^\dag \oa B^\dag)$, and $(A \oa B) \dashvv (A^\dag \ox B^\dag)$ are 
$\dagger$-duals.

A {\bf morphism of $\dagger$-duals} is a morphism of duals in one of the following forms:
\begin{itemize}  
\item For a right $\dagger$-dual it consists of a map and its dagger:  
\[ (f, f^\dagger): ((\eta, \epsilon): A \dagdual A^\dagger) \to ((\eta', \epsilon'): B \dagdual B^\dagger) \]
\item For a left $\dagger$-dual it is of the form:
\[ (f^\dag, f): ((\eta, \epsilon): A^\dagger \dagdual A) \to ((\eta', \epsilon'): B^\dagger \dagdual B) \]
\item  A pair $(f, g): ((\eta, \epsilon): A \dagdual B) \to ((\eta', \epsilon'): C \dagdual D)$ 
for general dagger duals such that one of the following (equivalent) diagrams commute:
\begin{equation} 
	(a)~~~\xymatrix{ 
	C^\dag \ar[r]^{f^\dag} \ar[d]_{q'} & A^\dag \ar[d]^{q} \\ 
	D \ar[r]_g & B } ~~~~~ (or)~~~~~ (b) ~~~ \xymatrix{ 
		A \ar[r]^{f} \ar[d]_{p} & C \ar[d]^{p'} \\ 
		B^\dag \ar[r]_{g^\dag} & D^\dag }
\end{equation}
where $(p,q)$ is the isomorphism between the dual $(\eta, \epsilon): A \dashvv B$ and its dagger  
$(\epsilon\dagger, \eta\dagger): B^\dag  \dashvv A^\dag$, and  
$(p',q')$ is the isomorphism between the dual $(\eta', \epsilon'): C \dashvv D$ and its dagger 
$(\epsilon'\dagger, \eta'\dagger): D^\dag  \dashvv C^\dag$, see Statement $(iii)$ of 
Lemma \ref{Lemma: dagger equiv}.
\end{itemize}

The commuting diagrams $(a)$ and $(b)$ imply that the morphism $(f,g)$ behaves coherently with 
the isomorphisms which determine $\dagger$-dual, see the below diagram:
\[ \begin{matrix}
	\xymatrixcolsep{2cm}
	\xymatrix{
	A \ar@{-||}[r]^{(\eta,  \epsilon)} \ar[d]_{p} \ar@/_{2pc}/[ddd]_{f} &  B \ar@{<-}[d]^{q}    \\
	B^\dag   \ar@{-||}[r]_{(\epsilon\dagger, \eta\dagger)} \ar@/_{2pc}/[ddd]_{g^\dag} & A^\dagger \\ 
	\\ 
	C \ar@{-||}[r]^{(\eta',  \epsilon')} \ar[d]_{p'}  &  D \ar@{<-}[d]^{q'} \ar@/_{2pc}/[uuu]_{g}   \\
	D^\dag   \ar@{-||}[r]_{(\epsilon'\dagger, \eta'\dagger)}  & C^\dagger \ar@/_{2pc}/[uuu]_{f^\dag}
	} \end{matrix} 
\]



\subsection{Linear monoids}
\label{Sec: linear monoids}
In a linear setting with two distinct tensor products, Frobenius Algebras are generalized by linear monoids 
\cite{Egg10, CKS00}.  The simplest way to describe a linear monoid is as a  $\ox$-monoid  on an object together 
with a dual for that object.  The similarity to Frobenius algebras becomes more apparent when one regards 
a linear monoid as a $\ox$-monoid and a $\oa$-comonoid with actions and coactions.  

\begin{definition}
	A {\bf linear monoid}, $A \linmonw B$, in an LDC consists of a monoid $(A,e: \top \to A,m: A \ox A \to A)$, 
	a left dual  $(\eta_L, \epsilon_L): A \dashvv B$, and  a right dual $(\eta_R, \epsilon_R): B \dashvv A$ such that: 
	\begin{equation}
		\label{eqn: linear monoid} 
		\begin{tikzpicture}
			\begin{pgfonlayer}{nodelayer}
				\node [style=none] (12) at (1.75, 0.25) {$B$};
				\node [style=none] (13) at (0.25, 0.25) {$B$};
				\node [style=circle] (17) at (1, 1.25) {};
				\node [style=none] (18) at (1, 2.5) {};
				\node [style=none] (19) at (0.5, 0) {};
				\node [style=none] (20) at (1.5, 0) {};
				\node [style=none] (21) at (1.25, 2.25) {$B$};
			\end{pgfonlayer}
			\begin{pgfonlayer}{edgelayer}
				\draw [in=-150, out=90, looseness=1.25] (19.center) to (17);
				\draw [in=90, out=-30, looseness=1.25] (17) to (20.center);
				\draw (17) to (18.center);
			\end{pgfonlayer}
		\end{tikzpicture} := 		
		 \begin{tikzpicture}
		\begin{pgfonlayer}{nodelayer}
			\node [style=circle] (0) at (-2.75, 0.75) {};
			\node [style=none] (1) at (-3.25, 1.25) {};
			\node [style=none] (2) at (-2.75, 0.5) {};
			\node [style=none] (3) at (-2.25, 1.25) {};
			\node [style=none] (4) at (-1.5, 1.25) {};
			\node [style=none] (5) at (-1.5, -0) {};
			\node [style=none] (6) at (-3.25, 1.5) {};
			\node [style=none] (7) at (-1, 1.5) {};
			\node [style=none] (8) at (-1, -0) {};
			\node [style=none] (9) at (-3.75, 0.5) {};
			\node [style=none] (10) at (-3.75, 2.25) {};
			\node [style=none] (11) at (-4, 2) {$B$};
			\node [style=none] (12) at (-1.75, 0.25) {$B$};
			\node [style=none] (13) at (-0.75, 0.25) {$B$};
			\node [style=none] (14) at (-2.25, 2.5) {$\eta_L$};
			\node [style=none] (15) at (-2, 1.75) {$\eta_L$};
			\node [style=none] (16) at (-3.25, -0.25) {$\epsilon_L$};
		\end{pgfonlayer}
		\begin{pgfonlayer}{edgelayer}
			\draw [in=150, out=-90, looseness=1.00] (1.center) to (0);
			\draw [in=-90, out=30, looseness=1.00] (0) to (3.center);
			\draw (0) to (2.center);
			\draw [bend left=90, looseness=1.50] (2.center) to (9.center);
			\draw (9.center) to (10.center);
			\draw (6.center) to (1.center);
			\draw [bend left=90, looseness=1.25] (6.center) to (7.center);
			\draw (7.center) to (8.center);
			\draw (4.center) to (5.center);
			\draw [bend right=90, looseness=1.25] (4.center) to (3.center);
		\end{pgfonlayer}
	\end{tikzpicture} = \begin{tikzpicture}
		\begin{pgfonlayer}{nodelayer}
			\node [style=circle] (0) at (-2, 0.75) {};
			\node [style=none] (1) at (-1.5, 1.25) {};
			\node [style=none] (2) at (-2, 0.5) {};
			\node [style=none] (3) at (-2.5, 1.25) {};
			\node [style=none] (4) at (-3.25, 1.25) {};
			\node [style=none] (5) at (-3.25, -0) {};
			\node [style=none] (6) at (-1.5, 1.5) {};
			\node [style=none] (7) at (-3.75, 1.5) {};
			\node [style=none] (8) at (-3.75, -0) {};
			\node [style=none] (9) at (-1, 0.5) {};
			\node [style=none] (10) at (-1, 2.25) {};
			\node [style=none] (11) at (-0.75, 2) {$B$};
			\node [style=none] (12) at (-3, 0.25) {$B$};
			\node [style=none] (13) at (-4, 0.25) {$B$};
			\node [style=none] (14) at (-2.5, 2.5) {$\eta_R$};
			\node [style=none] (15) at (-2.75, 1.75) {$\eta_R$};
			\node [style=none] (16) at (-1.5, -0.25) {$\epsilon_R$};
		\end{pgfonlayer}
		\begin{pgfonlayer}{edgelayer}
			\draw [in=30, out=-90, looseness=1.00] (1.center) to (0);
			\draw [in=-90, out=150, looseness=1.00] (0) to (3.center);
			\draw (0) to (2.center);
			\draw [bend right=90, looseness=1.50] (2.center) to (9.center);
			\draw (9.center) to (10.center);
			\draw (6.center) to (1.center);
			\draw [bend right=90, looseness=1.25] (6.center) to (7.center);
			\draw (7.center) to (8.center);
			\draw (4.center) to (5.center);
			\draw [bend left=90, looseness=1.25] (4.center) to (3.center);
		\end{pgfonlayer}
	\end{tikzpicture} 
	~~~~~~~~~~~~ 
	\begin{tikzpicture}
		\begin{pgfonlayer}{nodelayer}
			\node [style=circle] (17) at (1, 0.25) {};
			\node [style=none] (18) at (1, 2.5) {};
			\node [style=none] (21) at (1.25, 2.25) {$B$};
		\end{pgfonlayer}
		\begin{pgfonlayer}{edgelayer}
			\draw (17) to (18.center);
		\end{pgfonlayer}
	\end{tikzpicture} := 	
	\begin{tikzpicture}
		\begin{pgfonlayer}{nodelayer}
			\node [style=circle] (0) at (-0.75, 1.5) {};
			\node [style=none] (1) at (-0.75, 0.5) {};
			\node [style=none] (2) at (-1.75, 0.5) {};
			\node [style=none] (3) at (-1.75, 2.25) {};
			\node [style=none] (4) at (-2, 2) {$B$};
			\node [style=none] (5) at (-1.25, -0.25) {$\epsilon_L$};
		\end{pgfonlayer}
		\begin{pgfonlayer}{edgelayer}
			\draw (0) to (1.center);
			\draw [bend left=90, looseness=1.50] (1.center) to (2.center);
			\draw (2.center) to (3.center);
		\end{pgfonlayer}
	\end{tikzpicture} = 
	\begin{tikzpicture}
		\begin{pgfonlayer}{nodelayer}
			\node [style=circle] (0) at (-2, 1.5) {};
			\node [style=none] (1) at (-2, 0.5) {};
			\node [style=none] (2) at (-1, 0.5) {};
			\node [style=none] (3) at (-1, 2.25) {};
			\node [style=none] (4) at (-0.75, 2) {$B$};
			\node [style=none] (5) at (-1.5, -0.25) {$\epsilon_R$};
		\end{pgfonlayer}
		\begin{pgfonlayer}{edgelayer}
			\draw (0) to (1.center);
			\draw [bend right=90, looseness=1.50] (1.center) to (2.center);
			\draw (2.center) to (3.center);
		\end{pgfonlayer}
	\end{tikzpicture} 
\end{equation}
\end{definition}

Note that, in a symmetric LDC, given a dual $(\eta_L, \epsilon_L): A \dashvv B$, any linear monoid $A \linmonw B$, 
can be normalized to use the symmetric dual $(\eta_L c_\oa, c_\ox, \epsilon_L): B \dashvv A$. Suppose, in a symmetric LDC, 
$(\eta_L, \epsilon_L): A \dashvv B$, and $(\eta_R, \epsilon_R): A \dashvv B$ are cyclic duals, then:
\[  \begin{tikzpicture}
	\begin{pgfonlayer}{nodelayer}
		\node [style=none] (0) at (1.25, 2.5) {};
		\node [style=none] (1) at (1.25, 1) {};
		\node [style=none] (2) at (2.5, 1) {};
		\node [style=none] (3) at (2.5, 1.75) {};
		\node [style=none] (4) at (1.75, 1.75) {};
		\node [style=none] (5) at (1.75, -0.5) {};
		\node [style=none] (6) at (2, 2.5) {$\eta_R$};
		\node [style=none] (7) at (2, 1.5) {$B$};
		\node [style=none] (8) at (2.75, 1.5) {$A$};
		\node [style=none] (9) at (1, 2) {$B$};
		\node [style=none] (10) at (1.25, 0.25) {$\epsilon_L$};
		\node [style=none] (11) at (2, -0.25) {$B$};
	\end{pgfonlayer}
	\begin{pgfonlayer}{edgelayer}
		\draw (0.center) to (1.center);
		\draw [bend right=90, looseness=1.50] (1.center) to (2.center);
		\draw (2.center) to (3.center);
		\draw [bend left=270, looseness=1.75] (3.center) to (4.center);
		\draw (4.center) to (5.center);
	\end{pgfonlayer}
\end{tikzpicture} = \left( \begin{tikzpicture}
	\begin{pgfonlayer}{nodelayer}
		\node [style=none] (0) at (2.5, 2.5) {};
		\node [style=none] (1) at (2.5, 1) {};
		\node [style=none] (2) at (1.25, 1) {};
		\node [style=none] (3) at (1.25, 1.75) {};
		\node [style=none] (4) at (2, 1.75) {};
		\node [style=none] (5) at (2, -0.5) {};
		\node [style=none] (6) at (1.75, 2.5) {$\eta_L$};
		\node [style=none] (7) at (1.75, 1.5) {$B$};
		\node [style=none] (8) at (1, 1.5) {$A$};
		\node [style=none] (9) at (2.75, 2) {$B$};
		\node [style=none] (10) at (2.5, 0.25) {$\epsilon_R$};
		\node [style=none] (11) at (1.75, -0.25) {$B$};
	\end{pgfonlayer}
	\begin{pgfonlayer}{edgelayer}
		\draw (0.center) to (1.center);
		\draw [bend left=90, looseness=1.50] (1.center) to (2.center);
		\draw (2.center) to (3.center);
		\draw [bend right=270, looseness=1.75] (3.center) to (4.center);
		\draw (4.center) to (5.center);
	\end{pgfonlayer}
\end{tikzpicture} \right)^{-1} \]
which is an isomorphism of duals: \[ \begin{tikzpicture}
	\begin{pgfonlayer}{nodelayer}
		\node [style=none] (0) at (2.5, 2.25) {};
		\node [style=none] (1) at (2.5, 1) {};
		\node [style=none] (2) at (1.25, 1) {};
		\node [style=none] (3) at (1.25, 1.25) {};
		\node [style=none] (4) at (2, 1.25) {};
		\node [style=none] (5) at (2, -0.5) {};
		\node [style=none] (6) at (1.5, 2) {$\eta_L$};
		\node [style=none] (7) at (1.75, 1) {$B$};
		\node [style=none] (8) at (1, 1) {$A$};
		\node [style=none] (9) at (2.75, 1.5) {$B$};
		\node [style=none] (10) at (2.5, 0.25) {$\epsilon_R$};
		\node [style=none] (11) at (1.75, -0.25) {$B$};
		\node [style=none] (12) at (3.25, 2.25) {};
		\node [style=none] (13) at (3.25, -0.5) {};
		\node [style=none] (14) at (2.75, 3) {$\eta_R$};
	\end{pgfonlayer}
	\begin{pgfonlayer}{edgelayer}
		\draw (0.center) to (1.center);
		\draw [bend left=90, looseness=1.50] (1.center) to (2.center);
		\draw (2.center) to (3.center);
		\draw [bend right=270, looseness=1.75] (3.center) to (4.center);
		\draw (4.center) to (5.center);
		\draw (12.center) to (13.center);
		\draw [bend left=90, looseness=2.00] (0.center) to (12.center);
	\end{pgfonlayer}
\end{tikzpicture} = \begin{tikzpicture}
	\begin{pgfonlayer}{nodelayer}
		\node [style=none] (3) at (2.25, 2.75) {};
		\node [style=none] (4) at (2.75, 2.75) {};
		\node [style=none] (5) at (1.75, 1.25) {};
		\node [style=none] (6) at (2.5, 3.25) {$\eta_L$};
		\node [style=none] (7) at (3, 2.75) {$B$};
		\node [style=none] (8) at (2, 2.75) {$A$};
		\node [style=none] (11) at (1.5, 1) {$B$};
		\node [style=none] (12) at (3.25, 1.25) {};
		\node [style=none] (13) at (3.5, 1) {$A$};
		\node [style=none] (14) at (3.25, 0.75) {};
		\node [style=none] (15) at (1.75, 0.75) {};
	\end{pgfonlayer}
	\begin{pgfonlayer}{edgelayer}
		\draw [bend right=270, looseness=1.75] (3.center) to (4.center);
		\draw [in=90, out=-90] (4.center) to (5.center);
		\draw [in=90, out=-90] (3.center) to (12.center);
		\draw (14.center) to (12.center);
		\draw (15.center) to (5.center);
	\end{pgfonlayer}
\end{tikzpicture} ~~~~~~~~~ \begin{tikzpicture}
	\begin{pgfonlayer}{nodelayer}
		\node [style=none] (0) at (1.75, 0.25) {};
		\node [style=none] (1) at (1.75, 1.5) {};
		\node [style=none] (2) at (3, 1.5) {};
		\node [style=none] (3) at (3, 1.25) {};
		\node [style=none] (4) at (2.25, 1.25) {};
		\node [style=none] (5) at (2.25, 3) {};
		\node [style=none] (6) at (2.75, 0.5) {$\epsilon_L$};
		\node [style=none] (7) at (2.5, 1.5) {$B$};
		\node [style=none] (8) at (3.25, 1.5) {$A$};
		\node [style=none] (9) at (1.5, 1) {$B$};
		\node [style=none] (10) at (1.75, 2.25) {$\eta_R$};
		\node [style=none] (11) at (2.5, 2.75) {$B$};
		\node [style=none] (12) at (1, 0.25) {};
		\node [style=none] (13) at (1, 3) {};
		\node [style=none] (14) at (1.5, -0.5) {$\epsilon_R$};
		\node [style=none] (15) at (0.75, 2.5) {$A$};
	\end{pgfonlayer}
	\begin{pgfonlayer}{edgelayer}
		\draw (0.center) to (1.center);
		\draw [bend left=90, looseness=1.50] (1.center) to (2.center);
		\draw (2.center) to (3.center);
		\draw [bend right=270, looseness=1.75] (3.center) to (4.center);
		\draw (4.center) to (5.center);
		\draw (12.center) to (13.center);
		\draw [bend left=90, looseness=2.00] (0.center) to (12.center);
	\end{pgfonlayer}
\end{tikzpicture} = \begin{tikzpicture}
	\begin{pgfonlayer}{nodelayer}
		\node [style=none] (3) at (2.75, 1.25) {};
		\node [style=none] (4) at (2.25, 1.25) {};
		\node [style=none] (5) at (3.25, 2.75) {};
		\node [style=none] (6) at (2.5, 0.75) {$\epsilon_L$};
		\node [style=none] (7) at (2, 1.25) {$B$};
		\node [style=none] (8) at (3, 1.25) {$A$};
		\node [style=none] (11) at (3.5, 3) {$B$};
		\node [style=none] (12) at (1.75, 2.75) {};
		\node [style=none] (13) at (1.5, 3) {$A$};
		\node [style=none] (14) at (1.75, 3.25) {};
		\node [style=none] (15) at (3.25, 3.25) {};
	\end{pgfonlayer}
	\begin{pgfonlayer}{edgelayer}
		\draw [bend right=270, looseness=1.75] (3.center) to (4.center);
		\draw [in=-90, out=90] (4.center) to (5.center);
		\draw [in=-90, out=90] (3.center) to (12.center);
		\draw (14.center) to (12.center);
		\draw (15.center) to (5.center);
	\end{pgfonlayer}
\end{tikzpicture} \] 
Hence, in a symmetric LDC, for any linear monoid, one can assume that the duals are given using the symmetry map. 
A linear monoid is said to be {\bf symmetric} when its duals are symmetric. Section \ref{Sec: linear monoid examples} 
discusses different examples of linear monoids.

An alternate form for linear monoids, which displays their similarity to Frobenius algebras, is:

\begin{proposition}
	\label{Lemma: alternate presentation of linear monoids}
	A {\bf linear monoid}, $A \linmonw B$, in an LDC is equivalent to the following data:
	\begin{itemize}
	\item a monoid $(A, \mulmap{1.5}{white}: A \ox A \to A, \unitmap{1.5}{white}: \top \to A) $
	\item a comonoid $(B, \comulmap{1.5}{white}: B \to B \oa B, \counitmap{1.5}{white}: B \to \bot) $
	\item actions, $\leftaction{0.5}{white}: A \ox B \to B$, $\rightaction{0.5}{white}: B \ox A \to B$,
	and coactions $\leftcoaction{0.55}{white}: A \to B \oa A$, $\rightcoaction{0.5}{white}: A \to A \oa B$,
    \end{itemize}
	such that the following axioms (and their `op' and `co' symmetric forms) hold:  
    \[ (a) ~~
 \]
If the linear monoid is symmetric, then: 
$ %
$
\end{proposition}
\begin{proof}  Suppose $A \linmonw B$ is a linear monoid. Then, the comonoid on $B$
	is given by the left dual (or equivalently the right dual) of the monoid. The action and the 
	coaction maps are defined as follows:
	\[%
	 \]  
	Equations $(a)$-$(d)$ can be verified easily. 

	For the converse, the left and the right duals are defined as follows:
	\[ \eta_L := %
 \]
	\end{proof}

A linear monoid, $A \linmonw B$, in a monoidal category gives a Frobenius Algebra when it is a 
self-linear monoid, that is $A = B$, and the dualities coincide with the self-dual cup, and the cap. Note that 
while a Frobenius Algebra is always on a self-dual object, a linear monoid allows Frobenius interaction 
between distinct dual objects.

\begin{definition}
	\label{defn: linear monoid morphism}
A {\bf morphism of linear monoids} is a pair of maps, $(f,g): (A \linmonw B) \to (A' \linmonw B')$, such that 
$f: A \to A'$ is a monoid morphism (or equivalently $g: B' \to B$ is a comonoid morphism),
 and $(f, g)$ and $(g,f)$ preserves the left and the right duals respectively. 
 \end{definition}

Note that a morphism of Frobenius algebras is usually given by a single monoid and comonoid morphism which has the 
effect of giving an isomorphism. However, the morphisms of linear monoids as in Definition \ref{defn: linear monoid morphism}, 
are not restricted to isomorphisms. Here, the comonoid morphism $g: B' \to B$, is the cyclic mate of the monoid morphism, $f: A \to A'$. 

\begin{definition}
A {\bf $\dagger$-linear monoid}, $(A,  \mulmap{1.2}{white}, \unitmap{1.2}{white}) 
\dagmonw (A^\dagger, \comulmap{1.2}{white}, \counitmap{1.2}{white})$, 
in a $\dagger$-LDC, is a linear monoid such that $(\eta_L, \epsilon_L): A \dashvv A^\dagger$, 
and $(\eta_R, \epsilon_R): A^\dagger \dashvv A$ are $\dagger$-duals and:
\begin{equation}
	\label{MainEqn: rightdagmon}
	\begin{tikzpicture}
		\begin{pgfonlayer}{nodelayer}
			\node [style=circle] (17) at (1, 0.25) {};
			\node [style=none] (18) at (1, 2.5) {};
			\node [style=none] (21) at (1.25, 2.25) {$A^\dagger$};
		\end{pgfonlayer}
		\begin{pgfonlayer}{edgelayer}
			\draw (17) to (18.center);
		\end{pgfonlayer}
	\end{tikzpicture} := 	
	\begin{tikzpicture}
		\begin{pgfonlayer}{nodelayer}
			\node [style=circle] (0) at (-0.75, 3.5) {};
			\node [style=none] (1) at (-0.75, 2.5) {};
			\node [style=none] (2) at (-1.75, 2.5) {};
			\node [style=none] (3) at (-1.75, 4) {};
			\node [style=none] (4) at (-1.75, 2.5) {};
			\node [style=none] (5) at (-2.2, 3.75) {$A^\dagger$};
			\node [style=none] (6) at (-1.25, 1.75) {$\epsilon_L$};
		\end{pgfonlayer}
		\begin{pgfonlayer}{edgelayer}
			\draw (0) to (1.center);
			\draw [bend left=75, looseness=1.75] (1.center) to (2.center);
			\draw (2.center) to (3.center);
		\end{pgfonlayer}
	\end{tikzpicture} = \begin{tikzpicture}
		\begin{pgfonlayer}{nodelayer}
			\node [style=circle] (0) at (-1.75, 1.5) {};
			\node [style=none] (1) at (-1.75, 2.5) {};
			\node [style=none] (2) at (-1.4, 0.25) {$A^\dagger$};
			\node [style=none] (3) at (-2.5, 0.75) {};
			\node [style=none] (4) at (-1, 0.75) {};
			\node [style=none] (5) at (-1, 2.5) {};
			\node [style=none] (6) at (-2.5, 2.5) {};
			\node [style=none] (7) at (-1.75, 0.75) {};
			\node [style=none] (8) at (-1.75, -0) {};
		\end{pgfonlayer}
		\begin{pgfonlayer}{edgelayer}
			\draw (0) to (1.center);
			\draw (3.center) to (4.center);
			\draw (4.center) to (5.center);
			\draw (5.center) to (6.center);
			\draw (6.center) to (3.center);
			\draw (8.center) to (7.center);
		\end{pgfonlayer}
	\end{tikzpicture}   ~~~~~~~~
	\begin{tikzpicture}
		\begin{pgfonlayer}{nodelayer}
			\node [style=none] (12) at (2, 0.25) {$A^\dagger$};
			\node [style=none] (13) at (0, 0.25) {$A^\dagger$};
			\node [style=circle] (17) at (1, 1.25) {};
			\node [style=none] (18) at (1, 2.5) {};
			\node [style=none] (19) at (0.5, 0) {};
			\node [style=none] (20) at (1.5, 0) {};
			\node [style=none] (21) at (1.5, 2.25) {$A^\dagger$};
		\end{pgfonlayer}
		\begin{pgfonlayer}{edgelayer}
			\draw [in=-150, out=90, looseness=1.25] (19.center) to (17);
			\draw [in=90, out=-30, looseness=1.25] (17) to (20.center);
			\draw (17) to (18.center);
		\end{pgfonlayer}
	\end{tikzpicture} := \begin{tikzpicture}
		\begin{pgfonlayer}{nodelayer}
			\node [style=circle] (0) at (-2.75, 0.75) {};
			\node [style=none] (1) at (-3.25, 1.25) {};
			\node [style=none] (2) at (-2.75, 0.5) {};
			\node [style=none] (3) at (-2.25, 1.25) {};
			\node [style=none] (4) at (-1.5, 1.25) {};
			\node [style=none] (5) at (-1.5, -0) {};
			\node [style=none] (6) at (-3.25, 1.5) {};
			\node [style=none] (7) at (-1, 1.5) {};
			\node [style=none] (8) at (-1, -0) {};
			\node [style=none] (9) at (-3.75, 0.5) {};
			\node [style=none] (10) at (-3.75, 2.25) {};
			\node [style=none] (11) at (-4.1, 2) {$A^\dagger$};
			\node [style=none] (12) at (-1.85, 0.25) {$A^\dagger$};
			\node [style=none] (13) at (-0.65, 0.25) {$A^\dagger$};
			\node [style=none] (14) at (-2, 2.5) {$\eta_L$};
			\node [style=none] (15) at (-3.25, -0.25) {$\epsilon_L$};
		\end{pgfonlayer}
		\begin{pgfonlayer}{edgelayer}
			\draw [in=150, out=-90, looseness=1.00] (1.center) to (0);
			\draw [in=-90, out=30, looseness=1.00] (0) to (3.center);
			\draw (0) to (2.center);
			\draw [bend left=90, looseness=1.50] (2.center) to (9.center);
			\draw (9.center) to (10.center);
			\draw (6.center) to (1.center);
			\draw [bend left=90, looseness=1.00] (6.center) to (7.center);
			\draw (7.center) to (8.center);
			\draw (4.center) to (5.center);
			\draw [bend right=90, looseness=2.00] (4.center) to (3.center);
		\end{pgfonlayer}
	\end{tikzpicture} = \begin{tikzpicture}
		\begin{pgfonlayer}{nodelayer}
			\node [style=circle] (0) at (-2.75, 1.5) {};
			\node [style=none] (1) at (-3.25, 2.25) {};
			\node [style=none] (2) at (-2.75, 1) {};
			\node [style=none] (3) at (-2.25, 2.25) {};
			\node [style=none] (4) at (-4, 1) {};
			\node [style=none] (5) at (-1.5, 1) {};
			\node [style=none] (6) at (-1.5, 2.25) {};
			\node [style=none] (7) at (-4, 2.25) {};
			\node [style=none] (8) at (-2.75, 2.75) {};
			\node [style=none] (9) at (-2.75, 2.25) {};
			\node [style=none] (10) at (-3.5, 1) {};
			\node [style=none] (11) at (-3.5, 0.25) {};
			\node [style=none] (12) at (-2, 1) {};
			\node [style=none] (13) at (-2, 0.25) {};
			\node [style=none] (14) at (-2.35, 2.75) {$A^\dagger$};
			\node [style=none] (15) at (-1.4, 0.5) {$A^\dagger$};
			\node [style=none] (16) at (-4, 0.5) {$A^\dagger$};
		\end{pgfonlayer}
		\begin{pgfonlayer}{edgelayer}
			\draw [in=150, out=-90] (1.center) to (0);
			\draw [in=-90, out=30] (0) to (3.center);
			\draw (0) to (2.center);
			\draw (7.center) to (6.center);
			\draw (5.center) to (6.center);
			\draw (5.center) to (4.center);
			\draw (4.center) to (7.center);
			\draw (8.center) to (9.center);
			\draw (10.center) to (11.center);
			\draw (12.center) to (13.center);
		\end{pgfonlayer}
	\end{tikzpicture}	
\end{equation}
\end{definition}

In a $\dagger$-monoidal category, all $\dagger$-Frobenius algebras are 
also $\dagger$-linear monoids. However, the converse is not true. For example,  
in the category of complex Hilbert spaces and linear maps, the basic 
Weil Algebra $\C[x]/x^2=0$ is a commutative $\dagger$-linear monoid but not 
a $\dagger$-Frobenius Algebra. More examples are discussed in Section \ref{Sec: linear monoid examples}.

\begin{definition} 
	\label{Defn: twisted}
A $\dagger$-linear monoid is {\bf twisted} if the dual of the multiplication coincides with the dagger of the multiplication 
but with a twist: 
\begin{equation} 
	\label{eqn: twisted monoid}
	\begin{tikzpicture}
		\begin{pgfonlayer}{nodelayer}
			\node [style=none] (12) at (2, 0.25) {$A^\dagger$};
			\node [style=none] (13) at (0, 0.25) {$A^\dagger$};
			\node [style=circle] (17) at (1, 1.25) {};
			\node [style=none] (18) at (1, 2.5) {};
			\node [style=none] (19) at (0.5, 0) {};
			\node [style=none] (20) at (1.5, 0) {};
			\node [style=none] (21) at (1.5, 2.25) {$A^\dagger$};
		\end{pgfonlayer}
		\begin{pgfonlayer}{edgelayer}
			\draw [in=-150, out=90, looseness=1.25] (19.center) to (17);
			\draw [in=90, out=-30, looseness=1.25] (17) to (20.center);
			\draw (17) to (18.center);
		\end{pgfonlayer}
	\end{tikzpicture} := \begin{tikzpicture}
		\begin{pgfonlayer}{nodelayer}
			\node [style=circle] (0) at (-2.75, 0.75) {};
			\node [style=none] (1) at (-3.25, 1.25) {};
			\node [style=none] (2) at (-2.75, 0.5) {};
			\node [style=none] (3) at (-2.25, 1.25) {};
			\node [style=none] (4) at (-1.5, 1.25) {};
			\node [style=none] (5) at (-1.5, -0) {};
			\node [style=none] (6) at (-3.25, 1.5) {};
			\node [style=none] (7) at (-1, 1.5) {};
			\node [style=none] (8) at (-1, -0) {};
			\node [style=none] (9) at (-3.75, 0.5) {};
			\node [style=none] (10) at (-3.75, 2.25) {};
			\node [style=none] (11) at (-4.1, 2) {$A^\dagger$};
			\node [style=none] (12) at (-1.85, 0.25) {$A^\dagger$};
			\node [style=none] (13) at (-0.65, 0.25) {$A^\dagger$};
			\node [style=none] (14) at (-2, 2.5) {$\eta_L$};
			\node [style=none] (15) at (-3.25, -0.25) {$\epsilon_L$};
		\end{pgfonlayer}
		\begin{pgfonlayer}{edgelayer}
			\draw [in=150, out=-90, looseness=1.00] (1.center) to (0);
			\draw [in=-90, out=30, looseness=1.00] (0) to (3.center);
			\draw (0) to (2.center);
			\draw [bend left=90, looseness=1.50] (2.center) to (9.center);
			\draw (9.center) to (10.center);
			\draw (6.center) to (1.center);
			\draw [bend left=90, looseness=1.00] (6.center) to (7.center);
			\draw (7.center) to (8.center);
			\draw (4.center) to (5.center);
			\draw [bend right=90, looseness=2.00] (4.center) to (3.center);
		\end{pgfonlayer}
	\end{tikzpicture} = \begin{tikzpicture}
		\begin{pgfonlayer}{nodelayer}
			\node [style=circle] (0) at (-2.75, 1.5) {};
			\node [style=none] (1) at (-3.25, 2.25) {};
			\node [style=none] (2) at (-2.75, 1) {};
			\node [style=none] (3) at (-2.25, 2.25) {};
			\node [style=none] (4) at (-4, 1) {};
			\node [style=none] (5) at (-1.5, 1) {};
			\node [style=none] (6) at (-1.5, 2.25) {};
			\node [style=none] (7) at (-4, 2.25) {};
			\node [style=none] (8) at (-2.75, 2.75) {};
			\node [style=none] (9) at (-2.75, 2.25) {};
			\node [style=none] (10) at (-3.25, 1) {};
			\node [style=none] (11) at (-2.25, -0.25) {};
			\node [style=none] (12) at (-2.25, 1) {};
			\node [style=none] (13) at (-3.25, -0.25) {};
			\node [style=none] (14) at (-2.35, 2.75) {$A^\dagger$};
			\node [style=none] (15) at (-1.9, 0) {$A^\dagger$};
			\node [style=none] (16) at (-3.5, 0) {$A^\dagger$};
		\end{pgfonlayer}
		\begin{pgfonlayer}{edgelayer}
			\draw [in=150, out=-90] (1.center) to (0);
			\draw [in=-90, out=30] (0) to (3.center);
			\draw (0) to (2.center);
			\draw (7.center) to (6.center);
			\draw (5.center) to (6.center);
			\draw (5.center) to (4.center);
			\draw (4.center) to (7.center);
			\draw (8.center) to (9.center);
			\draw [in=90, out=-90] (10.center) to (11.center);
			\draw [in=90, out=-90] (12.center) to (13.center);
		\end{pgfonlayer}
	\end{tikzpicture}	
\end{equation}
\end{definition}

In symmetric $\dagger$-LDCs, endomorphism monoids (pants monoids) is a twisted $\dagger$-linear monoid, 
see Section \ref{Sec: pants linear monoid}. Note that, every commutative $\dagger$-linear monoid 
is a twisted $\dagger$-linear monoid.

Following the same pattern as $\dagger$-duals, the $\dagger$-linear monoids can be classified as follows:
\begin{itemize}
\item a {\bf (twisted) right $\dagger$-linear monoid} which is a linear monoid $A \linmonw A^\dag$ with 
left and right $\dagger$-duals such that the dual of the monoid is same as its (twisted) dagger,
\item a {\bf (twisted) left $\dagger$-linear monoid}, that is, a linear monoid $A^\dag \linmonw A$ with 
left and right $\dagger$-duals such that the dual of the comonoid is same as its (twisted) dagger, and 
\item a {\bf general $\dagger$-linear monoid} which is a linear monoid 
$(A, \mulmap{1.5}{white}, \unitmap{1.5}{white}) \linmonw (B, \comulmap{1.5}{white}, \counitmap{1.5}{white})$ satisfying 
one of the following three conditions: 
\begin{enumerate}[(i)]
\item  $A \linmonw B$ is isomorphic to its right $\dagger$-linear 
monoid via an isomorphism $(1_A, q): A \linmonw B \to A \dagmonw A^\dag$, 
\item $A \linmonw B$ is isomorphic to its left $\dagger$-
dual via an isomorphism $(p, 1_B): B^\dag \dagmonw B \to A \linmonw B$, or 
\item $A \linmonw B$ is isomorphic to $B^\dag \linmonb A^\dag$ via an isomorphism $(p, q)$ and 
$pq^\dagger = \iota_A$; 
\end{enumerate}
a general $\dagger$-linear monoid is {\bf twisted} if one of the following conditions hold:
\begin{itemize} 
     \item in condition $(i)$, $q: (A^\dag, \mulmap{1.5}{white}^\dag, 
\unitmap{1.5}{white}^\dagger) \to (B, \comulmap{1.5}{white}, \counitmap{1.5}{white})$ 
is a twisted comonoid morphism: 
\[ \begin{tikzpicture}
	\begin{pgfonlayer}{nodelayer}
		\node [style=none] (4) at (0.25, -2.5) {};
		\node [style=none] (5) at (1.25, -2.5) {};
		\node [style=none] (6) at (1.25, -1.75) {};
		\node [style=none] (7) at (0.25, -1.75) {};
		\node [style=circle] (8) at (0.75, -1.25) {};
		\node [style=none] (9) at (0.75, 0.25) {};
		\node [style=none] (10) at (1.5, -2.25) {$B$};
		\node [style=none] (11) at (0, -2.25) {$B$};
		\node [style=none] (12) at (1.25, 0.25) {$A^\dag$};
		\node [style=circle, scale=1.5] (13) at (0.75, -0.5) {};
		\node [style=none] (14) at (0.75, -0.5) {$q$};
	\end{pgfonlayer}
	\begin{pgfonlayer}{edgelayer}
		\draw [in=-45, out=90] (6.center) to (8);
		\draw [in=-135, out=90] (7.center) to (8);
		\draw (7.center) to (4.center);
		\draw (6.center) to (5.center);
		\draw (9.center) to (13);
		\draw (13) to (8);
	\end{pgfonlayer}
\end{tikzpicture} =\begin{tikzpicture}
	\begin{pgfonlayer}{nodelayer}
		\node [style=circle, scale=1.5] (0) at (0.25, -1.75) {};
		\node [style=circle, scale=1.5] (1) at (1.25, -1.75) {};
		\node [style=none] (2) at (0.25, -1.75) {$q$};
		\node [style=none] (3) at (1.25, -1.75) {$q$};
		\node [style=none] (4) at (0.25, -2.5) {};
		\node [style=none] (5) at (1.25, -2.5) {};
		\node [style=none] (6) at (1.25, -1) {};
		\node [style=none] (7) at (0.25, -1) {};
		\node [style=none] (9) at (0.75, 0.75) {};
		\node [style=none] (10) at (1.5, -2.25) {$B$};
		\node [style=none] (11) at (0, -2.25) {$B$};
		\node [style=none] (12) at (1.25, 0.5) {$A^\dag$};
		\node [style=none] (13) at (1.25, 0) {};
		\node [style=none] (14) at (0.25, 0) {};
		\node [style=circle] (15) at (0.75, -0.5) {};
		\node [style=none] (16) at (0.75, -1) {};
		\node [style=none] (17) at (-0.25, 0) {};
		\node [style=none] (18) at (1.75, 0) {};
		\node [style=none] (19) at (1.75, -1) {};
		\node [style=none] (20) at (-0.25, -1) {};
		\node [style=none] (21) at (0.75, 0) {};
	\end{pgfonlayer}
	\begin{pgfonlayer}{edgelayer}
		\draw (4.center) to (0);
		\draw (5.center) to (1);
		\draw [in=90, out=-90] (6.center) to (0);
		\draw [in=90, out=-90, looseness=0.75] (7.center) to (1);
		\draw [in=45, out=-90] (13.center) to (15);
		\draw [in=135, out=-90] (14.center) to (15);
		\draw (16.center) to (15);
		\draw (17.center) to (20.center);
		\draw (20.center) to (19.center);
		\draw (19.center) to (18.center);
		\draw (18.center) to (17.center);
		\draw (21.center) to (9.center);
	\end{pgfonlayer}
\end{tikzpicture}  \]
\item in condition  $(ii)$, $p: (A, \mulmap{1.5}{white}, \unitmap{1.5}{white}) 
\to (B^\dag, \comulmap{1.5}{white}^\dag, \counitmap{1.5}{white}^\dag)$
is a twisted monoid morphism:
\[  \begin{tikzpicture}
	\begin{pgfonlayer}{nodelayer}
		\node [style=none] (4) at (0.25, 0.25) {};
		\node [style=none] (5) at (1.25, 0.25) {};
		\node [style=none] (6) at (1.25, -0.25) {};
		\node [style=none] (7) at (0.25, -0.25) {};
		\node [style=circle] (8) at (0.75, -1) {};
		\node [style=none] (9) at (0.75, -2.5) {};
		\node [style=none] (10) at (1.5, 0) {$A$};
		\node [style=none] (11) at (0, 0) {$A$};
		\node [style=none] (12) at (1.25, -2.25) {$B^\dag$};
		\node [style=circle, scale=1.5] (13) at (0.75, -1.75) {};
		\node [style=none] (14) at (0.75, -1.75) {$p$};
	\end{pgfonlayer}
	\begin{pgfonlayer}{edgelayer}
		\draw [in=45, out=-90] (6.center) to (8);
		\draw [in=135, out=-90] (7.center) to (8);
		\draw (4.center) to (7.center);
		\draw (5.center) to (6.center);
		\draw (9.center) to (13);
		\draw (13) to (8);
	\end{pgfonlayer}
\end{tikzpicture}
 = \begin{tikzpicture}
	\begin{pgfonlayer}{nodelayer}
		\node [style=circle, scale=1.5] (0) at (0.25, 0) {};
		\node [style=circle, scale=1.5] (1) at (1.25, 0) {};
		\node [style=none] (2) at (0.25, 0) {$p$};
		\node [style=none] (3) at (1.25, 0) {$p$};
		\node [style=none] (4) at (0.25, 0.5) {};
		\node [style=none] (5) at (1.25, 0.5) {};
		\node [style=none] (9) at (0.75, -2.5) {};
		\node [style=none] (10) at (1.5, 0.5) {$A$};
		\node [style=none] (11) at (0, 0.5) {$A$};
		\node [style=none] (12) at (1.25, -2.25) {$B^\dag$};
		\node [style=none] (13) at (0, -1) {};
		\node [style=none] (14) at (0, -2) {};
		\node [style=none] (15) at (1.5, -2) {};
		\node [style=none] (16) at (1.5, -1) {};
		\node [style=none] (17) at (1.25, -2) {};
		\node [style=none] (18) at (0.25, -2) {};
		\node [style=circle] (19) at (0.75, -1.5) {};
		\node [style=none] (20) at (0.75, -1) {};
		\node [style=none] (21) at (0.75, -2) {};
		\node [style=none] (22) at (0.25, -1) {};
		\node [style=none] (23) at (1.25, -1) {};
	\end{pgfonlayer}
	\begin{pgfonlayer}{edgelayer}
		\draw (4.center) to (0);
		\draw (5.center) to (1);
		\draw [in=-45, out=90] (17.center) to (19);
		\draw [in=-135, out=90] (18.center) to (19);
		\draw (20.center) to (19);
		\draw (14.center) to (15.center);
		\draw (15.center) to (16.center);
		\draw (13.center) to (14.center);
		\draw (13.center) to (16.center);
		\draw (21.center) to (9.center);
		\draw [in=-90, out=90] (22.center) to (1);
		\draw [in=-105, out=105, looseness=0.75] (23.center) to (0);
	\end{pgfonlayer}
\end{tikzpicture}\]

\item  in condition $(iii)$, either $p$ is twisted monoid morphism, or $q$ is a twisted comonoid morphism. 
\end{itemize}
\end{itemize}

In condition $(i)$  for a general $\dagger$-linear monoid, the identity map completely determines $q$. Hence, 
we observe that the existence of a right $\dagger$-linear monoid is all that is needed 
for any linear monoid to be a  $\dagger$-linear monoid. By Lemma \ref{Lemma: twisted extension}, the observation extends to 
twisted $\dagger$-linear monoids too: the existence of a twisted right $\dagger$-linear monoid is all that 
is needed any linear monoid to be a twisted $\dagger$-linear monoid.

A linear monoid has a twisted isomorphism to its right (or left) $\dagger$-linear monoid 
if and only if the right (or the left) $\dagger$-linear monoid is twisted:
\begin{lemma}
	\label{Lemma: twisted extension}
In a $\dagger$-LDC, let $(\eta, \epsilon): (A, \mulmap{1.5}{white}, \unitmap{1.5}{white}) \linmonw \comonoid{B}$ 
be a linear monoid with an isomorphism its right $\dagger$-linear monoid: 
\[  \begin{matrix}
	\xymatrixcolsep{2cm}
	\xymatrix{
	\monoid{A} \ar@{-||}[r]^{(\eta,  \epsilon)~ \circ} \ar@{=}[d]_{} &  \comonoid{B} \ar@{<-}[d]^{q}    \\
	\monoid{A}   \ar@{-||}[r]_{\dagger (\tau, \gamma) ~\circ}  & (A^\dagger, \mulmap{1.5}{white}^\dagger, \unitmap{1.5}{white}^\dagger) } 
\end{matrix} \]
The isomorphism $(1,q)$ is twisted, that is, $q$ is a twisted comonoid morphism if and 
only if  $A \dagmonw A^\dag$ is a twisted $\dagger$-linear monoid.
\end{lemma}
\begin{proof}
Assume that $q$ is a twisted monoid morphism. Then, $A \dagmonw A^\dag$ is a twisted $\dagger$-linear monoid because: 
\[ \begin{tikzpicture}
	\begin{pgfonlayer}{nodelayer}
		\node [style=none] (9) at (-2.25, 0.5) {};
		\node [style=none] (10) at (1.75, -3.75) {$A^\dag$};
		\node [style=none] (11) at (-0.25, -3.75) {$A^\dag$};
		\node [style=none] (12) at (-1.75, 0.25) {$A^\dag$};
		\node [style=none] (17) at (-0.25, -1) {};
		\node [style=none] (18) at (-1.25, -1) {};
		\node [style=circle] (19) at (-0.75, -1.5) {};
		\node [style=none] (21) at (-2.25, 0) {};
		\node [style=none] (28) at (0.25, -4) {};
		\node [style=none] (29) at (1.25, -4) {};
		\node [style=none] (36) at (-0.75, -2) {};
		\node [style=none] (37) at (0.25, -1) {};
		\node [style=none] (38) at (1.25, -0.5) {};
		\node [style=none] (39) at (-1.25, -0.5) {};
		\node [style=none] (40) at (-2.25, -2) {};
		\node [style=none] (41) at (0, 0.25) {$\tau_L$};
		\node [style=none] (42) at (0, -0.5) {$\tau_L$};
		\node [style=none] (43) at (-1.5, -3) {$\gamma_L$};
	\end{pgfonlayer}
	\begin{pgfonlayer}{edgelayer}
		\draw [in=45, out=-90] (17.center) to (19);
		\draw [in=135, out=-90] (18.center) to (19);
		\draw (21.center) to (9.center);
		\draw (36.center) to (19);
		\draw [bend left=90, looseness=2.00] (17.center) to (37.center);
		\draw [bend left=90, looseness=0.75] (39.center) to (38.center);
		\draw (39.center) to (18.center);
		\draw [bend right=90, looseness=1.50] (40.center) to (36.center);
		\draw (21.center) to (40.center);
		\draw (28.center) to (37.center);
		\draw (38.center) to (29.center);
	\end{pgfonlayer}
\end{tikzpicture} \stackrel{(1)}{=} \begin{tikzpicture}
	\begin{pgfonlayer}{nodelayer}
		\node [style=none] (9) at (-2.25, 0.5) {};
		\node [style=none] (10) at (1.75, -3.75) {$A^\dag$};
		\node [style=none] (11) at (-0.25, -3.75) {$A^\dag$};
		\node [style=none] (12) at (-1.75, 0.25) {$A^\dag$};
		\node [style=none] (17) at (-0.25, -1) {};
		\node [style=none] (18) at (-1.25, -1) {};
		\node [style=circle] (19) at (-0.75, -1.5) {};
		\node [style=none] (21) at (-2.25, 0) {};
		\node [style=circle, scale=2] (26) at (0.25, -3.25) {};
		\node [style=circle, scale=2] (27) at (1.25, -3.25) {};
		\node [style=none] (28) at (0.25, -4) {};
		\node [style=none] (29) at (1.25, -4) {};
		\node [style=none] (32) at (0.25, -3.25) {$q^{-1}$};
		\node [style=none] (33) at (1.25, -3.25) {$q^{-1}$};
		\node [style=circle, scale=1.5] (34) at (-2.25, -0.25) {};
		\node [style=none] (35) at (-2.25, -0.25) {$q$};
		\node [style=none] (36) at (-0.75, -2) {};
		\node [style=none] (37) at (0.25, -1) {};
		\node [style=none] (38) at (1.25, -0.5) {};
		\node [style=none] (39) at (-1.25, -0.5) {};
		\node [style=none] (40) at (-2.25, -2) {};
		\node [style=none] (41) at (0, 0.25) {$\eta_L$};
		\node [style=none] (42) at (0, -0.5) {$\eta_L$};
		\node [style=none] (43) at (-1.5, -3) {$\epsilon_L$};
		\node [style=none] (44) at (-2, -1.25) {$B$};
		\node [style=none] (45) at (0, -2.5) {$B$};
		\node [style=none] (46) at (1.5, -2.5) {$B$};
	\end{pgfonlayer}
	\begin{pgfonlayer}{edgelayer}
		\draw [in=45, out=-90] (17.center) to (19);
		\draw [in=135, out=-90] (18.center) to (19);
		\draw (21.center) to (9.center);
		\draw (27) to (29.center);
		\draw (26) to (28.center);
		\draw (34) to (21.center);
		\draw (36.center) to (19);
		\draw [bend left=90, looseness=2.00] (17.center) to (37.center);
		\draw (26) to (37.center);
		\draw (38.center) to (27);
		\draw [bend left=90, looseness=0.75] (39.center) to (38.center);
		\draw (39.center) to (18.center);
		\draw [bend right=90, looseness=1.50] (40.center) to (36.center);
		\draw (34) to (40.center);
	\end{pgfonlayer}
\end{tikzpicture}
 = \begin{tikzpicture}
	\begin{pgfonlayer}{nodelayer}
		\node [style=none] (9) at (0.75, 1) {};
		\node [style=none] (10) at (1.75, -3.75) {$A^\dag$};
		\node [style=none] (11) at (-0.25, -3.75) {$A^\dag$};
		\node [style=none] (12) at (1.25, 0.75) {$A^\dag$};
		\node [style=none] (17) at (1.25, -2.25) {};
		\node [style=none] (18) at (0.25, -2.25) {};
		\node [style=circle] (19) at (0.75, -1.75) {};
		\node [style=none] (21) at (0.75, 0.25) {};
		\node [style=circle, scale=2] (26) at (0.25, -2.75) {};
		\node [style=circle, scale=2] (27) at (1.25, -2.75) {};
		\node [style=none] (28) at (0.25, -4) {};
		\node [style=none] (29) at (1.25, -4) {};
		\node [style=none] (32) at (0.25, -2.75) {$q^{-1}$};
		\node [style=none] (33) at (1.25, -2.75) {$q^{-1}$};
		\node [style=circle, scale=1.5] (34) at (0.75, -0.75) {};
		\node [style=none] (35) at (0.75, -0.75) {$q$};
	\end{pgfonlayer}
	\begin{pgfonlayer}{edgelayer}
		\draw [in=-45, out=90] (17.center) to (19);
		\draw [in=-135, out=90] (18.center) to (19);
		\draw (21.center) to (9.center);
		\draw (27) to (29.center);
		\draw (26) to (28.center);
		\draw (27) to (17.center);
		\draw (26) to (18.center);
		\draw (19) to (34);
		\draw (34) to (21.center);
	\end{pgfonlayer}
\end{tikzpicture} \stackrel{(2)}{=} \begin{tikzpicture}
	\begin{pgfonlayer}{nodelayer}
		\node [style=none] (0) at (0.25, -2) {};
		\node [style=none] (1) at (1.25, -2) {};
		\node [style=none] (9) at (0.75, 0.5) {};
		\node [style=none] (10) at (1.75, -3.75) {$A^\dag$};
		\node [style=none] (11) at (-0.25, -3.75) {$A^\dag$};
		\node [style=none] (12) at (1.25, 0.25) {$A^\dag$};
		\node [style=none] (13) at (0, -1) {};
		\node [style=none] (14) at (0, 0) {};
		\node [style=none] (15) at (1.5, 0) {};
		\node [style=none] (16) at (1.5, -1) {};
		\node [style=none] (17) at (1.25, 0) {};
		\node [style=none] (18) at (0.25, 0) {};
		\node [style=circle] (19) at (0.75, -0.5) {};
		\node [style=none] (20) at (0.75, -1) {};
		\node [style=none] (21) at (0.75, 0) {};
		\node [style=none] (22) at (0.25, -1) {};
		\node [style=none] (23) at (1.25, -1) {};
		\node [style=circle, scale=1.5] (24) at (0.25, -2.5) {};
		\node [style=circle, scale=1.5] (25) at (1.25, -2.5) {};
		\node [style=circle, scale=2] (26) at (0.25, -3.25) {};
		\node [style=circle, scale=2] (27) at (1.25, -3.25) {};
		\node [style=none] (28) at (0.25, -4) {};
		\node [style=none] (29) at (1.25, -4) {};
		\node [style=none] (30) at (0.25, -2.5) {$q$};
		\node [style=none] (31) at (1.25, -2.5) {$q$};
		\node [style=none] (32) at (0.25, -3.25) {$q^{-1}$};
		\node [style=none] (33) at (1.25, -3.25) {$q^{-1}$};
	\end{pgfonlayer}
	\begin{pgfonlayer}{edgelayer}
		\draw [in=45, out=-90] (17.center) to (19);
		\draw [in=135, out=-90] (18.center) to (19);
		\draw (20.center) to (19);
		\draw (14.center) to (15.center);
		\draw (15.center) to (16.center);
		\draw (13.center) to (14.center);
		\draw (13.center) to (16.center);
		\draw (21.center) to (9.center);
		\draw [in=90, out=-90] (22.center) to (1.center);
		\draw [in=105, out=-105, looseness=0.75] (23.center) to (0.center);
		\draw (27) to (29.center);
		\draw (25) to (27);
		\draw (1.center) to (25);
		\draw (0.center) to (24);
		\draw (24) to (26);
		\draw (26) to (28.center);
	\end{pgfonlayer}
\end{tikzpicture}
= \begin{tikzpicture}
	\begin{pgfonlayer}{nodelayer}
		\node [style=none] (0) at (0.25, -2) {};
		\node [style=none] (1) at (1.25, -2) {};
		\node [style=none] (9) at (0.75, 0.5) {};
		\node [style=none] (10) at (1.75, -3.75) {$A^\dag$};
		\node [style=none] (11) at (-0.25, -3.75) {$A^\dag$};
		\node [style=none] (12) at (1.25, 0.25) {$A^\dag$};
		\node [style=none] (13) at (0, -1) {};
		\node [style=none] (14) at (0, 0) {};
		\node [style=none] (15) at (1.5, 0) {};
		\node [style=none] (16) at (1.5, -1) {};
		\node [style=none] (17) at (1.25, 0) {};
		\node [style=none] (18) at (0.25, 0) {};
		\node [style=circle] (19) at (0.75, -0.5) {};
		\node [style=none] (20) at (0.75, -1) {};
		\node [style=none] (21) at (0.75, 0) {};
		\node [style=none] (22) at (0.25, -1) {};
		\node [style=none] (23) at (1.25, -1) {};
		\node [style=none] (28) at (0.25, -4) {};
		\node [style=none] (29) at (1.25, -4) {};
	\end{pgfonlayer}
	\begin{pgfonlayer}{edgelayer}
		\draw [in=45, out=-90] (17.center) to (19);
		\draw [in=135, out=-90] (18.center) to (19);
		\draw (20.center) to (19);
		\draw (14.center) to (15.center);
		\draw (15.center) to (16.center);
		\draw (13.center) to (14.center);
		\draw (13.center) to (16.center);
		\draw (21.center) to (9.center);
		\draw [in=90, out=-90] (22.center) to (1.center);
		\draw [in=105, out=-105, looseness=0.75] (23.center) to (0.center);
		\draw (0.center) to (28.center);
		\draw (29.center) to (1.center);
	\end{pgfonlayer}
\end{tikzpicture} \]

Step $(1)$ is true because $(1,q)$ is a morphism of the underlying duals. Step $(2)$ holds because $q$ is a twisted comonoid morphism.

For the converse, assume that $(\tau, \gamma): A \dagmonw A^\dag$ is a twisted right $\dagger$-linear monoid. 
Then, $q$ is a twisted comonoid morphism because: 
\[ \begin{tikzpicture}
	\begin{pgfonlayer}{nodelayer}
		\node [style=none] (0) at (0.25, -3) {};
		\node [style=none] (1) at (1.25, -3) {};
		\node [style=none] (12) at (0.25, 0) {$A^\dag$};
		\node [style=none] (17) at (1.25, -3) {};
		\node [style=none] (18) at (0.25, -3) {};
		\node [style=circle] (19) at (0.75, -2.5) {};
		\node [style=none] (20) at (0.75, -2) {};
		\node [style=none] (28) at (0.25, -4) {};
		\node [style=none] (29) at (1.25, -4) {};
		\node [style=none] (32) at (-0.25, -3.75) {$B$};
		\node [style=none] (33) at (1.75, -3.75) {$B$};
		\node [style=none] (39) at (0.75, -2) {};
		\node [style=none] (40) at (0.75, 0.5) {};
		\node [style=circle, scale=2] (43) at (0.75, -1) {};
		\node [style=none] (44) at (0.75, -1) {$q$};
	\end{pgfonlayer}
	\begin{pgfonlayer}{edgelayer}
		\draw [in=-45, out=90] (17.center) to (19);
		\draw [in=-135, out=90] (18.center) to (19);
		\draw (20.center) to (19);
		\draw (39.center) to (43);
		\draw (43) to (40.center);
		\draw (28.center) to (0.center);
		\draw (29.center) to (1.center);
	\end{pgfonlayer}
\end{tikzpicture} = \begin{tikzpicture}
	\begin{pgfonlayer}{nodelayer}
		\node [style=none] (0) at (0.25, -3) {};
		\node [style=none] (1) at (1.25, -2) {};
		\node [style=none] (9) at (0, 0.5) {$\eta_L$};
		\node [style=none] (12) at (-2.75, 0.25) {$A^\dag$};
		\node [style=none] (17) at (-0.25, -1.25) {};
		\node [style=none] (18) at (-1.25, -1.25) {};
		\node [style=circle] (19) at (-0.75, -1.75) {};
		\node [style=none] (20) at (-0.75, -2.25) {};
		\node [style=none] (28) at (0.25, -4) {};
		\node [style=none] (29) at (1.25, -4) {};
		\node [style=none] (32) at (-0.25, -3.75) {$B$};
		\node [style=none] (33) at (1.75, -3.75) {$B$};
		\node [style=none] (36) at (0.25, -1.25) {};
		\node [style=none] (37) at (-1.25, -0.5) {};
		\node [style=none] (38) at (1.25, -0.5) {};
		\node [style=none] (39) at (-2.25, -2.25) {};
		\node [style=none] (40) at (-2.25, 0.75) {};
		\node [style=none] (41) at (0, -0.75) {$\eta_L$};
		\node [style=none] (42) at (-1.5, -3) {$\epsilon_L$};
		\node [style=circle, scale=2] (43) at (-2.25, -0.75) {};
		\node [style=none] (44) at (-2.25, -0.75) {$q$};
		\node [style=none] (45) at (-2.5, -2) {$B$};
	\end{pgfonlayer}
	\begin{pgfonlayer}{edgelayer}
		\draw [in=45, out=-90] (17.center) to (19);
		\draw [in=135, out=-90] (18.center) to (19);
		\draw (20.center) to (19);
		\draw [bend left=90, looseness=1.50] (17.center) to (36.center);
		\draw (36.center) to (0.center);
		\draw (38.center) to (1.center);
		\draw [bend left=270] (38.center) to (37.center);
		\draw (37.center) to (18.center);
		\draw [bend left=90, looseness=1.25] (20.center) to (39.center);
		\draw (39.center) to (43);
		\draw (43) to (40.center);
		\draw (28.center) to (0.center);
		\draw (29.center) to (1.center);
	\end{pgfonlayer}
\end{tikzpicture} = \begin{tikzpicture}
	\begin{pgfonlayer}{nodelayer}
		\node [style=none] (0) at (0.25, -2) {};
		\node [style=none] (1) at (1.25, -2) {};
		\node [style=none] (9) at (0, 0.5) {$\tau_L$};
		\node [style=none] (10) at (1.75, -2.5) {$A^\dag$};
		\node [style=none] (11) at (-0.25, -2.5) {$A^\dag$};
		\node [style=none] (12) at (-2.75, 0.25) {$A^\dag$};
		\node [style=none] (17) at (-0.25, -1.25) {};
		\node [style=none] (18) at (-1.25, -1.25) {};
		\node [style=circle] (19) at (-0.75, -1.75) {};
		\node [style=none] (20) at (-0.75, -2.25) {};
		\node [style=none] (28) at (0.25, -4) {};
		\node [style=none] (29) at (1.25, -4) {};
		\node [style=circle, scale=2] (30) at (0.25, -3) {};
		\node [style=circle, scale=2] (31) at (1.25, -3) {};
		\node [style=none] (32) at (-0.25, -3.75) {$B$};
		\node [style=none] (33) at (1.75, -3.75) {$B$};
		\node [style=none] (34) at (0.25, -3) {$q$};
		\node [style=none] (35) at (1.25, -3) {$q$};
		\node [style=none] (36) at (0.25, -1.25) {};
		\node [style=none] (37) at (-1.25, -0.5) {};
		\node [style=none] (38) at (1.25, -0.5) {};
		\node [style=none] (39) at (-2.25, -2.25) {};
		\node [style=none] (40) at (-2.25, 0.75) {};
		\node [style=none] (41) at (0, -0.75) {$\tau_L$};
		\node [style=none] (42) at (-1.5, -3) {$\gamma_L$};
	\end{pgfonlayer}
	\begin{pgfonlayer}{edgelayer}
		\draw [in=45, out=-90] (17.center) to (19);
		\draw [in=135, out=-90] (18.center) to (19);
		\draw (20.center) to (19);
		\draw (0.center) to (30);
		\draw (30) to (28.center);
		\draw (31) to (29.center);
		\draw (1.center) to (31);
		\draw [bend left=90, looseness=1.50] (17.center) to (36.center);
		\draw (36.center) to (0.center);
		\draw (38.center) to (1.center);
		\draw [bend left=270] (38.center) to (37.center);
		\draw (37.center) to (18.center);
		\draw [bend left=90, looseness=1.25] (20.center) to (39.center);
		\draw (40.center) to (39.center);
	\end{pgfonlayer}
\end{tikzpicture}
\stackrel{(*)}{=} \begin{tikzpicture}
	\begin{pgfonlayer}{nodelayer}
		\node [style=none] (0) at (0.25, -2) {};
		\node [style=none] (1) at (1.25, -2) {};
		\node [style=none] (9) at (0.75, 0.5) {};
		\node [style=none] (10) at (1.75, -2.5) {$A^\dag$};
		\node [style=none] (11) at (-0.25, -2.5) {$A^\dag$};
		\node [style=none] (12) at (1.25, 0.25) {$A^\dag$};
		\node [style=none] (13) at (0, -1) {};
		\node [style=none] (14) at (0, 0) {};
		\node [style=none] (15) at (1.5, 0) {};
		\node [style=none] (16) at (1.5, -1) {};
		\node [style=none] (17) at (1.25, 0) {};
		\node [style=none] (18) at (0.25, 0) {};
		\node [style=circle] (19) at (0.75, -0.5) {};
		\node [style=none] (20) at (0.75, -1) {};
		\node [style=none] (21) at (0.75, 0) {};
		\node [style=none] (22) at (0.25, -1) {};
		\node [style=none] (23) at (1.25, -1) {};
		\node [style=none] (28) at (0.25, -4) {};
		\node [style=none] (29) at (1.25, -4) {};
		\node [style=circle, scale=2] (30) at (0.25, -3) {};
		\node [style=circle, scale=2] (31) at (1.25, -3) {};
		\node [style=none] (32) at (-0.25, -3.75) {$B$};
		\node [style=none] (33) at (1.75, -3.75) {$B$};
		\node [style=none] (34) at (0.25, -3) {$q$};
		\node [style=none] (35) at (1.25, -3) {$q$};
	\end{pgfonlayer}
	\begin{pgfonlayer}{edgelayer}
		\draw [in=45, out=-90] (17.center) to (19);
		\draw [in=135, out=-90] (18.center) to (19);
		\draw (20.center) to (19);
		\draw (14.center) to (15.center);
		\draw (15.center) to (16.center);
		\draw (13.center) to (14.center);
		\draw (13.center) to (16.center);
		\draw (21.center) to (9.center);
		\draw [in=90, out=-90] (22.center) to (1.center);
		\draw [in=105, out=-105, looseness=0.75] (23.center) to (0.center);
		\draw (0.center) to (30);
		\draw (30) to (28.center);
		\draw (31) to (29.center);
		\draw (1.center) to (31);
	\end{pgfonlayer}
\end{tikzpicture} \]
Step $(*)$ is true because $A$ is a twisted right $\dagger$-linear monoid.
\end{proof}

Next we note that $\dagger$-linear monoids are preserved by $\dagger$-linear functor:
\begin{lemma} 
	\label{Lemma: daglin dagmon}
	$\dagger$-linear functors preserve  (twisted) $\dagger$-linear monoids.
	\end{lemma}
	\begin{proof}
	Suppose $(F_\ox, F_\oa): \X \to \Y$ is a $\dagger$-linear functor, and $A \dagmonw A^\dagger$ 
	is a $dagger$-linear monoid in $\X$. Since, linear functors preserve linear monoids, $F_\ox(A) 
	\linmonw F_\oa(A^\dagger)$ is a linear monoid.
	
	To prove that $F_\ox(A)$ is a right $\dagger$-linear monoid, define the multiplication, unit, 
	comultiplication, counit, action and coaction maps are given as follows:
	\[ \top \to^{\unitmap{1.5}{black}} F_\ox(A) := 
	\]
	It remains to prove that:
	\begin{itemize}
	\item The counit is the dagger of the unit map:
	\[ %
 \]
	\item The comultiplication is the dagger of the multiplication map:
	\[%
 \]
	Observe that in step $(*)$, a twist is introduced if the $\dagger$-linear monoid is twisted.
	\end{itemize}
\end{proof}  

A morphism of (twisted) $\dagger$-linear monoids, follow the same pattern as a morphism 
of $\dagger$-duals:
\begin{definition} 
	\label{defn: morphism of dagger monoids}
	A {\bf morphism of (twisted) $\dagger$-linear monoids} is defined to be one of the following:
\begin{itemize}
\item For (twisted) {\bf right} $\dagger$-linear monoids, it is a pair of maps $(f, f^\dagger)$ which 
is a morphism of underlying linear monoids. 
\item For (twisted) {\bf left} $\dagger$-linear monoids, it is a pair of maps $(f^\dagger, f)$ 
which is a morphism of underlying linear monoids. 
\item For (twisted) {\bf general} $\dagger$-linear monoids, it is a pair $(f,g): A \dagmonw B \to C \dagmonb D$ 
such that $(f,g)$ is a morphism of underlying linear monoids such that one of the 
following equivalent diagrams commute:
\begin{equation} 
	(a)~~~\xymatrix{ 
	C^\dag \ar[r]^{f^\dag} \ar[d]_{q'} & A^\dag \ar[d]^{q} \\ 
	D \ar[r]_g & B } ~~~~~ (or)~~~~~ (b) ~~~ \xymatrix{ 
		A \ar[r]^{f} \ar[d]_{p} & C \ar[d]^{p'} \\ 
		B^\dag \ar[r]_{g^\dag} & D^\dag }
\end{equation}
where $(p,q)$ is the isomorphism between the linear monoid 
$A \linmonw B$ and its dagger $B^\dag \linmonw A^\dag$, and  
$(p',q')$ is an isomorphism between the linear monoid 
$C \linmonb D$ and its dagger $D^\dag  \linmonb C^\dag$ (as per the condition   
satisfied by a general $\dagger$-linear monoid).
\end{itemize}
\end{definition}

In the commuting diagram $(a)$, $f^\dagger: C^\dag \to A^\dag$ is a morphism between the dagger comonoids, 
and $g: D \to B$  is a morphism between the dual comonoids. If the given $\dagger$-linear monoids are 
twisted, then $q$ and $q'$ are twisted comonoid morphisms. Dually, in the commuting diagram $(b)$, 
$f: A \to C$ is a monoid morphism, and $g^\dag: B^\dag \to D^\dag$ is a morphism between the 
conjugate monoids. If the given $\dagger$-linear monoids are twisted, then the morphisms $p$ 
and $p'$ are twisted monoid morphisms. 

\subsection{The endomorphism  linear monoid}
\label{Sec: pants linear monoid}

In Section \ref{Sec: observables}, we discussed the pants $\dagger$-Frobenius algebras in $\dagger$-KCCs.
For an object $A$ in a KCC, the dual structure of $A$ induces a pants algebra a.k.a. an endomorphism monoid 
on the object $A \ox A^*$. 
If the KCC is symmetric, then pants algebras are Frobenius. Moreover, in a $\dagger$-KCC, 
these algebras are $\dagger$-Frobenius. In this section, we present analogue of a pants algebra in $\dagger$-LDCs.

\begin{lemma}
\label{Lemma: pants monoid}
In an LDC, every dual $(\eta_L, \epsilon_L): A \dashvv B : (\eta_R, \epsilon_R)$ induces 
a linear monoid:  $ (\tau, \gamma): (A \oa B) \linmonb (A \ox B)$: 
\[ \unitmap{1.5}{black} := \begin{tikzpicture}
	\begin{pgfonlayer}{nodelayer}
		\node [style=none] (6) at (-1.75, 2.25) {};
		\node [style=none] (7) at (-0.75, 2.25) {};
		\node [style=none] (8) at (-1.75, 1.75) {};
		\node [style=none] (9) at (-0.75, 1.75) {};
		\node [style=none] (16) at (-1.25, 0) {$A \oa B$};
		\node [style=none] (20) at (-0.75, 1.75) {};
		\node [style=oa] (21) at (-1.25, 1) {};
		\node [style=none] (22) at (-1.25, 0.25) {};
		\node [style=none] (23) at (-1.75, 1.75) {};
		\node [style=none] (25) at (-1.25, 3) {$\eta_L$};
	\end{pgfonlayer}
	\begin{pgfonlayer}{edgelayer}
		\draw (6.center) to (8.center);
		\draw (22.center) to (21);
		\draw [in=-90, out=165] (21) to (23.center);
		\draw (7.center) to (20.center);
		\draw [bend left=90, looseness=2.00] (6.center) to (7.center);
		\draw [in=15, out=-90] (20.center) to (21);
	\end{pgfonlayer}
\end{tikzpicture} ~~~~~~~~ \mulmap{1.5}{black} := \begin{tikzpicture}
	\begin{pgfonlayer}{nodelayer}
		\node [style=oa] (0) at (-2, 2) {};
		\node [style=oa] (1) at (0, 2) {};
		\node [style=oa] (3) at (-1, 0) {};
		\node [style=none] (4) at (-2, 2.75) {};
		\node [style=none] (5) at (0, 2.75) {};
		\node [style=none] (6) at (-1, -0.75) {};
		\node [style=none] (8) at (-2.5, -0.5) {$A \oa B$};
		\node [style=none] (9) at (-2, 3.15) {$A \oa B$};
		\node [style=none] (11) at (0, 3.15) {$A \oa B$};
		\node [style=none] (12) at (-1.5, 1.25) {};
		\node [style=none] (13) at (-0.5, 1.25) {};
		\node [style=none] (14) at (-1, 0.5) {$\epsilon_L$};
	\end{pgfonlayer}
	\begin{pgfonlayer}{edgelayer}
		\draw (4.center) to (0);
		\draw (5.center) to (1);
		\draw [bend left=60] (3) to (0);
		\draw [in=0, out=-45] (1) to (3);
		\draw (3) to (6.center);
		\draw [bend right=90, looseness=1.75] (12.center) to (13.center);
		\draw [in=90, out=-165] (1) to (13.center);
		\draw [in=90, out=-30] (0) to (12.center);
	\end{pgfonlayer}
\end{tikzpicture} \] \[ \tau_L := \begin{tikzpicture}
	\begin{pgfonlayer}{nodelayer}
		\node [style=none] (0) at (-0.75, 1.75) {};
		\node [style=none] (1) at (-1.75, 1.75) {};
		\node [style=none] (2) at (-0.75, 2) {};
		\node [style=none] (3) at (-1.75, 2) {};
		\node [style=none] (6) at (-2.25, 2.25) {};
		\node [style=none] (7) at (-0.25, 2.25) {};
		\node [style=none] (8) at (-2.25, 1.75) {};
		\node [style=none] (9) at (-0.25, 1.75) {};
		\node [style=none] (12) at (-1.25, 2.75) {$\eta_R$};
		\node [style=none] (13) at (-0.5, 0.25) {$A \ox B$};
		\node [style=none] (16) at (-2, 0.25) {$A \oa B$};
		\node [style=ox] (17) at (-0.5, 1.25) {};
		\node [style=none] (18) at (-0.5, 0.5) {};
		\node [style=none] (19) at (-0.75, 1.75) {};
		\node [style=none] (20) at (-0.25, 1.75) {};
		\node [style=oa] (21) at (-2, 1.25) {};
		\node [style=none] (22) at (-2, 0.5) {};
		\node [style=none] (23) at (-2.25, 1.75) {};
		\node [style=none] (24) at (-1.75, 1.75) {};
		\node [style=none] (25) at (-1.25, 3.25) {$\eta_L$};
	\end{pgfonlayer}
	\begin{pgfonlayer}{edgelayer}
		\draw (0.center) to (2.center);
		\draw [bend right=90, looseness=1.75] (2.center) to (3.center);
		\draw (3.center) to (1.center);
		\draw (6.center) to (8.center);
		\draw (18.center) to (17);
		\draw [in=-90, out=30, looseness=1.25] (17) to (20.center);
		\draw [in=-90, out=150, looseness=1.25] (17) to (19.center);
		\draw (22.center) to (21);
		\draw [in=-90, out=30, looseness=1.25] (21) to (24.center);
		\draw [in=-90, out=150, looseness=1.25] (21) to (23.center);
		\draw (7.center) to (20.center);
		\draw [bend left=90, looseness=1.25] (6.center) to (7.center);
	\end{pgfonlayer}
\end{tikzpicture}
~~~~~~~~ \gamma_L := \begin{tikzpicture}
	\begin{pgfonlayer}{nodelayer}
		\node [style=none] (0) at (-1.75, 1.75) {};
		\node [style=none] (1) at (-0.75, 1.75) {};
		\node [style=none] (2) at (-1.75, 1.5) {};
		\node [style=none] (3) at (-0.75, 1.5) {};
		\node [style=none] (6) at (-0.25, 1.25) {};
		\node [style=none] (7) at (-2.25, 1.25) {};
		\node [style=none] (8) at (-0.25, 1.75) {};
		\node [style=none] (9) at (-2.25, 1.75) {};
		\node [style=none] (12) at (-1.25, 0.75) {$\epsilon_L$};
		\node [style=none] (13) at (-2, 3.25) {$A \ox B$};
		\node [style=none] (16) at (-0.5, 3.25) {$A \oa B$};
		\node [style=oa] (17) at (-2, 2.25) {};
		\node [style=none] (18) at (-2, 3) {};
		\node [style=none] (19) at (-1.75, 1.75) {};
		\node [style=none] (20) at (-2.25, 1.75) {};
		\node [style=ox] (21) at (-0.5, 2.25) {};
		\node [style=none] (22) at (-0.5, 3) {};
		\node [style=none] (23) at (-0.25, 1.75) {};
		\node [style=none] (24) at (-0.75, 1.75) {};
		\node [style=none] (25) at (-1.25, 0.25) {$\epsilon_R$};
	\end{pgfonlayer}
	\begin{pgfonlayer}{edgelayer}
		\draw (0.center) to (2.center);
		\draw [bend right=90, looseness=1.75] (2.center) to (3.center);
		\draw (3.center) to (1.center);
		\draw (6.center) to (8.center);
		\draw (18.center) to (17);
		\draw [in=90, out=-150, looseness=1.25] (17) to (20.center);
		\draw [in=90, out=-30, looseness=1.25] (17) to (19.center);
		\draw (22.center) to (21);
		\draw [in=90, out=-150, looseness=1.25] (21) to (24.center);
		\draw [in=90, out=-30, looseness=1.25] (21) to (23.center);
		\draw (7.center) to (20.center);
		\draw [bend left=90, looseness=1.25] (6.center) to (7.center);
	\end{pgfonlayer}
\end{tikzpicture} ~~~~~~~~ \tau_R := \begin{tikzpicture}
	\begin{pgfonlayer}{nodelayer}
		\node [style=none] (0) at (-1.75, 1.75) {};
		\node [style=none] (1) at (-0.75, 1.75) {};
		\node [style=none] (2) at (-1.75, 2) {};
		\node [style=none] (3) at (-0.75, 2) {};
		\node [style=none] (6) at (-0.25, 2.25) {};
		\node [style=none] (7) at (-2.25, 2.25) {};
		\node [style=none] (8) at (-0.25, 1.75) {};
		\node [style=none] (9) at (-2.25, 1.75) {};
		\node [style=none] (12) at (-1.25, 2.75) {$\eta_R$};
		\node [style=none] (13) at (-2, 0.25) {$A \ox B$};
		\node [style=none] (16) at (-0.5, 0.25) {$A \oa B$};
		\node [style=oa] (17) at (-2, 1.25) {};
		\node [style=none] (18) at (-2, 0.5) {};
		\node [style=none] (19) at (-1.75, 1.75) {};
		\node [style=none] (20) at (-2.25, 1.75) {};
		\node [style=ox] (21) at (-0.5, 1.25) {};
		\node [style=none] (22) at (-0.5, 0.5) {};
		\node [style=none] (23) at (-0.25, 1.75) {};
		\node [style=none] (24) at (-0.75, 1.75) {};
		\node [style=none] (25) at (-1.25, 3.25) {$\eta_L$};
	\end{pgfonlayer}
	\begin{pgfonlayer}{edgelayer}
		\draw (0.center) to (2.center);
		\draw [bend left=90, looseness=1.75] (2.center) to (3.center);
		\draw (3.center) to (1.center);
		\draw (6.center) to (8.center);
		\draw (18.center) to (17);
		\draw [in=-90, out=150, looseness=1.25] (17) to (20.center);
		\draw [in=-90, out=30, looseness=1.25] (17) to (19.center);
		\draw (22.center) to (21);
		\draw [in=-90, out=150, looseness=1.25] (21) to (24.center);
		\draw [in=-90, out=30, looseness=1.25] (21) to (23.center);
		\draw (7.center) to (20.center);
		\draw [bend right=90, looseness=1.25] (6.center) to (7.center);
	\end{pgfonlayer}
\end{tikzpicture} ~~~~~~~~ \gamma_R := \begin{tikzpicture}
	\begin{pgfonlayer}{nodelayer}
		\node [style=none] (0) at (-0.75, 1.75) {};
		\node [style=none] (1) at (-1.75, 1.75) {};
		\node [style=none] (2) at (-0.75, 1.5) {};
		\node [style=none] (3) at (-1.75, 1.5) {};
		\node [style=none] (6) at (-2.25, 1.25) {};
		\node [style=none] (7) at (-0.25, 1.25) {};
		\node [style=none] (8) at (-2.25, 1.75) {};
		\node [style=none] (9) at (-0.25, 1.75) {};
		\node [style=none] (12) at (-1.25, 0.75) {$\epsilon_R$};
		\node [style=none] (13) at (-0.5, 3.25) {$A \ox B$};
		\node [style=none] (16) at (-2, 3.25) {$A \oa B$};
		\node [style=oa] (17) at (-0.5, 2.25) {};
		\node [style=none] (18) at (-0.5, 3) {};
		\node [style=none] (19) at (-0.75, 1.75) {};
		\node [style=none] (20) at (-0.25, 1.75) {};
		\node [style=ox] (21) at (-2, 2.25) {};
		\node [style=none] (22) at (-2, 3) {};
		\node [style=none] (23) at (-2.25, 1.75) {};
		\node [style=none] (24) at (-1.75, 1.75) {};
		\node [style=none] (25) at (-1.25, 0.25) {$\epsilon_L$};
	\end{pgfonlayer}
	\begin{pgfonlayer}{edgelayer}
		\draw (0.center) to (2.center);
		\draw [bend left=90, looseness=1.75] (2.center) to (3.center);
		\draw (3.center) to (1.center);
		\draw (6.center) to (8.center);
		\draw (18.center) to (17);
		\draw [in=90, out=-30, looseness=1.25] (17) to (20.center);
		\draw [in=90, out=-150, looseness=1.25] (17) to (19.center);
		\draw (22.center) to (21);
		\draw [in=90, out=-30, looseness=1.25] (21) to (24.center);
		\draw [in=90, out=-150, looseness=1.25] (21) to (23.center);
		\draw (7.center) to (20.center);
		\draw [bend right=90, looseness=1.25] (6.center) to (7.center);
	\end{pgfonlayer}
\end{tikzpicture}  \]
\end{lemma}

In a symmetric LDC, every dual automatically induces a pants linear monoid as shown 
in the previous lemma. Every linear monoid embeds into its pants linear monoids:

\begin{lemma} 
	\label{Lemma: pants embedding}
In an LDC, every linear monoid embeds into the pants linear monoids induced 
by its left and right duals. 
\end{lemma}
\begin{proof}
Suppose $A \linmonw B$ is a linear monoid containing a monoid $\monoid{A}$, a left dual  
$(\eta_L, \epsilon_L): A \dashvv B$, and  a right dual $(\eta_R, \epsilon_R): A \dashvv B$. 
The left and the right duals gives two pants linear monoids, 
$ (\tau, \gamma): (A \oa B) \linmonw (A \ox B)$, and 
$(\tau', \gamma'): (B \oa A) \linmonb (B \ox A)$  as shown in Lemma \ref{Lemma: pants monoid}.

In order to prove that $A \linmonw B$ embeds into its 
pants linear monoids, we must show that there exist morphisms of linear monoids, 
$(R,S)$ and $(R', S')$, and that these morphisms are left-invertible:
\[ 
\]

We first prove that $(\eta, \epsilon): A \linmonw B$ embeds into $ (\tau, \gamma): (A \oa B) \linmonw (A \ox B)$.

Define $R: A \to A \oa B$ and $S: A \ox B \to A$ as follows: 
\[ %
  \]

We prove that $R: \monoid{A} \to (A \oa B, \mulmap{1.5}{black}, \unitmap{1.5}{black})$ is a monoid morphism.

$R$ preserves the unit:
\[ 	%
  \] 

$R$ preserves the multiplication:
\[ %
 \]

Finally, $(R, S): ((\eta_L, \epsilon_L): A \dashvv B) \to ((\tau_L, \gamma_L): (A \oa B) \dashvv (A \ox B))$ is a morphism of duals:

\[ %
 \]
Similarly, we can prove that $(S,R): ((\tau_R, \gamma_R): (A \ox B) \dashvv (A \oa B)) \to ((\eta_R, \epsilon_R): B \dashvv A) $ 
is a morphism of duals. 

To prove that $R$ is an embedding, it remains to prove that $R$ is left invertible, that is, there exists $(M,N)$ such that 
$(R,S) (M,N) := (RM, NS) = (1_{A}, 1_B)$:
\[  %
 \]

It is easy to verify that $RM = 1_A$ and $NS = 1_B$. Thus, $(\eta, \epsilon): A \linmonw B$ embeds into $(\tau, \gamma): (A \oa B) \linmonb (A \ox B)$.
The proof that $(\eta, \epsilon): A \linmonw B$ embeds into $(\tau, \gamma): (A \oa B) \linmonb (A \ox B)$ is the mirror reflection 
of the above proof.
\end{proof}

In a symmetric $\dagger$-LDC, 
It turns out that the pants ``monoid'' is a  {\em twisted\/} $\dagger$-linear monoid, see \ref{Defn: twisted}: 
thus, the comultiplication is the dagger of the multiplication but with a twist (using 
the symmetry map). 

In a symmetric $\dagger$-LDC, a twisted $\dagger$-linear monoid embeds into the 
pants $\dagger$-linear monoid. 

\begin{lemma}
\label{Lemma: twisted pants}
	In a symmetric $\dagger$-LDC,  if $(\eta, \epsilon): A \dashvv A^\dag$ is a $\dagger$-dual 
    then the pants monoid $(\tau, \gamma): (A \oa A^\dag) \linmonw (A \ox A^\dag)$ is a 
    twisted $\dagger$-linear 
    monoid.
\end{lemma}
\begin{proof}
    In order to prove that $(\tau, \gamma): (A \oa A^\dag) \linmonw (A \ox A^\dag)$
    is a $\dagger$-linear monoid, we show that there exists a linear monoid morphism 
    $(p,q)$ as shown below such that $p$ and $q$ are isomorphisms, and $pq^\dag = \iota$:
    \[ \begin{matrix}
        \xymatrix{
            A \oa A^\dag  \ar@{-||}[r]^{( \tau,  \gamma)} \ar[d]_{p}  & A \ox A^\dag \ar@{<-}[d]^{q}   \\
            (A \ox A^\dag)^\dag  \ar@{-||}[r]_{(\gamma\dagger,\tau\dagger)} & (A \oa A^\dag)^\dag }
        \end{matrix} \]
    Moreover, to show that the $\dagger$-linear monoid is twisted, we must prove that 
    $p$ is a twisted monoid morphism (or equivalently $q$ is a twisted comonoid morphism).
    A monoid morphism $f: (A, m, u) \to (B, m', u')$ in an LDC is said to be twisted if $(f \ox f) m = c_\ox m' f$.

	The maps $p$ and $q$ as follows:
	\begin{align*}
	p &:= A \oa A^\dag \to^{\iota \oa 1} A^{\dag\dag} 
    \oa A^\dagger \to^{c_\ox} A^\dag \oa A^{\dag \dag} \to^{\lambda_\oa} (A \ox A^\dag)^\dagger \\
	q  &:= (A \oa A^\dag)^\dagger 
    \to^{\lambda_\ox^{-1}} A^\dagger \ox A^{\dag \dag} \to^{1 \ox i^{-1}} 
    A^\dag \ox A \to^{c_\ox} A \ox A^\dag
	\end{align*}

	Diagrammatically, 
	\[ p:= \begin{tikzpicture}
        \begin{pgfonlayer}{nodelayer}
            \node [style=oa] (0) at (-0.25, 2) {};
            \node [style=none] (2) at (0, 0.25) {};
            \node [style=none] (5) at (-1, -0.25) {};
            \node [style=none] (6) at (0, -0.25) {};
            \node [style=none] (7) at (-1, -1.25) {};
            \node [style=none] (8) at (0, -1.25) {};
            \node [style=ox] (9) at (-0.5, -0.75) {};
            \node [style=none] (10) at (-0.5, -0.25) {};
            \node [style=none] (11) at (-1.25, -1.25) {};
            \node [style=none] (12) at (-1.25, -0.25) {};
            \node [style=none] (13) at (0.25, -0.25) {};
            \node [style=none] (14) at (0.25, -1.25) {};
            \node [style=none] (15) at (-0.5, -1.25) {};
            \node [style=none] (16) at (-0.5, -1.75) {};
            \node [style=none] (17) at (-0.5, -2) {$(A \ox A^\dag)^\dagger$};
            \node [style=none] (18) at (-1.5, 0) {$A^\dagger$};
            \node [style=none] (19) at (0.5, 0) {$A^{\dag \dag}$};
            \node [style=none] (20) at (-0.25, 2.75) {};
            \node [style=none] (21) at (-1.25, 1) {$A^\dag$};
            \node [style=none] (22) at (0.25, 1) {$A$};
            \node [style=none] (23) at (-0.25, 3) {$A \oa A^\dag$};
            \node [style=none] (24) at (-1, 0.25) {};
            \node [style=onehalfcircle] (26) at (-0.75, 1.5) {};
            \node [style=none] (25) at (-0.75, 1.5) {$i$};
        \end{pgfonlayer}
        \begin{pgfonlayer}{edgelayer}
            \draw [in=270, out=90] (6.center) to (2.center);
            \draw [in=-165, out=75, looseness=1.25] (7.center) to (9);
            \draw (9) to (10.center);
            \draw [bend left, looseness=1.25] (9) to (8.center);
            \draw (11.center) to (12.center);
            \draw (12.center) to (13.center);
            \draw (13.center) to (14.center);
            \draw (14.center) to (11.center);
            \draw (16.center) to (15.center);
            \draw (0) to (20.center);
            \draw (24.center) to (5.center);
            \draw [in=105, out=-60, looseness=1.50] (0) to (24.center);
            \draw [in=90, out=-90] (26) to (2.center);
            \draw [in=-165, out=75, looseness=1.25] (26) to (0);
        \end{pgfonlayer}
    \end{tikzpicture}   ~~~~~~~~ 	q := \begin{tikzpicture}
        \begin{pgfonlayer}{nodelayer}
            \node [style=ox] (0) at (-0.5, -0.75) {};
            \node [style=none] (1) at (-1, 0.75) {};
            \node [style=circle, scale=2] (2) at (0, 0.75) {};
            \node [style=none] (4) at (0, 0.75) {$i^{-1}$};
            \node [style=none] (5) at (-1, 1.5) {};
            \node [style=none] (6) at (0, 1.5) {};
            \node [style=none] (7) at (-1, 2.5) {};
            \node [style=none] (8) at (0, 2.5) {};
            \node [style=oa] (9) at (-0.5, 2) {};
            \node [style=none] (10) at (-0.5, 1.5) {};
            \node [style=none] (11) at (-1.25, 2.5) {};
            \node [style=none] (12) at (-1.25, 1.5) {};
            \node [style=none] (13) at (0.25, 1.5) {};
            \node [style=none] (14) at (0.25, 2.5) {};
            \node [style=none] (15) at (-0.5, 2.5) {};
            \node [style=none] (16) at (-0.5, 3) {};
            \node [style=none] (17) at (-0.5, 3.5) {$(A \oa A^\dag)^\dagger$};
            \node [style=none] (18) at (-1.75, 1.25) {$A^\dagger$};
            \node [style=none] (19) at (0.75, 1.25) {$A^{\dag \dag}$};
            \node [style=none] (20) at (-0.5, -1.5) {};
            \node [style=none] (23) at (-0.5, -1.75) {$A \ox A^\dag$};
        \end{pgfonlayer}
        \begin{pgfonlayer}{edgelayer}
            \draw (5.center) to (1.center);
            \draw (6.center) to (2);
            \draw [in=30, out=-90, looseness=1.50] (1.center) to (0);
            \draw [in=135, out=-90, looseness=1.50] (2) to (0);
            \draw [in=165, out=-75] (7.center) to (9);
            \draw (9) to (10.center);
            \draw [bend right, looseness=1.25] (9) to (8.center);
            \draw (11.center) to (12.center);
            \draw (12.center) to (13.center);
            \draw (13.center) to (14.center);
            \draw (14.center) to (11.center);
            \draw (16.center) to (15.center);
            \draw (0) to (20.center);
        \end{pgfonlayer}
    \end{tikzpicture}  \]

    Note that $p$ and $q$ are isomorphisms. We must prove that  $pq^\dag = \iota$:
	\begin{align*}
	 pq ^\dagger &= (\iota_A \oa 1) (c_\oa)_{A^{\dag \dag}, A} (\lambda_\oa)_{A^\dag, A^{\dag \dag}} 
     ( (\lambda_\oa)^{-1}_{A^\dag, A^{\dag \dag}} (1 \ox \iota_A^{-1}) (c_\ox)_{A,A^\dag} )^\dag \\
	&=  (\iota_A \oa 1) (c_\oa)_{A^{\dag \dag}, A} (\lambda_\oa)_{A^\dag, A^{\dag \dag}} 
    (c_\ox)_{A,A^\dag}^{\dag}  (1 \ox \iota_A^{-1})^\dag (\lambda_\oa^{-1})_{A^\dag, A^{\dag \dag}}^\dag\\
	&=  (\iota_A \oa 1) (c_\oa)_{A^{\dag \dag}, A} (c_\oa)_{A^\dag,A^{\dag\dag}} 
    (\lambda_\oa)_{A^{\dag \dag}, A^\dag}  (1 \ox \iota_A^{-1})^\dag (\lambda_\oa^{-1})_{A^\dag, A^{\dag \dag}}^\dag \\
    &\stackrel{[\dagger\text{-LDC]-}7}{=} (\iota_A \oa 1) (\lambda_\oa)_{A^{\dag \dag}, A^\dag}  (1 \ox \iota_A^{-1})^\dag (\lambda_\oa^{-1})_{A^\dag, A^{\dag \dag}}^\dag \\
    &\stackrel{Nat.~\lambda_\oa}{=} (\iota_A \oa 1) (1 \oa i_A^{-1\dag}) (\lambda_\oa)( \lambda_\ox^{-1})^\dagger \\
    &\stackrel{\iota_A^{-1 \dag} = \iota_{A^\dag}}{=}(\iota_A \oa \iota_{A^\dag}) (\lambda_\oa)( \lambda_\ox^{-1})^\dagger     \\
	&\stackrel{[\dagger\text{-LDC]-}4}{=} \iota 
	\end{align*}
	
    It remains to show that $(p,q)$ is a twisted linear monoid morphism for which we prove that $q$ is a twisted comonoid morphism 
    and $(p,q)$ is a morphism of the underlying duals. 
	Recall that the linear monoid $(\tau, \gamma): (A \oa A^\dag) \linmonw (A \ox A^\dag)$ is given as follows:  
    Showing that the $q$ preserves the counit:
	\[ 	\mulmap{3}{white} := 
     \] 
    
    The comonoid $(A \ox A^\dag, \comulmap{1.5}{white}, \counitmap{1.5}{white})$ is the dual of the monoid $A \oa A^\dag$ 
    using the $\tau$ and the $\gamma$:
    \[ \comulmap{3}{white} = %
     \]

    The linear monoid $(\gamma \dag, \tau \dag): (A \ox A^\dag)^\dag \linmonb (A \oa A^\dag)^\dag$ is the dagger 
    of $(\tau, \gamma): (A \oa A^\dag) \linmonw (A \ox A^\dag)$.  The following is the proof that 
    $q: ((A \oa A^\dag)^\dag, \comulmap{1.5}{black}, \counitmap{1.5}{black}) \to ((A \ox A^\dag), \comulmap{1.5}{white}, \counitmap{1.5}{white}))$
    is twisted comonoid morphism.

    Showing that $q$ preserves the counit:
		\[ %
 \]

   Showing that $q$ preserves the comultiplication but with a twist:
	\[ %
  \]
It remains to prove that $(p,q)$ is a morphism of duals:
\[ %
 \]

The proof is as follows:

\medskip

$ \tau (p \oa 1) = %
 = (\gamma^\dagger) (1 \oa q) $
\end{proof}
	
It is simultaneously surprising and unsurprising that the pants $\dagger$-linear monoid 
is twisted in $\dagger$-LDCs. A twist is expected because bending the wires of a 
multiplication into a comultiplication changes the original order of the wires 
while daggering does not have this effect (equivalently conjugation reverses the 
tensor product: $\overline{A \ox B} \simeq \overline{B} ~\ox ~\overline{A}$).  
However, in $\dagger$-KCCs, the pants $\dagger$-FAs are not twisted, because the 
$\dagger$-duals are already equipped with a 
twist $(\eta^\dag = c_\ox \epsilon)$ which untwists the pants $\dagger$-FA. 

\begin{lemma} 
	In a symmetric $\dagger$-LDC, every twisted $\dagger$-linear monoid embeds into its pants $\dagger$-linear monoid.
\end{lemma}
\begin{proof} 
	In Lemma \ref{Lemma: pants embedding}, we showed that there exists a left invertible map, $(R,S)$,
	for every linear monoid into its pants monoids. Let $A \dagmonw A^\dag$ be a twisted $\dagger$-linear monoid. 
	\[ 
 \]
	where, $R: A \to A \oa A^\dag$ and $S: A \ox A^\dag \to A$ are defined as follows: 
\[ %
 \]

	To prove that every twisted $\dagger$-linear monoid embeds into its pants $\dagger$-linear monoid, 
	it suffices to prove that $(R,S)$ is a morphism of twisted $\dagger$-linear monoids.

	From Lemma \ref{Lemma: twisted pants}, we have that $(p,q)$ is an isomorphism of linear monoids:
	\[ %
 \]
	where the maps $p$ and $q$ are given diagrammatically as follows:
	\[ p:= %
  \]

	In order to prove that $(R, S)$ is a morphism of twisted $\dagger$-linear monoids, we must prove that $qS = R^\dagger$ or $Rp = S^\dag$ 
	(See Definition \ref{defn: morphism of dagger monoids}):
	\[    qS =  %
 =  R^\dagger 	\]

Step $(1)$ holds because $A$ is a twisted $\dagger$-linear monoid. Step $(2)$ holds because $A$ is a $\dagger$-dual.
\end{proof}



\subsection{Being Frobenius}
Linear monoids generalize Frobenius algebras from monoidal categories to LDCs.
In \ref{compact-mix-functor}, it has been proved that there exists a linear equivalence between compact LDCs and monoidal 
categories. Thus, it is useful to know the precise conditions under which a linear monoid in a 
compact LDC is a Frobenius algebra in the equivalent monoidal category. In this 
section we develop these conditions and extend them to show the correspondence between 
$\dagger$-Frobenius algebras and $\dagger$-linear monoids. 

\begin{lemma}
	\label{AppLemma: equivPresentations} In an LDC, the following conditions (and their 
	`op' symmetries) are equivalent for a self-linear monoid $A \linmonw A'$ with $\alpha: 
	A \to A'$ being an isomorphism. 
	\begin{enumerate}[(a)]
		\item The isomorphism $\alpha$ coincides with the following maps: 
		\begin{equation} 
			\label{eqn: unitary coincidence}
		
		\end{equation}
	\item  The coaction maps given by the comultiplication and the linear duality coincides with the comultiplication:
	\begin{equation}
		\label{eqn: action coincidence}
	%
	
	\end{equation}
	\item The self-dual cup and cap coincide with the cups and the caps of the duals given in the linear monoid:
	\begin{equation}
		\label{eqn: cup coincidence}
		%
	
	\end{equation}
	\end{enumerate}
	\end{lemma}
	\begin{proof}
		For $(a) \Rightarrow (b)$, substitute $\alpha$ in $(b)$  using \ref{eqn: unitary coincidence}. 
		\[  %
 \]
		$(b) \Rightarrow (c)$, and 	$(c) \Rightarrow (a)$ are  straightforward. 
	\end{proof}

In a compact LDC, a linear monoid satisfying one of the conditions in the previous Lemma, precisely 
corresponds to a Frobenius Algebra:
\begin{lemma}
	\label{Lemma: dagmondagFrob}
	In a compact LDC, a self-linear monoid, $A \linmonw A'$ with an isomorphism 
	$\alpha: A \to A'$ precisely corresponds to a Frobenius algebra under the linear equivalence, 
	${\sf Mx}_\downarrow$ if and only if the following equation holds. 
		\begin{equation} 
		\label{eqnn: unitary coincidence}
			%
	\end{equation}
	\end{lemma}
	\begin{proof} 
		Let $(\X, \ox, \oa)$ be a compact-LDC. We know that there exists a linear equivalence, ${\sf Mx}_\downarrow: 
		(X, \ox, \oa) \to (\X, \ox, \ox)$. Let $A \linmonw A'$ be a self-linear monoid in $\X$ 
		with isomorphism $\alpha: A \to A'$. The Frobenius Algebra  $(A, \mulmap{1.5}{white},
		 \unitmap{1.5}{white}, \comulmap{1.5}{black}, \counitmap{1.5}{black})$ in $(\X, \ox, \ox)$ is defined 
		 as follows: the monoid is same as that of the linear monoid, the comonoid is given the by the following data. 
		 \[ %
 \]
		
		It has been given that \ref{eqnn: unitary coincidence} holds. Then, equation \ref{eqn: action coincidence} holds by 
		 Lemma \ref{AppLemma: equivPresentations}. 
		
		 $%
		\stackrel{\ref{eqn: action coincidence}}{=} %
 $
		
		Now assume that there exists a Frobenius Algebra $(A, \mulmap{1.5}{white},
		\unitmap{1.5}{white}, \comulmap{1.5}{black}, \counitmap{1.5}{black})$ in $(\X, \ox, \ox)$. The 
		Frobenius Algebra induces a linear monoid, $A \linmonw A'$ in $(\X, \ox, \oa)$. The linear 
		monoid contains the monoid  $(A, \mulmap{1.5}{white}, \unitmap{1.5}{white})$ and the duals are given as follows:
		\[%
  \]
		It is easy to verify using the Frobenius Law that the left and the right duals are cyclic for the monoid.
		
		For the converse assume that a self-linear monoid $A \linmonw A'$ corresponds to a Frobenius Algebra 
		$(A, \mulmap{1.5}{white}, \unitmap{1.5}{white}, \comulmap{1.5}{black}, \counitmap{1.5}{black})$ 
		under the linear equivalence ${\sf Mx}_\downarrow$. We prove that statement $(iii)$ of Lemma 
		\ref{Lemma: alternate presentation of linear monoids} which is an equivalent form of the equation given in the current Lemma.
		\[  %
 \]
		This proves the converse.
	\end{proof}
\begin{corollary}
	In a unitary category, any $\dagger$-linear monoid $A$ 
	corresponds to a $\dagger$-Frobenius algebra under the equivalence, $\Mx_\downarrow$,
	if and only if the unitary structure map $\varphi_A: A \to A^\dagger$ 
		satisfies the equations: 
		\begin{equation}
            \label{eqnn: unitary coincidence-b}
            %
 
\end{equation}
\end{corollary}
\begin{proof}
		Suppose $A \dagmonw A^\dagger$ is a self-$\dagger$-linear monoid in a unitary category $(\X, \ox, \oa)$ 
		and the given equation holds. There 
		exists a $\dagger$-linear equivalence ${\sf Mx}_\downarrow: (\X, \ox, \oa) \to (\X, \ox, \ox)$ with the dagger 
		on $(\X, \ox, \ox)$ given as follows:
		\[ f^\ddagger = B \to^{\varphi_B} B^\dagger \to^{f^\dagger} A^\dagger \to^{\varphi_A^{-1}} A \]
		From Lemma \ref{Lemma: dagmondagFrob}, we know that the self-linear monoid $A \linmonw A^\dag$ in  $(\X, \ox, \oa)$ corresponds to a 
		Frobenius Algebra $(A, \mulmap{1.5}{white},
		\unitmap{1.5}{white}, \comulmap{1.5}{black}, \counitmap{1.5}{black})$ in $(\X, \ox, \ox)$.
		If $A$ is a $\dagger$-linear monoid, then $(A, \mulmap{1.5}{white},
		\unitmap{1.5}{white}, \comulmap{1.5}{black}, \counitmap{1.5}{black})$
		is a $\dagger$-Frobenius Algebra because:
		\[ 
 \]
		
		For the other way, the linear monoid $A \linmonw A^\dag$ given by the $\dagger$-Frobenius Algebra 
$(A, \mulmap{1.5}{white}, \unitmap{1.5}{white}, \comulmap{1.5}{black}, \counitmap{1.5}{black})$  is a 
$\dagger$-linear monoid because the duals of the linear monoid are $\dagger$-duals:

\vspace{1em}

$ %
 $

\medskip

The dual of the monoid is same as its dagger for $A \linmonw A^\dag$ because 
$(A, \mulmap{1.5}{white}, \unitmap{1.5}{white}, \comulmap{1.5}{black}, \counitmap{1.5}{black})$  is 
$\dagger$-Frobenius. 

For the converse assume that $A \dagmonw A^\dag$ corresponds to the $\dagger$-Frobenius Algebra 
$(A, \mulmap{1.5}{white}, \unitmap{1.5}{white}, \comulmap{1.5}{black}, \counitmap{1.5}{black})$ 
under the equivalence discussed above. Then it follows from the converse of the previous Lemma 
that the unitary structure map $\varphi_A: A \to A^\dag$ satisfies the given equation.
\end{proof}

The equation given in this Corollary should be reminiscent of involutive 
monoids \cite[Theorem 5.28]{HeV19} in $\dagger$-monoidal categories.  

\section{Linear comonoids}
\label{Sec: linear comonoid}
Our motive behind defining $\dagger$-linear monoids is to describe complementary systems in a $\dagger$-LDC setting. 
The bialgebra law is a central ingredient of complementary systems. The directionality of the linear distributors in an LDC  
makes a bialgebraic interaction between two linear monoids impossible.  A linear monoid, however, can interact bialgebraically 
with a linear comonoid.

\begin{definition}
	A {\bf linear comonoid}, $A \lincomonwtik B$, in an LDC consists of a $\ox$-comonoid, 
	$(A, \trianglecomult{0.55}, \trianglecounit{0.55})$, and a left and a right dual, $(\eta_L, \epsilon_L):A \dashvv B$, 
	and $(\eta_R, \epsilon_R): B \dashvv A$, such that:
	\begin{equation}
		\label{eqn: lin comon}
	(a)~~~ \begin{tikzpicture}
		\begin{pgfonlayer}{nodelayer}
			\node [style=none] (19) at (4.5, 1.25) {};
			\node [style=none] (20) at (4.75, 1.5) {$B$};
			\node [style=oa] (26) at (4.4, 3.5) {};
			\node [style=none] (27) at (4.4, 4.25) {};
			\node [style=none] (28) at (3.75, 4.1) {$B \oa B$};
			\node [style=circle] (29) at (4.5, 2.25) {};
		\end{pgfonlayer}
		\begin{pgfonlayer}{edgelayer}
			\draw (26) to (27.center);
			\draw [bend left=60] (26) to (29);
			\draw [bend right=60] (26) to (29);
			\draw (29) to (19.center);
		\end{pgfonlayer}
	\end{tikzpicture} :=\begin{tikzpicture}
		\begin{pgfonlayer}{nodelayer}
			\node [style=circle] (0) at (1.7, 2.75) {};
			\node [style=none] (1) at (1.2, 2) {};
			\node [style=none] (2) at (2.2, 2) {};
			\node [style=none] (3) at (1.7, 3.5) {};
			\node [style=none] (4) at (1.95, 3.25) {$A$};
			\node [style=none] (5) at (0.7, 2) {};
			\node [style=none] (6) at (-0.3, 2) {};
			\node [style=none] (7) at (2.7, 3.5) {};
			\node [style=none] (8) at (2.7, 1.25) {};
			\node [style=none] (9) at (2.95, 1.75) {$B$};
			\node [style=oa] (10) at (0.2, 3) {};
			\node [style=none] (11) at (0.2, 3.75) {};
			\node [style=none] (12) at (-0.5, 3.65) {$B \oa B$};
			\node [style=none] (13) at (2.2, 4.15) {$\eta_L$};
		\end{pgfonlayer}
		\begin{pgfonlayer}{edgelayer}
			\draw [in=-165, out=90, looseness=1.25] (1.center) to (0);
			\draw [in=90, out=-15, looseness=1.25] (0) to (2.center);
			\draw (0) to (3.center);
			\draw [bend left=90, looseness=1.75] (1.center) to (5.center);
			\draw [bend right=90] (6.center) to (2.center);
			\draw [bend right=90, looseness=1.25] (7.center) to (3.center);
			\draw (7.center) to (8.center);
			\draw [in=90, out=-45] (10) to (5.center);
			\draw [in=90, out=-135, looseness=1.25] (10) to (6.center);
			\draw (10) to (11.center);
		\end{pgfonlayer}
	\end{tikzpicture} = \begin{tikzpicture}
			\begin{pgfonlayer}{nodelayer}
				\node [style=circle] (0) at (0.75, 2.75) {};
				\node [style=none] (1) at (1.25, 2) {};
				\node [style=none] (2) at (0.25, 2) {};
				\node [style=none] (3) at (0.75, 3.5) {};
				\node [style=none] (4) at (0.5, 3.25) {$A$};
				\node [style=none] (5) at (1.75, 2) {};
				\node [style=none] (6) at (2.75, 2) {};
				\node [style=none] (7) at (-0.25, 3.5) {};
				\node [style=none] (8) at (-0.25, 1.25) {};
				\node [style=none] (9) at (-0.5, 1.75) {$B$};
				\node [style=oa] (10) at (2.25, 3) {};
				\node [style=none] (11) at (2.25, 3.75) {};
				\node [style=none] (12) at (2.95, 3.65) {$B \oa B$};
				\node [style=none] (13) at (0.25, 4.15) {$\eta_R$};
			\end{pgfonlayer}
			\begin{pgfonlayer}{edgelayer}
				\draw [in=-15, out=90, looseness=1.25] (1.center) to (0);
				\draw [in=90, out=-165, looseness=1.25] (0) to (2.center);
				\draw (0) to (3.center);
				\draw [bend right=90, looseness=1.75] (1.center) to (5.center);
				\draw [bend left=90] (6.center) to (2.center);
				\draw [bend left=90, looseness=1.25] (7.center) to (3.center);
				\draw (7.center) to (8.center);
				\draw [in=90, out=-135] (10) to (5.center);
				\draw [in=90, out=-45, looseness=1.25] (10) to (6.center);
				\draw (10) to (11.center);
			\end{pgfonlayer}
		\end{tikzpicture}		
~~~~~~~~
(b)~~~ \begin{tikzpicture}
	\begin{pgfonlayer}{nodelayer}
		\node [style=none] (19) at (4.5, 1.25) {};
		\node [style=none] (20) at (5, 1.75) {$B$};
		\node [style=none] (24) at (4.5, 4) {};
		\node [style=none] (26) at (5, 4) {$\bot$};
		\node [style=circle] (27) at (4.5, 3.5) {};
	\end{pgfonlayer}
	\begin{pgfonlayer}{edgelayer}
		\draw (19.center) to (27);
		\draw (27) to (24.center);
	\end{pgfonlayer}
\end{tikzpicture} := 
\begin{tikzpicture}
	\begin{pgfonlayer}{nodelayer}
		\node [style=none] (0) at (1.5, 3.5) {};
		\node [style=none] (1) at (1.25, 3.5) {$A$};
		\node [style=none] (2) at (2.5, 3.5) {};
		\node [style=none] (3) at (2.5, 1) {};
		\node [style=none] (4) at (2.75, 1.75) {$B$};
		\node [style=circle, scale=1.5] (5) at (0.75, 1.25) {};
		\node [style=none] (6) at (0.75, 1.25) {$\bot$};
		\node [style=none] (7) at (0.75, 4.25) {};
		\node [style=circle] (8) at (0.75, 2) {};
		\node [style=none] (9) at (1.5, 2.75) {};
		\node [style=none] (13) at (2, 4.25) {$\eta_L$};
		\node [style=circle] (14) at (1.5, 2.75) {};
	\end{pgfonlayer}
	\begin{pgfonlayer}{edgelayer}
		\draw [bend right=90, looseness=1.50] (2.center) to (0.center);
		\draw (2.center) to (3.center);
		\draw (5) to (7.center);
		\draw [dotted, in=-90, out=0, looseness=1.25] (8) to (9.center);
		\draw (0.center) to (14);
	\end{pgfonlayer}
\end{tikzpicture} = \begin{tikzpicture}
	\begin{pgfonlayer}{nodelayer}
		\node [style=none] (0) at (2, 3.5) {};
		\node [style=none] (1) at (2.25, 3.5) {$A$};
		\node [style=none] (2) at (1, 3.5) {};
		\node [style=none] (3) at (1, 1) {};
		\node [style=none] (4) at (0.75, 1.75) {$B$};
		\node [style=circle, scale=1.5] (5) at (2.75, 1.25) {};
		\node [style=none] (6) at (2.75, 1.25) {$\bot$};
		\node [style=none] (7) at (2.75, 4.25) {};
		\node [style=circle] (8) at (2.75, 2) {};
		\node [style=none] (9) at (2, 2.75) {};
		\node [style=none] (13) at (1.5, 4.25) {$\eta_R$};
		\node [style=circle] (14) at (2, 2.75) {};
	\end{pgfonlayer}
	\begin{pgfonlayer}{edgelayer}
		\draw [bend left=90, looseness=1.50] (2.center) to (0.center);
		\draw (2.center) to (3.center);
		\draw (5) to (7.center);
		\draw [dotted, in=-90, out=180, looseness=1.25] (8) to (9.center);
		\draw (0.center) to (14);
	\end{pgfonlayer}
\end{tikzpicture}
\end{equation}
\end{definition}

Note that while a linear monoid has a $\ox$-monoid and a $\oa$-comonoid, a linear comonoid 
has a $\ox$-comonoid and a $\oa$-monoid. The commuting diagram for Equation \ref{eqn: lin comon}(b) is as follows. 
Note that the diagram includes the unitors and linear distributors.
\[ \xymatrix{
    \bot  \ar[r]^{u_\ox^{-1}}   \ar[d]_{u_\ox^{-1}} \ar@{}[dddrrr]|{LHS = RHS}
    & \top \ox \bot \ar[r]^{\eta' \ox 1} 
    &  (A \oa B) \ox \bot \ar[r]^{\partial^r} 
    & B \oa (A \ox \bot) \ar[d]^{1 \oa ( e \ox 1)} \\
    \bot \ox \top  \ar[d]_{1 \ox \eta} & & &  B \oa (\top \ox \bot)  \ar[d]^{1 \oa u_\ox} \\
    \bot \ox (A \oa B) \ar[d]_{\partial^l} & & & B \oa \bot \ar[d]^{u_\oa} \\
    (\bot \ox A) \oa B \ar[r]_{(1 \ox e) \oa 1} 
    & ( \bot \ox \top) \oa B \ar[r]_{u_\ox \oa 1} 
    & \bot \oa B \ar[r]_{u_\oa} 
    & B  } \]

Similar to a linear monoid, however only in an isomix setting, a linear comonoid allows for Frobenius interaction between its $\ox$-monoid and $\oa$-comonoid:
\begin{lemma}
\label{Lemma: Frobenius linear comonoid}
A {\bf linear comonoid}, $A \lincomonwtik B$, in an isomix category is equivalent to the following data:
	\begin{itemize}
	\item a monoid $(B, \mulmap{1.5}{white}: B \oa B \to B, \unitmap{1.5}{white}: \bot \to B)$
	\item a comonoid $(A, \comulmap{1.5}{white}: A \to A \ox A, \counitmap{1.5}{white}: A \to \top) $
	\item actions, $\leftaction{0.5}{white}: B \ox A \to A$, $\rightaction{0.5}{white}: A \ox B \to A$,
	and 
	
	coactions $\leftcoaction{0.55}{white}: B \to A \oa B$, $\rightcoaction{0.5}{white}: B \to B \oa A$,
    \end{itemize}
	such that the following axioms (and their `op' and `co' symmetric forms) hold:  

	\[ (a) ~~
 \]
\end{lemma}
\begin{proof}
	Given that $A \lincomonw B$ is a linear comonoid. 
	Then, $(A, \comulmap{1.5}{white}, \counitmap{1.5}{white})$ is a $\ox$-comonoid and there exists  
	duals $(\eta_L, \epsilon_L): A \dashvv B$, and $(\eta_R, \epsilon_R): B \dashvv A$.

The $\oa$-monoid on $B$ is given as follows: 
\[ %
 \]

The left and the right actions are defined as follows:

\[  %
  \] 

The left and the right coactions are defined as follows:
\[ %
 \] 

The proof that the required axioms hold is as follows. We prove a representative set of axioms. The other axioms are 
horizontal and vertical reflections of the following axioms:

 \begin{enumerate}[(a)]
\item $%
 $
\end{enumerate} 

For the converse, assume the data given in the Lemma.  Define the linear dual $(\eta_L, \epsilon_L): B \dashvv A$ as follows:
	\[ \eta_L = %
  \] 
	Similarly, the cup and the cap of the linear dual $(\eta, \epsilon): A \dashvv B$ is defined as follows:
	\[ \eta_R = %
 	\]
	It is straightforward to check that the snake equations hold, and that $A \lincomonw B$ is a linear comonoid.
\end{proof}

The exponential modalities for LDCs provide an example of linear comonoids. We will prove 
later in Section \ref{Sec: exp modalities} while discussing exponential modalities in LDCs, 
that every dual in an LDC with exponential modalities is also a linear comonoid. 

Next, we explore the correspondence between linear monoids and linear comonoids in compact LDCs. 
A linear comonoid is same as a linear monoid except that a linear monoid has a $\ox$-monoid 
(and a $\oa$-comonoid) while a linear comonoid has a $\ox$-comonoid (and a $\oa$-monoid). 
This implies a correspondence between linear monoids and linear comonoids in a compact setting 
where the tensor and the par products are isomorphic.
In fact, we show that there exists a symmetry called {\bf compact reflection} which is an involution on the category of compact LDCs 
and linear functors: under this symmetry, a linear monoid translate into a linear comonoid and vice 
versa. The compact reflection symmetry is defined as follows: 

Let $\X$ be a compact LDC. Consider the category $\X^\op$ which is a compact LDC with the same tensor and par products as $\X$. 
The maps in $\X^\op$ are as follows:
\[ \infer{B \to^{\overline{f}} A \in \X^{\op}} { A \to^{f} B \in \X}  \] 
Here, the overline does not imply conjugation, and we use the notion to refer to the translation in the opposite category.
The coherence isomorphisms for the tensor and the par products are given as follows: 
\begin{align*}
    & A \ox (B \oa C) \to^{\widehat{\partial}} (A \ox B) \oa C := \overline{\partial}^{-1} \\
	& \bot \to^{\widehat{\m}} \top :=  \overline{\m}^{-1} \\ 
    & A \ox B \to^{\widehat{\mx}} A \oa B := \overline{\mx}^{-1}  
\end{align*}
Associators and unitors are defined similarly. 

\begin{lemma}
	\label{Lemma: monoid-comonoid-translation}
	Every linear monoid in a compact LDC, $\X$, corresponds to a linear comonoid in $\X^\op$ 
    under the compact reflection of $\X$, and vice versa. 
\end{lemma}
\begin{proof}
	Let $A \linmonw B$ be a linear monoid in the compact LDC, $\X$. We know that any 
	linear monoid, $A \linmonw B$, is precisely a monoid, $(A, m, u)$ with duals, 
	$(\eta_L, \epsilon_L): A \dashvv B$, and $(\eta_R, \epsilon_R): B \dashvv A$, and the 
	duals are cyclic on the monoid. The monoid $(A, m, u)$ gives a $\ox$-comonoid 
	$(A, \Delta := \overline{m} , e := \overline{u} )$ 
	in $\X^{\op}$ under compact reflection. The dual  $(\eta_L, \epsilon_L): A \dashvv B$ 
	in the linear monoid gives the following maps in $\X^\op$:
	\[   \infer{A \oa B \to^{\overline{\eta_L}} \top \in \X^{\op}} { \top \to^{\eta_L} A \oa B \in \X} 
	~~~~~~~~~~~~ 
	\infer{  \bot \to^{\overline{\epsilon_L}} B \ox A \in \X^{\op}} {  B \ox A \to^{\epsilon_L} \bot \in \X}   \]
	
	Similarly, the dual $(\eta_R, \epsilon_R): B \dashvv A$ gives maps $\overline \eta_R : 
	B \oa A \to \top$, and $\overline \epsilon_R : \bot \to A \ox B$ in $\X^\op$. Under this 
	translation of $A \linmonw B$, we get a linear comonoid, $A \lincomonw B$, in $\X^\op$ 
	with the $\ox$-comonoid $(A, \Delta, e)$, and duals $(\tau_R, \gamma_R): B \dashvv A$, 
	$(\tau_L, \gamma_L):  A \dashvv B$ where, 
	\begin{align*}
    \top \to^{\tau_L} A \oa B &:=   \top \to^{\overline{\m}} \bot \to^{\overline{\epsilon_R}} A \ox B \to^{\overline{\mx}^{-1}} A \oa B \\
    B \ox A \to^{\gamma_L} \bot &:=  B \ox A \to^{\overline{\mx}^{-1}} B \oa A \to^{\overline{\eta}_R} \top \to^{\overline{\m}} \bot 
	\end{align*}
	$\tau_R : \top \to B \oa A$, and $\gamma_R: A \ox B \to \bot$ are defined similarly.  Compact reflection of the resulting linear comonoid gives back the original linear monoid  in $(\X^{\op})^\op = \X$. 	
\end{proof}

It follows from the previous Lemma that the results pertaining to linear monoids in a compact LDC 
also applies for linear comonoids under compact reflection.

We move on to the discussion of morphisms of linear comonoids:

\begin{definition}
A {\bf morphism} of linear comonoids, $(f,g) : (A \lincomonw B) \to (A' \lincomonw B')$, consists of a pair of maps,
$f: A \to A'$ and $g: B' \to B$, such that $f$ is a comonoid morphism, and $(f,g)$ and $(g,f)$ are morphisms of 
the left and the right duals respectively.  
\end{definition}

\begin{definition}
A {\bf $\dagger$-linear comonoid} in a $\dagger$-LDC is $A \dagcomonwtik A^\dagger$ is a linear comonoid, 
$A \lincomonw A^\dagger$ such that $(\eta_L, \epsilon_L): A \dashvv A^\dagger$, and $(\eta_R, \epsilon_R): 
A^\dagger \dashvv A$ are $\dagger$-duals, and:
\begin{equation}
    \label{Eqn: leftdagcomon}
	\begin{tikzpicture}
		\begin{pgfonlayer}{nodelayer}
			\node [style=none] (19) at (4.5, 1.25) {};
			\node [style=none] (20) at (4.85, 1.75) {$A^\dagger$};
			\node [style=none] (21) at (4.5, 3.5) {};
			\node [style=none] (24) at (4.5, 4) {};
			\node [style=none] (25) at (4.5, 3.5) {};
			\node [style=none] (26) at (5, 4) {$\bot$};
			\node [style=circle] (27) at (4.5, 3.5) {};
		\end{pgfonlayer}
		\begin{pgfonlayer}{edgelayer}
			\draw (19.center) to (21.center);
			\draw (24.center) to (25.center);
		\end{pgfonlayer}
	\end{tikzpicture} := \begin{tikzpicture}
		\begin{pgfonlayer}{nodelayer}
			\node [style=none] (0) at (1.5, 3.5) {};
			\node [style=none] (1) at (1.25, 3.5) {$A$};
			\node [style=none] (2) at (2.5, 3.5) {};
			\node [style=none] (3) at (2.5, 1) {};
			\node [style=none] (4) at (2.75, 1.75) {$A^\dagger$};
			\node [style=circle, scale=1.5] (5) at (0.75, 1.25) {};
			\node [style=none] (6) at (0.75, 1.25) {$\bot$};
			\node [style=none] (7) at (0.75, 4.25) {};
			\node [style=circle] (8) at (0.75, 2) {};
			\node [style=none] (9) at (1.5, 2.75) {};
			\node [style=circle] (10) at (1.5, 2.75) {};
		\end{pgfonlayer}
		\begin{pgfonlayer}{edgelayer}
			\draw [bend right=90, looseness=1.25] (2.center) to (0.center);
			\draw (2.center) to (3.center);
			\draw (5) to (7.center);
			\draw [dotted, in=-90, out=0, looseness=1.25] (8) to (9.center);
			\draw (0.center) to (10);
		\end{pgfonlayer}
	\end{tikzpicture}  = \begin{tikzpicture}
	\begin{pgfonlayer}{nodelayer}
		\node [style=none] (1) at (0, 2.5) {};
		\node [style=none] (2) at (-0.75, 2.5) {};
		\node [style=none] (3) at (-0.75, 1.25) {};
		\node [style=none] (4) at (0.75, 1.25) {};
		\node [style=none] (5) at (0.75, 2.5) {};
		\node [style=none] (6) at (0, 0) {};
		\node [style=none] (7) at (0, 1.25) {};
		\node [style=none] (8) at (0.5, 0.25) {$A^\dagger$};
		\node [style=none] (11) at (0, 1.75) {};
		\node [style=none] (12) at (0.5, 1.25) {};
		\node [style=circle, scale=1.5] (13) at (1.5, 0.25) {};
		\node [style=none] (14) at (1.5, 3.75) {};
		\node [style=circle] (15) at (1.5, 3) {};
		\node [style=none] (16) at (0.5, 2.5) {};
		\node [style=none] (17) at (1.5, 0.25) {$\bot$};
		\node [style=circle] (18) at (0, 1.75) {};
	\end{pgfonlayer}
	\begin{pgfonlayer}{edgelayer}
		\draw (2.center) to (5.center);
		\draw (5.center) to (4.center);
		\draw (4.center) to (3.center);
		\draw (3.center) to (2.center);
		\draw (6.center) to (7.center);
		\draw [dotted, bend right=15] (11.center) to (12.center);
		\draw (14.center) to (13);
		\draw [dotted, in=180, out=90, looseness=1.25] (16.center) to (15);
		\draw (1.center) to (18);
	\end{pgfonlayer}
\end{tikzpicture}
~~~~~~~~  \begin{tikzpicture}
	\begin{pgfonlayer}{nodelayer}
		\node [style=none] (19) at (4.5, 1.25) {};
		\node [style=none] (20) at (4.75, 1.5) {$A^\dag$};
		\node [style=oa] (26) at (4.4, 3.5) {};
		\node [style=none] (27) at (4.4, 4.25) {};
		\node [style=none] (28) at (3.5, 4.1) {$A^\dag \oa A^\dag$};
		\node [style=circle] (29) at (4.5, 2.25) {};
	\end{pgfonlayer}
	\begin{pgfonlayer}{edgelayer}
		\draw (26) to (27.center);
		\draw [bend left=60] (26) to (29);
		\draw [bend right=60] (26) to (29);
		\draw (29) to (19.center);
	\end{pgfonlayer}
\end{tikzpicture} :=\begin{tikzpicture}
	\begin{pgfonlayer}{nodelayer}
		\node [style=circle] (0) at (1.7, 2.75) {};
		\node [style=none] (1) at (1.2, 2) {};
		\node [style=none] (2) at (2.2, 2) {};
		\node [style=none] (3) at (1.7, 3.5) {};
		\node [style=none] (4) at (1.95, 3.25) {$A$};
		\node [style=none] (5) at (0.7, 2) {};
		\node [style=none] (6) at (-0.3, 2) {};
		\node [style=none] (7) at (2.7, 3.5) {};
		\node [style=none] (8) at (2.7, 1.25) {};
		\node [style=none] (9) at (2.95, 1.75) {$A^\dag$};
		\node [style=oa] (10) at (0.2, 3) {};
		\node [style=none] (11) at (0.2, 3.75) {};
		\node [style=none] (12) at (-0.6, 3.65) {$A^\dag \oa A^\dag$};
	\end{pgfonlayer}
	\begin{pgfonlayer}{edgelayer}
		\draw [in=-165, out=90, looseness=1.25] (1.center) to (0);
		\draw [in=90, out=-15, looseness=1.25] (0) to (2.center);
		\draw (0) to (3.center);
		\draw [bend left=90, looseness=1.75] (1.center) to (5.center);
		\draw [bend right=90] (6.center) to (2.center);
		\draw [bend right=90, looseness=1.25] (7.center) to (3.center);
		\draw (7.center) to (8.center);
		\draw [in=90, out=-45] (10) to (5.center);
		\draw [in=90, out=-135, looseness=1.25] (10) to (6.center);
		\draw (10) to (11.center);
	\end{pgfonlayer}
\end{tikzpicture} = \begin{tikzpicture}
		\begin{pgfonlayer}{nodelayer}
			\node [style=none] (0) at (0.75, -0.25) {};
			\node [style=none] (1) at (0.25, -1) {$A^\dagger$};
			\node [style=none] (2) at (0, 1.5) {};
			\node [style=none] (4) at (-1, -0.25) {};
			\node [style=none] (5) at (-1, 1.5) {};
			\node [style=none] (6) at (1, 1.5) {};
			\node [style=none] (7) at (1, -0.25) {};
			\node [style=none] (8) at (0.5, 1.5) {};
			\node [style=none] (9) at (-0.5, 1.5) {};
			\node [style=none] (10) at (0, -0.25) {};
			\node [style=none] (11) at (0, -1) {};
			\node [style=oa] (12) at (0, 2.25) {};
			\node [style=none] (13) at (0, 3) {};
			\node [style=none] (14) at (1, 2.75) {$A^\dagger \oa A^\dagger$};
			\node [style=none] (18) at (-0.75, -0.25) {};
			\node [style=circle] (19) at (0, 0.75) {};
		\end{pgfonlayer}
		\begin{pgfonlayer}{edgelayer}
			\draw (5.center) to (4.center);
			\draw (4.center) to (7.center);
			\draw (6.center) to (7.center);
			\draw (6.center) to (5.center);
			\draw (11.center) to (10.center);
			\draw [bend left] (9.center) to (12);
			\draw [bend left] (12) to (8.center);
			\draw (12) to (13.center);
			\draw [bend right] (19) to (18.center);
			\draw [bend left] (19) to (0.center);
			\draw (19) to (2.center);
		\end{pgfonlayer}
	\end{tikzpicture}	
\end{equation}
A {\bf $\dagger$-self-linear comonoid} consists of an 
isomorphism $\alpha: A \to A^\dagger$ such that $\alpha \alpha^{-1 \dagger} = \iota$. 
\end{definition}

Next, we discuss the compact reflection of $\dagger$-linear comonoids in compact $\dagger$-LDCs. 
If $\X$ is a compact $\dagger$-LDC, then $\X^\op$ given by compact reflection of $\X$ has the same dagger functor. In this case, 
the laxors and involution natural isomorphisms are defined similar to the other coherence isomorphisms:
\begin{align*}
    & A^\dagger \ox B^\dagger \to^{\widehat{\lambda}_\ox} (A \oa B)^\dagger := \overline{\lambda_\ox}^{-1} \\
    & \top \to^{\widehat{\lambda}_\top} \bot^{\dagger} := \overline{\lambda_\top}^{-1} \\
    & A \to^{\widehat{\iota}} A^{\dagger \dagger} := \overline{\iota}^{-1} 
\end{align*}
If $\X$ is a unitary category, then $\X^\op$ is also a unitary category under compact reflection 
with the unitary structure map for an object $A$, given as 
$\widehat{\varphi}_A := \overline{\varphi_A}^{-1}$.

In Lemma \ref{Lemma: monoid-comonoid-translation}, we proved that the compact reflection of 
a linear monoid in a compact LDC gives a linear comonoid and vice versa. In the following lemma 
we extend the result to $\dagger$-linear monoids on compact $\dagger$-LDCs: 
\begin{lemma}
	If $\X$ is a compact $\dagger$-LDC, then every dagger linear monoid corresponds to a 
	dagger linear comonoid under the 
	compact reflection. 
\end{lemma}
\begin{proof}
Suppose $\X$ is a compact $\dagger$-LDC, and $A \dagmonw A^\dagger$ is a $\dagger$-linear monoid.
Then, we have a $\ox$-monoid $(A, m, u)$ with a right $\dagger$-dual, 
$A \dashvv A^\dagger$, and a left $\dagger$-dual $A^\dagger \dashvv A$ satisfying the equations for 
a $\dagger$-linear monoid. 
To show that under compact reflection, a $\dagger$-linear monoid translates to a $\dagger$-linear comonoid, 
it suffices to prove that the $\dagger$-duals in $\X$ translates to $\dagger$-duals in $\X^\op$. 
The dual, $(\tau_L, \gamma_L): A \dashvv A^\dagger$ defined as in Lemma 
\ref{Lemma: monoid-comonoid-translation} is a $\dagger$-dual because: 
\medskip

{ \centering $\begin{tikzpicture}
	\begin{pgfonlayer}{nodelayer}
		\node [style=none] (36) at (4.25, 1.75) {};
		\node [style=none] (37) at (2.75, 1.75) {};
		\node [style=none] (38) at (4.25, 0) {};
		\node [style=none] (39) at (2.75, 0) {};
		\node [style=circle, scale=2] (40) at (2.75, 1) {};
		\node [style=none] (41) at (2.75, 1) {$\widehat{\iota}$};
		\node [style=none] (42) at (4.5, 0.25) {$A^\dagger$};
		\node [style=none] (43) at (2.25, 1.5) {$A$};
		\node [style=none] (44) at (2.25, 0.25) {$A^{\dagger \dagger}$};
		\node [style=none, scale=1.5] (45) at (3.5, 3) {$\tau_L$};
	\end{pgfonlayer}
	\begin{pgfonlayer}{edgelayer}
		\draw (38.center) to (36.center);
		\draw (39.center) to (40);
		\draw (37.center) to (40);
		\draw [bend right=90, looseness=2.00] (36.center) to (37.center);
	\end{pgfonlayer}
\end{tikzpicture} = \begin{tikzpicture}
	\begin{pgfonlayer}{nodelayer}
		\node [style=none] (46) at (8.5, 0) {};
		\node [style=none] (48) at (8.5, 1.75) {};
		\node [style=none] (49) at (7, 1.75) {};
		\node [style=none, scale=1.5] (50) at (7.75, 3) {$\overline{\epsilon_R}$};
		\node [style=none] (51) at (9, 0.25) {$A^\dagger$};
		\node [style=none] (52) at (6.75, 1.5) {$A$};
		\node [style=circle, scale=2] (53) at (7, 1) {};
		\node [style=none] (54) at (7, 1) {$\widehat{\iota}$};
		\node [style=none] (55) at (7, 0) {};
		\node [style=none] (56) at (6.5, 0.25) {$A^{\dagger \dagger}$};
	\end{pgfonlayer}
	\begin{pgfonlayer}{edgelayer}
		\draw (48.center) to (46.center);
		\draw (55.center) to (53);
		\draw [bend left=270, looseness=1.75] (48.center) to (49.center);
		\draw (53) to (49.center);
	\end{pgfonlayer}
\end{tikzpicture} = \begin{tikzpicture}
	\begin{pgfonlayer}{nodelayer}
		\node [style=none] (57) at (12.25, 0) {};
		\node [style=none] (59) at (12.25, 1.75) {};
		\node [style=none] (60) at (10.75, 1.75) {};
		\node [style=none, scale=1.5] (61) at (11.5, 3) {$\overline{\epsilon_R}$};
		\node [style=none] (62) at (12.75, 0.25) {$A^\dagger$};
		\node [style=none] (63) at (10.25, 1.75) {$A$};
		\node [style=circle, scale=2] (64) at (10.75, 1) {};
		\node [style=none] (65) at (10.75, 1) {$\overline{\iota}^{-1}$};
		\node [style=none] (66) at (10.75, 0) {};
		\node [style=none] (67) at (10.25, 0.25) {$A^{\dagger \dagger}$};
	\end{pgfonlayer}
	\begin{pgfonlayer}{edgelayer}
		\draw (59.center) to (57.center);
		\draw (66.center) to (64);
		\draw [bend left=270, looseness=1.75] (59.center) to (60.center);
		\draw (64) to (60.center);
	\end{pgfonlayer}
\end{tikzpicture}
 \stackrel{(*)}{=}  \begin{tikzpicture}
	\begin{pgfonlayer}{nodelayer}
		\node [style=none] (68) at (17, 0.25) {$A^\dagger$};
		\node [style=none] (69) at (16.25, 2.75) {$A$};
		\node [style=none] (70) at (15.25, 0) {};
		\node [style=none] (71) at (14.75, 0.25) {$A^{\dagger \dagger}$};
		\node [style=none] (72) at (16.5, 3) {};
		\node [style=none, scale=1.5] (73) at (15.75, 2) {$\overline{\eta_R}$};
		\node [style=none] (74) at (16.75, 1.75) {};
		\node [style=none] (75) at (14.75, 1.75) {};
		\node [style=none] (76) at (14.75, 3) {};
		\node [style=none] (77) at (16.75, 3) {};
		\node [style=none] (78) at (16.5, 1.75) {};
		\node [style=none] (79) at (15.25, 1.75) {};
		\node [style=none] (80) at (16.5, 0) {};
		\node [style=none] (81) at (15, 3) {};
	\end{pgfonlayer}
	\begin{pgfonlayer}{edgelayer}
		\draw (74.center) to (75.center);
		\draw (75.center) to (76.center);
		\draw (77.center) to (74.center);
		\draw (80.center) to (78.center);
		\draw (70.center) to (79.center);
		\draw (77.center) to (76.center);
		\draw [bend left=90, looseness=1.50] (72.center) to (81.center);
	\end{pgfonlayer}
\end{tikzpicture}  =  \begin{tikzpicture}
	\begin{pgfonlayer}{nodelayer}
		\node [style=none] (82) at (6.75, 0.25) {$A^\dagger$};
		\node [style=none] (83) at (6, 2.75) {$A$};
		\node [style=none] (84) at (5, 0) {};
		\node [style=none] (85) at (4.5, 0.25) {$A^{\dagger \dagger}$};
		\node [style=none] (86) at (6.25, 3) {};
		\node [style=none, scale=1.3] (87) at (5.5, 2) {$\gamma_L$};
		\node [style=none] (88) at (6.5, 1.75) {};
		\node [style=none] (89) at (4.5, 1.75) {};
		\node [style=none] (90) at (4.5, 3) {};
		\node [style=none] (91) at (6.5, 3) {};
		\node [style=none] (92) at (6.25, 1.75) {};
		\node [style=none] (93) at (5, 1.75) {};
		\node [style=none] (94) at (6.25, 0) {};
		\node [style=none] (95) at (4.75, 3) {};
	\end{pgfonlayer}
	\begin{pgfonlayer}{edgelayer}
		\draw (88.center) to (89.center);
		\draw (89.center) to (90.center);
		\draw (91.center) to (88.center);
		\draw (94.center) to (92.center);
		\draw (84.center) to (93.center);
		\draw (91.center) to (90.center);
		\draw [bend left=90, looseness=1.50] (86.center) to (95.center);
	\end{pgfonlayer}
\end{tikzpicture}$ \par} 

\medskip

 The step $(*)$ is true because $(\eta_R, \epsilon_R): A^\dag \dagdual A$ is a $\dagger$-dual.
\end{proof}

\begin{definition}
A {\bf morphism of $\dagger$-linear comonoids} is a pair $(f, f^\dagger)$ such that $(f, f^\dagger)$ is a 
morphism of the underlying linear comonoids. 
\end{definition}


\section{Linear bialgebras}
\label{Sec: linear bialgebra}
Linear bialgebras provide the basis for defining complementary systems in isomix categories.
A linear bialgebra has a $\ox$-bialgebra and a $\oa$-bialgebra which are given by the bialgebraic interaction of a linear monoid 
and a linear comonoid.  Thus, in a linear bialgebra two distinct dualities, that of the linear monoid and that of the linear comonoid, 
are at play.   

All the results concerning bialgebras are necessarily set in symmetric LDCs, and we shall assume 
that the linear monoids and comonoids are symmetric.

\begin{definition}
	\label{Defn: linear bialg}
	A {\bf linear bialgebra}, $\frac{(a,b)}{(a',b')}\!\!:\! A \linbialgwtik B$, in an LDC consists of:
	\begin{enumerate}[(a)]
		\item a linear monoid, $(a,b)\!\!:\! A \linmonw B$, and
		\item a linear comonoid, $(a',b')\!\!:\! A \lincomonwtritik B$,
	\end{enumerate}
	such that:
	\begin{enumerate}[(i)]
		\item $(A, \mulmap{1.5}{white}, \unitmap{1.5}{white}, \trianglecomult{0.65} , \trianglecounit{0.65})$ is a $\ox$-bialgebra, and 
		\item $(B, \trianglemult{0.65}, \triangleunit{0.65}, \comulmap{1.5}{white} , \counitmap{1.5}{white})$ is a $\oa$-bialgebra. 
	\end{enumerate}
\end{definition}	

A linear bialgebra is {\bf commutative} if the $\oa$-monoid and $\ox$-monoid are commutative.
A {\bf self-linear bialgebra} is a linear bialgebra, in which there is an isomorphism $A \to^{\alpha} B$ 
(so essentially the algebra is on one object).

\begin{definition}
A {\bf morphism} of linear bialgebras $(f,g):  \frac{(a,b)}{(a',b')}\!\!:\! A \linbialgwtik B \to 
\frac{(c,d)}{(c',d')}\!\!:\! C \linbialgwtik D$ is a morphism the linear monoids 
$((f,g): (A \linmonw B) \to (C \linmonw B))$, and a morphism of the linear comonoids 
$((f,g): (A \lincomonwtritik B) \to (C \lincomonwtritik D))$.
\end{definition}

Isomorphisms transport linear bialgebras: 
\begin{lemma} In an LDC,  if $\frac{(a,b)}{(a',b')}\!\!:\! A \linbialgwtik B$ is a linear 
	bialgebra and $f: A' \to A$ and $g: B \to B'$ are isomorphisms, then 
	$A' \linbialgwtik B'$ 	is a linear bialgebra with the linear monoid 
	$(a(f^{-1} \oa g), (g^{-1} \ox f)b): A' \linmonw B'$, and the
	linear comonoid $(a'(f^{-1} \oa g),(g^{-1} \ox f)b'): A' \lincomonwtritik B'$.
\end{lemma}

A $\dagger$-linear bialgebra is defined as follows:

\begin{definition}
\begin{enumerate}[(i)]
	\item A {\bf $\dagger$-linear bialgebra}, $\frac{(a,b)}{(a',b')} \!:\!A \dagbialgwtik A^\dagger$, is a linear bialgebra with a $\dagger$-linear monoid, 
and a $\dagger$-linear comonoid. 
   \item A {\bf $\dagger$-self-linear bialgebra} is $\dagger$-linear bialgebra which is also a self-linear bialgebra  such 
that the isomorphism, $\alpha: A \to A^\dagger$, satisfies $\alpha \alpha^{-1 \dagger} = \iota$. 
\end{enumerate}
\end{definition}
Note that $A$ is a weak pre-unitary object, which if in the core, is a pre-unitary object. 


The next chapter develops complementary systems and measurements in MUCs 
using the structures developed in the current chapter.

%% file: chapter-conclusion.tex

\chapter{Conclusion and future work}
\label{Chap: conclusion}
\section{Conclusion}

This thesis is a product of our journey towards developing a suitable categorical semantics for quantum 
processes without the constraint of dimensionality. Since its inception, CQM has piqued the interest 
of researchers by its versatile and elegant approach towards studying quantum foundations and 
quantum processes. While CQM was built on compact closed categories - 
categorical semantics of linear logic with the pair of multiplicatives combined into one and the pair of 
additives combined into one - which forces the Hilbert space model to be finite dimensional, there have 
been efforts to address the dimensionality constraint of CQM. 

While one approach was to construct a category which is still compact closed but can suitably describe quantum processes of 
arbitrary dimensions \cite{GG17, HeR18}, the other approach was to construct suitable algebraic structures 
in $\dagger$-SMCs \cite{CoH16, AbH12}.  Attempts have also been made by adding properties 
to the underlying monoidal category so that it accommodates objects of arbitrary dimensional \cite{Heu08,Vic10}. 
In our journey, we stepped out of the well-trodden path of $\dagger$-monoidal categories in CQM, went back to the 
fundamentals, and as category theorists asked, {\em ``Can we consider the semantics of $\dagger$-linear logic that indeed is not degenerate?"} 
This turned our attention towards linearly distributive categories and $*$-autonomous categories (LDCs without negation).

Proceeding to define a $\dagger$ functor for LDCs, and a unitary structure for $\dagger$-isomix categories, we arrived at 
the notion of mixed unitary categories (MUCs) which are $\dagger$-isomix categories with a chosen 
unitary core. As revealed by the schematic diagram of MUC in Figure \ref{Fig: MUC}, the $\dagger$-isomix 
category is a larger space in which the {(smaller)} unitary category resides. For any MUC, its unitary 
core is equivalent to a $\dagger$-monoidal category, and in the presence of unitary duals, it is equivalent to a 
$\dagger$-compact closed category. Thus, one gets a neat description of the traditional CQM framework 
in terms of a larger and a more general framework. We have discussed quite a few examples of MUCs 
as we developed the structure. Significant among these examples are, the category of finiteness relations $\FRel$ 
which has a faithful functor into $\Rel$, and the category of finiteness matrices over of commutative rig $R$, $\FMat(R)$. 
In fact, the category of finite complex matrices $\Mat(\C)$ is isomorphic to the core of $\FMat(\C)$. 

Further in our journey, we tested our framework in its applicability to quantum mechanics by generalizing 
the algebraic structures fundamental to CQM - completely positive maps modeling quantum processes, 
special commutative $\dagger$-Frobenius algebras modeling quantum observables, 
and bialgebras and Hopf algebras modeling complementary observables - from $\dagger$-SMCs to $\dagger$-isomix categories. 
By generalizing the $\CP^\infty$ construction from $\dagger$-SMCs to MUCs, we noticed that, in MUCs, 
interestingly, a completely positive map always factors through the unitary core in its Kraus decomposition form.  
Moreover, in order to characterize the $\CP^\infty$ construction, it suffices that one can discard information 
for unitary objects. Hence, the notion of environment structure, in other words, discarding information, 
is a requirement only within the unitary core. 

Moving on, we turned our focus on {\em linear monoids} \cite{BCST96} in LDCs which are a general version of 
Frobenius algebras: in linear monoids, the monoid and the comonoid pair generally occur on different objects, 
the monoid is on the $\ox$-product, the comonoid is on the $\oa$-product, and the monoid and the comonoid 
objects are dual to one another. We introduced $\dagger$-linear monoids, which in a unitary category are 
equivalent to $\dagger$-Frobenius algebras when the unitary 
structure map satisfies a certain condition. Linear comonoids which are same as linear monoids except with  
the monoid on the $\oa$-product and the comonoid on the $\ox$-product, are a surprising consequence 
of our attempt to understand bialgebraic interaction of linear monoids. A bialgebraic interaction 
between two linear monoids implies a bialgebraic interaction between a $\ox$-monoid and 
a $\oa$-comonoid, which is not supported in an LDC. This led to the idea of linear comonoids. Indeed, 
linear comonoids are not mere algebraic constructs defined for convenience, but they are significant, since 
exponential modalities of linear logic provide a source of examples for linear comonoids. 

A linear monoid and linear comonoid may interact beautifully to produce a $\ox$-bialgebra and a $\oa$-bialgebra, 
we refer to this interaction as a linear bialgebra. At this stage, one can perceive the significance of 
keeping the tensor products (the multiplicatives of linear logic) distinct 
in the framework.  Complementary observables are intimately connected to the notion of measurement. 
In a MUC, we showed that a measurement takes place in two steps: first an arbitrary (non-unitary) type is 
compacted into the unitary core, while inside the unitary core, one can use the traditional machinery of CQM 
to perform a measurement. Interestingly, compaction produces a binary idempotent with certain properties. 
We called this a coring $\dagger$-binary idempotent. Conversely, a coring $\dagger$-binary idempotent 
induces a compaction. 

Finally, turning our attention to complementary systems, we considered such a system in 
$\dagger$-isomix categories as a self-linear bialgebra satisfying certain equations.
These equations implied that the pair of bialgebras are indeed Hopf with a particular anitpode.
We noted that when the $\dagger$-isomix category has free exponential modalities (the exponential operator of  
linear logic which provides for non-linear types with infinite duplication and discarding), every complementary system 
induces a linear bialgebra on the free exponential modalities (in the larger space). Indeed, in such a setting, 
every complementary system arises by splitting a coring $\dag$-binary idempotent 
on the induced $\dag$-linear bialgebra. This is perhaps, the most interesting result in this thesis in relation to 
quantum mechanics. 

Bohr's principle of complementarity \cite{Gri18} states that, due to the wave and particle nature of matter, 
physical properties occur in complementary pairs.  Our result connecting complementarity with exponential modalities  
displays a complementary system as a compaction of a $\dagger$-linear bialgebra in which the 
two dual structures (one pertaining to linear monoid and the other one of the linear comonoid) 
occur separately providing an interesting perspective on Bohr's principle.   

\section{Future work}

The current journey has many interesting further directions, a few of which are discussed below:

\subsection{Models of physical systems}

In this thesis, we have discussed multiple examples of MUCs such as 
 $\FRel$, $\FMat(R)$ and Chu spaces from a mathematical viewpoint.  
The role of these categories in the study of quantum foundations and the
other areas of quantum theory is yet to be explored. In particular, $\FMat(R)$ seems 
to be an interesting candidate to study quantum mechanics since its canonical unitary core is isomorphic to $\FHilb$ 
which is a well-studied model in CQM. Moreover, $\FMat(R)$ is a model of full linear logic and 
comes with the free exponential modalities.  

In \cite{Vic08}, Vicary uses the $!$ exponential modality and a differential map \cite{BCS06} 
to categorify the creation and annihilation operators for a Fock space in $\dagger$-monoidal categories 
with $\dagger$-biproducts. He assumes that the $\dagger$ commutes with the $!$ modality. 
Recall that, in a linear setting, we saw that applying $\dagger$ to $!$ modality gives the $?$ modality. 
Vicary applies the machinery to a category ${\sf Inner}$ with complex inner product spaces of countable 
dimensions and well-defined linear maps to study the state space of Quantum Harmonic Oscillators (QHO). 
However,  the claim that the adjunction between the cofree and the forgetful functor producing the 
$!$-modality is well-defined is left as a conjecture. It has been mentioned that the difficulty lies in proving 
that $\eta$ map of the adjunction is well-defined. It would be interesting to revisit Vicary's ideas 
on categorifying QHO \cite{Vic08}  in a linearly distributive setting, in particular, in $\FMat(R)$ with distinct 
$!$ and $?$  modalities. 

In \cite{Vic08}, Vicary points out that with exponential modality one can model the state space of a QHO, 
however, in order to achieve a categorical description of the dynamics of the system, one needs the 
ability to express differential equations categorically. Differential categories \cite{BCS06, BCL19} which 
are additive symmetric monoidal categories with a coalgebra comodality and a differential combinator 
seems to be a natural candidate satisfying the requirement. \cite{Lem20} emphasizes the relevance of differential 
categories to quantum foundations. In fact, $\FRel$ and $\FMat(R)$ are indeed differential 
categories due to the presence of the free exponential modalities.  This arises a question if can one obtain a complete 
categorical description of quantum harmonic oscillator in differential categories?

Our thesis proved that in a $\dagger$-isomix category, in the presence of free exponential modalities, 
every complementary system inside the canonical unitary core arises from the linear bialgebra induced on the 
free exponential modalities. Free exponential modalities imply the presence of a differential combinator \cite{Fio07, BCL19}
 and provide a categorical description for the state space of a QHO \cite{BPS94, Vic08}. This leads one to wonder 
if there is any interesting connection between quantum harmonic oscillators and complementary observables in physics. 

\subsection{Joyal's Bicompletion construction}

An interesting source of examples for MUC is Joyal's \cite{Joy95, Joy95b} bicompletion 
procedure on monoidal categories. Starting with a $\dagger$-monoidal category, or unitary category, $\C$, 
one can form a MUC $i : \C \to \Lambda(\C)$ by simply bicompleting 
(by adding arbitrary limits and colimits) to the $\dagger$-monoidal category. 
Furthermore, the bicompletion,  $\Lambda(\C)$, is a (non-compact) $\dagger$-isomix category which, 
when the starting point, $\C$, is $\dagger$-compact 
closed, is a $\dagger$-isomix $*$-autonomous category. Bicompleting a monoidal category, $\C$, 
causes its tensor to split into two linearly related tensors products $\ox$ and $\oa$ and 
induces a cofree functor $i: \C \to \Lambda(\C)$ on the category of bicomplete categories and 
bicontinuous functors. The free bicomplete category generated by a single object is a $*$-autonomous 
category \cite[Corollary - Theorem 3]{Joy95b}.

\subsection{Clifford algebras in linear settings}

Clifford algebras\footnote{Note that the term Clifford algebras and geometric algebra are used 
interchangeably in the physics literature. However, Clifford algebras are free geometric algebras 
which satisfy a universal property.}  \cite{Hes12, LuS09, DoL03} are regarded as a universal 
language for physics due to their intimate connection to geometry.  They neatly geometrify 
algebra and algebraize geometry.  Clifford algebras have been applied in many fields of physics including 
quantum gravity \cite{FrK04, DoL07}, quantum field theory \cite{VaR19} and quantum electrodynamics \cite{Bay04}. 

The Clifford algebra of space time captures the geometry of special relativity \cite{Hes12}. 
 In his discussion on categorification of a quantum harmonic oscillator \cite[Sec. Discussion]{Vic08}, Vicary  
notes that within the current CQM formalism,

{\em``\ldots an elegant categorical description of the other ingredients 
of the Schr\"odinger equation, such as Planck’s constant $h$ and the imaginary unit $i$, is far from apparent."}

Clifford algebra addresses the above concern by providing a geometric interpretation for the 
imaginary constant $i$ as rotations in space time. The gamma matrices which are $4 \times 4$ anti-commuting 
unitary matrices solving the Dirac equation  in quantum field theory gives a matrix representation for 
a Clifford algebra, $Cl_{1,3}(\R)$ \cite{VaR19}. The Pauli $I, X, Y, Z$ matrices which are quite significant for 
quantum mechanics and quantum computing provides a representation for the Clifford algebra 
$Cl_{0,3}(\R)$ \cite{DoL03}.  Schr\"odinger's equation can be represented as an element in the Clifford algebra $Cl_{0,1}(\R)$ \cite{HiC10}. 

Clifford algebras enable an interesting possibility of understanding and formalizing quantum computation 
via geometry rather than the other standard non-intuitive methods such as unitary matrices formalism or the 
circuit language method borrowed from classical computing.  This approach is quite different from the 
existing approaches and can provide a fresh  perspective on quantum computation. 
Efforts \cite{SCH98,  Vla99, Vla01, AeC07, Cza07, Lin2021, Lin2021b} to apply 
Clifford algebras to quantum computing in a non-categorical setting 
emphasize the conceptual clarity and computational 
advantages provided by Clifford algebras. In \cite{SCH98} the authors show that the Clifford  
algebra description of quantum computation operations has a direct correlation to NMR spectroscopy, hence
can be implemented in NMR quantum computing without further translation. 
The more recent works \cite{Lin2021, Lin2021b} interesting uses string diagrams of monoidal categories to describe Clifford operations. 

In the future, we would like to formulate abstract Clifford algebras in linear settings. It would also be interesting 
to adapt the string diagrams of CQM to the Clifford algebras.  The elements of a Clifford algebra 
always anti commute. This suggests a connection between these algebras to complementary 
measurements. Can Clifford algebras be used to provide a fresh perspective on quantum computation 
by taking an approach quite different from classical computation? Drawing inspiration from the ZX calculus \cite{CoD11}, can one build a
 universal Clifford calculus that would bridge multiple areas of quantum research?  
Can such a Clifford calculus provide a neat generalization to the ZX calculus to arbitrary dimensional systems?
It would be quite  interesting and highly useful to devise a diagrammatic calculus for the
complicated equations  areas in quantum chemistry and other branches of quantum mechanics. In ZX calculus, the bialgebraic interaction between 
the complementary $Z$ and the $X$ observables is exploited to construct gates for quantum computation. 
The Pauli $X, Y$ and $Z$ matrices along with the identity matrix is a representation for a Clifford algebra. 
This arises a question if there is any interesting connection between bialgebras and Clifford algebras?  
These questions are yet to be explored.